%% file: main.tex
\documentclass[10pt, oneside]{memoir}

\usepackage{mukherjee} 
\include{layout_thesis.tex}

\usepackage{breakcites}
\usepackage[utf8]{inputenc}
\usepackage{emptypage}
\usepackage[final]{pdfpages}

\author{\vspace{.7cm}\\ter verkrijging van de graad van doctor aan de\\
	Technische Universiteit Eindhoven, op gezag van de\\
	rector magnificus, prof.dr.ir. F.P.T. Baaijens, voor een\\
	commissie aangewezen door het College voor\\
 	Promoties in het openbaar te verdedigen\\
	op dinsdag 28 augustus 2018 om 16.00 uur\\
	\vspace{1cm}\\
door\\
\vspace{1cm}\\Debankur Mukherjee\\
\vspace{.7cm}\\
geboren te Hooghly, India}
\title{\textbf{Scalable Load Balancing Algorithms\\
 in Networked Systems}\\
\vspace{1.5cm}
{\LARGE\scshape proefschrift}}
\date{}

\usepackage{shorttoc,titletoc,libertine}
\newcommand\partialtocname{Contents}
\newcommand\ToCrule{\noindent\rule[5pt]{\textwidth}{1pt}}
\newcommand\ToCtitle{{\large\bfseries\partialtocname}\vskip2pt\ToCrule}
\makeatletter
\newcommand\Mprintcontents{%
  \ToCtitle
  \ttl@printlist[chapters]{toc}{}{1}{}\par\nobreak
  \ToCrule
  \vskip10pt}
\makeatother

\usepackage{diagbox}

\usepackage[breaklinks=true]{hyperref}
\hypersetup{colorlinks=true,linkcolor=black,citecolor=darkgray,pdfdisplaydoctitle=true,pdfpagemode=UseOutlines,bookmarksnumbered=true,bookmarksopen=true}

\begin{document} 

\frontmatter

\thispagestyle{empty}
\begin{vplace}[0.7]
\begin{center}
\Huge Scalable Load Balancing Algorithms \\in Networked Systems
\end{center}
\end{vplace}
\newpage

\thispagestyle{empty}%
\noindent This work was financially supported by The Netherlands Organization for Scientific Research (NWO) through the TOP-GO grant 613.001.012 and Gravitation Networks grant 024.002.003.\\\\
\includegraphics[width = 0.55\linewidth]{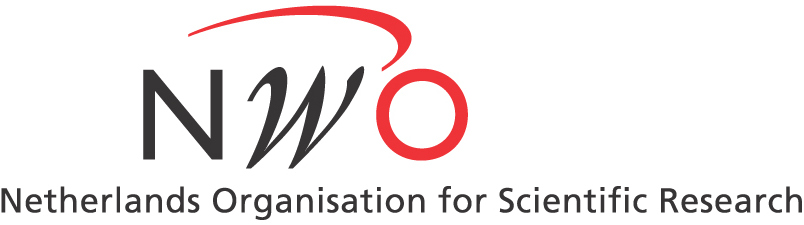}\hfill
\includegraphics[width = 0.22\linewidth]{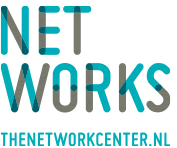}
\vfill

\noindent \copyright~Debankur Mukherjee, 2018
\vspace{1em}

\noindent Scalable Load Balancing Algorithms in Networked Systems
\vspace{1em}

\noindent A catalogue record is available from the Eindhoven University of Technology Library\\
\noindent ISBN: 978-90-386-4558-2
\vspace{1em}


\noindent Printed by Gildeprint Drukkerijen, Enschede

\maketitle
\aliaspagestyle{title}{empty}

\newpage
\noindent Dit proefschrift is goedgekeurd door de promotoren en de samenstelling\\van de promotiecommissie is als volgt:\\ \thispagestyle{empty}
\begin{flushleft}
\begin{tabular}{ l l }
voorzitter:   & prof.dr. M.T. de Berg\\
$1^\text{e}$ promotor: &prof.dr.ir. S.C. Borst\\
$2^\text{e}$ promotor: &prof.dr. J.S.H. van Leeuwaarden\\
leden: 
&prof.dr. S. Bhulai (Vrije Universiteit)\\
&prof.dr. R.W. van der Hofstad\\
&prof.dr. A.P. Zwart\\
&dr. A. Dieker (Columbia University)\\
&dr. D.A. Goldberg (Cornell University)
\end{tabular}
\end{flushleft}
\vfill
Het onderzoek dat in dit proefschrift wordt beschreven is uitgevoerd in\\ overeenstemming met de TU/e Gedragscode Wetenschapsbeoefening.

\chapter*[Acknowledgments]{Acknowledgments}
\input{acknowledgment}

\cleardoublepage

\begin{KeepFromToc}
\tableofcontents
\end{KeepFromToc}

\mainmatter

\chapter[Overview of Results]{Overview of Results}
\label{chap:introduction}
\emph{Based on:}
\begin{itemize}
\item[\cite{BBLM18}] {Van der Boor}, M., Borst, S.~C., {Van Leeuwaarden}, J. S.~H., and Mukherjee,
  D. (2018).
\newblock {Scalable load balancing in networked systems: A survey of recent advances}.
\newblock {\em	arXiv:1806.05444}. Extended abstract appeared in {\em Proc. ICM '18}.
\end{itemize}
\vfill
\startcontents[chapters]
\Mprintcontents
\include{intro} 
\stopcontents[chapters]
\cleardoublepage

\chapter[Universality of JSQ($d$) Policies]{Universality of JSQ($d$) Policies} 
\label{chap:univjsqd}
\emph{Based on:}
\begin{itemize}
\item[\cite{MBLW16-3}] Mukherjee, D., Borst, S.~C., {Van Leeuwaarden}, J. S.~H., and Whiting, P.~A.
  (2016b).
\newblock {Universality of power-of-d load balancing in many-server systems}.
\newblock {\em Stochastic Systems}, to appear. {\em arXiv:1612.00723}.
\end{itemize}
\vfill
\startcontents[chapters]
\Mprintcontents
 \include{univjsqd}
\stopcontents[chapters]
\cleardoublepage

\chapter[Universality of JIQ($d$) Policies]{Universality of JIQ($d$) Policies} 
\label{chap:jiq}
\emph{Based on:}
\begin{itemize}
\item[\cite{MBLW16-1}] Mukherjee, D., Borst, S.~C., {Van Leeuwaarden}, J. S.~H., and Whiting, P.~A.
  (2016a).
\newblock {Universality of load balancing schemes on the diffusion scale}.
\newblock {\em J. Appl. Probab.}, 53(4).
\end{itemize}
\vfill
\startcontents[chapters]
\Mprintcontents
 \include{jiq}
\stopcontents[chapters]
\cleardoublepage

\chapter[Steady-state Analysis of JSQ in the Diffusion Regime]{Steady-state Analysis of JSQ in the Diffusion Regime} 
\label{chap:jsqdiffusion}
\emph{Based on:}
\begin{itemize}
\item[\cite{BM18}] Banerjee, S. and Mukherjee, D. (2018).
\newblock {Join-the-shortest queue diffusion limit in Halfin-Whitt regime: Tail
  asymptotics and scaling of extrema}.
\newblock {\em Annals of Applied Probability}, minor revision.  {\em arXiv:1803.03306}.
\end{itemize}
\vfill
\startcontents[chapters]
\Mprintcontents
 \include{jsqdiffusion}
\stopcontents[chapters]
\cleardoublepage

\chapter[Asymptotic Optimality for Infinite-Server Dynamics]{Asymptotic Optimality for Infinite-Server Dynamics} 
\label{chap:asympjsqd}
\emph{Based on:}
\begin{itemize}
\item[\cite{MBLW16-4}] Mukherjee, D., Borst, S.~C., {Van Leeuwaarden}, J. S.~H., and Whiting, P.~A.
  (2016a).
\newblock {Asymptotic optimality of power-of-d load balancing in large-scale
  systems}.
\newblock {\em Mathematics of Operations Research}, under revision. {\em arXiv:1612.00722}.
\end{itemize}
\vfill
\startcontents[chapters]
\Mprintcontents
 \include{asympjsqd}
\stopcontents[chapters]
\cleardoublepage


\chapter[Optimal Service Elasticity]{Optimal Service Elasticity in Large-Scale Distributed Systems} 
\label{chap:energy1}
\emph{Based on:}
\begin{itemize}
\item[\cite{MDBL17}] Mukherjee, D., Dhara, S., Borst, S.~C., and {Van Leeuwaarden}, J.~S. (2017).
\newblock {Optimal service elasticity in large-scale distributed systems}.
\newblock {\em Proc. ACM Meas. Anal. Comput. Syst.}, 1(1):25.
\end{itemize}
\vfill
\startcontents[chapters]
\Mprintcontents
 \include{energy1}
\stopcontents[chapters]
\cleardoublepage

\chapter[Optimal Service Elasticity for Infinite Buffers]{Optimal Service Elasticity for Infinite Buffers: Large-Scale Asymptotics of a Non-monotone System} 
\label{chap:energy2}
\emph{Based on:}
\begin{itemize}
\item[\cite{MS18}] 
Mukherjee, D. and Stolyar, A. (2018).
\newblock {Join-Idle-Queue with service elasticity: Large-scale asymptotics of
  a non-monotone system}.
\newblock {\em Stochastic Systems}, minor revision. {\em arXiv:1803.07689}.
\end{itemize}
\vfill
\startcontents[chapters]
\Mprintcontents
 \include{energy2}
\stopcontents[chapters]
\cleardoublepage

\chapter[Load Balancing Topologies: JSQ on Graphs]{Load Balancing Topologies: \\JSQ on Graphs} 
\label{chap:networkjsq}
\emph{Based on:}
\begin{itemize}
\item[\cite{MBL17}] Mukherjee, D., Borst, S.~C., and {Van Leeuwaarden}, J. S.~H. (2018).
\newblock {Asymptotically optimal load balancing topologies}.
\newblock {\em Proc. ACM Meas. Anal. Comput. Syst.}, 2(1):1--29.
\end{itemize}
\vfill
\startcontents[chapters]
\Mprintcontents
 \include{networkjsq}
\stopcontents[chapters]
\cleardoublepage

\chapter[Load Balancing Topologies: JSQ(d) on Graphs]{Load Balancing Topologies:\\ JSQ(d) on Graphs} 
\label{chap:networkjsqd}
\emph{Based on:}
\begin{itemize}
\item[\cite{BMW17}] Budhiraja, A., Mukherjee, D., and Wu, R. (2017).
\newblock {Supermarket model on graphs}.
\newblock {\em Annals of Applied Probability}, minor revision. {\em arXiv:1712.07607}.
\end{itemize}
\vfill
\startcontents[chapters]
\Mprintcontents
 \include{networkjsqd}
\stopcontents[chapters]
\cleardoublepage


\cleardoublepage

\chapter*[Summary]{Summary}
\addcontentsline{toc}{chapter}{Summary}
\vfill
\input{TT_Summary}

\chapter*[About the author]{About the author}
\addcontentsline{toc}{chapter}{About the author}
\input{TT_cv}

\cleardoublepage

\end{document}

%% file: layout_thesis.tex
\setstocksize{240mm}{170mm}
\settrimmedsize{240mm}{170mm}{*}

\settrims{0mm}{0mm}

\settypeblocksize{168.96mm}{112mm}{*}

\setlrmargins{22mm}{*}{*}

\setulmargins{27mm}{*}{*}

\setheaderspaces{*}{7mm}{*}

\checkandfixthelayout
\fixpdflayout

\setsecnumdepth{subsection}

\numberwithin{equation}{chapter}
\counterwithin{figure}{chapter}

\maxtocdepth{section}

\setpnumwidth{3em}
\setrmarg{4em}

\usepackage{graphicx}
\chapterstyle{bianchi}

\aliaspagestyle{title}{empty}

\aliaspagestyle{part}{empty}

\usepackage{microtype}

%% file: acknowledgment.tex
Science is a collaborative effort, said twice Nobel prize winning physicist John Bardeen. From my experience in the last four years, I can only agree to this statement. This thesis has been possible because of the involvements of many devoted individuals, and I would like to take this opportunity to convey my gratitude to them.

First of all, I am immensely grateful to my supervisors Sem Borst and Johan van Leeuwaarden, who showed me the rudiments of independent research, and how to see worthy potential coming out of simple ideas. Honestly, I could not ask for a better supervision. Sem, your dedication, thoroughness, and deep insight will always keep inspiring me. I received your valuable and constructive suggestions whenever I so needed. You not only guided me academically, but also shaped me as a human being. Johan, thank you so much for your endless enthusiasm and constant encouragement to push my boundaries. Your frank and straightforward advice helped me get a clear perspective on things. 

I am particularly grateful to Jim Dai, David Goldberg, Maria Vlasiou, and Bert Zwart for providing countless valuable suggestions and honest feedback to prepare me for the job interview, and for helping me out during the difficult path of my career selection.

I would like to thank all my collaborators who shared their enthusiasm and brilliant ideas with me. A special thanks to Phil Whiting for our numerous insightful discussions in many joint projects. It was a unique opportunity for me to work with Alexander Stolyar. Sasha, I learned a lot from your passion for mathematical rigor and elegant style of approaching a problem. Visiting UNC Chapel Hill was a remarkable experience. 
Thanks Sayan Banerjee, Shankar Bhamidi, and Amarjit Budhiraja for making my visit memorable. Sayan, thanks for teaching me stochastic analysis in such depth. It was a wonderful experience spending so many hours thinking with you in front of the blackboard. A special thanks goes to Ruoyu Wu, whose tremendous effort made Chapter 9 of this thesis possible. It was a pleasure working with Subhabrata Sen. Thanks, Subhabrata, for sharing your enthusiasm and deep insight in our joint project. 

An internship with Mark de Berg and Bart Jansen enabled me to go beyond my research boundaries, and venture into a completely new realm. Thanks Bart and Mark for mentoring me during the NETWORKS internship. It was a very pleasant experience for me.
 
I would like to express my gratitude to Sandjai Bhulai, Ton Dieker, David Goldberg, Remco van der Hofstad, and Bert Zwart for agreeing to serve on my doctoral committee and for providing helpful comments on my thesis.

I am indebted to Antar Bandyopadhyay, Sreela Gangopadhyay, and Arup Pal from Indian Statistical Institute, whose lectures immensely influenced my decision for pursuing research in mathematics, and probability in particular. Also, I would not have decided to come to the Netherlands without the guidance from Krishanu Maulik. Thanks, Krishanu for sharing the opportunity and informing us about this place.

In the last four years I have visited many departments, and it is hard to find such a vibrant and academically rich environment as the Stochastics section at TU/e. The presence of EURANDOM definitely makes this department special. The seminars and workshops here maintain a constant flow of frontier researchers from all around the world throughout the year. I would like to express my very great appreciation to Onno Boxma and Remco van der Hofstad. It has been an honor to work under your leadership. I wish to acknowledge the generous financial support and wide academic exposure provided by the NETWORKS grant throughout my PhD.

My sincere thanks to all my colleagues in the department for maintaining a vibrant academic environment. Life in the department would have been much more difficult without such helpful, kind and efficient secretaries. Chantal Reemers and Petra Rozema-Hoekerd, thanks to both of you for taking care of so many things. A special mention goes to Enrico, my first office mate, for being a good friend, for the awesome road trip in Canada, and also for being our local guide in Italy. I would also like to thank Alessandro, Britt, Fabio, Gianmarco, Jori, Jaron, and Thomas for being so welcoming when I first arrived here, and for being so helpful ever since. Thanks to my current office mates Ellen, Kay, Marta, and Richard for creating the perfect work environment.

Moving from India to the Netherlands was a huge environmental and social change for me. I am eternally indebted to Soma Ray, who never let me feel that I am far away from home. My deepest gratitude goes to Souvik Dhara. Thanks Souvik, for sharing so many fantastic ideas, knowledge, and enthusiasm. Our everyday discussions taught me more than I could ever give you credit for here. I hope we will continue enriching each other like this. Finally, I would like to thank my parents, whose love and guidance are with me in whatever I pursue. I have always found them by my side in all my ups and downs, successes and failures.

%% file: intro.tex
\section{Introduction}

In this monograph we pursue scalable load balancing algorithms
(LBAs) that achieve excellent delay performance in large-scale
systems and yet only involve low implementation overhead.
LBAs play a critical role in distributing service requests or tasks
(e.g.~compute jobs, data base look-ups, file transfers) among servers
or distributed resources in parallel-processing systems.
The analysis and design of LBAs has attracted strong attention in recent
years, mainly spurred by scalability challenges arising
in cloud networks and data centers with massive numbers of servers.

LBAs can be broadly categorized as static, dynamic, or some intermediate
blend, depending on the amount of feedback or state information
(e.g.~congestion levels) that is used in allocating tasks.
The use of state information naturally allows dynamic policies
to achieve better delay performance, but also involves higher
implementation complexity and a substantial communication burden.
The latter issue is particularly pertinent in cloud networks and data
centers with immense numbers of servers handling a huge influx
of service requests.
In order to capture the large-scale context, we examine scalability
properties through the prism of asymptotic scalings where the system
size grows large, and identify LBAs which strike an optimal balance
between delay performance and implementation overhead in that regime.

The most basic load balancing scenario consists of $N$~identical
parallel servers and a dispatcher where tasks arrive that must
immediately be forwarded to one of the servers.
Tasks are assumed to have unit-mean exponentially distributed service
requirements, and the service discipline at each server is supposed
to be oblivious to the actual service requirements, i.e., the service time only gets revealed once a server begins processing the task.
In this canonical setup, the celebrated Join-the-Shortest-Queue (JSQ)
policy has several strong stochastic optimality properties.
In particular, the JSQ policy achieves the minimum mean overall
delay among all non-anticipating policies that do not have any
advance knowledge of the service requirements \cite{EVW80,Winston77}.
In order to implement the JSQ policy however, a dispatcher requires
instantaneous knowledge of all the queue lengths, which may involve
a prohibitive communication burden with a large number of servers~$N$.

This poor scalability has motivated consideration of JSQ($d$) policies,
where an incoming task is assigned to a server with the shortest queue
among $d \geq 2$ servers selected uniformly at random.
Note that this involves an exchange of $2 d$ messages per task,
irrespective of the number of servers~$N$.
Results in Mitzenmacher~\cite{Mitzenmacher01} and Vvedenskaya
{\em et al.}~\cite{VDK96} indicate that even sampling as few as $d = 2$
servers yields significant performance enhancements over purely
random assignment ($d = 1$) as $N$ grows large, which is commonly
referred to as the ``power-of-two'' or ``power-of-choice'' effect.
Specifically, when tasks arrive at rate $\lambda N$,
the queue length distribution at each individual server exhibits
super-exponential decay for any fixed $\lambda < 1$ as $N$ grows large,
a considerable improvement compared to exponential decay for purely random assignment.

As illustrated by the above, the diversity parameter~$d$ induces
a fundamental trade-off between the amount of communication overhead
and the delay performance.
Specifically, a random assignment policy does not entail any
communication burden, but the mean waiting time remains \emph{constant}
as $N$ grows large for any fixed $\lambda > 0$.
In contrast, a nominal implementation of the JSQ policy (without
maintaining state information at the dispatcher) involves $2 N$
messages per task, but the mean waiting time \emph{vanishes}
as $N$ grows large for any fixed $\lambda < 1$.
Although JSQ($d$) policies with $d \geq 2$ yield major performance
improvements over purely random assignment while reducing the
communication burden by a factor O($N$) compared to the JSQ policy,
the mean waiting time \emph{does not vanish} in the limit.
Thus, no fixed value of~$d$ will provide asymptotically optimal
delay performance.
This is evidenced by results of Gamarnik {\em et al.}~\cite{GTZ16}
indicating that in the absence of any memory at the dispatcher the
communication overhead per task \emph{must increase} with~$N$ in order
for any scheme to achieve a zero mean waiting time in the limit.

We will explore the intrinsic trade-off between delay performance
and communication overhead as governed by the diversity parameter~$d$,
in conjunction with the relative load~$\lambda$.
The latter trade-off is examined in an asymptotic regime where not
only the overall task arrival rate is assumed to grow with~$N$,
but also the diversity parameter is allowed to depend on~$N$.
We write $\lambda(N)$ and $d(N)$, respectively, to explicitly reflect
that, and investigate what growth rate of $d(N)$ is required,
depending on the scaling behavior of $\lambda(N)$, in order to achieve
a zero mean waiting time in the limit.
The analysis covers both fluid-scaled and diffusion-scaled versions of
the queue length process in regimes where $\lambda(N) / N \to \lambda < 1$
and $(N - \lambda(N)) / \sqrt{N} \to \beta > 0$ as $N \to \infty$,
respectively.
We establish that the limiting processes are insensitive to the exact
growth rate of $d(N)$, as long as the latter is sufficiently fast,
and in particular coincide with the limiting processes for the JSQ policy.
This reflects a remarkable universality property and demonstrates
that the optimality of the JSQ policy can asymptotically be preserved 
while dramatically lowering the communication overhead.

We will extend these universality properties to network
scenarios where the $N$~servers are assumed to be inter-connected
by some underlying graph topology $G_N$.
Tasks arrive at the various servers as independent Poisson processes
of rate~$\lambda$, and each incoming task is assigned to whichever
server has the shortest queue among the one where it appears and its
neighbors in $G_N$.
In case $G_N$ is a clique (fully connected graph), each incoming task is assigned to the
server with the shortest queue across the entire system,
and the behavior is equivalent to that under the JSQ policy.
The  stochastic optimality properties of the JSQ policy
thus imply that the queue length process in a clique will be `better'
than in an arbitrary graph~$G_N$.
We will establish sufficient conditions for the fluid-scaled
and diffusion-scaled versions of the queue length process in
an arbitrary graph to be equivalent to the limiting processes
in a clique as $N \to \infty$.
The conditions reflect similar universality properties as described above,
and in particular demonstrate that the optimality of a clique can
be asymptotically preserved while dramatically reducing the number
of connections, provided the graph $G_N$ is suitably random.

While a zero waiting time can be achieved in the limit by sampling
only $d(N) \ll N$ servers, the amount of communication overhead
in terms of $d(N)$ must still grow with~$N$.
This may be explained from the fact that a large number of servers
need to be sampled for each incoming task to ensure that at least one
of them is found idle with high probability.
As alluded to above, this can be avoided by introducing memory
at the dispatcher, in particular maintaining a record of vacant servers,
and assigning tasks to idle servers, if there are any.
This so-called Join-the-Idle-Queue (JIQ) scheme \cite{BB08,LXKGLG11}
has gained huge popularity recently, and can be implemented through
a simple token-based mechanism generating at most one message per task.
As shown by Stolyar~\cite{Stolyar15}, the fluid-scaled queue length
process under the JIQ scheme is equivalent to that under the JSQ policy
as $N \to \infty$, and we will extend this result to the
diffusion-scaled queue length process.
Thus, the use of memory allows the JIQ scheme to achieve asymptotically
optimal delay performance with minimal communication overhead.
In particular, ensuring that tasks are assigned to idle servers
whenever available is sufficient to achieve asymptotic optimality,
and using any additional queue length information yields no meaningful
performance benefits on the fluid or diffusion levels.

Stochastic coupling techniques play an instrumental role in the proofs
of the above-described universality and asymptotic optimality properties.
A direct analysis of the queue length processes under a JSQ($d(N)$)
policy, in a load balancing graph $G_N$, or under the JIQ scheme
is confronted with formidable obstacles, and does not seem tractable.
As an alternative route, we leverage novel stochastic coupling
constructions to relate the relevant queue length processes to the
corresponding processes under a JSQ policy, and show that the
deviation between these processes is asymptotically negligible
under suitable assumptions on $d(N)$ or $G_N$.

While the stochastic coupling schemes provide an effective
and overarching approach, they defy a systematic recipe and involve
some degree of ingenuity and customization.
Indeed, the specific coupling arguments that we develop are not only
different from those that were originally used in establishing the
stochastic optimality properties of the JSQ policy, but also differ
in critical ways between a JSQ($d(N)$) policy, a load balancing graph
$G_N$, and the JIQ scheme.
Yet different coupling constructions are devised for model variants
with infinite-server dynamics that we will discuss in Section~\ref{bloc}.

For readability, we occasionally use somewhat informal arguments and phrases in this introductory chapter, but completely rigorous statements and proofs can be found in the subsequent chapters.
In order for some of the sections and chapters to be mostly self-contained, we have also allowed for a certain degree of repetition in a few places.

The remainder of this introduction is organized as follows.
In Section~\ref{spec} we discuss various LBAs and evaluate
their scalability properties.
In Section~\ref{sec:powerofd} we introduce some useful preliminary concepts,
and then review fluid and diffusion limits for the JSQ policy as well
as JSQ($d$) policies with a fixed value of~$d$.
In Section~\ref{univ} we explore the trade-off between delay
performance and communication overhead as function of the diversity
parameter~$d$, in conjunction with the relative load.
In particular, we establish asymptotic universality properties
for JSQ($d$) policies, which are extended to systems with server pools
and network scenarios in Sections~\ref{bloc} and~\ref{networks},
respectively.
In Section~\ref{token} we establish asymptotic optimality properties
for the JIQ scheme.
We discuss somewhat related redundancy policies and alternative scaling
regimes and performance metrics in Section~\ref{miscellaneous}.
The chapter is concluded in Section~\ref{sec:ext} with a discussion of yet
further extensions and several open problems and emerging research
directions.

\section{Scalability spectrum}
\label{spec}

In this section we review a wide spectrum of LBAs and examine their
scalability properties in terms of the delay performance vis-a-vis
the associated implementation overhead in large-scale systems.

\subsection{Basic model}

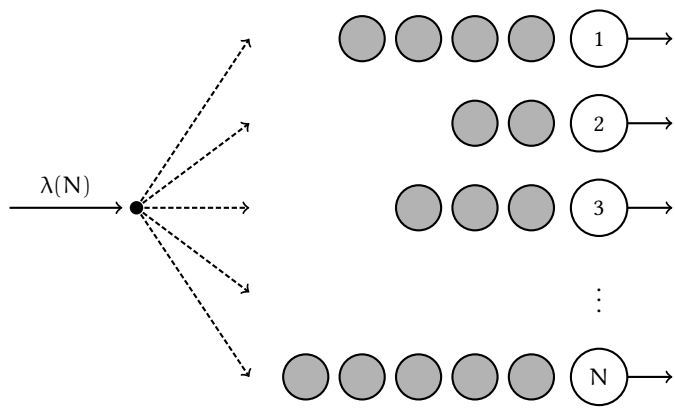
\begin{figure}
\begin{center}
\begin{tikzpicture}[scale=0.75]
\node[above] at (3.750,5) {\small$\lambda(N)$};
\draw[thick,black,->] (2.75,5)--(4.75,5);
\draw[thick,black,fill=black](5,5) circle [radius=0.1];
\draw[thick,black,->,dash pattern=on 2 off 1] (5,5)--(7,5);
\draw[thick,black,->,dash pattern=on 2 off 1] (5,5)--(7,6.5);
\draw[thick,black,->,dash pattern=on 2 off 1] (5,5)--(7,8);
\draw[thick,black,->,dash pattern=on 2 off 1] (5,5)--(7,3.5);
\draw[thick,black,->,dash pattern=on 2 off 1] (5,5)--(7,2);
\foreach \i in {1,...,4}
{
\draw[thick,black,fill=mygray](13-\i,8) circle [radius=0.4];
}
\foreach \i in {1,...,2}
{
\draw[thick,black,fill=mygray](13-\i,6.5) circle [radius=0.4];
}
\foreach \i in {1,...,3}
{
\draw[thick,black,fill=mygray](13-\i,5) circle [radius=0.4];
}
\foreach \i in {1,...,5}
{
\draw[thick,black,fill=mygray](13-\i,2) circle [radius=0.4];
}
\draw[thick,black](13.2,8) circle [radius=0.5];
\draw[thick,black](13.2,6.5) circle [radius=0.5];
\draw[thick,black](13.2,5) circle [radius=0.5];
\draw[thick,black](13.2,2) circle [radius=0.5];
\draw[thick,black,->] (13.7,8)--(14.5,8);
\draw[thick,black,->] (13.7,6.5)--(14.5,6.5);
\draw[thick,black,->] (13.7,5)--(14.5,5);
\draw[thick,black,->] (13.7,2)--(14.5,2);
\node at (13.2,8) {\small 1};
\node at (13.2,6.5) {\small 2};
\node at (13.2,5) {\small 3};
\node at (13.2,3.5) {$\vdots$};
\node at (13.2,2) {\small $N$};
\end{tikzpicture}
\end{center}
\caption{Tasks arrive at the dispatcher as a Poisson process
of rate $\lambda(N)$, and are forwarded to one of the $N$ servers
according to some specific load balancing algorithm.}
\label{figJSQ}
\end{figure}

Throughout this section and most of the chapter, we focus on a basic
scenario with $N$ parallel single-server infinite-buffer queues
and a single dispatcher where tasks arrive as a Poisson process
of rate~$\lambda(N)$, as depicted in Figure~\ref{figJSQ}.
Arriving tasks cannot be queued at the dispatcher,
and must immmediately be forwarded to one of the servers.
This canonical setup is commonly dubbed the \emph{supermarket model}.
Tasks are assumed to have unit-mean exponentially distributed service
requirements, and the service discipline at each server is supposed
to be oblivious to the actual service requirements.

When tasks do not get served and never depart but simply accumulate,
the above setup corresponds to a so-called balls-and-bins model,
and we will further elaborate on the connections and differences
with work in that domain in Subsection~\ref{ballsbins}.

\subsection{Asymptotic scaling regimes}
\label{asym}

An exact analysis of the delay performance is quite involved,
if not intractable, for all but the simplest LBAs.
A common approach is therefore to consider various limit regimes,
which not only provide mathematical tractability and illuminate the
fundamental behavior, but are also natural in view of the typical
conditions in which cloud networks and data centers operate.
One can distinguish several asymptotic scalings that have been used
for these purposes:

(i) In the classical heavy-traffic regime, $\lambda(N) = \lambda N$
with a {\em fixed} number of servers~$N$ and a relative load~$\lambda$
that tends to one (i.e., there is no asymptotics in~$N$).

(ii) In the conventional large-capacity or many-server regime,
the relative load $\lambda(N) / N$ approaches a constant $\lambda < 1$
as the number of servers~$N$ grows large.

(iii) The popular Halfin-Whitt regime~\cite{HW81} combines heavy traffic
with a large capacity, with
\begin{equation}
\label{eq:HW}
\frac{N - \lambda(N)}{\sqrt{N}} \to \beta > 0 \mbox{ as } N \to \infty,
\end{equation}
so the relative capacity slack behaves as $\beta / \sqrt{N}$
as the number of servers~$N$ grows large.

(iv) The so-called non-degenerate slow-down regime~\cite{Atar12}
involves $N - \lambda(N) \to \gamma > 0$, so the relative capacity
slack shrinks as $\gamma / N$ as the number of servers~$N$ grows large.

The term non-degenerate slow-down refers to the fact that in the
context of a centralized multi-server queue, the mean waiting time
in regime (iv) tends to a strictly positive constant as $N \to \infty$,
and is thus of similar magnitude as the mean service requirement.
In contrast, in regimes (ii) and (iii), the mean waiting time in
a multi-server queue decays exponentially fast in~$N$ or is of the
order~$1 / \sqrt{N}$, respectively, as $N \to \infty$,
while in regime (i) the mean waiting time grows arbitrarily large
relative to the mean service requirement.

In the context of a centralized M/M/N queue, scalings  (ii), (iii) and (iv) are commonly referred to as
Quality-Driven (QD), Quality-and-Efficiency-Driven (QED)
and Efficiency-Driven (ED) regimes.
These terms reflect that (ii) offers excellent service quality
(vanishing waiting time), (iv) provides high resource efficiency
(utilization approaching one), and (iii) achieves a combination
of these two, providing the best of both worlds.

In the present thesis, and in particular in the current chapter we will focus on scalings (ii) and (iii),
and occasionally also refer to these as fluid and diffusion scalings, 
since it is natural to analyze the relevant queue length process
on fluid scale ($1 / N$) and diffusion scale ($1 / \sqrt{N}$)
in these regimes, respectively.
We will not provide a detailed account of scalings (i) and (iv),
which do not capture the large-scale perspective and do not allow
for low delays, respectively, but we will briefly mention some results
for these regimes in Subsections~\ref{classical}
and~\ref{nondegenerate}.

An important issue in the context of scaling limits is the rate
of convergence and the accuracy for finite-size systems.
Some interesting results for the accuracy of mean-field approximations
for interacting-particle networks and in particular load balancing
models may be found in recent work of Gast~\cite{Gast17}, Gast \& Van Houdt~\cite{GH18}, and Ying \cite{Ying16,Ying17}.

\subsection{Random assignment: N independent M/M/1 queues}
\label{random}

One of the most basic LBAs is to assign each arriving task to a server
selected uniformly at random.
In that case, the various queues collectively behave as
$N$~independent M/M/1 queues, each with arrival rate $\lambda(N) / N$
and unit service rate.
In particular, at each of the queues, the total number of tasks in
stationarity has a geometric distribution with parameter $\lambda(N) / N$.
By virtue of the PASTA property, the probability that an arriving task
incurs a non-zero waiting time is $\lambda(N) / N$.
The mean number of waiting tasks (excluding the possible task in service)
at each of the queues is $\frac{\lambda(N)^2}{N (N - \lambda(N))}$,
so the total mean number of waiting tasks is
$\frac{\lambda(N)^2}{N - \lambda(N)}$, which by Little's law implies that
the mean waiting time of a task is $\frac{\lambda(N)}{N - \lambda(N)}$.
In particular, when $\lambda(N) = N \lambda$, the probability that
a task incurs a non-zero waiting time is $\lambda$,
and the mean waiting time of a task is $\frac{\lambda}{1 - \lambda}$,
independent of~$N$, reflecting the independence of the various queues.

As we will see later, a broad range of queue-aware LBAs can deliver
a probability of a non-zero waiting time and a mean waiting time
that vanish asymptotically.
While a random assignment policy is evidently not competitive
with such queue-aware LBAs, it still plays a relevant role due to
the strong degree of tractability inherited from its simplicity.
For example, the queue process under purely random assignment can
be shown to provide an upper bound (in a stochastic majorization sense)
for various more involved queue-aware LBAs for which even stability
may be difficult to establish directly, yielding conservative
performance bounds and stability guarantees.

A slightly better LBA is to assign tasks to the servers in
a Round-Robin manner, dispatching every $N$-th task to the same server.
In the fluid regime where $\lambda(N) = N \lambda$, the inter-arrival
time of tasks at each given queue will then converge to a constant
$1 / \lambda$ as $N \to \infty$.
Thus each of the queues will behave as a D/M/1 queue in the limit,
and the probability of a non-zero waiting time and the mean waiting
time will be somewhat lower than under purely random assignment.
However, both the probability of a non-zero waiting time and the mean
waiting time will still tend to strictly positive values and not vanish
as $N \to \infty$.

\subsection{Join-the-Shortest Queue (JSQ)}
\label{ssec:jsq}

Under the Join-the-Shortest-Queue (JSQ) policy, each arriving task is
assigned to the server with the currently shortest queue.
In the basic model described above, the JSQ policy has several strong
stochastic optimality properties, and yields the `most balanced
and smallest' queue process among all non-anticipating policies that
do not have any advance knowledge of the service requirements
\cite{EVW80,Winston77}.

\subsection{Join-the-Smallest-Workload (JSW): centralized M/M/N queue}
\label{ssec:jsw}

Under the Join-the-Smallest-Workload (JSW) policy, each arriving task
is assigned to the server with the currently smallest workload.
Note that this is an anticipating policy, since it requires advance
knowledge of the service requirements of all the tasks in the system.
Further observe that this policy (myopically) minimizes the waiting
time for each incoming task, and mimicks the operation of a centralized
$N$-server queue with a FCFS discipline.
The equivalence with a centralized $N$-server queue with a FCFS
discipline yields a strong optimality property of the JSW policy:
The vector of joint workloads at the various servers observed by each
incoming task is smaller in the Schur convex sense than under any
alternative admissible policy~\cite{FC01}.

It is worth observing that the above optimality properties in fact do
not rely on Poisson arrival processes or exponential service
requirement distributions.
Even though the JSW policy requires a similar excessive communication
overhead as the JSQ policy, aside from its anticipating nature, the equivalence
with a centralized FCFS queue means that there cannot be any idle
servers while tasks are waiting and that the total number of tasks behaves
as a birth-death process, which renders it far more tractable.
Specifically, given that all the servers are busy, the total number
of waiting tasks is geometrically distributed with parameter
$\lambda(N) / N$.
Thus the total mean number of waiting tasks is
$\Pi_W(N, \lambda(N)) \frac{\lambda(N)}{N - \lambda(N)}$,
and the mean waiting time is
$\Pi_W(N, \lambda(N)) \frac{1}{N - \lambda(N)}$,
with $\Pi_W(N, \lambda(N)$ denoting the probability of all servers
being occupied and a task incurring a non-zero waiting time.
This immediately shows that the mean waiting time is smaller
by at least a factor $\lambda(N)$ than for the random assignment
policy considered in Subsection~\ref{random}.

In the large-capacity regime $\lambda(N) = N \lambda$, it can be shown
that the probability $\Pi_W(N, \lambda(N))$ of a non-zero waiting time
decays exponentially fast in~$N$, and hence so does the mean waiting time.
In the Halfin-Whitt heavy-traffic regime~\eqref{eq:HW}, the probability
$\Pi_W(N, \lambda(N))$ of a non-zero waiting time converges to a finite
constant $\Pi_W^*(\beta)$, implying that the mean waiting time of a task
is of the order $1 / \sqrt{N}$, and thus vanishes as $N \to \infty$.

\subsection{Power-of-d load balancing (JSQ(d))}
\label{ssec:powerd}

We have seen that the achilles heel of the JSQ policy is its
excessive communication overhead in large-scale systems.
This poor scalability has motivated consideration of so-called
JSQ($d$) policies, where an incoming task is assigned to a server
with the shortest queue among $d$~servers selected uniformly at random.
Results in Mitzenmacher~\cite{Mitzenmacher01} and Vvedenskaya
{\em et al.}~\cite{VDK96} indicate that 
in the fluid regime where $\lambda(N) = \lambda N$,
the probability that there are $i$ or more tasks at a given queue is
proportional to $\lambda^{\frac{d^i - 1}{d - 1}}$ as $N \to \infty$,
and thus exhibits super-exponential decay as opposed to exponential decay
for the random assignment policy considered in Subsection~\ref{random}.

The diversity parameter~$d$ thus induces
a fundamental trade-off between the amount of communication overhead
and the performance in terms of queue lengths and delays.
A rudimentary implementation of the JSQ policy ($d = N$, without
replacement) involves $O(N)$ communication overhead per task,
but it can be shown that the probability of a non-zero waiting time
and the mean waiting \emph{vanish} as $N \to \infty$, just like
in a centralized queue.
Although JSQ($d$) policies with a fixed parameter $d \geq 2$ yield
major performance improvements, 
the probability of a non-zero waiting time
and the mean waiting time \emph{do not vanish} as $N \to \infty$.


\subsection{Token-based mechanisms: Join-the-Idle-Queue (JIQ)}
\label{ssec:jiq}

While a zero waiting time can be achieved in the limit by sampling
only $d(N) \ll N$ servers, the amount of communication overhead
in terms of $d(N)$ must still grow with~$N$.
This can be countered by introducing memory at the dispatcher,
in particular maintaining a record of vacant servers,
and assigning tasks to idle servers as long as there are any,
or to a uniformly at random selected server otherwise.
This so-called Join-the-Idle-Queue (JIQ) scheme \cite{BB08,LXKGLG11}
has received keen interest recently, and can be implemented through
a simple token-based mechanism.
Specifically, idle servers send tokens to the dispatcher to advertize
their availability, and when a task arrives and the dispatcher has
tokens available, it assigns the task to one of the corresponding
servers (and disposes of the token).
Note that a server only issues a token when a task completion leaves
its queue empty, thus generating at most one message per task.
Surprisingly, the mean waiting time and the probability of a non-zero
waiting time vanish under the JIQ scheme in both the fluid and diffusion
regimes, as we will further discuss in Section~\ref{token}.
Thus, the use of memory allows the JIQ scheme to achieve asymptotically
optimal delay performance with minimal communication overhead.

\subsection{Performance comparison}
\label{ssec:perfcomp}

We now present some simulation experiments 
to compare the above-described LBAs in terms of delay performance.
\begin{figure}
\includegraphics[width=\linewidth]{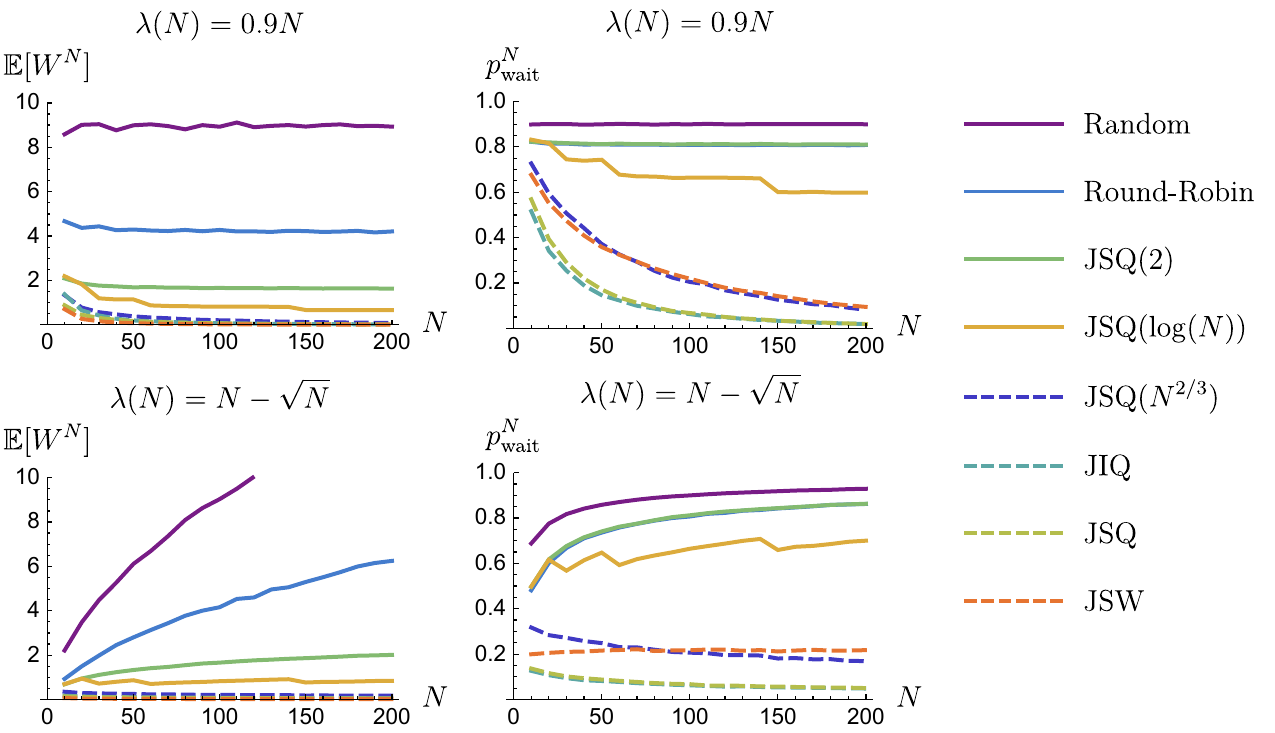}
\caption{Simulation results for mean waiting time $\mathbb{E}[W^N]$
and probability of a non-zero waiting time $p_{\textup{wait}}^N$,
for both a fluid regime and a diffusion regime.}
\label{differentschemes}
\end{figure}
Specifically, we evaluate the mean waiting time and the probability
of a non-zero waiting time in both a fluid regime ($\lambda(N) = 0.9 N$)
and a diffusion regime ($\lambda(N) = N - \sqrt{N}$).
The results are shown in Figure~\ref{differentschemes}.
An overview of the delay performance and overhead associated with
various LBAs is given in Table~\ref{table}.

We are specifically interested in distinguishing two classes of LBAs --
the ones delivering a mean waiting time and probability of a non-zero
waiting time that vanish asymptotically, and the ones that fail to do so
-- and relating that dichotomy to the associated communication
overhead and memory requirement at the dispatcher.
We give these classifications for both the fluid regime and the
diffusion regime.

\paragraph{JSQ, JIQ and JSW.}
Three schemes that clearly have vanishing waiting time are JSQ, JIQ and JSW.
The optimality of JSW is observed in the figures;
JSW has the smallest mean waiting time, and all three schemes have
vanishing waiting time in both the fluid and diffusion regime. 

However, there is a significant difference between JSW and JSQ/JIQ.
We observe that the probability of positive wait does not vanish for JSW,
while it does vanish for JSQ/JIQ.
This implies that the mean of all positive waiting times is an order
larger in JSQ/JIQ compared to JSW.
Intuitively, this is clear since in JSQ/JIQ, when a task is placed
in a queue, it waits for at least one specific other task.
In JSW, which is equivalent to the M/M/N queue, a task that cannot
start service immediately, can start service when one of the $N$~servers
becomes idle.

\paragraph{Random and Round-Robin.}
The mean waiting time does not vanish for Random and Round-Robin
in the fluid regime, as already mentioned in Subsection \ref{random}.
Moreover, the waiting time grows without bound in the diffusion regime
for these two schemes.
This is because the system can still be decomposed into single-server queues, and the loads
of the individual M/M/1 and D/M/1 queues tend to~1.

\paragraph{JSQ($d$) policies.}
Three versions of JSQ($d$) are included in Figure~\ref{differentschemes};
$d(N)=2\not\to \infty$, $d(N)=\lfloor\log(N)\rfloor\to \infty$
and $d(N)=N^{2/3}$ for which $\frac{d(N)}{\sqrt{N}\log(N)}\to \infty$.
Note that the graph for $d(N)=\lfloor \log(N) \rfloor$ shows sudden
jumps when $d(N)$ increases by~1.
As can be seen in Figure~\ref{differentschemes}, the variants
for which $d(N)\to\infty$ have vanishing wait in the fluid regime,
while $d = 2$ does not.
The latter could be readily observed, since JSQ($d$) uses no memory
and the overhead per task does not increase with~$N$,
as already mentioned in the introduction.
Furthermore, it follows that JSQ($d$) policies clearly outperform
Random and Round-Robin dispatching, while JSQ/JIQ/JSW are better
in terms of mean wait.

\begin{sidewaystable}\centering
\def\arraystretch{2}%
\begin{tabular}{|C{3cm}|C{3cm}|C{2.5cm}|C{3cm}|C{1.75cm}|}
\hline
Scheme & Queue length & Waiting time (fixed $\lambda < 1$) &
Waiting time ($1 - \lambda \sim 1 / \sqrt{N}$) & Overhead per task \\
\hline\hline
Random & $q_i^\star = \lambda^i$ &  $\frac{\lambda}{1 - \lambda}$ &
$\Theta(\sqrt{N})$ & 0 \\
\hline
JSQ($d$) & $q_i^\star = \lambda^{\frac{d^i - 1}{d - 1}}$ &
$\Theta$(1) & $\Omega(\log{N})$ & $2 d$ \\
\hline
\mbox{$d(N)$} \mbox{$\to \infty$} & same as JSQ & same as JSQ & ?? &
$2d(N)$ \\
\hline
$\frac{d(N)}{\sqrt{N} \log(N)}\to \infty$ &
same as JSQ & same as JSQ & same as JSQ & $2d(N)$ \\
\hline
JSQ & $q_1^\star = \lambda$, $q_2^\star =$ o(1) & o(1) &
$\Theta(1 / \sqrt{N})$ & $2 N$ \\
\hline\hline
JIQ & same as JSQ & same as JSQ & same as JSQ & $\leq 1$ \\
\hline
\end{tabular}
\caption{Queue length distribution, waiting times, and communication
overhead for various~LBAs.}
\label{table}
\end{sidewaystable}

\section{Preliminaries, JSQ policy, and power-of-d algorithms}
\label{sec:powerofd}

In this section we first introduce some useful notation and preliminary
concepts, and then review fluid and diffusion limits for the JSQ policy
as well as JSQ($d$) policies with a fixed value of~$d$.

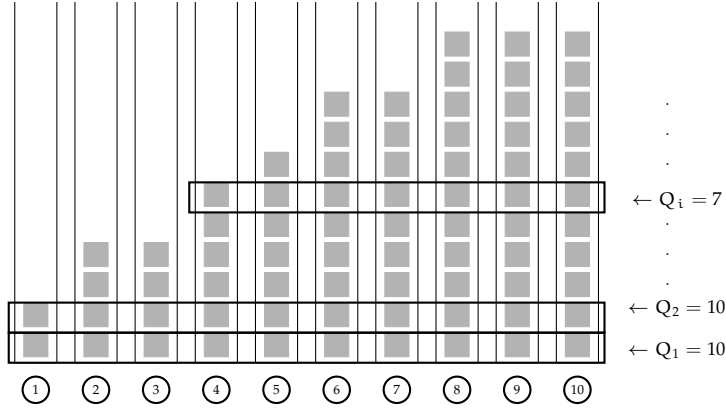
\begin{figure}
\begin{center}
\begin{tikzpicture}[scale=.80]
\foreach \x in {10, 9,...,1}
    \draw (\x,6)--(\x,0)--(\x+.7,0)--(\x+.7,6);
\foreach \x in {10, 9,...,1}
    \draw (\x+.35,-.45) node[circle,inner sep=0pt, minimum size=10pt,draw,thick] {{{\tiny $\mathsmaller{\x}$}}} ;
\foreach \y in {0, .5}
    \draw[fill=mygray,mygray] (1.15,.1+\y) rectangle (1.55,.5+\y);
\foreach \y in {0, .5, 1, 1.5}
    \draw[fill=mygray,mygray] (2.15,.1+\y) rectangle (2.55,.5+\y);
\foreach \y in {0, .5, 1, 1.5}
    \draw[fill=mygray,mygray] (3.15,.1+\y) rectangle (3.55,.5+\y);
\foreach \y in {0, .5, 1, 1.5, 2, 2.5}
    \draw[fill=mygray,mygray] (4.15,.1+\y) rectangle (4.55,.5+\y);
\foreach \y in {0, .5, 1, 1.5, 2, 2.5, 3}
    \draw[fill=mygray,mygray] (5.15,.1+\y) rectangle (5.55,.5+\y);
\foreach \y in {0, .5, 1, 1.5, 2, 2.5, 3, 3.5, 4}
    \draw[fill=mygray,mygray] (6.15,.1+\y) rectangle (6.55,.5+\y);
\foreach \y in {0, .5, 1, 1.5, 2, 2.5, 3, 3.5, 4}
    \draw[fill=mygray,mygray] (7.15,.1+\y) rectangle (7.55,.5+\y);
\foreach \y in {0, .5, 1, 1.5, 2, 2.5, 3, 3.5, 4, 4.5, 5}
    \draw[fill=mygray,mygray] (8.15,.1+\y) rectangle (8.55,.5+\y);
\foreach \y in {0, .5, 1, 1.5, 2, 2.5, 3, 3.5, 4, 4.5, 5}
    \draw[fill=mygray,mygray] (9.15,.1+\y) rectangle (9.55,.5+\y);
\foreach \y in {0, .5, 1, 1.5, 2, 2.5, 3, 3.5, 4, 4.5, 5}
    \draw[fill=mygray,mygray] (10.15,.1+\y) rectangle (10.55,.5+\y);

\draw[thick] (.9,0) rectangle (10.8,.5);
\draw[thick] (.9,.5) rectangle (10.8,1);
\draw[thick] (3.9,2.5) rectangle (10.8,3);

\draw  (12, .2) node {{\scriptsize $\leftarrow Q_1=10$}};
\draw  (12, .9) node {{\scriptsize $\leftarrow Q_2=10$}};

\draw  (11.85, 1.3) node {{\tiny $\cdot$}};
\draw  (11.85, 1.8) node {{\tiny $\cdot$}};
\draw  (11.85, 2.3) node {{\tiny $\cdot$}};
\draw  (12, 2.7) node {{\scriptsize $\leftarrow Q_i=7$}};
\draw  (11.85, 3.3) node {{\tiny $\cdot$}};
\draw  (11.85, 3.8) node {{\tiny $\cdot$}};
\draw  (11.85, 4.3) node {{\tiny $\cdot$}};
\end{tikzpicture}
\end{center}
\caption{The value of $Q_i$ represents the width of the
$i$-th row, when the servers are arranged in non-descending order
of their queue lengths.}
\label{figB}
\end{figure}

We keep focusing  on a basic scenario
where all the servers are homogeneous, the service requirements are
exponentially distributed, and the service discipline at each server
is oblivious of the actual service requirements.
In order to obtain a Markovian state description, it therefore
suffices to only track the number of tasks, and in fact we do not need
to keep record of the number of tasks at each individual server,
but only count the number of servers with a given number of tasks.
Specifically, we represent the state of the system by a vector
\begin{equation}
\label{sos}
\QQ(t) := \left(Q_1(t), Q_2(t), \dots\right),
\end{equation}
with $Q_i(t)$ denoting the number of servers with $i$ or more tasks
at time~$t$, including the possible task in service, $i = 1, 2 \dots$.
Note that if we represent the queues at the various servers as (vertical)
stacks, and arrange these from left to right in non-descending order,
then the value of $Q_i$ corresponds to the width of the $i$-th (horizontal)
row, as depicted in the schematic diagram in Figure~\ref{figB}.

In order to examine the fluid and diffusion limits in regimes
where the number of servers~$N$ grows large, we consider a sequence
of systems indexed by~$N$, and attach a superscript~$N$ to the
associated state variables.

The fluid-scaled occupancy state is denoted by
$\qq^N(t) := (q_1^N(t), q_2^N(t), \dots)$, with $q_i^N(t) = Q_i^N(t) / N$
representing the fraction of servers in the $N$-th system with $i$
or more tasks as time~$t$, $i = 1, 2, \dots$.
Let
$\cS = \{\qq \in [0, 1]^\infty: q_i \leq q_{i-1} \forall i = 2, 3,\dots\}$ 
be the set of all possible fluid-scaled states.
Whenever we consider fluid limits, we assume the sequence of initial
states is such that $\qq^N(0) \to \qq^\infty \in \cS$ as $N \to \infty$.

The diffusion-scaled occupancy state is defined as
$\bar{\QQ}^N(t) = (\bar{Q}_1^N(t), \bar{Q}_2^N(t), \dots)$, with
\begin{equation}\label{eq:diffscale}
\bar{Q}_1^N(t) = - \frac{N - Q_1^N(t)}{\sqrt{{N}}}, \qquad
\bar{Q}_i^N(t) = \frac{Q_i^N(t)}{\sqrt{{N}}}, \quad i = 2,3, \dots.
\end{equation}
Note that $-\bar{Q}_1^N(t)$ corresponds to the number of vacant servers,
normalized by $\sqrt{N}$.
The reason why $Q_1^N(t)$ is centered around~$N$ while $Q_i^N(t)$,
$i = 2,3, \dots$, are not, is that for the scalable LBAs that we
consider, the fraction of servers with exactly one task tends to one,
whereas the fraction of servers with two or more tasks tends to zero
as $N\to\infty$.
For convenience, we will assume that each server has an infinite-capacity
buffer, but all the results extend to the finite-buffer case.

\subsection{Fluid limit for JSQ(d) policies}

We first consider the fluid limit for JSQ($d$) policies
with an arbitrary but fixed value of~$d$ as characterized
by Mitzenmacher~\cite{Mitzenmacher01} and Vvedenskaya
{\em et al.}~\cite{VDK96}:

{\em The sequence of processes $\{\qq^N(t)\}_{t \geq 0}$ has a weak
limit $\{\qq(t)\}_{t \geq 0}$ that satisfies the system
of differential equations}
\begin{equation}
\label{fluid:standard}
\frac{\dif q_i(t)}{\dif t} =
\lambda (q_{i-1}^d(t) - q_i^d(t)) - (q_i(t) - q_{i+1}(t)).
\quad i = 1, 2, \dots.
\end{equation}
The fluid-limit equations may be interpreted as follows.
The first term represents the rate of increase in the fraction
of servers with $i$ or more tasks due to arriving tasks that are
assigned to a server with exactly $i - 1$ tasks.
Note that the latter occurs in fluid state $\qq \in \cS$ with probability
$q_{i-1}^d - q_i^d$, i.e., the probability that all $d$~sampled servers
have $i - 1$ or more tasks, but not all of them have $i$ or more tasks.
The second term corresponds to the rate of decrease in the fraction
of servers with $i$ or more tasks due to service completions from servers
with exactly $i$ tasks, and the latter rate is given by $q_i - q_{i+1}$.
While the system in~\eqref{fluid:standard} characterizes the functional law of large numbers (FLLN) behavior of systems under the JSQ($d$) scheme, weak convergence to a certain Ornstein-Ulenbeck process (both in the transient regime and in steady state) was shown in~\cite{Graham05}, establishing a functional central limit theorem (FCLT) result.
Strong approximations for systems under the JSQ($d$) scheme on any finite time interval by the deterministic system in~\eqref{fluid:standard}, a certain infinite-dimensional jump process, and a diffusion approximation were established in~\cite{LN05}.

When the derivatives in~\eqref{fluid:standard} are set equal to zero
for all~$i$, the unique fixed point for any $d \geq 2$ is obtained as
\begin{equation}
\label{eq:fixedpoint1}
q_i^* = \lambda^{\frac{d^i-1}{d-1}}.
\quad i = 1, 2, \dots.
\end{equation}
It can be shown that the fixed point is asymptotically stable in the
sense that $\qq(t) \to \qq^*$ as $t \to \infty$ for any initial fluid
state $\qq^\infty$ with $\sum_{i = 1}^{\infty} q_i^\infty < \infty$.

As mentioned earlier, the fixed point reveals that the stationary queue length distribution at
each individual server exhibits super-exponential decay as $N \to \infty$,
as opposed to exponential decay for a random assignment policy.

It is worth observing that this involves an interchange of the
many-server ($N \to \infty$) and stationary ($t \to \infty$) limits.
The justification is provided by the asymptotic stability of the fixed
point along with a few further technical conditions.

\subsection{Fluid limit for JSQ policy}
\label{ssec:jsqfluid}

We now turn to the fluid limit for the ordinary JSQ policy,
which rather surprisingly was not rigorously established until fairly
recently in~\cite{MBLW16-3}, leveraging martingale functional limit
theorems and time-scale separation arguments~\cite{HK94}.
A more detailed description of the fluid limit along with the proofs is presented in Chapter~\ref{chap:univjsqd}.

In order to state the fluid limit starting from an arbitrary
fluid-scaled occupancy state, we first introduce some additional notation.
For any fluid state $\qq \in \cS$,
denote by $m(\qq) = \min\{i: q_{i + 1} < 1\}$ the minimum queue length
among all servers.
Now if $m(\qq)=0$, then define $p_0(m(\qq))=1$ and $p_i(m(\qq))=0$
for all $i=1,2,\ldots$. 
Otherwise, in case $m(\qq)>0$, define
\begin{equation}
\label{eq:fluid-gen-intro}
p_i(\qq) =
\begin{cases}
\min\big\{(1 - q_{m(\qq) + 1})/\lambda,1\big\} & \quad\mbox{ for }\quad i=m(\qq)-1, \\
1 - p_{m(\qq) - 1}(\qq) & \quad\mbox{ for }\quad i=m(\qq), \\
0&\quad \mbox{ otherwise.}
\end{cases}
\end{equation}
{\em Any weak limit of the sequence of processes $\{\qq^N(t)\}_{t \geq 0}$
is given by the deterministic system $\{\qq(t)\}_{t \geq 0}$ satisfying
the system of differential equations
\begin{equation}
\label{eq:fluid-intro}
\frac{\dif^+ q_i(t)}{\dif t} =
\lambda p_{i-1}(\qq(t)) - (q_i(t) - q_{i+1}(t)),
\quad i = 1, 2, \dots,
\end{equation}
where $\dif^+/\dif t$ denotes the right-derivative.}
The reason why we have used the derivative in~\eqref{fluid:standard},
and the right-derivative in~\eqref{eq:fluid-intro} is that the limiting
trajectory for the JSQ policy may not be differentiable at all time points.
In fact, one of the major technical challenges in proving the fluid limit
for the JSQ policy is that the drift of the process is not continuous,
which leads to non-smooth limiting trajectories.

As in the case of the fluid-limit for JSQ($d$) policies
in~\eqref{fluid:standard}, the fluid-limit trajectory in~\eqref{eq:fluid-intro}
can be interpreted as follows.
The coefficient $p_i(\qq)$ represents the instantaneous fraction
of incoming tasks assigned to servers with a queue length of exactly~$i$
in the fluid state $\qq \in \mathcal{S}$.
Note that a strictly positive fraction $1 - q_{m(\qq) + 1}$
of the servers have a queue length of exactly~$m(\qq)$.
Clearly the fraction of incoming tasks that get assigned to servers
with a queue length of $m(\qq) + 1$ or larger is zero:
$p_i(\qq) = 0$ for all $i = m(\qq) + 1, \dots$.
Also, tasks at servers with a queue length of exactly~$i$ are completed
at (normalized) rate $q_i - q_{i + 1}$, which is zero for all
$i = 0, \dots, m(\qq) - 1$, and hence the fraction of incoming tasks
that get assigned to servers with a queue length of $m(\qq) - 2$ or less
is zero as well: $p_i(\qq) = 0$ for all $i = 0, \dots, m(\qq) - 2$.
This only leaves the fractions $p_{m(\qq) - 1}(\qq)$
and $p_{m(\qq)}(\qq)$ to be determined.
Now observe that the fraction of servers with a queue length of exactly
$m(\qq) - 1$ is zero.
If $m(\qq)=0$, then clearly the incoming tasks will join an empty queue, 
and thus, $p_{m(\qq)}=1$, and $p_i(\qq) = 0$ for all $i \neq m(\qq)$.
Furthermore, if $m(\qq) \geq 1$, since tasks at servers with a queue
length of exactly $m(\qq)$ are completed at (normalized) rate
$1 - q_{m(\qq) + 1} > 0$, incoming tasks can be assigned to servers
with a queue length of exactly $m(\qq) - 1$ at that rate.
We thus need to distinguish between two cases, depending on whether the
normalized arrival rate $\lambda$ is larger than $1 - q_{m(\qq) + 1}$
or not.
If $\lambda < 1 - q_{m(\qq) + 1}$, then all the incoming tasks can be
assigned to a server with a queue length of exactly $m(\qq) - 1$,
so that $p_{m(\qq) - 1}(\qq) = 1$ and $p_{m(\qq)}(\qq) = 0$.
On the other hand, if $\lambda > 1 - q_{m(\qq) + 1}$, then not all
incoming tasks can be assigned to servers with a queue length of
exactly $m(\qq) - 1$ active tasks, and a positive fraction will be
assigned to servers with a queue length of exactly $m(\qq)$:
$p_{m(\qq) - 1}(\qq) = (1 - q_{m(\qq) + 1}) / \lambda$
and $p_{m(\qq)}(\qq) = 1 - p_{m(\qq) - 1}(\qq)$.

The unique fixed point
$\qq^\star = (q_1^\star,q_2^\star,\ldots)$ of the dynamical system
in~\eqref{eq:fluid-intro} is given by
\begin{equation}
\label{eq:fpjsq}
q_i^* = \left\{\begin{array}{ll} \lambda, & i = 1, \\
0, & i = 2, 3,\dots. \end{array} \right.
\end{equation}
Note that the fixed point naturally emerges when $d \to \infty$ in the
fixed point expression~\eqref{eq:fixedpoint1} for fixed~$d$.
However, the process-level results in \cite{Mitzenmacher01,VDK96}
for fixed~$d$ cannot be readily used to handle joint scalings,
and do not yield the entire fluid-scaled sample path for arbitrary
initial states as given by~\eqref{eq:fluid-intro}.

The fixed point in~\eqref{eq:fpjsq}, in conjunction with an interchange
of limits argument, indicates that in stationarity the fraction
of servers with a queue length of two or larger under the JSQ policy
is negligible as $N \to \infty$.

\subsection{Diffusion limit for JSQ policy}\label{ssec:diffjsq}

We next describe the diffusion limit for the JSQ policy in the
Halfin-Whitt heavy-traffic regime \eqref{eq:HW}, as recently derived
by Eschenfeldt \& Gamarnik~\cite{EG15}.

\paragraph{Transient regime.}
Recall the centered and diffusion-scaled processes in~\eqref{eq:diffscale}.
{\em For suitable initial conditions, the sequence of processes
$\big\{\bar{\QQ}^N(t)\big\}_{t \geq 0}$ converges weakly to the limit
$\big\{\bar{\QQ}(t)\big\}_{t \geq 0}$, 
where 
$(\bar{Q}_1(t), \bar{Q}_2(t),\ldots)$
is the unique solution 
to the system of SDEs
\begin{equation}
\label{eq:diffusionjsq}
\begin{split}
\dif\bar{Q}_1(t) &= \sqrt{2}\dif W(t) - \beta\dif t - \bar{Q}_1(t) + \bar{Q}_2(t)-\dif U_1(t), \\
\dif\bar{Q}_2(t) &= \dif U_1(t) -  (\bar{Q}_2(t)-\bar{Q}_3(t)), \\
\dif\bQ_i(t) &= - (\bQ_i(t) - \bQ_{i+1}(t)), \quad i \geq 3,
\end{split}
\end{equation}
for $t \geq 0$, where $W$ is the standard Brownian motion and $U_1$ is
the unique non-decreasing non-negative process 
satisfying}
$\int_0^\infty \mathbbm{1}_{[\bar{Q}_1(t) < 0]} \dif U_1(t) = 0$.

Now introduce
\[
\bar{Q}_{tot}^N(t) = \frac{Q_{tot}^N(t) - N}{\sqrt{{N}}},
\]
as the centered and diffusion-scaled version of the total number of tasks
$Q_{tot}^N(t) = \sum_{i = 1}^{\infty} Q_i^N(t)$ in the $N$-th system
at time~$t$, and denote by $\bar{Q}_{vac}^N(t) = - \bar{Q}_1^N(t)$ the
diffusion-scaled number of vacant servers in the $N$-th system at time~$t$.
Summing the equations in~\eqref{eq:diffusionjsq} over $i = 1, 2, \dots$,
and rewriting the top equation in terms of $\bar{Q}_{vac}^N(t)$,
we obtain that for suitable initial conditions, the sequence of processes
$\{(\bar{Q}_{tot}^N(t), \bar{Q}_{vac}^N(t))\}_{t \geq 0}$ converges
weakly to the limit $\{(\bar{Q}_{tot}(t), \bar{Q}_{vac}(t))\}_{t \geq 0}$,
as the unique solution to the system of SDEs
\begin{equation}
\begin{split}
\dif\bar{Q}_{tot}(t) &= \sqrt{2} \dif W(t) - \beta\dif t +
\bar{Q}_{vac}(t), \\
\dif\bar{Q}_{vac}(t) &= \sqrt{2}\dif W(t) + \beta\dif t -
(\bar{Q}_{vac}(t) + \bar{Q}_2(t)) +\dif U_1(t),
\end{split}
\end{equation}
for $t \geq 0$, where $W$ is the standard Brownian motion and $U_1$ is
the unique non-decreasing non-negative process satisfying
$\int_0^\infty \mathbbm{1}_{[\bar{Q}_{vac}(t) > 0]} \dif U_1(t) = 0$.

Strikingly, the top equation has the exact same form as in the
corresponding centralized M/M/$N$ queue, while the bottom equation
is nearly identical, except for the term $\bar{Q}_2(t)$.
As it turns out, despite the differences in the dynamics between the JSQ policy and the M/M/$N$ system, there are surprising similarities in terms of the qualitative behavior of the total number of tasks in the system.
We will reflect more on the behavior of the JSQ policy and the M/M/$N$ system in Remark~\ref{rem:comp2-intro} below.

\paragraph{Interchange of limits.}
In~\cite{EG15} the convergence of the scaled occupancy measure was established only in the transient regime on any finite time interval.
The tightness of the diffusion-scaled occupancy measure and the interchange of limits were open until Braverman~\cite{Braverman18} recently further established that the weak-convergence result extends to the steady state as well, i.e., $\bar{\QQ}^N(\infty)$ converges weakly to the random variable $(Q_1(\infty), Q_2(\infty), 0, 0,\ldots)$ as $N\to\infty$, where $(Q_1(\infty), Q_2(\infty))$ has the stationary distribution of the process $(Q_1, Q_2)$.
 Thus, the steady state of the diffusion process in~\eqref{eq:diffusionjsq} is proved to capture the asymptotic behavior of large-scale systems under the JSQ policy.

Although the above interchange of limits result~\cite{Braverman18} establishes that the mean steady-state waiting time under the JSQ policy is of a similar
order $O(1 / \sqrt{N})$ as in the M/M/$N$ queue, it is important to
observe a subtle but fundamental difference in the distributional
properties due to the distributed versus centralized queueing operation.
In the ordinary M/M/$N$ queue a fraction $\Pi_W^*(\beta)$ of the
tasks incur a non-zero waiting time as $N \to \infty$, but a non-zero
waiting time is only of length $1 / (\beta \sqrt{N})$ in expectation.
In contrast, under the JSQ policy, the fraction of tasks that experience
a non-zero waiting time is only of the order $O(1 / \sqrt{N})$.
However, such tasks will have to wait for the duration of a residual
service time, yielding a waiting time of the order $O(1)$.

\paragraph{Tail asymptotics of the steady state.}
In Chapter~\ref{chap:jsqdiffusion} the tail asymptotics of the steady-state distribution $\pi$ of the diffusion in~\eqref{eq:diffusionjsq} will be studied.
In particular, using a classical regenerative process construction of the diffusion process in~\eqref{eq:diffusionjsq}, Theorem~\ref{th:statail} in Chapter~\ref{chap:jsqdiffusion} establishes that $\bQ_1(\infty)$ has a Gaussian tail, and the tail exponent is uniformly bounded by constants which do not depend on $\beta$, whereas $\bQ_2(\infty)$ has an exponentially decaying tail, and the coefficient in the exponent is linear in $\beta$.
More precisely, for any $\beta>0$ there exist positive constants $C_1, C_2, D_1, D_2$ not depending on $\beta$ and positive constants $C^{l}(\beta)$, $C^{u}(\beta)$, $D^{l}(\beta)$, $D^{u}(\beta)$, $C_R(\beta)$, $D_R(\beta)$ depending only on $\beta$ such that 
\begin{equation}
\begin{split}
C^{l}(\beta)e^{-C_1x^2} \le \pi(\bQ_1(\infty) < -x) \le C^{u}(\beta)e^{-C_2x^2}, \ \ x \ge C_R(\beta)\\
D^{l}(\beta)e^{-D_1\beta y} \le \pi(\bQ_2(\infty) > y) \le D^{u}(\beta)e^{-D_2\beta y}, \ \ y \ge D_R(\beta).
\end{split}
\end{equation}
It is further shown in Theorem~\ref{th:lil} that there exists a positive constant $\mathcal{C^*}$ not depending on $\beta$ such that almost surely along any sample path
\begin{equation}\label{eq:fluc-width}
\begin{split}
-2\sqrt{2} &\le \liminf_{t \rightarrow \infty} \frac{\bQ_1(t)}{\sqrt{\log t}} \le -1,\\
\frac{1}{\beta} &\le \limsup_{t \rightarrow \infty} \frac{\bQ_2(t)}{\log t} \le \frac{2}{\mathcal{C^*} \beta}.
\end{split}
\end{equation}
Equation~\eqref{eq:fluc-width} captures the explicit dependence on $\beta$ of the width of the fluctuation window of $\bQ_1$ and $\bQ_2$. 
Specifically, note that the width of fluctuation of $\bQ_1$ does not depend on the value of $\beta$, whereas that of $\bQ_2$ 
is linear in $\beta^{-1}$. 
\begin{remark}\label{rem:comp2-intro}\normalfont
It is worth mentioning that in case of M/M/N systems in the Halfin-Whitt heavy-traffic regime~\cite[Theorem 2]{HW81}, the centered and scaled total number of tasks in the system $(\bar{S}^N(t)-N)/\sqrt{N}$ converges weakly to a diffusion process $\{\bar{S}(t)\}_{t\geq 0}$ having infinitesimal generator 
$A = (\sigma^2(x)/2)(\dif^2/\dif x^2)+m(x)(\dif/\dif x)$
with
\[m(x) = \begin{cases}
-\beta & \mbox{ if } x >0\\
-(x+\beta) & \mbox{ if } x \le 0
\end{cases}\qquad \text{and}\qquad\sigma^2(x) = 2.\]
Note that since this is a simple combination of a Brownian motion with a negative drift (when all servers are fully occupied) and an Ornstein Uhlenbeck process (when there are idle servers), the steady-state distribution $\bar{S}(\infty)$ can be {\em computed explicitly}, and is a combination of an exponential distribution (from the Brownian motion with a negative drift) and a Gaussian distribution (from the OU process).
Although in terms of tail asymptotics, $S(\infty) = \bQ_1(\infty) + \bQ_2(\infty)$ behaves somewhat similarly to that for the centered and scaled total number of tasks in the corresponding M/M/$N$ system, there are some fundamental differences between the two processes, which not only make the analysis of the JSQ diffusion much harder, but also lead to several completely different qualitative properties.
\begin{enumerate}[{\normalfont (i)}]
\item Observe that in case of M/M/N systems, whenever there are some waiting tasks (equivalent to $Q_2$ being positive in our case), the queue length has a constant negative drift towards zero. This leads to the exponential upper tail of $\bar{S}(\infty)$, by comparing with the stationary distribution of a reflected Brownian motion with constant negative drift. In our case, the rate of decrease of $Q_2$ is always proportional to itself, which makes it somewhat counter-intuitive that its stationary distribution has an exponential tail.
\item Further, from \eqref{eq:diffusionjsq}, $Q_2$ {\em never} hits zero. 
Thus, in the steady state, there is no mass at $Q_2=0$, and the system always has waiting tasks. 
This is in sharp contrast to the M/M/N case, where the system has no waiting tasks with positive probability in steady state. 
\item In the M/M/N system, given that a task faces a non-zero wait, \emph{the steady-state waiting time is of order $1/\sqrt{N}$} whereas in the JSQ case it is of \emph{constant order} (the time till the service of the task ahead of it in its queue finishes). 
Moreover, in the JSQ case, it is easy to see that $Q_1$ (the limit of the scaled number of idle servers) spends zero time at the origin, i.e., in steady state  the fraction of arriving tasks that find all servers busy vanishes in the large-$N$ limit.
Consequently, JSQ achieves an  asymptotically vanishing steady-state probability of non-zero wait (in fact, this is of order $1/\sqrt{N}$, see~\cite{Braverman18}). 
This is another sharp contrast with the M/M/N case, where the asymptotic steady-state probability of non-zero wait is strictly positive.
\item In the M/M/N system, the number of idle servers can be non-zero only when the number of waiting tasks is zero.
Thus, the dynamics of both the number of idle servers and the number of waiting tasks are completely captured by the one-dimensional process $\bar{S}^N$ and by the one-dimensional diffusion  $\bar{S}$ in the limit. 
But in the JSQ case, $Q_2$ is never zero, and the dynamics of $(Q_1, Q_2)$ are truly two-dimensional (although the diffusion is non-elliptic) with $Q_1$ and $Q_2$ interacting with each other in an intricate manner.
\end{enumerate}
\end{remark}

\subsection{JSQ(d) policies in heavy-traffic regime}

Finally, we briefly discuss the behavior of JSQ($d$) policies
with a fixed value of~$d$ in the Halfin-Whitt heavy-traffic
regime~\eqref{eq:HW}.
While a complete characterization of the occupancy process
for fixed~$d$ has remained elusive so far, significant partial results
were recently obtained by Eschenfeldt \& Gamarnik~\cite{EG16}.
In order to describe the transient asymptotics,
introduce the following rescaled processes
\begin{equation}
\label{eq:HTtransient}
\bQ_i^N(t) := \frac{N-Q_i^N(t)}{\sqrt{N}}, \quad i = 1, 2,\ldots.
\end{equation}
Note that in contrast with~\eqref{eq:diffscale}, in~\eqref{eq:HTtransient}
\emph{all} components are centered by $N$. 
We also note that in~\cite{EG16} a considerably more general class of heavy-traffic regimes
have been considered (not just the Halfin-Whitt regime).
Then {\em for suitable initial states, \cite[Theorem~2]{EG16}
establishes that on any finite time interval, $\bar{\QQ}^N(\cdot)$
converges weakly to a deterministic system $\bar{\QQ}(\cdot)$ that
satisfies the following system of ODEs
\begin{equation}
\dif \bQ_i(t) = - d (\bQ_i(t) - \bQ_{i-1}(t)) + \bQ_{i+1}(t) - \bQ_i(t),
\quad i = 1, 2, \ldots,
\end{equation}
with the convention that $\bQ_0(t) \equiv 0$.}
It is noteworthy that the scaled occupancy process loses its
diffusive behavior for fixed~$d$.
It is further shown in~\cite{EG16} that with high probability the
steady-state fraction of queues with length at least
$\log_d(\sqrt{N} / \beta) - \omega(1)$ tasks approaches unity,
which in turn implies that with high probability the steady-state delay
is {\em at least} $\log_d(\sqrt{N} / \beta) - O(1)$ as $N \to \infty$.
The diffusion approximation of the JSQ($d$) policy in the Halfin-Whitt regime~\eqref{eq:HW}, starting from a different initial scaling, has been studied by Budhiraja \& Friedlander~\cite{BF17}.

In the work of Ying~\cite{Ying17} a broad framework involving Stein's
method was introduced to analyze the rate of convergence of the
stationary distribution in a heavy-traffic regime where
$\frac{N - \lambda(N)}{\eta(N)} \to \beta > 0$ as $N \to \infty$,
with $\eta(N)$ a positive function diverging to infinity as $N \to \infty$.
Note that the case $\eta(N)= \sqrt{N}$ corresponds to the Halfin-Whitt
heavy-traffic regime~\eqref{eq:HW}.
Using this framework, it was proved that when $\eta(N) = N^\alpha$
with some $\alpha>0.8$, 
\begin{equation}
\label{eq:po2ht}
\mathbb{E}\Big(\sum_{i=1}^\infty \Big|q_i^N(\infty) - q_i^\star\Big|\Big)\leq\frac{1}{N^{2\alpha -1-\xi}},\qquad\mbox{where}\qquad q_i^\star = \Big(\frac{\lambda(N)}{N}\Big)^{2^k-1},
\end{equation}
and $\xi>0$ is an arbitrarily small constant.
Equation~\eqref{eq:po2ht} not only shows that the
stationary occupancy measure asymptotically concentrates at $\qq^\star$,
but also provides the rate of convergence.

\section{Universality of JSQ(d) policies}
\label{univ}


In this section we will further explore the trade-off between delay
performance and communication overhead as a function of the diversity
parameter~$d$, in conjunction with the relative load.
The latter trade-off will be examined in an asymptotic regime where
not only the total task arrival rate $\lambda(N)$ grows with~$N$,
but also the diversity parameter depends on~$N$,
and we write $d(N)$ to explicitly reflect that.
We will specifically investigate what growth rate of $d(N)$ is required,
depending on the scaling behavior of $\lambda(N)$,
in order to asymptotically match the optimal performance of the JSQ
policy and achieve a zero mean waiting time in the limit.
The results presented in the remainder of the section are discussed in greater detail in Chapter~\ref{chap:univjsqd}.

\begin{theorem}{\normalfont (Fluid limit for JSQ($d(N)$))}
\label{fluidjsqd}
If $d(N)\to\infty$ as $N\to\infty$, then the fluid limit of the
JSQ$(d(N))$ scheme coincides with that of the ordinary JSQ policy,
and in particular, is given by the dynamical system in~\eqref{eq:fluid-intro}.
Consequently, the stationary occupancy states converge to the unique
fixed point as in~\eqref{eq:fpjsq}.
\end{theorem}

\begin{theorem}{\normalfont (Diffusion limit for JSQ($d(N)$))}
\label{diffusionjsqd}
If $d(N) /( \sqrt{N} \log N)\to\infty$, then for suitable initial
conditions the weak limit of the sequence of processes
$\big\{\bar{\QQ}^{ d(N)}(t)\big\}_{t \geq 0}$ coincides with that
of the ordinary JSQ policy, and in particular, is given by the system
of SDEs in~\eqref{eq:diffusionjsq}.
\end{theorem}

The above universality properties indicate that the JSQ overhead can
be lowered by almost a factor O($N$) and O($\sqrt{N} / \log N$)
while retaining fluid- and diffusion-level optimality, respectively.
In other words, Theorems~\ref{fluidjsqd} and~\ref{diffusionjsqd}
 reveal that it is sufficient for $d(N)$ to grow at any rate
and faster than $\sqrt{N} \log N$ in order to observe similar scaling
benefits as in a pooled system with $N$~parallel single-server queues
on fluid scale and diffusion scale, respectively.
The stated conditions are in fact close to necessary, in the sense
that if $d(N)$ is uniformly bounded and $d(N) /( \sqrt{N} \log N) \to 0$
as $N \to \infty$, then the fluid-limit and diffusion-limit paths
of the system occupancy process under the JSQ($d(N)$) scheme differ
from those under the ordinary JSQ policy.
In particular, if $d(N)$ is uniformly bounded, the mean steady-state
delay does not vanish asymptotically as $N \to \infty$.

It is worth mentioning that from a high level, conceptually related scaling limits were examined using quite different
techniques by Dieker and Suk~\cite{DS15} in a dynamic
scheduling framework (as opposed to the load balancing context).

\begin{remark}\normalfont
One implication of Theorem~\ref{fluidjsqd} is that in the subcritical regime any growth rate of $d(N)$ is enough to achieve an asymptotically vanishing steady-state probability of wait. 
This result is complemented by recent results of Liu and Ying~\cite{LY18} and Brightwell et al.~\cite{BFL18}, where the steady-state analysis is extended to the heavy-traffic regime.
Specifically, it is established in~\cite{LY18} that when the system load of the $N$-th system scales as $N-N^{\alpha}$ with $\alpha\in (0,1/2)$ (i.e., the system is in heavy traffic, but the load is lighter than that in the Halfin-Whitt regime), the steady-state probability of wait for the JSQ($d(N)$) policy with $d(N)\geq N^{1-\alpha}\log N$ vanishes as $N\to\infty$. 
The results of~\cite{BFL18} imply that when $\lambda(N) = N - N^{\alpha}$ and $d(N) = \lfloor N^{\beta}\rfloor$ with $\alpha,\beta\in (0,1]$, $k = \lceil (1-\alpha)/\beta\rceil$, and $2\alpha+\beta (k-1)>1$, with probability tending to 1 as $N\to\infty$, the proportion of queues with queue length equal to $k$ is at least $1-2N^{-1+\alpha + (k-1)\beta}$ and there are no longer queues.
It is important to note that in contrast to the latter papers, the result stated in Theorem~\ref{diffusionjsqd} considers behavior of the system on diffusion scale (and described in terms of a limiting diffusion process).
\end{remark}

\paragraph{High-level proof idea.}
The proofs of both Theorems~\ref{fluidjsqd} and~\ref{diffusionjsqd}
rely on a stochastic coupling construction to bound the difference
in the queue length processes between the JSQ policy and a scheme
with an arbitrary value of $d(N)$.  
This coupling is then exploited to obtain the fluid and diffusion
limits of the JSQ($d(N)$) policy, along with the associated fixed point,
under the conditions stated in Theorems~\ref{fluidjsqd}
and~\ref{diffusionjsqd}.

\begin{figure}
\begin{center}
\includegraphics[scale=.50]{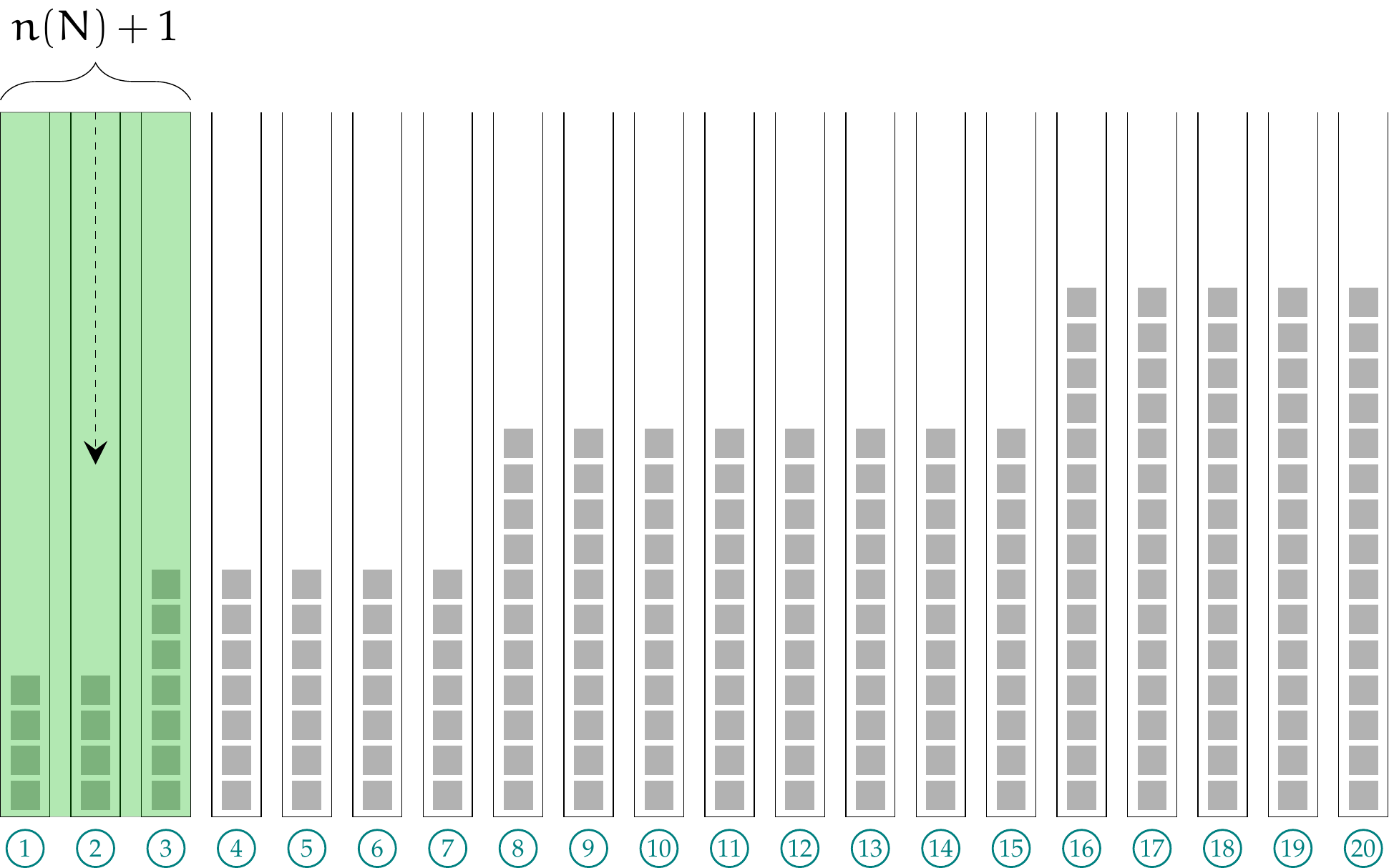}
  \caption{High-level view of the $\CJSQ(n(N))$ class of schemes, where as in Figure~\ref{figB}, the servers are arranged in nondecreasing order of their queue lengths, and the arrival must be assigned through the green left tunnel.}
  \label{fig:sfigCJSQ}
\end{center}
\end{figure}

\begin{figure}
\begin{center}
\includegraphics[scale=1.3]{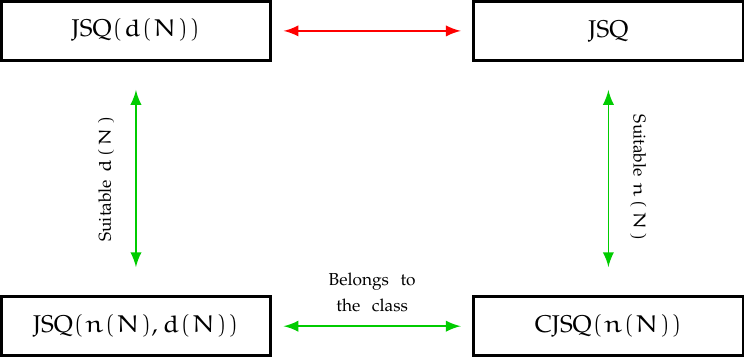}
  \caption{The asymptotic equivalence structure is depicted for various intermediate load balancing schemes to facilitate the comparison between the JSQ$(d(N))$ scheme and the ordinary JSQ policy.}
  \label{fig:sfigRelation}
\end{center}
\end{figure}

A direct comparison between the JSQ$(d(N))$ scheme and the ordinary JSQ policy is not straightforward, which is why the $\CJSQ(n(N))$ class of schemes is introduced as an intermediate scenario to establish the universality result.
Just like the JSQ$(d(N))$ scheme, the schemes in the class $\CJSQ(n(N))$ may be thought of as ``sloppy'' versions of the JSQ policy, in the sense that tasks are not necessarily assigned to a server with the shortest queue length but to one of the $n(N)+1$
lowest ordered servers, as graphically illustrated in Figure~\ref{fig:sfigCJSQ}.
In particular, for $n(N)=0$, the class only includes the ordinary JSQ policy. 
Note that the JSQ$(d(N))$ scheme is guaranteed to identify the lowest ordered server, but only among a randomly sampled subset of $d(N)$ servers.
In contrast, a scheme in the $\CJSQ(n(N))$ class only guarantees that
one of the $n(N)+1$ lowest ordered servers is selected, but 
across the entire pool of $N$ servers. 
We will show that for sufficiently small $n(N)$, any scheme from the class $\CJSQ(n(N))$ is still `close' to the ordinary JSQ policy. 
We will further prove that for sufficiently large $d(N)$ relative to $n(N)$ we can construct a scheme
called JSQ$(n(N),d(N))$, belonging to the $\CJSQ(n(N))$ class, which differs `negligibly' from the JSQ$(d(N))$ scheme. 
Therefore,  for a `suitable' choice of $d(N)$ the idea is to produce a `suitable' $n(N)$.
This proof strategy is schematically represented in Figure~\ref{fig:sfigRelation}.

In order to prove the stochastic comparisons among the various schemes,
the many-server system is described as an ensemble of stacks,
in a way that two different ensembles can be ordered.
This stack formulation has also been considered in the literature
for establishing the stochastic optimality properties of the JSQ
policy \cite{towsley,Towsley95,Towsley1992}.
In Remark~\ref{rem:literatureJSQ} we will compare and contrast the
various stochastic comparison techniques.
In this formulation, at each step, items are added or removed
(corresponding to an arrival or departure) according to some rule. 
From a high level, it is then shown that if two systems follow some
specific rules, then at any step, the two ensembles maintain some kind
of deterministic ordering. 
This deterministic ordering turns into an almost sure ordering in the
probability space constructed by a specific coupling.
In what follows, each server along with its queue is thought of as
a stack of items, and the stacks are always considered to be arranged
in non-decreasing order of their heights. 
The ensemble of stacks then represents the empirical CDF of the queue
length distribution, and the $i^{\mathrm{th}}$ horizontal bar
corresponds to $Q_i^\Pi$ (for some task assignment scheme~$\Pi$),
as depicted in Figure~\ref{figB}.
For the sake of full exposure, we will describe the coupling
construction in the scenario when the buffer capacity~$B$ at each
stack can possibly be finite.
If $B < \infty$ and an arriving item happens to land on a stack which
already contains $B$ items, then the item is discarded, and is added
to a special stack $L^\Pi$ of discarded items, where it stays forever.

Any two ensembles $\cA$ and $\cB$, each having $N$~stacks
and a maximum height $B$ per stack, are said to follow
Rule($n_{\cA},n_{\cB},k$) at some step, if either an item is removed
from the $k^{\mathrm{th}}$ stack in both ensembles (if nonempty), 
or an item is added to the $n_{\cA}^{\mathrm{th}}$ stack in ensemble
$\cA$ and to the $n_{\cB}^{\mathrm{th}}$ stack in ensemble $\cB$. 

\begin{proposition}
\label{prop:det-ord-intro}
For any two ensembles of stacks $\cA$ and $\cB$, if \emph{Rule}$(n_{\cA},n_{\cB},k)$ is followed at each step for some value of $n_{\cA}$, $n_{\cB}$, and $k$, with $n_{\cA}\leq n_{\cB}$ (the value of $n_{\cA}$, $n_{\cB}$, and $k$ might differ from step to step), then the following ordering is always preserved: for all $m\leq B$,
\begin{equation}\label{eq:det ord-intro}
\sum_{i=m}^B Q_i^{\cA} +L^{\cA}\leq \sum_{i=m}^B Q_i^{\cB} +L^{\cB}.
\end{equation}
\end{proposition} 

This proposition says that, while adding the items to the ordered stacks,
if we ensure that in ensemble $\cA$ the item is always placed to the
left of that in ensemble $\cB$, and if the items are removed from the
same ordered stack in both ensembles, then the aggregate size of the
$B-m+1$ highest horizontal bars as depicted in Figure~\ref{figB}
plus the cumulative number of discarded items is no larger in $\cA$
than in $\cB$ throughout.

\paragraph{Another type of sloppiness.} 
Recall that $\CJSQ(n(N))$ contains all schemes that assign
incoming tasks by some rule to any of the $n(N)+1$ lowest ordered servers. 
Let MJSQ$(n(N))$ be a particular scheme that always assigns incoming
tasks to precisely the $(n(N)+1)^{\mathrm{th}}$ ordered server. 
Notice that this scheme is effectively the JSQ policy when the system
always maintains $n(N)$ idle servers, or equivalently, uses only
$N-n(N)$ servers, and $\MJSQ(n(N)) \in \CJSQ(n(N))$. 
For brevity, we will often suppress $n(N)$ in the notation where it is
clear from the context.
We call any two systems \emph{S-coupled}, if they have synchronized
arrival clocks and departure clocks of the $k^{\mathrm{th}}$ longest queue,
for $1 \leq k \leq N$ (`S' in the name of the coupling stands for `Server').
Consider three S-coupled systems following respectively the JSQ policy,
any scheme from the class $\CJSQ$, and the $\MJSQ$ scheme. 
Recall that $Q_i^\Pi(t)$ is the number of servers with at least $i$~tasks
at time~$t$ and $L^\Pi(t)$ is the total number of lost tasks up to time~$t$,
for the schemes $\Pi=$ JSQ, $\CJSQ$, $\MJSQ$. 
The following proposition provides a stochastic ordering for any scheme
in the class CJSQ with respect to the ordinary JSQ policy and the MJSQ scheme. 

\begin{proposition}
\label{prop:stoch-ord-intro}
For any fixed $m \geq 1$,
\begin{enumerate}[{\normalfont(i)}] 
\item\label{item:jsq-cjsq-intro} $
\left\{\sum_{i=m}^{B} Q_i^{\JSQ}(t) + L^{\JSQ}(t)\right\}_{t \geq 0} \leq_{\mathrm{st}}
\left\{\sum_{i=m}^{B} Q_i^{\CJSQ}(t) + L^{\CJSQ}(t)\right\}_{t \geq 0},
$
\item\label{item:cjsq-mjsq-intro} $\left\{\sum_{i=m}^{B} Q_i^{\CJSQ}(t) +
L^{\CJSQ}(t)\right\}_{t \geq 0}\leq_{\mathrm{st}}
\left\{\sum_{i=m}^{B} Q_i^{\MJSQ}(t) + L^{\MJSQ}(t)\right\}_{t \geq 0},$
\end{enumerate}
provided the inequalities hold at time $t=0$.
\end{proposition}

The above proposition has the following immediate corollary,
which will be used to prove bounds on the fluid and the diffusion scale.

\begin{corollary}
\label{cor:bound-intro}
In the joint probability space constructed by the S-coupling of the three systems under respectively \emph{JSQ}, \emph{MJSQ}, and any scheme from the class \emph{CJSQ}, the following ordering is preserved almost surely throughout the sample path: for any fixed $m \geq 1$
\begin{enumerate}[{\normalfont(i)}]
\item $Q_m^{\CJSQ}(t) \geq \sum_{i=m}^{B} Q_i^{\JSQ}(t) -
\sum_{i=m+1}^{B} Q_i^{\MJSQ}(t) + L^{\JSQ}(t) - L^{\MJSQ}(t)$,
\item $Q_m^{\CJSQ}(t) \leq \sum_{i=m}^{B} Q_i^{\MJSQ}(t) -
\sum_{i=m+1}^{B} Q_i^{\JSQ}(t) + L^{\MJSQ}(t)-L^{\JSQ}(t),$
\end{enumerate}
provided the inequalities hold at time $t=0$.
\end{corollary}

\begin{remark}
\label{rem:literatureJSQ}
\normalfont
Note that $\sum_{i=1}^{B} \min\big\{Q_i, k\big\}$ represents the
aggregate size of the rightmost $k$~stacks, i.e., the $k$~longest queues.
Using this observation, the stochastic majorization property of the
JSQ policy as stated in \cite{towsley,Towsley95,Towsley1992} can be
shown following similar arguments as in the proof
of Proposition~\ref{prop:stoch-ord-intro}.
Conversely, the stochastic ordering between the JSQ policy and the
MJSQ scheme presented in Proposition~\ref{prop:stoch-ord-intro} can also be
derived from the weak majorization arguments developed
in \cite{towsley,Towsley95,Towsley1992}. 
But it is only through the stack arguments developed 
in Chapter~\ref{chap:univjsqd}
as described above, that the results could be extended to compare
any scheme from the class CJSQ with the scheme MJSQ as stated in Proposition~\ref{prop:stoch-ord-intro}~\eqref{item:cjsq-mjsq-intro}.
\end{remark}

\paragraph{Comparing two arbitrary schemes.}
To analyze the JSQ$(d(N))$ scheme, we need a further stochastic
comparison argument.
Consider two S-coupled systems following schemes $\Pi_1$ and $\Pi_2$.
Fix a specific arrival epoch, and let the arriving task join the
$n_{\Pi_i}^{\mathrm{th}}$ ordered server in the $i^{\mathrm{th}}$
system following scheme $\Pi_i$, $i = 1, 2$ (ties can be broken
arbitrarily in both systems). 
We say that at a specific arrival epoch the two systems
\emph{differ in decision}, if $n_{\Pi_1} \neq n_{\Pi_2}$,
and denote by $\Delta_{\Pi_1,\Pi_2}(t)$ the cumulative number of times
the two systems differ in decision up to time~$t$.

\begin{proposition}
\label{prop:stoch-ord2-intro}
For two S-coupled systems under schemes $\Pi_1$ and $\Pi_2$ the
following inequality is preserved almost surely:
\begin{equation}
\label{eq:stoch-ord2-intro}
\sum_{i=1}^{B} |Q_i^{\Pi_1}(t) - Q_i^{\Pi_2}(t)| \leq
2 \Delta_{\Pi_1,\Pi_2}(t) \qquad \forall\ t \geq 0,
\end{equation}
provided the two systems start from the same occupancy state at $t=0$,
i.e., $Q_i^{\Pi_1}(0) = Q_i^{\Pi_2}(0)$ for all $i = 1, 2, \ldots, B$.
\end{proposition}

\paragraph{A bridge between two types of sloppiness.}
We will now introduce the JSQ$(n,d)$ scheme with $n, d \leq N$,
which is an intermediate blend between the CJSQ$(n)$ schemes
and the JSQ$(d)$ scheme.
We now specify the JSQ$(d,n)$ scheme.
At its first step, just as in the JSQ$(d)$ scheme, it first chooses
the shortest of $d$~random candidates but only sends the arriving task
to that server's queue if it is one of the $n+1$ shortest queues. 
If it is not, then at the second step it picks any of the $n+1$ shortest
queues uniformly at random and then sends the task to that server's queue. 
Note that by construction, JSQ$(d,n)$ is a scheme
in CJSQ$(n)$.
Consider two S-coupled systems with a JSQ$(d)$ and a JSQ$(n,d)$ scheme.
Assume that at some specific arrival epoch, the incoming task is
dispatched to the $k^{\mathrm{th}}$ ordered server in the system
under the JSQ($d$) scheme.
If $k \in \{1, 2, \ldots, n+1\}$, then the system under the JSQ$(n,d)$
scheme also assigns the arriving task to the $k^{\mathrm{th}}$ ordered
server. 
Otherwise, it dispatches the arriving task uniformly at random among
the first $(n+1)$ ordered servers.

The next proposition provides a bound on the number of times these two
systems differ in decision on any finite time interval.
For any $T \geq 0$, let $A(T)$ and $\Delta(T)$ be the total number
of arrivals to the system and the cumulative number of times that the
JSQ($d$) scheme and JSQ$(n,d)$ scheme differ in decision up to time~$T$.

\begin{proposition}
\label{prop:differ-intro}
For any $T\geq 0$, and $M>0,$
\[
\Pro{\Delta(T) \geq
M \given A(T)} \leq \frac{A(T)}{M} \left(1-\frac{n}{N}\right)^d.
\]
\end{proposition}

\noindent
\emph{Proof sketch of  Theorem~\ref{fluidjsqd}.}
The proof of Theorem~\ref{fluidjsqd} uses the S-coupling
and consists of three main steps:
\begin{enumerate}[{\normalfont (i)}]
\item First it is shown that if $n(N) / N \to 0$ as $N \to \infty$,
then the MJSQ$(n(N))$ scheme has the same fluid limit as the ordinary
JSQ policy.
\item Then application of Corollary~\ref{cor:bound-intro} proves that
as long as $n(N) / N \to 0$, \emph{any} scheme from the class
$\CJSQ(n(N))$ has the same fluid limit as the ordinary JSQ policy.
\item Next, Propositions~\ref{prop:stoch-ord2-intro} and~\ref{prop:differ-intro}
are used to establish that if $d(N) \to \infty$, then for \emph{some}
$n(N)$ with $n(N) / N \to 0$, both the JSQ$(d(N))$ scheme and the
JSQ$(n(N),d(N))$ scheme have the same fluid limit.
The proposition then follows by observing that the JSQ$(n(N),d(N))$
scheme belongs to the class $\CJSQ(n(N))$.\hfill\qed
\end{enumerate}

The proof of Theorem~\ref{diffusionjsqd} follows the same arguments, but
uses the condition $n(N)/\sqrt{N}\to 0$ (instead of $n(N)/N\to 0$) in Steps (i) and (ii),
and the condition $d(N)/(\sqrt{N}\log(N))\to\infty$ (instead of $d(N)\to\infty$) in Step (iii).


\paragraph*{Extension to batch arrivals.}
We now consider an extension of the model in which tasks arrive in batches.
We assume that the batches arrive as a Poisson process of rate
$\lambda(N) / \ell(N)$, and have fixed size $\ell(N) > 0$,
so that the effective total task arrival rate remains $\lambda(N)$. 
We will show that even for arbitrarily slowly growing batch size,
fluid-level optimality can be achieved with $O(1)$ communication
overhead per task. 
For that, we define the JSQ($d(N)$) scheme adapted to batch arrivals:
When a batch of size $\ell(N)$ arrives, the dispatcher samples
$d(N) \geq \ell(N)$ servers without replacement, and assigns the
$\ell(N)$ tasks to the $\ell(N)$ servers with the smallest queue
length among the sampled servers.

\begin{theorem}{\normalfont (Batch arrivals)}
\label{th:batch-intro}
Consider the batch arrival scenario with growing batch size $\ell(N)\to\infty$ and $\lambda(N)/N\to\lambda<1$ as $N\to\infty$. For the JSQ$(d(N))$ scheme with $d(N)\geq \ell(N)/(1-\lambda-\varepsilon)$ for any fixed $\varepsilon>0$, if $q^{\sss d(N)}_1(0)\to q_1(0)\leq \lambda$, and $q_i^{\sss d(N)}(0)\to 0$ for all $i\geq 2$, then the weak limit of the sequence of processes $\big\{\qq^{\sss d(N)}(t)\big\}_{t\geq 0}$ coincides with that of the ordinary JSQ policy, and in particular, is given by the system in~\eqref{eq:fluid-intro}.  
\end{theorem}

Observe that for a fixed $\varepsilon>0$, the communication overhead
per task is on average given by $(1-\lambda-\varepsilon)^{-1}$
which is $O(1)$.
Thus Theorem~\ref{th:batch-intro} ensures that in case of batch arrivals
with growing batch size, fluid-level optimality can be achieved
with $O(1)$ communication overhead per task.
The result for the fluid-level optimality in stationarity can also be
obtained indirectly by exploiting the fluid-limit result in~\cite{YSK15}.
Specifically, it can be deduced from the result in~\cite{YSK15} that
for batch arrivals with growing batch size, the JSQ$(d(N))$ scheme
with suitably growing $d(N)$ yields the same fixed point of the fluid
limit as described in~\eqref{eq:fpjsq}.

\section{Blocking and infinite-server dynamics}
\label{bloc}

The basic scenario that we have focused on so far involved
single-server queues.
In this section we turn attention to a system with parallel server
pools, each with $B$~servers, where $B$ can possibly be infinite.
As before, tasks arrive at a single dispatcher and must immediately
be forwarded to one of the server pools, but also directly start
execution or be discarded otherwise.
The execution times are assumed to be exponentially distributed,
and do not depend on the number of other tasks receiving service
simultaneously, but the experienced performance (e.g.~in terms
of received throughput or packet-level delay) does degrade
in a convex manner with an increasing number of concurrent tasks.
In order to distinguish it from the single-server queueing dynamics
as considered earlier, the current scenario will henceforth be referred
to as the `infinite-server dynamics'.
These characteristics pertain for instance to video streaming sessions
and various interactive applications.
In contrast to elastic data transfers or computing-intensive jobs,
the duration of such sessions is hardly affected by the number
of contending service requests.
The perceived performance in terms of video quality or packet-level
latency however strongly varies with the number of concurrent tasks,
creating an incentive to distribute the incoming tasks across the
various server pools as evenly as possible.

As it turns out, the JSQ policy has similar stochastic optimality
properties as in the case of single-server queues, and in particular
stochastically minimizes the cumulative number of discarded tasks
\cite{STC93,J89,M87,MS91}.
However, the JSQ policy also suffers from a similar scalability issue
due to the excessive communication overhead in large-scale systems,
which can be mitigated through JSQ($d$) policies.
Results of Turner~\cite{T98} and recent papers by Mukhopadhyay
{\em et al.}~\cite{MKMG15,MMG15}, Karthik {\em et al.}~\cite{KMM17},
and Xie {\em et al.}~\cite{XDLS15} indicate that JSQ($d$) policies
provide similar ``power-of-choice'' gains for loss probabilities.
It may be shown though that the optimal performance of the JSQ policy
cannot be matched for any fixed value of~$d$.

Motivated by these observations, we explore the trade-off between
performance and communication overhead for infinite-server dynamics.
We will demonstrate that the optimal performance of the JSQ policy
can be asymptotically retained while drastically reducing the
communication burden, mirroring the universality properties described
in Section~\ref{univ} for single-server queues.
The results presented in the remainder of the section 
along with their full proofs are contained in Chapter~\ref{chap:asympjsqd}.

\subsection{Fluid limit for JSQ policy}
\label{ssec:jsqfluid-infinite}

As in Subsection~\ref{ssec:jsqfluid}, for any fluid state $\qq \in \cS$,
denote by $m(\qq) = \min\{i: q_{i + 1} < 1\}$ the minimum queue length
among all servers.
Now if $m(\qq)=0$, then define $p_0(m(\qq))=1$ and $p_i(m(\qq))=0$
for all $i=1,2,\ldots$. 
Otherwise, in case $m(\qq)>0$, define
\begin{equation}
\label{eq:fluid-prob-infinite}
p_{i}(\qq) =
\begin{cases}
\min\big\{m(\qq)(1 - q_{m(\qq) + 1})/\lambda,1\big\} & \quad \mbox{ for }
\quad i=m(\qq)-1, \\
1 - p_{ m(\qq) - 1}(\qq) & \quad \mbox{ for } \quad i=m(\qq), \\
0 & \quad \mbox{ otherwise.}
\end{cases}
\end{equation}
{\em Any weak limit of the sequence of processes $\{\qq^N(t)\}_{t \geq 0}$
is given by the deterministic system $\{\qq(t)\}_{t \geq 0}$ satisfying
the following of differential equations
\begin{equation}
\label{eq:fluid-infinite}
\frac{\dif^+ q_i(t)}{\dif t} =
\lambda p_{i-1}(\qq(t)) - i (q_i(t) - q_{i+1}(t)),
\quad i = 1, 2, \dots,
\end{equation}
where $\dif^+/\dif t$ denotes the right-derivative.}

Equations~\eqref{eq:fluid-prob-infinite} and \eqref{eq:fluid-infinite}
are to be contrasted with Equations~\eqref{eq:fluid-gen-intro}
and~\eqref{eq:fluid-intro}.
While the form of the evolution equations~\eqref{eq:fluid-infinite}
of the limiting dynamical system remains similar to~\eqref{eq:fluid-intro},
the rate of decrease of $q_i$ is now $i (q_i - q_{i+1})$,
reflecting the infinite-server dynamics.

Let $K := \lfloor \lambda \rfloor$ and $f := \lambda - K$ denote the
integral and fractional parts of~$\lambda$, respectively.
It is easily verified that, assuming $\lambda<B$, the unique fixed point
of the dynamical system in~\eqref{eq:fluid-infinite} is given by
\begin{equation}
\label{eq:fixed-point-infinite}
q_i^\star = \left\{\begin{array}{ll} 1 & i = 1, \dots, K \\
f & i = K + 1 \\
0 & i = K + 2, \dots, B, \end{array} \right.
\end{equation}
and thus $\sum_{i=1}^{B} q_i^\star = \lambda$.
This is consistent with the results in Mukhopadhyay
{\em et al.}~\cite{MKMG15,MMG15} and Xie {\em et al.}~\cite{XDLS15}
for fixed~$d$, where taking $d \to \infty$ yields the same fixed point.
However, the results in \cite{MKMG15,MMG15,XDLS15} for fixed~$d$
cannot be directly used to handle joint scalings, and do not yield the
universality of the entire fluid-scaled sample path for arbitrary
initial states as stated in~\eqref{eq:fluid-infinite}.

The fixed point in~\eqref{eq:fixed-point-infinite}, in conjunction
with an interchange of limits argument, indicates that in stationarity
the fraction of server pools with at least $K+2$ and at most $K-1$
active tasks is negligible as $N \to \infty$.

\subsection{Diffusion limit for JSQ policy}
\label{ssec:jsq-diffusion-infinite}

As it turns out, the diffusion-limit results may be qualitatively
different, depending on whether $f = 0$ or $f > 0$,
and we will distinguish between these two cases accordingly.
Observe that for any assignment scheme, in the absence of overflow events,
the total number of active tasks evolves as the number of jobs
in an M/M/$\infty$ system with arrival rate $\lambda(N)$ and unit
service rate, for which the diffusion limit is well-known~\cite{Robert03}.
For the JSQ policy we can establish, for suitable initial conditions,
that the total number of server pools with $K - 2$ or less and $K + 2$
or more tasks is negligible on the diffusion scale.
If $f > 0$, the number of server pools with $K - 1$ tasks is negligible
as well, and the dynamics of the number of server pools with $K$
or $K + 1$ tasks can then be derived from the known diffusion limit
of the total number of tasks mentioned above.
In contrast, if $f = 0$, the number of server pools with $K - 1$ tasks
is not negligible on the diffusion scale, and the limiting behavior is
qualitatively different, but can still be characterized.

\subsubsection{Diffusion-limit results for non-integral
\texorpdfstring{$\boldsymbol{\lambda}$}{lambda}}

We first consider the case $f > 0$, and define $f(N) := \lambda(N) - K N$.
Based on the above observations,
we define the following centered and scaled processes:
\begin{equation}
\begin{split}
\bar{Q}^N_i(t)&=N-Q^N_i(t)\geq 0\quad \mathrm{for}\quad i\leq K-1,\\
\bar{Q}_{K}^N(t)&:=\frac{N-Q_K^N(t)}{\log (N)}\geq 0,\\
\bar{Q}_{K+1}^N(t)&:=\frac{Q^N_{K+1}(t)-f(N)}{\sqrt{N}}\in\R,\\
\bar{Q}^N_i(t)&:=Q^N_i(t)\geq 0\quad \mathrm{for}\quad i\geq K+2.
\end{split}
\end{equation}
\begin{theorem}[{Diffusion limit for JSQ policy; $f>0$}]
\label{th:diffusion-infinite}
Assume $\bar{Q}^N_i(0) $ converges to $\bar{Q}_i(0)$ in probability,
and $\lambda(N)/N \to \lambda > 0$ as $N \to \infty$, then
\begin{enumerate}[{\normalfont(i)}]
\item $\lim\limits_{N\to\infty}\Pro{\sup_{t\in[0,T]}\bar{Q}_{K-1}^N(t)\leq 1}=1$, and $\big\{\bar{Q}^N_i(t)\big\}_{t\geq 0}$ converges weakly to $\big\{\bar{Q}_i(t)\big\}_{t\geq 0}$, where $\bar{Q}_i(t)\equiv 0$, provided 
$\lim_{N\to\infty}\Pro{\bar{Q}_{K-1}^N(0)\leq 1}=1$, and
$\bar{Q}_i^N(0)$ converges to 0 in probability, for $i\leq K-2$.
\item $\big\{\bar{Q}^N_K(t)\big\}_{t\geq 0}$ is a stochastically bounded sequence of processes.
\item $\big\{\bar{Q}^N_{K+1}(t)\big\}_{t\geq 0}$ converges weakly to $\big\{\bar{Q}_{K+1}(t)\big\}_{t\geq 0}$, where $\bar{Q}_{K+1}(t)$ is given by the Ornstein-Uhlenbeck process satisfying the following stochastic differential equation:
$$d\bar{Q}_{K+1}(t)=-\bar{Q}_{K+1}(t)dt+\sqrt{2\lambda}dW(t),$$
where $W(t)$ is the standard Brownian motion,
provided $\bar{Q}_{K+1}^N(0)$  converges to $\bar{Q}_{K+1}(0)$ in probability.
\item For $i\geq K+2$, $\big\{\bar{Q}^N_i(t)\big\}_{t\geq 0}$ converges weakly to $\big\{\bar{Q}_i(t)\big\}_{t\geq 0}$, where $\bar{Q}_i(t)\equiv 0$, provided $\bar{Q}_i^N(0)$ converges to 0 in probability.\\
\end{enumerate}
\end{theorem}

Theorem~\ref{th:diffusion-infinite} implies that for suitable initial
states, for large $N$, there will be almost no server pool with $K-2$
or less tasks and $K+2$ or more tasks on any finite time interval. 
Also, the number of server pools having fewer than $K$ tasks is
of order $\log(N)$, and there are $f(N)+O_P(\sqrt{N})$ server pools
with precisely $K + 1$ active tasks.

\paragraph{High-level proof idea.}
Informally speaking, the proof of Theorem~\ref{th:diffusion-infinite}
proceeds along the following lines of arguments.
Observe that $\sum_{i=1}^{K} (N-Q_i^N(\cdot))$ increases by one at rate 
\[
\sum_{i=1}^{K} i (Q_i(t)-Q_{i+1}(t)) =
\sum_{i=1}^{K} (Q_i(t)-Q_{K+1}(t)) \approx K (1 - f) N,
\] 
which is when there is a departure from some server pool with at most
$K$~active tasks, and if positive, decreases by one at constant rate
$\lambda(N) = (K + f) N + o(N)$, which is whenever there is an arrival.
Thus, $\sum_{i=1}^{K} (N - Q_i^N(\cdot))$ roughly behaves as
a birth-and-death process with birth rate $K (1 - f) N$ and death rate
$(K + f) N$. 
Since $f > 0$, we have $K + f > K (1 - f)$, and on any finite time
interval the maximum of such a birth-and-death process scales as $\log(N)$.

Similar to the argument above, the process
$\sum_{i=1}^{K-1} \bQ_i^N(\cdot)$ increases by one at rate 
\begin{align*}
\sum_{i=1}^{K-1} i (Q^N_i(t)-Q^N_{i+1}(t)) &=
\sum_{i=1}^{K-1} Q^N_i(t) - (K-1) Q_K^N(t) \\
&\leq (K-1) (N-Q_K^N(t)) = O(\log(N)),
\end{align*}
which is when there is a departure from some server pool with at most
$K - 1$ active tasks, and if positive, decreases by one at rate
$\lambda(N)$, which is whenever there is an arrival.
Thus, $\sum_{i=1}^{K-1} \bQ_i^N(\cdot)$ roughly behaves as
a birth-and-death process with birth rate $O(\log(N))$ and death rate $O(N)$.
This leads to the asymptotic result for $\sum_{i=1}^{K-1} \bQ_i^N(\cdot)$,
and in particular for $\bQ_{K-1}^N(\cdot)$.
This completes the proof of Parts~(i) and~(ii)
of Theorem~\ref{th:diffusion-infinite}.

Furthermore, since $\lambda < K + 1$, the number of tasks that are
assigned to server pools with at least $K + 1$ tasks converges to zero
in probability 
This completes the proof of Part (iv)
of Theorem~\ref{th:diffusion-infinite}.
 
Finally, all the above combined also means that on any finite time
interval the total number of tasks in the system behaves with high
probability as the total number of jobs in an M/M/$\infty$ system. 
Therefore with the help of the following diffusion limit result
for the M/M/$\infty$ system in Theorem~\ref{th:robert-book-mmn},
we conclude the proof of Part (iii) of Theorem~\ref{th:diffusion-infinite}.

\begin{theorem}[{\cite[Theorem 6.14]{Robert03}}]
\label{th:robert-book-mmn}
Let $\big\{Y_\infty^N(t)\big\}_{t\geq 0}$ be the total number of jobs
in an M/M/$\infty$ system with arrival rate $\lambda (N)$
and unit-mean service time. 
If $(Y_\infty^N(0) - \lambda(N))/\sqrt{N} \to v \in \R$,
then the process $\big\{\bar{Y}_\infty^N(t)\big\}_{t\geq 0}$, with
$$\bar{Y}^N_{\infty}(t)=\frac{Y^N_\infty(t)-\lambda( N)}{\sqrt{N}},$$
converges weakly to an Ornstein-Uhlenbeck process
$\big\{X(t)\big\}_{t\geq 0}$ described by the SDE
\begin{align*}
X(0) = v, \qquad \dif X(t) &= - X(t) \dif t + \sqrt{2 \lambda} \dif W(t).
\end{align*}
\end{theorem}

\subsubsection{Diffusion-limit results for integral
\texorpdfstring{$\boldsymbol{\lambda}$}{lambda}}
We now turn to the case $f = 0$, and assume that
\begin{equation}
\label{eq:f=0-intro}
\frac{K N - \lambda(N)}{\sqrt{N}} \to \beta \in \R \quad \mbox{ as }
\quad N \to \infty,
\end{equation} 
which can be thought of as an analog of the Halfin-Whitt regime.
As mentioned above, the limiting behavior in this case is
qualitatively different from the case $f > 0$.
Hence, we now consider the following scaled quantities:
\begin{equation}
\begin{split}
\zeta_1^N(t) := \frac{1}{\sqrt{N}}\sum_{i=1}^{K} (N-Q_i^N(t)), \qquad
\zeta_2^N(t) := \frac{Q_{K+1}^N(t)}{\sqrt{N}}.
\end{split}
\end{equation}

\begin{theorem}
\label{th: f=0 diffusion}
Assuming the convergence of initial states, on any finite time interval the process
$\big\{(\zeta_1^N(t),\zeta_2^N(t))\big\}_{t \geq 0}$ converges weakly
to the process $\big\{(\zeta_1(t),\zeta_2(t))\big\}_{t \geq 0}$
governed by the following system of SDEs:
\begin{align*}
\dif\zeta_1(t) &= \sqrt{2K} \dif W(t) - (\zeta_1(t) + K \zeta_2(t)) +
\beta \dif t + \dif V_1(t), \\
\dif\zeta_2(t) &= \dif V_1(t) - (K + 1) \zeta_2(t),
\end{align*}
where $W$ is the standard Brownian motion, and $V_1(t)$ is the unique
non-decreasing process satisfying
\begin{align*}
\int_0^t \ind{\zeta_1(s) \geq 0}\dif V_1(s)=0.
\end{align*}
\end{theorem}

Unlike the $f > 0$ case, the above theorem says that if $f = 0$,
then over any finite time horizon, there will be $O_P(\sqrt{N})$
server pools with fewer than $K$ or more than $K$~active tasks,
and hence most of the server pools have precisely $K$~active tasks.
The proof of Theorem~\ref{th: f=0 diffusion} uses the reflection
argument developed in~\cite{EG15}.

\begin{remark}
\normalfont
Let $Y^N(t)$ denote the total number of tasks in the system at time~$t$.
Note that $Y^N(t) - K N = Z_2^N(t) - Z_1^N(t)$.
Thus, under the assumption in~\eqref{eq:f=0-intro}, the diffusion limit
in Theorem~\ref{th: f=0 diffusion} implies that 
\begin{align*}
\frac{Y^N(\cdot) - \lambda(N)}{\sqrt{N}} =
\frac{Y^N(\cdot) - K N}{\sqrt{N}} + \frac{K N - \lambda(N)}{\sqrt{N}} \dto
\zeta_2(\cdot) - \zeta_1(\cdot) + \beta.
\end{align*}
Writing $X(t) = \zeta_2(t) - \zeta_1(t) - \beta$,
from Theorem~\ref{th: f=0 diffusion}, we see that the process
$\big\{X(t)\big\}_{t \geq 0}$ satisfies 
$$\dif X(t) = - X(t) \dif t - \sqrt{2 K} \dif W(t),$$
which is consistent with the diffusion-level behavior of $Y^N(\cdot)$
stated in Theorem~\ref{th:robert-book-mmn}. 
\end{remark}

\subsection{Universality of JSQ(d) policies in infinite-server dynamics}
\label{ssec:univ-infinite}

As in Section~\ref{univ}, we will now further explore the trade-off
between performance and communication overhead as a function
of the diversity parameter~$d(N)$, in conjunction with the relative load.
We will specifically investigate what growth rate of $d(N)$ is required,
depending on the scaling behavior of $\lambda(N)$, in order to
asymptotically match the optimal performance of the JSQ policy.

\begin{theorem}[ Fluid limit for JSQ($d(N)$) in infinite-server dynamics]
\label{fluidjsqd-infinite}
If $d(N)\to\infty$ as $N\to\infty$, then the fluid limit of the
JSQ$(d(N))$ scheme coincides with that of the ordinary JSQ policy,
and in particular, is given by the dynamical system
in~\eqref{eq:fluid-infinite}. 
Consequently, the stationary occupancy states converge to the unique
fixed point as in~\eqref{eq:fixed-point-infinite}.
\end{theorem}

In order to state the universality result on diffusion scale,
define in case $f > 0$,
\begin{equation}
\begin{split}
\bar{Q}_i^{d(N)}(t) & := \dfrac{N - Q_i^{d(N)}(t)}{\sqrt{N}} \geq 0,
\quad i \leq K, \\  
\bar{Q}_{K+1}^{d(N)}(t) & := \dfrac{Q_{K+1}^{d(N)}(t) - f(N)}{\sqrt{N}}
\in \R, \\ 
\bar{Q}_i^{d(N)}(t) & := \frac{Q_i^{d(N)}(t)}{\sqrt{N}}\geq 0,
\quad \text{ for } \quad i \geq K + 2,
\end{split}
\end{equation}
and otherwise, if $f = 0$,
\begin{equation}
\begin{split}
\hQ_{K-1}^{d(N)}(t) & :=
\sum_{i=1}^{K-1} \dfrac{N - Q_i^{d(N)}(t)}{\sqrt{N}} \geq 0, \\
\hQ_K^{d(N)}(t) & := \dfrac{N - Q_K^{d(N)}(t)}{\sqrt{N}} \geq 0, \\
\hQ_i^{d(N)}(t) & := \dfrac{Q_i^{d(N)}(t)}{\sqrt{N}} \geq 0, \quad
\text{ for } \quad i \geq K + 1.
\end{split}
\end{equation}

\begin{theorem}[\normalfont Diffusion limit for JSQ($d(N)$) in infinite-server dynamics]
\label{diffusionjsqd-infinite}
Assume $d(N) / (\sqrt{N} \log N) \to \infty$. 
Under suitable initial conditions, the following hold.
\begin{enumerate}[{\normalfont (i)}]
\item If $f>0$, then on any finite time interval the process $\bQ_i^{d(N)}(\cdot)$ converges to the zero process for $i\neq K+1$, and the process $\bQ^{d(N)}_{K+1}(\cdot)$ converges weakly to the Ornstein-Uhlenbeck process satisfying the following SDE:
$$d\bar{Q}_{K+1}(t)=-\bar{Q}_{K+1}(t)dt+\sqrt{2\lambda}dW(t),$$
where $W(t)$ is the standard Brownian motion.
\item If $f=0$, then on any finite time interval the process $\hQ_{K-1}^{d(N)}(\cdot)$ converges weakly to the zero process, and the process $(\hQ_{K}^{d(N)}(\cdot), \hQ_{K+1}^{d(N)}(\cdot))$ converges weakly to $(\hQ_{K}(\cdot), \hQ_{K+1}(\cdot))$, described by the unique solution of the following system of SDEs:
\begin{align*}
\dif\hQ_{K}(t) &= \sqrt{2 K} \dif W(t) - (\hQ_K(t) + K \hQ_{K+1}(t)) +
\beta \dif t + \dif V_1(t) \\
\dif\hQ_{K+1}(t) &= \dif V_1(t) - (K + 1) \hQ_{K+1}(t),
\end{align*}
where $W$ is the standard Brownian motion, and $V_1(t)$ is the unique
non-decreasing process satisfying
\begin{align*}
\int_0^t \ind{\hQ_K(s)\geq 0} \dif V_1(s) = 0.
\end{align*}
\end{enumerate}
\end{theorem}

Having established the asymptotic results for the JSQ policy
in Subsections~\ref{ssec:jsqfluid-infinite}
and~\ref{ssec:jsq-diffusion-infinite}, the proofs of the asymptotic
results for the JSQ$(d(N))$ scheme in Theorems~\ref{fluidjsqd-infinite}
and~\ref{diffusionjsqd-infinite} involve establishing a universality
result which shows that the limiting processes for the JSQ$(d(N))$
scheme are `asymptotically equivalent' to those for the ordinary JSQ
policy for suitably large~$d(N)$.
The notion of asymptotic equivalence between different schemes is now
formalized in the next definition. 
\begin{definition}
Let $\Pi_1$ and $\Pi_2$ be two schemes parameterized by the number of server pools $N$. For any positive function $g:\N\to\R_+$, we say that $\Pi_1$ and $\Pi_2$ are `$g(N)$-alike' if there exists a common probability space, such that for any fixed $T\geq 0$, for all $i\geq 1$,
$$\sup_{t\in[0,T]}(g(N))^{-1}|Q_i^{\Pi_1}(t)-Q_i^{\Pi_2}(t)|\pto 0\quad \mathrm{as}\quad N\to\infty.$$
\end{definition}

\noindent
Intuitively speaking, if two schemes are $g(N)$-alike, then in some sense,
the associated system occupancy states are indistinguishable
on $g(N)$-scale. 
For brevity, for two schemes $\Pi_1$ and $\Pi_2$ that are $g(N)$-alike,
we will often say that $\Pi_1$ and $\Pi_2$ have the same process-level
limits on $g(N)$-scale.
The next theorem states a sufficient criterion for the JSQ$(d(N))$
scheme and the ordinary JSQ policy to be $g(N)$-alike, and thus,
provides the key vehicle in establishing the universality result.

\begin{theorem}
\label{th:pwr of d-intro}
Let $g: \N \to \R_+$ be a function diverging to infinity.
Then the JSQ policy and the JSQ$(d(N))$ scheme are $g(N)$-alike,
with $g(N) \leq N$, if 
\begin{align}
\label{eq:fNalike cond1-intro}
\mathrm{(i)}&\quad d(N) \to \infty, \quad \text{for} \quad g(N) = O(N), \\
\mathrm{(ii)}&\quad d(N) \left(\frac{N}{g(N)}\log\left(\frac{N}{g(N)}\right)\right)^{-1} \to \infty, \quad \text{for} \quad g(N)=o(N).
\label{eq:fNalike cond2-intro}
\end{align}
\end{theorem}

Theorem~\ref{th:pwr of d-intro} yields the next two immediate corollaries.
\begin{corollary}
\label{cor:fluid-intro}
If $d(N) \to \infty$ as $N \to \infty$, then the JSQ$(d(N))$ scheme
and the ordinary JSQ policy are $N$-alike.
\end{corollary}

\begin{corollary}
\label{cor-diff-intro}
If $d(N)/(\sqrt{N} \log(N)) \to \infty$ as $N \to \infty$, then the
JSQ$(d(N))$ scheme and the ordinary JSQ policy are $\sqrt{N}$-alike.
\end{corollary}

Observe that Corollaries~\ref{cor:fluid-intro} and~\ref{cor-diff-intro} together
with the asymptotic results for the JSQ policy
in Subsections~\ref{ssec:jsqfluid-infinite}
and~\ref{ssec:jsq-diffusion-infinite} imply
Theorems~\ref{fluidjsqd-infinite} and~\ref{diffusionjsqd-infinite}.
The rest of the section will be devoted to the proof
of Theorem~\ref{th:pwr of d-intro}.
The proof crucially relies on a novel coupling construction,
which will be used to (lower and upper) bound the difference
of occupancy states of two arbitrary schemes.

\paragraph{The coupling construction.}
Throughout the description of the coupling, we fix $N$, and suppress
the superscript $N$ in the notation. 
Let $Q_i^{\Pi_1}(t)$ and $Q_i^{\Pi_2}(t)$ denote the number of server
pools with at least $i$~active tasks at time~$t$ in two systems
following schemes $\Pi_1$ and $\Pi_2$, respectively.
With a slight abuse of terminology, we occasionally use $\Pi_1$
and $\Pi_2$ to refer to systems following schemes $\Pi_1$ and $\Pi_2$,
respectively.
To couple the two systems, we synchronize the arrival epochs 
and maintain a single exponential departure clock with instantaneous rate
at time $t$ given by $M(t) := \max\left\{\sum_{i=1}^{B} Q_i^{\Pi_1}(t),
\sum_{i=1}^{B} Q_i^{\Pi_2}(t)\right\}$.
We couple the arrivals and departures in the various server pools
as follows:

(1) \emph{Arrival:}
At each arrival epoch, assign the incoming task in each system to one
of the server pools according to the respective schemes.

(2) \emph{Departure:} 
Define $$H(t) :=
\sum_{i=1}^{B} \min\left\{Q_i^{\Pi_1}(t), Q_i^{\Pi_2}(t)\right\}$$
and $$p(t):=
\begin{cases}
\dfrac{H(t)}{M(t)}, & \quad \text{if} \quad M(t) > 0, \\
0, & \quad \text{otherwise.}
\end{cases}
$$
At each departure epoch $t_k$ (say), draw a uniform$[0,1]$ random
variable $U(t_k)$. 
The departures occur in a coupled way based upon the value of $U(t_k)$.
In either of the systems, assign a task index $(i,j)$, if that task is
at the $j^{\mathrm{th}}$ position of the $i^{\mathrm{th}}$ ordered
server pool. 
Let $\mathcal{A}_1(t)$ and $\mathcal{A}_2(t)$ denote the set of all
task indices present at time~$t$ in systems $\Pi_1$ and $\Pi_2$,
respectively. 
Color the indices (or tasks) in $\mathcal{A}_1 \cap \mathcal{A}_2$,
$\mathcal{A}_1 \setminus \mathcal{A}_2$
and $\mathcal{A}_2 \setminus \mathcal{A}_1$, green, blue and red,
respectively, and note that $|\mathcal{A}_1 \cap \mathcal{A}_2| = H(t)$. 
Define a total order on the set of indices as follows:
$(i_1,j_1)<(i_2,j_2)$ if $i_1<i_2$, or $i_1=i_2$ and $j_1<j_2$.
Now, if $U(t_k)\leq p(t_k-)$, then select one green index uniformly
at random and remove the corresponding tasks from both systems.
Otherwise, if $U(t_k)> p(t_k-)$, then choose one integer~$m$,
uniformly at random from all the integers between~$1$
and $M(t) - H(t) = M(t) (1 - p(t))$, and remove the tasks corresponding
to the $m^{\mathrm{th}}$ smallest (according to the order defined above)
red and blue indices in the corresponding systems.
If the number of red (or blue) tasks is less than~$m$, then do nothing.

\begin{figure}
\begin{center}
\includegraphics[scale=1]{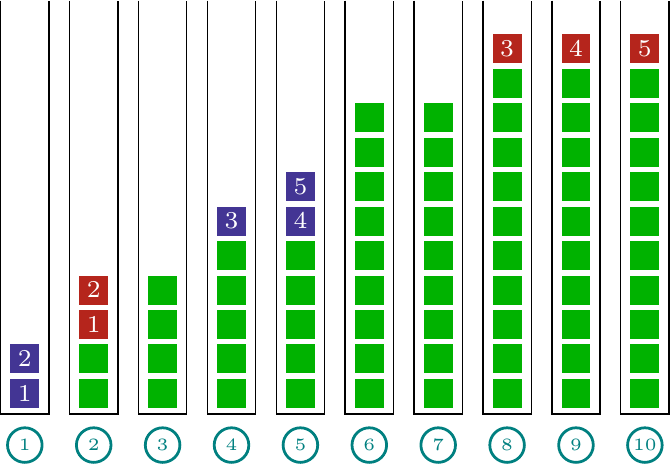}
\caption{T-coupling: Superposition of the occupancy states at some particular time instant, of  schemes $\Pi_1$ and $\Pi_2$ when the server pools in both systems are arranged in nondecreasing order of the number of active tasks. The $\Pi_1$ system is the union of the green and blue tasks, and the $\Pi_2$ system is the union of the green and red tasks.}
\label{fig:tcoupling}
\end{center}
\end{figure}

\begin{figure}
\begin{center}
\begin{tikzpicture}[scale = .8]
  
  \draw(14,1) node[mybox, text width = 3cm,text centered]  {%
   { JSQ$(n(N),d(N))$}};
   
  \draw(22,1) node[mybox, text width = 3cm,text centered]  {%
   {  CJSQ$(n(N))$}};
   
   \draw(14,6) node[mybox, text width = 3cm,text centered]  {%
   {  JSQ$(d(N))$}};
   
   \draw(22,6) node[mybox, text width = 3cm,text centered]  {%
   {  JSQ}};

\draw[doublearr] (16.5,6) to (19.5,6);
\draw[doublearr2] (14,2) to (14,5);
\draw[doublearr2] (16.5,1) to (19.5,1);
\draw[doublearr2] (22,2) to (22,5);

\node  at (18,6.5) {\small Theorem~\ref{th:pwr of d-intro}};

\node [rotate=270] at (21.5,3.5) {\small Proposition~\ref{prop: modified JSL-intro}};
\node [rotate=270] at (22.5,3.5) {\small Suitable $n(N)$};

\node [rotate=90] at (14.5,3.5) {\small Proposition~\ref{prop: power of d-intro}};
\node [rotate=90] at (13.5,3.5) {\small Suitable $d(N)$};

\draw[text width=3cm, align=center] (18,1.75) node {\small Belongs to};
\draw[text width=2cm, align=center] (18,1.35) node {\small the class};
\end{tikzpicture}
\caption{The asymptotic equivalence structure is depicted for various intermediate load balancing schemes to facilitate the comparison between the JSQ$(d(N))$ scheme and the ordinary JSQ policy.}
\label{fig:sfigrelation-inf}
\end{center}
\end{figure}
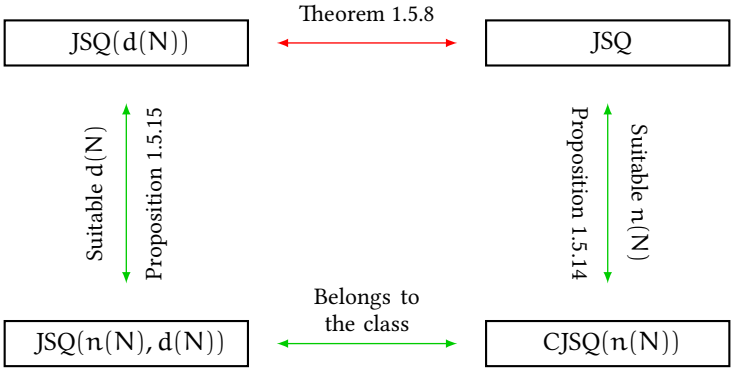

The above coupling has been schematically represented
in Figure~\ref{fig:tcoupling}, and will henceforth be referred to as
T-coupling, where T stands for `task-based'. 
Now we need to show that, under the T-coupling, the two systems,
considered independently, evolve according to their own statistical laws. 
This can be seen in several steps.
Indeed, the T-coupling basically uniformizes the departure rate
by the maximum number of tasks present in either of the two systems. 
Then informally speaking, the green region signifies the common portion of
tasks, and the red and blue regions represent the separate contributions. 
Now observe that
\begin{enumerate}[{\normalfont (i)}]
\item The total departure rate from $\Pi_i$ is 
\begin{align*}
& M(t) \left[p(t)+(1-p(t)) \frac{|\mathcal{A}_i\setminus\mathcal{A}_{3-i}|}{M(t)-H(t)}\right]
= |\mathcal{A}_1\cap\mathcal{A}_2|+|\mathcal{A}_i\setminus\mathcal{A}_{3-i}|
=|\mathcal{A}_i|, 
\end{align*}
for $i= 1,2.$
\item Assuming without loss of generality
$|\mathcal{A}_1| \geq |\mathcal{A}_2|$, each task in $\Pi_1$ is equally
likely to depart.
\item Each task in $\Pi_2$ within $\mathcal{A}_1\cap \mathcal{A}_2$
and each task within $\mathcal{A}_2 \setminus \mathcal{A}_1$ is equally
likely to depart, and the probabilities of departures are proportional
to $|\mathcal{A}_1 \cap \mathcal{A}_2|$
and $|\mathcal{A}_2 \setminus \mathcal{A}_1|$, respectively.
\end{enumerate}

\begin{remark}[Comparison of T-coupling and S-coupling]
\label{rem:contrast}
\normalfont
As briefly mentioned earlier, in the current infinite-server scenario,
the departures of the ordered server pools cannot be coupled,
mainly since the departure rate at the $m^\mathrm{ th}$ ordered server pool,
for some $m = 1, 2, \ldots, N$, depends on its number of active tasks.
It is worthwhile to mention that the T-coupling in the current section is
stronger than the S-coupling used in Section~\ref{univ} in the single-server queueing scenario.
Observe that due to Lemma~\ref{lem:majorization}, the absolute
difference of the occupancy states of the JSQ policy and any scheme
from the CJSQ class at any time point can be bounded deterministically
(without any terms involving the cumulative number of lost tasks).
It is worth emphasizing that the universality result on some specific
scale, stated in Theorem~\ref{th:pwr of d-intro}, does not depend on the
behavior of the JSQ policy on that scale, whereas in the single-server queueing scenario
it does, mainly because the upper and lower bounds
in Corollary~\ref{cor:bound-intro} involve tail sums of two different policies.
More specifically, in the single-server queueing scenario the fluid and diffusion limit results of CJSQ($n(N)$) class crucially use those of the MJSQ($n(N)$) scheme, while in the current scenario it does not -- the results for the MJSQ($n(N)$) scheme comes as a consequence of those for the CJSQ($n(N)$) class of schemes.
Also, the bounds in Lemma~\ref{lem:majorization} 
do not depend on~$t$, and hence, apply in the steady state as well.
Moreover, the S-coupling compares the $k$ \emph{highest} horizontal bars, whereas the T-coupling
in the current section compares the $k$ \emph{lowest}
horizontal bars.
As a result, the bounds on the occupancy states established
in Corollary~\ref{cor:bound-intro} 
involve tail sums of the occupancy states of the ordinary JSQ policy,
which necessitates proving the convergence of tail sums of the
occupancy states of the ordinary JSQ policy. 
In contrast, as we will see in the proof of Proposition~\ref{prop: modified JSL-intro}, the bound in the infinite-server scenario involves only a single component
(see Equations \eqref{eq:upperocc} and \eqref{eq:lowerocc}),
and thus, proving convergence of each component suffices.
\end{remark}

The T-coupling can be used to derive several stochastic inequality results
that will play an instrumental role in proving Theorem~\ref{th:pwr of d-intro}.

In order to compare the JSQ policy with the CJSQ($n(N)$) schemes,
denote by $Q_i^{\Pi_1}(t)$ and $Q_i^{\Pi_2}(t)$ the number of server
pools with at least $i$~tasks under the JSQ policy and $\CJSQ(n(N))$
scheme, respectively.

\begin{lemma}
\label{lem:majorization}
For any $k \in \big\{1, 2, \ldots, B\big\}$,
\begin{equation}
\label{eq:majorization}
\left\{\sum_{i=1}^{k} Q_i^{\Pi_1}(t) - k n(N)\right\}_{t \geq 0} \leq_{st}
\left\{\sum_{i=1}^{k} Q_i^{\Pi_2}(t)\right\}_{t \geq 0} \leq_{st}
\left\{\sum_{i=1}^{k} Q_i^{\Pi_1}(t)\right\}_{t \geq 0},
\end{equation}
provided at $t=0$ the two systems start from the same occupancy states.
\end{lemma}

In the next remark we comment on the contrast
of Lemma~\ref{lem:majorization} with stochastic dominance properties
for the ordinary JSQ policy in the existing literature.

\begin{remark}
\label{rem:novelty}
{\normalfont
The stochastic ordering in Lemma~\ref{lem:majorization} is to be
contrasted with the weak majorization results
in~\cite{Winston77,towsley,Towsley95,Towsley1992,W78} in the context
of the ordinary JSQ policy in the single-server queueing scenario,
and in~\cite{STC93,J89,M87,MS91} in the scenario of state-dependent
service rates, non-decreasing with the number of active tasks.
In the current infinite-server scenario, the results
in~\cite{STC93,J89,M87,MS91} imply that for any non-anticipating
scheme~$\Pi$ taking assignment decisions based on the number of active
tasks only, for all $t \geq 0$,
\begin{align}
\label{eq: towsley}
\sum_{m=1}^{\ell} X_{(m)}^\mathrm{ JSQ}(t) & \leq_{st}
\sum_{m=1}^\ell X_{(m)}^{\Pi}(t), \mbox{ for } \ell = 1, 2, \ldots, N, \\
\left\{L^\mathrm{ JSQ}(t)\right\}_{t \geq 0} &\leq_{st}
\left\{L^{\Pi}(t)\right\}_{t \geq 0},
\end{align}
where $X_{(m)}^\Pi(t)$ is the number of tasks in the $m^{\mathrm{th}}$
ordered server pool at time~$t$ in the system following scheme~$\Pi$
and $L^{\Pi}(t)$ is the total number of overflow events under policy~$\Pi$
up to time~$t$. 
Observe that $X_{(m)}^\Pi$ can be visualized as the $m^{\mathrm{th}}$
largest (rightmost) vertical bar (or stack) in Figure~\ref{figB}.
Thus~\eqref{eq: towsley} says that the sum of the lengths of the
$\ell$~largest \emph{vertical} stacks in a system following any
scheme~$\Pi$ is stochastically larger than or equal to that following
the ordinary JSQ policy for any $\ell = 1, 2, \ldots, N$.
Mathematically, this ordering can be equivalently written as
\begin{equation}
\label{eq:equiv-ord}
\sum_{i = 1}^{B} \min\big\{\ell, Q_i^\mathrm{ JSQ}(t)\big\} \leq_{st}
\sum_{i = 1}^{B} \min\big\{\ell, Q_i^{\Pi}(t)\big\},
\end{equation}
for all $\ell = 1, \dots, N$.
In contrast, in order to show asymptotic equivalence on various scales,
we need to both upper and lower bound the occupancy states of the
$\CJSQ(n(N))$ schemes in terms of the JSQ policy, and therefore need
a much stronger hold on the departure process.
The T-coupling provides us just that, and has several useful properties
that are crucial for our proof technique.
For example, Proposition~\ref{prop:stoch-ord2-intro} uses the fact that
if two systems are T-coupled, then departures cannot increase the sum
of the absolute differences of the $Q_i$-values, which is not true
for the coupling considered in the above-mentioned literature.
The left stochastic ordering in~\eqref{eq:majorization} also does not
remain valid in those cases.
Furthermore, observe that the right inequality in~\eqref{eq:majorization}
(i.e., $Q_i$'s) implies the stochastic inequality is \emph{reversed}
in~\eqref{eq:equiv-ord}, which is counter-intuitive in view of the
well-established optimality properties of the ordinary JSQ policy,
as mentioned above.
The fundamental distinction between the two coupling techniques is
also reflected by the fact that the T-coupling does not allow
for arbitrary nondecreasing state-dependent departure rate functions,
unlike the couplings in~\cite{STC93,J89,M87,MS91}.}
\end{remark}

\begin{proposition}
\label{prop: modified JSL-intro}
For any function $g:\N\to\R_+$ diverging to infinity,
if $$n(N)/ g(N)\to 0\quad \mbox{as}\quad N \to \infty,$$ 
then the JSQ policy and the
$\CJSQ(n(N))$ schemes are $g(N)$-alike.
\end{proposition}

\begin{proof}[Proof of Proposition~\ref{prop: modified JSL-intro}]
Using Lemma~\ref{lem:majorization}, there exists a common probability
space such that for any $k \geq 1$ we can write
\begin{equation}\label{eq:upperocc}
\begin{split}
Q_k^{\Pi_2}(t) &= \sum_{i=1}^{k} Q_i^{\Pi_2}(t) - \sum_{i=1}^{k-1} Q_i^{\Pi_2}(t) \\
&\leq \sum_{i=1}^{k} Q_i^{\Pi_1}(t) - \sum_{i=1}^{k-1} Q_i^{\Pi_1}(t) + k n(N) \\
&= Q_k^{\Pi_1}(t) + k n(N).
\end{split}
\end{equation}
Similarly, we can write
\begin{equation}
\label{eq:lowerocc}
\begin{split}
Q_k^{\Pi_2}(t) &= \sum_{i=1}^{k} Q_i^{\Pi_2}(t) - \sum_{i=1}^{k-1} Q_i^{\Pi_2}(t) \\
&\geq \sum_{i=1}^{k} Q_i^{\Pi_1}(t) - k n(N) - \sum_{i=1}^{k-1} Q_i^{\Pi_1}(t) \\
&= Q_k^{\Pi_1}(t) - k n(N).
\end{split}
\end{equation}
Therefore, for all $k \geq 1$, we have
$\sup_t |Q_k^{\Pi_2}(t) - Q_k^{\Pi_1}(t)| \leq k n(N)$.
Since we know $n(N)/ g(N) \to 0$ as $N \to \infty$, the proof is complete.
\end{proof}

Next we compare schemes from the CJSQ($n(N)$) class with the JSQ($d(N)$) scheme. 
The comparison follows a somewhat similar line of argument as 
in Section~\ref{univ}, and involves a JSQ$(n(N),d(N))$ scheme
which is an intermediate blend between the $\CJSQ(n(N))$ schemes
and the JSQ$(d(N))$ scheme.
Specifically, the JSQ$(n(N),d(N))$ scheme selects a candidate server
pool in the exact same way as the JSQ$(d(N))$ scheme.
However, it only assigns the task to that server pool if it belongs
to the $n(N)+1$ lowest ordered ones,
and to a randomly selected server pool among these otherwise.
Note that by construction, the JSQ$(n(N),d(N))$ scheme belongs to the
$\CJSQ(n(N))$ class of schemes.

We now consider two T-coupled systems: one with a JSQ$(d(N))$
scheme and another with a JSQ$(n(N),d(N))$ scheme.
Assume that at some specific arrival epoch, the incoming task is
assigned to the $k^{\mathrm{th}}$ ordered server pool in the system
under the JSQ($d(N)$) scheme. 
If $k \in \big\{1, 2, \ldots, n(N)+1\big\}$, then the scheme
JSQ$(n(N),d(N))$ also assigns the arriving task to the $k^{\mathrm{th}}$
ordered server pool. 
Otherwise it dispatches the arriving task uniformly at random
among the first $n(N)+1$ ordered server pools.

We will establish a sufficient criterion on $d(N)$ in order for the
JSQ$(d(N))$ scheme and JSQ$(n(N),d(N))$ scheme to be close in terms
of $g(N)$-alikeness, as stated in the next proposition.

\begin{proposition}
\label{prop: power of d-intro}
Assume, $n(N) / g(N)\to 0$ as $N \to \infty$ for some function
$g: \N\to \R_+$ diverging to infinity.
The JSQ$(n(N),d(N))$ scheme and the JSQ($d(N)$) scheme are $g(N)$-alike
if the following condition holds:
\begin{equation}
\label{eq:condition-same-intro}
\frac{n(N)}{N}d(N) - \log\frac{N}{g(N)} \to \infty, \quad \text{as}
\quad N \to \infty.
\end{equation}
\end{proposition}

Finally, Proposition~\ref{prop: power of d-intro} in conjunction
with Proposition~\ref{prop: modified JSL-intro} yields Theorem~\ref{th:pwr of d-intro}.
The overall proof strategy as described above, is schematically
represented in Figure~\ref{fig:sfigrelation-inf}.

\section{Universality of load balancing in networks}
\label{networks}

In this section we return to the single-server queueing dynamics,
and extend the universality properties to network scenarios,
where the $N$~servers are assumed to be inter-connected by some
underlying graph topology~$G_N$.
An extensive treatment of the model considered in this section can be found in Chapters~\ref{chap:networkjsq} and~~\ref{chap:networkjsqd}.

Tasks arrive at the various servers as independent Poisson processes
of rate~$\lambda$, and each incoming task is assigned to whichever
server has the smallest number of tasks among the one where it arrives
and its neighbors in~$G_N$.  
Thus, in case $G_N$ is a clique, each incoming task is assigned to the
server with the shortest queue across the entire system,
and the behavior is equivalent to that under the JSQ policy.
The stochastic optimality properties of the JSQ policy thus imply that
the queue length process in a clique will be better balanced
and smaller (in a majorization sense) than in an arbitrary graph~$G_N$.

Besides the prohibitive communication overhead discussed earlier,
a further scalability issue of the JSQ policy arises when executing
a task involves the use of some data.
Storing such data for all possible tasks on all servers will typically
require an excessive amount of storage capacity.
These two burdens can be effectively mitigated in sparser graph
topologies where tasks that arrive at a specific server~$i$ are only
allowed to be forwarded to a subset of the servers ${\mathcal N}_i$.
For the tasks that arrive at server~$i$, queue length information
then only needs to be obtained from servers in ${\mathcal N}_i$,
and it suffices to store replicas of the required data on the
servers in ${\mathcal N}_i$.
The subset ${\mathcal N}_i$ containing the peers of server~$i$ can
be naturally viewed as its neighbors in some graph topology~$G_N$.
Here we consider the case of undirected graphs,
but most of the analysis can be extended to directed graphs.

While sparser graph topologies relieve the scalability issues
associated with a clique, the queue length process will be worse
(in the majorization sense) because of the limited connectivity.
Surprisingly, however, even quite sparse graphs can asymptotically
match the optimal performance of a clique, provided they are
suitably random, as we will further describe below.

The above model has been studied in~\cite{G15,T98}, focusing
on certain fixed-degree graphs and in particular ring topologies.
The results demonstrate that the flexibility to forward tasks
to a few neighbors, or even just one, with possibly shorter queues
significantly improves the performance in terms of the waiting time
and tail distribution of the queue length.
This resembles the ``power-of-choice'' gains observed for JSQ($d$)
policies in complete graphs.

However, the results in \cite{G15,T98} also establish that
the performance sensitively depends on the underlying graph topology,
and that selecting from a fixed set of $d - 1$ neighbors typically
does not match the performance of re-sampling $d - 1$ alternate
servers for each incoming task from the entire population,
as in the power-of-$d$ scheme in a complete graph.
Further related problems have been investigated
in \cite{M96focs,ACM01,KM00,MPS02}.

If tasks do not get served and never depart but simply accumulate,
then the scenario described above amounts to a so-called balls-and-bins
problem on a graph.
Viewed from that angle, a close counterpart of our setup is studied
in Kenthapadi \& Panigrahy~\cite{KP06}, where in our terminology each
arriving task is routed to the shortest of $d \geq 2$ randomly
selected neighboring queues.
In this setup \cite{KP06} show that if the underlying graph is almost
regular with degree $N^{\varepsilon}$, where $\varepsilon$ is not too
small, the maximum number of balls in a bin scales as
$\log(\log(N)) / \log(d) + O(1)$, just like when the underlying graph
is a clique.
There are fundamental differences between the ball-and-bins and the
queueing scenarios, however, and an inherently different approach is
required in the current setup than what was developed in~\cite{KP06}.
Moreover, \cite{KP06} considers only the scaling of the maximum
queue length, whereas we analyze a more detailed time-varying
evolution of the entire system along with its stationary behavior.
We will further elaborate on the connections and differences
with balls-and-bins problems in Subsection~\ref{ballsbins}

When each arriving task is routed to the shortest of $d \geq 2$ randomly
selected neighboring queues, the process-level convergence in the transient regime is established in Chapter~\ref{chap:networkjsqd}.
In this work, we analyze the evolution of the queue length process at an arbitrary tagged server as the system size becomes large.
The main ingredient is a careful analysis of local occupancy measures associated with the neighborhood of each server and to argue that under suitable conditions their asymptotic behavior is the same for all servers.
Under mild conditions on the graph topology $G_N$ (diverging minimum degree and the ratio between minimum degree and maximum degree in each connected component converges to 1), for a suitable initial occupancy measure, Theorem~\ref{th:deterministic-d} in Chapter~\ref{chap:networkjsqd} establishes that for any fixed $d\geq 2$, the global occupancy state process for the JSQ($d$) scheme on $G_N$ has the same weak limit in~\eqref{fluid:standard} as that on a clique, as the number of vertices $N$ becomes large.
Also, the propagation of chaos property was shown to hold for this system, in the sense that the  queue lengths at any finite collection of tagged servers are asymptotically independent, and the 
 queue length process for each server converges in distribution (in the path space) to the corresponding McKean-Vlasov process, see Theorem~\ref{th:tagged-d} in Chapter~\ref{chap:networkjsqd}.
 Furthermore, when the graph sequence is random, with the $N$-th graph given as an Erd\H{o}s-R\'enyi random graph (ERRG) on $N$ vertices with average degree $c(N)$, annealed convergence of the occupancy process to the same deterministic limit as above, is established under the condition $c(N)\to\infty$, and under a stronger condition $c(N)/ \log N\to\infty$, convergence (in probability) is shown for almost every realization of the random graph.

As mentioned above, the queue length process in a clique will be
better balanced and smaller (in a majorization sense) than
in an arbitrary graph~$G_N$.
Accordingly, a graph $G_N$ is said to be $N$-optimal or
$\sqrt{N}$-optimal when the queue length process on $G_N$ is equivalent
to that on a clique on an $N$-scale or $\sqrt{N}$-scale, respectively.
Roughly speaking, a graph is $N$-optimal if the \emph{fraction}
of nodes with $i$~tasks, for $i = 0, 1, \ldots$, behaves as in a clique
as $N \to \infty$.
The fluid-limit results for the JSQ policy discussed
in Subsection~\ref{ssec:jsqfluid} imply that the latter fraction is
zero in the limit for all $i \geq 2$ in a clique in stationarity,
i.e., the fraction of servers with two or more tasks vanishes in any
graph that is $N$-optimal, and consequently the mean waiting time
vanishes as well as $N \to \infty$.
Furthermore, the diffusion-limit results of~\cite{EG15} for the JSQ
policy discussed in Subsection~\ref{ssec:diffjsq} imply that the
number of nodes with zero tasks and that with two tasks both scale as
$\sqrt{N}$ as $N \to \infty$.
Again loosely speaking, a graph is $\sqrt{N}$-optimal if in the
heavy-traffic regime the number of nodes with zero tasks and that
with two tasks when scaled by $\sqrt{N}$ both evolve as in a clique
as $N \to \infty$.
Formal definitions of asymptotic optimality on an $N$-scale or
$\sqrt{N}$-scale will be introduced in Definition~\ref{def:opt} below.

As one of the main results, we will demonstrate that, remarkably,
asymptotic optimality can be achieved in quite sparse
ERRGs.
We prove that a sequence of ERRGs indexed by the number of vertices~$N$
with $d(N) \to \infty$ as $N \to \infty$, is $N$-optimal.
We further establish that the latter growth condition for the average
degree is in fact necessary in the sense that any graph sequence that
contains $\Theta(N)$ bounded-degree vertices cannot be $N$-optimal.
This implies that a sequence of ERRGs with finite average degree
cannot be $N$-optimal.
The growth rate condition is more stringent for optimality
on $\sqrt{N}$-scale in the heavy-traffic regime.
Specifically, we prove that a sequence of ERRGs indexed by the number
of vertices $N$ with $d(N) / (\sqrt{N} \log(N)) \to \infty$
as $N \to \infty$, is $\sqrt{N}$-optimal.

The above results demonstrate that the asymptotic optimality of cliques
on an $N$-scale and $\sqrt{N}$-scale can be achieved in far sparser graphs,
where the number of connections is reduced by nearly a factor $N$
and $\sqrt{N}/\log(N)$, respectively, provided the topologies are
suitably random in the ERRG sense.
This translates into equally significant reductions in communication
overhead and storage capacity, since both are roughly proportional
to the number of connections.

While quite sparse graphs can achieve asymptotic optimality in the
presence of randomness, the worst-case graph instance may even
in very dense regimes (high average degree) not be optimal.
In particular, we prove that any graph sequence with minimum degree
$N-o(N)$ is $N$-optimal, but that for any $0 < c < 1/2$ one can
construct graphs with minimum degree $c N + o(N)$ which are not
$N$-optimal for some $\lambda < 1$.

The key challenge in the analysis of load balancing on arbitrary graph
topologies is that one needs to keep track
of the evolution of the number of tasks at each vertex along with their
corresponding neighborhood relationship.
This creates a major problem in constructing a tractable Markovian
state descriptor, and renders a direct analysis of such processes
highly intractable.
Consequently, even asymptotic results for load balancing processes
on an arbitrary graph have remained scarce so far.
We take a radically different approach and aim to compare the load
balancing process on an arbitrary graph with that on a clique.
Specifically, rather than analyze the behavior for a given class
of graphs or degree value, we explore for what types of topologies
and degree properties the performance is asymptotically similar
to that in a clique.

Our proof arguments build on the stochastic coupling constructions
developed in Section~\ref{univ} for JSQ($d$) policies.
Specifically, we view the load balancing process on an arbitrary
graph as a `sloppy' version of that on a clique, and thus construct
several other intermediate sloppy versions.
By constructing novel couplings, we develop a method of comparing
the load balancing process on an arbitrary graph and that on a clique. 
In particular, we bound the difference between the fraction of vertices
with $i$ or more tasks in the two systems for $i = 1, 2, \dots$,
to obtain asymptotic optimality results.
From a high level, conceptually related graph conditions for
asymptotic optimality were examined using quite different
techniques by Tsitsiklis and Xu \cite{TX11,TX13} in a dynamic
scheduling framework (as opposed to the load balancing context).

For $k = 1, \ldots, N$, denote by $X_k(G_N, t)$ the queue length
at the $k$-th server at time~$t$ (including the task possibly in service),
and by $X_{(k)}(G_N, t)$ the queue length at the $k$-th ordered server
at time~$t$ when the servers are arranged in  non-decreasing order
of their queue lengths (ties can be broken in some way that will be
evident from the context).
Let $Q_i(G_N, t)$ denote the number of servers with queue length
at least~$i$ at time~$t$, $i = 1, 2, \ldots, B$.
It is important to note that $\{(q_i(G_N, t))_{i \geq 1}\}_{t \geq 0}$
is itself \emph{not} a Markov process, but the joint process
$\{(q_i(G_N, t))_{i \geq 1}, (X_{k}(G_N, t))_{k=1}^{N}\}_{t \geq 0}$
is Markov.
Also, in the Halfin-Whitt heavy-traffic regime~\eqref{eq:HW},
define the centered and scaled processes 
\begin{equation}
\label{eq:HWOcc}
\bar{Q}_1(G_N,t) = - \frac{N - Q_1(G_N,t)}{\sqrt{N}}, \qquad
\bar{Q}_i(G_N,t) = \frac{ Q_i(G_N,t)}{\sqrt{N}},
\end{equation}
analogous to~\eqref{eq:diffscale}. \\


As stated before, a clique is an optimal load balancing topology,
as the occupancy process is better balanced and smaller
(in a majorization sense) than in any other graph topology.
In general the optimality is strict, but it turns out that
near-optimality can be achieved asymptotically in a broad class
of other graph topologies.
Therefore, we now introduce two notions of \emph{asymptotic optimality},
which will be useful to characterize the performance in large-scale
systems. 

\begin{definition}[{Asymptotic optimality}]
\label{def:opt}
A graph sequence $\GG = \{G_N\}_{N \geq 1}$ is called `asymptotically
optimal on $N$-scale' or `$N$-optimal', if for any $\lambda < 1$,
the scaled occupancy process $(q_1(G_N, \cdot), q_2(G_N, \cdot), \ldots)$
converges weakly, on any finite time interval, to the process
$(q_1(\cdot), q_2(\cdot),\ldots)$ given by~\eqref{eq:fluid-intro}.

Moreover, a graph sequence $\GG = \{G_N\}_{N \geq 1}$ is called
`asymptotically optimal on $\sqrt{N}$-scale' or `$\sqrt{N}$-optimal',
if in the Halfin-Whitt heavy-traffic regime~\eqref{eq:HW},
on any finite time interval, the process
$(\bQ_1(G_N, \cdot), \bQ_2(G_N, \cdot), \ldots)$ as in~\eqref{eq:HWOcc}
converges weakly to the process $(\bQ_1(\cdot), \bQ_2(\cdot), \ldots)$
given by~\eqref{eq:diffusionjsq}.
\end{definition}

Intuitively speaking, if a graph sequence is $N$-optimal
or $\sqrt{N}$-optimal, then in some sense, the associated occupancy
processes are indistinguishable from those of the sequence of cliques
on $N$-scale or $\sqrt{N}$-scale.
In other words, on any finite time interval their occupancy processes
can differ from those in cliques by at most $o(N)$ or $o(\sqrt{N})$,
respectively. 
For brevity, $N$-scale and $\sqrt{N}$-scale will henceforth be referred
to as \emph{fluid scale} and \emph{diffusion scale}, respectively.
In particular, exploiting interchange of the stationary ($t \to \infty$)
and many-server ($N \to \infty$) limits, we obtain that for any
$N$-optimal graph sequence $\{G_N\}_{N \geq 1}$,
\begin{equation}
q_1(G_N, \infty) \to \lambda \quad \mbox{ and} \quad
q_i(G_N, \infty) \to 0 \quad \mbox{ for all } i = 2, \dots, B,
\end{equation}
as  $N \to \infty$,
implying that the stationary fraction of servers with queue length two
or larger and the mean waiting time vanish.



\subsection{Asymptotic optimality criteria for deterministic graph sequences}

We now proceed to develop a criterion for asymptotic optimality
of an arbitrary deterministic graph sequence on different scales.
Next this criterion will be leveraged to establish optimality
of a sequence of random graphs.
We start by introducing some useful notation, and two measures
of \emph{well-connectedness}.
Let $G = (V, E)$ be any graph.
For a subset $U \subseteq V$, define $\com(U) := |V\setminus N[U]|$
to be the cardinality of the set of all vertices that do not share an edge with any vertex in $U$, where
$N[U] := U\cup \{v \in V:\ \exists\ u \in U \mbox{ with } (u, v) \in E\}$.
For any fixed $\varepsilon > 0$ define
\begin{equation}
\begin{split}
\label{def:dis}
\dis_1(G,\varepsilon) &:= \sup_{U\subseteq V, |U|\geq \varepsilon |V|}\com(U),
\\
\dis_2(G,\varepsilon) &:= \sup_{U\subseteq V, |U|\geq \varepsilon \sqrt{|V|}}\com(U).
\end{split}
\end{equation}

The next theorem provides sufficient conditions for asymptotic
optimality on $N$-scale and $\sqrt{N}$-scale in terms of the above
two well-connectedness measures.

\begin{theorem}
\label{th:det-seq}
For any graph sequence $\GG = \{G_N\}_{N \geq 1}$,
\begin{enumerate}[{\normalfont (i)}]
\item $\GG$ is $N$-optimal if for any $\varepsilon > 0$, 
${\normalfont\dis}_1(G_N, \varepsilon) / N \to 0$ as $N \to \infty$.
\item $\GG$ is $\sqrt{N}$-optimal if for any $\varepsilon > 0$, 
${\normalfont\dis}_2(G_N, \varepsilon) / \sqrt{N} \to 0$ as $N \to \infty$.
\end{enumerate}
\end{theorem}

The next corollary is an immediate consequence
of Theorem~\ref{th:det-seq}.

\begin{corollary}
Let $\GG = \{G_N\}_{N \geq 1}$ be any graph sequence.
Then \textrm{\normalfont(i)} If $d_{\min}(G_N) = N - o(N)$, then $\GG$ is $N$-optimal,
and 
\textrm{\normalfont(ii)} If $d_{\min}(G_N) = N - o(\sqrt{N})$, then $\GG$ is
$\sqrt{N}$-optimal.
\end{corollary}

The rest of the subsection is devoted to a discussion of the main
proof arguments for Theorem~\ref{th:det-seq}, focusing on the proof
of $N$-optimality. 
The proof of $\sqrt{N}$-optimality follows along similar lines.
We establish in Proposition~\ref{prop:cjsq} that if a system is able
to assign each task to a server in the set $\cS^N(n(N))$ of the $n(N)$
nodes with shortest queues, where $n(N)$ is $o(N)$,
then it is $N$-optimal. 
Since the underlying graph is not a clique however (otherwise there
is nothing to prove), for any $n(N)$ not every arriving task can be
assigned to a server in $\cS^N(n(N))$.
Hence we further prove in Proposition~\ref{prop:stoch-ord-new}
a stochastic comparison property implying that if on any finite time
interval of length~$t$, the number of tasks $\Delta^N(t)$ that are not
assigned to a server in $\cS^N(n(N))$ is $o_P(N)$, then the system is
$N$-optimal as well.
The $N$-optimality can then be concluded when $\Delta^N(t)$ is
$o_P(N)$, which we establish in Proposition~\ref{prop:dis-new}
under the condition that $\dis_1(G_N, \varepsilon) / N \to 0$ as
$N \to \infty$ as stated in Theorem~\ref{th:det-seq}.

To further explain the idea described in the above proof outline,
it is useful to adopt a slightly different point of view towards
load balancing processes on graphs.
From a high level, a load balancing process can be thought of as follows:
there are $N$~servers, which are assigned incoming tasks by some scheme.
The assignment scheme can arise from some topological structure,
in which case we will call it \emph{topological load balancing},
or it can arise from some other property of the occupancy process,
in which case we will call it \emph{non-topological load balancing}.
As mentioned earlier, the JSQ policy or the clique is optimal
among the set of all non-anticipating schemes, irrespective of being
topological or non-topological.
Also, load balancing on graph topologies other than a clique can be
thought of as a `sloppy' version of that on a clique, when each server
only has access to partial information on the occupancy state.
Below we first introduce a different type of sloppiness in the task
assignment scheme, and show that under a limited amount of sloppiness
optimality is retained on a suitable scale.
Next we will construct a scheme which is a hybrid of topological
and non-topological schemes, whose behavior is simultaneously close
to both the load balancing process on a suitable graph and that
on a clique.

\paragraph{A class of sloppy load balancing schemes.}

Fix some function $n: \N \to \N$, and recall the set $\cS^N(n(N))$
as before as well as the class $\CJSQ(n(N))$ where each arriving task
is assigned to one of the servers in $\cS^N(n(N))$. 
It should be emphasized that for any scheme in $\CJSQ(n(N))$,
we are not imposing any restrictions on how the incoming task should
be assigned to a server in $\cS^N(n(N))$.
The scheme only needs to ensure that the arriving task is assigned
to some server in $\cS^N(n(N))$ with respect to \emph{some} tie-breaking mechanism.
Observe that using Corollary~\ref{cor:bound-intro} and following the
arguments as in the proof of Theorems~\ref{fluidjsqd}
and~\ref{diffusionjsqd}, we obtain the next proposition,
which provides a sufficient criterion for asymptotic optimality
of any scheme in $\CJSQ(n(N))$.

\begin{proposition}
\label{prop:cjsq}
For $0 \leq n(N) < N$, let $\Pi \in \CJSQ(n(N))$ be any scheme.
\textrm{\normalfont(i)} If $n(N) / N \to 0$ as $N \to \infty$, then $\Pi$ is $N$-optimal,
and
\textrm{\normalfont(ii)} If $n(N) / \sqrt{N} \to 0$ as $N \to \infty$, then $\Pi$ is
$\sqrt{N}$-optimal.
\end{proposition}

\paragraph{A bridge between topological and non-topological load balancing.}
For any graph $G_N$ and $n \leq N$, we first construct a scheme
called $I(G_N, n)$, which is an intermediate blend between the
topological load balancing process on $G_N$ and some kind
of non-topological load balancing on $N$~servers.
The choice of $n = n(N)$ will be clear from the context.

To describe the scheme $I(G_N, n)$, first synchronize the arrival
epochs at server~$v$ in both systems, $v = 1, 2, \ldots, N$.
Further, synchronize the departure epochs at the $k$-th ordered
server with the $k$-th smallest number of tasks in the two systems,
$k = 1, 2, \ldots, N$.
When a task arrives at server~$v$ at time~$t$ say, it is assigned
in the graph $G_N$ to a server $v' \in N[v]$ according to its own
statistical law.
For the assignment under the scheme $I(G_N,n)$, first observe that if
\begin{equation}
\label{eq:criteria}
\min_{u \in N[v]} X_u(G_N,t) \leq \max_{u \in \cS(n)} X_u(G_N,t),
\end{equation}
then there exists \emph{some} tie-breaking mechanism for which
$v' \in N[v]$ belongs to $\cS(n)$ under $G_N$.
Pick such an ordering of the servers, and assume that $v'$ is the
$k$-th ordered server in that ordering, for some $k \leq n+1$.
Under $I(G_N,n)$ assign the arriving task to the $k$-th ordered
server (breaking ties arbitrarily in this case).
Otherwise, if \eqref{eq:criteria} does not hold, then the task is
assigned to one of the $n+1$ servers with minimum queue lengths
under $G_N$ uniformly at random.

Denote by $\Delta^N(I(G_N, n), T)$ the cumulative number of arriving
tasks up to time $T \geq 0$ for which Equation~\eqref{eq:criteria} is
violated under the above coupling.
The next proposition shows that the load balancing process under the
scheme $I(G_N,n)$ is close to that on the graph $G_N$ in terms of the
random variable $\Delta^N(I(G_N, n), T)$.

\begin{proposition}
\label{prop:stoch-ord-new}
The following inequality is preserved almost surely:
\begin{equation}
\label{eq:stoch-ord-new}
\sum_{i=1}^{B} |Q_i(G_N, t) - Q_i(I(G_N, n), t)| \leq
2 \Delta^N(I(G_N, n), t) \qquad \forall\ t \geq 0,
\end{equation}
provided the two systems start from the same occupancy state at $t=0$.
\end{proposition}

In order to conclude optimality on $N$-scale or $\sqrt{N}$-scale,
it remains to be shown that $\Delta^N(I(G_N, n), T)$ is sufficiently small.
The next proposition provides suitable asymptotic bounds for
$\Delta^N(I(G_N, n), T)$ under the conditions on $\dis_1(G_N, \varepsilon)$
and $\dis_2(G_N, \varepsilon)$ stated in Theorem~\ref{th:det-seq}.

\begin{proposition}
\label{prop:dis-new}
The following properties hold for any graph sequence:
\begin{enumerate}[{\normalfont (i)}]
\item For any $\varepsilon>0$, there exists $\varepsilon'>0$ and $n_{\varepsilon'}(N)$ with $n_{\varepsilon'}(N)/N\to 0$ as $N\to\infty$, such that if $\dis_1(G_N,\varepsilon')/N\to 0$ as $N\to\infty$, then for all $T>0$, $\Pro{\Delta^N(I(G_N,n_{\varepsilon'}),T)/N>\varepsilon}\to 0$. 
\item For any $\varepsilon>0$, there exists $\varepsilon'>0$ and $m_{\varepsilon'}(N)$ with $m_{\varepsilon'}(N)/\sqrt{N}\to 0$ as $N\to\infty$, such that if $\dis_2(G_N,\varepsilon')/\sqrt{N}\to 0$ as $N\to\infty$, then for all $T>0$, $\Pro{\Delta^N(I(G_N,m_{\varepsilon'}),T)/\sqrt{N}>\varepsilon}\to 0$.
\end{enumerate}
\end{proposition}

The proof of Theorem~\ref{th:det-seq} then readily follows
by combining Propositions~\ref{prop:cjsq}-\ref{prop:dis-new}
and observing that the scheme $I(G_N,n)$ belongs to the class
$\CJSQ(n)$ by construction.

From the conditions of Theorem~\ref{th:det-seq} it follows that
if for all $\varepsilon > 0$, $\dis_1(G_N, \varepsilon)$
and $\dis_2(G_N, \varepsilon)$ are $o(N)$ and $o(\sqrt{N})$,
respectively, then the total number of edges in $G_N$ must be
$\omega(N)$ and $\omega(N\sqrt{N})$, respectively.
Theorem~\ref{th:bdd-deg} below states that the \emph{super-linear}
growth rate of the total number of edges is not only sufficient,
but also necessary in the sense that any graph with $O(N)$ edges is
asymptotically sub-optimal on $N$-scale.

\begin{theorem}
\label{th:bdd-deg}
Let $\GG = \{G_N\}_{N \geq 1}$ be any graph sequence, such that there
exists a fixed integer $M < \infty$ with 
\begin{equation}
\label{eq:bdd-deg}
\limsup_{N \to \infty} \dfrac{\#\big\{v \in V_N: d_v \leq M\big\}}{N} > 0,
\end{equation}
where $d_v$ is the degree of the vertex~$v$.
Then $\GG$ is sub-optimal on $N$-scale.
\end{theorem}

To prove Theorem~\ref{th:bdd-deg}, we show that starting
from an all-empty state, in finite time, a positive fraction
of servers in $G_N$ will have at least two tasks. 
This will prove that the occupancy processes when scaled by~$N$ cannot
agree with those in the sequence of cliques, and hence $\{G_N\}_{N \geq 1}$
cannot be $N$-optimal.
The idea of the proof can be explained as follows: If a system contains
$\Theta(N)$ bounded-degree vertices, then starting from an all-empty state,
in any finite time interval there will be $\Theta(N)$ servers~$u$ say,
for which all the servers in $N[u]$ have at least one task.
For all such servers an arrival at~$u$ must produce a server
with queue length two.
Thus, it shows that the instantaneous rate at which servers of queue
length two are formed is bounded away from zero, and hence $\Theta(N)$
servers of queue length two are produced in finite time. 

\paragraph{Worst-case scenario.}
Next we consider the worst-case scenario.
Theorem~\ref{th:min-deg-negative} below asserts that a graph sequence
can be sub-optimal for some $\lambda < 1$ even when the minimum degree
$d_{\min}(G_N)$ is~$\Theta(N)$.

\begin{theorem}
\label{th:min-deg-negative}
For any $\big\{d(N)\big\}_{N \geq 1}$, such that $d(N) / N\to c$
with $0 < c < 1/6$, there exists $\lambda < 1$, and a graph sequence
$\GG = \{G_N\}_{N \geq 1}$ with $d_{\min}(G_N) = d(N)$, such that
$\GG$ is sub-optimal on $N$-scale.
\end{theorem}

To construct such a sub-optimal graph sequence, consider a sequence
of complete bipartite graphs $G_N = (V_N, E_N)$,
with $V_N = A_N \sqcup B_N$ and $|A_N| / N \to c \in (0, 1/2)$
as $N \to \infty$.
If this sequence were $N$-optimal, then starting from an all-empty state,
asymptotically the fraction of servers with queue length one would
converge to~$\lambda$, and the fraction of servers with queue length
two or larger should remain zero throughout.
Now note that for large $N$ the rate at which tasks join the empty
servers in $A_N$ is given by $(1-c) \lambda$, whereas the rate
of empty server generation in $A_N$ is at most~$c$.
Choosing $\lambda > c / (1-c)$, one can see that in finite time each
server in $A_N$ will have at least one task.
From that time onward with at least instantaneous rate
$\lambda (\lambda - c) - c$, servers with queue length two start forming.
The range for $c$ stated in Theorem~\ref{th:min-deg-negative} is only to
ensure that there exists $\lambda<1$ with $\lambda (\lambda - c) - c > 0$. 

\subsection{Asymptotic optimality of random graph sequences}

Next we investigate how the load balancing process behaves on random
graph topologies. 
Specifically, we aim to understand what types of graphs are
asymptotically optimal in the presence of randomness
(i.e., in an average-case sense).
Theorem~\ref{th:inhom} below establishes sufficient conditions
for asymptotic optimality of a sequence of inhomogeneous random graphs.
Recall that a graph $G' = (V', E')$ is called a supergraph
of $G = (V, E)$ if $V = V'$ and $E \subseteq E'$.

\begin{theorem}
\label{th:inhom}
Let $\GG= \{G_N\}_{N \geq 1}$ be a graph sequence such that for each~$N$,
$G_N = (V_N, E_N)$ is a supergraph of the inhomogeneous random graph
$G_N'$ where any two vertices $u, v \in V_N$ share an edge
with probability $p_{uv}^N$, independent to each other.
\begin{enumerate}[{\normalfont (i)}]
\item If $\inf\ \{p^N_{uv}: u, v\in V_N\}$ is $\omega(1/N)$, then $\GG$ is $N$-optimal.
\item If $\inf\ \{p^N_{uv}: u, v\in V_N\}$ is $\omega(\log(N)/\sqrt{N})$, then $\GG$ is $\sqrt{N}$-optimal.
\end{enumerate}
\end{theorem}

The proof of Theorem~\ref{th:inhom} relies on Theorem~\ref{th:det-seq}.
Specifically, if $G_N$ satisfies conditions~(i) and~(ii) in
Theorem~\ref{th:inhom}, then the corresponding conditions~(i) and~(ii)
in Theorem~\ref{th:det-seq} hold.

As an immediate corollary of Theorem~\ref{th:inhom} we obtain
an optimality result for the sequence of ERRGs.
Let $\ER_N(p(N))$ denote a graph on $N$ vertices, such that any pair of vertices share an edge with probability $p(N)$.

\begin{corollary}
\label{cor:errg}
Let $\GG = \{G_N\}_{N \geq 1}$ be a graph sequence such that for each~$N$,
$G_N$ is a super-graph of $\ER_N(p(N))$, and $d(N) = (N-1) p(N)$.
{\normalfont (i)}
If $d(N)\to\infty$ as $N\to\infty$, then $\GG$ is $N$-optimal.
{\normalfont (ii)}
If $d(N)/(\sqrt{N}\log(N))\to\infty$ as $N\to\infty$, then $\GG$ is
$\sqrt{N}$-optimal.
\end{corollary}

Theorem~\ref{th:det-seq} can be further leveraged to establish the
optimality of the following sequence of random graphs.
For any $N \geq 1$ and $d(N) \leq N-1$ such that $N d(N)$ is even,
construct the \emph{erased random regular} graph on $N$~vertices
as follows:
Initially, attach $d(N)$ \emph{half-edges} to each vertex. 
Call all such half-edges \emph{unpaired}.
At each step, pick one half-edge arbitrarily, and pair it to another
half-edge uniformly at random among all unpaired half-edges
to form an edge, until all the half-edges have been paired.
This results in a uniform random regular multi-graph with degree
$d(N)$~\cite[Proposition 7.7]{remco-book-1}. 
Now the erased random regular graph is formed by erasing all the
self-loops and multiple edges, which then produces a simple graph.

\begin{theorem}
\label{th:reg}
Let $\GG = \{G_N\}_{N \geq 1}$ be a sequence of erased random regular
graphs with degree $d(N)$.
Then
{\normalfont (i)}
If $d(N)\to\infty$ as $N\to\infty$, then $\GG$ is $N$-optimal.
{\normalfont (ii)}
If $d(N)/(\sqrt{N}\log(N))\to\infty$ as $N\to\infty$, then $\GG$ is $\sqrt{N}$-optimal.
\end{theorem}

Note that due to Theorem~\ref{th:bdd-deg}, we can conclude that the
growth rate condition for $N$-optimality in Corollary~\ref{cor:errg}~(i)
and Theorem~\ref{th:reg}~(i) is not only sufficient, but necessary as well.
Thus informally speaking, $N$-optimality is achieved under the minimum
condition required as long as the underlying topology is suitably random.

\section{Token-based load balancing}
\label{token}

While a zero waiting time can be achieved in the limit by sampling
only $d(N) =o(N)$ servers as Sections~\ref{univ} and~\ref{networks} showed,
even in network scenarios, the amount of communication overhead
in terms of $d(N)$ must still grow with~$N$.
As mentioned earlier, this can be avoided by introducing memory at the
dispatcher, in particular maintaining a record of only vacant servers,
and assigning tasks to idle servers, if there are any,
or to a uniformly at random selected server otherwise.
This so-called Join-the-Idle-Queue (JIQ) scheme \cite{BB08,LXKGLG11}
can be implemented through a simple token-based mechanism generating
at most one message per task.
Remarkably enough, even with such low communication overhead,
the mean waiting time and the probability of a non-zero waiting time
vanish under the JIQ scheme in both the fluid and diffusion regimes,
as we will discuss in the next two subsections.

\subsection{Fluid-level optimality of JIQ scheme}
\label{ssec:fluidjiq}

We first consider the fluid limit of the JIQ policy.
Recall that $q_i^N(\infty)$ denotes a random variable denoting the
process $q_i^N(\cdot)$ in steady state.
Under significantly more general conditions (in the presence
of finitely many heterogeneous server pools and for general service
time distributions with decreasing hazard rate) it was proved
in~\cite{Stolyar15} that under the JIQ scheme
\begin{equation}
\label{eq:fpjiq}
q_1^N(\infty) \to \lambda, \qquad q_i^N(\infty) \to 0 \quad
\mbox{ for all } i \geq 2, \qquad \mbox{ as } \quad N \to \infty.
\end{equation}
The above equation in conjunction with the PASTA property yields that
the steady-state probability of a non-zero wait vanishes as $N \to \infty$,
thus exhibiting asymptotic optimality of the JIQ policy on fluid scale.

\paragraph{High-level proof idea.}
Loosely speaking, the proof of~\eqref{eq:fpjiq} consists of three
principal components:
\begin{enumerate}[{\normalfont (i)}]
\item Starting from an all-empty state, observe that the asymptotic
rate of increase of~$q_1$ is given by the arrival rate~$\lambda$. 
Also, the rate of decrease is~$q_1$. 
Thus, on a small time interval $\dif t$, the rate of change of~$q_1$
is given by 
\begin{equation}\label{eq:fluidjiq}
\frac{\dif q_1(t)}{\dif t} = \lambda - q_1(t).
\end{equation}
Under the above dynamics, the system occupancy states converge to the
unique fixed point of the above ODE, given by $(\lambda, 0, 0, \ldots)$.
\item The occupancy process is monotone, in the sense that
(a) Starting from an all-empty state, the occupancy process is
componentwise stochastically nondecreasing in time, and 
(b) The occupancy process at any fixed time~$t$ starting
from an arbitrary state is componentwise stochastically dominated
by the occupancy process at time~$t$ starting from an all-empty state.
\item Under the JIQ scheme, the system is stable, and hence the
occupancy process is ergodic.
Since $q_1(t)$ is the instantaneous rate of departure from the system,
ergodicity implies that in steady state there can be {\em at most}
a $\lambda$ fraction of busy servers (containing at least one task).
In fact, it further establishes that the steady-state fraction
of servers with more than one tasks vanishes asymptotically.
\end{enumerate}

Points~(i) and~(ii) above imply that starting from any state the
system must have at least a $\lambda$ fraction of busy servers,
and finally this along with Point (iii) establishes that the
steady-state occupancy process must converge to $(\lambda, 0, 0, \ldots)$. 
 
\subsection{Diffusion-level optimality of JIQ scheme}

We now turn to the diffusion limit of the JIQ scheme.
Recall the centered and scaled occupancy process
as in~\eqref{eq:diffscale}, and the Halfin-Whitt heavy-traffic regime
in~\eqref{eq:HW}.

\begin{theorem}{\normalfont (Diffusion limit for JIQ)}
\label{diffusionjiq}
Assume that $\lambda(N)$ satisfies~\eqref{eq:HW}.
Under suitable initial conditions the weak limit of the sequence of centered and diffusion-scaled occupancy process in~\eqref{eq:diffscale} coincides with that of the ordinary JSQ policy, and in particular, is given by the system of SDEs in~\eqref{eq:diffusionjsq}.
\end{theorem}

The above theorem implies that for suitable states,
on any finite time interval, the occupancy process of a system
under the JIQ policy is indistinguishable from that under the JSQ policy.

\paragraph{High-level proof idea.}
A rigorous proof of Theorem~\ref{diffusionjiq} is presented in Chapter~\ref{chap:jiq}.
The proof relies on a novel coupling
construction as described below in detail.
The idea is to compare the occupancy processes of two systems
following JIQ and JSQ policies, respectively. 
Comparing the JIQ and JSQ policies is facilitated when viewed as follows:
(i) If there is an idle server in the system, both JIQ and JSQ perform
 similarly,
(ii)~Also, when there is no idle server and only $O(\sqrt{N})$ servers
with queue length two, JSQ assigns the arriving task to a server
with queue length one. 
In that case, since JIQ assigns at random, the probability that the
task will land on a server with queue length two and thus JIQ acts
differently than JSQ is $O(1/\sqrt{N})$.
Since on any finite time interval the number of times an arrival finds
all servers busy is at most $O(\sqrt{N})$, all the arrivals except
an $O(1)$ of them are assigned in exactly the same manner in both JIQ
and JSQ, which then leads to the same scaling limit for both policies. \\

The diffusion limit result in Theorem~\ref{diffusionjiq} is in fact
true for an even broader class of load balancing schemes.
Recall that $B$ denotes the buffer capacity (possibly infinite)
of each server, and in case $B < \infty$, if a task is assigned
to a server with $B$~outstanding tasks, it is instantly discarded.
Define the class of schemes
$$\Pi^{(N)} := \{\Pi(d_0, d_1, \ldots, d_{B-1}): d_0 = N, 1 \leq d_i \leq N,
1 \leq i \leq B-1, B \geq 2\},$$ where 
in the scheme $\Pi(d_0, d_1, \ldots, d_{B-1}),$ the dispatcher assigns
an incoming task to the server with the minimum queue length among $d_k$
(possibly depending on~$N$) servers selected uniformly at random when the
minimum queue length across the system is~$k$, $k = 0, 1, \ldots, B - 1$.
The system analyzed in~\cite{EG15} (JSQ with $B = 2$) can be written
as $\Pi(N,N)$, JIQ can be expressed as $\Pi(N, 1, 1,\ldots)$,
and JIQ with a buffer capacity $B = 2$ is $\Pi(N, 1)$.

The crux of the argument in proving diffusion-level optimality for any
scheme in $\Pi^{(N)}$ goes as follows: 
First the scheme $\Pi(N, d_1, \ldots, d_{B-1})$ is sandwiched
between $\Pi(N, 1)$ and $\Pi(N, d_1)$. 
More specifically, the gap between $\Pi(N, d_1, \ldots, d_{B-1})$
and $\Pi(N, 1)$ is bounded by the number of items lost due to full buffers.
Next, this loss is bounded using the number of servers with queue length
$2$ in $\Pi(N,N)$.
This allows the use of the results in~\cite{EG15}, and yields that
on any finite time interval with high probability an $O(1)$ number of items
are lost due to full buffers, which is negligible on $\sqrt{N}$ scale.
Specifically, this shows that for suitable initial states,
the schemes $\Pi(N, 1)$ and $\Pi(N, d_1)$, along with any scheme
in the class $\Pi^{(N)}$ have the same diffusion limits in the Halfin-Whitt
heavy-traffic regime.
We conclude this subsection by describing the coupling construction
stating the stochastic inequalities, and a brief proof sketch
for Theorem~\ref{diffusionjiq}.

\paragraph{The coupling construction.}
We now construct a stochastic coupling between two systems following any two schemes, say
 $\Pi_1 = \Pi(l_0, l_1, \ldots, l_{B-1})$
and $\Pi_2 = \Pi(d_0, d_1, \ldots, d_{B'-1})$ in $\Pi^{(N)}$, respectively,
to establish the desired stochastic ordering results.
With slight abuse of notation we will denote by $\Pi_i$ the system
following scheme $\Pi_i$, $i = 1, 2$.

For the arrival process we couple the two systems as follows. First we synchronize the  arrival epochs of the two systems. Now assume that in the systems $A$ and $B$, the minimum queue lengths are $k$ and $m$, respectively, $k\leq B-1$, $m\leq B'-1$. Therefore, when a task arrives, the dispatchers in $\Pi_1$ and $\Pi_2$ have to select $l_k$ and $d_m$ servers, respectively, and then have to send the task to the one having the minimum queue length among the respectively selected servers. Since the servers are being selected uniformly at random we can assume without loss of generality, as in the stack construction, that the servers are arranged in non-decreasing order of their queue lengths and are indexed in increasing order. 
Hence, observe that when a few server indices are selected, the server having the minimum of those indices will be the server with the minimum queue length among these. Hence, in this case the dispatchers in $\Pi_1$ and $\Pi_2$ select $l_k$ and $d_m$ random numbers (without replacement) from $\{1,2,\ldots,N\}$ and then send the incoming task to the servers having indices to be the minimum of those selected numbers. To couple the decisions of the two systems, at each arrival epoch a single random permutation of $\{1,2,\ldots,N\}$ is drawn, denoted by $\boldsymbol{\Sigma}^{(N)}:=(\sigma_1, \sigma_2,\ldots,\sigma_N)$. Define $\sigma_{(i)}:= \min_{j\leq i}\sigma_j$. Then observe that system $\Pi_1$ sends the task to the server with the index $\sigma_{(l_k)}$ and system $\Pi_2$ sends the task to the server with the index $\sigma_{(d_m)}$. Since at each arrival epoch both systems use a common random permutation, they take decisions in a coupled manner.

For the potential departure process, couple the service completion times of the $k^{th}$ queue in both scenarios, $k= 1,2,\ldots,N$. More precisely, for the potential departure process assume that we have a single synchronized exp($N$) clock independent of the arrival epochs for both systems. Now when this clock rings, a number $k$ is uniformly selected from $\{1,2,\ldots,N\}$ and a potential departure occurs from the $k^{th}$ queue in both systems. If at a potential departure epoch an empty queue is selected, then we do nothing. In this way the two schemes, considered independently, still evolve according to their appropriate statistical laws.

\begin{proposition}
\label{prop: stoch_ord}
For any two schemes say, $\Pi_1=\Pi(l_0,l_1,\ldots,l_{B-1})$ and $\Pi_2=\Pi(d_0, d_1,\ldots, d_{B'-1})$ with $B\leq B'$ assume $l_0=\ldots=l_{B-2}=d_0=\ldots=d_{B-2}=d$, $l_{B-1}\leq d_{B-1}$ and either $d=N$ or $d\leq d_{B-1}$. Then the following holds:
\begin{enumerate}[{\normalfont (i)}]
\item\label{component_ordering} $\{Q^{ \Pi_1}_i(t)\}_{t\geq 0}\leq_{st}\{Q^{ \Pi_2}_i(t)\}_{t\geq 0}$ for $i=1,2,\ldots,B$
\item\label{upper bound} $\{\sum_{i=1}^B Q^{ \Pi_1}_i(t)+L^{ \Pi_1}(t)\}_{t\geq 0}\geq_{st} \{\sum_{i=1}^{B'} Q^{ \Pi_2}_i(t)+L^{ \Pi_2}(t)\}_{t\geq 0}$
\item\label{delta_ineq} $\{\Delta(t)\}_{t\geq 0}\geq \{\sum_{i=B+1}^{B'}Q_i^{ \Pi_2}(t)\}_{t\geq 0}$ almost surely under the coupling defined above,
\end{enumerate}
for any fixed $N\in\mathbbm{N}$ where $\Delta(t):=L^{ \Pi_1}(t)-L^{ \Pi_2}(t)$, provided that at time $t=0$ the above ordering holds.
\end{proposition}

\begin{proof}[Proof of Theorem~\ref{diffusionjiq}]
Let $\Pi=\Pi(N,d_1,\ldots,d_{B-1})$ be a load balancing scheme in the class $\Pi^{(N)}$. Denote by $\Pi_1$ the scheme $\Pi(N,d_1)$ with buffer size $B=2$ and let $\Pi_2$ denote the JIQ policy $\Pi(N,1)$ with buffer size $B=2$.

Observe that from Proposition~\ref{prop: stoch_ord} we have under the coupling defined above,
\begin{equation}
\label{eq: bound}
\begin{split}
|Q_i^\Pi(t) - Q_i^{\Pi_2}(t)| & \leq |Q_i^{\Pi}(t) - Q_i^{\Pi_1}(t)| +
|Q_i^{\Pi_1}(t) - Q_i^{\Pi_2}(t)| \\ & \leq
|L^{\Pi_1}(t) - L^{\Pi}(t)| + |L^{\Pi_2}(t) - L^{\Pi_1}(t)| \\
& \leq 2 L^{\Pi_2}(t),
\end{split}
\end{equation}
for all $i \geq 1$ and $t \geq 0$ with the understanding that $Q_j(t)=0$
for all $j > B$, for a scheme with buffer capacity~$B$.
The third inequality above is
due to Proposition~\ref{prop: stoch_ord}\eqref{delta_ineq},
which in particular says that
$$\{L^{\Pi_2}(t)\}_{t \geq 0} \geq \{L^{\Pi_1}(t)\}_{t \geq 0} \geq
\{L^{\Pi}(t)\}_{t \geq 0}$$ 
almost surely under the coupling.
Now we have the following lemma which we will prove below.

\begin{lemma}
\label{lem: tight}
For all $t\geq 0$, under the assumptions of Theorem~\ref{diffusionjiq},
$\{L^{\Pi_2}(t)\}_{N \geq 1}$ forms a tight sequence.
\end{lemma}

Since $L^{\Pi_2}(t)$ is non-decreasing in~$t$, the above lemma
in particular implies that
\begin{equation}
\label{eq: conv 0}
\sup_{t\in[0,T]}\frac{L^{\Pi_2}(t)}{\sqrt{N}}\pto 0.
\end{equation}
For any scheme $\Pi\in\Pi^{(N)}$, from \eqref{eq: bound} we know that 
$$\{Q_i^{\Pi_2}(t) - 2 L^{\Pi_2}(t)\}_{t \geq 0} \leq
\{Q_i^{\Pi}(t)\}_{t \geq 0} \leq \{Q_i^{\Pi_2}(t) +
2 L^{\Pi_2}(t)\}_{t \geq 0}.$$
Combining~\eqref{eq: bound} and~\eqref{eq: conv 0} shows that if the
weak limits under the $\sqrt{N}$ scaling exist, they must be the same
for all the schemes in the class $\Pi^{(N)}$. 
Also, as described in Section~\ref{sec:powerofd}, the weak limit
for $\Pi(N, N)$ exists and the common weak limit can be described
by the unique solution of the SDEs in~\eqref{eq:diffusionjsq}.
Hence, the proof of Theorem~\ref{diffusionjiq} is complete.
\end{proof}

\begin{remark}\normalfont
We have only focused on the scenario with a single dispatcher, but it is not uncommon for LBAs to operate across multiple dispatchers.
While the presence of multiple dispatchers does not affect the
queueing dynamics of JSQ($d$) policies, it does matter for the JIQ
scheme which uses memory at the dispatcher.
Scenarios with multiple dispatchers have received limited attention
in the literature, and the scant papers that exist
\cite{LXKGLG11,Mitzenmacher16,Stolyar17} almost exclusively assume
that the loads at the various dispatchers are strictly equal.
The results in~\cite{Stolyar17} in fact show that the JIQ scheme
remains asymptotically optimal even when the servers are heterogeneous,
while it is readily seen that JSQ($d$) policies cannot even be
maximally stable in that case for any fixed value of~$d$.
The case when the arrival rates at the various dispatchers are not
perfectly equal, is more delicate, and has been considered by Van der Boor {\em et al.}~\cite{BBL17a}.
\end{remark}

\subsection{Joint load balancing and auto-scaling}

Besides delay performance and implementation overhead, a further key
attribute in the context of large-scale cloud networks and data centers
is energy consumption.
So-called auto-scaling algorithms have emerged as a popular mechanism
for adjusting service capacity in response to varying demand levels
so as to minimize energy consumption while meeting performance targets,
but have mostly been investigated in settings with a centralized queue,
and queue-driven auto-scaling techniques have been widely investigated
in the literature \cite{ALW10,GDHS13,LCBWGWMH12,LLWLA11a,LLWLA11b,LLWA12,LWAT13,PP16,UKIN10,WLT12}.
In systems with a centralized queue it is very common to put servers
to `sleep' while the demand is low, since servers in sleep mode
consume much less energy than active servers.
Under Markovian assumptions, the behavior of these mechanisms can
be described in terms of various incarnations of M/M/$N$ queues
with setup times.
There are several further recent papers which examine on-demand server
addition/removal in a somewhat different vein~\cite{PS16,NS16}. 
Unfortunately, data centers and cloud networks with massive numbers
of servers are too complex to maintain any centralized queue,
as it involves a prohibitively high communication burden to obtain
instantaneous state information.

Motivated by these observations, a joint load balancing and auto-scaling strategy is proposed in Chapter~\ref{chap:energy1},
which retains the excellent delay performance and low implementation
overhead of the ordinary JIQ scheme,
and at the same time minimizes the energy consumption.
The strategy is referred to as TABS (Token-Based Auto-Balance Scaling)
and operates as follows:
\begin{itemize}
\item When a server becomes idle, it sends a `green' message to the dispatcher, waits for an $\exp(\mu)$ time (standby period), and turns itself off by sending a `red' message to the dispatcher (the corresponding green message is destroyed).
\item When a task arrives, the dispatcher selects a green message at random if there are any, and assigns the task to the corresponding server (the corresponding green message is replaced by a `yellow' message). 
Otherwise, the task is assigned to an arbitrary busy  server, and if at that arrival epoch there is a red message at the dispatcher, then it selects one at random, and the setup procedure of the corresponding server is initiated, replacing its red message by an `orange' message.
Setup procedure takes $\exp(\nu)$ time after which the server becomes active. 
\item Any server which activates due to the latter event, sends a green message to the dispatcher (the corresponding orange message is replaced), waits for an $\exp(\mu)$ time for a possible assignment of a task, and again turns itself off by sending a red message to the dispatcher.
\end{itemize}

\begin{figure}
\begin{center}
\includegraphics[scale=1]{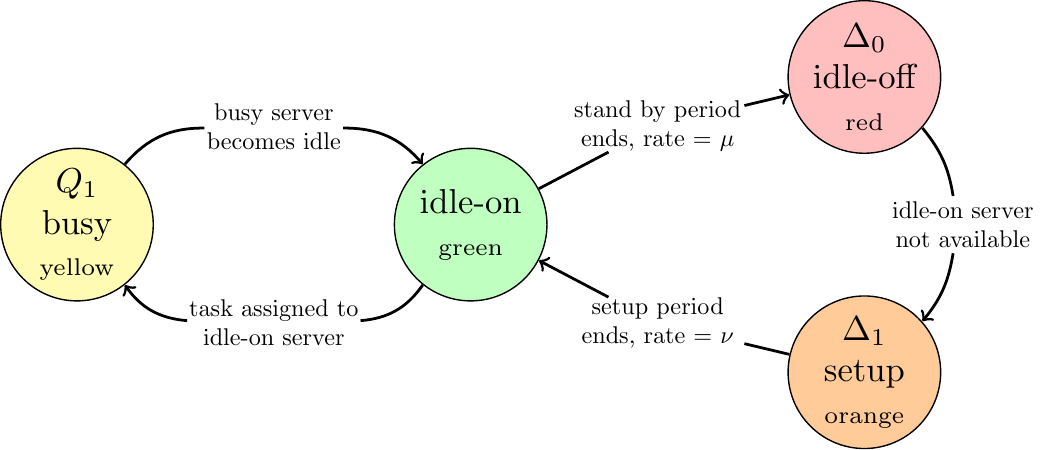}
\end{center}
\caption{Illustration of server on-off decision rules in the TABS scheme,
along with message colors and state variables.}
\label{fig:scheme}
\end{figure}

The TABS scheme gives rise to a distributed operation in which servers
are in one of four states (busy, idle-on, idle-off or standby),
and advertize their state to the dispatcher via exchange of tokens.
Figure~\ref{fig:scheme} illustrates this token-based exchange protocol.
Note that setup procedures are never aborted and continued even when
idle-on servers do become available.
Very recently dynamic scaling and load balancing with variable service
capacity and on-demand agents has been further examined in~\cite{GFPJ17}.

To describe systems under the TABS scheme,
we use the notation $\mathbf{Q}^N(t) := (Q_1^N(t), Q_2^N(t), \dots, Q_B^N(t))$
to denote the system occupancy state at time~$t$ as before.
Also, let $\Delta_0^N(t)$ and $\Delta_1^N(t)$ denote the number
of idle-off servers and servers  in setup mode at time $t$, respectively. 
The fluid-scaled quantities are denoted by the respective small letters,
\emph{viz.} $q_i^{N}(t) := Q_i^{N}(t) / N$,
$\delta_0^N(t) = \Delta_0^N(t) / N$,
and $\delta_1^N(t) = \Delta_1^N(t)/N$.
For brevity in notation, we will write
$\mathbf{q}^N(t) = (q_1^N(t),\dots,q_B^N(t))$
and $\boldsymbol{\delta}^N(t) = (\delta_0^N(t), \delta_1^N(t))$. 


\paragraph{Fluid limit.}
Under suitable initial conditions, on any finite time interval,
with probability~$1$, any sequence $\{N\}$ has a further subsequence
along which the sequence of processes $(\qq^N(\cdot), \dd^N(\cdot))$
converges to a deterministic limit $(\qq(\cdot), \dd(\cdot))$
that satisfies the following system of ODEs
\begin{equation}
\label{eq:tabsfluid}
\begin{split}
\frac{\dif^+ q_i(t)}{\dif t} &= \lambda(t) p_{i-1}(\qq(t),\dd(t),\lambda(t))
 - (q_i(t) - q_{i+1}(t)),\ i = 1, \ldots, B, \\
\frac{\dif^+\delta_0(t)}{\dif t} &= u(t) - \frac{\dif^+\xi(t)}{\dif t}, \qquad
\frac{\dif^+\delta_1(t)}{\dif t} = \frac{\dif^+\xi(t)}{\dif t} - \nu \delta_1(t),
\end{split}
\end{equation}
where by convention $q_{B+1}(\cdot) \equiv 0$, and
\begin{align*}
u(t) &= 1 - q_1(t) - \delta_0(t) - \delta_1(t), \\
\frac{\dif^+\xi(t)}{\dif t} &=
\lambda(t) (1 - p_0(\qq(t), \dd(t), \lambda(t))) \ind{\delta_0(t) > 0}.
\end{align*}
For any $(\qq, \dd)$ and $\lambda > 0$,
$(p_i(\qq, \dd, \lambda))_{i \geq 0}$ are given by 
\begin{align*}
p_0(\qq,\dd,\lambda) &= 
\begin{cases}
& 1 \qquad \text{ if } \qquad u = 1 - q_1 - \delta_0 - \delta_1 > 0, \\
& \min\{\lambda^{-1}(\delta_1 \nu + q_1 - q_2), 1\},\quad \text{otherwise,}
\end{cases} \\
\quad p_i(\qq, \dd, \lambda) &=
(1 - p_0(\qq, \dd, \lambda)) (q_i - q_{i+1})q_1^{-1},\ i = 1, \ldots, B.
\end{align*}

We now provide an intuitive explanation of the fluid limit stated above.
The term $u(t)$ corresponds to the asymptotic fraction of idle-on
servers in the system at time~$t$, and $\xi(t)$ represents the
asymptotic cumulative number of server setups (scaled by~$N$)
that have been initiated during $[0,t]$.
The coefficient $p_i(\qq, \dd, \lambda)$ can be interpreted as the
instantaneous fraction of incoming tasks that are assigned to some
server with queue length~$i$, when the fluid-scaled occupancy state is
$(\qq, \dd)$ and the scaled instantaneous arrival rate is~$\lambda$.
Observe that as long as $u > 0$, there are idle-on servers,
and hence all the arriving tasks will join idle servers. 
This explains that if $u > 0$, $p_0(\qq, \dd, \lambda) = 1$
and $p_i(\qq, \dd, \lambda)=0$ for $i = 1, \ldots, B-1$.
If $u = 0$, then observe that servers become idle at rate $q_1 - q_2$,
and servers in setup mode turn on at rate $\delta_1 \nu$.
Thus the  idle-on servers are created at a total rate
$\delta_1 \nu + q_1 - q_2$.
If this rate is larger than the arrival rate~$\lambda$, then almost
all the arriving tasks can be assigned to idle servers.
Otherwise, only a fraction $(\delta_1 \nu + q_1 - q_2) / \lambda$
of arriving tasks join idle servers. 
The rest of the tasks are distributed uniformly among busy servers,
so a proportion $(q_i - q_{i+1}) q_1^{-1}$ are assigned to servers
having queue length~$i$.
For any $i = 1, \ldots, B$, $q_i$ increases when there is an arrival
to some server with queue length $i - 1$, which occurs at rate
$\lambda p_{i-1}(\qq, \dd, \lambda)$, and it decreases when there is
a departure from some server with queue length~$i$, which occurs
at rate $q_i - q_{i-1}$. 
Since each idle-on server turns off at rate~$\mu$, the fraction
of servers in the off mode increases at rate $\mu u$.
Observe that if $\delta_0 > 0$, for each task that cannot be assigned
to an idle server, a setup procedure is initiated  at one idle-off server. 
As noted above, $\xi(t)$ captures the (scaled) cumulative number
of setup procedures initiated up to time~$t$.
Therefore the fraction of idle-off servers and the fraction of servers
in setup mode decreases and increases by $\xi(t)$, respectively,
during $[0,t]$.
Finally, since each server in setup mode becomes idle-on at rate~$\nu$,
the fraction of servers in setup mode decreases at rate $\nu \delta_1$.

\paragraph{Fixed point and global stability.}
In case of a constant arrival rate $\lambda(t) \equiv \lambda < 1$,
any fluid sample path in~\eqref{eq:tabsfluid} has a unique fixed point:
\begin{equation}
\label{eq:fixed point-intro}
\delta_0^* = 1 - \lambda, \qquad \delta_1^* = 0, \qquad
q_1^* = \lambda \quad \mbox{ and } \quad q_i^*=0,
\end{equation}
for $i = 2, \ldots, B$.
Indeed, it can be verified that $p_0(\qq^*, \dd^*, \lambda) = 1$
and $u^* = 0$ for $(\qq^*, \dd^*)$ given by~\eqref{eq:fixed point-intro}
so that the derivatives of $q_i$, $i = 1, \dots, B$, $\delta_0$,
and $\delta_1$ become zero, and that these cannot be zero at any other
fluid-scaled occupancy state.
Note that, at the fixed point, a fraction $\lambda$ of the servers
have exactly one task while the remaining fraction have zero tasks,
independently of the values of the parameters~$\mu$ and~$\nu$.

In order to establish the convergence of the sequence of steady states,
we need the global stability of the fluid limit, i.e., starting
from any fluid-scaled occupancy state, any fluid sample path described
by~\eqref{eq:tabsfluid} converges to the unique fixed
point~\eqref{eq:fixed point-intro} as $t \to \infty$.
More specifically, irrespective of the starting state,
\begin{equation}
\label{eq:globalstab}
(\qq(t), \dd(t)) \to (\qq^*, \dd^*), \quad \mbox{ as } \quad t \to \infty,
\end{equation}
where $(\qq^*, \dd^*)$ is as defined in~\eqref{eq:fixed point-intro}.

\paragraph{Interchange of limits.}
The global stability can be leveraged to show that the steady-state
distribution of the $N^{\mathrm{th}}$ system, for large $N$, can be
well approximated by the fixed point of the fluid limit
in~\eqref{eq:fixed point-intro}. 
Specifically, it justifies the interchange of the many-server
($N \to \infty$) and stationary ($t \to \infty$) limits.
Since the buffer capacity~$B$ at each server is supposed to be finite,
for every~$N$, the Markov process $(\QQ^N(t), \Delta_0^N(t), \Delta_1^N(t))$ is irreducible, has a finite state space, and thus has a unique
steady-state distribution.
Let $\pi^N$ denote the steady-state distribution of the
$N^{\mathrm{th}}$ system, i.e.,
$$\pi^{N}(\cdot) = \lim_{t \to \infty}
\mathbb{P}\ \big(\qq^{N}(t) = \cdot, \dd^N(t) = \cdot\big).$$ 
The fluid limit result and the global stability thus yield that
$\pi^N$ converges weakly to~$\pi$ as $N \to \infty$, where $\pi$ is
given by the Dirac mass concentrated at $(\qq^*,\dd^*)$ defined
in~\eqref{eq:fixed point-intro}.
\begin{remark}\normalfont
Note that the above interchange of limits result was obtained \emph{under the assumption that the queues have finite buffers}, and analysis of the infinite-buffer scenario was left open.
The key challenge in the latter case stems from the fact that the system stability under the usual subcritical load assumption is not automatic.
In fact as explained in Chapter~\ref{chap:energy2}, when the number of servers $N$ is fixed, the stability may {\em not} hold even under a subcritical load assumption.
In Chapter~\ref{chap:energy2} the stability issue of the TABS scheme is addressed and the convergence of the sequence of steady states will be shown for the infinite-buffer scenario.
In particular, it will be established that for a fixed choice of parameters $\lambda<1$, $\mu>0$, and $\nu>0$, the system with $N$ servers under the TABS scheme is stable for large enough $N$.
There we introduce an induction-based approach which uses both the conventional fluid limit (in the sense of a large starting state) and the mean-field fluid limit (when $N\to\infty$) in an intricate fashion to prove the large-$N$ stability of the system.
\end{remark}

\paragraph{Performance metrics.}

As mentioned earlier, two key performance metrics are the expected
waiting time of tasks $\expt[W^N]$ and energy consumption $\expt[P^N]$
for the $N^{\mathrm{th}}$ system in steady state.
In order to quantify the energy consumption, we assume that the energy
usage of a server is $P_{\full}$ when busy or in set-up mode,
$P_{\idle}$ when idle-on, and zero when turned off.
Evidently, for any value of~$N$, at least a fraction $\lambda$ of the
servers must be busy in order for the system to be stable,
and hence $\lambda P_{\full}$ is the minimum mean energy usage per
server needed for stability.
We will define $\expt[Z^N] = \expt[P^N] - \lambda P_{\full}$ as the
relative energy wastage accordingly.
The interchange of limits result can be leveraged to obtain that
asymptotically the expected waiting time and energy consumption
for the TABS scheme vanish in the limit, for any strictly positive
values of~$\mu$ and~$\nu$.
More specifically, for a constant arrival rate 
$\lambda(t) \equiv \lambda < 1$, for any $\mu > 0$, $\nu > 0$,
as $N \to \infty$,
\begin{enumerate}[{\normalfont (a)}]
\item Zero mean waiting time: $\expt[W^N] \to 0$,
\item Zero energy wastage: $\expt[Z^N] \to 0$.
\end{enumerate}
The key implication is that the TABS scheme, while only involving
constant communication overhead per task, provides performance
in a distributed setting that is as good at the fluid level as can
possibly be achieved, even in a centralized queue, or with unlimited
information exchange.

\paragraph{Comparison to ordinary JIQ policy.}
Consider a constant arrival rate $\lambda(t) \equiv \lambda$. 
It is worthwhile to observe that the component $\qq$ of the fluid limit
as in~\eqref{eq:tabsfluid} coincides with that for the ordinary JIQ
policy where servers always remain on, when the system following the
TABS scheme starts with all the servers being idle-on,
and $\lambda + \mu < 1$. 
To see this, observe that the component $\qq$ depends on~$\dd$ only
through $(p_{i-1}(\qq, \dd))_{i \geq 1}$. 
Now, $p_0 = 1$, $p_i = 0$, for all $i \geq 1$,
whenever $q_1 + \delta_0 + \delta_1 < 1$, irrespective of the precise
values of $(\qq, \dd)$. 
Moreover, starting from the above initial state, $\delta_1$ can
increase only when $q_1 + \delta_0 = 1$. 
Therefore, the fluid limit of $\qq$ in~\eqref{eq:tabsfluid} and the
ordinary JIQ scheme are identical if the system parameters
$(\lambda, \mu, \nu)$ are such that $q_1(t) + \delta_0(t) < 1$,
for all $t \geq 0$.
Let $y(t) = 1 - q_1(t) - \delta_0(t)$.
The solutions to the differential equations
\[
\frac{\dif q_1(t)}{\dif t} = \lambda - q_1(t), \quad
\frac{\dif y(t)}{\dif t} = q_1(t) - \lambda - \mu y(t),
\]
$y(0) = 1$, $q_1(0) = 0$ are given by 
\[
q_1(t) = \lambda (1 - \e^{- t}), \quad
y(t) = \frac{\e^{- (1+\mu) t}}{\mu-1} \big(\e^t(\lambda+\mu-1) -
\lambda \e^{\mu t}\big).
\]
Notice that if $\lambda + \mu<1$, then $y(t) > 0$ for all $t \geq 0$
and thus, $q_1(t) + \delta_0(t) < 1$, for all $t \geq 0$.
The fluid-level optimality of the JIQ scheme was described
in Subsection~\ref{ssec:fluidjiq}. 
This observation thus establishes the optimality of the fluid-limit
trajectory under the TABS scheme for suitable parameter values
in terms of response time performance.
From the energy usage perspective, under the ordinary JIQ policy,
since the asymptotic steady-state fraction of busy servers ($q_1^*$)
and idle-on servers are given by $\lambda$ and $1-\lambda$, respectively,
the asymptotic steady-state (scaled) energy usage is given by 
\begin{align*}
\expt[P^{\mathrm{JIQ}}] = \lambda P_{\full} + (1-\lambda) P_{\idle} =
\lambda P_{\full}(1 + (\lambda^{-1}-1)f),
\end{align*}
where $f = P_{\idle}/P_{\full}$ is the relative energy consumption
of an idle server.
As described earlier, the asymptotic steady-state (scaled) energy
usage under the TABS scheme is $\lambda P_{\full}$.
Thus the TABS scheme reduces the asymptotic steady-state energy usage
by $\lambda P_{\full}(\lambda^{-1}-1)f = (1-\lambda)P_{\idle}$, which
amounts to a relative saving of $(\lambda^{-1}-1)f/(1+(\lambda^{-1}-1)f)$.
In summary, the TABS scheme performs as good as the ordinary JIQ policy
in terms of the waiting time and communication overhead while providing
a significant energy saving.

\section[Redundancy and alternative scaling regimes]{Redundancy policies and alternative scaling regimes}
\label{miscellaneous}

In this section we discuss somewhat related redundancy policies,
 alternative scaling regimes, and some additional performance metrics of interest.

\subsection{Redundancy-d policies}

So-called redundancy-$d$ policies involve a somewhat similar operation
as JSQ($d$) policies, and also share the primary objective of ensuring
low delays \cite{AGSS13,VGMSRS13}.
In a redundancy-$d$ policy, $d \geq 2$ candidate servers are selected
uniformly at random (with or without replacement) for each arriving task,
just like in a JSQ($d$) policy.
Rather than forwarding the task to the server with the shortest queue
however, replicas are dispatched to all sampled servers.
Note that the initial replication to $d$ servers selected uniformly
at random does not entail any communication burden, but the abortion
of redundant copies at a later stage does involve a significant amount
of information exchange and complexity.

Two common options can be distinguished for abortion of redundant clones.
In the first variant, as soon as the first replica starts service,
the other clones are abandoned.
In this case, a task gets executed by the server which had the smallest
workload at the time of arrival (and which may or may not have had the
shortest queue length) among the sampled servers.
This may be interpreted as a power-of-$d$ version of the Join-the-Smallest
Workload (JSW) policy discussed in Subsection~\ref{ssec:jsw}.
The optimality properties of the JSW policy mentioned in that subsection
suggest that redundancy-$d$ policies should outperform JSQ($d$) policies,
which appears to be supported by simulation experiments,
but has not been established by analytical comparisons so far.

In the second option the other clones of the task are not aborted
until the first replica has completed service (which may or may not
have been the first replica to start service).
While a task is only handled by one of the servers in the former case, 
it may be processed by several servers in the latter case.
When the service times are exponentially distributed and independent for
the various clones, the aggregate amount of time spent by all the servers
until completion remains exponentially distributed with the same mean.
An exact analysis of the delay distribution in systems with $N = 2$
or $N = 3$ servers is provided in \cite{GZDHHS15,GZDHHS16},
and exact expressions for the mean delay with an arbitrary number
of servers are established in~\cite{GZHS16}.
The limiting delay distribution in a fluid regime with $N \to \infty$
is derived in \cite{GZVHS16,GHSVZ17}
based on an asymptotic independence assumption among the servers.
In general, the mean aggregate amount of time devoted to a task
and the resulting delay may be larger or smaller for less or more
variable service time distributions, also depending on the number
of replicas per task \cite{PC16,SLR16,WJW14,WJW15}.
In particular, for heavy-tailed service time distributions, the mean
aggregate time spent on a task may be considerably reduced by virtue
of the redundancy.
Indeed, even if the first replica to start service has an extremely long
service time, that is not likely to be case for the other clones as well.
In spite of the extremely long service time of the first replica,
it is therefore unlikely for the aggregate amount of time spent on the
task or its waiting time to be large.
This provides a significant performance benefit to redundancy-$d$
policies over JSQ($d$) policies, and has also motivated a strong
interest in adaptive replication schemes \cite{APS17a,Joshi17a, Joshi17b}.

A further closely related model is where $k$ of the replicas need to
complete service, $1 \leq k \leq d$, in order for the task to finish 
which is relevant in the context of storage systems with coding
and MapReduce tasks \cite{JSW15a,JSW17}.
The special case where $k = d = N$ corresponds to a classical
fork-join system.

\subsection{Conventional heavy traffic}
\label{classical}

In this subsection we briefly discuss a few asymptotic results for LBAs
in the classical heavy-traffic regime as described in Subsection~\ref{asym}
where the number of servers~$N$ is fixed and the relative load tends
to one in the limit.

The papers \cite{Foschini77,FS78,Reiman84,ZHW95} establish diffusion
limits for the JSQ policy in a sequence of systems with Markovian
characteristics as in our basic model set-up, but where in the $K$-th
system the arrival rate is $K \lambda + \hat\lambda \sqrt{K}$, while
the service rate of the $i$-th server is $K \mu_i + \hat\mu_i \sqrt{K}$,
$i = 1, \dots, N$, with $\lambda = \sum_{i = 1}^{N} \mu_i$,
inducing critical load as $K \to \infty$.
It is proved that for suitable initial conditions the queue lengths are
of the order O($\sqrt{K}$) over any finite time interval and exhibit
a state-space collapse property.
In particular, a properly scaled version of the joint queue length
process lives in a one-dimensional rather than $N$-dimensional space,
reflecting that the various queue lengths evolve in lock-step,
with the relative proportions remaining virtually identical in the limit,
while the aggregate queue length varies.

Atar {\em et al.}~\cite{AKM17} investigate a similar scenario,
and establish diffusion limits for three policies: the JSQ($d$) policy,
the redundancy-$d$ policy (where the redundant clones are abandoned
as soon as the first replica starts service), and a combined policy
called Replicate-to-Shortest-Queues (RSQ) where $d$ replicas are
dispatched to the $d$-shortest queues.
Note that the latter policy requires instantaneous knowledge of all
the queue lengths, and hence involves a similar excessive communication
overhead as the ordinary JSQ policy, besides the substantial information
exchange associated with the abortion of redundant copies.
Conditions are derived for the values of the relative service rates
$\mu_i$, $i = 1, \dots, N$, in conjunction with the diversity
parameter~$d$, in order for the queue lengths under the JSQ($d$)
and redundancy-$d$ policies to be of the order O($\sqrt{K}$) over any
finite time interval and exhibit state-space collapse.
The conditions for the two policies are distinct, but in both cases
they are weaker for larger values of~$d$, as intuitively expected.
While the conditions for the values of~$\mu_i$ depend on~$d$,
whenever they are met, the actual diffusion-scaled queue length
processes do not depend on the exact value of~$d$ in the limit,
showing a certain resemblance with the universality property
as identified in Subsection~\ref{ssec:powerd} for the large-capacity
and Halfin-Whitt regimes.

Zhou {\em et al.} \cite{ZWTSS17} consider a slightly different model
set-up with a time-slotted operation, and identify a class~$\Pi$ of LBAs
that not only provide throughput-optimality (or maximum stability,
i.e., keep the queues stable in a suitable sense whenever feasible to
do so at all), but also achieve heavy-traffic delay optimality, in
the sense that the properly scaled aggregate queue length is the same
as that in a centralized queue where all the resources are pooled
as the load tends to one.
As it turns out, the class~$\Pi$ includes JSQ($d$) policies with $d \geq 2$,
but does \emph{not} include the JIQ scheme, which tends to degenerate 
into a random assignment policy when idle servers are rarely available.
The authors further propose a threshold-based policy which has low
implementation complexity like the JIQ scheme, but \emph{does} belong
to the class~$\Pi$, and hence achieves heavy-traffic delay optimality.

\subsection{Non-degenerate slowdown}
\label{nondegenerate}

In this subsection we briefly discuss a few of the scarce asymptotic
results for LBAs in the so-called non-degenerate slow-down regime
described in Subsection~\ref{asym} where $N - \lambda(N) \to \gamma > 0$,
as the number of servers~$N$ grows large.
We note that in a centralized queue the process tracking the evolution
of the number of waiting tasks, suitably accelerated and normalized
by~$N$, converges in this regime to a Brownian motion with drift
$- \gamma$ reflected at zero as $N \to \infty$.
In stationarity, the number of waiting tasks, normalized by~$N$,
converges in this regime to an exponentially distributed random
variable with parameter~$\gamma$ as $N \to \infty$.
Hence, the mean number of waiting tasks must be at least of the order
$N / \gamma$, and the waiting time cannot vanish as $N \to \infty$
under any policy.

Gupta \& Walton~\cite{GW17} characterize the diffusion-scaled queue
length process under the JSQ policy in this asymptotic regime.
They further compare the diffusion limit for the JSQ policy with that
for a centralized queue as described above as well as several LBAs
such as the JIQ scheme and a refined version called Idle-One-First (I1F),
where a task is assigned to a server with exactly one task if no idle
server is available and to a randomly selected server otherwise.

It is proved that the diffusion limit for the JIQ scheme is no longer
asymptotically equivalent to that for the JSQ policy in this asymptotic
regime, and the JIQ scheme fails to achieve asymptotic optimality
in that respect, as opposed to the behavior in the large-capacity
and Halfin-Whitt regimes discussed in Subsection~\ref{ssec:jiq}.
In contrast, the I1F scheme does preserve the asymptotic equivalence
with the JSQ policy in terms of the diffusion-scaled queue length process,
and thus retains asymptotic optimality in that sense.

These results provide further indication that the amount and accuracy
of queue length information needed to achieve asymptotic equivalence
with the JSQ policy depend not only on the scale dimension (e.g.~fluid
or diffusion), but also on the load regime.
Put differently, the finer the scale and the higher the load,
the more strictly one can distinguish various LBAs in terms
of the relative performance compared to the JSQ policy.

\subsection{Scaling of maximum queue length}
\label{ballsbins}

So far we have focused on the asymptotic behavior of LBAs in terms
of the number of servers with a certain queue length, either on fluid
scale or diffusion scale, in various regimes as $N \to \infty$.
A related but different performance metric is the maximum queue length
$M(N)$ among all servers as $N \to \infty$.
Luczak \& McDiarmid~\cite{LM06} showed that 
for a fixed $d\geq 2$ the steady-state maximum queue length 
$M(N)$ in a system under JSQ($d$) policy is given by $\log(\log(N)) / \log(d) + O(1)$ and is concentrated on at most two adjacent values, whereas 
for purely random assignment ($d=1$), it scales as $\log(N) / \log(1/\lambda)$ and does not concentrate on a bounded range of values.
This is yet a further manifestation of the power-of-choice effect.

The maximum queue length $M(N)$ is in fact the central performance
metric in balls-and-bins models where arriving items (balls) do not
get served and never depart but simply accumulate in bins,
and (stationary) queue lengths are not meaningful.
In fact, the very notion of randomized load balancing and power-of-$d$
strategies was introduced in a balls-and-bins setting in the seminal
paper by Azar {\em et al.}~\cite{ABKU94}.
Several further variations and extensions in that context have been considered in \cite{V99,ACMR95,BCSV00,BCSV06, CMS95,DM93,FPSS05,PR04, P05, G81}.

As alluded to above, there are natural parallels between the
balls-and-bins setup and the queueing scenario that we have focused
on so far.
These commonalities are for example reflected in the fact that
power-of-$d$ strategies yield similar dramatic performance improvements
over purely random assignment in both settings.

However, there are also quite fundamental differences between the
balls-and-bins setup and the queueing scenario, besides the obvious
contrasts in the performance metrics.
The distinction is for example evidenced by the fact that a simple
Round-Robin strategy produces a perfectly balanced allocation in
a balls-and-bins setup but is far from optimal in a queueing scenario
as observed in Subsection~\ref{random}.
In particular, the stationary fraction of servers with two or more
tasks under a Round-Robin strategy remains positive in the limit
as $N \to \infty$, whereas it vanishes under the JSQ policy.
Furthermore, it should also be noted~\cite{LM05} that the maximum number of balls in a bin under the purely random assignment policy scales as $\log(N)/\log(\log(N))$ and is concentrated on two adjacent values, which is again in contrast with the queueing scenario.
On a related account, since tasks get served and eventually depart
in a queueing scenario, less balanced allocations with a large
portion of vacant servers will generate fewer service completions
and result in a larger total number of tasks.
Thus different schemes yield not only various degrees of balance,
but also variations in the aggregate number of tasks in the system,
which is not the case in a balls-and-bins set-up.

\section{Extensions}\label{sec:ext}

Throughout most of the chapter we have focused on the supermarket model
as a canonical setup and adopted several common assumptions in that context:
(i) all servers are identical;
(ii) the service requirements are exponentially distributed;
(iii) no advance knowledge of the service requirements is available;
(iv) in particular, the service discipline at each server is oblivious
to the actual service requirements.
As mentioned earlier, the stochastic optimality of the JSQ policy,
and hence its central role as an ideal performance benchmark,
critically rely on these assumptions.
The latter also broadly applies to the stochastic coupling techniques
and asymptotic universality properties that we have considered
in the previous sections.
In this section however we review some results for scenarios
where these assumptions are relaxed, in particular allowing for general
service requirement distributions and possibly heterogeneous servers,
along with some broader methodological issues.
In Subsection~\ref{ssec:powerofd-gen} we focus on the behavior of JSQ($d$)
policies in such scenarios, mainly in the large-$N$ limit,
while also briefly commenting on the JIQ policy.
In Subsection~\ref{anticipating} we discuss strategies which specifically exploit
knowledge of server speeds or service requirements of arriving tasks
in making task assignment decisions, and may not necessarily use
queue length information, mostly in a finite-$N$ regime.

\subsection{JSQ(d) policies with general service requirement
distributions}
\label{ssec:powerofd-gen}

Foss \& Chernova~\cite{FC91,FC98} use direct probabilistic methods
and fluid limits to obtain stability conditions
for finite-size systems with a renewal arrival process, a FCFS discipline
at each server, various state-dependent routing policies, including JSQ,
and general service requirement distributions, which may depend on the
task type, the server or both.
Using fluid limits as well as Lyapunov functions,
Bramson~\cite{Bramson98,Bramson11} shows that JSQ($d$) policies achieve
stability for any subcritical load in finite-size systems
with a renewal arrival process, identical servers, non-idling local
service disciplines and general service requirement distributions.
In addition, he derives uniform bounds on the tails of the marginal
queue length distributions, and uses these to prove relative
compactness of these distributions.

Bramson {\em et al.}~\cite{BLP10,BLP12} examine mean-field limits
for JSQ($d$) policies with generally distributed service requirements,
leveraging the above-mentioned tail bounds and relative compactness.
They establish that similar ``power-of-choice'' benefits occur
as originally demonstrated for exponentially distributed service
requirements in the work of Mitzenmacher~\cite{Mitzenmacher96}
and Vvedenskaya {\em et al.}~\cite{VDK96}, provided a certain `ansatz'
holds asserting that finite subsets of queues become independent in the
large-$N$ limit.
The latter `propagation of chaos' property is shown to hold in several
settings, e.g. when the service requirement distribution has
a decreasing hazard rate and the discipline at each server is FCFS
or when the service requirement distribution has a finite second moment
and the load is sufficiently low.
The ansatz also always holds for the power-of-$d$ version of the JSW
rather than JSQ policy.

It is further shown in~\cite{BLP10,BLP12} that the arrival process
at any given server tends to a state-dependent Poisson process in the
large-$N$ limit, and that the queue length distribution becomes insensitive
with respect to the service requirement distribution when the service
discipline is either Processor Sharing or LCFS with preemptive resume.
This may be explained from the insensitivity property of queues
with state-dependent Poisson arrivals and symmetric service disciplines.

There are strong plausibility arguments that a similar asymptotic
insensitivity property should hold for the JIQ policy in a queueing
scenario, even if the discipline at each server is not symmetric
but FCFS for example.
So far, however, this has only been rigorously established for service
requirement distributions with decreasing hazard rate in~\cite{Stolyar15}.
This result was in fact proved for systems with heterogeneous server
pools, and was further extended in~\cite{Stolyar17} to systems
with multiple symmetric dispatchers.
As it turns out, general service requirement distributions
with an increasing hazard rate give rise to major technical challenges
due to a lack of certain monotonicity properties.
This has only allowed a proof of the asymptotic zero-wait property
for the JIQ policy for load values strictly below $1/2$ so far~\cite{FS17}.

A fundamental technical issue associated with any general service
requirement distribution is that the joint queue length no longer
provides a suitable state description, and that the state space
required for a Markovian description is no longer countable.
Aghajani \& Ramanan~\cite{AR17} and Aghajani {\em et al.}~\cite{ALR18}
introduce a particle representation for the state of the system
and describe the state dynamics for a JSQ($d$) policy via a sequence
of interacting measure-valued processes.
They prove that as $N$ grows large, a suitably scaled sequence
of state processes converges to a hydrodynamic limit which is
characterized as the unique solution of a countable system of coupled
deterministic measure-valued equations, i.e., a system of PDE rather
than the usual ODE equations.
They also establish a `propagation of chaos' result, meaning that
finite collections of queues are asymptotically independent.

Mukhopadhyay \& Mazumdar~\cite{MM14,MM16}
and Mukhopadhyay {\em et al.}~\cite{MKM16} analyzed the performance
and stability of static probabilistic routing strategies and power-of-$d$
policies in the large-$N$ limit in systems with exponential service
requirement distributions, but heterogeneous server pools
and a Processor-Sharing discipline at each server.
They also considered variants of the JSQ($d$) policy which account
for the server speed in the selection criterion as well as hybrid
combinations of the JSQ($d$) policy with static probabilistic routing.
Related results for heterogeneous loss systems rather than queueing
scenarios are presented in \cite{KMM17,MMG15,MKMG15}.
As the results in~\cite{MM14,MM16} reflect, ordinary JSQ($d$) policies
may fail to sample the faster servers sufficiently often in such scenarios,
and therefore fail to achieve maximum stability, let alone
asymptotic optimality.
In~\cite{MKM16} a weighted version of JSQ($d$) policies is presented that
does provide maximum stability, without requiring any specific knowledge
of the underlying system parameters and server speeds in particular.

Vasantam {\em et al.} \cite{VMM17a,VMM17b} examine mean-field limits
for power-of-$d$ policies in many-server loss systems with phase-type
service requirement distributions.
They observe that the fixed point suggests a similar insensitivity
property of the stationary occupancy distribution as mentioned above.
In view of the insensitivity of loss systems with possibly
state-dependent Poisson arrivals, this may be interpreted as an indirect
indication that the arrival process at any given server pool tends
to a state-dependent Poisson arrival process in the large-$N$ limit.
In a somewhat different strand of work, Jonckheere \& Prabhu~\cite{JP16}
investigate the behavior of blocking probabilities in various load
regimes in systems with many single-server finite-buffer queues,
a Processor-Sharing discipline at each server, and an \emph{insensitive}
routing policy.

\subsection{Heterogeneous servers and knowledge of service requirements}
\label{anticipating}

The bulk of the literature has focused on systems
with identical servers, and scenarios with non-identical server speeds
have received relatively limited attention.
A natural extension of the JSQ policy is to assign jobs to the server
with the normalized shortest queue length, or equivalently,
assuming exponentially distributed service requirements, the shortest
expected delay.
While such a Generalized JSQ (GJSQ) or Shortest Expected Delay (SED)
strategy tends to perform well~\cite{BZ89}, it is not strictly
optimal in general~\cite{EVW80}, and the true optimal strategy may in
fact have a highly complicated structure.
Selen {\em et al.}~\cite{SAK16} present approximations for the
performance of GJSQ policies in a finite-$N$ regime with generally
distributed service requirements and a Processor-Sharing discipline
at each server, extending the analysis in Gupta {\em et al.}~\cite{GHSW07}
for the ordinary JSQ policy with homogeneous servers.

In a separate line of work, Feng {\em et al.}~\cite{FMR05} consider
static dispatching policies in a finite-$N$ regime with heterogeneous
servers and a FCFS or Processor-Sharing discipline at each server.
The assignment decision may depend on the service requirement of the
arriving task, but not on the actual queue lengths or any other
state information.
In case of FCFS the optimal routing policy is shown to have a nested
size interval structure, generalizing the strict size interval
structure of the task assignment strategies in Harchol-Balter
{\em et al.}~\cite{HCM99} which are optimal for homogeneous servers.
In case of Processor Sharing, the knowledge of the service requirements
of arriving tasks is irrelevant, in the absence of any state information.

Altman {\em et al.}~\cite{AAP11} consider static probabilistic routing
policies in a somewhat similar setup of a finite-$N$ regime
with multiple task types, servers with heterogeneous speeds,
and a Processor-Sharing discipline at each server.
The routing probabilities are selected so as to either minimize the
global weighted holding cost or the expected holding cost for
an individual task, and may depend on the type of the task and its
service requirement, but not on any other state information.

When knowledge of the service requirements of arriving tasks is available,
it is natural to exploit that for the purpose of local scheduling
at the various servers, and for example use size-based disciplines.
The impact of the local scheduling discipline and server heterogeneity
on the performance and degree of efficiency of load balancing strategies
is examined in~\cite{CMW09}.
An interesting broader issue concerns the relative benefits provided
by exploiting knowledge of service requirements of arriving tasks versus
using information on queue lengths or workloads at the various servers,
which strongly depend on the service requirement distribution~\cite{HSY09}.

%% file: univjsqd.tex
\begin{abstract}
We consider a system of $N$~parallel single-server queues with unit-mean exponential service requirements and a single dispatcher where tasks arrive as a Poisson process of rate $\lambda(N)$. When a task arrives, the dispatcher assigns it to a server with the shortest queue among $d(N)$ randomly selected servers ($1 \leq d(N) \leq N$). This load balancing strategy is referred to as a JSQ($d(N)$) scheme, marking that it subsumes the celebrated Join-the-Shortest Queue (JSQ) policy as a crucial special case for $d(N) = N$.

We construct a stochastic coupling to bound the difference in the queue length processes between the JSQ policy and a JSQ($d(N)$) scheme with an arbitrary value of~$d(N)$. We use the coupling to derive the fluid limit in the regime where $\lambda(N) / N \to \lambda < 1$ as $N \to \infty$ with $d(N) \to\infty$, along with the associated fixed point. The fluid limit turns out not to depend on the exact growth rate of $d(N)$, and in particular coincides with that for the JSQ policy. We further leverage the coupling to establish that the diffusion limit in the critical regime where $(N - \lambda(N)) / \sqrt{N} \to \beta > 0$ as $N \to \infty$ with $d(N)/(\sqrt{N} \log (N))\to\infty$ corresponds to that for the JSQ policy. These results indicate that the optimality of the JSQ policy can be preserved at the fluid-level and diffusion-level while reducing the overhead by nearly a factor~O($N$) and O($\sqrt{N}/\log(N)$), respectively.
\end{abstract}



\section{Model description and main results}\label{sec: model descr}
In this chapter we establish a universality property for a broad
class of randomized load balancing schemes in many-server systems as  described in Section~\ref{univ}.

Consider a system with $N$~parallel single-server queues with identical servers
and a single dispatcher.
Tasks with unit-mean exponential service requirements arrive at the
dispatcher as a Poisson process of rate $\lambda(N)$,
and are instantaneously forwarded to one of the servers.
Specifically, when a task arrives, the dispatcher assigns it
to a server with the shortest queue among $d(N)$ randomly selected
servers ($1 \leq d(N) \leq N$).
This load balancing strategy will
be referred to as the JSQ($d(N)$) scheme, marking that it subsumes
the ordinary JSQ policy as a crucial special case for $d(N) = N$.
The buffer capacity at each of the servers is~$b$ (possibly infinite),
and when a task is assigned to a server with $b$~pending tasks,
it is permanently discarded.

For any $d(N)$ ($1 \leq d(N) \leq N$), let 
$$\QQ^{\sss d(N)}(t): =
\left(Q_1^{\sss d(N)}(t), Q_2^{\sss d(N)}(t), \dots, Q_b^{\sss d(N)}(t)\right)$$ 
denote the
system occupancy state, where $Q_i^{\sss d(N)}(t)$ is the number of servers
under the JSQ($d(N)$) scheme with a queue length of~$i$ or larger,
at time~$t$, including the possible task in service, $i = 1, \dots, b$, recall Figure~\ref{figB} in Chapter~\ref{chap:introduction}. 
Throughout we assume that at each arrival epoch the servers are ordered
in  nondecreasing order of their queue lengths (ties can be broken arbitrarily), and whenever we refer
to some ordered server, it should be understood with respect to this prior ordering.

We occasionally omit the superscript $d(N)$, and replace it by~$N$, to refer to
the $N^{\mathrm{th}}$ system, when the value of $d(N)$ is clear from the context.
When a task is discarded, in case of a finite buffer size,
we call it an \emph{overflow} event, and we denote by $L^{\sss d(N)}(t)$ the
total number of overflow events under the JSQ($d(N)$) scheme up to time~$t$.

A sequence of random variables $\big\{X_N\big\}_{N\geq 1}$, for some function $f:\R\to\R_+$, is said to be $\Op(f(N))$, if the sequence of scaled random variables $\big\{X_N/f(N)\big\}_{N\geq 1}$ is  tight, or said to be $\op(f(N))$, if $\big\{X_N/f(N)\big\}_{N\geq 1}$ converges to zero in probability.
Boldfaced letters are used to denote vectors.
We denote by $\ell_1$ the space of all summable sequences. 
For any set $K$, the closure is denoted by $\overline{K}$. 
We denote by $D_E[0,\infty)$ the set of all \emph{c\`adl\`ag} (right continuous left limits exist) functions from $[0,\infty)$ to a complete separable metric space $E$, and
by `$\dto$'  convergence in distribution for real-valued random variables and with respect to the Skorohod $J_1$ topology for c\'adl\'ag processes.

\subsection{Fluid-limit results}\label{sec:fluidresult}

In the fluid-level analysis, we consider the subcritical regime
where $\lambda(N) / N \to \lambda < 1$ as $N \to \infty$.
In order to state the results, we first introduce some useful notation.
Denote the fluid-scaled system occupancy state by
$\qq^{\sss d(N)}(t) := \QQ^{\sss d(N)}(t) / N$, i.e., $q^{\sss d(N)}_i(t)=Q^{\sss d(N)}_i(t)/N$, and define 
$$\mathcal{S} =
\left\{\qq \in [0, 1]^b: q_i \leq q_{i-1} \mbox{ for all } i = 2, \dots, b, \mbox{ and } \sum_{i=1}^b q_i < \infty\right\}$$
as the set of all possible fluid-scaled occupancy states equipped with the $\ell_1$ topology.
For any $\qq \in \mathcal{S}$, denote $m(\qq) = \min\{i: q_{i + 1} < 1\}$,
with the convention that
$q_{b+1} = 0$ if $b < \infty$.
Note that $m(\qq)<\infty$, since $\qq\in\ell_1$.
If $m(\qq)=0$, then define $p_{0}(\qq)=1$ and $p_i(\qq) = 0$ for all $i \geq 1$.
If $m(\qq)>0$, distinguish two cases, depending on whether the normalized arrival
rate $\lambda$ is larger than $1 - q_{m(\qq) + 1}$ or not.
If $\lambda < 1 - q_{m(\qq) + 1}$, then define $p_{m(\qq) - 1}(\qq) = 1$
and $p_i(\qq) = 0$ for all $i \neq m(\qq) - 1$.
On the other hand, if $\lambda > 1 - q_{m(\qq) + 1}$,
then $p_{m(\qq) - 1}(\qq) = (1 - q_{m(\qq) + 1}) / \lambda$,
$p_{m(\qq)}(\qq) = 1 - p_{m(\qq) - 1}(\qq)$,
and $p_i(\qq) = 0$ for all $i \neq m(\qq) - 1, m(\qq)$.
Note that the assumption $\lambda < 1$ ensures that the latter case
cannot occur when $m(\qq) = b<\infty$.

\begin{theorem}{\normalfont (Universality of fluid limit for JSQ($d(N)$) scheme)}
\label{fluidjsqd-ssy}
Assume that the initial occupancy state $\qq^{\sss d(N)}(0)$ converges to $\qq^\infty$ in~$\mathcal{S}$ and $\lambda(N)/N\to\lambda<1$ as $N \to \infty$.
For the JSQ$(d(N))$ scheme with $d(N) \to\infty$, any subsequence of the sequence of processes
$\big\{\qq^{\sss d(N)}(t)\big\}_{t \geq 0}$ has a further subsequence that converges weakly with respect to the Skorohod $J_1$ topology, to the limit $\big\{\qq(t)\big\}_{t \geq 0}$ 
satisfying the following system of integral equations
\begin{equation}\label{eq:fluid}
 q_i(t) = q_i^\infty+ \lambda \int_0^t p_{i-1}(\qq(s))\dif s - \int_0^t(q_i(s) - q_{i+1}(s))\dif s,\quad i=1,\ldots,b,
\end{equation}
where the coefficients $p_i(\cdot)$ are as
defined above.
\end{theorem}

The above theorem shows that the fluid-level dynamics do not depend
on the specific growth rate of $d(N)$ as long as $d(N) \to \infty$
as $N \to \infty$.
In particular, the JSQ($d(N)$) scheme with $d(N)\to\infty$ exhibits the
same behavior as the ordinary JSQ policy  in the limit, and thus achieves fluid-level
optimality.

The coefficient $p_i(\qq)$ represents the instantaneous fraction
of incoming tasks assigned to servers with a queue length of exactly~$i$
in the fluid-level state $\qq \in \mathcal{S}$.
Assuming $m(\qq) < b$, a strictly positive fraction $1 - q_{m(\qq) + 1}$
of the servers have a queue length of exactly $m(\qq)$.
Since $d(N) \to \infty$, the fraction of incoming tasks that get assigned
to servers with a queue length of $m(\qq) + 1$ or larger is zero:
$p_i(\qq) = 0$ for all $i = m(\qq) + 1, \dots, b - 1$.
Also, tasks at servers with a queue length of exactly~$i$ are completed
at (normalized) rate $q_i - q_{i + 1}$, which is zero for all
$i = 0, \dots, m(\qq) - 1$, and hence the fraction of incoming tasks
that get assigned to servers with a queue length of $m(\qq) - 2$ or less
is zero as well: $p_i(\qq) = 0$ for all $i = 0, \dots, m(\qq) - 2$.
This only leaves the fractions $p_{m(\qq) - 1}(\qq)$
and $p_{m(\qq)}(\qq)$ to be determined.
Now observe that the fraction of servers with a queue length of exactly
$m(\qq) - 1$ is zero.
If $m(\qq)=0$, then clearly the incoming tasks will join the empty queue, 
and thus, $p_{m(\qq)}=1$, and $p_i(\qq) = 0$ for all $i \neq m(\qq)$.
Furthermore, if $m(\qq)\geq 1$, since tasks at servers with a queue length of exactly $m(\qq)$
are completed at (normalized) rate $1 - q_{m(\qq) + 1} > 0$,
incoming tasks can be assigned to servers with a queue length of exactly
$m(\qq) - 1$ at that rate.
We thus need to distinguish between two cases, depending on whether the
normalized arrival rate $\lambda$ is larger than $1 - q_{m(\qq) + 1}$ or not.
If $\lambda < 1 - q_{m(\qq) + 1}$, then all the incoming tasks can be
assigned to a server with a queue length of exactly $m(\qq) - 1$,
so that $p_{m(\qq) - 1}(\qq) = 1$ and $p_{m(\qq)}(\qq) = 0$.
On the other hand, if $\lambda > 1 - q_{m(\qq) + 1}$, then not all
incoming tasks can be assigned to servers with a queue length of
exactly $m(\qq) - 1$ active tasks, and a positive fraction will be
assigned to servers with a queue length of exactly $m(\qq)$:
$p_{m(\qq) - 1}(\qq) = (1 - q_{m(\qq) + 1}) / \lambda$
and $p_{m(\qq)}(\qq) = 1 - p_{m(\qq) - 1}(\qq)$.

It is easily verified that the unique fixed point $\qq^\star = (q_1^\star,q_2^\star,\ldots, q_b^\star)$ of the system of differential
equations in~\eqref{eq:fluid} is given by
\begin{equation}
\label{eq:fixed point-ssy}
q_i^* = \left\{\begin{array}{ll} \lambda, & i = 1, \\
0, & i =  2, \dots, b. \end{array} \right.
\end{equation}
Note that the fixed point in~\eqref{eq:fixed point-ssy} is consistent with the results in \cite{Mitzenmacher01,VDK96,YSK15}
for fixed~$d$, where taking $d \to \infty$ yields the same fixed point.
However, the results in \cite{Mitzenmacher01,VDK96,YSK15} for fixed~$d$
cannot be directly used to handle joint scalings, and do not yield the
universality of the entire fluid-scaled sample path
for arbitrary initial states as established in Theorem~\ref{fluidjsqd-ssy}.

The fixed point in~\eqref{eq:fixed point-ssy} in conjunction with the interchange of limits result in Proposition~\ref{th:interchange} below indicates that in stationarity the fraction of servers with a queue
length of two or larger is negligible.
Let 
$$\pi^{\sss d(N)}(\cdot)=\lim_{t\to\infty}\Pro{\qq^{\sss d(N)}(t)=\cdot}$$ 
be the stationary measure of the occupancy states of the $N^{\mathrm{th}}$ system under the JSQ($d(N)$) scheme.  
\begin{proposition}{{\normalfont (Interchange of limits)}}
\label{th:interchange}
For the JSQ$(d(N))$ scheme 
let $\pi^{\sss d(N)}$ be the stationary measure of the occupancy states of the $N^{\mathrm{th}}$ system. 
Then $\pi^{\sss d(N)}\dto\pi^\star$ as $N\to\infty$ with $d(N)\to\infty$, where $\pi^\star=\delta_{\qq^\star}$ with $\delta_x$ being the Dirac measure concentrated upon $x$, and $\qq^\star$ as in~\eqref{eq:fixed point-ssy}.
\end{proposition}
The above proposition relies on tightness of $\big\{\pi^{d(N)}\big\}_{N\geq 1}$ and the global stability of the fixed point, and is proved in Subsection~\ref{ssec:globstab}.
\\

We now consider an extension of the model in which tasks arrive in batches. We assume that the batches arrive as a Poisson process with rate  $\lambda(N)/\ell(N)$, and have fixed size $\ell(N)>0$, so that the effective total task arrival rate remains $\lambda(N)$. 
We will show that even for arbitrarily slowly growing batch size,  fluid-level optimality can be achieved with $O(1)$ communication overhead per task. 
For that, we define the JSQ($d(N)$) scheme adapted for batch arrivals. When a batch of size $\ell(N)$ arrives, the dispatcher samples $d(N)\geq \ell(N)$ servers without replacement, and assigns the $\ell(N)$ tasks to the $\ell(N)$ servers with the smallest queue length among the sampled servers.

\begin{theorem}{\normalfont (Batch arrivals)}
\label{th:batch-ssy}
Consider the batch arrival scenario with growing batch size $\ell(N)\to\infty$ and $\lambda(N)/N\to\lambda<1$ as $N\to\infty$. For the JSQ$(d(N))$ scheme with $d(N)\geq \ell(N)/(1-\lambda-\varepsilon)$ for any fixed $\varepsilon>0$, if $q^{\sss d(N)}_1(0)\pto q_1^\infty\leq \lambda$, and $q_i^{\sss d(N)}(0)\pto 0$ for all $i\geq 2$, then the sequence of processes
$\big\{\qq^{\sss d(N)}(t)\big\}_{t\geq 0}$ converges weakly to the limit $\big\{\qq(t)\big\}_{t\geq 0}$, described as follows: 
\begin{equation}\label{eq:batch}
q_1(t) = \lambda + (q_1^\infty-\lambda)\e^{-t},\quad
q_i(t)\equiv 0\quad \mathrm{for\ all}\quad i= 2,\ldots,b.
\end{equation}
\end{theorem}
The fluid limit in~\eqref{eq:batch} agrees with the fluid limit of the JSQ$(d(N))$ scheme if the initial state is taken as in Theorem~\ref{th:batch-ssy}. 
Further observe that the fixed point also coincides with that of the JSQ policy, as given by \eqref{eq:fixed point-ssy}.
Also, for a fixed $\varepsilon>0$, the communication overhead per task is on average given by $(1-\lambda-\varepsilon)^{-1}$ which is $O(1)$.
Thus Theorem~\ref{th:batch-ssy} ensures that in case of batch arrivals with growing batch size, fluid-level optimality can be achieved with $O(1)$ communication overhead per task.
The result for the fluid-level optimality in stationarity can also be obtained indirectly by exploiting the fluid-limit result in~\cite{YSK15}.
Specifically, it can be deduced from the result in~\cite{YSK15} that for batch arrivals with growing batch size, the JSQ$(d(N))$ scheme with suitably growing $d(N)$ yields the same fixed point of the fluid limit as described in~\eqref{eq:fixed point-ssy}.

\subsection{Diffusion-limit results}

In the diffusion-limit analysis, we consider the Halfin-Whitt
regime where 
$$\frac{N - \lambda(N)}{ \sqrt{N}} \to \beta\quad\text{as}\quad N \to \infty$$
for some positive coefficient $\beta > 0$.
In order to state the results, we first introduce some useful notation.
Let $\bar{\QQ}^{\sss d(N)}(t) =
\big(\bar{Q}_1^{\sss d(N)}(t), \bar{Q}_2^{\sss d(N)}(t), \dots, \bar{Q}_b^{\sss d(N)}(t)\big)$
be a properly centered and scaled version of the system occupancy state
$\QQ^{\sss d(N)}(t)$, with 
$$\bar{Q}_1^{\sss d(N)}(t) = - \frac{N-Q_1^{\sss d(N)}(t)}{ \sqrt{N}},\qquad \bar{Q}_i^{\sss d(N)}(t) =\frac{ Q_i^{\sss d(N)}(t)}{\sqrt{N}},\quad i = 2, \dots, b.$$
The reason why $Q_1^{\sss d(N)}(t)$ is centered around~$N$ while $Q_i^{\sss d(N)}(t)$,
$i = 2, \dots, b$, are not, is because the fraction of servers with a queue
length of exactly one tends to one, whereas the fraction of servers with a queue length of two or more tends to zero as $N\to\infty$.

\begin{theorem}{\normalfont (Universality of diffusion limit for JSQ($d(N)$) scheme)}
\label{diffusionjsqd-ssy}
Assume that the initial occupancy state $\bQ_i^{\sss d(N)}(0)$ converges to $\bQ_i(0)$ in $\mathbbm{R}$ as $N \to \infty$, the
buffer capacity $b\geq 2$ (possibly infinite),
and there exists some $k\geq 2$ such that $\bQ_{k+1}^N(0)=0$ for all sufficiently large $N$.
For $d(N) /( \sqrt{N} \log N)\to\infty$,
the sequence of processes $\big\{\bar{\QQ}^{\sss d(N)}(t)\big\}_{t \geq 0}$ 
converges weakly to the limit
 $\big\{\bar{\QQ}(t)\big\}_{t \geq 0}$ in~$D_{\mathcal{S}}[0,\infty)$,
where $\bar{Q}_i(t) \equiv 0$ for $i \geq k+1$
and $(\bar{Q}_1(t), \bar{Q}_2(t),\ldots,\bQ_k(t))$
are the unique solutions in $D_{\R^k}[0,\infty)$ of the stochastic integral equations
\begin{equation}\label{eq:diffusionjsqd-ssy}
\begin{split}
\bar{Q}_1(t) &= \bar{Q}_1(0) + \sqrt{2} W(t) - \beta t +
\int_0^t (- \bar{Q}_1(s) + \bar{Q}_2(s)) \dif s - U_1(t), \\
\bar{Q}_2(t) &= \bar{Q}_2(0) + U_1(t) - \int_0^t (\bar{Q}_2(s)-\bar{Q}_3(s) )\dif s, \\
\bQ_i(t)&=\bQ_i(0)-\int_0^t (\bQ_i(s) - \bQ_{i+1}(s))\dif s,\quad i= 3,\ldots,k,
\end{split}
\end{equation}
for $t \geq 0$, where $W$ is the standard Brownian motion and $U_1$ is
the unique nondecreasing nonnegative process in~$D_\R[0,\infty)$ satisfying
$\int_0^\infty \mathbbm{1}_{[\bar{Q}_1(t) < 0]} \dif U_1(t) = 0$.

\end{theorem}
Although~\eqref{eq:diffusionjsqd-ssy} differs from the diffusion limit obtained for the fully pooled M/M/N queue in the Halfin-Whitt regime~\cite{HW81, LK11, LK12}, it shares similar favorable properties.
Observe that $-\bQ_1^{\sss d(N)}$ is the scaled number of vacant servers. 
Thus, Theorem~\ref{diffusionjsqd-ssy} shows that over any finite time horizon,
there will be $O_P(\sqrt{N})$ servers with queue length zero
and $O_P(\sqrt{N})$ servers with a queue length larger than two,
and hence all but $O_P(\sqrt{N})$ servers have a queue length of exactly one.
This diffusion limit is proved in~\cite{EG15} for the ordinary JSQ policy.
Our contribution is to construct a stochastic coupling
and establish that, somewhat remarkably, the diffusion limit is the
same for any JSQ($d(N)$) scheme, as long as $d(N) /( \sqrt{N} \log(N))\to\infty$.
In particular, the JSQ($d(N)$) scheme with $d(N) /( \sqrt{N} \log(N))\to\infty$ exhibits the
same behavior as the ordinary JSQ policy in the limit, and thus achieves
diffusion-level optimality.
This growth condition for $d(N)$ is not only sufficient, but also nearly necessary,
as indicated by the next theorem.

\begin{theorem}{\normalfont (Almost necessary condition)}
\label{th:diff necessary}
Assume $\bQ_i^{\sss d(N)}(0) \dto \bQ_i(0)$ in $\mathbbm{R}$ as $N \to \infty$. If $d(N)/(\sqrt{N}\log N)\to 0$, then the diffusion limit of the JSQ$(d(N))$ scheme differs from that of the JSQ policy.
\end{theorem}

Theorem~\ref{th:diff necessary}, in conjunction with Theorem~\ref{diffusionjsqd-ssy}, shows that $\sqrt{N}\log N$ is the minimal order of $d(N)$ for the JSQ$(d(N))$ scheme to achieve diffusion-level optimality.


\subsection{Proof strategy}\label{subsec:strategy}
The idea behind the proofs of the asymptotic results for the JSQ$(d(N))$ scheme in Theorems~\ref{fluidjsqd-ssy} and~\ref{diffusionjsqd-ssy} is to 
(i)~prove the fluid limit and exploit the existing diffusion limit result for the ordinary JSQ policy, and then (ii)~prove a universality result by establishing that the ordinary JSQ policy and the JSQ$(d(N))$ scheme coincide under some suitable conditions on $d(N)$. 
For the ordinary JSQ policy the fluid limit in the subcritical regime is established in~Subsection~\ref{ssec:fluidjsq}, and the diffusion limit in the Halfin-Whitt heavy-traffic regime in~\cite[Theorem 2]{EG15}.
A direct comparison between the JSQ$(d(N))$ scheme and the ordinary JSQ policy is not straightforward, which is why we introduce the CJSQ$(n(N))$ class of schemes as an intermediate scenario to establish the universality result.

Just like the JSQ$(d(N))$ scheme, the schemes in the class CJSQ$(n(N))$ may be thought of as ``sloppy'' versions of the JSQ policy, in the sense that tasks are not necessarily assigned to a server with the shortest queue length but to one of the $n(N)+1$
lowest ordered servers, as graphically illustrated in Figure~\ref{fig:sfigCJSQ}.
In particular, for $n(N)=0$, the class only includes the ordinary JSQ policy. 
Note that the JSQ$(d(N))$ scheme is guaranteed to identify the lowest ordered server, but only among a randomly sampled subset of $d(N)$ servers.
In contrast, a scheme in the CJSQ$(n(N))$ class only guarantees that
one of the $n(N)+1$ lowest ordered servers is selected, but 
across the entire pool of $N$ servers. 
We will show that for sufficiently small $n(N)$, any scheme from the class CJSQ$(n(N))$ is still `close' to the ordinary JSQ policy. 
We will further prove that for sufficiently large $d(N)$ relative to $n(N)$ we can construct a scheme
called JSQ$(n(N),d(N))$, belonging to the CJSQ$(n(N))$ class, which differs `negligibly' from the JSQ$(d(N))$ scheme. 
Therefore,  for a `suitable' choice of $d(N)$ the idea is to produce a `suitable' $n(N)$.
This proof strategy has been schematically represented in Figure~\ref{fig:sfigRelation}.

In the next section we construct a stochastic coupling called S-coupling, which will be the key vehicle in establishing the universality result mentioned above.
\begin{remark} \normalfont
Observe that sampling \emph{without} replacement polls more servers than \emph{with} replacement, and hence the minimum number of active tasks among the selected servers is stochastically smaller in the case without replacement.
As a result, for sufficient conditions as in Theorems~\ref{fluidjsqd-ssy} and~\ref{th:diff necessary}, it is enough to consider sampling with replacement.
Also, for notational convenience, in the proof of the almost necessary condition stated in Theorem~\ref{th:diff necessary} we will assume sampling with replacement, although the proof technique and the result is valid if the servers are chosen without replacement.
\end{remark} 

The remainder of this chapter is organized as follows.
In Section~\ref{sec:coupling} we construct a stochastic coupling called S-coupling, and establish the stochastic ordering relations which will be the key vehicle in establishing the universality result mentioned above.
Sections~\ref{sec:fluid} and~\ref{sec:diffusion} contain  the proofs of the fluid and diffusion limit results, respectively. 
Finally in Section~\ref{sec:conclusion} we
make some concluding remarks and briefly comment on topics for further research.

\section{Coupling and stochastic ordering}\label{sec:coupling}

In this section, we construct a path-wise coupling between any scheme from the class CJSQ($n(N)$) and the ordinary JSQ policy, which ensures that for sufficiently small $n(N)$, on any finite time interval, the two schemes differ negligibly. This plays an instrumental role in establishing  the universality results in Theorems~\ref{fluidjsqd-ssy} and~\ref{diffusionjsqd-ssy}.
All the statements in this section should be understood to apply to the $N^{\mathrm{th}}$ system with $N$ servers.

\subsection{Stack formation and deterministic ordering}
In order to prove the stochastic comparisons among the various schemes, as in~\cite{MBLW16-1}, we describe the many-server system as an ensemble of stacks, in a way that two different ensembles can be ordered.
In this formulation, at each step, items are added or removed according to some rule. From a high level,  we then show that if two systems follow some specific rules, then at any step, the two ensembles maintain some kind of deterministic ordering. 
This deterministic ordering turns into an almost sure ordering in the next subsection, when we construct the S-coupling.

Each server along with its queue is thought of as a stack of items, and we always consider the stacks to be arranged in nondecreasing order of their heights. 
The ensemble of stacks then represents the empirical CDF of the queue length distribution, and the $i^{\mathrm{th}}$ horizontal bar corresponds to $Q_i^{\Pi}$ (for some task assignment scheme $\Pi$), as depicted in Figure~\ref{figB}.
 If an arriving item happens to land on a stack which already contains $b$ items, then the item is discarded, and is added to a special stack $L^\Pi$ of discarded items, where it stays forever.

Any two ensembles $A$ and $B$, each having $N$ stacks and a maximum height $b$ per stack, are said to follow Rule($n_A,n_B,k$) at some step, 
if either an item is removed from the $k^{\mathrm{th}}$ stack in both ensembles (if nonempty), 
or an item is added to the $n_A^{\mathrm{th}}$ stack in ensemble $A$ and to the $n_B^{\mathrm{th}}$ stack in ensemble $B$. 

\begin{proposition}\label{prop:det ord}
For any two ensembles of stacks $A$ and $B$, as described above, if at any step \emph{Rule}$(n_A,n_B,k)$ is followed for some value of $n_A$, $n_B$, and $k$, with $n_A\leq n_B$, then the following ordering is always preserved: for all $m\leq b$,
\begin{equation}\label{eq:det ord-ssy}
\sum_{i=m}^b Q_i^A +L^A\leq \sum_{i=m}^b Q_i^B +L^B.
\end{equation}
\end{proposition} 

This proposition says that, while adding the items to the ordered stacks, if we ensure that in ensemble $A$ the item is always placed to the left of that in ensemble $B$, and if the items are removed from the same ordered stack in both ensembles, then the 
aggregate size of the $b-m+1$ highest horizontal bars as depicted in Figure~\ref{figB} plus the cumulative number of discarded items is no larger in $A$ than in $B$ throughout.

\begin{proof}[Proof of Proposition~\ref{prop:det ord}]
We prove the ordering by forward induction on the time-steps, i.e., we assume that at some step the ordering holds, and show that in the next step it will be preserved. 
In ensemble $\Pi$, where $\Pi=A$, $B$, after applying Rule($n_A,n_B,k$), the updated lengths of the horizontal bars are denoted by $\tQ^\Pi_i$, $i\geq 1$. Also, define $I_\Pi(c):=\max\big\{i\geq 0:Q_i^\Pi\geq N-c+1\big\}$, $c=1,\ldots,N$, with the convention that $Q_0^{\Pi}\equiv N$.

Now if the rule prescribes removal of an item from the $k^\mathrm{th}$ stack,  then the updated ensemble will have the values
\begin{equation}\label{eq:removal}
\tilde{Q}^{\Pi}_i=
\begin{cases}
Q^{\Pi}_i-1, &\mbox{ for }i=I_{\Pi}(k),\\
Q^{\Pi}_j,&\mbox{ otherwise, }
\end{cases}
\end{equation}
if $I_{\Pi}(k)\geq 1$; otherwise all the $Q^{\Pi}_i$-values remain unchanged. On the other hand, if the rule produces the addition of an item to stack $n_{\sss\Pi}$, then the values will be updated as
\begin{equation}\label{eq:addition}
\tilde{Q}^{\Pi}_i=
\begin{cases}
Q^{\Pi}_i+1, &\mbox{ for }i=I_{\Pi}(n_{\sss\Pi})+1,\\
Q^{\Pi}_j,&\mbox{ otherwise, }
\end{cases}
\end{equation}
if $I_{\Pi}(n_\Pi)<b$; otherwise all values remain unchanged.

Fix any $m\leq b$. Observe that in any event the $Q_i$-values change by at most one at any step, and hence it suffices to prove the preservation of the ordering in the case when \eqref{eq:det ord-ssy} holds with equality:
\begin{equation}\label{eq:equality}
\sum_{i=m}^b Q_i^A +L^A= \sum_{i=m}^b Q_i^B +L^B.
\end{equation}

We distinguish between two cases depending on whether an item is removed or added.
First suppose that the rule prescribes removal of an item from the $k-$th stack from both ensembles. Observe from \eqref{eq:removal} that the value of $\sum_{i=m}^b Q_i^\Pi +L^\Pi$ changes if and only if $I_\Pi(k)\geq m$. Also, since removal of an item can only decrease the sum, without loss of generality we may assume that $I_B(k)\geq m$, otherwise the right side of \eqref{eq:equality} remains unchanged, and the ordering is trivially preserved. From our initial hypothesis,
\begin{equation}
\sum_{i=m+1}^b Q_i^{A}+L^{A}\leq
\sum_{i=m+1}^b Q_i^{B}+L^{B}.
\end{equation}
This implies
\begin{equation}
\begin{split}
Q_{m}^{A}&=\sum_{i=m}^b Q_i^{A}-\sum_{i=m+1}^b Q_i^{A}\geq \sum_{i=m}^b Q_i^{B}-\sum_{i=m+1}^b Q_i^{B}= Q_m^{B}.
\end{split}
\end{equation}
Also, 
\begin{equation}
\begin{split}
I_B(k)\geq m &\iff Q_m^B\geq N-k+1 \\
&\implies Q_m^A\geq N-k+1 \iff I_A(k)\geq m.
\end{split}
\end{equation}
Therefore the sum $\sum_{i=m}^b Q_i^A +L^A$ also decreases, and the ordering is preserved.

Now suppose that the rule prescribes addition of an item to the respective stacks in both ensembles.
From \eqref{eq:addition} we get that after adding an item to the ensemble, the value of $\sum_{i=m}^b Q_i^\Pi +L^\Pi$ increases only if $I_\Pi(n_{\sss\Pi})\geq~m-1$. As in the previous case, we assume \eqref{eq:equality}, and since adding an item can only increase the concerned sums, we assume that $I_A(n_{A})\geq~m-1$, because otherwise the left side of~\eqref{eq:equality} remains unchanged, and the ordering is trivially preserved. Now from our initial hypothesis we have
\begin{equation}\label{eq:initial}
\sum_{i=m-1}^b Q_i^A +L^A\leq \sum_{i=m-1}^b Q_i^B +L^B.
\end{equation}
Combining \eqref{eq:equality} with \eqref{eq:initial} gives
\begin{equation}\label{eq:prop-stoch-assump3}
\begin{split}
Q_{m-1}^{A}&=\left(\sum_{i=m-1}^b Q_i^{A}+L^A\right)-\left(\sum_{i=m}^b Q_i^{A}+L^A\right)\\
&\leq \left(\sum_{i=m-1}^b Q_i^{B}+L^B\right)-\left(\sum_{i=m}^b Q_i^{B}+L^B\right)= Q_{m-1}^{B}.
\end{split}
\end{equation}
Observe that
\begin{equation*}
\begin{split}
I_A(n_{A})\geq m-1 &\iff Q_{m-1}^A\geq N-n_A+1\implies Q_{m-1}^A\geq N-n_B+1\\
&\implies Q_{m-1}^B\geq N-n_B+1\iff I_B(n_{B})\geq m-1.
\end{split}
\end{equation*}
Hence, $\sum_{i=m}^b Q_i^B +L^B$ also increases, and the ordering is preserved. 
\end{proof}

\subsection{Stochastic ordering}\label{ssec:stoch-ord}
We now use the deterministic ordering established in Proposition~\ref{prop:det ord} in conjunction with the S-coupling construction to prove a stochastic comparison between the JSQ$(d(N))$ scheme, a specific scheme from the class CJSQ$(n(N))$ and the ordinary JSQ policy.
As described earlier, the class CJSQ$(n(N))$ contains all schemes that assign incoming tasks by some rule to any of the $n(N)+1$ lowest ordered servers. 
Observe that when $n(N)=0$, the class contains only the ordinary  JSQ policy. Also, if $n^{(1)}(N)< n^{(2)}(N)$, then CJSQ$(n^{(1)}(N))\subset$ CJSQ$(n^{(2)}(N)).$ 
Let MJSQ$(n(N))$ be a particular scheme that always assigns incoming tasks to precisely the $(n(N)+1)^{\mathrm{th}}$ ordered server. 
Notice that this scheme is effectively the JSQ policy when the system always maintains $n(N)$ idle servers, or equivalently, uses only $N-n(N)$ servers, and MJSQ($n(N)$) $\in$ CJSQ($n(N)$). 
For brevity, we suppress $n(N)$ in the notation for the remainder of this subsection.

We call any two systems \emph{S-coupled}, if they have synchronized arrival clocks and departure clocks of the $k^{\mathrm{th}}$ longest queue, for $1\leq k\leq N$ (`S' in the name of the coupling stands for `Server').
Consider three S-coupled systems following respectively the JSQ policy, any scheme from the class CJSQ, and the MJSQ scheme. 
Recall that $Q^\Pi_i(t)$ is the number of servers with at least $i$ tasks at time $t$ and $L^\Pi(t)$ is the total number of lost tasks up to time $t$, for the schemes $\Pi=$ JSQ, CJSQ, MJSQ. 
The following proposition provides a stochastic ordering for any scheme in the class CJSQ with respect to the ordinary JSQ policy and the MJSQ scheme.
\begin{proposition}\label{prop:stoch-ord}
For any fixed $m\geq 1$,
\begin{enumerate}[{\normalfont(i)}] 
\item\label{item:jsq-cjsq-ssy} $
\left\{\sum_{i=m}^bQ_i^{\JSQ}(t)+L^{\JSQ}(t)\right\}_{t\geq 0}\leq_{\mathrm{st}}
\left\{\sum_{i=m}^bQ_i^{\CJSQ}(t)+L^{\CJSQ}(t)\right\}_{t\geq 0},
$
\item\label{item:cjsq-mjsq-ssy} $\left\{\sum_{i=m}^bQ_i^{\CJSQ}(t)+L^{\CJSQ}(t)\right\}_{t\geq 0}\leq_{\mathrm{st}}
\left\{\sum_{i=m}^bQ_i^{\MJSQ}(t)+L^{\MJSQ}(t)\right\}_{t\geq 0},$
\end{enumerate}
provided the inequalities hold at time $t=0$.
\end{proposition}
The above proposition has the following immediate corollary, which will be used to prove bounds on the fluid and the diffusion scale.
\begin{corollary}\label{cor:bound-ssy}
In the joint probability space constructed by the S-coupling of the three systems under respectively \emph{JSQ}, \emph{MJSQ}, and any scheme from the class \emph{CJSQ}, the following ordering is preserved almost surely throughout the sample path: for any fixed $m\geq 1$
\begin{enumerate}[{\normalfont(i)}]
\item $Q_m^{\CJSQ}(t)\geq \sum_{i=m}^{b}Q_i^{\JSQ}(t)-\sum_{i=m+1}^{b}Q_i^{\MJSQ}(t)+L^{\JSQ}(t)-L^{\MJSQ}(t)$,
\item $Q_m^{\CJSQ}(t)
\leq \sum_{i=m}^{b}Q_i^{\MJSQ}(t)-\sum_{i=m+1}^{b}Q_i^{\JSQ}(t)+L^{\MJSQ}(t)-L^{\JSQ}(t),$
\end{enumerate}
provided the inequalities hold at time $t=0$.
\end{corollary}

\begin{proof}[Proof of Proposition~\ref{prop:stoch-ord}]
We first S-couple the concerned systems. 
Let us say that an incoming task is assigned to the $n_\Pi^{\mathrm{th}}$ ordered server under scheme $\Pi$, $\Pi$= JSQ, CJSQ, MJSQ. Then observe that, under the S-coupling, almost surely, 
$n_\JSQ\leq n_\CJSQ\leq n_\MJSQ.$ Therefore, Proposition~\ref{prop:det ord} ensures that in the probability space constructed through the S-coupling, the ordering is preserved almost surely throughout the sample path. 
\end{proof}
\begin{remark}
\normalfont
Note that $\sum_{i=1}^b\min\big\{Q_i,k\big\}$ represents the aggregate size of the rightmost $k$ stacks, i.e., the $k$ longest queues.
Using this observation, the stochastic majorization property
of the JSQ policy as stated in \cite{towsley, Towsley95,Towsley1992}
can be shown following similar arguments as in the proof of Proposition~\ref{prop:stoch-ord}.
Conversely, the stochastic ordering between the JSQ policy and the MJSQ scheme presented in Proposition~\ref{prop:stoch-ord} can also be derived from the weak majorization arguments developed in \cite{towsley, Towsley95,Towsley1992}. But it is only through the stack arguments developed in the previous subsection that we could extend the results to compare any scheme from the class CJSQ with the scheme MJSQ as well,
as stated in Proposition~\ref{prop:stoch-ord}~\eqref{item:cjsq-mjsq-ssy}.
\end{remark}

To analyze the JSQ$(d(N))$ scheme, we need a further stochastic comparison argument.
Consider two S-coupled systems following schemes $\Pi_1$ and $\Pi_2$.
Fix a specific arrival epoch, and let the arriving task join the $n_{\Pi_i}^{\mathrm{th}}$ ordered server in the $i^{\mathrm{th}}$ system following scheme $\Pi_i$, $i=1,2$ (ties can be broken arbitrarily in both systems). 
We say that at a specific arrival epoch the two systems \emph{differ in decision} if $n_{\Pi_1}\neq n_{\Pi_2}$,
and denote by $\Delta_{\Pi_1,\Pi_2}(t)$ the cumulative number of times the two systems differ in decision up to time~$t$.

\begin{proposition}\label{prop:stoch-ord2}
For two S-coupled systems under schemes $\Pi_1$ and $\Pi_2$ the following inequality is preserved almost surely
\begin{equation}\label{eq:stoch-ord2}
\sum_{i=1}^b |Q_i^{\Pi_1}(t)-Q_i^{\Pi_2}(t)|\leq 2\Delta_{\Pi_1,\Pi_2}(t)\qquad\forall\ t\geq 0,
\end{equation}
provided the two systems start from the same occupancy state at $t=0$, i.e., $Q_i^{\Pi_1}(0)= Q_i^{\Pi_2}(0)$ for all $i=1,2,\ldots, b$.
\end{proposition}

\begin{proof}
We will again use forward induction on the event times of arrivals and departures. Let the inequality~\eqref{eq:stoch-ord2} hold at time epoch $t_0$, and let $t_1$ be the next event time.
We distinguish between two cases, depending on whether $t_1$ is an arrival epoch or a departure epoch. 

If $t_1$ is an arrival epoch and the systems differ in decision, then observe that the left side of \eqref{eq:stoch-ord2} can only increase by two. In this case, the right side also increases by two, and the inequality is preserved.
Therefore, it is enough to prove that the left side of \eqref{eq:stoch-ord2} remains unchanged if the two systems do not differ in decision. In that case, assume that both $\Pi_1$ and $\Pi_2$ assign the arriving task to the $k^{\mathrm{th}}$ ordered server.
Recall from the proof of Proposition~\ref{prop:det ord} the definition of $I_\Pi$ for some scheme $\Pi$. 
If $I_{\sss\Pi_1}(k)=I_{\sss\Pi_2}(k)$, then the left side of \eqref{eq:stoch-ord2} clearly remains unchanged. Now, without loss of generality, assume $I_{\sss\Pi_1}(k)<I_{\sss\Pi_2}(k)$. Therefore, 
$$Q_{\sss I_{\Pi_1}(k)+1}^{\Pi_1}(t_0)< Q_{\sss I_{\Pi_1}(k)+1}^{\Pi_2}(t_0).$$
After an arrival, the $ (I_{\Pi_1}(k)+1)$-th term in the left side of \eqref{eq:stoch-ord2} decreases by one, and the $ (I_{\Pi_2}(k)+1)$-th term may increase by at most one. Thus the inequality is preserved.

If $t_1$ is a departure epoch, then due to the S-coupling, without loss of generality, assume that a potential departure occurs from the $k^{\mathrm{th}}$ ordered server. Also note that a departure in either of the two systems can change at most one of the $Q_i$-values. 
If at time epoch $t_0$, $I_{\Pi_1}(k)=I_{\Pi_2}(k)=i$, then both $Q_{\sss i}^{\Pi_1}$ and $Q_{\sss i}^{\Pi_2}$ decrease by one, and hence the left side of \eqref{eq:stoch-ord2} does not change. 
Otherwise, without loss of generality assume $I_{\Pi_1}(k)<I_{\Pi_2}(k).$ Then observe that 
$$Q_{\sss I_{\Pi_2}(k)}^{\Pi_1}(t_0)<Q_{\sss I_{\Pi_2}(k)}^{\Pi_2}(t_0).$$ 
Furthermore, after the departure, $Q_{\sss I_{\Pi_1}(k)}^{\Pi_1}$ may decrease by at most  one. Therefore $|Q_{\sss I_{\Pi_1}(k)}^{\Pi_1}- Q_{\sss I_{\Pi_1}(k)}^{\Pi_2}|$ may increase by at most one, and $Q_{\sss I_{\Pi_2}(k)}^{\Pi_2}$ decreases by one, thus $|Q_{\sss I_{\Pi_2}(k)}^{\Pi_1}- Q_{I_{\sss \Pi_2}(k)}^{\Pi_2}|$ decreases by one. Hence, in total, the left side of \eqref{eq:stoch-ord2} either remains the same or decreases by one.
\end{proof}

\subsection{Comparing the JSQ(d) and CJSQ(n) schemes}
We will now introduce the JSQ$(n,d)$ scheme with $n,d\leq N$, which is an intermediate blend between 
the CJSQ$(n)$ schemes and the JSQ$(d)$ scheme.
The JSQ$(n,d)$ scheme will be seen in a moment to  be a scheme in the CJSQ$(n)$ class. 
It will also be seen to approximate the JSQ$(d)$ scheme closely.
We now specify the JSQ$(d,n)$ scheme. At its first step, just as in the JSQ$(d)$ scheme, it first chooses the shortest of $d$ random candidates but only sends this to that server's queue if it is one of the $n+1$ shortest queues. 
If it is not, then at the second step it picks any of the $n+1$ shortest queues uniformly at random and then sends to that server's queue. 
As was mentioned earlier, by construction, JSQ$(d,n)$ is a scheme in CJSQ$(n)$.

We now consider two S-coupled systems with a JSQ$(d)$ and a JSQ$(n,d)$ scheme.
Assume that at some specific arrival epoch, the incoming task is dispatched to the $k^{\mathrm{th}}$ ordered server in the system under the JSQ($d$) scheme. If $k\in\{1,2,\ldots,n+1\}$, then the system under JSQ$(n,d)$ scheme also assigns the arriving task to the $k^{\mathrm{th}}$ ordered server. 
Otherwise, it dispatches the arriving task uniformly at random among the first $(n+1)$ ordered servers.

In the next proposition we will bound the number of times these two systems differ in decision on any finite time interval.
For any $T\geq 0$, let $A(T)$ and $\Delta(T)$ be the total number of arrivals to the system and the cumulative number of times that the JSQ($d$) scheme and JSQ$(n,d)$ scheme differ in decision up to time $T$.
\begin{proposition}\label{prop:differ}
For any $T\geq 0$, and $M>0,$
\begin{equation}
\Pro{\Delta(T)\geq M\given A(T)}\leq \frac{A(T)}{M}\left(1-\frac{n}{N}\right)^{d}.
\end{equation}
\end{proposition}
\begin{proof}
Observe that at any arrival epoch, the systems under the JSQ$(d)$ scheme and the JSQ$(n,d)$ scheme will differ in decision 
only if none of the $n$ lowest ordered servers gets selected by the JSQ$(d)$ scheme.
Now, at any arrival epoch, the probability that the JSQ($d$) scheme does not select any of the $n$ lowest ordered servers, is given by 
$$p=\left(1-\frac{n}{N}\right)^{d}.$$
Since at each arrival epoch, $d$ servers are selected independently, given $A(T)$, 
$$\Delta(T)\sim \mbox{Bin}(A(T),p).$$
Therefore, for $T\geq 0$, Markov's inequality yields, for any fixed $M>0$,
$$\Pro{\Delta(T)\geq M\given A(T)}\leq \frac{\expect{\Delta(T)\given A(T)}}{M}=\frac{A(T)}{M}\left(1-\frac{n}{N}\right)^{d}.$$
\end{proof}

\section{Fluid-limit proofs}\label{sec:fluid}
In this section we prove the fluid-limit results for the JSQ$(d(N))$ scheme stated in Theorems~\ref{fluidjsqd-ssy} and~\ref{th:batch-ssy}.
The fluid limit for the ordinary JSQ policy is provided in Subsection~\ref{ssec:fluidjsq}, and in Subsection~\ref{ssec:equivjsq} we prove a universality result establishing that under the condition that $d(N)\to\infty$ as $N\to\infty$, the fluid limit for the JSQ$(d(N))$ scheme coincides with that for the ordinary JSQ policy.

\subsection{Fluid limit of JSQ}\label{ssec:fluidjsq}
We now prove Proposition~\ref{th:interchange} using the time scale separation technique developed in~\cite{HK94}, suitably extended to an infinite-dimensional space. 
As mentioned in the introduction, to the best of our knowledge, this is the first time the transient fluid limit of the ordinary JSQ policy is rigorously established.
We also observe that in order to exploit the coupling framework in Section~\ref{ssec:stoch-ord} and in particular Proposition~\ref{prop:stoch-ord}, we need convergence of tail-sums.
Thus we need to establish the fluid convergence result with respect to the $\ell_1$ topology, which makes the analysis technically challenging.

To leverage the time scale separation technique, note that the rate at which incoming tasks join a server
with $i$ active tasks is determined only by the process $\ZZ^N(\cdot)=(Z_1^N(\cdot),\ldots,Z_b^N(\cdot))$, where $Z_i^N(t)=N-Q_i^N(t)$, $i=1,\ldots,b$, represents the number of servers with fewer than $i$ tasks at time $t$.
Furthermore, the dynamics of the $\ZZ^N(\cdot)$ process can be described as
\begin{equation}\label{eq:prelimit-slowprocess}
\ZZ^N \rightarrow 
\begin{cases}
\ZZ^N+e_i& \quad\mbox{ at rate }\quad N(q_i-q_{i+1}),\\
\ZZ^N-e_i& \quad\mbox{ at rate }\quad N\lambda \ind{\ZZ_\qq\in\mathcal{R}_{i}},
\end{cases}
\end{equation}
where $e_i$ is the $i^{\mathrm{th}}$ unit vector, and 
\begin{equation}\label{eq:partition}
\mathcal{R}_i := \big\{(z_1,z_2,\ldots, z_b): z_1=\ldots=z_{i-1}=0<z_i\leq z_{i+1}\leq\ldots\leq z_b\big\}\in\mathcal{G},
\end{equation}
$i=1,2,\ldots,b$,
with the convention that $Q^N_{b+1}$ is always taken to be zero, if $b<\infty$.
Observe that  in any time interval $[t,t+\varepsilon]$ of length $\varepsilon>0$, the $\ZZ^N(\cdot)$ process experiences $O(\varepsilon N)$ events (arrivals and departures), 
while the $\qq^N(\cdot)$ process can change by only an $O(\varepsilon)$ amount.
In other words, loosely speaking, around a `small' neighborhood of time $t$, the $q_i(t)$'s are constants, while as $N\to \infty$, the process $\ZZ^N(\cdot)$ behaves as a {\em time-scaled} version of the following process:
\begin{equation}
\ZZ_{\qq(t)} \rightarrow 
\begin{cases}
\ZZ_{\qq(t)}+e_i& \quad\mbox{ at rate }\quad q_i(t)-q_{i+1}(t),\\
\ZZ_{\qq(t)}-e_i& \quad\mbox{ at rate }\quad \lambda \ind{\ZZ_{\qq(t)}\in\mathcal{R}_{i}}.
\end{cases}
\end{equation}
Therefore, the $\ZZ^N(\cdot)$ process evolves on a much faster time scale than
the $\qq^N(\cdot)$ process.
As a result, in the limit as $N\to\infty$, at each time point $t$, the $\ZZ^N(\cdot)$ process
achieves stationarity depending on the instantaneous value of the $\qq^N(\cdot)$ process, i.e., a separation of time scales takes place. 
In order to establish the time-scale separation and the fluid limit results, we first write the evolution of the occupancy states in terms of a suitable random measure (see~\eqref{eq:martingale rep assumption 2-2}) and establish in Proposition~\ref{prop:rel compactness} that the sequence of joint occupancy process and the random measure is relatively compact. 
We also characterize the limit of any convergent subsequence, where we invoke analogous arguments as used in the proofs of~\cite[Lemma 2]{HK94} and \cite[Theorem 3]{HK94} to complete the proof of the separation of time scales.
The proof of the fluid limit result is then completed by establishing uniqueness of the instantaneous stationary distribution achieved by the fast process, given any fluid-scaled occupancy state.

Denote by $\bZ_+$ the one-point compactification of the set of nonnegative integers~$\Z_+$, i.e., $\bZ_+=\Z_+\cup\{\infty\}$.
Equip $\bZ_+$ with the order topology. Denote $G=\bZ_+^b$ equipped with product topology, and with the Borel $\sigma$-algebra $\mathcal{G}$.
Let us consider the $G$-valued process $\ZZ^N(s):=\big(Z_i^N(s)\big)_{i\geq 1}$ as introduced above.
Note that for the ordinary JSQ policy, the probability that a task arriving at (say) $t_k$ is assigned to 
some server with $i$ active tasks is given by $p_{i-1}^N(\QQ^N(t_k-))=\ind{\ZZ^N(t_k-)\in\mathcal{R}_{i}}$, where 
$\mathcal{R}_i$ is as in~\eqref{eq:partition}.
We prove the following fluid-limit result for the ordinary JSQ policy.
Recall the definition of $m(\qq)$ in Subsection~\ref{sec:fluidresult}. If $m(\qq)>0$, then define
\begin{equation}\label{eq:fluid-gen}
p_{i}(\qq)=
\begin{cases}
\min\big\{(1-q_{m(\qq)+1})/\lambda,1\big\} & \quad\mbox{ for }\quad i=m(\qq)-1,\\
1 - p_{\sss m(\qq) - 1}(\qq) & \quad\mbox{ for }\quad i=m(\qq),\\
0&\quad \mbox{ otherwise,}
\end{cases}
\end{equation}
and else, define $p_0(\qq) = 1$ and $p_i(\qq) = 0$ for all $i = 1,\ldots,b$.
\begin{theorem}[{Fluid limit of JSQ}]
\label{th:genfluid1}
Assume $\qq^N(0)\pto \qq^\infty$ in $\mathcal{S}$ and $\lambda(N)/N\to\lambda>0$ as $N\to\infty$. Then any subsequence of the sequence of processes $\big\{\qq^N(t)\big\}_{t\geq 0}$ for the ordinary JSQ policy has a further subsequence that converges weakly with respect to the Skorohod $J_1$ topology to the limit $\{\qq(t)\}_{t\geq 0}$ satisfying the following system of integral equations
\begin{equation}\label{eq:fluidfinal}
q_i(t) = q_i(0)+\lambda \int_0^t p_{i-1}(\qq(s))\dif s - \int_0^t (q_i(s)-q_{i+1}(s))\dif s, \quad i=1,2,\ldots,b,
\end{equation}
where $\qq(0) = \qq^\infty$ and the coefficients $p_i(\cdot)$ are as
defined in~\eqref{eq:fluid-gen}.
\end{theorem}
 The rest of this section will be devoted to the proof of Theorem~\ref{th:genfluid1}.
First we construct the martingale representation of the occupancy state process $\QQ^N(\cdot)$.
Note that the component $Q_i^N(t)$, satisfies the identity relation
\begin{align}\label{eq:recursion}
Q_i^N(t)=Q_i^N(0)+A_i^N(t)-D_i^N(t),&\quad\mbox{ for }\quad i=1,\ldots, b,
\end{align}
where
\begin{align*}
A_i^N(t)&=\mbox{ number of arrivals during $[0,t]$ to some server with }i-1\mbox{ active tasks,} \\
D_i^N(t)&=\mbox{ number of departures during $[0,t]$ from some server with }i\mbox{ active tasks}.
\end{align*}
We can express $A^N_i(t)$ and $D_i^N(t)$ as
\begin{align*}
A^N_i(t) &=  \mathcal{N}_{A,i}\left(\lambda(N)\int_0^t p_{i-1}^N(\QQ^N(s))\dif s\right),\\
D_i^N(t) &=  \mathcal{N}_{D,i}\left(\int_0^t(Q^N_i(s)-Q^N_{i+1}(s))\dif s\right),
\end{align*}
where $\mathcal{N}_{A,i}$ and $\mathcal{N}_{D,i}$ are mutually independent unit-rate Poisson processes, $i=1,2,\ldots,b$.
Define the sigma fields
\begin{align*}
\mathcal{A}^N_i(t)&:= \sigma\left(A^N_i(s): 0\leq s\leq t\right),\\
\mathcal{D}_i^N(t)&:= \sigma\left(D_i^N(s): 0\leq s\leq t\right),\mbox{ for }i = 1,\ldots,b,
\end{align*}
and the filtration $\mathbf{F}^N\equiv\big\{\mathcal{F}^N_t:t\geq 0\big\}$ with
\begin{equation}\label{eq:filtration}
\mathcal{F}^N_t:=\bigvee_{i=1}^{\infty} [\mathcal{A}_i^N(t)\vee \mathcal{D}_i^N(t)]
\end{equation}
augmented by all the null sets. 
Now we have the following martingale decomposition from the random time change of a unit-rate Poisson process result in \cite[Lemma~3.2]{PTRW07}.

\begin{proposition}[Martingale decomposition]
\label{prop:mart-rep}
The following are $\mathbf{F}^N$-martingales, for $i\geq 1$:
\begin{equation}\label{eq:martingales}
\begin{split}
M_{A,i}^N(t)&:=  \mathcal{N}_{A,i}\left(\lambda(N)\int_0^t p_{i-1}^N(\QQ^N(s))\dif s\right)-\lambda(N)\int_0^t p_{i-1}^N(\mathbf{Q}^N(s)) \dif s,\\
M_{D,i}^N(t)&:=\mathcal{N}_{D,i}\left(\int_0^t (Q^N_i(s)-Q^N_{i+1}(s))\dif s\right)-\int_0^t (Q^N_i(s)-Q^N_{i+1}(s))\dif s,
\end{split}
\end{equation}
with respective compensator and predictable quadratic variation processes given by
\begin{align*}
\langle M_{A,i}^N\rangle(t)&:= \lambda(N)\int_0^t p_{i-1}^N(\mathbf{Q}^N(s-))\dif s,\\
\langle M_{D,i}^N\rangle(t)&:=\int_0^t (Q^N_i(s)-Q^N_{i+1}(s))\dif s.
\end{align*}
\end{proposition}
Therefore, finally we have the following martingale representation of the $N^{\mathrm{th}}$ process:
\begin{equation}\label{eq:mart-unscaled}
\begin{split}
Q_i^N(t)&=Q_i^N(0)+\lambda(N)\int_0^t p_{i-1}^N(\mathbf{Q}^N(s))\dif s
-\int_0^t (Q^N_i(s)-Q^N_{i+1}(s))\dif s \\
&\hspace{2cm}+(M_{A,i}^N(t)-M_{D,i}^N(t)),\quad t\geq 0,\quad i= 1,\ldots,b.
\end{split}
\end{equation}
In the proposition below, we prove that the martingale part vanishes in $\ell_1$ when scaled by~$N$. 

\begin{proposition}[Convergence of martingales]
\label{prop:mart zero1}
$$\left\{\frac{1}{N}\sum_{i\geq 1}(|M_{A,i}^N(t)|+|M_{D,i}^N(t)|)\right\}_{t\geq 0}\dto \big\{m(t)\big\}_{t\geq 0}\equiv 0.$$
\end{proposition}
\begin{proof}
The proof follows using the same line of arguments as in the proof of~\cite[Theorem 3.13]{Mitzenmacher96}, and hence is sketched only briefly for the sake of completeness.
Fix any $T\geq 0$, and observe that 
\begin{align}
&\lim_{N\to\infty}\sup_{t\in[0,T]}\frac{1}{N}\sum_{i\geq 1}|M_{A,i}^N(t)|\label{subeq:4-8} \\
&=
\lim_{N\to\infty}\sup_{t\in[0,T]}\frac{1}{N}\left(\sum_{i\geq 1} \left|\mathcal{N}_{A,i}\left(\lambda(N)\int_0^t p_{i-1}^N(\QQ^N(s))\dif s\right)\right.\right.\nonumber\\
&\hspace{6cm}\left.\left.-\lambda(N)\int_0^t p_{i-1}^N(\QQ^N(s))\dif s\right| \right)\nonumber\\
&\leq \lim_{N\to\infty}\frac{1}{N}\sum_{i\geq 1} \mathcal{N}_{A,i}\left(\lambda(N)\int_0^T p_{i-1}^N(\QQ^N(s))\dif s\right)+\lambda T\label{subeq:4-10}.
\end{align} 
Since $N^{-1}\sum_{i\geq 1} \lambda(N)\int_0^t p_{i-1}^N(\QQ^N(s))\dif s\to \lambda t<\infty,$ the $\lim_{N\to\infty}$ and $\sum_{i\geq 1}$ above can be interchanged in~\eqref{subeq:4-10}, and hence in~\eqref{subeq:4-8}.
Now for each $i\geq 1$, from Doob's inequality~\cite[Theorem 1.9.1.3]{LS89}, we have for any $\epsilon>0,$
\begin{align*}
&\Pro{\sup_{t\in[0,T]}\frac{1}{N}|M_{A,i}^N(t)|\geq \epsilon}=\Pro{\sup_{t\in[0,T]}|M_{A,i}^N(t)|\geq N\epsilon}\\
&\leq \frac{1}{N^2\epsilon^2}\expect{\langle M_{A,i}^N\rangle (T)}
\leq \frac{1}{N\epsilon^2}\int_0^T p_{i-1}(\mathbf{Q}^N(s-))\lambda(N)\dif s
\leq \frac{\lambda T}{N\epsilon^2}\to 0,
\end{align*}
as $N\to\infty$.
Thus $\sup_{t\in[0,T]}N^{-1}M_{A,i}^N(t)\pto 0$, and hence
$$\sup_{t\in[0,T]}N^{-1}\sum_{i\geq 1}|M_{A,i}^N(t)|\pto 0.$$
Using similar arguments as above, we can also show that 
$$\sup_{t\in[0,T]}N^{-1}\sum_{i\geq 1}|M_{D,i}^N(t)|\pto 0,$$ and the proof is complete.
\end{proof}

Now we prove the relative compactness of the sequence of fluid-scaled processes.
Recall that we denote all the fluid-scaled quantities by their respective small letters, e.g.~$\mathbf{q}^N(t):=\mathbf{Q}^N(t)/N$, componentwise, i.e., $q_i^N(t):=Q_i^N(t)/N$ for $i\geq 1$. Therefore the martingale representation in~\eqref{eq:mart-unscaled} can be written as
\begin{equation}\label{eq:mart1}
\begin{split}
q_i^N(t)&=q_i^N(0)+\frac{\lambda(N)}{N}\int_0^t p_{i-1}^N(\mathbf{Q}^N(s))\dif s
-\int_0^t (q^N_i(s)-q^N_{i+1}(s))\dif s\\
&\hspace{3cm} +\frac{1}{N}(M_{A,i}^N(t)-M_{D,i}^N(t)),\quad i=1,2,\ldots, b,
\end{split}
\end{equation}
or equivalently,
\begin{equation}\label{eq:martingale rep assumption 2}
\begin{split}
q_i^N(t)&=q_i^N(0)+\frac{\lambda(N)}{N}\int_0^t \ind{\ZZ^N(s)\in\mathcal{R}_{i}}\dif s
-\int_0^t (q^N_i(s)-q^N_{i+1}(s))\dif s \\
&\hspace{3cm}+\frac{1}{N}(M_{A,i}^N(t)-M_{D,i}^N(t)),\quad i=1,2,\ldots, b.
\end{split}
\end{equation}
Now, we consider the Markov process $(\qq^N,\ZZ^N)(\cdot)$ defined on $\mathcal{S}\times G$. 
Define a random measure $\alpha^N$ on the measurable space $([0,\infty)\times G, \mathcal{C}\otimes\mathcal{G})$, when $[0,\infty)$ is endowed with the Borel sigma algebra $\mathcal{C}$, by
\begin{equation}
\alpha^N(A_1\times A_2):=\int_{A_1} \ind{\ZZ^N(s)\in A_2}\dif s,
\end{equation}
for $A_1\in\mathcal{C}$ and $A_2\in\mathcal{G}$. 
Then the representation in \eqref{eq:martingale rep assumption 2} can be written in terms of the random measure as
\begin{equation}\label{eq:martingale rep assumption 2-2}
\begin{split}
q_i^N(t)&=q_i^N(0)+\lambda\alpha^N([0,t]\times\mathcal{R}_i)
-\int_0^t (q^N_i(s)-q^N_{i+1}(s))\dif s \\
&\hspace{3cm} +\frac{1}{N}(M_{A,i}^N(t)-M_{D,i}^N(t)),\quad i=1,2,\ldots, b.
\end{split}
\end{equation}
Let $\mathfrak{L}$ denote the space of all measures on $[0,\infty)\times G$ satisfying $\gamma([0,t],G)= t$, endowed with the topology corresponding to weak convergence of measures restricted to $[0,t]\times G$ for each $t$.
\begin{proposition}[Relative compactness]
\label{prop:rel compactness}
Assume $\mathbf{q}^N(0)\dto\qq^\infty\in \mathcal{S}$ as $N\to\infty$.
Then $\big\{(\mathbf{q}^N(\cdot),\alpha^N)\big\}_{N\geq 1}$ is a relatively compact sequence in $D_{\mathcal{S}}[0,\infty)\times\mathfrak{L}$ and the limit $(\mathbf{q}(\cdot),\alpha)$ of any convergent subsequence satisfies
\begin{equation}\label{eq:rel compact}
q_i(t)=q_i^\infty+\lambda \alpha([0,t]\times\mathcal{R}_i) -\int_0^t (q_i(s)-q_{i+1}(s))\dif s,\quad i=1,2,\ldots, b.
\end{equation}
\end{proposition}

To prove Proposition~\ref{prop:rel compactness}, we will verify the  relative compactness conditions given in~\cite{EK2009}. 
Let $(E,r)$ be a complete and separable metric space. For any $x\in D_E[0,\infty)$, $\delta>0$ and $T>0$, define
\begin{equation}\label{eq:mod-continuity}
w'(x,\delta,T)=\inf_{\{t_i\}}\max_i\sup_{s,t\in[t_{i-1},t_i)}r(x(s),x(t)),
\end{equation}
where $\{t_i\}$ ranges over all partitions of the form $0=t_0<t_1<\ldots<t_{n-1}<T\leq t_n$ with $\min_{1\leq i\leq n}(t_i-t_{i-1})>\delta$ and $n\geq 1$.
 Below we state the conditions for the sake of completeness.
\begin{theorem}\label{th:from EK}
\begin{normalfont}
\cite[Corollary~3.7.4]{EK2009}
\end{normalfont}
Let $(E,r)$ be complete and separable, and let $\big\{X_n\big\}_{n\geq 1}$ be a family of processes with sample paths in $D_E[0,\infty)$. Then $\big\{X_n\big\}_{n\geq 1}$ is relatively compact if and only if the following two conditions hold:
\begin{enumerate}[{\normalfont (a)}]
\item For every $\eta>0$ and rational $t\geq 0$, there exists a compact set $\Gamma_{\eta, t}\subset E$ such that $$\varliminf_{n\to\infty}\Pro{X_n(t)\in\Gamma_{\eta, t}}\geq 1-\eta.$$
\item For every $\eta>0$ and $T>0$, there exists $\delta>0$ such that
$$\varlimsup_{n\to\infty}\Pro{w'(X_n,\delta, T)\geq\eta}\leq\eta.$$
\end{enumerate}
\end{theorem}

In order to prove the relative compactness, we will need the next three lemmas: Lemma~\ref{lem:technicalball} characterizes the relatively compact subsets of $\mathcal{S}$, Lemma~\ref{lem:tightcond} provides a necessary and sufficient criterion for a sequence of $\ell_1$-valued random variables to be tight, and Lemma~\ref{lem:trivialbound} is needed to ensure that at all finite times $t$, the occupancy state process lies in some compact set (possibly depending upon $t$).
\begin{lemma}[Compact subsets of $\mathcal{S}$]
\label{lem:technicalball}
Assume the buffer $b=\infty.$
A set $K\subseteq \mathcal{S}$ is relatively compact in $\mathcal{S}$ with respect to the $\ell_1$ topology if and only if
\begin{equation}
\label{eq:compactcond}
\lim_{k\to\infty} \sup_{\xx\in K} \sum_{i=k}^\infty x_i = 0.
\end{equation}
\end{lemma} 
\begin{proof}
For the if part, fix any $K\subseteq \mathcal{S}$ satisfying~\eqref{eq:compactcond}.
We will show that 
%
 any sequence $\big\{\xx^n\big\}_{n\geq 1}$ in $K$ 
has a Cauchy subsequence.
Since the $\ell_1$ space is complete, this will then imply that $\big\{\xx^n\big\}_{n\geq 1}$ has a convergent subsequence with the limit in $\overline{K}$,
which will complete the proof.

To show the existence of a Cauchy sequence, fix any $\varepsilon>0$, and choose $k\geq 1$ (depending on $\varepsilon$) such that
\begin{equation}\label{eq:uppertailtech}
\sum_{i\geq k}|x_i^n|<\frac{\varepsilon}{4}\qquad\forall\ n\geq 1.
\end{equation}
Now observe that the set of first coordinates $\big\{x_1^n\big\}_{n\geq 1}$ is a sequence in $[0,1]$, and hence has a convergent subsequence. Along that subsequence, the set of the second coordinates has a further convergent subsequence. Proceeding this way, we can get a subsequence along which the first $k-1$ coordinates converge. Therefore, depending upon $\varepsilon$, an $N'\in\N$ can be chosen, such that 
\begin{equation}\label{eq:lowertailtech}
\sum_{i<k}|x_i^n-x_i^m|<\frac{\varepsilon}{2}\qquad\forall\ m,n \geq N'.
\end{equation}
Therefore,~\eqref{eq:uppertailtech} and~\eqref{eq:lowertailtech} yields for all $n\geq \max\big\{N,N'\big\}$,
\begin{align*}
\norm{\xx^n-\xx^m}
&=\sum_{i\geq 1}|x_i^n-x_i^m|
\leq \sum_{i< k}|x_i^n-x_i^m| +\sum_{i\geq k}|x_i^n-x_i^m|\\
&\leq \sum_{i< k}|x_i^n-x_i^m| +\sum_{i\geq k}x_i^n+\sum_{i\geq k}x_i^m
<\varepsilon
\end{align*}
along the above suitably constructed subsequence. 
Now that the limit point is in $\mathcal{S}$ follows from the completeness of $\ell_1$ space and the fact that $\mathcal{S}$ is a closed subset of $\ell_1$.
Indeed, since the $\ell_1$ topology is finer than the product topology, any set that is closed with respect to the product topology is closed with respect to the $\ell_1$ topology, and observe that $\mathcal{S}$ is closed with respect to the product topology.

For the only if part, let $K\subseteq \mathcal{S}$ be relatively compact, and on the contrary, assume that there exists an $\varepsilon>0,$ such that
\begin{equation}
\lim_{k\to\infty} \sup_{\xx\in K} \sum_{i=k}^\infty x_i \geq \varepsilon.
\end{equation}
Therefore, for each $k\geq 1$, there exists $\xx^{(k)}\in K$, such that $\sum_{i=k}^\infty x^{(k)}_i \geq \varepsilon/2$.
Consider any limit point $\xx^*$ of the sequence $\big\{\xx^{(k)}\big\}_{k\geq 1}$, and note that $\sum_{i=j}^\infty x^*_i \geq \varepsilon/2$ for all $j\geq 1.$ This contradicts that $\xx^*\in\ell_1$, and the proof is complete.
\end{proof}

\begin{lemma}[Criterion for $\ell_1$-tightness]
\label{lem:tightcond}
Let $\big\{\XX^N\big\}_{N\geq 1}$ be a sequence of random variables in $\mathcal{S}$. Then the following are equivalent: 
\begin{enumerate}[{\normalfont (i)}]
\item $\big\{\XX^N\big\}_{N\geq 1}$ is tight with respect to product topology, and
for all $\varepsilon>0,$
\begin{equation}\label{eq:smalltail}
\lim_{k\to\infty}\varlimsup_{N\to\infty}\mathbbm{P}\Big(\sum_{i\geq k}X_i^N>\varepsilon\Big) = 0.
\end{equation}
\item $\big\{\XX^N\big\}_{N\geq 1}$ is tight with respect to $\ell_1$ topology.
\end{enumerate}
\end{lemma}
\begin{proof}
To prove (i)$\implies$(ii),
for any $\varepsilon>0$, we will construct a relatively compact set compact set $K(\varepsilon)$ such that
$$\Pro{\XX^N\notin \overline{K(\varepsilon)}}<\varepsilon\quad\mbox{for all } N.$$
Observe from~\eqref{eq:smalltail} that for all $\varepsilon>0$, there exists 
an $r(\varepsilon)\geq 1$, such that
$$\varlimsup_{N\to\infty}\mathbbm{P}\Big(\sum_{i\geq r(\varepsilon)}X_i^N>\varepsilon\Big) <\varepsilon,$$
and with it an $N(\varepsilon)\geq 1$, such that 
$$\mathbbm{P}\Big(\sum_{i\geq k(\varepsilon)}X_i^N>\varepsilon\Big)<\varepsilon\quad\mbox{for all } N> N(\varepsilon).$$
Furthermore, since $\big\{\XX^1, \XX^2,\ldots, \XX^{N(\varepsilon)}\big\}$ is a finite set of $\ell_1$-valued random variables, there exists 
$k(\varepsilon)\geq r(\varepsilon)$, such that 
$$\mathbbm{P}\Big(\sum_{i\geq k(\varepsilon)}X_i^N>\varepsilon\Big)<\varepsilon\quad\mbox{for all } N.$$
Thus, there exists an increasing sequence $\big\{k(n)\big\}_{n\geq 1}$
such that 
$$\mathbbm{P}\Big(\sum_{i\geq k(n)}X_i^N>\frac{\varepsilon}{2^n}\Big)<\frac{\varepsilon}{2^n}\quad\mbox{for all } N.$$
Define the set $K(\varepsilon)$ as
$$K(\varepsilon):= \Big\{\xx\in \mathcal{S}:\sum_{i\geq k(n)}x_i\leq \frac{\varepsilon}{2^n}\quad\mbox{for all}\quad n\geq 1\Big\}.$$
Due to Lemma~\ref{lem:technicalball}, we know that $K(\varepsilon)$ is relatively compact in $\ell_1$.
Also,
\begin{align*}
\mathbbm{P}\big(\XX^N\notin \overline{K(\varepsilon)}\big)&= \mathbbm{P}\Big(\bigcup_{n\geq 1}\Big\{\sum_{i\geq k(n)}X_i^N>\frac{\varepsilon}{2^n}\Big\}\Big)\\
&\leq \sum_{n\geq 1}\mathbbm{P}\Big(\sum_{i\geq k(n)}X_i^N>\frac{\varepsilon}{2^n}\Big)<\varepsilon.
\end{align*}
To prove (ii)$\implies$(i), first observe the fact that a sequence of random variables is tight with respect to the $\ell_1$ topology implies that it must be tight with respect to the product topology. 
Now assume on the contrary to~\eqref{eq:smalltail}, that there exists $\varepsilon>0$, such that 
\begin{equation}\label{eq:heavytail}
\lim_{k\to\infty}\varlimsup_{N\to\infty}\mathbbm{P}\Big(\sum_{i\geq k}X_i^N>\varepsilon\Big)>\varepsilon.
\end{equation}
Since $\big\{\XX^N\big\}_{N\geq 1}$ is tight with respect to the $\ell_1$ topology, take any convergent subsequence $\big\{\XX^{N(n)}\big\}_{n\geq 1}$ with $\XX^*$ being a random variable following the limiting measure.
In that case, observe that~\eqref{eq:heavytail} implies $\mathbbm{P}\big(\sum_{i\geq k}X_i^*>\varepsilon/2\big)>\varepsilon$ for all $k\geq 1$, which leads to a contradiction since $\XX^*$ is an $\ell_1$-valued random variable.
\end{proof}

\begin{lemma}
\label{lem:trivialbound}
For any $\qq\in \mathcal{S}$, assume that $\qq^N(0)\dto \qq^\infty$, as $N\to\infty$.
Then for any $t\geq 0$, there exists $M(t, \qq^\infty)\geq 1$, such that under the JSQ policy, with probability tending to one as $N\to\infty$, no arriving task is assigned to a server with $M(t,\qq^\infty)-1$ active tasks up to time $t$.
\end{lemma}
\begin{proof}
Let $A^N(t)$ be the cumulative number of tasks arriving up to time $t$.
Since the arrival rate is $\lambda(N)$, and $\lambda(N)/N\to\lambda$, as $N\to\infty,$ for any $\varepsilon>0$, 
$$\Pro{A^N(t)\geq (\lambda t + \varepsilon)N}\to 0\qquad\mbox{as}\qquad N\to\infty.$$
Define $M(t,\qq^\infty):=\min\big\{k\geq 1: \sum_{i=1}^{k-1}(1-q_i^\infty)>\lambda t\big\},$
and choose 
$$\varepsilon = \sum_{i=1}^{M(t,\qq^\infty)-1}(1-q_i^\infty)-\lambda t >0.$$
Note that since $\qq^\infty\in \mathcal{S}\subset \ell_1$, $M(t,\qq^\infty)$ exists and is finite for all $t\geq 0.$
We now claim that the probability that in the interval $[0,t]$ a task is assigned to some server with $M(t,\qq^\infty)$ active tasks tends to 0, as $N\to\infty.$
Indeed, in order for a task to be assigned to some server with $M(t,\qq^\infty)-1$ active tasks, all the servers must have at least $M(t,\qq^\infty)-1$ active tasks. 
Now, the minimum number of tasks required for this, is given by 
$\sum_{i=1}^{M(t,\qq^\infty)-1}(N-Q_i^N(0))$.
Therefore, the proof is complete by observing that 
\begin{align*}
&\mathbbm{P}\Big(A^N(t)\geq \sum_{i=1}^{M(t,\qq^\infty)-1}(N-Q_i^N(0))\Big) \\
&\hspace{2cm}= \Pro{A^N(t)\geq \Big(\lambda t +\frac{\varepsilon}{2}\Big) N}\to 0, \quad \mbox{as}\quad N\to\infty.
\end{align*}
\end{proof}

\begin{proof}[Proof of Proposition~\ref{prop:rel compactness}]
The proof goes in two steps. We first prove the relative compactness, and then show that the limit satisfies~\eqref{eq:rel compact}.

Observe from \cite[Proposition 3.2.4]{EK2009} that, to prove the relative compactness of the sequence of processes $\big\{(\mathbf{q}^N(\cdot),\alpha^N)\big\}_{N\geq 1}$, it is enough to prove relative compactness of the individual components.
Note that from Prohorov's theorem~\cite[Theorem 3.2.2]{EK2009}, $\mathfrak{L}$ is compact since $G$ is compact. Now, relative compactness of $\big\{\alpha^N\big\}_{N\geq 1}$ follows from the compactness of $\mathfrak{L}$ under the topology of weak convergence of measures and Prohorov's theorem.
To claim the relative compactness of $\big\{\mathbf{q}^N(\cdot)\big\}_{N\geq 1}$, we will verify the conditions of Theorem~\ref{th:from EK}. 

Observe that in order to verify Theorem~\ref{th:from EK} (a), we need to show tightness of the sequence $\big\{\mathbf{q}^N(t)\big\}_{N\geq 1}$ for each fixed (rational) $t\geq 0$.
Fix any $t\geq 0.$
Due to Lemma~\ref{lem:trivialbound}, we know  
$$\lim_{N\to\infty}\Pro{q_i^N(t)\leq q_i^N(0),\quad \forall\ i\geq M(t,\qq^\infty)}= 1.$$
Also, $\qq^N(0)\dto \qq^\infty$ with respect to the $\ell_1$ topology. In particular, $\big\{\qq^N(0)\big\}_{N\geq 1}$ is tight in $\ell_1$.
Therefore, using (ii)$\implies$(i) in Lemma~\ref{lem:tightcond} we obtain, for any $\varepsilon>0$,
\begin{align*}
\lim_{k\to\infty}\varlimsup_{N\to\infty} \mathbbm{P}\Big(\sum_{i\geq k}q_i^N(t)>\varepsilon\Big)
&\leq \lim_{k\to\infty}\varlimsup_{N\to\infty} \mathbbm{P}\Big(\sum_{i\geq k}q_i^N(0)>\varepsilon\Big)=0.
\end{align*}
Also, since $\qq^N(t)\in \mathcal{S}\subseteq [0,1]^b$, which is compact with respect to the product topology,  $\big\{\qq^N(t)\big\}_{N\geq 1}$ is tight with respect to the product topology.
Hence using (i)$\implies$(ii) in Lemma~\ref{lem:tightcond} we conclude that the sequence $\big\{\mathbf{q}^N(t)\big\}_{N\geq 1}$ is tight in $\ell_1.$
For condition (b),  first note that 
for all $i = 1,\ldots, b$.
\begin{align*}
&|q_i^N(t_1)-q_i^N(t_2)|
\leq \lambda \alpha^N([t_1,t_2]\times\mathcal{R}_i)+\int_{t_1}^{t_2} (q^N_i(s)-q^N_{i+1}(s))\dif s \\
&\hspace{3cm}+\frac{1}{N}\big|M_{A,i}^N(t_1)-M_{D,i}^N(t_1)-M_{A,i}^N(t_2)+M_{D,i}^N(t_2)\big| + o(1).
\end{align*}
Thus, 
\begin{equation}\label{mart-norm-ub}
\begin{split}
&\norm{\qq^N(t_1)-\qq^N(t_2)}\\
&\leq \lambda \sum_{i=1}^b\alpha^N([t_1,t_2]\times\mathcal{R}_i)+\int_{t_1}^{t_2} \sum_{i=1}^b (q^N_i(s)-q^N_{i+1}(s))\dif s \\
&\hspace{1cm}+\frac{1}{N}\sum_{i=1}^b\big|M_{A,i}^N(t_1)-M_{D,i}^N(t_1)-M_{A,i}^N(t_2)+M_{D,i}^N(t_2)\big|+o(1)\\
&\leq \lambda (t_1-t_2)+\int_{t_1}^{t_2} q^N_1(s)\dif s +\frac{1}{N}\sum_{i=1}^b\big|M_{A,i}^N(t_1)-M_{D,i}^N(t_1)\\
&\hspace{5.5cm}-M_{A,i}^N(t_2)+M_{D,i}^N(t_2)\big|+o(1)\\
&\leq (\lambda+1) (t_1-t_2)+\frac{1}{N}\sum_{i=1}^b\big|M_{A,i}^N(t_1)-M_{D,i}^N(t_1)\\
&\hspace{5.5cm}-M_{A,i}^N(t_2)+M_{D,i}^N(t_2)\big|+o(1).
\end{split}
\end{equation}
From the $\ell_1$ convergence of scaled martingales in Proposition~\ref{prop:mart zero1}, we get, for any $T\geq 0$,
$$\sup_{t\in[0,T]}\frac{1}{N}\sum_{i=1}^b|M_{A,i}^N(t_1)-M_{D,i}^N(t_1)-M_{A,i}^N(t_2)+M_{D,i}^N(t_2)|\pto 0.$$
Observe that the proof of the relative compactness of $\big\{\qq^N(t)\big\}_{t\geq 0}$ is complete if we show that for any $\eta>0$, there exists a $\delta >0$ and a finite partition $(t_j)_{i=1}^n$ of $[0,T]$ with $\min_j|t_{j}-t_{j-1}|>\delta$ such that 
\begin{equation}
\varlimsup_{N\to\infty}\mathbbm{P}\Big(\max_j \sup_{s,t\in [t_{j-1},t_j)}\norm{\qq^N(s)-\qq(t)} \geq \eta \Big) < \eta.
\end{equation}
Now, \eqref{mart-norm-ub} implies that, for any finite partition $(t_j)_{j= 1}^n$ of $[0,T]$,
\begin{align*}
\max_j \sup_{s,t\in [t_{j-1},t_j)} \norm{\qq^N(s)-\qq^N(t)} &\leq (\lambda+1) \max_j (t_{j}-t_{j-1})+\zeta_N,
\end{align*}
where $\Pro{\zeta_N>\eta/2}<\eta$ for all sufficiently large $N$. Now take $\delta = \eta/(4(\lambda+1))$ and any partition with $\max_j(t_j-t_{j-1})< \eta/(2(\lambda+1))$ and $\min_j(t_j-t_{j-1})>\delta$. 
On the event $\big\{\zeta_N\leq \eta/2\big\}$,  
$$\max_i \sup_{s,t\in [t_{i-1},t_i)}\norm{\qq^N(s)-\qq^N(t)} \leq \eta.$$
Therefore, for all sufficiently large $N$,
\begin{align*}
&\mathbbm{P}\Big(\max_j \sup_{s,t\in [t_{j-1},t_j)}\norm{\qq^N(s)-\qq^N(t)} \geq \eta \Big)
\leq \Pro{\zeta_N>\eta/2}\leq \eta.
\end{align*}

To prove that the limit $(\qq(\cdot),\alpha)$ of any convergent subsequence satisfies~\eqref{eq:rel compact}, we will use the continuous-mapping theorem~\cite[Theorem~3.4.1]{W02}.
Specifically, we will show that the right side of~\eqref{eq:martingale rep assumption 2-2} is a continuous map of suitable arguments.
Let $\big\{\qq(t)\big\}_{t\geq 0}$ and $\big\{\yy(t)\big\}_{t\geq 0}$ be an $\mathcal{S}$-valued and an $\ell_1$-valued c\`adl\`ag function, respectively. 
Also, let $\alpha$ be a measure on the measurable space $([0,\infty)\times G, \mathcal{C}\otimes\mathcal{G})$. Then for $\qq^0\in \mathcal{S}$, define for $i\geq 1$,
$$F_i(\qq,\alpha,\qq^0,\yy)(t):=q_i^0+y_i(t)+\lambda \alpha([0,t]\times\mathcal{R}_i)-\int_0^t(q_i(s) - q_{i+1}(s))\dif s.$$
Observe that it is enough to show that $\FF=(F_1,\ldots,F_b)$ is a continuous
operator. 
Indeed, in that case the right side of~\eqref{eq:martingale rep assumption 2-2} can be written as $\FF(\qq^N,\alpha^N,\qq^N(0),\yy^N)$, where $\yy^N=(y_1^N,\ldots,y_b^N)$ with $y_i^N= (M_{A,i}^N-M_{D,i}^N)/N$, and since each argument converges, we will get the convergence to the right side of~\eqref{eq:rel compact}.
Therefore, we now prove the continuity of $\FF$ below. 
In particular assume that (a)~the sequence of processes $\big\{(\qq^N,\yy^N)\big\}_{N\geq 1}$ converges to $(\qq,\yy)$ with respect to the $\ell_1$ topology, (b)~for any fixed $t\geq 0$, the sequence $\big\{\big(\alpha^N([0,t],\mathcal{R}_i)\big)_{i\geq 1}\big\}_{N\geq 1}$ in $\ell_1$ converges to $\big(\alpha([0,t],\mathcal{R}_i)\big)_{i\geq 1}$, and (c)~the sequence of $\mathcal{S}$-valued random  variables $\qq^N(0)$
 converges to $\qq(0)$ with respect tothe  $\ell_1$ topology.
 
 Fix any $T\geq 0$ and $\varepsilon>0$.
 \begin{enumerate}[{\normalfont (i)}]
 \item  Due to (a) above, choose $N_1\in\N$, such that for all $N\geq N_1$ 
 $$\sup_{t\in[0,T]}\norm{\qq^N(t)-\qq(t)}<\varepsilon/(4T).$$ 
 In that case, observe that
 \begin{align*}
 \sup_{t\in [0,T]}\int_0^t|q_1^N(t) - q_1(t)|\dif s& \leq T\sup_{t\in [0,T]}\norm{\qq^N(t))-\qq(t))} <\frac{\varepsilon}{4}.
 \end{align*}
 \item Again, due to (a), choose $N_2\in\N$, such that for all $N\geq N_2$ 
 $$\sup_{t\in[0,T]}\norm{\yy^N(t)-\yy(t)}<\varepsilon/4.$$
 \item   We now claim that for the $\epsilon > 0$ given above there is an $N_3 \in {\mathbb N}$
such that for all $N\geq N_3$
 \begin{equation}\label{eq:measureconv}
 \lambda\sum_{i\geq 1}  \left|\alpha^N([0,T]\times\mathcal{R}_i)-\alpha([0,T]\times\mathcal{R}_i)\right|<\frac{\varepsilon}{4}.
 \end{equation}
 Observe that we only know the weak convergence of the sequence of measures~$\alpha^N$, and therefore we cannot directly make assumption (b) above.
We are therefore about to  show that assumption (b) is valid in our case and that it follows from weak convergence.
Indeed, since $\qq^\infty\in \mathcal{S}\subseteq \ell_1$, there exists $\hat{M}(\qq^\infty)$, such that 
$q_{\hat{M}(\qq^\infty)}^\infty<1$, and consequently $q_i^\infty<1$ for all $i\geq \hat{M}(\qq^\infty)$.
Also, due to Lemma~\ref{lem:trivialbound}, 
$$\lim_{N\to\infty}\mathbbm{P}\Big(\sup_{t\in[0,T]}q_i^N(t)\leq q_i^N(0)\quad\mbox{for all}\quad i\geq M(T,\qq^\infty)\Big)= 1.$$
Thus, if $N_0 := \max\big\{\hat{M}(\qq^\infty), M(T,\qq^\infty)\big\}$, then 
$$\lim_{N\to\infty}\mathbbm{P}\Big(\sup_{t\in[0,T]}q_i^N(t)<1\quad\mbox{for all}\quad i\geq N_0\Big)= 1.$$
This implies
$$\sum_{i\geq N_0}\alpha^N([0,T]\times\mathcal{R}_i)\pto \sum_{i\geq N_0}\alpha([0,T]\times\mathcal{R}_i)=0.$$
Also, due to weak convergence of $\alpha^N$, 
$$\sum_{i< N_0}\alpha^N([0,T]\times\mathcal{R}_i)\pto \sum_{i< N_0}\alpha([0,T]\times\mathcal{R}_i).$$ 
 \item  Finally, due to (c), choose $N_4\in\N$, such that for all $N\geq N_4$ 
 $$\norm{\qq^N(0)-\qq(0)}<\varepsilon/4.$$
 \end{enumerate}
Let $\hat{N}=\max\big\{N_1,N_2,N_3,N_4\big\}$, then for $N\geq \hat{N}$,
\begin{align*}
&\sup_{t\in [0,T]}\norm{\FF(\qq^N,\alpha^N,\qq^N(0),\yy^N)-\FF(\qq,\alpha,\qq(0),\yy)}(t)<\varepsilon.
\end{align*}
Thus the proof of continuity of $\FF$ is complete.
\end{proof}

To characterize the limit in~\eqref{eq:rel compact}, for any $\qq\in \mathcal{S}$, define the Markov process  $\ZZ_{\qq}$ on $G$ as
\begin{equation}\label{eq:slowprocess}
\ZZ_{\qq} \rightarrow 
\begin{cases}
\ZZ_{\qq}+e_i& \quad\mbox{ at rate }\quad q_i-q_{i+1},\\
\ZZ_{\qq}-e_i& \quad\mbox{ at rate }\quad \lambda \ind{\ZZ_\qq\in\mathcal{R}_{i}},
\end{cases}
\end{equation}
where $e_i$ is the $i^{\mathrm{th}}$ unit vector, $i=1,\ldots,b$.

\begin{proof}[{Proof of Theorem~\ref{th:genfluid1}}]
Having proved the relative compactness in Proposition~\ref{prop:rel compactness},  it follows from analogous arguments as used in the proofs of~\cite[Lemma 2]{HK94} and \cite[Theorem 3]{HK94}, that the limit of any convergent subsequence of the sequence of processes $\big\{\qq^N(t)\big\}_{t\geq 0}$ satisfies
\begin{equation}
q_i(t) = q_i(0)+\lambda \int_0^t \pi_{\qq(s)}(\mathcal{R}_i)\dif s - \int_0^t (q_i(s)-q_{i+1}(s))\dif s, \quad i=1,2,\ldots,b,
\end{equation}
for \emph{some} stationary measure $\pi_{\qq(t)}$ of the Markov process  $\ZZ_{\qq(t)}$ described in~\eqref{eq:slowprocess} satisfying $\pi_{\qq}\big\{\ZZ: Z_i=\infty\big\}=1$ if $q_i<1$. 

Now it remains to show that $\qq(t)$ \emph{uniquely} determines the measure $\pi_{\qq(t)}$, and that $\pi_{\qq(s)}(\mathcal{R}_i)=p_{i-1}(\qq(s))$ described in~\eqref{eq:fluid-gen}. 
As mentioned earlier, in this proof we will now assume the specific assignment probabilities in~\eqref{eq:partition}, corresponding to the ordinary JSQ policy.
To see this, fix any $\qq=(q_1,\ldots,q_b)\in \mathcal{S}$.  
Observe that due to summability of the components of~$\qq$, there exists $0\leq m<\infty$, such that $q_{m+1}<1$ and $q_1=\ldots=q_m=1$,
with the convention that $q_0\equiv 1$ and $q_{b+1}\equiv 0$ if $b<\infty$. In that case,
$$\pi_{\qq}\big(\big\{Z_{m+1}=\infty, Z_{m+2}=\infty,\ldots,Z_b=\infty\big\}\big)=1.$$
Also, 
note that $q_i = 1$ forces $\dif q_i/\dif t \leq 0$, i.e., $\lambda \pi_{\qq}(\mathcal{R}_i) \leq q_i-q_{i+1}$ for all $i = 1, \ldots, m$, and in particular $\pi_{\qq}(\mathcal{R}_i) = 0$ for all $i = 1, \ldots, m - 1.$ Thus,
$$\pi_{\qq}\big(\big\{Z_1=0,Z_2=0,\ldots,Z_{m-1}=0\big\}\big)=1.$$

Therefore, $\pi_\qq$ is determined only by the stationary distribution of the $m^{\mathrm{th}}$ component, which can be described as a birth-death process
\begin{equation}\label{eq:bdprocess}
Z \rightarrow 
\begin{cases}
Z+1& \quad\mbox{ at rate }\quad q_m-q_{m+1},\\
Z-1& \quad\mbox{ at rate }\quad \lambda\ind{Z>0},
\end{cases}
\end{equation}
and let $\pi^{(m)}$ be its stationary distribution. 
Now it is enough to show that $\pi^{(m)}$ is uniquely determined by $\qq$. 
First observe that the process on $\bZ$ described in~\eqref{eq:bdprocess} is reducible, and can be decomposed into
two irreducible classes given by $\mathbbm{Z}$ and $\{\infty\}$, respectively.
Therefore, if $\pi^{(m)}(Z=\infty)=0$ or $1$, then it is unique. 
Indeed, if $\pi^{(m)}(Z=\infty)=0$, then $Z$ is a birth-death process on $\mathbbm{Z}$ only, and hence it has a unique stationary distribution. 
Otherwise, if $\pi^{(m)}(Z=\infty)=1$, then it is trivially unique. 
Now we distinguish between two cases depending on whether $q_m-q_{m+1}\geq \lambda$ or not.

Note that if $q_m-q_{m+1}\geq\lambda$, then $\pi^{(m)}(Z\geq k)=1$ for all $k\geq 0$. 
On $\bZ$ this shows that $\pi^{(m)}(Z=\infty)=1$.
Furthermore, if $q_m-q_{m+1}<\lambda$, we will show that $\pi^{(m)}(Z=\infty)=0$.
On the contrary, assume that $\pi^{(m)}(Z=\infty)=\varepsilon\in (0,1]$.
Also, let $\hat{\pi}^{(m)}$ be the unique stationary distribution of the birth-death process in~\eqref{eq:bdprocess} on $\mathbbm{Z}$.
Therefore, 
$$\pi_\qq(\mathcal{R}_m)=\pi^{(m)}(Z>0)=(1-\varepsilon)\hat{\pi}^{(m)}(Z>0)+\varepsilon = (1-\varepsilon)\frac{q_m-q_{m+1}}{\lambda}+\varepsilon.$$
Substituting into the differential form of the fluid equation~\eqref{eq:fluidfinal} at the given time $t$, we obtain that 
\begin{align*}
\frac{\dif q_m(t)}{\dif t} &= \lambda \Big[(1-\varepsilon)\frac{q_m-q_{m+1}}{\lambda}+\varepsilon \Big] - (q_m - q_{m+1})\\
 &= -\varepsilon(q_m-q_{m+1}) +\lambda \varepsilon >0,
\end{align*}
where the last inequality follows since we are considering the case when $q_m - q_{m+1}<\lambda$.
Now since $q_m(t)=1$, this leads to a contradiction for any $\varepsilon>0$, and hence it must be the case that $\pi^{(m)}(Z=\infty)=0$. 

Therefore, for all $\qq\in \mathcal{S}$, $\pi_\qq$ is uniquely determined by $\qq$. 
Furthermore, we can identify the expression for $\pi_q(\mathcal{R}_i)$ as
\begin{equation}
\pi_\qq(\mathcal{R}_i)=
\begin{cases}
\min\big\{(q_m-q_{m+1})/\lambda,1\big\}& \quad\mbox{ for }\quad i=m,\\
1- \min\big\{(q_m-q_{m+1})/\lambda,1\big\} & \quad\mbox{ for }\quad i=m+1,\\
0&\quad \mbox{ otherwise,}
\end{cases}
\end{equation}
and hence $\pi_{\qq(s)}(\mathcal{R}_i)=p_{i-1}(\qq(s))$ as claimed.
\end{proof}

\subsection{Equivalence on fluid scale}\label{ssec:equivjsq}
Having proved Theorem~\ref{th:genfluid1}, it suffices to prove the universality property stated in the next proposition. This will complete the proof of Theorem~\ref{fluidjsqd-ssy}.
\begin{proposition}\label{prop:samefluid}
If $d(N)\to\infty$ as $N\to\infty$, then the JSQ$(d(N))$ scheme and the ordinary JSQ policy have the same fluid limit.
\end{proposition}
The proof of the above proposition uses the S-coupling 
 results from Section~\ref{sec:coupling}, and consists of three steps:
\begin{enumerate}[{\normalfont (i)}]
\item First we show that if $n(N)/N\to 0$ as $N\to\infty$, then the MJSQ$(n(N))$ scheme has the same fluid limit as the ordinary JSQ policy.
\item Then we apply Corollary~\ref{cor:bound-ssy} to prove that as long as $n(N)/N\to 0$, \emph{any} scheme from the class CJSQ$(n(N))$ has the same fluid limit as the ordinary JSQ policy.
\item Next, using Propositions~\ref{prop:stoch-ord2} and~\ref{prop:differ} we establish that if $d(N)\to\infty$, then for \emph{some} $n(N)$ with $n(N)/N\to 0$, the  JSQ$(d(N))$ scheme and the JSQ$(n(N),d(N))$ scheme have the same fluid limit. The proposition then follows by observing that the JSQ$(n(N),d(N))$ scheme belongs to the class CJSQ$(n(N))$.
\end{enumerate}

\begin{proof}[Proof of Proposition~\ref{prop:samefluid}]
First, to show Claim~(i) above, define $\bar{N}=N-n(N)$ and $\bar{\lambda}(\bar{N})=\lambda(N)$.
Observe that the MJSQ$(n(N))$ scheme with $N$ servers can be thought of as the ordinary JSQ policy with $\bar{N}$ servers and arrival rate $\bar{\lambda}(\bar{N})$.
Also, since $n(N)/N\to 0$,
\begin{align*}
\frac{\bar{\lambda}(\bar{N})}{\bar{N}}=\frac{\lambda(N)}{N-n(N)}\to \lambda\quad \text{as}\quad \bar{N}\to\infty.
\end{align*}
Furthermore, observe that the fluid limit of the JSQ policy in Theorem~\ref{th:genfluid1} as given by~\eqref{eq:fluidfinal} is characterized by the parameter $\lambda$ only, and hence the fluid limit of the MJSQ$(n(N))$ scheme is the same as that of the ordinary JSQ policy.

Second, observe from the fluid limit of the JSQ policy that if $\lambda< 1$, then for any buffer capacity $b\geq 1$, and any starting state, the fluid-scaled cumulative overflow is negligible, i.e., for any $t\geq 0$, $L^N(t)/N\pto 0$.
Since the above fact is induced by the fluid limit only, the same holds for the MJSQ$(n(N))$ scheme.
Therefore, using the lower and upper bounds in Corollary~\ref{cor:bound-ssy} and the tail bound in Proposition~\ref{prop:stoch-ord}, we obtain Claim~(ii) above.

Finally, choose $n(N)=N/\sqrt{d(N)}$, and consider the JSQ$(n(N),d(N))$ scheme. 
Since $d(N)\to\infty$, it is clear that $n(N)/N\to 0$ as $N\to\infty$.
Also, if $\Delta^N(T)$ denotes the cumulative number of times that the JSQ($d(N)$) scheme and JSQ$(n(N),d(N))$ scheme differ in decision up to time $T$, then Proposition~\ref{prop:differ} yields
\begin{align*}
\Pro{\Delta^N(T)\geq \varepsilon N\given A^N(T)}&\leq \frac{A^N(T)}{\varepsilon N}\left(1-\frac{n(N)}{N}\right)^{d(N)}\\
&=\frac{A^N(T)}{\varepsilon N}\left(1-\frac{1}{\sqrt{d(N)}}\right)^{d(N)}.
\end{align*}
Since $\big\{A^N(T)/N\big\}_{N\geq 1}$ is a tight sequence of random variables, we have
\begin{align*}
\frac{A^N(T)}{\varepsilon N}\left(1-\frac{1}{\sqrt{d(N)}}\right)^{d(N)}\pto 0
\quad\text{as}\quad N\to\infty,
\end{align*}
and hence, $\Delta^N(T)/N\pto 0$. 
Therefore, applying the $\ell_1$ distance bound stated in Proposition~\ref{prop:stoch-ord2}, we obtain Claim~(iii).
The proof is then completed by observing that the JSQ$(n(N),d(N))$ scheme belongs to the class CJSQ$(n(N))$.
\end{proof}

\begin{proof}[Proof of Theorem~\ref{th:batch-ssy}]
For any $\varepsilon>0$, define 
$$T^N_\varepsilon:=\inf\big\{t:Q_1^{\sss d(N)}(t)>(\lambda+\varepsilon)N\big\}.$$
Now the proof consists of two main steps. First we show that if $d(N)\geq \ell(N)/(1-\lambda-\varepsilon)$ for some $\varepsilon>0$, then there exists an $\varepsilon'>0$, such that if for some $T>0$, $\Pro{T^N_{\varepsilon'}>T}\to 1$ as $N\to\infty$,
 then the number of times that the JSQ($d(N)$) scheme and the ordinary JSQ policy differ in decision in $[0,T]$ is $\op(N)$. This then implies that up to such a time $T$, it is enough to consider the fluid limit of the ordinary JSQ policy with batch arrivals. Second, we show that if the conditions stated in Theorem~\ref{th:batch-ssy} hold, then for any finite time $T> 0$, $\Pro{T^N_{\varepsilon'}>T}\to 1$ as $N\to\infty$. This will complete the proof.

To prove the first part, consider the JSQ($d(N)$) scheme in case of batch arrivals. Choose $\varepsilon'=\varepsilon/2$, and assume that $T>0$ is such that $\Pro{T^N_{\varepsilon'}>T}\to 1$ as $N\to\infty$.
Let $I_i$ denote the number of idle servers among $d(N)$ randomly chosen servers for the $i^{\mathrm{th}}$ batch arrival, and define $W^N(t)$ to be the cumulative number of tasks that have not been assigned to some idle server, up to time $t$. If $A^N(t)$ denotes the number of batch arrivals that occurred up to time $t$, then 
$$W^N(t)=\sum_{i=1}^{A^N(t)}[\ell(N)-I_i]^+.$$
We show that $W^N(t)/N\pto 0$ for all $t\leq T$ for $d(N)=\ell(N)/(1-\lambda-\varepsilon)$. Observe that $I_i$ follows a Hypergeometric distribution with sample size $d(N)$, and population size $N$ containing  $N-Q_1^N(t)\geq (1-\lambda-\varepsilon/2) N$ successes. Define $J_i$ to be distributed as $d(N)-I_i$. Then
$$[\ell(N)-I_i]^+=k\iff J_i=d(N)-\ell(N)+k.$$
Therefore, for $c=1-\lambda-\varepsilon/2$ we have,
$$\expect{[\ell(N)-I_i]^+}=\sum_{k\geq 1}k\Pro{J_i=(1-c)d(N)+k}\leq d(N)\Pro{J_i\geq (1-c)d(N)}.$$
Now, from \cite{H63, LP14}, we know
$$\Pro{J_i\geq (1-c)d(N)}\leq \exp(-d(N)H(\lambda,c)),$$
where $$H(\lambda,c)=(1-c)\log\left(\frac{1-c}{\lambda}\right)+c\log\left(\frac{c}{1-\lambda}\right)>0,$$
since $c<1-\lambda.$
Therefore, 
\begin{equation}
\begin{split}
\Pro{W^N(t)>\varepsilon N}
&\leq\frac{\expect{W^N(t)}}{\varepsilon N}\\
&\leq \frac{d(N)}{\varepsilon N}\times  \frac{\lambda(N) t}{\ell(N)}\times \exp(-d(N)H(\lambda,c))\\
&=O(\exp(-d(N)H(\lambda,c))).
\end{split}
\end{equation}
This implies that whenever $\ell(N)\to\infty$, if $d(N)=\ell(N)/(1-\lambda-\varepsilon/2)$, then $W^N(t)$ is $o_P(N)$ for all $t\leq T$. 
Now the analysis of the batch arrivals with ordinary JSQ policy in Theorem~\ref{th:batchjsq} below, up to time $T$, shows that the process
$\big\{\qq^N(t)\big\}_{0\leq t\leq T}$ converges to the deterministic limit $\big\{\qq(t)\big\}_{0\leq t\leq T}$, described by \eqref{eq:batch}.

Therefore, it is enough to show that any $T>0$ satisfies the required criterion. This can be seen by observing that for any $T\geq0$, and any $\varepsilon'>0$,
\begin{align*}
&\Pro{T_{\varepsilon'}^N\leq T}\leq \Pro{T_{\varepsilon'/2}^N< T}\\
&\leq\Pro{\sup_{t\in[0,T]}Q_1^{\sss d(N)}(t)>(\lambda+\varepsilon'/2)N}\\
&\leq\Pro{\sup_{t\in[0,T]}Q_1^\JSQ(t)>\Big(\lambda+\frac{\varepsilon'}{4}\Big)N}\Pro{\sup_{t\in[0,T]}|Q_1^\JSQ(t)-Q_1^{\sss d(N)}(t)|\leq \frac{N\varepsilon'}{4}}\\
&\hspace{2cm}+\Pro{\sup_{t\in[0,T]}|Q_1^\JSQ(t)-Q_1^{\sss d(N)}(t)|> \frac{N\varepsilon'}{4}}\longrightarrow 0\quad\mathrm{as}\quad N\to\infty.
\end{align*}
Therefore the proof is complete.
\end{proof}

\begin{theorem}{\normalfont (Batch arrivals JSQ)}
\label{th:batchjsq}
Consider the batch arrival scenario with growing batch size $\ell(N)\to\infty$ and $\lambda(N)/N\to\lambda<1$ as $N\to\infty$. For the JSQ policy, if $q^{\sss d(N)}_1(0)\pto q_1^\infty\leq \lambda$, and $q_i^{\sss d(N)}(0)\pto 0$ for all $i\geq 2$, then the sequence of processes
$\big\{\qq^{\sss d(N)}(t)\big\}_{t\geq 0}$ converges weakly to the limit $\big\{\qq(t)\big\}_{t\geq 0}$, described as follows: 
\begin{equation}\label{eq:batchjsq}
q_1(t) = \lambda + (q_1^\infty-\lambda)\e^{-t},\quad
q_i(t)\equiv 0\quad \mathrm{for\ all}\quad i= 2,\ldots,b.
\end{equation}
\end{theorem}
\begin{proof}
Fix any finite time $T\geq 0$.
To analyze the JSQ policy with batch arrivals, observe that before time $T$, all the arriving tasks join idle servers. Therefore, assuming $Q_2^N(0)=0$, for all $t\leq T$, the evolution for $Q_1^N$ can be written as
\begin{equation}
Q_1^N(t)=Q_1^N(0)+\ell(N)A\left(t\lambda(N)/\ell(N)\right)-D\left(\int_0^tQ_1^N(s)ds\right),
\end{equation}
where $A$ and $D$ are independent unit-rate Poisson processes.
Using the random time change of unit-rate Poisson processes \cite[Lemma~3.2]{PTRW07}, and applying the arguments in \cite[Lemma~3.4]{PTRW07}, the above process scaled by $N$, then admits the martingale decomposition
\begin{equation}\label{eq:fluid-batch}
q_1^N(t)=q_1^N(0)+\frac{M^N_1(t)}{N}+\lambda t-\frac{M^N_2(t)}{N}-\int_0^tq_1^N(s)ds,
\end{equation}
where 
\begin{align*}
M^N_1(t)&=\ell(N)A\left(t\lambda(N)/\ell(N)\right)-t\lambda(N),\\
M^N_2(t)&=D\left(\int_0^tQ_1^N(s)ds\right)-\int_0^tQ_1^N(s)ds,
\end{align*}
are square integrable martingales with respective quadratic variation processes given by
\begin{align*}
\langle M^N_1\rangle(t)&=t\lambda(N),\\
\langle M^N_2\rangle(t)&=\int_0^tQ_1^N(s)ds.
\end{align*}
Now, since for any $T\geq 0$, $\langle M^N_1\rangle(T)/N^2\to 0$, and $\langle M^N_2\rangle(T)/N^2\pto 0$, from the stochastic boundedness criterion for square integrable martingales \cite[Lemma~5.8]{PTRW07}, we get that both $\big\{M_1^N(t)/N\big\}_{t\geq 0}\dto 0$ and $\big\{M_2^N(t)/N\big\}_{t\geq 0}\dto 0$. Therefore, from the continuous mapping theorem and \eqref{eq:fluid-batch}, it follows that $\big\{q_1^N(t)\big\}_{t\geq 0}$ as $N\to\infty$ converges weakly to a deterministic limit described by the integral equation
\begin{equation}
q_1(t)=q_1^\infty+\lambda t-\int_0^tq_1(s)ds
\end{equation}
having~\eqref{eq:batch} as the unique solution.
This completes the proof of the fluid limit of JSQ with batch arrivals.
\end{proof}

\subsection{Global stability and interchange of limits}\label{ssec:globstab}

To prove the interchange of limits result stated in Proposition~\ref{th:interchange}, we will establish the global stability of the fixed point, i.e., all fluid paths converge to the fixed point in~\eqref{eq:fixed point-ssy} as $t\to\infty$. 
This is formally stated in the following lemma.
\begin{lemma}\label{lem:global-stab}
Let $\qq(t)$ be the fluid limit, i.e., the solution of the dynamical system described by the system of integral equations in~\eqref{eq:fluid}.
For any $\qq^\infty\in \mathcal{S}$, if $\qq(0) = \qq^\infty$, then $\qq(t)\to
\qq^*$ as $t\to\infty$, where $\qq^*$ is defined as in \eqref{eq:fixed point-ssy}.
\end{lemma}
In case of the JSQ$(d)$ scheme with fixed $d$, the global stability is proved by constructing a Lyapunov function that measures the `distance' (in terms of a weighted $L_1$-norm) between the trajectory and the fixed point, and that strictly decreases everywhere except at the fixed point, see~\cite[Theorem 3.6]{Mitzenmacher96}.
In case of the ordinary JSQ policy however, we can exploit a more direct method to establish the global stability, as further detailed below.
\begin{proof}[Proof of Lemma~\ref{lem:global-stab}]
The proof follows in two steps: we will first establish that as $t\to\infty$, $q_1(t)\to\lambda<1$, and then show that $q_2(t)\to 0$.

Observe that the rate of change of $q_1(t)$ is $\lambda p_0(\qq(t))-(q_1(t)-q_2(t))$.
For any $\varepsilon\geq 0$, if $q_1(t)\leq \lambda-\varepsilon$, then $p_0(\qq(t))=1$, so that the rate of change is  $\lambda -(q_1(t)-q_2(t))\geq \varepsilon$, i.e., positive and bounded away from zero when $\varepsilon>0$.
Also, $q_1(t)$ cannot decrease if $q_1(t)\leq\lambda$.
This shows that for all $\varepsilon>0$, there exists a time $t_0 = t_0(\varepsilon, \qq^\infty)$, such that, $q_1(t)\geq \lambda -\varepsilon$ for all $t\geq t_0$.
Thus, $\liminf_{t\to\infty} q_1(t)\geq\lambda$.

On the other hand, we claim that $\limsup_{t\to\infty} q_1(t)\leq \lambda$.
Suppose not, i.e., assume $\limsup_{t\to\infty}q_1(t) = \lambda+\varepsilon$ for some $\varepsilon>0$.
Because $q_1(t)$ is non-decreasing when $q_1(t)\leq \lambda$, there must exist a $t_0$ such that $q_1(t)\geq \lambda$ $\forall\ t\geq t_0$.
The high-level idea behind the claim is as follows.
If $q_1(t)$ were to remain above $\lambda$ by a non-vanishing margin, then the cumulative number of departures would exceed the cumulative number of arrivals by an infinite amount, which cannot occur since the initial number of tasks is bounded.
More formally,
\begin{align*}
\sum_{i=1}^b q_i(t) &= \sum_{i=1}^bq_i(t_0) + \lambda\int_{t_0}^t\sum_{i=1}^b p_{i-1}(\qq(s))\dif s - \int_{t_0}^t q_1(s)\dif s\\
& \leq \sum_{i=1}^bq_i(t_0)  - \int_{t_0}^t [q_1(s) - \lambda]^+\dif s,
\end{align*}
and thus,
$$\int_{t_0}^t [q_1(s) - \lambda]^+\dif s\leq \sum_{i=1}^b q_i(t) - \sum_{i=1}^bq_i(t_0)<\infty.$$
This provides a contradiction with $\limsup_{t\to\infty} q_1(t) = \lambda +\varepsilon$, since the rate of decrease of $q_1(t)$ is at most~1.
Therefore, $q_1(t)\to \lambda$ as $t\to\infty$.  

Consequently, for any $\qq^\infty\in\mathcal{S}$ and $\varepsilon>0$, if $\qq(0) = \qq^\infty$, then there exists a time $t_2 = t_2(\qq^\infty, \varepsilon)<\infty$, such that $q_1(t)\leq \lambda+\varepsilon$ for all $t\geq t_2$.
Thus choosing $\varepsilon = (1-\lambda)/2$ say, for all $t\geq t_2$, $q_1(t)<1$, and thus $p_0(\qq(t))=1$, i.e., $\sum_{i=2}^b p_{i-1}(\qq(t))=0$.
Define $q_{2+}(t) := \sum_{i = 2}^b q_i(t)$.
Observe that
\begin{align*}
q_{2+}(t)&= q_{2+}(t_2)+\lambda\int_{t_2}^t\sum_{i=2}^b p_{i-1}(\qq(s))\dif s- \int_{t_2}^tq_2(s)\dif s\\
&=  q_{2+}(t_2)- \int_{t_2}^tq_2(s)\dif s\qquad \mbox{for all}\quad t\geq t_2,
\end{align*}
which implies $q_2(t)\leq q_{2+}(t_2)\e^{-(t-t_2)}$.
Thus, $q_2(t)$ and consequently, $q_{2+}(t)$ converges to 0 as $t\to\infty$.
This completes the proof of global stability of the fixed point.
\end{proof}

\begin{proof}[Proof of Proposition~\ref{th:interchange}]
The proof follows in two steps: (i) we first establish that the sequence of stationary measures $\big\{\pi^{\sss d(N)}\big\}_{N\geq 1}$ is tight, and then (ii) show the interchange of limits.

(i) Observe that if $b<\infty,$ then the space $[0,1]^b$ is compact, and hence Prohorov's theorem implies that $\big\{\pi^{\sss d(N)}\big\}_{N\geq 1}$ is tight. 
Now assume $b=\infty.$
For any two positive integers $d_1\leq d_2$, note that at each arrival, the JSQ$(d_2)$ scheme polls more servers than the JSQ$(d_1)$ scheme.
Thus using the S-coupling and Proposition~\ref{prop:det ord}, we can conclude for every~$N$,
$$\sum_{i\geq m} Q_i^{d_2}\leq_{st}\sum_{i\geq m} Q_i^{d_1},\quad \mbox{for all}\quad m \geq 1.$$
In particular, putting $d_1=1$ and $d_2 = d(N)$,
\begin{equation}\label{eq:1vsdN}
\sum_{i\geq m} Q_i^{d(N)}\leq_{st}\sum_{i\geq m} Q_i^{1},\quad \mbox{for all}\quad m \geq 1.
\end{equation}
Let $\XX^N$ and $\YY^N$ denote random variables following the stationary distribution of two systems with $N$ servers under the JSQ$(d(N))$ and JSQ$(1)$ schemes, respectively. 
We will verify the tightness criterion stated in Lemma~\ref{lem:tightcond}.
Note that since $\XX^N$ takes value in $\mathcal{S}\subset [0,1]^\infty$, which is compact with respect to the product topology, Prohorov's theorem implies that $\big\{\XX^N\big\}_{N\geq 1}$ is tight with respect to the product topology.
To verify the condition in~\eqref{eq:smalltail}, note that the system under the JSQ$(1)$ scheme is essentially a collection of $N$ independent M/M/1 systems. Therefore, for each $k\geq 1$,
\begin{align*}
\varlimsup_{N\to\infty}\mathbbm{P}\Big(\sum_{i\geq k}X_i^N>\varepsilon\Big)
\leq \varlimsup_{N\to\infty}\mathbbm{P}\Big(\sum_{i\geq k}Y_i^N>\varepsilon\Big)
= (1-\lambda)\sum_{i\geq k}\lambda^i.
\end{align*}
Since $\lambda<1$, taking the limit $k\to\infty$, the right side of the above inequality tends to zero, and hence, the condition in~\eqref{eq:smalltail} is verified.

(ii) Now observe that since $\big\{\pi^{\sss d(N)}\big\}_{N\geq 1}$ is tight, any subsequence has a convergent further subsequence. 
Let $\big\{\pi^{\sss d(N_n)}\big\}_{n\geq 1}$ be any such convergent subsequence, with $\big\{N_n\big\}_{n\geq 1}\subseteq\N$, such that $\pi^{\sss d(N_n)}\dto\hat{\pi}$ as $n\to\infty$. We will show that $\hat{\pi}$ is unique and equals the measure $\pi^\star$, as defined in the statement of Proposition~\ref{th:interchange}.
Notice that if $\qq^{\sss d(N_n)}(0)\sim\pi^{\sss d(N_n)}$, then $\qq^{\sss d(N_n)}(t)\sim\pi^{\sss d(N_n)}$ for all $t\geq 0$.
Thus, $\hat{\pi}$ is an invariant distribution of the deterministic process $\big\{\qq(t)\big\}_{t\geq 0}$.
 This in conjunction with the global stability in Lemma~\ref{lem:global-stab} implies that $\hat{\pi}$ must be the fixed point of the fluid limit.  
Thus, we have shown the convergence of the stationary measure.
\end{proof}

\section{Diffusion-limit proofs}\label{sec:diffusion} 
In this section we prove the diffusion-limit results for the JSQ$(d(N))$ scheme stated in Theorem~\ref{diffusionjsqd-ssy}, and the almost necessity condition for diffusion-level optimality stated in Theorem~\ref{th:diff necessary}.
As noted in Subsection~\ref{subsec:strategy}, the diffusion limit 
for the ordinary JSQ policy is obtained in~\cite[Theorem~2]{EG15}, and characterized by~\eqref{eq:diffusionjsqd-ssy}.
Therefore it suffices to prove the universality property stated in the next proposition. 
\begin{proposition}\label{prop:samedif}
If $d(N)/(\sqrt{N}\log(N))\to\infty$ as $N\to\infty$, then the JSQ$(d(N))$ scheme and the ordinary JSQ policy have the same diffusion limit.
\end{proposition}
The proof of the above proposition follows similar lines as that of Proposition~\ref{prop:samefluid}, leveraging again the S-coupling results from Section~\ref{sec:coupling}, and
involves three steps:
\begin{enumerate}[{\normalfont (i)}]
\item First we show that if $n(N)/\sqrt{N}\to 0$ as $N\to\infty$, then the MJSQ$(n(N))$ scheme has the same diffusion limit as the ordinary JSQ policy.
\item Then we use Corollary~\ref{cor:bound-ssy} to prove that as long as $n(N)/\sqrt{N}\to 0$, \emph{any} scheme from the class CJSQ$(n(N))$ has the same diffusion limit as the ordinary JSQ policy.
\item Next we establish using Propositions~\ref{prop:stoch-ord2} and~\ref{prop:differ} that if $d(N)$ is such that $d(N)/(\sqrt{N}\log(N))\to\infty$ as $N\to\infty$, then for \emph{some} $n(N)$ with $n(N)/\sqrt{N}\to 0$, the  JSQ$(d(N))$ scheme and the JSQ$(n(N),d(N))$ scheme have the same diffusion limit. The proposition then follows by observing that the JSQ$(n(N),d(N))$ scheme belongs to the class CJSQ$(n(N))$.
\end{enumerate}

\begin{proof}[Proof of Proposition~\ref{prop:samedif}]
To show Claim~(i) above, define $\bar{N}=N-n(N)$ and $\bar{\lambda}(\bar{N})=\lambda(N)$.
As mentioned earlier, the MJSQ$(n(N))$ scheme with $N$ servers can be thought of as the ordinary JSQ policy with $\bar{N}$ servers and arrival rate $\bar{\lambda}(\bar{N})$.
Also, since $n(N)/\sqrt{N}\to 0$,
\begin{align*}
\frac{\bar{N}-\bar{\lambda}(\sqrt{\bar{N}})}{\bar{N}}=\frac{N-n(N)-\lambda(N)}{\sqrt{N-n(N)}}\to \beta>0\quad \text{as}\quad \bar{N}\to\infty.
\end{align*}
Furthermore, observe that the diffusion limit of the JSQ policy in \cite[Theorem~2]{EG15} as given in~\eqref{eq:diffusionjsqd-ssy} is characterized by the parameter $\beta>0$, and hence the diffusion limit of the MJSQ$(n(N))$ scheme is the same as that of the ordinary JSQ policy.

Observe from the diffusion limit of the JSQ policy that if $\beta>0$, then for any buffer capacity $b\geq 2$, and suitable initial state as described in~Theorem~\ref{diffusionjsqd-ssy}, the cumulative overflow is negligible, i.e., for any $t\geq 0$, $L^N(t)\pto 0$.
Indeed observe that if $b\geq 2$, and $\big\{\bQ_2^N(0)\big\}_{N\geq 1}$ is a tight sequence, then the sequence of processes $\big\{\bQ_2^N(t)\big\}_{t\geq 0}$ is stochastically bounded.
Therefore, on any finite time interval, there will be only $\Op(\sqrt{N})$ servers with queue length more than one, whereas, for an overflow event to occur all the $N$ servers must have at least two pending tasks.
Therefore, for any $t\geq 0$,
\begin{align*}
\limsup_{N\to\infty}\Pro{L^N(t)>0}&\leq \limsup_{N\to\infty}\Pro{\sup_{s\in[0,t]}Q^N_2(s)=N}\\
&\leq \limsup_{N\to\infty}\Pro{\sup_{s\in[0,t]}\bQ^N_2(s)=\sqrt{N}}=0.
\end{align*} 
Since the above fact is implied by the diffusion limit only, the same holds for the MJSQ$(n(N))$ scheme.
Therefore, using the lower and upper bounds in Corollary~\ref{cor:bound-ssy} we arrive at Claim~(ii).

Finally, choose
$$n(N)=\frac{N\log N}{d(N)},$$
 and consider the JSQ$(n(N),d(N))$ scheme. 
Since $d(N)/( \sqrt{N}\log N)\to\infty$, it is clear that $n(N)/\sqrt{N}\to 0$ as $N\to\infty$.
Again, if $\Delta^N(T)$ denotes the cumulative number of times that the JSQ($d(N)$) scheme and JSQ$(n(N),d(N))$ scheme differ in decision up to time $T$, then Proposition~\ref{prop:differ} yields
\begin{equation}\label{eq:propdif}
\begin{split}
\Pro{\Delta^N(T)\geq \varepsilon \sqrt{N}\given A^N(T)}&\leq \frac{A^N(T)}{\varepsilon \sqrt{N}}\left(1-\frac{n(N)}{N}\right)^{ d(N)}\\
&\leq \frac{A^N(T)}{\varepsilon \sqrt{N}}\left(1-\frac{\log(N)}{d(N)}\right)^{ d(N)}\\
&\leq \frac{A^N(T)}{\varepsilon N} \sqrt{N}\left(1-\frac{\log(N)}{d(N)}\right)^{ d(N)}.
\end{split}
\end{equation}
Since $\big\{A^N(T)/N\big\}_{N\geq 1}$ is a tight sequence of random variables, and
\begin{align*}
&\sqrt{N}\left(1-\frac{\log(N)}{d(N)}\right)^{ d(N)}\to 0,
\quad\text{as}\quad N\to\infty,\\
\iff &\frac{1}{2}\log(N)+d(N)\log\left(1-\frac{\log(N)}{d(N)}\right)\to-\infty,\quad\text{as}\quad N\to\infty,\\
\Longleftarrow\hspace{.15cm}& \frac{1}{2}\log N - \frac{\log(N)}{d(N)}\times d(N)\to-\infty,\quad\text{as}\quad N\to\infty,
\end{align*}
from~\eqref{eq:propdif}, $\Delta^N(T)/N\pto 0$. 
Therefore, by invoking Proposition~\ref{prop:stoch-ord2}, we obtain Claim~(iii).
The proof is then completed by observing that the JSQ$(n(N),d(N))$ scheme belongs to the class CJSQ$(n(N))$.
\end{proof}

We next prove that the growth condition $d(N)/( \sqrt{N}\log N)\to\infty$ is nearly necessary: for any $d(N)$ such that $d(N)/(\sqrt{N}\log N)\to 0$ as $N\to\infty$, the diffusion limit of the JSQ$(d(N))$ scheme differs from that of the ordinary JSQ policy. 
Note that it is enough to consider the truncated system where any arrival to a server with at least two tasks is discarded, since the truncated system and the original system have the same diffusion limit~\cite{MBLW16-1}. 

Now consider the JSQ($d(N)$) scheme for some $d(N)$ with $d(N)/(\sqrt{N}\log N)\to 0$ as $N\to\infty$, and assume on the contrary, the hypothesis that the process 
$$\big\{(Q_1^{\sss d(N)}(t)-N)/\sqrt{N}, Q_2^{\sss d(N)}(t)/\sqrt{N}\big\}_{t\geq 0}$$ converges to the diffusion limit corresponding that of the JSQ policy.
From a high level, the idea is to show that if the processes $(N-Q^{\sss d(N)}_1(\cdot))$ and $Q^{\sss d(N)}_2(\cdot)$ are $\Op(\sqrt{N})$, then in any finite time interval the number of tasks assigned to a server with queue length at least one, by the JSQ$(d(N))$ scheme with $d(N)/(\sqrt{N}\log N)\to 0$ does not scale with $\sqrt{N}$, which then immediately proves that the diffusion limit cannot coincide with that of the ordinary JSQ policy.

To formalize the above idea, we first define an artificial scheme below, which will serve as an asymptotic  lower bound to the number of servers with queue length two in a system following the JSQ$(d(N))$ scheme, under the hypothesis that the diffusion limit of the JSQ$(d(N))$ coincides with that of the ordinary JSQ policy.
For any nonnegative sequence $c(N)$, define a scheme $\Pi(c(N))$ which 
\begin{enumerate}[{\normalfont (i)}]
\item At each external arrival, assigns the task to a server having queue length one with probability $(1-c(N)/N)^{\sss d(N)}$, and else discards it (ties can be broken randomly),
\item If a departure occurs from a server with queue length one, then it immediately makes the server busy with a dummy arrival, i.e., essentially $\Pi(c(N))$ prohibits any server to remain idle.
\end{enumerate}
We use a coupling argument to show the following:
\begin{lemma}\label{lem:upper}
For any nonnegative sequence $c(N)$ with $c(N)/ \sqrt{N}\to \infty$ as $N\to\infty$, there exists a common probability space, such that for any $T>0$, 
$$\Pro{\sup_{t\in [0,T]}\left\{Q_2^{\sss d(N)}(t)-Q_2^{\sss \Pi(c(N))}(t)\right\}\geq 0}\longrightarrow 1\quad\mathrm{as}\quad N\to\infty,$$
provided $Q_2^{\sss d(N)}(0)\geq Q_2^{\sss \Pi(c(N))}(0)$ for all sufficiently large $N$,
and the hypothesis that the sequences of processes $\big\{(N-Q_1^{\sss d(N)}(t))/\sqrt{N}\big\}_{t\geq 0}$ and $\big\{Q_2^{\sss d(N)}(t)/\sqrt{N}\big\}_{t\geq 0}$ are stochastically bounded.
\end{lemma}
In order to prove Lemma~\ref{lem:upper}, we first S-couple the two systems under schemes $\Pi(c(N))$ and JSQ($d(N)$) respectively. 
Now at each external arrival, to assign the task in the two systems in a coupled way, draw a single uniform$[0,1]$ random variable $U$, independent of any other processes. 
\begin{itemize}
\item Under the JSQ($d(N)$) scheme, if $u<(Q_1^{\sss d(N)}/N)^{\sss d(N)}-(Q_2^{\sss d(N)}/N)^{\sss d(N)}$, assign the task to a server with queue length one, if 
\begin{equation}\label{eq:idleserver}
(Q_1^{\sss d(N)}/N)^{\sss d(N)}-(Q_2^{\sss d(N)}/N)^{\sss d(N)}<U<1-(Q_2^{\sss d(N)}/N)^{\sss d(N)},
\end{equation}
then assign the task to an idle server, and otherwise discard it. 
This preserves the statistical law of the JSQ($d(N)$) scheme with a buffer size $b=2$.
Indeed note that according to the above rule the probability that an incoming task will be assigned
to some server with queue length zero, one, and two, are respectively given by $(1-(Q_1^{\sss d(N)}/N)^{\sss d(N)})$, $(Q_1^{\sss d(N)}/N)^{\sss d(N)}-(Q_2^{\sss d(N)}/N)^{\sss d(N)}$, and $(Q_2^{\sss d(N)}/N)^{\sss d(N)}$.
\item Under the scheme $\Pi(c(N))$, if $U<(1-c(N)/N)^{\sss d(N)}$, assign the incoming task to a server with queue length one, otherwise discard it.
Clearly, the statistical law of the $\Pi(c(N))$ scheme is preserved by this rule.
\end{itemize}
\begin{proof}[Proof of Lemma~\ref{lem:upper}]
Fix any $T\geq 0$. Now the proof follows in two steps:

(i) First assume that at each external arrival up to time $T$, whenever an incoming task joins a server with queue length one, under the $\Pi(c(N))$ scheme, then so does the incoming task under the JSQ($d(N)$) scheme. 
In that case, since the two systems are S-coupled, by forward induction on event times, it can be seen that $Q_2^{\sss d(N)}(t)\geq Q_2^{\sss \Pi(c(N))}(t)$ for all $0\leq t\leq T$, provided $Q_2^{\sss d(N)}(0)\geq Q_2^{\sss \Pi(c(N))}(0)$.

(ii) 
Now, for any $T\geq0$, 
according to the hypothesis, both $\sup_{t\in [0,T]}Q_2^{\sss d(N)}(t)$ and $\sup_{t\in[0,T]}\big\{N-Q_1^{\sss d(N)}(t)\big\}$  are $\Op(\sqrt{N})$. 
Also, since $c(N)/\sqrt{N}\to\infty$, it is straightforward to check that
\begin{equation}\label{eq:pi(c(N))}
\begin{split}
&\liminf_{N\to\infty}\P\Big(\sup_{t\in [0,T]}\left(Q_1^{\sss d(N)}(t)/N\right)^{\sss d(N)}-\left(Q_2^{\sss d(N)}(t)/N\right)^{\sss d(N)}\\
&\hspace{6cm}\geq \left(1-c(N)/N\right)^{\sss d(N)}\Big)=1.
\end{split}
\end{equation}

Note that the probabilities that an incoming task joins a server with queue length one are given by the quantities
$(Q_1^{\sss d(N)}(t)/N)^{\sss d(N)}-(Q_2^{\sss d(N)}(t)/N)^{\sss d(N)}$ and $(1-c(N)/N)^{\sss d(N)}$ for the JSQ($d(N)$) and the $\Pi(c(N))$ scheme, respectively. 
Informally speaking, due to the above coupling,~\eqref{eq:pi(c(N))} then implies that with high probability, on any finite time interval, whenever an external incoming task joins a server with queue length one under the $\Pi(c(N))$ scheme, then so does the incoming task under the JSQ($d(N)$) scheme. 
 Therefore, from Part~(i) above, we can say
 \begin{align*}
&\Pro{\sup_{t\in [0,T]}\big\{Q_2^{\sss d(N)}(t)-Q_2^{\sss \Pi(c(N))}(t)\big\}\geq 0}\\
 &\geq \Pro{\sup_{t\in [0,T]}\left(\frac{Q_1^{\sss d(N)}(t)}{N}\right)^{\sss d(N)}-\left(\frac{Q_2^{\sss d(N)}(t)}{N}\right)^{\sss d(N)}\geq \left(1-\frac{c(N)}{N}\right)^{\sss d(N)}}\\
 &\longrightarrow 1 \quad\text{as}\quad N\to\infty.
 \end{align*}

Thus for any $T\geq 0$, if $\sup_{t\in [0,T]}Q_2^{\sss d(N)}(t)$ and $\sup_{t\in[0,T]}\left\{N-Q_1^{\sss d(N)}(t)\right\}$  are  $\Op(\sqrt{N})$, then with probability tending to one as $N\to\infty$, up to time $t$, the process $\big\{Q_2^{\sss \Pi(c(N))}(t)\big\}_{0\leq t\leq T}$ is indeed a lower bound for  $\big\{Q_2^{\sss d(N)}(t)\big\}_{0\leq t\leq T}$, and hence by our hypothesis, the proof is complete. 
\end{proof}

\begin{proof}[Proof of Theorem~\ref{th:diff necessary}]
Fix any sequence $d(N)$ such that $d(N)/(\sqrt{N}\log N)\to 0$ as $N\to\infty$.
Assume the hypothesis that for the JSQ($d(N)$) scheme, the process
$$\big\{(Q_1^{\sss d(N)}(t)-N)/\sqrt{N}, Q_2^{\sss d(N)}(t)/\sqrt{N}\big\}_{t\geq 0}$$ 
converges to the appropriate diffusion limit corresponding to that of the ordinary JSQ policy.
We will show that under this hypothesis, the process $\big\{Q_2^{\sss d(N)}(t)/\sqrt{N}\big\}_{t\geq 0}$ is not stochastically bounded, which will then lead to a contradiction. 

In order to show this, we will choose an appropriate $c(N)$ such that $c(N)/\sqrt{N}\to\infty$ as $N\to\infty$, and the process $\big\{Q_2^{\sss \Pi(c(N))}(t)/\sqrt{N}\big\}_{t\geq 0}$ is not stochastically bounded.
The conclusion then follows by the application of Lemma~\ref{lem:upper}.

Observe that the martingale decomposition of the scaled $Q_2^{\sss \Pi(c(N))}(\cdot)$ process can be written as
\begin{equation}\label{eq:mart-modif}
\begin{split}
\bQ_2^{\sss \Pi(c(N))}(t)=\bQ_2^{\sss \Pi(c(N))}(0)+\frac{M^N(t)}{\sqrt{N}}+\frac{\lambda(N) t}{\sqrt{N}}&\left(1-\frac{c(N)}{N}\right)^{\sss d(N)}\\
&-\int_0^t\bQ_2^{\sss \Pi(c(N))}(s)ds,
\end{split}
\end{equation}
where $\bQ^{\sss \Pi(c(N))}_2(t)=Q^{\sss \Pi(c(N))}_2(t)/\sqrt{N}.$ 
Now write $c(N)=g(N)\sqrt{N}$, for some $g(N)\to\infty$ (to be chosen later),  and $d(N)=\sqrt{N}\log(N)/\omega(N)$, where $\omega(N)=\sqrt{N}\log(N)/d(N)\to\infty$ as $N\to\infty$. Therefore, we write \eqref{eq:mart-modif} as
\begin{equation}
\begin{split}
\bQ_2^{\sss \Pi(c(N))}(t)=\bQ_2^{\sss \Pi(c(N))}(0)+\frac{M^N(t)}{\sqrt{N}}+\frac{\lambda(N) t}{\sqrt{N}}&\left(1-\frac{g(N)}{\sqrt{N}}\right)^{\frac{\sqrt{N}\log N}{\omega(N)}}\\
&-\int_0^t\bQ_2^{\sss \Pi(c(N))}(s)ds.
\end{split}
\end{equation}
Observe that for any $t\geq 0$,
\begin{align*}
&\lim_{N\to\infty}\frac{\lambda(N) t}{\sqrt{N}}\left(1-\frac{g(N)}{\sqrt{N}}\right)^{\frac{\sqrt{N}\log N}{\omega(N)}}\\
&=t\lim_{N\to\infty}\exp\left[\log (\sqrt{N}-\beta)+\frac{\sqrt{N}\log N}{\omega(N)}\log\left(1-\frac{g(N)}{\sqrt{N}}\right)\right]\\
&=t\lim_{N\to\infty}\exp\left[\log (\sqrt{N}-\beta)-\frac{g(N)\log N}{\omega(N)}-o\left(\frac{g(N)\log N}{\omega(N)}\right)\right]
\end{align*}
Choosing $g(N)$ such that $g(N)/\omega(N)\to 0$ implies 
$$\frac{\lambda(N) t}{\sqrt{N}}\left(1-\frac{g(N)}{\sqrt{N}}\right)^{\frac{\sqrt{N}\log N}{\omega(N)}}\to\infty,\quad \text{as}\quad N\to\infty.$$
Note that for any $\omega(N)$, this choice of $g(N)$ is feasible (choose $g(N)=\sqrt{\omega(N)}$, say).
Furthermore, the process $\big\{M^N(t)/\sqrt{N}\big\}_{t\geq 0}$ in \eqref{eq:mart-modif} is stochastically bounded due to the martingale FCLT \cite[Theorem~7.1]{EK2009} and our hypothesis. 
Now we can conclude that for the above choices of $g(N)$ and $\omega(N)$, the process $\big\{\bQ_2^{\sss \Pi(c(N))}(t)\big\}_{t\geq 0}$, and hence the process $\big\{\bQ_2^{\sss d(N)}(t)\big\}_{t\geq 0}$ (due to Lemma~\ref{lem:upper}) is not stochastically bounded. 
Therefore, the limit does not coincide with the limit of the scaled $Q_2^{\sss \jsq}$-process.
\end{proof}

\section{Conclusion}\label{sec:conclusion}
In this chapter we have established universality properties
for power-of-$d$ load balancing schemes in many-server systems.
Specifically, we considered a system of $N$~parallel exponential
servers and a single dispatcher which assigns arriving tasks to the
server with the shortest queue among $d(N)$ randomly selected servers.
We developed a novel stochastic coupling construction to bound the
difference in the queue length processes between the JSQ policy
($d = N$) and a scheme with an arbitrary value of~$d$.
As it turns out, a direct comparison between the JSQ policy
and a JSQ($d$) scheme is a significant challenge.
Hence, we adopted a two-stage approach based on a novel class
of schemes which always assign the incoming task to one of the
servers with the $n(N) + 1$ smallest number of tasks.
Just like the JSQ($d(N)$) scheme, these schemes may be thought
of as `sloppy' versions of the JSQ policy.
Indeed, the JSQ($d(N)$) scheme is guaranteed to identify the
server with the minimum number of tasks, but only among
a randomly sampled subset of $d(N)$ servers.
In contrast, the schemes in the above class only guarantee that
one of the $n(N) + 1$ servers with the smallest number of tasks
is selected, but across the entire system of $N$ servers.
We showed that the system occupancy processes for an intermediate
blend of these schemes are simultaneously close on a $g(N)$ scale
($g(N) = N$ or $g(N) = \sqrt{N}$) to both the JSQ policy
and the JSQ($d(N)$) scheme for suitably chosen values of $d(N)$
and $n(N)$ as function of $g(N)$.
Based on the latter asymptotic universality, it then sufficed to
establish the fluid and diffusion limits for the ordinary JSQ policy.
Thus deriving the fluid limit of the ordinary JSQ policy, and using the above coupling argument we establish the fluid limit of the JSQ$(d(N))$ scheme in a regime with
$d(N) \to \infty$ as $N \to \infty$, along with the corresponding
fixed point.
The fluid limit turns out not to depend on the exact growth rate
of $d(N)$, and in particular coincides with that for the ordinary JSQ policy.
We further leveraged the coupling to prove that the diffusion limit
in the Halfin-Whitt regime with $d(N)/(\sqrt{N} \log (N)) \to \infty$
as $N \to \infty$ corresponds to that for the JSQ policy.
These results indicate that the optimality of the JSQ policy can be
preserved at the fluid-level and diffusion-level while reducing the
overhead by nearly a factor~O($N$) and O($\sqrt{N} / \log(N)$),
respectively.
In future work we plan to extend the results to heterogeneous
servers and non-exponential service requirement distributions.


%% file: jiq.tex
\begin{abstract}
We consider a system of $N$~parallel queues with identical exponential service rates and a single dispatcher where tasks arrive as a Poisson process. When a task arrives, the dispatcher always assigns it to an idle server, if there is any, and to a server with the shortest queue among $d$~randomly selected servers otherwise ($1 \leq d \leq N$). This load balancing scheme subsumes the so-called Join-the-Idle Queue (JIQ) policy ($d = 1$) and the celebrated Join-the-Shortest Queue (JSQ) policy ($d = N$) as two crucial special cases.
We develop a stochastic coupling construction to obtain the diffusion limit of the queue process in the Halfin-Whitt heavy-traffic regime, and establish that it does not depend on the value of~$d$, implying that assigning tasks to idle servers is sufficient for diffusion level optimality.
\end{abstract}

\section{Introduction}

In this chapter we establish a universality property for a broad class of load balancing schemes in a many-server Halfin-Whitt heavy-traffic regime, as described in Section~\ref{token}. 
Specifically, we consider a family of load balancing schemes termed JIQ($d$), where the dispatcher always assigns an incoming task to an idle server, if there is any, and to a server with the shortest queue among $d$~uniformly at random selected servers otherwise. Observe that the JIQ($N$) scheme coincides with the ordinary JSQ policy, while the JIQ($1$) scheme corresponds to the so-called Join-the-Idle-Queue
(JIQ) policy considered in~\cite{BB08,LXKGLG11,Stolyar15}.


We exploit a stochastic coupling construction to extend the weak convergence result for the JSQ policy as established by Eschenfeldt and Gamarnik~\cite{EG15} to the entire class of JIQ($d$) policies. 
We specifically establish that the diffusion limit, rather surprisingly, does not depend on the value of~$d$ at all, so that in particular the JIQ and JSQ policies yield the same diffusion limit. 
The latter property implies that in a many-server heavy-traffic regime, ensuring that tasks are assigned to idle servers whenever possible, suffices to achieve optimality at the diffusion level, and not just at the fluid level as proved by Stolyar~\cite{Stolyar15} for the under-loaded scenario. It further suggests that using any additional queue length information beyond the knowledge of empty queues yields only limited performance gains in large-scale systems in the Halfin-Whitt heavy-traffic regime.

A coupling method was used in Chapter~\ref{chap:univjsqd} to establish fluid and diffusion-level optimality of JSQ($d(N)$) policies.
There the idea pivots on two key observations: 
(i)~For any scheme, if each arrival is assigned to \emph{approximately} the shortest queue, then the scheme can still retain its optimality on various scales, and
(ii)~For any two schemes, if on any finite time interval not \emph{too many} arrivals are assigned to different ordered servers, then they can have the same scaling limits. 
Combination of the above two ideas provided a coupling framework involving an intermediate class of schemes that enabled us to establish the asymptotic optimality results.
In the current chapter the stochastic comparison framework is inherently different.
Comparing the JIQ and JSQ policies in the Halfin-Whitt regime will be facilitated when viewed as follows:
(i) If there is an idle server in the system, both JIQ and JSQ perform similarly,
(ii)~Also, when there is no idle server and only $O(\sqrt{N})$ servers with queue length two, JSQ assigns the arriving task to a server with queue length one. 
In that case, since JIQ assigns at random, the probability that the task will land on a server with queue length two and thus acts differently than JSQ is $O(1/\sqrt{N})$.
First we show that on any finite time interval the number of times an arrival finds all servers busy is at most $O(\sqrt{N})$.
Hence, all the arrivals except an $O(1)$ of them are assigned in exactly the same manner in both JIQ and JSQ, which then leads to the same scaling limit for both policies with the same initial state condition.

The chapter is organized as follows.
In Section~\ref{sec: model descr-jap} we present a detailed model description and formulate the main result.
In Section~\ref{sec: coupling} we develop a stochastic coupling
construction to compare the system occupancy state under various task assignment policies.
We then combine in Section~\ref{sec: conv} the stochastic comparison results with some of the derivations in~\cite{EG15} to obtain the common diffusion limit and finally make a few concluding remarks in Section~\ref{sec:conclusion-jap}.

\section{Model description and main results}
\label{sec: model descr-jap}
Consider a system with $N$~parallel queues with independent and identical servers having unit-exponential service rates and a single \emph{dispatcher}.
Tasks arrive at the dispatcher as a Poisson process of rate $\lambda(N)$, and are instantaneously forwarded to one of the servers. Tasks can be queued at the various servers, possibly subject to a buffer capacity limit as further described below, but \emph{cannot} be queued at the dispatcher. The dispatcher always assigns an incoming task to an idle server, if there is any, and to a server with the shortest queue among $d$~uniformly at random selected servers otherwise ($1 \leq d \leq N$), ties being broken arbitrarily. The buffer capacity at each of the servers is~$b\geq 2$ (possibly infinite), and when a task is assigned to a server with $b$~pending tasks, it is instantly discarded.
As mentioned earlier, the above-described scheme coincides with
the ordinary JSQ policy when $d = N$,
and corresponds to the JIQ policy
considered in~\cite{BB08,LXKGLG11,Stolyar15} when $d = 1$.

We consider the Halfin-Whitt heavy-traffic regime where the arrival rate increases with the number of servers as $\lambda(N) = N-\beta\sqrt{N}$ for some $\beta>0$. We denote the class of above-described policies by $\Pi^{(N)}(d)$, where the superscript~$N$ indicates that the diversity parameter~$d$ is allowed to depend on the number of servers. For any policy $\Pi \in \Pi^{(N)}(d)$ and buffer size~$b$, let $\mathbf{Q}^\Pi = (Q_1^\Pi, Q_2^\Pi, \ldots, Q_b^\Pi)$, where $Q_i^\Pi$ is the number of servers with a queue length greater than or equal to $i = 1, \ldots, b$, including the possible task in service. Also, let $\mathbf{X}^\Pi = (X_1^\Pi, X_2^\Pi, \ldots, X_b^\Pi)$ be a properly centered and scaled version of the vector $\mathbf{Q}^{\Pi}$, with $X_1^\Pi = (Q_1^\Pi-N)/\sqrt{N}$ and $X_i^\Pi = Q_i^\Pi/\sqrt{N}$ for $i = 2, \dots, b$. The reason why $Q_1^\Pi$ is centered around~$N$ while $Q_i^\Pi$, $i = 2, \dots, b$, are not, is because the fraction of servers with exactly one task tends to one as $N$ grows large as we will see. In case of a finite buffer size $b < \infty$, when a task is discarded, we call it an \emph{overflow} event, and we denote by $L^\Pi(t)$ the total number of overflow events under policy~$\Pi$ up to time~$t$. 

The next theorem states our main result. In the rest of the chapter let $D$ be the set of all right continuous functions from $[0,\infty)$ to $\mathbbm{R}$ having left limits and let `$\dist$' denote  convergence in distribution.


\begin{theorem}
\label{th: main}

For any policy $\Pi \in \Pi^{(N)}(d)$,
if for $i=1,2,\ldots$, $X_i^\Pi(0) \dist X_i(0)$ in $\mathbbm{R}$
as $N \to \infty$ with $X_i(0)=0$ for $i\geq 3$, then the processes
$\{X_i^\Pi(t)\}_{t \geq 0} \dist \{X_i(t)\}_{t \geq 0}$
in~$D$, where $X_i(t) \equiv 0$ for $i \geq 3$ and $(X_1(t), X_2(t))$
are unique solutions in $D \times D$ of the stochastic
integral equations
\begin{equation}\label{eq: main theorem}
\begin{split}
X_1(t) &=
X_1(0) + \sqrt{2} W(t) - \beta t + \int_0^t (-X_1(s)+X_2(s)) ds - U_1(t), \\
X_2(t) &= X_2(0) + U_1(t) + \int_0^t(-X_2(s)) ds,
\end{split}
\end{equation}
where $W$ is a standard Brownian motion and $U_1$ is
the unique non-decreasing non-negative process in~$D$ satisfying
$\int_0^\infty \mathbbm{1}_{[X_1(t)<0]} dU_1(t) = 0$.
\end{theorem}

The above result is proved in~\cite{EG15} for the ordinary JSQ
policy.
Our contribution is to develop a stochastic ordering construction
and establish that, somewhat remarkably, the diffusion limit is
the same for any policy in $\Pi^{(N)}(d)$.
In particular, the JIQ and JSQ policies yield the same diffusion limit.

\begin{remark}\textnormal{
We note that as in~\cite{EG15} we assume the convergence of the
initial state, which implies that the process has to start from
a state in which the number of vacant servers as well as the number
of servers with two tasks scale with $\sqrt{N}$,
and the number of servers with three or more tasks is $o(\sqrt{N})$.}
\end{remark}


\section{Coupling and stochastic ordering}
\label{sec: coupling}

In this section we prove several stochastic comparison results for the
system occupancy state under various load balancing schemes for a fixed
number of queues~$N$ (and hence we shall often omit the superscript~$N$
in this section).
These stochastic ordering results will be leveraged in the next section
to prove the main result stated in Theorem~\ref{th: main}.

In order to bring out the full strength of the stochastic comparison results,
we will in fact consider a broader class of load balancing schemes
$$\Pi^{(N)} := \big\{\Pi(d_0, d_1, \ldots, d_{b-1}): d_0 = N, 1 \leq d_i \leq N,
1 \leq i \leq b-1, b \geq 2\big\},$$ 
and show that Theorem~\ref{th: main} actually
holds for this entire class of schemes.
In the scheme $\Pi(d_0, d_1, \ldots, d_{b-1})$, the dispatcher assigns
an incoming task to the server with the minimum queue length among $d_k$
(possibly depending on~$N$) servers selected uniformly at random when the
minimum queue length across the system is~$k$, $k = 0, 1, \ldots, b - 1$.
As before, $b$ represents the buffer size, and when a task is assigned
to a server with $b$~outstanding tasks, it is instantly discarded.

\subsection{Stack formation and deterministic ordering}
\label{subsec: det_ord}

Let us consider the servers arranged in non-decreasing order of their queue lengths. Each server along with its queue can be thought of as a stack of items. The ensemble of stacks then represent the empirical CDF of the queue length distribution, and the $i^{\mathrm{th}}$ horizontal bar corresponds to $Q_i^{\Pi}$ (for the concerned policy $\Pi$). The items are added to and removed from the various stacks according to some rule. Before proceeding to the coupling argument, we first state and prove a deterministic comparison result under the above setting.

Consider two ensembles $A$ and $B$ with the same total number of stacks. The stacks in ensemble $A$ have a maximum capacity of $b$ items and those in ensemble $B$ have a maximum capacity of $b'$ items with $b\leq b'$. For two such ensembles a step is said to follow $Rule(k,l,l_A,l_B)$ if either addition or removal of an item in both ensembles is done in that step as follows:
\begin{enumerate}[(i)]
\item Removal: An item is removed (if any) from the $k^{\mathrm{th}}$   stack from both ensembles or an item is removed from some stack in ensemble $A$ but no removal is done in ensemble $B$.
\item Addition: 
\begin{itemize}
\item[(ii.a)] System A: If the minimum stack height is less than $b-1$, then the item is added to the $l^{\mathrm{th}}$  stack. Else, the item is added to the $l_A^{\mathrm{th}}$   stack. If the item lands on a stack with height $b$, then it is dropped.
\item[(ii.b)] System B: If the minimum stack height is less than $b-1$, then the item is added to the $l^{\mathrm{th}}$  stack. Otherwise if the minimum stack height is precisely equal to $b-1$, the item is added to the $l_B^{\mathrm{th}}$  stack. When the minimum stack height in the system is at least $b$, the item can be sent to any stack. If the item lands on a stack with height $b'$, then it is dropped.
\end{itemize}
\end{enumerate}
Then we have the following result.
\begin{proposition}\label{prop: det_ord}
Consider two ensembles $A$ and $B$ as described above with the total number of stacks being $N$, stack capacities being $b$ and $b'$ respectively, with $b\leq b'$ and with $\mathbf{Q}^A\leq \mathbf{Q}^B$ component-wise \emph{i.e,} $Q^A_i\leq Q^B_i$ for all $i\geq 1$. The component-wise ordering is preserved if at any step $Rule(k,l,l_A,l_B)$ is followed with $l_A\geq l_B$ and either $l=1$ or $l\geq l_B$.
\end{proposition}
Before diving deeper into the proof of this proposition, let us discuss the high-level intuition behind it.
First observe that, if $\mathbf{Q}^A\leq \mathbf{Q}^B$, and an item is added (removed) to (from) the stack with the same index in both ensembles, then the component-wise ordering will be preserved. Hence, the preservation of ordering at the time of removal, and at the time of addition when, in both ensembles, the minimum stack height is less than $b-1$, is fairly straightforward.

Now, in other cases of addition, since in ensemble $A$ the stack capacity is $b\ (\leq b')$, if the minimum stack height in ensemble $B$ is at least $b$,  the ordering is preserved trivially. This leaves us with only the case when the minimum stack height in ensemble $B$ is precisely equal to $b-1$. 
In this case, when the minimum stack height in ensemble $A$ is also precisely equal to $b-1$, the preservation of the ordering follows from the assumption that $l_A\geq l_B$, which ensures that if in ensemble $A$, the item is added to some stack with $b-1$ items (and hence increases $Q^A_{b}$), then the same will be done in ensemble $B$ whenever $Q^A_b=Q^B_b$. 
Otherwise if the minimum stack height in ensemble $A$ is less than $b$, then assuming either $l=1$ (i.e.~the item will be sent to the minimum queue) or $l\geq l_B$ (i.e.~an increase in $Q^A_b$ implies an increase in $Q_b^B$) ensures the preservation of ordering.
\begin{proof}[Proof of Proposition~\ref{prop: det_ord}]
Suppose after following $Rule(k,l,l_A,l_B)$ the updated stack heights of ensemble $\Pi$ are denoted by $(\tilde{Q}_1^{\Pi},\tilde{Q}_2^{\Pi},\ldots)$, $\Pi= A,B$. We need to show $\tilde{Q}_i^A\leq\tilde{Q}_i^B$ for all $i\geq 1$.

For ensemble $\Pi$ let us define $I_{\Pi}(c):=\max\{i: Q^{\Pi}_i\geq N-c+1\}$, $c=1,\ldots,N$, $\Pi= A,B$. Define $I_{\Pi}(c)$ to be 0 if $Q_1^{\Pi}$ is (and hence all the $Q^{\Pi}_i$ values are) less than $N-c+1$. Note that $I_A(c)\leq I_B(c)$ for all $c= 1,2,\ldots N$ because of the initial ordering.

Now if the rule produces a removal of an item,  then the updated ensemble will have the values
\begin{equation}
\tilde{Q}^{\Pi}_i=
\begin{cases}
Q^{\Pi}_i-1, &\mbox{ for }i=I_{\Pi}(k),\\
Q^{\Pi}_i,&\mbox{ otherwise, }
\end{cases}
\end{equation}
if $I_{\Pi}(k)\geq 1$; otherwise all the $Q^{\Pi}_i$ values remain unchanged.
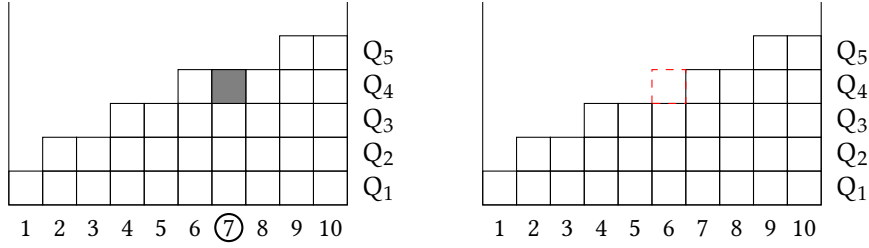
\begin{figure}
\begin{center}
\begin{tikzpicture}[scale=.45]
\draw (1,6)--(1,0)--(11,0)--(11,6);
\foreach \x in {10, 9,...,1}
	\draw (\x,1) rectangle (\x+1,0)
	(\x+.5,-.15) node [black,below] {\x} ;
\foreach \x in {10, 9,...,2}
	\draw (\x,2) rectangle (\x+1,1);
\foreach \x in {10, 9,...,4}
	\draw (\x,3) rectangle (\x+1,2);
\foreach \x in {10, 9,...,6}
	\draw (\x,4) rectangle (\x+1,3);
\foreach \x in {10, 9,...,9}
	\draw (\x,5) rectangle (\x+1,4);
\foreach \y in {1,2,...,5}
	\draw (11.15,\y-.5) node [black, right] {$Q_{\y}$};
\draw[thick] (7.5,-.67) circle [radius=.4];
\draw [fill=gray] (7,3) rectangle (8,4);

\draw (15,6)--(15,0)--(25,0)--(25,6);
\foreach \x in {10, 9,...,1}
	\draw (14+\x,1) rectangle (14+\x+1,0)
	(14+\x+.5,-.15) node [black,below] {\x} ;
\foreach \x in {24, 23,...,16}
	\draw (\x,2) rectangle (\x+1,1);
\foreach \x in {24, 23,...,18}
	\draw (\x,3) rectangle (\x+1,2);
\foreach \x in {24, 23,...,21}
	\draw (\x,4) rectangle (\x+1,3);
\foreach \x in {24, 23}
	\draw (\x,5) rectangle (\x+1,4);
\draw[dashed,thin, red] (20,3) rectangle (21,4);
\foreach \y in {1,2,...,5}
	\draw (25.15,\y-.5) node [black, right] {$Q_{\y}$};
\end{tikzpicture}
\caption{Removal of an item from the ensemble}\label{fig:removal}
\end{center}
\end{figure}
For example, in Figure~\ref{fig:removal}, $b=5$, $N=10$, and at the time of removal $k=7$. For this configuration $I_{\Pi}(7)=4$ since $Q^{\Pi}_4=5\geq 10-7+1=4$ but $Q^{\Pi}_5=2<4$. Hence, $Q_4^{\Pi}$ is reduced and all the other values remain unchanged. Note that the specific label of the servers does not matter here. So after the removal/addition of an item we consider the configuration as a whole by rearranging it again in non-decreasing order of the queue lengths.

Since in both $A$ and $B$ the values of $Q_i$ remain unchanged except for $i=I_A(k)$ and $I_B(k)$, it suffices to prove the preservation of the ordering for these two specific values of $i$. Now for $i=I_A(k)$, 
$$\tilde{Q}_i^A=Q_i^A-1\leq Q_i^B-1\leq\tilde{Q}_i^B.$$
If $I_B(k)=I_A(k)$, then we are done by the previous step. If $I_B(k)>I_A(k)$, then from the definition of $I_A(k)$ observe that $I_B(k)\notin\{i: Q_i^A\geq N-k+1\}$ and hence $Q_i^A<N-k+1$, for $i=I_B(k)$. Therefore, for $i=I_B(k)$,
$$\tilde{Q}_i^A\leq N-k\leq Q_i^B-1=\tilde{Q}_i^B.$$
On the other hand, if the rule produces the addition of an item to stack $l$, then the values will be updated as
\begin{equation}
\tilde{Q}^{\Pi}_i=
\begin{cases}
Q^{\Pi}_i+1, &\mbox{ for }i=I_{\Pi}(l)+1,\\
Q^{\Pi}_i,&\mbox{ otherwise, }
\end{cases}
\end{equation}
if $I_{\Pi}(l)<b_{\Pi}$, with $b_{\Pi}$ the stack-capacity of the corresponding system; otherwise the values remain unchanged.
\begin{figure}
\begin{center}
\begin{tikzpicture}[scale=.45]
\draw (1,6)--(1,0)--(11,0)--(11,6);
\foreach \x in {10, 9,...,1}
	\draw (\x,1) rectangle (\x+1,0)
	(\x+.5,-.15) node [black,below] {\x} ;
\foreach \x in {10, 9,...,2}
	\draw (\x,2) rectangle (\x+1,1);
\foreach \x in {10, 9,...,4}
	\draw (\x,3) rectangle (\x+1,2);
\foreach \x in {10, 9,...,6}
	\draw (\x,4) rectangle (\x+1,3);
\foreach \x in {10, 9,...,9}
	\draw (\x,5) rectangle (\x+1,4);
\foreach \y in {1,2,...,5}
	\draw (11.15,\y-.5) node [black, right] {$Q_{\y}$};
\draw[thick] (2.5,-.65) circle [radius=.4];
\draw [fill=gray] (2,2) rectangle (3,3);

\draw (15,6)--(15,0)--(25,0)--(25,6);
\foreach \x in {10, 9,...,1}
	\draw (14+\x,1) rectangle (14+\x+1,0)
	(14+\x+.5,-.15) node [black,below] {\x} ;
\foreach \x in {24, 23,...,16}
	\draw (\x,2) rectangle (\x+1,1);
\foreach \x in {24, 23,...,18}
	\draw (\x,3) rectangle (\x+1,2);
\foreach \x in {24, 23,...,21}
	\draw (\x,4) rectangle (\x+1,3);
\foreach \x in {24, 23}
	\draw (\x,5) rectangle (\x+1,4);
\draw[red,thin] (17,2) rectangle (18,3);
\foreach \y in {1,2,...,5}
	\draw (25.15,\y-.5) node [black, right] {$Q_{\y}$};
\end{tikzpicture}
\caption{Addition of an item to the ensemble}\label{fig:addition}
\end{center}
\end{figure}
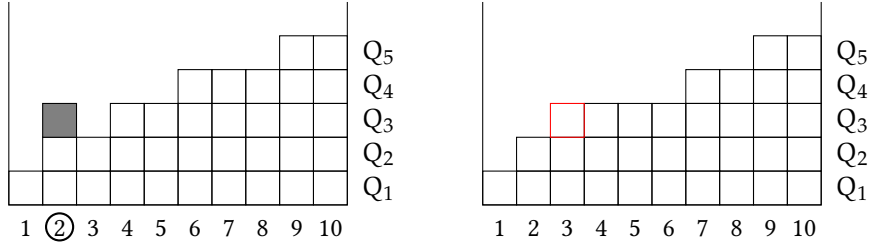
In Figure~\ref{fig:addition}, we have $l=2$ and for that particular configuration $I_{\Pi}(2)=2$. Hence, $Q_3^{\Pi}$ is incremented by one and the other variables remain fixed. 

Therefore, it is enough to consider the $i^{\mathrm{th}}$ horizontal bars for $i=(I_A(l)+1), (I_B(l)+1)$ when $I_A(l)<b$. According to the addition rule there are several cases which we now consider one by one:
\begin{enumerate}
\item First we consider the case when in both ensembles the minimum stack height is less than $b-1$. Then by part (ii) of the rule both incoming items are added to the $l^{\mathrm{th}}$   stack. When considering ensemble $B$ we may neglect the case $I_B(l)\geq b$ since then the value at $I_B(l)+1$ does not matter. Thus assume $I_B(l)\leq b-1$ and set $i=I_B(l)+1$ so that
$$\tilde{Q}_i^B=Q_i^B+1\geq Q_i^A+1\geq\tilde{Q}_i^A.$$
If $I_A(l)=I_B(l)$, then we are done by the previous case. If $I_A(l)+1\leq I_B(l)$, then it follows from the definition that $Q_i^A<N-l+1$ and $Q_i^B\geq N-l+1$, for $i=I_A(l)+1$. Hence,
$$\tilde{Q}_i^A=Q_i^A+1\leq N-l+1\leq Q_i^B\leq\tilde{Q}_i^B.$$
\item If the minimum stack height in $A$ is less than $b-1$ and that in $B$ is precisely $b-1$, then according to the rule the incoming item is added to the $l^{\mathrm{th}}$   stack in $A$ and the $l_B^{\mathrm{th}}$   stack in $B$. We here show that the component-wise ordering will be preserved if either $l=1$ or $l\geq l_B$. Observe that if $l=1$, then $I_A(l)<b-1$ which implies $I_A(l)+1\leq b-1$. But since the minimum stack height in $B$ is $b-1$, for all $i\leq b-1$ and in particular for $i=I_A(l)+1$, $\tilde{Q}^B_i=N\geq\tilde{Q}^A_i$. Now we consider the case when $l\geq l_B$. Also observe that the fact that the minimum stack height in $B$ is $b-1$, implies $I_B(l_B)\geq b-1\geq I_A(l_A)$ (since if $I_A(l)=b$, then nothing will be changed and so we do not need to consider this case). Then again if $I_A(l)=I_B(l_B)$, we are done.  Therefore, suppose $I_A(l)<I_B(l_B)$, which implies $I_A(l)+1\leq I_B(l_B)$. By definition, for $i=I_A(l)+1$, we have $Q_i^A<N-l+1$ and $Q_i^B\geq N-l_B+1\geq N-l+1$. Combining these two inequalities yields
$$\tilde{Q}_i^A=Q_i^A+1\leq N-l+1\leq Q_i^B=\tilde{Q}_i^B.$$
\item If the minimum stack height in both ensembles is $b-1$, then recall that the incoming item is added to the $l_A^{\mathrm{th}}$   stack in $A$ and to the $l_B^{\mathrm{th}}$   stack in $B$ with $l_A\geq l_B$. Arguing similarly as in the previous case we can conclude that the inequality is preserved.
\item Finally, if the minimum stack height in $B$ is larger than or equal to $b$, then the preservation of the inequality is trivial.
\end{enumerate}
Hence, the proof of the proposition is complete.
\end{proof}
\subsection{The coupling construction}\label{subsec: coupling}
We now construct a coupling between two systems $A$ and $B$ following any two schemes, say, $\Pi_A=\Pi(l_0,l_1,\ldots,l_{b-1})$ and $\Pi_B=\Pi(d_0, d_1,\ldots, d_{b'-1})$ in $\Pi^{(N)}$ respectively and combine it with Proposition~\ref{prop: det_ord} to get the desired stochastic ordering results.

For the arrival process we couple the two systems as follows. First we synchronize the  arrival epochs of the two systems. Now assume that in the systems $A$ and $B$, the minimum queue lengths are $k$ and $m$, respectively, $k\leq b-1$, $m\leq b'-1$. Therefore, when a task arrives, the dispatchers in $A$ and $B$ have to select $l_k$ and $d_m$ servers, respectively, and then have to send the task to the one having the minimum queue length among the respectively selected servers. Since the servers are being selected uniformly at random we can assume without loss of generality, as in the stack construction, that the servers are arranged in non-decreasing order of their queue lengths and are indexed in increasing order. 
Hence, observe that when a few server indices are selected, the server having the minimum of those indices will be the server with the minimum queue length among these. In this case the dispatchers in $A$ and $B$ select $l_k$ and $d_m$ random numbers (without replacement) from $\{1,2,\ldots,N\}$ and then send the incoming task to the servers having indices to be the minimum of those selected numbers. To couple the decisions of the two systems, at each arrival epoch a single random permutation of $\{1,2,\ldots,N\}$ is drawn, denoted by $\boldsymbol{\Sigma}^{(N)}:=(\sigma_1, \sigma_2,\ldots,\sigma_N)$. Define $\sigma_{(i)}:= \min_{j\leq i}\sigma_j$. Then observe that system $A$ sends the task to the server with the index $\sigma_{(l_k)}$ and system $B$ sends the task to the server with the index $\sigma_{(d_m)}$. Since at each arrival epoch both systems use a common random permutation, they take decisions in a coupled manner.

For the potential departure process, couple the service completion times of the $k^{\mathrm{th}}$ queue in both scenarios, $k= 1,2,\ldots,N$. More precisely, for the potential departure process assume that we have a single synchronized exp($N$) clock independent of arrival epochs for both systems. Now when this clock rings, a number $k$ is uniformly selected from $\{1,2,\ldots,N\}$ and a potential departure occurs from the $k^{\mathrm{th}}$ queue in both systems. If at a potential departure epoch an empty queue is selected, then we do nothing. In this way the two schemes, considered independently, still evolve according to their appropriate statistical laws.

Loosely speaking, our next result is based upon the following intuition: Suppose we have two systems $A$ and $B$ with two different schemes $\Pi_A$ and $\Pi_B$ having buffer sizes $b$ and $b'$ ($b\leq b'$) respectively. 
Also, for these two systems, initially, $Q^A_i\leq Q^B_i$ for all $i=1,\ldots,b$. Below we develop some intuition as to under what conditions the initial ordering of the $Q_i$-values will be preserved after one arrival or departure.

For the departure process if we ensure that departures will occur from the $k^{\mathrm{th}}$ largest queue in both systems for some $k\in\{1,2,\ldots,N\}$ (ties are broken in any way), then observe that the ordering will be preserved after one departure.

In case of the arrival process, assume that when the minimum queue length in both systems is less than $b-1$, the incoming task is sent to the server with the same index. In that case it can be seen that the $Q_i$-values in $A$ and $B$ will preserve their ordering after the arrival as well. Next consider the case when the minimum queue length in both systems is precisely $b-1$. Now, in $A$, an incoming task can either be rejected (and will not change the $Q$-values at all) or be accepted (and $Q^{\Pi_A}_b$ will increase by 1). Here we ensure that if the incoming task is accepted in $A$, then it is accepted in $B$ as well unless $Q_b^{\Pi_A}<Q_b^{\Pi_B}$, in which case it is clear that the initial ordering will be preserved after the arrival. Finally, if the minimum queue length in $A$ is less than $b-1$ and that in $B$ is precisely $b-1$, then the way to ensure the inequality is either by making the scheme $\Pi_A$ send the incoming task to the server with minimum queue length (and hence, it will only increase the value of $Q_i^{\Pi_A}$ for some $i<b$, leaving other values unchanged) or by letting the selected server in $\Pi_A$ have a smaller queue length than the selected server in $\Pi_B$. The former case corresponds to the condition $d=N$ and the latter corresponds to the condition $d\leq d_{b-1}$, either of which has to be satisfied, in order to ensure the preservation of the ordering. This whole idea is formalized below.

\begin{proposition}\label{prop: stoch_ord-jap}
For two schemes 
$$\Pi_A=\Pi(l_0,l_1,\ldots,l_{b-1})\quad \mbox{and} \quad\Pi_B=\Pi(d_0, d_1,\ldots, d_{b'-1})$$ 
with $b\leq b'$ assume $l_0=\ldots=l_{b-2}=d_0=\ldots=d_{b-2}=d$, $l_{b-1}\leq d_{b-1}$ and either $d=N$ or $d\leq d_{b-1}$. Then the following holds:
\begin{enumerate}[{\normalfont (i)}]
\item\label{component_ordering-jap} $\{Q^{ \Pi_A}_i(t)\}_{t\geq 0}\leq_{st}\{Q^{ \Pi_B}_i(t)\}_{t\geq 0}$ for $i=1,2,\ldots,b$,
\item\label{upper bound-jap} $\{\sum_{i=1}^b Q^{ \Pi_A}_i(t)+L^{ \Pi_A}(t)\}_{t\geq 0}\geq_{st} \{\sum_{i=1}^{b'} Q^{ \Pi_B}_i(t)+L^{ \Pi_B}(t)\}_{t\geq 0}$,
\item\label{delta_ineq-jap} $\{\Delta(t)\}_{t\geq 0}\geq \{\sum_{i=b+1}^{b'}Q_i^{ \Pi_B}(t)\}_{t\geq 0}$ almost surely under the coupling defined above,
\end{enumerate}
for any fixed $N\in\mathbbm{N}$ where $\Delta(t):=L^{ \Pi_A}(t)-L^{ \Pi_B}(t)$, provided that at time $t=0$ the above ordering holds.
\end{proposition}
\begin{proof}
To prove the stochastic ordering we use the coupling of the schemes as described above and show that the ordering holds for the entire sample path. That is, the two processes arising from the above pair of schemes will be defined on a common probability space and it will then be shown that the ordering is maintained almost surely over all time.

Note that we shall consider only the event times $0=t_0<t_1<\ldots$, i.e.~the time epochs when arrivals or potential service completions occur and apply forward induction to show that the ordering is preserved. By assumption the orderings hold at time $t_0 = 0$.\\

(i) The main idea of the proof is to use the coupling and show that at each event time the joint process of the two schemes follows a rule $Rule(k,l,l_A,l_B)$ described in Subsection~\ref{subsec: det_ord}, with some random $k$, $l$, $l_A$ and $l_B$ such that $l_A\geq l_B$ and either $l=1$ or $l\geq l_B$, and apply Proposition~\ref{prop: det_ord}. 
We now identify the rule at event time $t_1$ and verify that the conditions of Proposition~\ref{prop: det_ord} hold. If the event time $t_1$ is a potential departure epoch, then according to the coupling similarly as in the stack formation a random $k\in\{1, 2,\ldots, N\}$ will be chosen in both systems for a potential departure. Now assume that $t_1$ is an arrival epoch. 
In that case if the minimum queue length in both systems is less than $b-1$, then both schemes $ \Pi_A$ and $ \Pi_B$ will send the arriving task to the $\sigma_{(d)}^{\mathrm{th}}$   queue. If the minimum queue length in scheme $ \Pi_A$ is $b-1$, then the incoming task is sent to the $\sigma_{(l_{b-1})}^{\mathrm{th}}$ queue and if in scheme $ \Pi_B$ the minimum queue length is $b-1$, then the incoming task is sent to $\sigma_{(d_{b-1})}^{\mathrm{th}}$ queue where we recall that  $(\sigma_1, \sigma_2,\ldots,\sigma_N)$ is a random permutation of $\{1,2,\ldots,N\}$. Therefore, observe that at each step $Rule(\sigma_{(d)},k,\sigma_{(l_{b-1})}, \sigma_{(d_{b-1})})$ is followed.

Now to check the conditions, first observe that 
$$\sigma_{(l_{b-1})}= \min_{i\leq l_{b-1}}\sigma_i\geq\min_{i\leq d_{b-1}}\sigma_i=\sigma_{(d_{b-1})},$$
where the second inequality is due to the assumption $l_{b-1}\leq d_{b-1}$. In addition, we have assumed either $d=N$ or $d\leq d_{b-1}$. If $d=N$, then the dispatcher sends the incoming task to the server with the minimum queue length which is the same as sending to stack 1 as in Proposition $\ref{prop: det_ord}$. On the other hand, $d\leq d_{b-1}$ implies
$$\sigma_{(d)}= \min_{i\leq d}\sigma_i\geq \min_{i\leq d_{b-1}}\sigma_i=\sigma_{(d_{b-1})}.$$
Therefore, assertion~\eqref{component_ordering-jap} follows from Proposition~\ref{prop: det_ord}.\\

(ii) We again apply forward induction. Assume that the ordering holds at time $t_0$. If the next event time is an arrival epoch, then observe that both sides of the inequality in \eqref{upper bound-jap} will increase, since if the incoming task is accepted, then the $Q$-values will increase and if it is rejected, then the $L$-value will increase.\\
On the other hand, if the next event time is a potential departure epoch, then it suffices to show that,  if the left-hand-side decreases, then the right-hand-side decreases as well. Indeed, from assertion~\eqref{component_ordering-jap} we know that $Q^{ \Pi_A}_1\leq Q^{ \Pi_B}_1$ and hence we can see that if there is a departure from $ \Pi_A$ (i.e.~the $k^{\mathrm{th}}$ queue of $\Pi_A$ is non-empty), then there will be a departure from $ \Pi_B$ (i.e.~the $k^{\mathrm{th}}$ queue of $\Pi_B$ will be non-empty) as well.\\

(iii) Assertion~\eqref{delta_ineq-jap} follows directly from~\eqref{component_ordering-jap} and  \eqref{upper bound-jap}.
\end{proof}

\subsection{Discussion}
It is worth emphasizing that Proposition~\ref{prop: stoch_ord-jap}\eqref{component_ordering-jap} is fundamentally different from the stochastic majorization results
for the ordinary JSQ policy, and below we contrast our methodology with some existing literature.
As noted earlier, the ensemble of stacks, arranged in non-decreasing
order, represents the empirical CDF of the queue length distribution
at the various servers.
Specifically, if we randomly select one of the servers, then the
probability that the queue length at that server is greater than
or equal to~$i$ at time~$t$ under policy~$\Pi$ equals
$\frac{1}{N} \mathbbm{E} Q_i^\Pi(t)$.
Thus assertion~\eqref{component_ordering-jap} of Proposition~\ref{prop: stoch_ord-jap} implies that if we select one
of the servers at random, then its queue length is stochastically
larger under policy~$\Pi_B$ than under policy~$\Pi_A$.

The latter property does generally \emph{not} hold when we compare
the ordinary JSQ policy with an alternative load balancing policy.
Indeed,  the class of load balancing schemes $\tilde{\Pi}^{(N)}$ (for the $N^{\mathrm{th}}$ system say) considered in \cite{Towsley1992} consists of all the schemes that have instantaneous queue length information of all the servers and that have to send an incoming task to some server if there is at least some place available anywhere in the whole system. This means that a scheme can only discard an incoming task if the system is completely full. Observe that \emph{only} the JSQ policy lies both in the class $\Pi^{(N)}$ (defined in Section~\ref{sec: coupling}) and the class $\tilde{\Pi}^{(N)}$, because any scheme in $\Pi^{(N)}$ other than JSQ may reject an incoming task in some situations, where there might be some place available in the system. In this setup \cite{Towsley1992} shows that for any scheme $\Pi\in\tilde{\Pi}^{(N)}$, and for all $t\geq 0$,
\begin{align}\label{eq: towsley-jap}
\sum_{i=1}^k Y_{(i)}^{\JSQ}(t)&\leq_{st}\sum_{i=1}^k Y_{(i)}^{\Pi}(t),\mbox{ for } k=1,2,\ldots, N,\\
\{L^{\JSQ}(t)\}_{t\geq 0}&\leq_{st}\{L^{\Pi}(t)\}_{t\geq 0},
\end{align}
where $Y^{\Pi}_{(i)}(t)$ is the $i^{\mathrm{th}}$ largest queue length at time $t$ in the system following scheme $\Pi$ and $L^{\Pi}(t)$ is the total number of overflow events under policy $\Pi$ up to time $t$, as defined in Section~\ref{sec: model descr-jap}. Observe that $Y_{(i)}^{\Pi}$ can be visualized as the $i^{\mathrm{th}}$ largest vertical bar (or stack) as described in Subsection~\ref{subsec: det_ord}. Thus~\eqref{eq: towsley-jap} says that the sum of the lengths of the $k$ largest \emph{vertical} stacks in a system following any scheme $\Pi\in\tilde{\Pi}^{(N)}$ is stochastically larger than or equal to that following the scheme JSQ for any $k=1,2,\ldots,N$. Mathematically, this ordering can be written as
$$\sum_{i = 1}^{b} \min\{k, Q_i^{\JSQ}(t)\}  \leq_{st}
\sum_{i = 1}^{b} \min\{k, Q_i^{\Pi}(t)\},$$
for all $k = 1, \dots, N$.
In contrast, Proposition~\ref{prop: stoch_ord-jap} shows that the length of the $i^{\mathrm{th}}$ largest \emph{horizontal} bar in the system following some scheme $\Pi_A$ is stochastically smaller than that following some other scheme $\Pi_B$ if some conditions are satisfied. Also observe that the ordering between each of the horizontal bars (i.e.~$Q_i$'s) implies the ordering between the sums of the $k$ largest vertical stacks, but not the other way around. Further it should be stressed that, in crude terms, JSQ in our class $\Pi^{(N)}$, plays the role of upper bound, whereas what Equation~\eqref{eq: towsley-jap} implies is almost the opposite in nature to the conditions we require.

While in \cite{Towsley1992} no policies with admission control (where the dispatcher can discard an incoming task even if the system is not full) were considered, in a later paper \cite{towsley} and also in \cite{Towsley95} the class was extended to a class $\hat{\Pi}^{(N)}$ consisting of all the policies that have information about instantaneous queue lengths available and that can either send an incoming task to some server with available space or can reject an incoming task even if the system is not full. One can see that $\hat{\Pi}^{(N)}$ contains both $\tilde{\Pi}^{(N)}$ and $\Pi^{(N)}$ as subclasses. But then for such a class with admission control, \cite{towsley} notes that a stochastic ordering result like \eqref{eq: towsley-jap} cannot possibly hold. Instead, what was shown in \cite{Towsley95} is that for all $t\geq 0$,
\begin{align}\label{eq:ordering total jobs}
\sum_{i=1}^{k}Y_{(i)}^{\JSQ}(t)+L^{\JSQ}(t)\leq_{st}\sum_{i=1}^k Y_{(i)}^{\Pi}(t)+L^{\Pi}(t)\mbox{ for all }k\in\{1,2,\ldots,N\}
\end{align}
Note that the ordering in \eqref{eq:ordering total jobs} is the same in spirit as the ordering stated in Proposition~\ref{prop: stoch_ord-jap}\eqref{upper bound-jap} and the inequalities in \eqref{eq:ordering total jobs} are what in the language of \cite[Def.~14.4]{Towsley95} known as the \emph{weak sub-majorization by $p$}, where $p=L^{\Pi}(t)-L^{\JSQ}(t)$. 
But in this case also our inequalities in Proposition~\ref{prop: stoch_ord-jap}\eqref{component_ordering-jap} imply something completely orthogonal to what is implied by \eqref{eq:ordering total jobs}. 
In other words, the stochastic ordering results in Proposition~\ref{prop: stoch_ord-jap} provide both upper and lower bounds for the occupancy state of one scheme w.r.t.~another and 
are stronger than the stochastic majorization properties for the JSQ
policy existing in the literature. Hence we also needed to exploit a different proof methodology
than the majorization framework developed in~\cite{towsley, Towsley95, Towsley1992}.

\section{Convergence on diffusion scale}\label{sec: conv}
In this section we leverage the stochastic ordering established in Proposition~\ref{prop: stoch_ord-jap} to prove the main result stated in Theorem~\ref{th: main}. All the inequalities below are stated as  almost sure statements with respect to the common probability space constructed under the associated coupling. We shall use this joint probability space to make the probability statements about the marginals.
\begin{proof}[Proof of Theorem~\ref{th: main}]
Let $\Pi=\Pi(N,d_1,\ldots,d_{b-1})$ be a load balancing scheme in the class $\Pi^{(N)}$. Denote by $\Pi_1$ the scheme $\Pi(N,d_1)$ with buffer size $b=2$ and let $\Pi_2$ denote the JIQ policy $\Pi(N,1)$ with buffer size $b=2$.

Observe that from Proposition~\ref{prop: stoch_ord-jap} we have under the coupling defined in Subsection~\ref{subsec: coupling},
\begin{equation}\label{eq: bound-jap}
\begin{split}
|Q^{\Pi}_i(t)-Q^{\Pi_2}_i(t)|&\leq |Q^{\Pi}_i(t)-Q^{\Pi_1}_i(t)|+|Q^{\Pi_1}_i(t)-Q^{\Pi_2}_i(t)|\\
&\leq |L^{\Pi_1}(t)-L^{\Pi}(t)|+|L^{\Pi_2}(t)-L^{\Pi_1}(t)|\\
&\leq 2L^{\Pi_2}(t),
\end{split}
\end{equation}
for all $i\geq 1$ and $t\geq 0$ with the understanding that $Q_j(t)=0$ for all $j>b$, for a scheme with buffer $b$. The third inequality above is due to Proposition~\ref{prop: stoch_ord-jap}\eqref{delta_ineq-jap}, which in particular says that $\{L^{\Pi_2}(t)\}_{t\geq 0}\geq \{L^{\Pi_1}(t)\}_{t\geq 0}\geq\{ L^{\Pi}(t)\}_{t\geq 0}$ almost surely under the coupling. Now we have the following lemma which we will prove below.

\begin{lemma}\label{lem: tight-jap}
For all $t\geq 0$, under the assumption of Theorem~\ref{th: main}, $\{L^{\Pi_2}(t)\}_{N\geq 1}$ forms a tight sequence.
\end{lemma}

Since $L^{\Pi_2}(t)$ is non-decreasing in $t$, the above lemma in particular implies that
\begin{equation}\label{eq: conv 0-jap}
\sup_{t\in[0,T]}\frac{L^{\Pi_2}(t)}{\sqrt{N}}\pto 0.
\end{equation}
For any scheme $\Pi\in\Pi^{(N)}$, from \eqref{eq: bound-jap} we know that 
$$\{Q^{\Pi_2}_i(t)-2L^{\Pi_2}(t)\}_{t\geq 0}\leq\{Q^{\Pi}_i(t)\}_{t\geq 0}\leq\{Q^{\Pi_2}_i(t)+2L^{\Pi_2}(t)\}_{t\geq 0}.$$
Combining \eqref{eq: bound-jap} and \eqref{eq: conv 0-jap} shows that if the weak limits under the $\sqrt{N}$ scaling exist with respect to the Skorohod $J_1$-topology, they must be the same for all the schemes in the class $\Pi^{(N)}$. Also from Theorem 2 in \cite{EG15} we know that the weak limit for $\Pi(N,N)$ exists and the common weak limit for the first two components can be described by the unique solution in $D\times D$ of the stochastic differential equations in \eqref{eq: main theorem}. Hence the proof of Theorem~\ref{th: main} is complete.
\end{proof}

\begin{proof}[Proof of Lemma~\ref{lem: tight-jap}]
First we consider the evolution of $L^{\Pi_2}(t)$ as the following unit jump counting process. A task arrival occurs at rate $\lambda(N)$ at the dispatcher, and if $Q_1^{\Pi_1}=N$, then it sends it to a server chosen uniformly at random. If the chosen server has queue length 2, then $L^{\Pi_2}$ is increased by 1. It is easy to observe that this evolution can be equivalently described as follows. If $Q^{\Pi_2}_1(t)=N$, then each of the servers having queue length 2 starts increasing $L^{\Pi_2}$ by 1 at rate $\lambda(N)/N$. From this description we have
\begin{equation}\label{eq:L_rep}
L^{\Pi_2}(t)=A\left(\int_0^t\frac{\lambda(N)}{N}Q^{\Pi_2}_2(s)\mathbbm{1}[Q^{\Pi_2}_1(s)=N]ds\right)
\end{equation}
with $A(\cdot)$ being a unit rate Poisson process. Now using Proposition~\ref{prop: stoch_ord-jap} it follows that $\mathbbm{1}[Q^{\Pi_2}_1(s)=N]\leq \mathbbm{1}[Q^{\Pi_3}_1(s)=N]$ and $Q^{\Pi_2}_2(s)\leq Q^{\Pi_3}_2(s)$ where $\Pi_3=\Pi(N,N)$. Therefore, it is enough to prove the stochastic boundedness \cite[Def.~5.4]{PTRW07} of the sequence 
\begin{equation}
\Gamma^{(N)}(t):=A\left(\int_0^t\frac{\lambda(N)}{N}Q^{\Pi_3}_2(s)\mathbbm{1}[Q^{\Pi_3}_1(s)=N]ds\right).
\end{equation}
To prove this we shall use the martingale techniques described for instance in \cite{PTRW07}. 
Define the filtration $\mathbf{F}\equiv\{\mathcal{F}_t:t\geq 0\}$, where  for $t\geq 0$,
\begin{align*}
\mathcal{F}_t&:=\sigma\Big(Q^{\Pi_3}(0),A\left(\int_0^t\frac{\lambda(N)}{N}Q^{\Pi_3}_2(s)\mathbbm{1}[Q^{\Pi_3}_1(s)=N]ds\right),\\
&\hspace{5cm} Q^{\Pi_3}_1(s), Q^{\Pi_3}_2(s): 0\leq s\leq t\Big).
\end{align*}
Then using a random time change of unit rate Poisson process \cite[Lemma 3.2]{PTRW07} and similar arguments to those in \cite[Lemma 3.4]{PTRW07}, we have the next lemma.
\begin{lemma}
With respect to the filtration $\mathbf{F}$,
\begin{align*}
M^{(N)}(t)&:= A\left(\int_0^t\frac{\lambda(N)}{N}Q^{\Pi_3}_2(s)\mathbbm{1}[Q^{\Pi_3}_1(s)=N]ds\right)\\
&\hspace{5cm}-\int_0^t\frac{\lambda(N)}{N}Q^{\Pi_3}_2(s)\mathbbm{1}[Q^{\Pi_3}_1(s)=N]ds
\end{align*}
is a square-integrable martingale with $\mathbf{F}$-compensator $$I(t)=\int_0^t\frac{\lambda(N)}{N}Q^{\Pi_3}_2(s)\mathbbm{1}[Q^{\Pi_3}_1(s)=N]ds.$$ Moreover, the predictable quadratic variation process is given by $\langle M^{(N)}\rangle(t)=I(t).$
\end{lemma}
Now we apply Lemma 5.8 in \cite{PTRW07} which gives a stochastic boundedness criterion for square-integrable martingales. 
\begin{lemma}
\begin{normalfont}
\cite[Lemma 5.8]{PTRW07}
\end{normalfont}
Suppose that, for each $N\geq 1$, $M^{(N)}\equiv \{M^{(N)}(t):t\geq 0\}$ is a square-integrable martingale (with respect to a specified filtration) with predictable quadratic variation process $\langle M^{(N)}\rangle\equiv\{\langle M^{(N)}\rangle(t):t\geq 0\}$. If the sequence of random variables $\{\langle M^{(N)}\rangle(T): N\geq 1\}$ is stochastically bounded in $\mathbbm{R}$ for each $T>0$, then the sequence of stochastic processes $\{M^{(N)}: N\geq 1\}$ is stochastically bounded in $D$.
\end{lemma}
Therefore, it only remains to show the stochastic boundedness of $\{\langle M^{(N)}\rangle(T):N\geq 1\}$ for each $T>0$. Fix a $T>0$ and observe that
\begin{equation}\label{eq: qvp_boundedness}
\begin{split}
\langle M^{(N)}\rangle(T)&=\frac{\lambda(N)}{N}\int_0^T\frac{Q^{\Pi_3}_2(s)}{\sqrt{N}}\mathbbm{1}[Q^{\Pi_3}_1(s)=N]ds\\
&\leq \left[\sup_{t\in [0,T]}\frac{Q^{\Pi_3}_2(s)}{\sqrt{N}}\right]\times\left[\int_0^{T}\frac{1}{\sqrt{N}}\mathbbm{1}[Q^{\Pi_3}_1(s)=N]\lambda(N) ds\right].
\end{split}
\end{equation}
From \cite{EG15} we know that for any $T\geq 0$,  $\int_0^{T}1/\sqrt{N}\mathbbm{1}[Q^{\Pi_3}_1(s)=N]dA(\lambda(N)s)$  and $\sup_{t\in [0,T]}Q^{\Pi_3}_2(t)/\sqrt{N}$  are both tight. Moreover, since $\int_0^{T}1/\sqrt{N}\mathbbm{1}[Q^{\Pi_3}_1(s)=N]\lambda(N) ds$ is the intensity function of the stochastic integral $\int_0^{T}1/\sqrt{N}\mathbbm{1}[Q^{\Pi_3}_1(s)=N]dA(\lambda(N)s)$, which is a tight sequence, we have the following lemma.
\begin{lemma}
For all fixed $T\geq 0$, 
$\int_0^{T}\frac{1}{\sqrt{N}}\mathbbm{1}[Q^{\Pi_3}_1(s)=N]\lambda(N) ds$ is tight as a sequence in $N$.
\end{lemma}
Hence, both terms on the right-hand side of \eqref{eq: qvp_boundedness} are stochastically bounded and the resulting stochastic bound on $\langle M^{(N)}\rangle(T)$ completes the proof.
\end{proof}

\section{Conclusion}\label{sec:conclusion-jap}
In this chapter we have considered a system with symmetric
Markovian parallel queues and a single dispatcher.
We established the diffusion limit of the queue process in the
Halfin-Whitt regime for a wide class of load balancing schemes which
always assign an incoming task to an idle server, if there is any.
The results imply that assigning tasks to idle servers whenever
possible is sufficient to achieve diffusion level optimality.
Thus, using more fine-grained queue state information will increase the
communication burden and potentially impact the scalability in
large-scale deployments without significantly improving the performance.

In ongoing work we are aiming to extend the analysis to the stationary
distribution of the queue process, and in particular to quantify the
performance deviation from a system with a single centralized queue.
It would also be interesting to generalize the results to scenarios
where the individual nodes have general state-dependent service
rates rather than constant service rates.


%% file: jsqdiffusion.tex
\begin{abstract}
Consider a system of $N$~parallel single-server queues with unit-exponential service time distribution and a single dispatcher where tasks arrive as a Poisson process of rate $\lambda(N)$. When a task arrives, the dispatcher assigns it to one of the servers according to the Join-the-Shortest Queue (JSQ) policy. Eschenfeldt and Gamarnik~\cite{EG15} established that in the Halfin-Whitt regime where $(N - \lambda(N)) / \sqrt{N} \to \beta > 0$ as $N \to \infty$, an appropriately scaled occupancy measure of the system under the JSQ policy converges weakly on any finite time interval to a certain diffusion process as $N\to\infty$. Recently, it was further established by Braverman~\cite{Braverman18} that the convergence result extends to the steady state as well, i.e., the stationary occupancy measure of the system converges weakly to the steady state of the diffusion process as $N\to\infty$, proving the interchange of limits result.
In this chapter we perform a detailed analysis of the steady state of the above diffusion process, and obtain precise tail-asymptotics of the stationary distribution and scaling of extrema on large time intervals. 
\end{abstract}

\section{Introduction}
For any $\beta>0$, consider the following diffusion process
\begin{equation}\label{eq:diffusionjsq-aap}
\begin{split}
Q_1(t) &= Q_1(0) + \sqrt{2} W(t) - \beta t +
\int_0^t (- Q_1(s) + Q_2(s)) \dif s - L(t), \\
Q_2(t) &= Q_2(0) + L(t) - \int_0^t Q_2(s)\dif s.
\end{split}
\end{equation}
for $t \geq 0$, where $W$ is the standard Brownian motion, $L$ is
the unique nondecreasing nonnegative process in~$D_\R[0,\infty)$ satisfying
$\int_0^\infty \mathbbm{1}_{[Q_1(t) < 0]} \dif L(t) = 0$,
and $(Q_1(0), Q_2(0)) \in (-\infty, 0] \times [0, \infty)$.
In this chapter we establish tail asymptotics of the stationary distribution of the above diffusion process and identify the scaling behavior of $\inf_{0\leq s\leq t}Q_1(s)$ and $\sup_{0\leq s\leq t}Q_2(s)$ for large $t$.
Recall from Section~\ref{sec:powerofd} that under the Halfin-Whitt scaling for the arrival rate $\lambda(N)$ as in~\eqref{eq:HW}, the diffusion process in~\eqref{eq:diffusionjsq-aap} arises as the weak limit of the sequence of the scaled occupancy measure $\bar{\QQ}^N(t) = \left(\bar{Q}_1^N(t), \bar{Q}_2^N(t), \dots\right)$ (see~\eqref{eq:diffscale}) of systems under the Join-the-Shortest Queue (JSQ) policy, as the system size $N$ (number of servers in the system) becomes large~\cite{EG15}. 
Furthermore,  Braverman~\cite{Braverman18} recently established that the weak-convergence result extends to the steady state as well, i.e., $\bar{\QQ}^N(\infty)$ converges weakly to $(Q_1(\infty), Q_2(\infty), 0, 0,\ldots)$ as $N\to\infty$, where $(Q_1(\infty), Q_2(\infty))$ is distributed as the stationary distribution of the process $(Q_1, Q_2)$.
Thus, the steady state of the diffusion process in~\eqref{eq:diffusionjsq-aap} captures the asymptotic behavior of large-scale systems under the JSQ policy.
 
The steady state of the diffusion process in~\eqref{eq:diffusionjsq-aap} is technically hard to analyze.
In fact, even establishing its ergodicity is non-trivial. The standard method employed in studying steady-state behavior of diffusions \cite{ABD01, BL07,DW94, HM09} is to construct a suitable Lyapunov function which shows that the diffusion has a strong drift towards a compact set. Inside the compact set, some irreducibility condition, like uniform ellipticity (as in \cite{ABD01,BL07,DW94}) or hypoellipticity (as in \cite{HM09}), is used to show positive recurrence, and consequently, existence and uniqueness of the stationary distribution and ergodicity of the diffusion process. The construction of the Lyapunov function usually involves establishing stability of the associated noiseless dynamical system and having tractable bounds on hitting times for this deterministic system. In our setup, even the noiseless system requires non-trivial analysis (see Section 4.1 of \cite{Braverman18}).
In~\cite{Braverman18} a Lyapunov function is obtained via a generator expansion framework using Stein's method that establishes exponential ergodicity of $(Q_1,Q_2)$.
Although this approach gives a good handle on the rate of convergence to stationarity, the non-trivial dynamics of the noiseless system result in a complicated form for the Lyapunov function which sheds little light on the form of the stationary distribution itself. 
Moreover, the diffusion in~\eqref{eq:diffusionjsq-aap} (without the reflection term) is not hypoelliptic and this complicates things even further. It is also worth pointing out here that we obtain different tail behavior for $Q_1$ and $Q_2$ (Gaussian and Exponential, respectively) and get explicit dependence of $\beta$ in the exponents, which is hard to obtain using the Lyapunov function methods known in the literature.

This asks for a fundamentally different characterization of the stationary distribution, and we take resort to the theory of regenerative processes (see Chapter 10 of \cite{Thorisson}) to obtain a tractable representation of the steady state.
A variant of this method was first used in~\cite{BBD15} to study a diffusion process with inert drift, although the stationary distribution in that case had an explicit product form that facilitated the analysis, as opposed to the current scenario.
First, we show that the diffusion $(Q_1, Q_2)$ can be decomposed into i.i.d.~renewal cycles between carefully constructed regeneration times having good moment bounds. 
This decomposition gives an alternative, more transparent proof of ergodicity, and also shows that the diffusion falls in the category of classical regenerative processes. 
Loosely speaking, regeneration times are random times when the process \emph{starts afresh}, and the theory of classical regenerative processes can be used to conclude that the stationary behavior of a process is the same as the behavior within one renewal cycle (i.e., between two successive regeneration times).
The regenerative process representation enables us to obtain a form for the stationary distribution that is amenable to analysis (see Theorem \ref{th:stationary}).
Tail estimates for the stationary measure are then obtained by analyzing this form and are presented in Theorem \ref{th:statail}. 
Moreover, in Theorem~\ref{th:lil}, we obtain the precise almost sure
scaling behavior of the extrema of the process sample paths.

The regenerative structure of the diffusion process and the intermediate results might be of independent interest.
In fact, they might also be used to provide a detailed result for the behavior of the stationary measure near the center (bulk behavior) and produce sharp estimates on the stationary mean of $Q_2$. \\
 
The rest of the chapter is arranged as follows.
 In Section~\ref{sec:main}, we describe the two main results of this chapter. 
 In Section~\ref{sec:reg}, we establish $(Q_1, Q_2)$ as a classical regenerative process and state several crucial hitting-time estimates that are required to prove the main results.
 In Section~\ref{sec:renewal}, we obtain a tail estimate for the regeneration time which, in particular, implies that it has a finite first moment. This, in turn, implies the ergodicity of the diffusion process and gives a tractable form for the stationary distribution. 
 In Section~\ref{sec:fluc}, we obtain fluctuation estimates of the paths of $Q_1$ and $Q_2$ between two successive regeneration times, which are used in the proofs of Theorems~\ref{th:statail} and~\ref{th:lil}. 
 In Section~\ref{sec:proofmain}, we combine the results in Sections~\ref{sec:reg}, \ref{sec:renewal} and~\ref{sec:fluc} to prove Theorems~\ref{th:statail} and~\ref{th:lil}.

\section{Main results}\label{sec:main}

In this section we will state the main results, and discuss their ramifications. 
Recall the diffusion process $(Q_1, Q_2)$ as defined by Equation~\eqref{eq:diffusionjsq-aap}.
As mentioned in the introduction, it is known~\cite{Braverman18} that for any $\beta>0$, $(Q_1, Q_2)$ is an ergodic continuous-time Markov process.
Let $(Q_1(\infty), Q_2(\infty))$ denote a random variable distributed as the unique stationary distribution~$\pi$ of the process.
Then the next theorem gives a precise characterization of the tail of the stationary distribution.
\begin{theorem}\label{th:statail}
For any $\beta>0$ there exist positive constants $C_1, C_2, D_1, D_2$ \emph{not} depending on $\beta$ and positive constants $C^{l}(\beta), C^{u}(\beta), D^{l}(\beta), D^{u}(\beta), C_R(\beta), D_R(\beta)$ depending \emph{only} on $\beta$ such that 
\begin{equation}
\begin{split}
C^{l}(\beta)\e^{-C_1x^2} \le \pi(Q_1(\infty) < -x) \le C^{u}(\beta)\e^{-C_2x^2}, \ \ x \ge C_R(\beta)\\
D^{l}(\beta)\e^{-D_1\beta y} \le \pi(Q_2(\infty) > y) \le D^{u}(\beta)\e^{-D_2\beta y}, \ \ y \ge D_R(\beta).
\end{split}
\end{equation}
\end{theorem}
The dependence on $\beta$ of the tail-exponents is precisely captured in the above theorem.
Note that $Q_1(\infty)$ has a Gaussian tail, and the tail exponent is uniformly bounded by constants which do not depend on $\beta$, whereas $Q_2(\infty)$ has an exponentially decaying tail, and the coefficient in the exponent is linear in $\beta$.
\begin{remark}\label{rem:comp1}\normalfont
Let us now discuss a further implication of Theorem~\ref{th:statail}.
Recall that $Q_i^N(t)$ denotes the number of servers in the $N$-th system with queue length $i$ or larger at time $t$.
Let $S^N(t):=\sum_{i\geq 1}Q_i^N(t)$ denote the total number of tasks in the system.
Then \cite[Theorem 5]{Braverman18} implies that $(S^N(\infty)-N)/\sqrt{N}$ converges weakly to $S(\infty)\disteq Q_1(\infty)+Q_2(\infty)$.
In that case, Theorem~\ref{th:statail} implies that $S(\infty)$ has an exponential upper tail (large positive deviation) and a Gaussian lower tail (large negative deviation).
Although in terms of tail asymptotics, $S(\infty)$ behaves somewhat similarly to that for the centered and scaled total number of tasks in the corresponding M/M/$N$ system, there are some fundamental differences between the two processes that not only make the analysis of the JSQ policy much harder, but also lead to several completely different qualitative properties.
This has been discussed in Remark~\ref{rem:comp2-intro} in detail.
\end{remark}

The next theorem establishes the scaling behavior of the extrema of the process $\{(Q_1(t), Q_2(t))\}_{t \ge 0}$ on large time intervals.
\begin{theorem}\label{th:lil}
There exists a positive constant $\mathcal{C^*}$ not depending on $\beta$ such that the following hold almost surely along any sample path:
\begin{align*}
-2\sqrt{2} &\le \liminf_{t \rightarrow \infty} \frac{Q_1(t)}{\sqrt{\log t}} \le -1,\\
\frac{1}{\beta} &\le \limsup_{t \rightarrow \infty} \frac{Q_2(t)}{\log t} \le \frac{2}{\mathcal{C^*} \beta}.
\end{align*}
\end{theorem}
Again, Theorem~\ref{th:lil} captures the explicit dependence on $\beta$ of the width of the fluctuation window of $Q_1$ and $Q_2$. 
Specifically, note that the width of fluctuation of $Q_1$ does not depend on the value of $\beta$, whereas that of $Q_2$  is linear in $\beta^{-1}$. 
\begin{remark}\normalfont
Our proof of Theorem~\ref{th:statail} provides explicit values of the constants $C_1, C_2, D_1, D_2,\mathcal{C^*}$. 
We are not explicit about them in the statements of the theorems since these estimates are not sharp in the constants.
\end{remark}

\section{Regenerative process view of the diffusion}\label{sec:reg}
As mentioned in the introduction, the key challenge in analyzing the steady state of the diffusion process in~\eqref{eq:diffusionjsq-aap} stems from its lack of explicit characterization. 
In order to obtain sharp estimates for the stationary distribution we take resort to the theory of regenerative processes.
Loosely speaking, a stochastic process is called \emph{classical regenerative} if it starts anew at random times (called \emph{regeneration times}), independent of the past. See~\cite[Chapter 10]{Thorisson} for a rigorous treatment of regenerative processes.
The regeneration times split the process into renewal cycles that are independent and identically distributed, possibly except the first cycle.
Consequently, the behavior inside a specific renewal cycle characterizes the steady-state behavior.

In case of recurrent discrete state-space Markov chains regeneration times can be defined as hitting times of a fixed state.
Although the diffusion process in~\eqref{eq:diffusionjsq-aap} is two-dimensional, we will show that it actually \textit{exhibits point recurrence} and we can define regeneration times in terms of hitting times as follows.

First we introduce the following notations.
\begin{align*}
\tau_i(z)&:=\inf\{t \ge 0: Q_i(t)=z\}, \ i=1,2.\quad\text{and}\quad
\sigma(t):= \inf\{s \ge t: Q_1(s)=0\}.
\end{align*}
We now define the renewal cycles as follows.
Fix any $B>0$. For $k\geq 0$, define the stopping times
\begin{equation}\label{rendef}
\begin{split}
\alpha_{2k+1} &:= \inf \Big\{t\geq \alpha_{2k} : Q_2(t) = B\Big\},\\
\alpha_{2k+2} &:= \inf \left\{t>\alpha_{2k+1}:Q_2(t) = 2B\right\},\quad
\Xi_k := \alpha_{2k+2},
\end{split}
\end{equation}
with the convention that $\alpha_0=0$ and $\Xi_{-1}=0$.
The dependence of $B$ in the above stopping times is suppressed for convenience in notation. 
Hereafter we will assume $B>0$ to be fixed unless mentioned otherwise.
The next lemma describes the diffusion process as an appropriate classical regenerative process.
\begin{lemma}\label{lem:regeneration}
The process $\{Q_1(t), Q_2(t)\}_{t\geq 0}$ is a classical regenerative process with regeneration times given by $\{\Xi_k\}_{k\geq 0}$.
\end{lemma}
\begin{proof}
Note that it is enough to prove that 
$Q_1(\alpha_{2k}) = 0$ for all $k\geq 1$.
Indeed, this ensures that for all $k\geq 0$, $(Q_1(\Xi_k), Q_2(\Xi_k)) = (0, 2B)$, and the Markov process naturally regenerates at time $\Xi_k$.

Fix any $k\geq 1$.
Assume, if possible, $Q_1(\alpha_{2k}) < 0$. 
In that case, the path-continuity of $Q_1$ implies that the local time $L$ is constant in a small neighborhood of $\alpha_{2k}$.
Consequently, $Q_2$ must be strictly decreasing in an open time interval containing $\alpha_{2k}$.
This contradicts the fact that $\alpha_{2k}$ is the hitting time of a level \emph{from below} by the process $Q_2$.
\end{proof}
The above lemma implies that the regenerative cycles given by 
$$\big\{(Q_1(t), Q_2(t))\big\}_{\Xi_k\leq t<\Xi_{k+1}}$$ 
form an i.i.d.~sequence for $k\geq 0$.
The time intervals $\{\Xi_{k+1}-\Xi_k\}_{k\geq 0}$ are called the \emph{inter-regeneration times}.
In order to characterize the steady-state distribution using a regenerative approach, we first show that the initial \emph{delay length} $\Xi_0$ (time to enter into the regenerative cycles starting from an arbitrary state) as well as inter-regeneration times have finite expectations.
In fact, the next proposition establishes detailed tail asymptotics for the delay length $\Xi_0$ and thus, in particular, for the inter-regeneration times.
\begin{proposition}\label{prop:Xi}
Let $(Q_1(0), Q_2(0))= (x,y)$ with $x \le 0, y >0$. 
There exist constants $ c^{(1)}_\Xi, c^{(2)}_\Xi, t_{\Xi}>0$, possibly depending on $x,y,B,\beta$, such that for all $t\geq t_{\Xi}$,
$$\P_{(x,y)}(\Xi_0>t)\leq c^{(1)}_\Xi\exp(-c^{(2)}_\Xi t^{1/6}).$$
In particular, $\expt_{(x,y)}\Xi_0<\infty.$
\end{proposition}
Proposition~\ref{prop:Xi} is proved in Section~\ref{sec:renewal} and
yields the existence and uniqueness of the stationary distribution and ergodicity of the process as stated in Theorem~\ref{th:stationary} below. 
We note that the geometric ergodicity has already been proved in~\cite{Braverman18}.
The principal importance of Theorem~\ref{th:stationary} lies in the fact that it provides an explicit form of the stationary measure which will be the key vehicle in the study of the tail asymptotics and the fluctuation window, as stated in Theorems~\ref{th:statail} and~\ref{th:lil}. 
\begin{theorem}\label{th:stationary}
Fix any $B>0$. The process described by Equation~\eqref{eq:diffusionjsq-aap} has a unique stationary distribution $\pi$ which can be represented as
$$
\pi(A) = \dfrac{\mathbb{E}_{(0, 2B)}\left(\int_{0}^{\Xi_0}\mathbbm{1}_{[(Q_1(s),Q_2(s)) \in A]}ds\right)}{\mathbb{E}_{(0, 2B)}\left(\Xi_0\right)}
$$
for any measurable set $A \subseteq (-\infty, 0] \times (0, \infty)$. Moreover, the process is ergodic in the sense that for any measurable function $f$ satisfying 
$$\mathbb{E}_{(0, 2B)}\left(\int_{0}^{\Xi_0}f((Q_1(s),Q_2(s)))ds\right) < \infty,$$ the following holds:
\begin{equation}\label{ergodicity}
\frac{1}{t}\int\limits_0^t f((Q_1(s),Q_2(s)))ds \to \int\limits_{(-\infty, 0] \times (0, \infty)} f d\pi = \frac{\mathbb{E}_{(0, 2B)}\left(\int_{0}^{\Xi_0}f((Q_1(s),Q_2(s))ds\right)}{\mathbb{E}_{(0, 2B)}\left(\Xi_0\right)}
\end{equation}
almost surely as $t \rightarrow \infty$.
\end{theorem}
The above theorem follows using \cite[Chapter 10, Theorem 2.1]{Thorisson}, details of which are deferred till Section~\ref{sec:renewal}.
\begin{remark}\normalfont
We note that it can be shown by soft arguments involving Girsanov's theorem and the theory of L\'evy processes that the distribution of $\Xi_1 - \Xi_0$ has a density with respect to the Lebesgue measure, see the proof of Lemma 7.1 in \cite{BBD15}. 
This implies that the inter-regeneration time $\Xi_{k+1}-\Xi_k$ is \emph{spread-out} (see Section 3.5 of Chapter 10 in~\cite{Thorisson}).
Consequently, the total variation convergence of the diffusion process at time $t$ to the stationary distribution as $t\to\infty$, can be obtained using Theorem 3.3 of Chapter 10 in~\cite{Thorisson}. 
However, we skip this argument, since geometric ergodicity has already been established in \cite[Theorem 3]{Braverman18}. 
\end{remark}
In light of Theorem~\ref{th:stationary}, observe that establishing tail asymptotics of the stationary distribution reduces to studying the amount of time spent by the diffusion in a certain region in one particular renewal cycle.
The next theorem provides several important hitting-time estimates that will play a crucial role in the proofs of Theorems~\ref{th:statail} and~\ref{th:lil}.
Define
\begin{equation}\label{lzerodef}
l_0(\beta) := \max\left\lbrace \beta, \beta^{-1}, \frac{1}{\beta}\log\frac{1}{\beta}\right\rbrace.
\end{equation}
\begin{theorem}\label{th:excren}
There exists a positive constant $R_0$ such that with $B = R_0l_0(\beta)$ in~\eqref{rendef},
the following hold:
\begin{itemize}
\item[{\normalfont(i)}] There exist constants $C^*_1, C^*_2 >0$ that do not depend on $\beta$ such that for all $y \ge 4B$,
\begin{equation*}
\P_{(0, 2B)}\left(\tau_2(y) \le \Xi_0 \right) \le C^*_1 \e^{-C^*_2 \beta (y-\beta)/2}.
\end{equation*}
\item[{\normalfont(ii)}] 
For all $y \ge 2B$,
\begin{equation*}
\P_{(0,2B)}\left(\tau_2(y) \le \Xi_0\right) \ge (1-\e^{-\beta R_0l_0(\beta)})\e^{-\beta(y-2R_0l_0(\beta))}.
\end{equation*}
\item[{\normalfont(iii)}] There exists a constant $C^*(\beta)>0$ depending on $\beta$ such that for any $x \ge 18B$,
\begin{equation*}
\P_{(0,2B)}\left(\tau_1(-x) \le \Xi_0\right) \le C^*(\beta) \e^{-(x-2\beta)^2/8}.
\end{equation*}
\item[{\normalfont(iv)}] There exists a constant $C^{**}(\beta) >0$ depending on $\beta$ such that for any $x \ge \beta$,
\begin{equation*}
\P_{(0,2B)}\left(\inf_{t \le \Xi_0} Q_1(t) < -x\right) \ge C^{**}(\beta) \e^{-x^2}.
\end{equation*}
\end{itemize}
\end{theorem}
Theorem~\ref{th:excren} is proved in Section~\ref{sec:fluc} where we analyze the behavior of the process $(Q_1, Q_2)$ between two successive regeneration times. 
Results in Theorem~\ref{th:excren} in conjunction with Proposition~\ref{prop:Xi} and Theorem~\ref{th:stationary} are used to prove Theorems~\ref{th:statail} and~\ref{th:lil}, which is presented in Section~\ref{sec:proofmain}.

\section{Analysis of regeneration times}\label{sec:renewal}
In this section we will prove Proposition~\ref{prop:Xi} and Theorem~\ref{th:stationary}.
The proof of Proposition~\ref{prop:Xi} consists of several steps.
First, we analyze the down-crossings of $Q_2$, where 
we establish various hitting-time estimates in the time interval $[\alpha_{2k}, \alpha_{2k+1}]$, $k\geq 0$.
In particular, we prove the following lemma.
\begin{lemma}\label{lem:alphadowntail}
Fix $(Q_1(0), Q_2(0))=(x,y)$ with $x \le 0, y>0$. 
There exist positive constants $c_{\alpha_1}, c_{\alpha_1}', t_{\alpha_1}$ possibly depending on $(x,y)$, $B$, and $\beta$, such that for all $t \ge t_{\alpha_1}$,
$$
\P_{(x,y)}(\alpha_1 >t) \le c_{\alpha_1}'\exp(-c_{\alpha_1}t^{1/6}).
$$
\end{lemma}
As before, note that setting $(x,y) = (0,2B)$ furnishes the corresponding probabilities when $\alpha_1$ is replaced by $\alpha_{2k+1} - \alpha_{2k}$.
Lemma~\ref{lem:alphadowntail} is proved in Subsection~\ref{ssec:down}.
Next we consider the up-crossings of $Q_2$, where 
we establish various hitting-time estimates in the time interval $[\alpha_{2k+1}, \alpha_{2k+2}]$, $k \ge 0$. 
Specifically, we establish the following.
\begin{lemma}\label{lem:alphauptail}
Fix $(Q_1(0), Q_2(0))=(x,y)$ with $x \le 0, y>0$. 
There exist positive constants $c_{\alpha_2}, c_{\alpha_2}', t_{\alpha_2}$ possibly depending on $(x,y)$, $B$, and $\beta$, such that for all $t \ge t_{\alpha_2}$,
$$
\P_{(x,y)}(\alpha_2 - \alpha_1 >t) \le c_{\alpha_2}'\exp(-c_{\alpha_2}t^{1/6}).
$$
\end{lemma}
Lemma~\ref{lem:alphauptail} is proved in Section~\ref{ssec:up}.
Now observe that Lemmas~\ref{lem:alphadowntail} and~\ref{lem:alphauptail} together complete the proof of Proposition~\ref{prop:Xi}.\hfill $\qed$

\begin{proof}[Proof of Theorem~\ref{th:stationary}]
Due to Proposition~\ref{prop:Xi}, the fact that $\pi$ defined in the theorem is stationary follows from \cite[Chapter 10, Theorem 2.1]{Thorisson}. Now, we will prove the ergodicity result \eqref{ergodicity} which will also yield uniqueness. Take any starting point $(x,y)$ with $x \le 0$ and $y>0$ and recall $\Xi_{-1}=0$. Take any measurable function $f$ satisfying $\mathbb{E}_{(0, 2B)}\left(\int_{0}^{\Xi_0}f((Q_1(s),Q_2(s)))ds\right) < \infty$. Let $N_t = \sup\{k \ge -1: \Xi_{k} \le t\}$. Assume without loss of generality that $f$ is non-negative (for general $f$, consider the positive and negative parts of $f$ separately). We can write
\begin{align*}
&\int_0^{\Xi_0 \wedge t}f((Q_1(s),Q_2(s)))ds + \mathbbm{1}_{[\Xi_1 \le t]}\sum_{k=1}^{N_t}\int_{\Xi_{k-1}}^{\Xi_{k}} f((Q_1(s),Q_2(s)))ds \\
&\le \int_0^t f((Q_1(s),Q_2(s)))ds\\
&\le \int_0^{\Xi_0}f((Q_1(s),Q_2(s)))ds + \sum_{k=1}^{N_t+1}\int_{\Xi_{k-1}}^{\Xi_{k}} f((Q_1(s),Q_2(s)))ds.
\end{align*}
Clearly, $t^{-1}\int_0^{\Xi_0}f((Q_1(s),Q_2(s)))ds \rightarrow 0$ as $t \rightarrow \infty$. By Proposition 7.3 of \cite{Ross10},
$$
t^{-1} \sum_{k=1}^{N_t}\int_{\Xi_k}^{\Xi_{k+1}} f((Q_1(s),Q_2(s)))ds \rightarrow  \frac{\mathbb{E}_{(0, 2B)}\left(\int_{0}^{\Xi_0}f((Q_1(s),Q_2(s))ds\right)}{\mathbb{E}_{(0, 2B)}\left(\Xi_0\right)}
$$
and
 $$
t^{-1} \sum_{k=1}^{N_t+1}\int_{\Xi_k}^{\Xi_{k+1}} f((Q_1(s),Q_2(s)))ds \rightarrow  \frac{\mathbb{E}_{(0, 2B)}\left(\int_{0}^{\Xi_0}f((Q_1(s),Q_2(s))ds\right)}{\mathbb{E}_{(0, 2B)}\left(\Xi_0\right)}
$$
almost surely as $t \rightarrow \infty$. This proves \eqref{ergodicity}, and consequently uniqueness of the stationary distribution.
\end{proof}

\subsection{Down-crossings of \texorpdfstring{$\mathbf{Q_2}$}{Q2} and tightness estimates}\label{ssec:down}
In this subsection, we will prove tail asymptotics for the distribution of $\alpha_1$ as stated in Lemma~\ref{lem:alphadowntail}. 
This will require a crucial tightness estimate for the process $Q_2$, which is given in Lemma \ref{lem:q2regeneration} below. 
Loosely speaking, we need to have sharp estimates for the time  $Q_2$ takes to hit the level $B$ starting from a large initial state.
This, in turn, amounts to estimating the time integral of the $Q_1$ process when $Q_2$ is large, 
which is furnished by Lemma~\ref{lem:q1integral}.
The tail estimates presented in Lemmas~\ref{lem:integrated} and~\ref{lem:integrated2} will be used in the proof of Lemma~\ref{lem:q1integral}.\\

Fix any $M>0$ and $\varepsilon>0$. 
Observe that if $\inf_{0\leq s\leq t}Q_2(s)>M + \beta$, then the process $\{Q_1(s)\}_{0\leq s\leq t}$ is bounded below by the process $\{\eta(s)\}_{0\leq s\leq t}$, where 
$$\eta(t) = Q_1(0) + \sqrt{2}W(t) + Mt - L_\eta(t),$$ 
with $L_\eta$ being
the local time of $\eta$ given by $L_{\eta}(t) = \sup_{s \le t}\{Q_1(0) + \sqrt{2}W(s) + Ms\}^+$ (where $x^+ = \max\{x,0\}$ for any $x \in \mathbb{R}$), and $W$ being the standard Brownian motion.
Note that the dependence of $M$ in $\eta$ is suppressed for convenience in notation.
 For $i\geq 1$ define 
\begin{align*}
&T_{2i-1}:= \inf\ \{t>T_{2i-2}: \eta(t) = -\varepsilon\},
&T_{2i}:= \inf\ \{t>T_{2i-1}: \eta(t) = -\varepsilon/2\},\\
&\xi_i:= T_{2i}-T_{2i-1},\quad	\zeta_i := T_{2i+1}-T_{2i},
&u_i:= \sup_{T_{2i-1}\leq t\leq T_{2i}}(-\eta(t)),\\
&N_t = \inf\ \{n\geq 1: T_{2n}\geq t\}.
\end{align*}
with the convention that $T_{0}\equiv 0.$
Further, for $i\geq 1$, let $T_i^W$ denote the corresponding stopping
times when the process $\eta$ is replaced by the process $W_R$ described as 
$$W_R(t) = Q_1(0) + \sqrt{2}W(t) - L_W(t)$$
with $L_W$ being
the local time of $W_R$ given by $L_W(t) = \sup_{s \le t}\{Q_1(0) + \sqrt{2}W(s)\}^+$.
Also, similarly denote $\xi_i^W:= T_{2i}^W-T_{2i-1}^W$ and $\zeta_i^W := T_{2i+1}^W-T_{2i}^W$. 
\begin{lemma}\label{lem:integrated}
Assume that $Q_1(0) \in [-\varepsilon, 0]$.
Then the following hold:
\begin{enumerate}[{\normalfont (i)}]
\item For $i\geq 1$, $\zeta_i^W\leq_{st}\zeta_i$.
\item \label{subexp} There exist constants $c_\zeta, c_W>0$ not depending on $M, \varepsilon$ such that for $t \ge \varepsilon^2$
\begin{equation}
\begin{split}
&\mathrm{(a)}~\Pro{\zeta_1>t}\geq \exp(-c_\zeta t/\varepsilon^2)\\ 
&\mathrm{(b)}~\Pro{\zeta_1^W>t}\leq \exp(-c_Wt/\varepsilon^2).
\end{split}
\end{equation}
\item\label{lem:tails}
For all $x\geq \varepsilon$,
$\Pro{u_1>x}\leq \exp(-M(x-\varepsilon)),$
\item For all $t\geq \varepsilon/M$,
$\Pro{\xi_1>t}\leq \frac{2}{\sqrt{\pi}M\sqrt{t}}\exp(-M^2 t/16)$.
\item \label{lem:Nt-upper}
There exist constants $b,c_N^{(1)}>0$ not depending on $M, \varepsilon$, such that for $t \ge \varepsilon^2/b$
$$\Pro{N_t > b\varepsilon^{-2}t}\leq 2\exp(-c_N^{(1)}t/\varepsilon^2).$$ 
\end{enumerate}
\end{lemma}
\begin{proof}
(i) This is an immediate consequence of the fact that 
$$\{\eta(s)\}_{0\leq s\leq t}\geq_{st} \{W_R(s)\}_{0\leq s\leq t}.$$

\noindent
(ii)
Take $\varepsilon=1$. Using the Markov property for reflected Brownian motion, it is easy to see that there exist constants $c_\zeta, c_W>0$ such that $\exp(-c_Wt)\geq\Pro{\zeta_1^W>t}\geq \exp(-c_\zeta t)$ for $t \ge 1$.
(ii.a) now follows from (i) and Brownian scaling.
(ii.b) is also an immediate consequence of Brownian scaling. 

\noindent
(iii) Observe that
\begin{align*}
\Pro{u_1>x} 
\leq \Pro{\inf_{s<\infty}(-\varepsilon + \sqrt{2}W(s) + Ms)<-x} = \exp(-M(x-\varepsilon)),
\end{align*}
since $-\inf_{s < \infty}(\sqrt{2}W(s) + Ms)$ follows an exponential random variable with mean $1/M$.

\noindent
(iv) Note that
\begin{align*}
\Pro{\xi_1>t} &= \Pro{\sup_{s\leq t}(-\varepsilon+\sqrt{2}W(s) +Ms) \le -\varepsilon/2}\\
&\leq \Pro{\sqrt{2}W(t) + Mt \leq \varepsilon/2}
\leq \frac{2}{\sqrt{\pi}M\sqrt{t}}\exp(-M^2 t/16)\ \forall\ t\geq \varepsilon/M.
\end{align*}

\noindent
(v)
Observe that
\begin{align*}
\Pro{N_t>b\varepsilon^{-2}t} &\leq \Pro{\sum_{i=1}^{\lfloor b\varepsilon^{-2}t \rfloor}\zeta_i\leq t}\leq \Pro{\sum_{i=1}^{\lfloor b\varepsilon^{-2}t \rfloor}\zeta_i^W\leq t}, \ \text{ by part (i),}\\
&\leq \P\Big(\sum_{i=1}^{\lfloor b\varepsilon^{-2}t \rfloor}\frac{1}{\varepsilon^{2}}\left(\zeta_i^W - \expt(\zeta_i^W)\right)\leq - \left(\frac{b}{2\varepsilon^{2}}\expt \zeta_1^W - 1\right)\frac{t}{\varepsilon^{2}}\Big)\\
&\leq 2\exp(-c_N^{(1)}t/\varepsilon^2) \qquad\mbox{[choosing }b = 4\varepsilon^2/\expt(\zeta_1^W)\mbox{]},
\end{align*}
where the last step follows from part (ii), which shows that $\varepsilon^{-2}\left(\zeta_i^W - \expt(\zeta_i^W)\right)$ are sub-exponential random variables, and then using the Chernoff's inequality (see \cite[Pg.~16, Equation (2.2)]{Massart07}) to the sum $\sum_{i=1}^{\lfloor b\varepsilon^{-2}t \rfloor}\varepsilon^{-2}\left(\zeta_i^W - \expt(\zeta_i^W)\right)$.
Here, note that by Brownian scaling, $b$ chosen above does not depend on $\varepsilon$.
\end{proof}

The next technical lemma establishes a useful concentration inequality that will be crucial in obtaining tail probabilities for $\sum_{i=1}^{N_t}u_i\xi_i$.
\begin{lemma}\label{lem:BM-concentration}
Fix $\varepsilon>0$ and $M \ge \frac{1}{\varepsilon}$. 
Let $\Phi_i$'s be i.i.d.~nonnegative random variables with 
$$\Pro{\Phi_1>z}\leq \exp(-c'M^{3/2}\sqrt{z})\quad \text{for all}\quad z\geq 4\varepsilon^2/M,$$ 
and $\E{\Phi_1}\leq c_{11}\varepsilon^2/M$ where $c', c_{11}$ are positive constants not depending on $M, \varepsilon$. Then 
$$\P\Big(\sum_{i=1}^n\Phi_i\geq 4c_{11}n\frac{\varepsilon^2}{M}\Big)\leq \Big(1 + c_1\frac{1}{n^{2/5}\left(\varepsilon M\right)^{8/5}}\Big)\exp\Big(-c_2 (\varepsilon M)^{4/5}n^{1/5}\Big),$$
for $n \ge c_3\varepsilon M$, where $c_1,c_2, c_3$ are positive constants not depending on $M, \varepsilon$.
\end{lemma}
\begin{proof}
For some $A \ge 4\varepsilon^2/M$ to be chosen later, define 
$$\Phi_i^*:= \Phi_i\ind{\Phi_i\geq A} \qquad\mbox{and}\qquad \Phi_i^{**}:= \Phi_i\ind{\Phi_i < A}.$$
Thus, $\Phi_i = \Phi_i^* + \Phi_i^{**}$. Note that
\begin{align*}
\E{\Phi_i^*}^2 &= \int_{A^2}^{\infty}\P\left(\Phi_i>\sqrt{z}\right)dz = \int_{A}^{\infty}2z\P\left(\Phi_i>z\right)dz \\
&\le  \int_{A}^{\infty}2z\exp(-c'M^{3/2}\sqrt{z})dz
= \int_{\sqrt{A}}^{\infty}4z^3\exp(-c'M^{3/2}z)dz\\
& =  \frac{4}{M^6}\int_{M^{3/2}\sqrt{A}}^{\infty}z^3\exp(-c'z)dz
\leq c'' \frac{A^{3/2}}{M^{3/2}} \exp(-c'M^{3/2}\sqrt{A}),
\end{align*}
where the constant $c''$ does not depend on $M,A$.
Thus, using Chebyshev's inequality 
\begin{equation}\label{eq:*bound}
\Pro{\sum_{i=1}^n\Phi^*>2c_{11}n\frac{\varepsilon^2}{M}}\leq \frac{c'' M^{1/2} A^{3/2}\exp(-c'M^{3/2}\sqrt{A})}{4nc_{11}^2 \varepsilon^4}.
\end{equation}
Further note that $\Phi_i^{**}$'s are bounded random variables.
Therefore using the Azuma-Hoeffding inequality we obtain, 
\begin{equation}\label{eq:**bound}
\begin{split}
&\Pro{\sum_{i=1}^n\Phi_i^{**}>2c_{11}n\frac{\varepsilon^2}{M}}
= \Pro{\sum_{i=1}^n(\Phi_i^{**} - \E{\Phi_i^{**}})>c_{11}n\frac{\varepsilon^2}{M}}\\
&\hspace{1cm}\leq \exp\Big(-\big(\frac{c_{11}n\varepsilon^2}{M}\big)^2/(8A^2n)\Big)= \exp(-c_{11}^2n\varepsilon^4/(8A^2M^2)).
\end{split}
\end{equation}
Equating the exponents of equations~\eqref{eq:*bound} and~\eqref{eq:**bound}, and solving for $A$, we get
$$A = \left(\frac{c_{11}^2}{8c'}\right)^{2/5}\left(\frac{\varepsilon^{8/5}n^{2/5}}{M^{7/5}}\right).$$ The condition $A \ge 4\varepsilon^2/M$ implies $n \ge 2^5\left(\frac{8c'}{c_{11}^2}\right)\varepsilon M$. This choice for $A$ yields the bound claimed in the lemma. 
\end{proof}

\begin{lemma}\label{lem:integrated2}
Fix any $\varepsilon>0$ and $M \ge \frac{1}{\varepsilon}$.
\begin{enumerate}[{\normalfont (i)}]
\item \label{lem:u1xi1} There exist positive constants $c', c_{11}$ not depending on $M, \varepsilon$, such that
\begin{align*}
\mbox{\normalfont (a)}&\quad \Pro{u_1\xi_1>x} \le \exp(-c'M^{3/2}\sqrt{x})\qquad \forall \ x \ge 4\varepsilon^2/M,\\
\mbox{\normalfont (b)}&\quad \E{u_1\xi_1}\leq c_{11}\frac{\varepsilon^2}{M}.
\end{align*}
\item \label{fact:sum-uixi}
Let $b,c_{11}$ be the constants in Lemma~\ref{lem:integrated} \eqref{lem:Nt-upper} and Lemma~\ref{lem:integrated2} \eqref{lem:u1xi1} respectively. 
There exist constants $c_1, c_2, c_3$ not depending on $\varepsilon, M$ such that
$$
\Pro{\sum_{i=1}^{N_t}u_i\xi_i> 4 \frac{bc_{11}t}{M}} \leq c_1\exp(-c_2(\varepsilon M)^{4/5}(t/\varepsilon^2)^{1/5})
$$
for $t \ge c_3 \varepsilon^3M$.
\end{enumerate}
\end{lemma}
\begin{proof}
(i.a)~Recall that $M \ge \frac{1}{\varepsilon}$.
By Lemma~\ref{lem:integrated}~\eqref{lem:tails}, we obtain for $x \ge 4\varepsilon^2/M,$
\begin{align*}
\Pro{u_1\xi_1>x}&\leq \Pro{u_1>\sqrt{Mx}}+\Pro{u_1\xi_1>x, u_1\leq \sqrt{Mx}}\\
&\leq \Pro{u_1>\sqrt{Mx}}+\Pro{\xi_1>\frac{\sqrt{x}}{\sqrt{M}}}\\
&\leq \exp(-M(\sqrt{Mx}-\varepsilon))+\frac{2}{\sqrt{\pi}M^{3/4}x^{1/4}}\exp(-M^{3/2} \sqrt{x}/16)\\
&\leq \exp(-M^{3/2}\sqrt{x}/2)+\frac{2}{\sqrt{\pi}M^{3/4}x^{1/4}}\exp(-M^{3/2} \sqrt{x}/16)\\
&\leq \exp(-c'M^{3/2}\sqrt{x}),
\end{align*}
where the last line is a consequence of the fact that for $x \ge 4\varepsilon^2/M$ and $M \ge \frac{1}{\varepsilon}$, $M^{3/4}x^{1/4} \ge \sqrt{2M\varepsilon} > 1$.

\noindent
(i.b)~As a consequence of part (i.a) we obtain
\begin{align*}
\E{u_1\xi_1}&\leq \int_0^{4\varepsilon^2/M}\dif x+
+\frac{1}{M^3}\int_{4\varepsilon^2/M}^\infty \exp(-c'M^{3/2}\sqrt{x})M^3\dif x\\
&\leq \frac{4\varepsilon^2}{M} + \frac{c'''}{M^3} \leq c_{11}\frac{\varepsilon^2}{M},
\end{align*}
where we again used $M \ge \frac{1}{\varepsilon}$ to obtain $\frac{1}{M^3} \le \frac{\varepsilon^2}{M}$.

\noindent
(ii)
Observe that due to Lemma~\ref{lem:integrated} \eqref{lem:Nt-upper} and Lemma~\ref{lem:BM-concentration},
\begin{align*}
\P\Big(\sum_{i=1}^{N_t}u_i\xi_i>4 \frac{bc_{11}t}{M}\Big)
&\leq\Pro{N_t>b\varepsilon^{-2}t} + \P\ \Big(\sum_{i=1}^{\lfloor b\varepsilon^{-2}t\rfloor}u_i\xi_i>4 \frac{bc_{11}t}{M}\Big)\\
&\leq \exp(-c_N^{(1)}t/\varepsilon^2) + C_1\exp(-C_2(\varepsilon M)^{4/5}(t/\varepsilon^2)^{1/5})
\end{align*}
for $t \ge C_3 \varepsilon^3M$, where $C_1, C_2, C_3$ can be chosen to be independent of $M, \varepsilon$. This completes the proof.
\end{proof}

We are now in a position to state and prove Lemma~\ref{lem:q1integral} that provides us with a crucial estimate for the time-integral of the $Q_1$ process when $Q_2$ is large.
\begin{lemma}\label{lem:q1integral}
There exist $c'_1, c'_2, c'_3>0$, not depending on $\beta$ such that for any $y > c'_1\left(\beta \vee \beta^{-1}\right)+ \beta$,
\begin{align*}
&\P_{(0, y)}\Bigg(\int_0^t(-Q_1(s))\dif s>\left(\beta \wedge \beta^{-1}\right)\frac{t}{2},\ \inf_{s\leq t}Q_2(s) \geq c'_1\left(\beta \vee \beta^{-1}\right)+ \beta\Bigg)\\
&\leq \exp\Big(-c'_2t^{1/5}\left(\beta \vee \beta^{-1}\right)^{2/5}\Big)\qquad \text{for}\quad t \ge c'_3 \left(\beta \wedge \beta^{-1}\right)^2.
\end{align*} 
\end{lemma}
\begin{proof}
Recall the constants $b$ and $c_{11}$ from Lemma~\ref{lem:integrated}~\eqref{lem:Nt-upper} and Lemma~\ref{lem:integrated2}~\eqref{lem:u1xi1} respectively. 
As $c_{11}$ appears in the upper bound of $\mathbb{E}\left(u_1\xi_1\right)$ in Lemma~\ref{lem:integrated2}~\eqref{lem:u1xi1}, we can take $c_{11}> b^{-1} \vee 1$. First we consider the case $\beta \in (0,1)$. Take $\varepsilon = \beta/4$. Choose $M = 16c_{11}b/\beta$, since in that case 
$$\varepsilon = \frac{\beta}{4} = \frac{4c_{11}b}{M}.$$
Observe that 
\begin{align*}
&\P_{(0,y)}\Big(\int_0^t(-Q_1(s))\dif s>\frac{\beta t}{2},\ \inf_{s\leq t}Q_2(s) \geq M+\beta\Big)\\
&\leq\P_{(0,\ y)}\Big(\sum_{i=1}^{N_t}\int_{T_{2i-1}}^{T_{2i}}(-Q_1(s))\dif s > \frac{4c_{11}b}{M}t, \inf_{s\leq t}Q_2(s) \geq M+\beta\Big)\\
&\leq\Pro{\sum_{i=1}^{N_t}u_i\xi_i>\frac{4c_{11}b}{M}t}
\leq \exp\Big(-c''_2(\beta M)^{4/5}(t/\beta^2)^{1/5}\Big)\\
&\leq \exp\Big(-c'_2(t/\beta^2)^{1/5}\Big) \quad \text{for}\quad t \ge c''_3 \beta^3 M= c'_3\beta^2, \quad\mbox{due to Lemma~\ref{lem:integrated2}~\eqref{fact:sum-uixi}},
\end{align*}
where the constants $c'_2, c''_2 c'_3, c''_3$ do not depend on $\beta,M$. 
Next, for the case $\beta>1$, we take $\varepsilon = \frac{1}{4\beta}$ and $M= 16c_{11}\beta$ so that
$$
\varepsilon = \frac{1}{4\beta} = \frac{4c_{11}b}{M},
$$
and then apply the same argument. This completes the proof.
\end{proof}

\begin{lemma}\label{lem:q2regeneration}
There exist positive constants $c'_1, c'_2, c'_3, c'_4$ not depending on $\beta$ such that the following hold: 
\begin{enumerate}[{\normalfont (i)}]
\item For $\beta \ge 1$ and any $y \ge 1$, for all $ t\geq c'_4 y/\beta$
\begin{align*}
\P_{(0,\ y + c'_1\beta)}\big(\inf_{s\leq t}Q_2(s) > c'_1\beta \big)\leq c'_3\exp(-c'_2\beta^{2/5}t^{1/5}).
\end{align*}
\item For $\beta \in (0,1)$ and any $y \ge 1$, for all $t\geq c'_4 \big(y\beta^{-1} \vee \beta^{-2}\big)$
\begin{align*}
&\P_{(0,\ y + c'_1\beta^{-1})}\Big(\inf_{s\leq t}Q_2(s) > \frac{c'_1}{\beta}\Big)\\
&\hspace{1cm}\leq c'_3\left(\exp(-c'_2\beta^{-2/5}t^{1/5} + \exp(-c'_2\beta^2 t) + \beta^{-2}\exp(-c'_2 t)\right).
\end{align*}
\end{enumerate}
\end{lemma}
\begin{proof}
Let us denote the following events
\begin{align*}
\cE_t&:= \Big[\inf_{s\leq t} Q_2(s) > c'_1\left(\beta \vee \beta^{-1}\right)+\beta\Big],\\
\cE_t^1&:= \Bigg[\int_0^t(-Q_1(s))\dif s>\frac{\beta t}{2}, \quad\inf_{s\leq t} Q_2(s) > c'_1\left(\beta \vee \beta^{-1}\right)+ \beta\Bigg],\\
\cE_t^2&:=\Bigg[\int_0^t(-Q_1(s))\dif s\leq\frac{\beta t}{2},\quad \inf_{s\leq t} Q_2(s) > c'_1\left(\beta \vee \beta^{-1}\right)+ \beta\Bigg].
\end{align*}
Note that if $(Q_1(0), Q_2(0)) = (0,\ y + c'_1\left(\beta \vee \beta^{-1}\right)+ \beta)$, then from the evolution equation of the diffusion in~\eqref{eq:diffusionjsq-aap}, the event $\cE_t^2$ implies the event
\begin{align*}
\tilde{\cE}_t^2&:=\Bigg[Q_1(t) + Q_2(t) \leq y+ c'_1\left(\beta \vee \beta^{-1}\right) + \beta +\sqrt{2}W(t) - \frac{\beta t}{2}, \\ 
&\hspace{5cm}\inf_{s\leq t} Q_2(s) > c'_1\left(\beta \vee \beta^{-1}\right)+ \beta\Bigg].
\end{align*}
Therefore,
\begin{equation}\label{eq:lem3.7-1}
\begin{split}
\P_{(0,\ y + c'_1\left(\beta \vee \beta^{-1}\right)+ \beta)}\big(\cE_t\big) 
&\leq \P_{(0,\ y + c'_1\left(\beta \vee \beta^{-1}\right)+ \beta)}\Big(\cE_t^1\Big)\\
&\hspace{1cm}+\P_{(0,\ y + c'_1\left(\beta \vee \beta^{-1}\right)+ \beta)}\Big(\tilde{\cE}_t^2\Big).
\end{split}
\end{equation}
Now, choose $c'_1, c'_2$ as in Lemma \ref{lem:q1integral}. Then for any $y \ge 1$,
\begin{equation}\label{eq:lem3.7-2}
\P_{(0,\ y + c'_1\left(\beta \vee \beta^{-1}\right)+ \beta)}\Big(\cE_t^1\Big) \leq \exp(-c'_2t^{1/5}\left(\beta \vee \beta^{-1}\right)^{2/5}).
\end{equation}
Also, note that 
\begin{equation}\label{eq:lem3.7-3}
\begin{split}
&\P_{(0,\ y + c'_1\left(\beta \vee \beta^{-1}\right)+ \beta)}\Big(\tilde{\cE}_t^2\Big)\\
&\leq \P_{(0,\ y + c'_1\left(\beta \vee \beta^{-1}\right)+ \beta)}\big(Q_1(t) \leq y+\sqrt{2}W(t) - \frac{\beta t}{2},\\
&\hspace{5.5cm}  \inf_{s\leq t} Q_2(s) > c'_1\left(\beta \vee \beta^{-1}\right)+ \beta\big)\\
&\leq \Pro{\sqrt{2}W(t)>\frac{\beta t}{4}} + \P_{(0,\ y + c'_1\left(\beta \vee \beta^{-1}\right)+ \beta)}\Big(Q_1(t) \leq y - \frac{\beta t}{4}, \\
&\hspace{5.5cm} \inf_{s\leq t} Q_2(s) > c'_1\left(\beta \vee \beta^{-1}\right)+ \beta\Big).
\end{split}
\end{equation}
Due to Brownian scaling we have
\begin{equation}
\Pro{\sqrt{2}W(t)>\frac{\beta t}{4}}\leq c\exp(-c'\beta^2 t) \ \text{ for } \ t \ge \beta^{-2},
\end{equation}
where $c,c'$ do not depend on $\beta$.
Moreover, choosing $t > 8y/\beta$, and applying Lemma~\ref{lem:integrated}~\eqref{lem:tails} and Lemma~\ref{lem:integrated}~\eqref{lem:Nt-upper} with $\varepsilon = (\beta \wedge \beta^{-1})/4$ and $M=c'_1\left(\beta \vee \beta^{-1}\right)$,
\begin{equation}\label{eq:lem3.7-last}
\begin{split}
&\P_{(0,\ y + c'_1\left(\beta \vee \beta^{-1}\right)+ \beta)}\Big(Q_1(t) \leq y - \frac{\beta t}{4},  \inf_{s\leq t} Q_2(s) > c'_1\left(\beta \vee \beta^{-1}\right)+ \beta\Big)\\
&\leq \P_{(0,\ y + c'_1\left(\beta \vee \beta^{-1}\right)+ \beta)}\Big(Q_1(t) \leq - \frac{\beta t}{8},  \inf_{s\leq t} Q_2(s) > c'_1\left(\beta \vee \beta^{-1}\right)+ \beta\Big)\\
&\leq \P\Big(\sup_{1\leq i\leq N_t}u_i > \frac{\beta t}{8}\Big)\\
&\leq \Pro{N_t>16b\left(\beta \vee \beta^{-1}\right)^2t} + 16b\left(\beta \vee \beta^{-1}\right)^2t \Pro{u_1>\frac{\beta t}{8}}\\
&\leq \exp\Big(-c\left(\beta \vee \beta^{-1}\right)^2t\Big) + 16b\left(\beta \vee \beta^{-1}\right)^2t\\
&\hspace{5cm}\times \exp\left(-\left(\beta \vee \beta^{-1}\right)\left(\frac{\beta t}{8}-\frac{\beta \wedge \beta^{-1}}{4}\right)\right)\\
&\leq \exp\Big(-c\left(\beta \vee \beta^{-1}\right)^2t\Big) + 16b\left(\beta \vee \beta^{-1}\right)^2t\ \exp\left(-\left(\beta \vee \beta^{-1}\right)\left(\frac{\beta t}{16}\right)\right),
\end{split}
\end{equation}
where $b,c$ do not depend on $\beta$. 
Combining Equations~\eqref{eq:lem3.7-1} -- \eqref{eq:lem3.7-last} completes the proof of the lemma.
\end{proof}

We now have all the necessary results to prove Lemma~\ref{lem:alphadowntail}.
\begin{proof}[Proof of Lemma~\ref{lem:alphadowntail}]
From Lemma \ref{lem:q2regeneration}, for any $\beta>0$, we obtain $M^*>2B, t^*>0$ such that for all $t \ge t^*$,
\begin{align}\label{down}
\P_{(0,2M^*)}(\tau_2(M^*) >t) \le C_1\exp(-C_2t^{1/5}),
\end{align}
where the constants $C_1, C_2>0$ depend on $\beta,M^*$.  
Set the starting state to be $(Q_1(0), Q_2(0))=(0,y)$ where $M^* \ge y \ge B$. 
It will be clear from the proof that the same argument works for general starting points $(x,y)$ with $x \le 0, y >0$. For $k \ge 0$, define the following stopping times:
\begin{align*}
\alpha^*_{2k+1} & = \inf \Big\{t\geq \alpha^*_{2k} : Q_2(t) = 2M^* \text{ or } Q_2(t) = B\Big\},\\
\alpha^*_{2k+2} & = \inf \left\{t>\alpha^*_{2k+1}:Q_2(t) = M^* \text{ or } Q_2(t) = B\right\},
\end{align*}
where by convention, we take $\alpha^*_0=0$. 
Let $\mathcal{N}' := \inf\ \{k \ge 0: Q_2(\alpha^*_{2k}) = B\}$.

We will first prove the following: for some positive constant $p(M^*)$ that depends only on $M^*$,
\begin{equation}\label{alphadown}
\inf_{z \in [B, M^*]}\P_{(0,z)}(\tau_2(B) < \tau_2(2M^*))
 \ge p(M^*)>0.
\end{equation}
To see this, recall $S(t) = Q_1(t) + Q_2(t)$ and note that for $t \le \tau_1(-\beta/2)$,
$$
S(0) + \sqrt{2}W(t) - \beta t/2 \ge S(t) \ge Q_1(t).
$$
Further, note that $S(t) \le Q_2(t)$. 
Moreover, due to arguments similar to Lemma~\ref{lem:regeneration}, we know $Q_1(\tau_2(2M^*))=0$, and hence, $S(\tau_2(2M^*)) = Q_2(\tau_2(2M^*))$. 
Combining these facts, we obtain for any $z \in [B, M^*]$,
\begin{equation}\label{eq:lem3.1-p1m}
\begin{split}
\P_{(0,z)}(\tau_2(2M^*) &\le \tau_1(-\beta/2)) \\
&\le \P(S(t) \text{ hits } 2M^* \text{ before } -\beta/2)\\
&\le \P(S(0) + \sqrt{2}W(t) - \beta t/2 \text{ hits } 2M^* \text{ before } -\beta/2)\\
&\le \P(\sqrt{2}W(t) - \beta t/2 \text{ hits } M^* \text{ before } -(M^* + \beta/2))\\
&\le \e^{-\beta M^*/2} <1,
\end{split}
\end{equation}
where we have used the fact that the scale function (see \cite[V.46]{RW2000_2}) for $\sqrt{2}W(t) - \beta t/2$ is $s(x)= \exp(\beta x/2)$.

Now we will show that if the process $(Q_1, Q_2)$ starts with the initial state $(-\beta/2,z)$ with $z \le 2M^*$, 
then with positive probability $Q_1(t) <0$ for all $t \le \log (2M^* B^{-1})$.  
This in turn implies that $Q_2$ hits the level $B$ before time $\log (2M^* B^{-1})$, since for $t \le \tau_1(0)$, $(\dif/\dif t) Q_2(t) = - Q_2(t)$. 

Construct the Ornstein-Uhlenbeck process $Q_1^+$ on the same probability space as $Q_1$ as follow
$$
Q_1^+(t) = Q_1(0) + \sqrt{2}W(t) + \int_0^t(-Q_1^+(s) + (2M^*-\beta))ds,
$$
where the driving Brownian motion $W$ is the same as that for $Q_1$. 
By~\cite[Proposition 2.18]{Karatzas}, $Q_1(t) \le Q_1^+(t)$ for all $t \le \tau_1(0)$.
Now define the following event
$$
\cE(M^*) := \Big\{Q_1^+(t) < 0 \text{ for all } t \le \log (2M^* B^{-1})\Big\}.
$$
Note that $\cE(M^*)$ does not depend on $z$. It follows from the Doob representation for Ornstein-Uhlenbeck processes that $\P(\cE(M^*))>0$.
Thus,
\begin{equation}\label{timelength}
\begin{split}
&\inf_{z \in [B, 2M^*)}\P_{(-\beta/2,z)}\Big(\tau_2(B) \le  \log (2M^* B^{-1}) < \tau_2(2M^*)\Big)\\
&\hspace{.5cm}\ge \inf_{z \in [B, 2M^*)}\P_{(-\beta/2,z)}(\tau_1(0) \ge \log (2M^* B^{-1})) \ge \P(\mathcal{E}(M^*))>0.
\end{split}
\end{equation}
The strong Markov property in combination with \eqref{eq:lem3.1-p1m} and~\eqref{timelength} now produces the bound
\begin{align*}
\inf_{z \in [B, M^*]}\P_{(0,z)}(\tau_2(B) < \tau_2(2M^*))
 \ge (1-\e^{-\beta M^*/2})\P(\mathcal{E}(M^*))>0
\end{align*}
which proves \eqref{alphadown}. By virtue of \eqref{alphadown}, we have the following for $n \ge 1$,
\begin{align}\label{down1}
\P(\mathcal{N}'>n) \le (1-p(M^*))^n.
\end{align}
Now, let $T(M^*)$ be a number large enough such that
\begin{equation}\label{downadd}
\P\left(\sqrt{2} W(T(M^*)) \ge \beta T(M^*)/2 - (2M^* + \beta/2)\right) \le \P(\mathcal{E}(M^*))/2.
\end{equation}
Then,
\begin{equation}\label{downadd1}
\begin{split}
&\sup_{z \in [B, 2M^*)} \P_{(0,z)}\big(\tau_2(B) \wedge \tau_2(2M^*) > T(M^*) + \log (2M^* B^{-1})\big)\\
&\le \sup_{z \in [B, 2M^*)} \P_{(0,z)}\big(\tau_2(B) \wedge \tau_2(2M^*) > T(M^*)+ \log (2M^* B^{-1}),\\
&\hspace{1cm}\tau_1(-\beta/2) < T(M^*)\big)
+ \sup_{z \in [B, 2M^*)} \P_{(0,z)}(\tau_1(-\beta/2) \ge T(M^*))\\ 
&\le \sup_{z \in [B, 2M^*)}\P_{(-\beta/2,z)}(\tau_2(B) \wedge \tau_2(2M^*) > \log (2M^* B^{-1}))\\
&\hspace{4cm} + \sup_{z \in [B, 2M^*)} \P_{(0,z)}(\tau_1(-\beta/2) \ge T(M^*)),
\end{split}
\end{equation}
where we have used the strong Markov property in the last step.
By \eqref{timelength},
$$
 \sup_{z \in [B, 2M^*)}\P_{(-\beta/2,z)}\Big(\tau_2(B) \wedge \tau_2(2M^*) > \log (2M^* B^{-1})\Big) \le 1-\P(\mathcal{E}(M^*)).
$$
By using $S(0) + \sqrt{2}W(t) - \beta t/2 \ge S(t)$ for $t \le \tau_1(-\beta/2)$ and $Q_1(t) \le S(t) \le Q_2(t)$ for $t \ge 0$,
\begin{align*}
&\sup_{z \in [B, 2M^*)} \P_{(0,z)}\left(\tau_1(-\beta/2) \ge T(M^*)\right) \\
&\hspace{1.5cm}\le \P\left(\inf_{t \le T(M^*)}\left(\sqrt{2}W(t) - \beta t/2 \right) \ge  -(2M^* + \beta/2)\right)\\
&\hspace{1.5cm}\le \P\left(\sqrt{2} W(T(M^*)) \ge \beta T(M^*)/2 - (2M^* + \beta/2)\right) \le \P(\mathcal{E}(M^*))/2.
\end{align*}
Using these bounds in \eqref{downadd1}, we obtain
\begin{equation}\label{timefin}
\begin{split}
&\sup_{z \in [B, 2M^*)} \P_{(0,z)}\left(\tau_2(B) \wedge \tau_2(2M^*) > T(M^*) + \log (2M^* B^{-1})\right) \\
&\hspace{3cm}\le 1 - \frac{\P(\mathcal{E}(M^*))}{2} <1.
\end{split}
\end{equation}
Thus, using the strong Markov property and \eqref{timefin}, we obtain for any $k \ge 0$,
\begin{align}\label{down11}
\P_{(0,y)}(\alpha^*_{2k+1} - \alpha^*_{2k} > n\left(T(M^*) + \log (2M^* B^{-1})\right)) \le \left(1 - \frac{\P(\mathcal{E}(M^*))}{2}\right)^n.
\end{align}
Furthermore, by \eqref{down} we have constants $C_1$ and $C_2$, such that for $k \ge 1$ and for all $t \ge t^*$,
\begin{align}\label{down2}
\P_{(0,y)}(\alpha^*_{2k} - \alpha^*_{2k-1}>t) \le C_1\exp(-C_2t^{1/5}).
\end{align}
Writing $\alpha_1 = \sum_{j=0}^{2\mathcal{N}'} (\alpha^*_{j+1} - \alpha^*_j)$ and using \eqref{down11} and \eqref{down2}, we get positive constants $C, C', C''$ and $t_{\alpha}^{(2)}>0$, depending on $\beta, B, M^*$, such that for all $t \ge t_{\alpha}^{(2)}$,
\begin{align*}
\P_{(0,y)}(\alpha_1>t) &\le \P(\mathcal{N}'>n) + \P\Big(\sum_{j=0}^{2n} (\alpha^*_{j+1} - \alpha^*_j) > t\Big) \\
&\le \e^{-Cn} + C'n \e^{-C(t/n)^{1/5}} \le C'\e^{-C''t^{1/6}},
\end{align*}
where the last step is obtained by taking $n = \lfloor t^{1/6} \rfloor$.
\end{proof}

\subsection{Up-crossings of \texorpdfstring{$\mathbf{Q_2}$}{Q2}}\label{ssec:up}

In this subsection, we will prove tail asymptotics for the distribution of $\alpha_2-\alpha_1$ as stated in Lemma~\ref{lem:alphauptail}.
The proof consists of the following two major parts: (i) First we establish in Lemma~\ref{lem:alpha} the tail probability of the hitting time of $Q_2$ to level $2B$ starting below level $B$ when $Q_1(0)$ is not too small.
(ii) Then in Lemma~\ref{qonealphaone} we show that at time $\alpha_1$, $Q_1(\alpha_1)$ cannot be too small.
Lemmas~\ref{lem:alpha} and~\ref{qonealphaone} are combined to prove Lemma~\ref{lem:alphauptail}.
\begin{lemma}\label{lem:alpha}
For any fixed $B>0$ and $M > 8B + 6\beta$, there exists $c_\alpha^{(2)}>0$ (depending on $M,B,\beta$) such that for all $t\geq 9$,
$$\sup_{x \in [-M/2,0], \ y \in (0,B]}\P_{(x, y)}(\tau_2(2B)>t)\leq \exp(-c_\alpha^{(2)}\sqrt{t}).$$
\end{lemma}

In order to prove Lemma~\ref{lem:alpha}, set $M>0$ to be a fixed large number to be chosen later and $(Q_1(0), Q_2(0))=(x,y)$ for some $x \in [-M/2,0], y \in (0,B]$.
For $i\geq 1$ define the stopping times
\begin{align*}
\tau_{2,2i-1}&:= \inf\Big\{t\geq 0: Q_2(t)= 2B\quad\mbox{or}\quad Q_1(t)=-M\Big\},\\ 
\tau_{2,2i}&:= \inf\Big\{t\geq 0: Q_2(t)= 2B \quad\mbox{or}\quad Q_1(t)=-\frac{M}{2}\Big\},
\end{align*}
where by convention we take $\tau_{2,0}\equiv 0$.
Also define 
$$N^*:= \inf \Big\{k\geq 0: Q_2(\tau_{2,2k+1})= 2B\Big\}.$$
Therefore, note that 
\begin{equation}\label{eq:tau2B}
\tau_2(2B) = \sum_{j=1}^{2N^*+1}(\tau_{2,j}-\tau_{2,j-1}).
\end{equation}
The proof of Lemma \ref{lem:alpha} consists of three parts:
\begin{enumerate}[(i)]
\item Lemma~\ref{lem:q2upcross} contains the required probability estimate to analyze the time interval $\tau_{2, 2i-1}- \tau_{2,2i-2}$,
\item Lemma~\ref{lem:evencycle} contains estimates of the tail probabilities for the time interval $\tau_{2, 2i}- \tau_{2,2i-1}$,
and 
\item Lemma~\ref{lem:N-star-tail} provides estimates of the tail probabilities for the random variable $N^*$. Lemma~\ref{lem:exit-by-q2} is used in the proof of Lemma~\ref{lem:N-star-tail}.
\end{enumerate}
Combining Equation~\eqref{eq:tau2B} and Lemmas \ref{lem:q2upcross}, \ref{lem:evencycle}, and \ref{lem:N-star-tail}, we will complete the proof of Lemma~\ref{lem:alpha}.

\begin{lemma}\label{lem:q2upcross}
For any fixed $B,M> 0$, 
$$\inf_{\substack{x\in[-M,\ 0],\\ y \in (0, 2B]}}\P_{(x,y)}\Big(\sup_{0\leq s\leq 1}Q_2(s)>2B\Big)\geq p^{(1)}(M,B)>0.$$
\end{lemma}
\begin{proof}
Recall that 
$$Q_1(t) = Q_1(0) +\sqrt{2}W(t) -\beta t +\int_0^t (-Q_1(s) + Q_2(s))\dif s - L(t),$$
where 
\begin{equation}\label{eq:lt}
\begin{split}
L(t) &=\ \sup_{s\leq t} \Big(Q_1(0)+\sqrt{2}W(s)-\beta s + \int_0^s(-Q_1(u)+Q_2(u))\dif u\Big)^+\\
&\geq\ \sup_{s\leq t} (Q_1(0)+\sqrt{2}W(s)-\beta s)^+.
\end{split}
\end{equation}
Thus, $\Pro{L(1)> 4B}\geq \Pro{\sqrt{2}W(1)>\beta + 4B -Q_1(0)}$.
Observe that for any $Q_2(0) = y \le 2B$,
$$\big\{L(1)>4B\big\}\implies \big\{sup_{s\leq 1}Q_2(s)>2B\big\}.$$
To see this, suppose $L(1)>4B$. If $\sup_{s\leq 1} Q_2(s) \leq 2B$, then
$$Q_2(1) = y + L(1) - \int_0^1 Q_2(s) \dif s \ge L(1) - 2B > 2B,$$
which is a contradiction.
Therefore, 
\begin{align*}
&\inf_{\substack{x\in[-M,\ 0],\\ y \in (0,2B]}}\P_{(x,y)}\Big(\sup_{0\leq s\leq 1}Q_2(s)>2B\Big)
\geq
\inf_{\substack{x\in[-M,\ 0],\\ y \in (0, 2B]}}\P_{(x,y)}\Big(L(1)>4B\Big)\\
&\geq
\inf_{\substack{x\in[-M,\ 0],\\ y \in (0, 2B]}}\P_{(x,y)}\Big(\sqrt{2}W(1)>\beta + 4B - x\Big)\\
&\geq \Pro{\sqrt{2}W(1)>\beta + 4B+M} = p^{(1)}(M,B)>0.
\end{align*}
This completes the proof of Lemma~\ref{lem:q2upcross}.
\end{proof}

\begin{lemma}\label{lem:evencycle}
For any $j \ge 0$ and any fixed $M\ge 6\beta$, there exists $c_\tau^{(1)}>0$ such that for all $t\geq 2$,
$$\sup_{\substack{x\in[-M/2,\ 0],\\ y \in (0, B]}}\P_{(x,y)}\Big(\tau_{2,2j+2} - \tau_{2, 2j+1}>t\ \Big|\ N^*>j\Big)\leq \exp(-c_\tau^{(1)} t).$$
\end{lemma}
\begin{proof}
Let us denote $Q_1^* = Q_1+\beta$. Since $ N^*>j$, we know $Q_2(\tau_{2, 2j+1}) < 2B$.
In that case, for $t>\tau_{2,2j+1},$
\begin{align*}
Q_1^*(t) &= Q_1^*(\tau_{2, 2j+1}) + \sqrt{2}W(t) + \int_{\tau_{2,2j+1}}^t (-Q_1^*(s)+Q_2(s))\dif s\\
&\geq Q_1^*(\tau_{2, 2j+1}) +\sqrt{2}W(t) - \int_{\tau_{2,2j+1}}^tQ_1^*(s)\dif s\\
&= -M+\beta + \sqrt{2}W(t) - \int_{\tau_{2,2j+1}}^t Q_1^*(s)\dif s.
\end{align*}
Thus, we obtain
\begin{align*}
\P_{(x,y)}\Big(\tau_{2,2j+2} - \tau_{2, 2j+1}&>t\ \Big|\ N^*>j\Big)
\\
&\leq \P\Big(\sup_{s\leq t}(\sqrt{2}W(s) - (-M/2+\beta)s)\leq M/2\Big),
\end{align*}
since for $t\in (\tau_{2,2j+1}, \tau_{2,2j+2})$, $Q_1^*(s)\leq -M/2+\beta$.
Therefore, as $M \ge 6\beta$, for all $t\geq 2$,
\begin{align*}
\P_{(x,y)}\Big(\tau_{2,2j+2}-\tau_{2,2j+1}>t\ \Big|\ &N^*>j\Big)
\leq \Pro{\sqrt{2}W(t) \leq M/2- (M/2-\beta)t}\\
&\le \Pro{\sqrt{2}W(t) \leq - (M/2-\beta)t/4}\\
&\leq \exp(-c_\tau^{(1)}(M/2-\beta)^2 t)
 \le \exp(-c_\tau^{(1)} t),
\end{align*}
where $c_\tau^{(1)}$ does not depend on $x,y$.
\end{proof}

\begin{lemma}\label{lem:exit-by-q2}
For any fixed $B>0$ and $M>8B + 2\beta$, there exists positive $p^{(2)}=p^{(2)}(M,B)$ such that
\begin{align*}
&\inf_{\substack{ x\in [-M/2,\ 0], \\ y \in (0, B]}}\P_{(x,y)}\Big(\exists\ t^*\in [0,1],\ \mbox{ such that } \sup_{0\leq t\leq t^*}Q_2(t) \ge 2B,\\
&\hspace{7cm} \inf_{0\leq t\leq t^*}Q_1(t)>-M\Big)
\geq p^{(2)}.
\end{align*}
\end{lemma}
\begin{proof}
For fixed $B>0$ and $M>8B + 2\beta$, consider the event 
$$\cE(\beta,M):= \Big\{\sqrt{2}W(1)>\beta + 4B + \frac{M}{2}, \quad \inf_{t\in [0,1]}\sqrt{2}W(t)>\beta + 4B-\frac{M}{2}\Big\}.$$
From the representation \eqref{eq:lt}, note that the event $\cE(M,B)$ implies the event 
$\{L(1) > 4B\}$, which in turn implies that there exists $t^*\in [0,1]$ such that  $L(t^*)=4B$
and $\forall\ t\leq t^*$,
\begin{align*}
Q_1(t) &\geq -\frac{M}{2} +\sqrt{2}W(t) -\beta -4B
>-\frac{M}{2} -\beta -4B + \big(\beta + 4B -\frac{M}{2}\big)=-M.
\end{align*}
Therefore, $\inf_{0\leq t\leq t^*}Q_1(t) > -M$.
Furthermore, we claim that $$\sup_{0\leq t\leq t^*}Q_2(t) \ge 2B.$$
Indeed, if $\sup_{0\leq t\leq t^*}Q_2(t) < 2B$, then
$$
Q_2(t^*)\geq L(t^*) - \int_0^{t^*}Q_2(s)\dif s >4B - 2B t^*\geq 2B,
$$
since $0\leq t^*\leq 1$, which leads to a contradiction.
Finally, 
\begin{multline*}
\inf_{\substack{ x\in [-M/2,\ 0], \\ y \in (0, \beta^{-1}]}}\P_{(x,y)}\Big(\exists\ t^*\in [0,1],\ \mbox{ such that } \sup_{0\leq t\leq t^*}Q_2(t)>2B,\\ \inf_{0\leq t\leq t^*}Q_1(t)>-M\Big)
\geq \Pro{\cE(M,B)} >0.
\end{multline*}
This completes the proof of the lemma.
\end{proof}

\begin{lemma}\label{lem:N-star-tail}
For any fixed $B>0$ and $M> 8B + 2\beta$, there exist $c_N^{(2)}, n_N>0$ such that for all $n\geq n_N$, 
$$\sup_{\substack{x\in[-M/2,\ 0],\\ y \in (0, B]}}\P_{(x,y)}(N^*>n)\leq \exp(-c_N^{(2)}n).$$
\end{lemma}
\begin{proof}
Observe that
\begin{align*}
\P_{(x,y)}(N^*>n)&\leq \P_{(x,y)}(Q_1(\tau_{2,2k+1})=-M \mbox{ and }Q_2(\tau_{2,2k+1})< 2B \text{ for all } k \le n)\\
&\leq (1-p^*)^n,
\end{align*}
using the strong Markov property, where 
\begin{align*}
p^* &:= \inf_{\substack{x\in[-M/2,\ 0]\\ y\in(0, B]}}\P_{(x,y)}(Q_2\mbox{ hits } 2B \mbox{ before }Q_1\mbox{ hits }-M)\\
&\geq \inf_{\substack{x\in[-M/2,\ 0]\\ y\in(0, B]}}\P_{(x,y)}(\exists\ t^*\in[0,1]\mbox{ such that }\sup_{0\leq t\leq t^*}Q_2(t) >2B, \\
&\hspace{7cm}\inf_{0\leq t\leq t^*}Q_1(t)>-M)\\
&\geq p^{(2)}(M,B)>0,
\end{align*}
by Lemma \ref{lem:exit-by-q2}, choosing $M> 8B + 2\beta$.
\end{proof}

Now, we have all the necessary results to prove Lemma \ref{lem:alpha}.
\begin{proof}[Proof of Lemma~\ref{lem:alpha}]
Recall that $\tau_2(2B) = \sum_{j=1}^{2N^*+1}(\tau_{2,j}-\tau_{2,j-1}).$
From Lemma~\ref{lem:q2upcross} observe that for any fixed $M> 0$ and any $x \in [-M/2,0], \ y \in (0,B]$,
\begin{align*}
\P_{(x, y)}(\tau_{2,1}>n) &=\mathbb{E}_{(x, y)}\left(\mathbbm{1}_{[\tau_{2,1}>n-1]}\P_{(Q_1(n-1), Q_2(n-1))}(\tau_{2,1}>1)\right)\\
&\leq (1-p^{(1)}(M,B))\Pro{\tau_{2,1}>n-1},
\end{align*}
which implies $\P_{(x,y)}(\tau_{2,1}>n) \leq (1-p^{(1)}(M,B))^n$.
Furthermore, following the same argument as above, we can claim that for all $j\geq 1$,
\begin{equation}\label{eq:oddcycle}
\P_{(x,y)}(\tau_{2,2j-1}-\tau_{2,2j-2}\geq n)\leq (1-p^{(1)}(M,B))^n.
\end{equation}
Therefore for $t\geq 9$, choosing $M> 8B + 6\beta$, we can write for any $x \in [-M/2,0]$ and $y \in (0,B]$,
\begin{align*}
&\P_{(x,y)}(\tau_2(2B)>t)
\leq \P_{(x,y)}(N^*>n) + \P_{(x,y)}\Big(\sum_{j=1}^{2n+1}(\tau_{2,j}-\tau_{2,j-1})>t\Big) \\
&\leq \exp(-c_N^{(2)}n) + (2n+1)\exp(-ct/(2n+1)),\hspace{.25cm} \mbox{Due to Lemmas~\ref{lem:evencycle} \&~\ref{lem:N-star-tail},}\\
&\hspace{7.34cm} \mbox{and~\eqref{eq:oddcycle}}\\
&\leq c'\sqrt{t}\e^{-c\sqrt{t}}
\leq \e^{c_\alpha^{(2)}\sqrt{t}},\hspace{3.47cm}\mbox{[choosing }n = \lfloor(\sqrt{t}-1)/2\rfloor\mbox{]}
\end{align*} 
where $c_\alpha^{(2)}$ does not depend on $(x,y)$.
\end{proof}

As mentioned earlier, the next lemma gives a tail estimate on the distribution of $Q_1(\alpha_1)$.
\begin{lemma}\label{qonealphaone}
Fix $(Q_1(0), Q_2(0))=(x,y)$ with $x \le 0$, $y>0$. 
Recall the constant $t_{\alpha}^{(1)}$ obtained in Lemma~\ref{lem:alphadowntail}. There exist constants $C_1, C_2>0$ possibly depending on $(x,y)$, $B$, and $\beta$, such that for all $A \ge \max\{8 \beta t_{\alpha}^{(1)}, -4x\}$,
$$
\P_{(x,y)}(Q_1(\alpha_1) < -A) \le C_1 \e^{-C_2 A^{1/6}}.
$$
\end{lemma}
\begin{proof}
In the proof, $C, C'$ will denote generic positive constants depending on $\beta, x, y$ whose values change from line to line. Observe that for $t>0$,
$$
Q_1(t) \ge Q_1(0) + \sqrt{2}W(t) - \beta t - L^*(t),
$$
where $L^*(t) = \sup\limits_{s \le t}(Q_1(0) + \sqrt{2}W(s) - \beta s)^+$.
Thus, for any $A \ge \max\{8 \beta t_{\alpha}^{(1)}, -4x\}$,
\begin{equation}\label{Qal1}
\begin{split}
&\P_{(x,y)}(Q_1(\alpha_1) < -A) \le \P_{(x,y)}(\alpha_1 > A/(8\beta)) + \P_{(x,y)}\left(\inf_{s \le A/(8\beta)}Q_1(s) < -A\right)\\
 &\le \P_{(x,y)}\left(\inf_{s \le A/(8\beta)}\left(Q_1(0) + \sqrt{2}W(s) - \beta s - L^*(s)\right) < -A\right) \\
 &\hspace{4cm}+\P_{(x,y)}(\alpha_1 > A/(8\beta)) \\
& \le \P_{(x,y)}(\alpha_1 > A/(8\beta)) + \P_{(x,y)}\left(L^*(A/(8\beta)) > A/2\right) \\
&\hspace{4cm}+ \P_{(x,y)}\left(\inf_{s \le A/(8\beta)}\left(\sqrt{2}W(s) - \beta s\right) < -A/2 - x\right)\\
& \le \P_{(x,y)}(\alpha_1 > A/(8\beta)) + \P_{(x,y)}\left(\sup_{s \le A/(8\beta)}(\sqrt{2}W(s) - \beta s) > A/2\right)\\
  &\hspace{4cm}+ \P_{(x,y)}\left(\inf_{s \le A/(8\beta)}\left(\sqrt{2}W(s) - \beta s\right) < -A/4\right).
\end{split}
\end{equation}
By Lemma \ref{lem:alphadowntail},
$$
\P_{(x,y)}(\alpha_1 > A/(8\beta)) \le C\e^{-C'A^{1/6}}.
$$ 
Using the fact that the scale function (see \cite[V.46]{RW2000_2}) for $\sqrt{2}W(t) - \beta t$ is $s(z)= \exp(\beta z)$,
\begin{align*}
&\P_{(x,y)}\left(\sup_{s \le A/(8\beta)}(\sqrt{2}W(s) - \beta s) > A/2\right) \\
&\hspace{3cm}\le \P_{(x,y)}\left(\sup_{s < \infty}(\sqrt{2}W(s) - \beta s) > A/2\right) = \e^{-\beta A/2}.
\end{align*}
Moreover, by standard estimates on normal distribution functions,
\begin{align*}
&\P_{(x,y)}\left(\inf_{s \le A/(8\beta)}\left(\sqrt{2}W(s) - \beta s\right) < -A/4\right) \\
&\hspace{3cm}\le \P_{(x,y)}\left(\inf_{s \le A/(8\beta)}\left(\sqrt{2}W(s)\right) < -A/8\right) \le C\e^{-C'A}.
\end{align*}
Using the above bounds in \eqref{Qal1}, we obtain
$$
\P_{(x,y)}(Q_1(\alpha_1) < -A) \le C\e^{-C'A^{1/6}}
$$
for any $A \ge \max\{8 \beta t_{\alpha}^{(1)}, -4x\}$, proving the lemma.
\end{proof}

\begin{proof}[Proof of Lemma~\ref{lem:alphauptail}]
In the proof, $C, C'$ will denote generic positive constants depending on $\beta, x, y$ whose values change from line to line. Fix $M > 8B + 6\beta + 2$. Let $s_{\alpha}' = \inf\{ t \ge \alpha_1: Q_1(t) = -M/2\}$. Take $t_{\alpha}' = 4\max\{9, 8 \beta t_{\alpha}^{(1)}, -4x\}$. Then for $t \ge t_{\alpha}'$,
\begin{equation}\label{upcross}
\begin{split}
\P_{(x,y)}(\alpha_2 - \alpha_1 >t) &\le \P_{(x,y)}(Q_1(\alpha_1) < -t/4) \\
&\hspace{1cm}+ \sup_{u \ge -t/4, v>0}\P_{(u,v)}(\tau_1(-M/2) > t/4)\\
&\hspace{2cm} + \sup_{y \in (0,B]}\P_{(-M/2, y)}(\tau_2(2B) > t/2).
\end{split}
\end{equation}
By Lemma \ref{qonealphaone},
$$
\P_{(x,y)}(Q_1(\alpha_1) < -t/4) \le C\e^{-C't^{1/6}}.
$$
By computations similar to Lemma \ref{lem:evencycle},
\begin{align*}
&\sup_{u \ge -t/4, v>0}\P_{(u,v)}(\tau_1(-M/2) > t/4) \\
&\hspace{1cm}\le \P\left(-\frac{t}{4} + \sqrt{2}W(t/2) +\left(\frac{M}{2} - \beta\right)\frac{t}{2} < - \frac{M}{2}\right)\\
&\hspace{1cm}\le \P\left(\sqrt{2}W(t/2) < - \frac{t}{4}\right) \le C\e^{-C't}.
\end{align*}
By Lemma \ref{lem:alpha},
$$
\sup_{y \in (0,B]}\P_{(-M/2, y)}(\tau_2(2B) > t/2) \le \e^{-C'\sqrt{t}}.
$$
Using these bounds in \eqref{upcross}, we obtain for all $t \ge t_{\alpha}'$,
$$
\P_{(x,y)}(\alpha_2 - \alpha_1 >t) \le C\e^{-C't^{1/6}}
$$
proving the lemma.
\end{proof}

\section{Analysis of fluctuations within a renewal cycle}\label{sec:fluc}
In this section we prove Theorem~\ref{th:excren}.
Specifically,  we derive sharp estimates for the fluctuations of excursions of $Q_1$ and $Q_2$ between two successive regeneration times defined in \eqref{rendef}.
This will eventually furnish tail estimates for the stationary distribution of $Q_1$ and $Q_2$ and the scaling of extrema in large time intervals that are described in Theorems~\ref{th:statail} and~\ref{th:lil}.
First we state and prove Lemmas~\ref{lem:Q2hit1} -- \ref{lowqone}, which provide all the necessary results for proving Theorem~\ref{th:excren} at the end of this section.\\

Denote the Brownian motion with drift $b$ and and its corresponding reflected analogue by
\begin{align*}
W^{(b)}(t) &:= \sqrt{2} W(t) + bt,\\
W^{(b)}_R(t) &:= \sqrt{2} W(t) + bt - \sup_{s \le t}\left(\sqrt{2} W(s) + bs\right),
\end{align*}
where $W$ denotes the standard Brownian motion.
Also, denote the local time of the reflected Brownian motion $W^{(b)}_R$ and its hitting time of level $z$ by $L^{(b)}$ and $\tau^{(b)}(z)$ respectively.

\begin{lemma}\label{lem:Q2hit1}
There exist positive constants $C_1, C_2>0$ that do not depend on $\beta$ such that
\begin{equation*}
\P_{(0, y +\beta)}\left(\tau_1\left(-\frac{\beta}{2}\right) \le \tau_2\left(\frac{y}{2} + \beta\right)\right) \le C_1 \e^{-C_2\beta y}
\end{equation*}
for $y \ge \frac{1}{4\beta}$ if $\beta \ge 1$ and $y \ge \frac{64}{\beta}\log \frac{1}{\beta}$ if $\beta <1$.
\end{lemma}
\begin{proof}
From the evolution equation of $Q_1$ in~\eqref{eq:diffusionjsq-aap}, note that for $y>0$, $W^{(y/2)}_R$ can be constructed on the same probability space as $(Q_1, Q_2)$, such that starting from $(Q_1(0), Q_2(0))= (0, y +\beta)$, almost surely $Q_1(t) \ge W^{(y/2)}_R(t)$ for all $t \le  \tau_2\left(\frac{y}{2} + \beta\right)$. 
The scale function $s$ for $W^{(y/2)}_R(t)$ is obtained by solving the equation $\frac{y}{2}s'(z) + s"(z)=0$ (see \cite[V.46]{RW2000_2}) and one candidate is
\begin{equation}\label{eq:scalf}
s(z) = \frac{2}{y} \left(1- \e^{-yz/2}\right).
\end{equation}
We will estimate the time taken by $W^{(y/2)}_R$ to hit the level $-\beta/2$. Define stopping times for the process $W^{(y/2)}$ as follows: For $i \ge 0$
\begin{align*}
\gamma_{i+1} &= \inf\big\{t \ge \gamma_{i} : W^{(y/2)}(t) - W^{(y/2)}(\gamma_i) \text{ hits } \beta/4 \text{ or } -\beta/4\big\},
\end{align*}
with the convention that $\gamma_0=0$.
From the explicit form of the scale function $s$ in~\eqref{eq:scalf}, observe that for $i \ge 0$,
\begin{equation}\label{scalebound}
\P\left(W^{(y/2)}(\gamma_{i+1}) - W^{(y/2)}(\gamma_i) = -\beta/4\right) = \frac{1-\e^{-\beta y/8}}{\e^{\beta y/8} - \e^{-\beta y/8}} \le \e^{-\beta y /8}.
\end{equation}
Define
$$
\mathcal{N} := \inf\big\{ i \ge 1: W^{(y/2)}(\gamma_{i+1}) - W^{(y/2)}(\gamma_i) = -\beta/4\big\}.
$$
Then for any $n \ge 1$, by \eqref{scalebound}, $\P(\mathcal{N} \le n) \le n \e^{-\beta y /8}$. Note that for $t < \gamma_{\mathcal{N}}$, $W^{(y/2)}_R(t) > - \beta/2$. Thus, $\tau^{(y/2)}\left(-\frac{\beta}{2}\right) \ge \gamma_{\mathcal{N}}$. 
Consequently,
$$
L^{(y/2)}\left(\tau^{(y/2)}\left(-\frac{\beta}{2}\right)\right) \ge \sup_{t \le \gamma_{\mathcal{N}}}\left(W^{(y/2)}(t)\right) \ge \mathcal{N}\beta/4.
$$
Therefore, for any $n \ge 1$,
\begin{equation}\label{locbound}
\P\left(L^{(y/2)}\left(\tau^{(y/2)}\left(-\frac{\beta}{2}\right)\right) \le n \beta \right) \le \P\left(\mathcal{N} \le 4n\right) \le 4n \e^{-\beta y /8}.
\end{equation}
Further, on the event $\left[ \tau_1\left(-\frac{\beta}{2}\right) \le \tau_2\left(\frac{y}{2} + \beta\right)\right]$, $\tau_1\left(-\frac{\beta}{2}\right) \ge \tau^{(y/2)}\left(-\frac{\beta}{2}\right)$. Therefore, for $n \ge 1$,
\begin{equation}\label{tau1}
\begin{split}
&\P_{(0, y+ \beta)}\left(\tau_1\left(-\frac{\beta}{2}\right) \le n\beta/y, \tau_1\left(-\frac{\beta}{2}\right) \le \tau_2\left(\frac{y}{2} + \beta\right)\right)\\
&\le \P\left(\tau^{(y/2)}\left(-\frac{\beta}{2}\right) \le n\beta/y\right)\\
&\le \P\left(\tau^{(y/2)}\left(-\frac{\beta}{2}\right) \le n\beta/y, L^{(y/2)}\left(\tau^{(y/2)}\left(-\frac{\beta}{2}\right)\right) > n \beta\right)\\
&\hspace{3cm} + \P\left(L^{(y/2)}\left(\tau^{(y/2)}\left(-\frac{\beta}{2}\right)\right) \le n \beta\right).
\end{split}
\end{equation}
An upper bound for the second probability in the right side of~\eqref{tau1} has been obtained in \eqref{locbound}. To estimate the first probability, observe that
\begin{equation}\label{tau2}
\begin{split}
&\P\left(\tau^{(y/2)}\left(-\frac{\beta}{2}\right) \le n\beta/y, L^{(y/2)}\left(\tau^{(y/2)}\left(-\frac{\beta}{2}\right)\right) > n \beta\right)\\
&\hspace{1.5cm}\le \P\left(\sup_{t \le n\beta/y}\left(\sqrt{2}W(t) + yt/2\right) > n\beta\right)\\
& \hspace{1.5cm}\le \P\left(\sup_{t \le n\beta/y}\sqrt{2}W(t) > n\beta/2\right) \le \frac{4}{\sqrt{\pi n\beta y}} \e^{-n\beta y/16}.
\end{split}
\end{equation}
Using \eqref{locbound} and \eqref{tau2} in \eqref{tau1}, we obtain
\begin{equation}\label{imp1}
\begin{split}
&\P_{(0, y+ \beta)}\left(\tau_1\left(-\frac{\beta}{2}\right) \le n\beta/y, \tau_1\left(-\frac{\beta}{2}\right) \le \tau_2\left(\frac{y}{2} + \beta\right)\right)\\
&\hspace{3cm} \le 4n \e^{-\beta y /8} + \frac{4}{\sqrt{\pi n\beta y}} \e^{-n\beta y/16},
\end{split}
\end{equation}
where an appropriate choice of $n \ge 1$ (depending on $y$ and $\beta$) will be made later. Now, we want to estimate the probability 
$$\P_{(0, y+ \beta)}\left(n\beta/y < \tau_1\left(-\frac{\beta}{2}\right) \le \tau_2\left(\frac{y}{2} + \beta\right)\right).$$ 
Towards this end, recall that $S(t) = Q_1(t) + Q_2(t)$ has the representation
$$
S(t) = S(0) + \sqrt{2}W(t) - \beta t + \int_0^t (-Q_1(s))ds.
$$
Thus, for $t \le \tau_1\left(-\frac{\beta}{2}\right)$,
$$
S(t) \le S(0) + \sqrt{2}W(t) - \frac{\beta}{2} t.
$$
Therefore, if $n$ is chosen such that $y \le \sqrt{n} \beta/4$,
\begin{equation}\label{imp2}
\begin{split}
&\P_{(0, y+ \beta)}\left(n\beta/y < \tau_1\left(-\frac{\beta}{2}\right) \le \tau_2\left(\frac{y}{2} + \beta\right)\right)\\
 &\le \P\left(y+\beta + \sqrt{2}W(t) - \frac{\beta}{2} t \ge y/2 + \beta/2, \text{ for all } t \le n\beta/y\right)\\
 &\le \P\left(\sqrt{2}W(n\beta/y) - \frac{n\beta^2}{2y} \ge -y/2 - \beta/2\right)\\
& \le \P\left(\sqrt{2}W(n\beta/y) \ge \frac{n\beta^2}{8y}\right) \qquad \text{since } y \le \sqrt{n} \beta/4\\
 &\le \frac{8\sqrt{y}}{\sqrt{\pi n} \beta^{3/2}} \e^{-\frac{n\beta^3}{256 y}}.
\end{split}
\end{equation}
From \eqref{imp1} and \eqref{imp2} we obtain
\begin{equation}\label{imp3}
\begin{split}
&\P_{(0, y+ \beta)}\left(\tau_1\left(-\frac{\beta}{2}\right) \le \tau_2\left(\frac{y}{2} + \beta\right)\right)\\
&\hspace{2cm} \le 4n \e^{-\beta y /8} + \frac{4}{\sqrt{\pi n\beta y}} \e^{-n\beta y/16} + \frac{8\sqrt{y}}{\sqrt{\pi n} \beta^{3/2}} \e^{-\frac{n\beta^3}{256 y}}.
\end{split}
\end{equation}
Now, if $\beta \ge 1$, choose $n = 16 y^2\beta^2$. Then, clearly $y \le \sqrt{n} \beta/4$. With this choice of $n$, the above expression yields the bound
\begin{equation}\label{imp4}
\begin{split}
&\P_{(0, y+ \beta)}\left(\tau_1\left(-\frac{\beta}{2}\right) \le \tau_2\left(\frac{y}{2} + \beta\right)\right)\\
&\hspace{2cm}  \le 64(\beta y)^2 \e^{-\beta y /8} + \frac{1}{\sqrt{\pi}(\beta y)^{3/2}} \e^{-(\beta y)^3} + \frac{2}{\sqrt{\pi \beta y}} \e^{-\frac{\beta y}{16}}
\end{split}
\end{equation}
for $y \ge \frac{1}{4\beta}$ (this ensures $n \ge 1$).\\

If $\beta < 1$, choose $n=y^4$. 
Then $y \le \sqrt{n} \beta/4$ is satisfied if $y \ge 4/\beta$. Some routine calculations reveal that for $y \ge 4/\beta$ the second and third terms appearing on the right side of~\eqref{imp3} can be estimated by
$$
\frac{4}{\sqrt{\pi n\beta y}} \e^{-n\beta y/16} \le \frac{1}{8\sqrt{\pi}} \e^{-(\beta y)^5/16}
$$
and
$$
\frac{8\sqrt{y}}{\sqrt{\pi n} \beta^{3/2}} \e^{-\frac{n\beta^3}{256 y}} \le \frac{1}{\sqrt{\pi}}\e^{-(\beta y)^3/256}.
$$
To estimate the first term on the right side of~\eqref{imp3}, rewrite it as
$$
4n \e^{-\beta y /8} = \left[4(\beta y)^4 \e^{-(\beta y)/16}\right] \left[\beta^{-4}\e^{-(\beta y)/16}\right].
$$
Observe that $\beta^{-4}\e^{-(\beta y)/16} \le 1$ for $y \ge \frac{64}{\beta}\log \frac{1}{\beta}$. Therefore, for $\beta<1$ and $y \ge \frac{64}{\beta}\log \frac{1}{\beta}$, we have the following bound:
\begin{equation}\label{imp5}
\begin{split}
&\P_{(0, y+ \beta)}\left(\tau_1\left(-\frac{\beta}{2}\right) \le \tau_2\left(\frac{y}{2} + \beta\right)\right)\\
&\hspace{2cm}  \le \frac{1}{8\sqrt{\pi}} \e^{-(\beta y)^5/16} + \frac{1}{\sqrt{\pi}}\e^{-(\beta y)^3/256} + 4(\beta y)^4 \e^{-(\beta y)/16}.
\end{split}
\end{equation}
The lemma follows from \eqref{imp4} and \eqref{imp5}.
\end{proof}
The above lemma can be used to deduce the following hitting-time estimate for $Q_2$.
\begin{lemma}\label{Q2gebeta1}
There exist constants $\widetilde{C}_1, \widetilde{C}_2 >0$ that do not depend on $\beta$ such that
\begin{equation*}
\P_{(0, y+ \beta)}\left(\tau_2\left(2y + \beta\right) \le \tau_2\left(\frac{y}{2} + \beta\right)\right) \le \widetilde{C}_1 \e^{-\widetilde{C}_2 \beta y}
\end{equation*}
for $y \ge \frac{1}{4\beta}$ if $\beta \ge 1$ and $y \ge \frac{64}{\beta}\log \frac{1}{\beta}$ if $\beta <1$.
\end{lemma}
\begin{proof}
We can write for any $y>0$,
\begin{eq}\label{Q21}
\P_{(0, y+ \beta)}\left(\tau_2\left(2y + \beta\right) \le \tau_2\left(\frac{y}{2} + \beta\right)\right)
&\le \P_{(0, y+ \beta)}\left(\tau_1\left(-\frac{\beta}{2}\right) \le \tau_2\left(\frac{y}{2} + \beta\right)\right) \\
&+ \P_{(0, y+ \beta)}\left(\tau_2\left(2y + \beta\right) < \tau_1\left(-\frac{\beta}{2}\right)\right).
\end{eq}
By Lemma \ref{lem:Q2hit1},
\begin{equation}\label{Q22}
\P_{(0, y+ \beta)}\left(\tau_1\left(-\frac{\beta}{2}\right) \le \tau_2\left(\frac{y}{2} + \beta\right)\right) \le C_1\e^{-C_2\beta y}
\end{equation}
for $y \ge \frac{1}{4\beta}$ if $\beta \ge 1$ and $y \ge \frac{64}{\beta}\log \frac{1}{\beta}$ if $\beta <1$.
To estimate the second probability in \eqref{Q21}, recall that for $t \le \tau_1\left(-\frac{\beta}{2}\right)$, $S(t) = Q_1(t) + Q_2(t)$ satisfies
$$
S(t) \le S(0) + \sqrt{2}W(t) - \frac{\beta}{2} t.
$$
Therefore,
\begin{eq}\label{Q23}
&\P_{(0, y+ \beta)}\left(\tau_2\left(2y + \beta\right) < \tau_1\left(-\frac{\beta}{2}\right)\right)\\
&\hspace{3cm}
\le \P\left( \sup_{t< \infty}\left(\sqrt{2}W(t) - \frac{\beta}{2} t\right) \ge y\right) = \e^{-\frac{\beta y}{2}}
\end{eq}
for $y >0$. The first inequality above follows from the fact that points of time where $Q_2$ increases are precisely those where $Q_1$ equals zero: hence $Q_1(\tau_2(2y+\beta))=0$.

The lemma now follows by using \eqref{Q22} and \eqref{Q23} in \eqref{Q21}.
\end{proof}
The above estimate can be strengthened to the following tail estimate which will be used to study fluctuations of $Q_2$ between successive regeneration times.
\begin{lemma}\label{Q2gebeta2}
Recall the constants $\widetilde{C}_1, \widetilde{C}_2$ in the statement of Lemma \ref{Q2gebeta1}. There exist constants $C^*_1, C^*_2 >0$ that do not depend on $\beta$ such that
\begin{equation*}
\P_{(0, y+ \beta)}\left(\tau_2\left(2y + \beta\right) \le \tau_2\left(y_0 + \beta\right)\right) \le C^*_1 \e^{-C^*_2 \beta y}\quad\mbox{for all}\quad y \ge y_0,
\end{equation*}
 where $y_0 = \max\left\lbrace \frac{1}{4\beta}, \frac{\log(4\widetilde{C}_1)}{\widetilde{C}_2\beta}\right\rbrace$ if $\beta \ge 1$ and $\max\left\lbrace\frac{64}{\beta}\log \frac{1}{\beta}, \frac{\log(4\widetilde{C}_1)}{\widetilde{C}_2\beta}\right\rbrace$ if $\beta <1$.
\end{lemma}
\begin{proof}
Define stopping times:
\begin{align*}
T_{2k+1}&=\inf\left\lbrace t \ge T_{2k}: Q_2(t)=2y + \beta \text{ or } Q_2(t)=\frac{y}{2} + \beta \text{ or } Q_2(t)=y_0 + \beta \right\rbrace;\\
T_{2k+2}&=\inf\left\lbrace t \ge T_{2k+1}: Q_2(t)=y + \beta \text{ or } Q_2(t)=y_0 + \beta \right\rbrace,
\end{align*}
for $k \ge 0$, with the convention that $T_0=0$. Let
$$
\mathcal{N}^{0} = \inf\{ k \ge 1: Q_2\left(T_{2k}\right) = y_0 + \beta\}.
$$
Define $\hat{Q}_2(t) = \log_2(Q_2(t) - \beta)$. By Lemma \ref{Q2gebeta1} and our choice of $y_0$, for any $z \ge \log_2(y_0)$,
\begin{multline*}
\P(\hat{Q}_2 \text{ hits } z+1 \text{ before } z-1 \ \vert \ \hat{Q}_2(0) = z, Q_1(0)=0)\\
=\P_{(0,2^z + \beta)}(Q_2 \text{ hits } 2^{z+1}+\beta \text{ before } 2^{z-1} + \beta) \le 1/4.
\end{multline*}
Thus, $\hat{Q}_2$ starting from any $z \ge \log_2(y_0)$ and observed at the stopping times where the increments are $\pm1$ until the first time it crosses the level $\log_2(y_0)$ (i.e., strictly less than $\log_2(y_0)$) is stochastically dominated by a random walk $\left(\mathcal{S}_n\right)_{n \ge 0}$ where 
$$\P(S_{n+1} - S_n =1) = 1- \P(S_{n+1} - S_n =-1) = 1/4.$$ Therefore,
\begin{multline*}
\sup_{z \ge \log_2(y_0)}\P(\hat{Q}_2 \text{ hits } z+1 \text{ before it crosses } \log_2(y_0) \ \vert \ \hat{Q}_2(0) = z, Q_1(0)=0)\\
\le \sup_{z \ge \log_2(y_0)}\P(S_n \text{ hits } z+1 \ \vert \ S_0=z) = p^{(S)} <1,
\end{multline*}
which, in turn, implies that for any $y \ge y_0$,
\begin{align*}
\P_{(0, y+ \beta)}\left(\tau_2\left(2y + \beta\right) \le \tau_2\left(y_0 + \beta\right)\right) \le p^{(S)} <1.
\end{align*}
Thus, for any $k \ge 1$,
\begin{equation}\label{Nzero}
\P_{(0, y+ \beta)}\left(\mathcal{N}^{0} \ge k+1\right) \le (p^{(S)})^k.
\end{equation}
Finally, for any $y \ge y_0$,
\begin{align*}
&\P_{(0, y+ \beta)}\left(\tau_2\left(2y + \beta\right) \le \tau_2\left(y_0 + \beta\right)\right) 
=\P_{(0, y+ \beta)}\Big(\sup_{0 \le t \le T_{2\mathcal{N}^0}}Q_2(t) > 2y + \beta\Big)\\
&\le \sum_{k=1}^{\infty}\P_{(0, y+ \beta)}\left(\sup_{T_{2k-2} \le t \le T_{2k}}Q_2(t) > 2y + \beta, \mathcal{N}^0 \ge k\right)\\
 &=  \sum_{k=1}^{\infty}\mathbb{E}_{(0, y+ \beta)}\mathbb{I}(\mathcal{N}^0 \ge k) \P_{(0, y+ \beta)}\left(\tau_2\left(2y + \beta\right) \le \tau_2\left(\frac{y}{2} + \beta\right)\right), \\
 &\hspace{5.8cm} \text{by the strong Markov property at } T_{2k-2}\\
& \le \P_{(0, y+ \beta)}\left(\tau_2\left(2y + \beta\right) \le \tau_2\left(\frac{y}{2} + \beta\right)\right) \sum_{k=1}^{\infty} (p^{(S)})^{k-1}, \qquad \text{by } \eqref{Nzero}\\
 &\le (1-p^{(S)})^{-1} \widetilde{C}_1 \e^{-\widetilde{C}_2 \beta y},  \qquad \text{by Lemma } \ref{Q2gebeta1},
\end{align*}
which completes the proof of the lemma.
\end{proof}
The lower bound on the tail probabilities is achieved for all $\beta>0$ in the following lemma.
\begin{lemma}\label{Q2lb}
For any $\beta>0$ and any $B>0$,
\begin{equation*}
\P_{(0,2B)}\left(\tau_2(y) < \tau_2(B)\right) \ge (1-\e^{-\beta B})\e^{-\beta(y-2B)}
\end{equation*}
for all $y \ge 2B$.
\end{lemma}
\begin{proof}
Note that $Q_2(t) \ge Q_1(t) + Q_2(t) = S(t)$ for all $t\ge 0$. Further, recall that
$$
S(t) = S(0) + \sqrt{2}W(t) - \beta t + \int_0^t(-Q_1(s))ds \ge S(0) + \sqrt{2}W(t) - \beta t, \ t \ge 0.
$$
Therefore, for all $y \ge 2B$,
\begin{align*}
\P_{(0,2B)}\left(\tau_2(y) < \tau_2(B)\right) &\ge \P_{(0,2B)}\left(S(t) \text{ hits level } y \text{ before level } B\right)\\
&\ge \P\left(2B + \sqrt{2}W(t) - \beta t \text{ hits level } y \text{ before level }B\right)\\
&=\P\left(\sqrt{2}W(t) - \beta t \text{ hits level } y-2B \text{ before level } -B\right)\\
&= \frac{1-\e^{-\beta B}}{\e^{\beta(y-2B)} - \e^{-\beta B}}, \qquad \text{by scale function arguments}\\
&\ge (1-\e^{-\beta B})\e^{-\beta(y-2B)},
\end{align*}
proving the lemma.
\end{proof}
Now, we will study fluctuations of $Q_1$ within one renewal cycle. 
Recall 
$l_0(\beta)$ from~\eqref{lzerodef} and the notation
$$
\sigma(t) =\inf\{ s \ge t: Q_1(s)=0\}, \ \ t \ge 0.
$$
\begin{lemma}\label{unifprob}
There exist constants $R_1>0$ not depending on $\beta$ and $p^{**}(\beta) \in (0,1)$ such that for all $R \ge R_1$,
\begin{equation}
\sup_{y \ge Rl_0(\beta)} \P_{(0,y)}\left(\tau_1(-\beta) < \tau_2(Rl_0(\beta))\right) = p^*(\beta,R)\leq p^{**}(\beta).
\end{equation}
\end{lemma}
\begin{proof}
In the proof $C,C',C_1, C_2, \dots$ will denote generic positive constants not depending on $\beta$ and $R$ whose values might change from line to line. For any $y \ge Rl_0(\beta) - \beta$,
\begin{eq}\label{R11}
\P_{(0,y+\beta)}\left(\tau_1(-\beta) < \sigma\left(\tau_2\left(\frac{y}{2} + \beta\right)\right)\right)
 \le \P_{(0, y+\beta)} \left(\tau_1(-\beta/2) \le \tau_2\left(\frac{y}{2} + \beta\right)\right)\\
  + \P_{(0,y+\beta)}\left(\tau_2\left(\frac{y}{2} + \beta\right) < \tau_1(-\beta/2) < \tau_1(-\beta) < \sigma\left(\tau_2\left(\frac{y}{2} + \beta\right)\right)\right) \\
  \le \P_{(0, y+\beta)} \left(\tau_1(-\beta/2) \le \tau_2\left(\frac{y}{2} + \beta\right)\right) + \sup_{x \in [-\beta/2, 0]}\P_{\left(x,\frac{y}{2} + \beta\right)}\left(\tau_1(-\beta) < \tau_1(0)\right)
\end{eq}
where the last step is a consequence of the strong Markov property applied at the time $\tau_2\left(\frac{y}{2} + \beta\right)$.
From Lemma \ref{lem:Q2hit1}, for $R \ge 65$,
\begin{equation}\label{Rone}
\P_{(0, y+\beta)} \left(\tau_1(-\beta/2) \le \tau_2\left(\frac{y}{2} + \beta\right)\right) \le C_1\e^{-C_2\beta y}, \ \ y \ge Rl_0(\beta) - \beta.
\end{equation}
Now let us take the starting configuration to be $(Q_1(0), Q_2(0)) = \left(x,\frac{y}{2} + \beta\right)$ with $y \ge Rl_0(\beta) - \beta$ and $R \ge 5$.
In that case, since $(\dif/\dif t)Q_2(t)\geq -Q_2(t)$, therefore $Q_2(t)\geq (y/2 +\beta)/2$ for all $t \le \log 2$.
Consequently, for any $t \le \log 2$,
\begin{align*}
Q_1(t) &= Q_1(0) + \sqrt{2} W(t) - \beta t +
\int_0^t (- Q_1(s) + Q_2(s)) \dif s - L(t)\\
&\geq x + \sqrt{2} W(t) - \beta t +
\int_0^t Q_2(s) \dif s\\
&\ge x + \sqrt{2}W(t) + (y-2\beta)t/4 \ge x + \sqrt{2}W(t) + yt/8.
\end{align*}
Therefore,
\begin{equation}\label{R12}
\begin{split}
&\sup_{x \in [-\beta/2, 0]}\P_{\left(x,\frac{y}{2} + \beta\right)}\left(\tau_1(-\beta) < \tau_1(0) \le \log 2\right)\\
&\hspace{2cm} \le \sup_{x \in [-\beta/2, 0]}\P_{\left(x,\frac{y}{2} + \beta\right)}\left(x + \sqrt{2}W(t) + yt/8 \text{ hits } -\beta \text{ before } 0\right)\\
&\hspace{2cm}  \le \P\left(\sqrt{2}W(t) + yt/8 \text{ hits } -\beta/2 \text{ before } \beta/2\right) \le \e^{-\beta y/16}
\end{split}
\end{equation}
where the last step follows from standard scale function arguments. Moreover, for $y \ge Rl_0(\beta) - \beta$ with $R \ge 65$,
\begin{equation}\label{R13}
\begin{split}
&\sup_{x \in [-\beta/2, 0]}\P_{\left(x,\frac{y}{2} + \beta\right)}\left(\tau_1(0) > \log 2\right)\\
 &\le  \sup_{x \in [-\beta/2, 0]}\P_{\left(x,\frac{y}{2} + \beta\right)}\left(\sup_{t \le \log 2}(x + \sqrt{2}W(t) + yt/8) < 0\right)\\
&\le\P\left(\sup_{t \le \log 2}( \sqrt{2}W(t) + yt/8) < \beta/2\right)
 \le \P\left(\sqrt{2}W(\log 2)  < -y/32\right)\\
 &\le \e^{-y^2/(4(32^2) \log 2)} \le \e^{-\beta y/(64 \log 2)}.
\end{split}
\end{equation}
Using \eqref{Rone}--\eqref{R13} in \eqref{R11}, we obtain for $R \ge 65$, there exist positive constants $C,C'$ not depending on $\beta$ and $R$ such that for all $y \ge Rl_0(\beta) - \beta$
\begin{equation}\label{Rsigma}
\P_{(0,y+\beta)}\left(\tau_1(-\beta) < \sigma\left(\tau_2\left(\frac{y}{2} + \beta\right)\right)\right) \le C\e^{-C'\beta y}.
\end{equation}
Denote $\Gamma(y, \beta) = \left\lfloor\log_2\left(\frac{y}{Rl_0(\beta)-\beta}\right) + 2\right\rfloor$.
Now, for any $y \ge Rl_0(\beta) - \beta$, observe that the event  $[\tau_1(-\beta) \le \tau_2(Rl_0(\beta))]$ can be written as a subset of
\[ \bigcup_{k=1}^{\Gamma(y,\beta)}\Big[\sigma\left(\tau_2\left(\frac{y}{2^{k-1}} + \beta\right)\right) < \tau_1\left(-\beta\right) < \sigma\left(\tau_2\left(\frac{y}{2^k} + \beta\right)\right)\Big],\]
and therefore,
\begin{eq}\label{Rtwo}
&\P_{(0,y+\beta)}\left(\tau_1(-\beta) \le \tau_2(Rl_0(\beta))\right)\\
& \leq\sum_{k=1}^{\Gamma(y,\beta)}\P_{(0,y+\beta)}\left(\sigma\left(\tau_2\left(\frac{y}{2^{k-1}} + \beta\right)\right) < \tau_1(-\beta) < \sigma\left(\tau_2\left(\frac{y}{2^k} + \beta\right)\right)\right).
\end{eq}
Take any $R \ge 260$. 
For each $k \le \Gamma(y,\beta)$, by the strong Markov property,
\begin{align*}
&\P_{(0,y+\beta)}\left(\sigma\left(\tau_2\left(\frac{y}{2^{k-1}} + \beta\right)\right) < \tau_1(-\beta) < \sigma\left(\tau_2\left(\frac{y}{2^k} + \beta\right)\right)\right)\\
&\hspace{3cm} \le \sup_{z \in [y/2^k, y/2^{k-1}]} \P_{(0,z+\beta)}\left(\tau_1(-\beta) < \sigma\left(\tau_2\left(\frac{y}{2^k} + \beta\right)\right)\right)\\
&\hspace{3cm}\le \sup_{z \in [y/2^k, y/2^{k-1}]} \P_{(0,z+\beta)}\left(\tau_1(-\beta) < \sigma\left(\tau_2\left(\frac{z}{2} + \beta\right)\right)\right) \\
&\hspace{3cm}\le C\e^{-C'\beta y/2^k},
\end{align*}
where the last inequality follows from \eqref{Rsigma} as for $k \le \Gamma(y,\beta)$, $\frac{y}{2^k} \ge \frac{Rl_0(\beta) - \beta}{4} \ge \frac{R}{4}l_0(\beta) - \beta$ and $\frac{R}{4} \ge 65$.

Writing $p(\beta,R) = C_1\e^{-C_2\beta (Rl_0(\beta)-\beta)/4}$ and using the above bound in \eqref{Rtwo}, we obtain $R_1 > 0$ such that for any $R \ge R_1$ and any $y \ge Rl_0(\beta) - \beta$,
\begin{align}
\P_{(0,y+\beta)}(\tau_1(-\beta) &\le \tau_2(Rl_0(\beta))) \le \sum_{k=1}^{\Gamma(y,\beta)}C_1\e^{-C_2\beta y/2^k}\label{eq:lem5.5sum1}\\
& \le \sum_{k=0}^{\infty}p(\beta,R)^{2^k} \le \sum_{k=0}^{\infty}p(\beta,R_1)^{2^k} =: p^{**}(\beta) <1,\label{eq:lem5.5sum2}
\end{align}
where the second inequality can be seen as follows: For any $y\ge Rl_0(\beta) - \beta$, the last term in the sum in~\eqref{eq:lem5.5sum1} is bounded above by $p(\beta, R)$.
Also, starting from the last term and counting backwards in $k$, observe that each next term is the square of the previous term, which provides the $2^k$ in the exponent of $p(\beta, R)$ in~\eqref{eq:lem5.5sum2}.
Now, it is straightforward to see that for a fixed $\beta$ the first sum in \eqref{eq:lem5.5sum2} is a decreasing function in $R$, and is bounded away from 1 for all large enough $R$.
This proves the lemma.
\end{proof}
\begin{lemma}\label{alphaone}
There exists a constant $R_2>0$ not depending on $\beta$ such that for any $R \ge R_2$, there is a constant $C_2(\beta,R)>0$ (depending on $\beta, R$) satisfying
\begin{equation}
\sup_{z \in[-\beta, 0],y \ge 2R l_0(\beta)}\P_{(z,y)}\left(\tau_1(-x) < \tau_2(Rl_0(\beta))\right) \le C_2(\beta,R)\e^{-(x-\beta)^2/2}, 
\end{equation}
for all  $x \ge \beta +1$.
\end{lemma}
\begin{proof}
Take any $R>0$. Let $(Q_1(0), Q_2(0))=(z,y)$ where $z \in [-\beta,0]$ and $y \ge 2Rl_0(\beta)$. Define the stopping times: $\sigma^{(0)} =0$ and for $k \ge 0$,
\begin{align*}
\sigma^{(2k+1)} &= \inf\{t \ge \sigma^{(2k)}: Q_1(t) = -\beta-1 \text{ or } Q_2(t) \le Rl_0(\beta)\},\\
\sigma^{(2k+2)} &= \inf\{t \ge \sigma^{(2k+1)}: Q_1(t) = -\beta \text{ or } Q_2(t) \le Rl_0(\beta)\}.
\end{align*}
Define $\mathcal{N}^{\sigma} = \inf\{n \ge 1: Q_2(\sigma^{(n)}) \le Rl_0(\beta)\}$. Observe that for any $z \in[-\beta, 0]$, by the strong Markov property, we obtain
\begin{equation}\label{eq:lem5.6-1}
\begin{split}
&\sup_{y \ge R l_0(\beta)}\P_{(z,y)}\left(\tau_1(-\beta-1) < \tau_2(Rl_0(\beta))\right)\\
& \hspace{1cm}\le \sup_{y \ge R l_0(\beta)}\P_{(z,y)}\left(\tau_1(0) < \tau_1(-\beta-1) < \tau_2(Rl_0(\beta))\right)\\
& \hspace{3cm} + \sup_{y \ge R l_0(\beta)}\P_{(z,y)}\left(\tau_1(-\beta-1) < \tau_1(0) \wedge  \tau_2(Rl_0(\beta))\right)\\
&\hspace{1cm}\le \sup_{y \ge R l_0(\beta)}\P_{(0,y)}\left(\tau_1(-\beta) < \tau_2(Rl_0(\beta))\right) \\
& \hspace{3cm}+ \sup_{y \ge R l_0(\beta)}\P_{(z,y)}\left(\tau_1(-\beta-1) < \tau_1(0)\wedge  \tau_2(Rl_0(\beta))\right).
\end{split}
\end{equation}
By Lemma \ref{unifprob}, for large enough $R$,
\begin{equation}\label{lal}
\sup_{y \ge R l_0(\beta)}\P_{(0,y)}\left(\tau_1(-\beta) < \tau_2(Rl_0(\beta))\right) \le p^{**}(\beta) <1.
\end{equation}
Further, observe that for $t \le \tau_1(0)\wedge  \tau_2(Rl_0(\beta))$, 
$$Q_1(t) \ge z + \sqrt{2}W(t) + (Rl_0(\beta) - \beta)t \ge -\beta + \sqrt{2}W(t) + (Rl_0(\beta) - \beta)t.$$ 
Therefore,
\begin{multline}\label{nil}
 \sup_{y \ge R l_0(\beta)}\P_{(z,y)}\left(\tau_1(-\beta-1) < \tau_1(0)\wedge  \tau_2(Rl_0(\beta))\right)\\
 \le \P(-\beta + \sqrt{2}W(t) + (Rl_0(\beta) - \beta)t \text{ hits } -\beta -1 \text{ before } 0) \le \e^{-(Rl_0(\beta)-\beta)}.
\end{multline}
Using~\eqref{lal} and~\eqref{nil} in~\eqref{eq:lem5.6-1}, we conclude that there is $R_2>0$ such that for all $R \ge R_2$,
\begin{equation}\label{holud}
\sup_{z \in[-\beta, 0],y \ge R l_0(\beta)}\P_{(z,y)}\left(\tau_1(-\beta-1) < \tau_2(Rl_0(\beta))\right) \le p'(\beta,R)<1.
\end{equation}
Using \eqref{holud} and the strong Markov property, there exists a constant $C(\beta,R)>0$ depending on $\beta,R$ such that
\begin{equation}\label{Nexpo}
\begin{split}
\sup_{z \in[-\beta, 0],y \ge 2R l_0(\beta)}\mathbb{E}_{(z,y)}(\mathcal{N}^{\sigma}) \  &\le 2\sum_{n=0}^{\infty}\P(\mathcal{N}^{\sigma} > 2n)  \\
& \le 2\sum_{n=0}^{\infty}p'(\beta, R)^n \ \le C(\beta,R) \ <\infty.
\end{split}
\end{equation}
For $(Q_1(0), Q_2(0)) = (-\beta-u, y)$ for any $u\ge 1,y > 0$, by \cite[Proposition 2.18]{Karatzas}, a process $Z$ can be constructed on the same probability space as $(Q_1,Q_2)$, such that $Q_1(t) + \beta \ge Z(t)$ for $t \le \tau_1(0)$, where $Z$ is an Ornstein-Uhlenbeck process which solves the SDE:
$$
dZ(t) = \sqrt{2}dW(t) - Z(t)dt, \ \ Z(0)=-u.
$$
The scale function for $Z$ is given by $s_Z(z) = \int_0^z\e^{w^2/2}dw$. From this observation and elementary estimates on $s_Z$, we have for any $x \ge \beta+u$,
\begin{multline}\label{OUcom}
\sup_{y > 0}\P_{(-\beta-u, y)}\left(\tau_1(-x) < \tau_1(-\beta)\right) \le \P\left(Z(t) \text{ hits } -x+\beta \text{ before } 0\right)\\
= \frac{s_Z(0) - s_Z(-u)}{s_Z(0)-s_Z(-x+\beta)}\le \sqrt{9\pi/2}\e^{u^2/2} \e^{-(x-\beta)^2/2}.
\end{multline}
Finally, using \eqref{Nexpo} and \eqref{OUcom} along with the strong Markov property, for any $R \ge R_2$ and any $x \ge \beta+1$,
\begin{align*}
&\sup_{\substack{z \in[-\beta, 0],\\y \ge 2R l_0(\beta)}}\P_{(z,y)}\left(\tau_1(-x) < \tau_2(Rl_0(\beta))\right)\\
& \hspace{1cm}= \sup_{\substack{z \in[-\beta, 0],\\y \ge 2R l_0(\beta)}}\P_{(z,y)}\left(\inf_{t \le \sigma^{(\mathcal{N}^{\sigma})}}Q_1(t) < -x\right)\\
&\hspace{1cm}\le \sup_{\substack{z \in[-\beta, 0],\\y \ge 2R l_0(\beta)}}\sum_{k=0}^{\infty}\P_{(z,y)}\left(\inf_{t \in [\sigma^{(2k+1)}, \sigma^{(2k+2)}]}Q_1(t) < -x, \mathcal{N}^{\sigma} \ge 2k+2\right)\\
&\hspace{1cm}\le \sup_{\substack{z \in[-\beta, 0],\\y \ge 2R l_0(\beta)}}\sum_{k=0}^{\infty}\mathbb{E}_{(z,y)}\mathbbm{1}_{[\mathcal{N}^{\sigma} \ge 2k+2]}\sup_{y > 0}\P_{(-\beta-1, y)}\left(\tau_1(-x) < \tau_1(-\beta)\right)\\
&\hspace{1cm}\le \sup_{\substack{z \in[-\beta, 0],\\y \ge 2R l_0(\beta)}}\mathbb{E}_{(z,y)}(\mathcal{N}^{\sigma}) \ \sup_{y > 0}\P_{(-\beta-1, y)}\left(\tau_1(-x) < \tau_1(-\beta)\right)\\
&\hspace{1cm} \le C_2(\beta,R)\e^{-(x-\beta)^2/2}
\end{align*}
where $C_2(\beta,R)>0$ is a constant depending on $\beta, R$. This proves the lemma.
\end{proof}
\begin{lemma}\label{alphatwo}
For any $R > 1$ and any $x \ge 18Rl_0(\beta)$, there exists a positive constant $C_3(\beta,R)$ (depending on $\beta,R$) such that
\begin{equation*}
\sup_{z \in [-9Rl_0(\beta),0], \ y \le 2Rl_0(\beta)} \P_{(z,y)}\left(\tau_1(-x) < \tau_2(2Rl_0(\beta))\right) \le C_3(\beta,R)\e^{-(x-\beta)^2/2}.
\end{equation*}
\end{lemma}
\begin{proof}
Fix any $R>1$, $Q_1(0) = z \ge -9Rl_0(\beta)$ and $Q_2(0)=y \le 2Rl_0(\beta)$. Define the stopping times: $\gamma^{(0)} =0$ and for $k \ge 0$,
\begin{align*}
\gamma^{(2k+1)} &= \inf\{t \ge \gamma^{(2k)}: Q_1(t) = -18Rl_0(\beta) \text{ or } Q_2(t) = 2Rl_0(\beta)\},\\
\gamma^{(2k+2)} &= \inf\{t \ge \gamma^{(2k+1)}: Q_1(t) = -9Rl_0(\beta) \text{ or } Q_2(t) = 2Rl_0(\beta)\}.
\end{align*}
Define $\mathcal{N}^{\gamma} = \inf\{n \ge 1: Q_2(\gamma^{(n)}) = 2Rl_0(\beta)\}$. Taking $B=2Rl_0(\beta)$ and $M=18Rl_0(\beta)$ in Lemma~\ref{lem:exit-by-q2}, we know there exists $q(\beta,R)$ such that
\begin{multline}\label{Nexpo2}
\inf_{\substack{z \in [-9Rl_0(\beta),0], \\ y \le 2Rl_0(\beta)}} \P_{(z,y)}\left(\tau_2(2Rl_0(\beta)) < \tau_1(-18Rl_0(\beta))\right)\\
\ge \inf_{\substack{z \in [-9Rl_0(\beta),0], \\ y \le 2Rl_0(\beta)}} \P_{(z,y)}\left(\tau_2(4Rl_0(\beta)) < \tau_1(-18Rl_0(\beta))\right) \ge q(\beta,R)>0.
\end{multline}
Using \eqref{Nexpo2} and the strong Markov property, there exists a constant $C(\beta,R)>0$ depending on $\beta,R$ such that
\begin{multline}\label{Nexpo3}
\sup_{z \in [-9Rl_0(\beta),0], \ y \le 2Rl_0(\beta)} \mathbb{E}_{(z,y)}(\mathcal{N}^{\gamma}) \le 2\sum_{n=0}^{\infty}\P(\mathcal{N}^{\gamma} > 2n)\\
\le 2\sum_{n=0}^{\infty}(1-q(\beta,R))^n \le C(\beta,R) <\infty.
\end{multline}
Using \eqref{Nexpo3} and \eqref{OUcom} along with the strong Markov property, we obtain for any $x \ge 18Rl_0(\beta)$, 
\begin{align*}
&\sup_{\substack{z \in [-9Rl_0(\beta),0], \\ y \le 2Rl_0(\beta)}} \P_{(z,y)}\left(\tau_1(-x) < \tau_2(2Rl_0(\beta))\right) \\
&= \sup_{\substack{z \in [-9Rl_0(\beta),0], \\ y \le 2Rl_0(\beta)}} \P_{(z,y)}\left(\inf_{t \le \gamma^{(\mathcal{N}^{\gamma})}}Q_1(t) < -x\right)\\
&\le \sup_{\substack{z \in [-9Rl_0(\beta),0], \\ y \le 2Rl_0(\beta)}} \sum_{k=0}^{\infty}\P_{(z,y)}\left(\inf_{t \in [\gamma^{(2k+1)}, \gamma^{(2k+2)}]}Q_1(t) < -x, \mathcal{N}^{\gamma} \ge 2k+2\right)\\
&\le \sup_{\substack{z \in [-9Rl_0(\beta),0], \\ y \le 2Rl_0(\beta)}} \sum_{k=0}^{\infty}\mathbb{E}_{(z,y)}\mathbbm{1}_{[\mathcal{N}^{\gamma} \ge 2k+2]}\sup_{y > 0}\P_{(-18Rl_0(\beta), y)}\left(\tau_1(-x) < \tau_1(-\beta)\right)
\end{align*}
\begin{align*}
&\le \sup_{\substack{z \in [-9Rl_0(\beta),0], \\ y \le 2Rl_0(\beta)}} \mathbb{E}_{(z,y)}(\mathcal{N}^{\gamma}) \ \sup_{y > 0}\P_{(-18Rl_0(\beta), y)}\left(\tau_1(-x) < \tau_1(-\beta)\right)\\
&\le C_3(\beta,R)\e^{-(x-\beta)^2/2}
\end{align*}
for some constant $C_3(\beta,R)>0$ depending on $\beta,R$. This proves the lemma.
\end{proof}
Now, we are in a position to give an upper bound for the fluctuations of $Q_1$ between two successive regeneration times $\Xi_k$ and $\Xi_{k+1}$, $k \ge 0$, defined in~\eqref{rendef} taking $B = Rl_0(\beta)$ for sufficiently large fixed $R$.
\begin{lemma}\label{Q1ub}
Fix any $R \ge \max\{2,R_1,R_2\}$, where $R_1$ and $R_2$ are obtained from Lemmas \ref{unifprob} and \ref{alphaone} respectively. 
Let $(Q_1(0), Q_2(0))=(0,2Rl_0(\beta))$ and take $B=Rl_0(\beta)$ in~\eqref{rendef}. 
There exists a constant $C^*(\beta,R)>0$ depending on $\beta,R$ such that for any $x \ge 18Rl_0(\beta)$,
\begin{equation*}
\P_{(0,2Rl_0(\beta))}\left(\inf_{t \le \Xi_0} Q_1(t) < -x\right) \le C^*(\beta,R) \e^{-(x-2\beta)^2/8}.
\end{equation*}
\end{lemma} 
\begin{proof}
Choose and fix $R \ge \max\{2,R_1,R_2\}$. Define
$$
\Xi^* = \inf\big\{t \ge \tau_2(Rl_0(\beta)): Q_1(t) \ge -\beta -1\big\}.
$$
Then for any $x \ge 2(\beta+1)$,  by Lemma \ref{alphaone} and \eqref{OUcom} along with the strong Markov property,
\begin{equation}\label{interest}
\begin{split}
&\P_{(0,2Rl_0(\beta))}\left(\inf_{\tau_2(Rl_0(\beta)) \le t \le \Xi^*} Q_1(t) < -x\right)\\
& \leq \P_{(0,2Rl_0(\beta))}\left(\inf_{\tau_2(Rl_0(\beta)) \le t \le \Xi^*} Q_1(t) < -x, Q_1(\tau_2(Rl_0(\beta)))  \ge -x/2\right)\\
&\hspace{5cm}+ \P_{(0,2Rl_0(\beta))}\left(\tau_1(-x/2) < \tau_2(Rl_0(\beta))\right)\\
&\le  \sup_{u \in [1, \frac{x}{2} - \beta]}\P_{(-\beta-u, Rl_0(\beta))}\left(\tau_1(-x) < \tau_1(-\beta)\right)\\
&\hspace{5cm}+\P_{(0,2Rl_0(\beta))}\left(\tau_1(-x/2) < \tau_2(Rl_0(\beta))\right)\\
&\le  \sqrt{9\pi/2}\e^{(\frac{x}{2} - \beta)^2/2} \e^{-(x-\beta)^2/2}+C_2(\beta,R)\e^{-(\frac{x}{2}-\beta)^2/2}\\
&\le ( \sqrt{9\pi/2}+C_2(\beta,R))\e^{-(x-2\beta)^2/8}.
\end{split}
\end{equation}
Therefore, for any $x \ge 18Rl_0(\beta)$, using \eqref{interest} along with Lemmas~\ref{alphaone} and~\ref{alphatwo},
\begin{align*}
&\P_{(0,2Rl_0(\beta))}\left(\inf_{t \le \Xi_0} Q_1(t) < -x\right) \\
&\le \P_{(0,2Rl_0(\beta))}\left(\inf_{t \le \tau_2(Rl_0(\beta))} Q_1(t) < -x\right)\\
&\hspace{2cm} + \P_{(0,2Rl_0(\beta))}\left(\inf_{\tau_2(Rl_0(\beta)) \le t \le \Xi^*} Q_1(t) < -x\right) \\
&\hspace{4cm}+ \P_{(0,2Rl_0(\beta))}\left(\inf_{\Xi^* \le t \le \Xi_0} Q_1(t) < -x\right)\\
& \le \P_{(0,2Rl_0(\beta))}\left(\tau_1(-x) < \tau_2(Rl_0(\beta))\right)\\
&\hspace{2cm} + \sup_{z \in [-9Rl_0(\beta),0],\ y \le 2Rl_0(\beta)} \P_{(z,y)}\left(\tau_1(-x) < \tau_2(2Rl_0(\beta))\right)\\
&\hspace{4cm} + \P_{(0,2Rl_0(\beta))}\left(\inf_{\tau_2(Rl_0(\beta)) \le t \le \Xi^*} Q_1(t) < -x\right)\\
& \le C_2(\beta,R)\e^{-\frac{(x-\beta)^2}{2}} + \big( C_2(\beta,R) + \sqrt{9\pi/2}\big)\e^{-\frac{(x-2\beta)^2}{8}} + C_3(\beta,R)\e^{-\frac{(x-\beta)^2}{2}}\\
& \le C^*(\beta,R) \e^{-(x-2\beta)^2/8}
\end{align*}
which proves the lemma.
\end{proof}
Now, we prove a lower bound for the fluctuation of $Q_1$.
\begin{lemma}\label{lowqone}
Let $(Q_1(0), Q_2(0))=(0,2Rl_0(\beta))$ and and take $B=Rl_0(\beta)$ in~\eqref{rendef}. There exist constants $R^{**} >0$ not depending on $\beta$ such that for any $R \ge R^{**}$ and any $x \ge \beta$,
\begin{equation*}
\P_{(0,2Rl_0(\beta))}\left(\inf_{t \le \Xi_0} Q_1(t) < -x\right) \ge C^{**}(\beta, R) \e^{-x^2},
\end{equation*}
where the positive constant $C^{**}(\beta, R)$ depends on both $\beta$ and $R$.
\end{lemma} 
\begin{proof}
Using $y = 2Rl_0(\beta) -\beta$ in Lemma \ref{lem:Q2hit1}, we observe that there exists $R^{**}>0$ such that for all $R \ge R^{**}$, there is a constant $q_1(\beta,R)>0$ (depending on $\beta,R$) for which
\begin{equation}\label{low1}
\P_{(0,2Rl_0(\beta))}\left(\tau_1(-\beta/2) > \tau_2(Rl_0(\beta) + \beta/2\right) \ge q_1(\beta,R) >0.
\end{equation}
Recall $S(t) = Q_1(t) + Q_2(t)$. 
Recall that $Q_1(t) \le S(t) \le Q_2(t)$ for every $t$, and
when $Q_1(0)\in [0,\beta/2]$,
$$
S(t) = S(0) + \sqrt{2}W(t) - \beta t + \int_0^t(-Q_1(s))ds \le S(0) + \sqrt{2}W(t) - \frac{\beta}{2} t
$$
for $t \le \tau_1(-\beta/2)$.
Moreover, observe that if $Q_2(0)\leq 2Rl_0(\beta)$, then we have $Q_1(\tau_2(2Rl_0(\beta))) =0$ and consequently, 
$$S(\tau_2(2Rl_0(\beta))) = Q_2(\tau_2(2Rl_0(\beta))) = 2Rl_0(\beta).$$ 
Thus,
\begin{equation}\label{low2}
\begin{split}
&\sup_{z \in [-\beta/2,0]} \P_{(z,Rl_0(\beta) + \beta/2)}\left(\tau_2(2Rl_0(\beta)) < \tau_1(-\beta/2) \right)\\
& \le \sup_{z \in [-\beta/2,0]} \P_{(z,Rl_0(\beta) + \beta/2)}\left(S(t) \text{ hits } 2Rl_0(\beta) \text{ before } -\beta/2\right)\\
& \le \sup_{z \in [-\beta/2,0]} \P_{(z,Rl_0(\beta) + \beta/2)}\Big(z + Rl_0(\beta) + \beta/2 + \sqrt{2}W(t) -\frac{\beta}{2}t\\
&\hspace{7cm} \text{ hits } 2Rl_0(\beta) \text{ before } -\beta/2\Big)\\
& \le \P\left(\sqrt{2}W(t) -\frac{\beta}{2}t \text{ hits } Rl_0(\beta) -\beta/2\right) \\
&\le \e^{-\beta(Rl_0(\beta) -\beta/2)/2}=: 1-q_2(\beta, R) <1.
\end{split}
\end{equation}
For $y \le 2Rl_0(\beta)$ and $(Q_1(0), Q_2(0)) = (-\beta/2, y)$, by \cite[Proposition 2.18]{Karatzas}, a process $U$ can be constructed on the same probability space as $ (Q_1, Q_2)$ such that almost surely $Q_1(t) +\beta \le U(t)$ for all $t \le \tau_1(0)$, where $U$ is an Ornstein-Uhlenbeck process which solves the SDE:
$$
dU(t) = \sqrt{2}dW(t) + (2Rl_0(\beta)-U(t))dt, \ \ U(0)=\beta/2.
$$
The scale function for $U$ is given by $s_U(u) = \int_0^u\e^{(w-2Rl_0(\beta))^2/2}dw$. Therefore, by elementary estimates on $s_U$, there exists a constant $C(\beta,R)>0$ (depending on $\beta,R$) such that for any $x \ge \beta$,
\begin{multline}\label{low3}
\inf_{y \le 2Rl_0(\beta)}\P_{(-\beta/2, y)}\left(\tau_1(-x) < \tau_1(0)\right) \ge \P\left(U(t) \text{ hits } -(x-\beta) \text{ before } \beta\right)\\
= \frac{s_U(\beta) - s_U(\beta/2)}{s_U(\beta) - s_U(-(x-\beta))} \ge C(\beta,R)\e^{-x^2}.
\end{multline}
Recall the notation $\sigma(t) = \inf\{s \ge t: Q_1(s)=0\}$ and define the stopping time 
$$\sigma_R = \inf\{t > \tau_2(Rl_0(\beta)+\beta/2): Q_2(t)=2Rl_0(\beta)\}.$$ 
From \eqref{low1}--\eqref{low3} and the strong Markov property, for any $R \ge R^{**}$ and any $x \ge \beta$,
\begin{align*}
&\P_{(0,2Rl_0(\beta))}\left(\inf_{t \le \Xi_0} Q_1(t) < -x\right)\\
&\ge \P_{(0,2Rl_0(\beta))}(\tau_2(Rl_0(\beta)+\beta/2) < \tau_1(-\beta/2) < \sigma_R, \\
&\hspace{5cm} \tau_1(-x) \in (\tau_1(-\beta/2), \sigma(\tau_1(-\beta/2))))\\
&\ge \P_{(0,2Rl_0(\beta))}\left(\tau_1(-\beta/2) > \tau_2(Rl_0(\beta) + \beta/2)\right) \\
&\hspace{3cm}\times \inf_{z \in [-\beta/2,0]} \P_{(z,Rl_0(\beta) + \beta/2)}\left(\tau_1(-\beta/2) < \tau_2(2Rl_0(\beta)) \right)\\
&\hspace{5cm}\times \inf_{y \le 2Rl_0(\beta)}\P_{(-\beta/2, y)}\left(\tau_1(-x) < \tau_1(0)\right)\\
& \ge q_1(\beta,R)q_2(\beta, R)C(\beta,R)\e^{-x^2}.
\end{align*}
This proves the lemma.
\end{proof}

\begin{proof}[Proof of Theorem~\ref{th:excren}]
Fix any 
\begin{equation}\label{Rzerodef}
R_0 \ge 4\max\{64,\log(4\widetilde{C}_1)/\widetilde{C}_2, R_1,R_2, R^{**}\},
\end{equation}
where where $R_1, R_2$ and $R^{**}$ are obtained from Lemmas \ref{unifprob}, \ref{alphaone} and \ref{lowqone} respectively and $\widetilde{C}_1, \widetilde{C}_2$ are the constants defined in the statement of Lemma \ref{Q2gebeta1}. Choose $B = R_0l_0(\beta)$ in \eqref{rendef}.\\

\noindent
To prove (i), note that $y_0$ defined in Lemma \ref{Q2gebeta2} satisfies $y_0 + \beta < R_0l_0(\beta)$ for our specific choice of $R_0$. Therefore, taking $z=\frac{y-\beta}{2}$ in place of $y$ in  Lemma \ref{Q2gebeta2} and applying the strong Markov property at $\tau_2(z + \beta)$, we have for any $y \ge 4R_0l_0(\beta)$,
\begin{align*}
\P_{(0, 2R_0l_0(\beta))}\left(\tau_2(y) \le \Xi_0 \right) &= \P_{(0, 2R_0l_0(\beta))}\left(\tau_2(y) \le \tau_2(R_0l_0(\beta)) \right)\\
& \le \P_{(0, z+\beta)}\left(\tau_2(2z+\beta) \le \tau_2(R_0l_0(\beta)) \right)\\
& \le \P_{(0, z+\beta)}\left(\tau_2(2z+\beta) \le \tau_2(y_0 + \beta) \right)\le C^*_1 \e^{-C^*_2 \beta z}.
\end{align*}
Part (ii) follows from Lemma \ref{Q2lb} by taking $B= R_0l_0(\beta)$. Parts (iii) and (iv) are direct consequences of Lemmas \ref{Q1ub} and \ref{lowqone} respectively.
\end{proof}

\section{Proofs of the main results}\label{sec:proofmain}

\begin{proof}[Proof of Theorem~\ref{th:statail}]
We will show that the tail bounds stated in the theorem hold with $C_R(\beta) = 18R_0l_0(\beta)$ and $D_R(\beta) = 4R_0l_0(\beta)$, where $R_0$ is defined in \eqref{Rzerodef} and $l_0(\beta)$ was defined in \eqref{lzerodef}. Taking $B= R_0l_0(\beta)$ in Theorem \ref{th:stationary}, note that for any $x\ge 0, y>0$,
\begin{align}\label{re}
\pi(Q_1(\infty) < -x) &= \frac{\mathbb{E}_{(0, 2R_0l_0(\beta))}\left(\int_{0}^{\Xi_0}\mathbbm{1}_{[Q_1(s) < -x]}ds\right)}{\mathbb{E}_{(0, 2R_0l_0(\beta))}\left(\Xi_0\right)},\nonumber\\
\pi(Q_2(\infty) > y)& = \frac{\mathbb{E}_{(0, 2R_0l_0(\beta))}\left(\int_{0}^{\Xi_0}\mathbbm{1}_{[Q_2(s) > y]}ds\right)}{\mathbb{E}_{(0, 2R_0l_0(\beta))}\left(\Xi_0\right)}.
\end{align}
To prove the theorem, we only need to estimate the numerators in the above representation. By the Cauchy-Schwarz inequality, for $x \ge 18R_0l_0(\beta)$,
\begin{align*}
&\mathbb{E}_{(0, 2R_0l_0(\beta))}\left(\int_{0}^{\Xi_0}\mathbbm{1}_{[Q_1(s) < -x]}ds\right) \\
&\le \mathbb{E}_{(0, 2R_0l_0(\beta))}\left(\mathbbm{1}_{[\tau_1(-x)]< \Xi_0]} (\Xi_0 - \tau_1(-x))\right)\\
& \le \sqrt{\P_{(0,2R_0l_0(\beta))}(\tau_1(-x)< \Xi_0)} \sqrt{\mathbb{E}_{(0, 2R_0l_0(\beta))}(\Xi_0)^2} \\
&\le \sqrt{C^*(\beta)} \e^{-(x-2\beta)^2/16}\sqrt{\mathbb{E}_{(0, 2R_0l_0(\beta))}(\Xi_0^2)},
\end{align*}
where the last inequality is a consequence of Part (iii) of Theorem~\ref{th:excren}. By Proposition \ref{prop:Xi}, $\mathbb{E}_{(0,2R_0l_0(\beta))}\left(\Xi^2_0\right) < \infty$. 
Now, using this in the above bound, we obtain the upper bound on $\pi(Q_1(\infty) < -x)$ claimed in the theorem. The upper bound for $\pi(Q_2(\infty)>y)$ is obtained similarly using part (i) of Theorem \ref{th:excren}.

To obtain the lower bound on $\pi(Q_1(\infty) < -x)$, we proceed along the same line of arguments as in the proof of Lemma \ref{lowqone}. Recall the stopping time 
$$\sigma_R = \inf\{t > \tau_2(Rl_0(\beta)+\beta/2): Q_2(t)=2Rl_0(\beta)\}.$$ Observe that for $x \ge \beta$,
\begin{align*}
&\mathbb{E}_{(0, 2R_0l_0(\beta))}\left(\int_{0}^{\Xi_0}\mathbbm{1}_{[Q_1(s) < -x]}ds\right)\\
& \ge \mathbb{E}_{(0, 2R_0l_0(\beta))}\Big(\mathbbm{1}_{[\tau_2(R_0l_0(\beta) + \beta/2) < \tau_1(-\beta/2) < \sigma_{R_0}]}\int\limits_{\tau_1(-\beta/2)}^{\sigma(\tau_1(-\beta/2))}\mathbbm{1}_{[Q_1(s) < - x]}ds\Big)\\
& \ge \P_{(0, 2R_0l_0(\beta))}\left(\tau_2(R_0l_0(\beta) + \beta/2) < \tau_1(-\beta/2)\right)\\
&\hspace{2.5cm} \times\inf_{y \le 2R_0l_0(\beta)}\mathbb{E}_{(-\beta/2,y)}\left(\int_0^{\tau_1(0)}\mathbbm{1}_{[Q_1(s) < - x]}ds\right) \\
&\hspace{2.5cm} \times \inf_{z \in [-\beta/2,0]} \P_{(z,R_0l_0(\beta) + \beta/2)}\left(\tau_1(-\beta/2) < \tau_2(2R_0l_0(\beta)) \right)\\
& \ge q_1(\beta,R_0)q_2(\beta, R_0)\inf_{y \le 2R_0l_0(\beta)}\mathbb{E}_{(-\beta/2,y)}\left(\int_0^{\tau_1(0)}\mathbbm{1}_{[Q_1(s) < - x]}ds\right)
\end{align*}
where $q_1(\beta,R_0)>0, q_2(\beta,R_0)>0$ are obtained in \eqref{low1} and \eqref{low2} respectively with $R_0$ in place of $R$.

Recall that for $y \le 2R_0l_0(\beta)$ and $(Q_1(0), Q_2(0)) = (-\beta/2, y)$, by  \cite[Proposition 2.18]{Karatzas}, a process $U_{\beta/2}$ can be constructed on the same probability space as $(Q_1,Q_2)$, such that $Q_1(t) +\beta \le U_{\beta/2}(t)$ for $t \le \tau_1(0)$, where $U_z$ is an Ornstein-Uhlenbeck process which solves the SDE:
$$
dU_z(t) = \sqrt{2}dW(t) + (2Rl_0(\beta)-U(t))dt, \ \ U_z(0)=z,
$$
where the scale function for $U_z$ is given by $s_U(u) = \int_0^u\e^{(w-2Rl_0(\beta))^2/2}dw$.

Define $\tau_z^U(w) = \inf\{t \ge 0: U_z(t)=w\}$ and write the law of $U_z$ and the corresponding expectation as $\P^U_z$ and $\mathbb{E}^U_z$ respectively. Then, for $x \ge \beta$,
\begin{equation}\label{lowstatone}
\begin{split}
&\inf_{y \le 2R_0l_0(\beta)}\mathbb{E}_{(-\beta/2,y)}\left(\int_0^{\tau_1(0)}\mathbbm{1}_{[Q_1(s) < - x]}ds\right) \\
&\ge \mathbb{E}^U_{\beta/2}\left(\int_0^{\tau_{\beta/2}^U(\beta)}\mathbbm{1}_{[U_{\beta/2}(s) < - x+\beta]}ds\right), \quad\text{hence, by strong Markov property},\\
&\ge \P^U_{\beta/2}\left(\tau_{\beta/2}^U(-2x + \beta) < \tau_{\beta/2}^U(\beta)\right)\mathbb{E}^U_{-2x+\beta}\left(\tau_{-2x+\beta}^{U}(-x+\beta)\right)\\
&= \frac{s_U(\beta) - s_U(\beta/2)}{s_U(\beta) - s_U(-(2x-\beta))} \ \mathbb{E}^U_{-2x+\beta}\left(\tau_{-2x+\beta}^{U}(-x+\beta)\right)\\
& \ge C(\beta) \e^{-4x^2}\mathbb{E}^U_{-2x+\beta}\left(\tau_{-2x+\beta}^{U}(-x+\beta)\right),
\end{split}
\end{equation}
where $C(\beta)$ is a positive constant that only depends on $\beta$. 
Now, from the Doob representation of the Ornstein-Uhlenbeck process,
$$
U_{-2x+\beta}(t) = (-2x+\beta)\e^{-t} + 2R_0l_0(\beta)(1-\e^{-t}) + \e^{-t}\widetilde{W}(\e^{2t}-1)
$$
for a standard Brownian motion $\widetilde{W}$.
Therefore, taking $T=\log(5/4)$, for $x \ge 4R_0 l_0(\beta)$,
\begin{align*}
&\P^U_{-2x+\beta}\left(\tau_{-2x+\beta}^U(-x + \beta)  \le T\right)\\
& \le \P\left(\sup_{t \le T}\left((-2x+\beta)\e^{-t} + 2R_0l_0(\beta)(1-\e^{-t}) + \e^{-t}\widetilde{W}(\e^{2t}-1)\right) > -x+\beta\right)\\
& \le \P\left((-2x+\beta)\e^{-T} + 2R_0l_0(\beta)(1-\e^{-T}) + \sup_{t \le T}\left(\widetilde{W}(\e^{2t}-1)\right) > -x+\beta\right)
\end{align*}
\begin{flalign*}
& \le \P\left(\sup_{t \le T}\left(\widetilde{W}(\e^{2t}-1)\right) > x/2\right) \ \ \text{ by our choice of $T$}&&\\
& = \P\left(\sup_{t \le 1}\widetilde{W}(t) > \frac{x}{2\sqrt{\exp(2T)-1}}\right) \ \text{ by Brownian scaling }&&\\
& \le \frac{4\sqrt{\exp(2T)-1}}{\sqrt{2\pi}x} < \frac{1}{2}.&&
\end{flalign*}
Thus,
\begin{align*}
\mathbb{E}^U_{-2x+\beta}\left(\tau_{-2x+\beta}^{U}(-x+\beta)\right)&=\int_0^{\infty}\P^U_{-2x+\beta}\left(\tau_{-2x+\beta}^U(-x + \beta) > t\right)dt \\
&\ge \frac{1}{2}\log(5/4).
\end{align*}
Using this in \eqref{lowstatone} gives us the lower bound on $\pi(Q_1(\infty) < -x)$ claimed in the theorem.

Finally, we prove the lower bound on $\pi(Q_2(\infty) > y)$.
Note that by the strong Markov property, for any $y \ge R_0l_0(\beta)$,
\begin{equation}\label{lowtwo}
\begin{split}
&\mathbb{E}_{(0, 2R_0l_0(\beta))}\left(\int_{0}^{\Xi_0}\mathbbm{1}_{[Q_2(s) > y]}ds\right)\\
&\qquad \ge \P_{(0, 2R_0l_0(\beta))}\left(\tau_2(2y) \le \Xi_0\right) \times \mathbb{E}_{(0,2y)}\left(\tau_2(y)\right)\\
 &\qquad\ge (1-\e^{-\beta R_0l_0(\beta)})\e^{-\beta(2y-2R_0l_0(\beta))}\mathbb{E}_{(0,2y)}\left(\tau_2(y)\right)
\end{split}
\end{equation}
where the last step follows from Part (ii) of Theorem~\ref{th:excren}. Recall that
$$
Q_2(t) \ge S(t) \ge S(0) + \sqrt{2}W(t) - \beta t, \ \ t \ge 0,
$$
where $S(t) = Q_1(t) + Q_2(t)$. Therefore, starting with $(Q_1(0),Q_2(0)) = (0,2y)$, the hitting time of level $y$ of $Q_2$ is stochastically bounded below by the hitting time of $y$ by $S(0) + \sqrt{2}W(t) - \beta t$. Denoting the latter hitting time by $\tau^S(y)$, we obtain $\mathbb{E}_{(0,2y)}\left(\tau_2(y)\right) \ge \mathbb{E}_{(0,2y)}\left(\tau^S(y)\right)$. For $y \ge R_0 l_0(\beta)$,
\begin{align*}
&\P_{(0,2y)}\left(\tau^S(y) \le \frac{y}{2\beta}\right) = \P\left(\inf_{t \le \frac{y}{2\beta}}\left(2y + \sqrt{2}W(t) - \beta t\right) < y\right) \\
&\hspace{2cm}\le \P\left(\inf_{t \le \frac{y}{2\beta}}\left(\sqrt{2}W(t)\right) < -y/2\right)\\
&\hspace{2cm}= \P\left(\inf_{t \le 1}\left(W(t)\right) < -\sqrt{\beta y}/2\right) \le \frac{4}{\sqrt{2\pi \beta y}} \le \frac{4}{\sqrt{2\pi R_0}} < \frac{1}{2},
\end{align*}
for our choice of $R_0$. This gives
$$
\mathbb{E}_{(0,2y)}\left(\tau^S(y)\right) = \int_0^{\infty} \P_{(0,2y)}\left(\tau^S(y) > t\right) dt \ge \frac{y}{4\beta}.
$$
Using this in \eqref{lowtwo} gives us the lower bound on $\pi(Q_2(\infty) > y)$ claimed in the theorem.
\end{proof}  

\begin{proof}[Proof of Theorem~\ref{th:lil}]
Below we provide the proof of the fluctuation result for $Q_2$.
The proof for $Q_1$ follows using analogous arguments.

Take $\mathcal{C^*}$ in the theorem to be the positive constant $C^*_2$ not depending on $\beta$ that was obtained in Part (i) of Theorem~\ref{th:excren}. 
Fix $\epsilon \in (0,1/2)$. Fix any starting point $(Q_1(0), Q_2(0))=(x,y)$. Then by Parts (i) and (ii) of Theorem~\ref{th:excren}, we obtain constants $D_1(\beta)$ and $D_2(\beta)$ and an integer $N(\beta)>0$ depending only on $\beta$ and such that for all $n \ge N(\beta)$,
\begin{eqnarray*}
\P_{(x,y)}\left(\sup_{t \in [\Xi_{n}, \Xi_{n+1}]} Q_2(t) > \frac{2(1+\epsilon)\log n}{C^*_2 \beta}\right) &\le \frac{D_1(\beta)}{n^{1+\epsilon}},\\
\P_{(x,y)}\left(\sup_{t \in [\Xi_{n}, \Xi_{n+1}]} Q_2(t) > \frac{(1-\epsilon)\log n}{\beta}\right) &\ge \frac{D_2(\beta)}{n^{1-\epsilon}}.
\end{eqnarray*}
Therefore, by the Borel-Cantelli Lemma,
\begin{equation}\label{BC}
\frac{1-\epsilon}{\beta} \le \limsup_{n \rightarrow \infty}\frac{\sup_{t \in [\Xi_{n}, \Xi_{n+1}]} Q_2(t)}{\log n} \le \frac{2(1+\epsilon)}{C^*_2 \beta}, \ \ a.s.
\end{equation}
By Proposition \ref{prop:Xi}, $\mathbb{E}_{(0,2R_0l_0(\beta))}\left(\Xi_0\right) < \infty$ and as $\{\Xi_{n+1} - \Xi_{n}\}_{n \ge 0}$ are i.i.d., therefore by the Strong Law of Large Numbers,
\begin{equation}\label{SLLN}
\lim_{n \rightarrow \infty} \frac{\Xi_n}{n} \rightarrow \mathbb{E}_{(0,2R_0l_0(\beta))}\left(\Xi_0\right), \ \ a.s.
\end{equation}
From the lower bound in \eqref{BC}, with probability one, there exists a subsequence $\{n_k\}\subseteq\{n\}$ and $t_{n_k} \in [\Xi_{n_k}, \Xi_{n_k+1}]$ such that 
$$
Q_2(t_{n_k}) \ge (1-2\epsilon)\frac{\log n_k}{\beta}
$$ 
for all sufficiently large $k$. Moreover, by \eqref{SLLN}, almost surely,
$$
\log t_{n_k} \le \log \Xi_{n_k + 1} = \log \left(\frac{\Xi_{n_k+1}}{n_k+1}\right) + \log (n_k+1) \le (1+\epsilon)\log n_k
$$
for all sufficiently large $k$. Therefore, almost surely, for all sufficiently large $k$,
$$
\frac{Q_2(t_{n_k})}{\log t_{n_k}} \ge \frac{1-2\epsilon}{(1+\epsilon)\beta}.
$$
Since this holds for every $\epsilon \in (0,1/2)$, we obtain
$$
\limsup_{t \rightarrow \infty} \frac{Q_2(t)}{\log t} \ge \frac{1}{\beta}, \ \ a.s.
$$
From the upper bound in \eqref{BC} and \eqref{SLLN}, we obtain $n_0$ such that for all $n \ge n_0$ 
$$
\frac{\sup_{t \in [\Xi_{n}, \Xi_{n+1}]} Q_2(t)}{\log n} \le \frac{2(1+\epsilon)}{C^*_2 \beta}, \ \ \text{ and }\ \log t \ge (1-\epsilon) \log n.
$$
Therefore,
$$
\frac{Q_2(t)}{\log t} \le \frac{2(1+\epsilon)}{(1-\epsilon)C^*_2 \beta}, \ \ \text{ for all } t \ge \Xi_{n_0}
$$
and hence,
$$
\limsup_{t \rightarrow \infty} \frac{Q_2(t)}{\log t} \le \frac{2}{C^*_2 \beta}, \ \ a.s.
$$
The fluctuation result for $Q_1$ is obtained similarly using Parts (iii) and (iv) of Theorem \ref{th:excren}.
\end{proof}


%% file: asympjsqd.tex
\begin{abstract}
{We consider a system of $N$ identical server pools and a single
dispatcher where tasks with unit-exponential service requirements
arrive at rate $\lambda(N)$.
In order to optimize the experienced performance, the dispatcher aims
to evenly distribute the tasks across the various server pools.
Specifically, when a task arrives, the dispatcher assigns it to the
server pool with the minimum number of tasks among $d(N)$ randomly
selected server pools.

We construct a stochastic coupling to bound the difference in the
system occupancy processes between the JSQ policy and a scheme with
an arbitrary value of $d(N)$.
We use the coupling to derive the fluid limit in case $d(N) \to \infty$
and $\lambda(N)/N \to \lambda$ as $N \to \infty$, along with the
associated fixed point.
The fluid limit turns out to be insensitive to the exact growth rate
of $d(N)$, and coincides with that for the JSQ policy.
We further establish that the diffusion limit corresponds to that
for the JSQ policy as well, as long as $d(N)/\sqrt{N} \log(N) \to \infty$,
and characterize the common limiting diffusion process.
These results indicate that the JSQ optimality can be preserved
at the fluid-level and diffusion-level while reducing the overhead
by nearly a factor O($N$) and O($\sqrt{N}/\log(N)$), respectively.
}
\end{abstract}

\section{Introduction}\label{intro}

In this chapter, we establish asymptotic optimality for a broad
class of randomized load balancing strategies as described in Section~\ref{bloc}.
Specifically, we focus on a basic scenario of $N$~identical
parallel server pools and a single dispatcher where tasks arrive
as a Poisson process.
Incoming tasks cannot be queued, and must immediately be dispatched
to one of the server pools to start execution, or discarded.
Specifically, when a task arrives, the dispatcher assigns it
to a server with the shortest queue among $d(N)$ randomly selected
servers ($1 \leq d(N) \leq N$).
The execution times are assumed to be exponentially distributed,
and do not depend on the number of other tasks receiving service,
but the experienced performance (e.g.~in terms of received throughput
or packet-level delay) does degrade in a convex manner with
an increasing number of concurrent tasks.

The results in this chapter mirror the fluid-level and diffusion-level optimality properties reported in Chapter~\ref{chap:univjsqd}
for power-of-d($N$) strategies in a scenario with single-server queues (not server pools as in this chapter).
Another important difference between the current system and the one studied in Chapter~\ref{chap:univjsqd} is that arriving tasks that are not immediately placed into service can queue in the system of Chapter~\ref{chap:univjsqd}, whereas tasks cannot queue in the system of this chapter (so in case of finite server pools, arriving tasks that are not immediately served are discarded).
Consequently, here the load per server pool $\lambda$ can possibly be larger than~1, in contrast to the $\lambda<1$ assumption in Chapter~\ref{chap:univjsqd}.
As it turns out, due to these differences in the dynamics, a fundamentally different coupling argument is required in the present chapter to establish asymptotic equivalence.
In particular, for the single-server dynamics, first the servers are ordered according to the number of active tasks, and the departures at the ordered servers under two different policies are then coupled.
In contrast, for the infinite-server dynamics, the departure rate at the ordered server pools can vary depending on the exact number of active tasks.
Therefore, the departure processes under two different policies cannot be coupled as in Chapter~\ref{chap:univjsqd}, which necessitates the construction of a novel stochastic coupling.
Specifically, one can think of the coupling for the single-server dynamics as one-dimensional (depending only upon the ordering of the servers), while the coupling we introduce in this chapter is two-dimensional, with the server ordering as one coordinate and the number of tasks as the other, as will be explained in greater detail later.
We further elaborate on the necessity and novelty of the coupling methodology developed in the current chapter, and reflect on the contrast with the stochastic optimality results for the JSQ policy in the
existing literature and the coupling technique in Chapter~\ref{chap:univjsqd} in Remarks~\ref{rem:novelty2} and \ref{rem:contrast2}.
In addition, 
the fluid- and diffusion-limit results in the infinite-server scenario are also notably different from those in Chapter~\ref{chap:univjsqd}.
More specifically, we extend the fluid-limit result in Theorem~\ref{th:genfluid1} in Chapter~\ref{chap:univjsqd} to a more general class of assignment probabilities and departure rate functions, and depending on whether the scaled arrival rate converges to an integer or not,  obtain a qualitatively different behavior of the occupancy state process on diffusion scale.
Furthermore, the diffusion-limit result in Theorem~\ref{diffusionjsqd} characterizes the diffusion-scale behavior only in the transient regime, whereas in the current chapter, since tightness of the diffusion-scaled occupancy process is not an issue (due to the infinite-server dynamics), we are able to analyze the steady-state behavior as well.

The remainder of the chapter is organized as follows.
In Section~\ref{sec: model descr-mor} we present a detailed model description, and provide an overview of the main results. 
In Section~\ref{sec:eqiv} we explain the proof outline and introduce a notion of asymptotic equivalence of two assignment schemes. 
Section~\ref{sec:proof-equiv} introduces a stochastic coupling between any two schemes, and proves the asymptotic equivalence results.
Sections~\ref{sec:fluid-mor}--\ref{sec:integral} contain the proofs of the main results, and in Section~\ref{sec:performance} we reflect upon various performance implications. We conclude in Section~\ref{sec:conclusion-mor} with topics for further research.

\section{Main results}\label{sec: model descr-mor}
\label{sec:main-mor}

\subsection{Model description and notation}
Consider a system with $N$~parallel identical server pools and a single
dispatcher where tasks arrive as a Poisson process of rate~$\lambda(N)$.
Arriving tasks cannot be queued, and must immediately be assigned
to one of the server pools to start execution.
The execution times are assumed to be exponentially distributed with unit
mean, and do not depend on the number of other tasks receiving service.
Each server pool is however only able to accommodate a maximum
of $B$~simultaneous tasks (possibly $B = \infty$),
and when a task is allocated to a server pool that is already handling
$B$~active tasks, it gets permanently discarded.

Specifically, when a task arrives, the dispatcher assigns it to the
server pool with the minimum number of active tasks among $d(N)$
randomly selected server pools ($1 \leq d(N) \leq N$).
As mentioned earlier, this assignment strategy is called a JSQ$(d(N))$
scheme, as it closely resembles the power-of-$d$ version
of the  Join-the-Shortest-Queue (JSQ) policy, and will also
concisely be referred to as such in the special case $d(N) = N$.
We will consider an asymptotic regime where the number of server pools~$N$
and the task arrival rate $\lambda(N)$ grow large in proportion,
with $\lambda(N) / N \to \lambda\leq B$ as $N \to \infty$.
For convenience, we denote $K = \lfloor \lambda \rfloor$
and $f = \lambda - K \in [0, 1)$.

For any $d(N)$ ($1 \leq d(N) \leq N$), let 
$$\QQ^{ d(N)}(t) =
(Q_1^{ d(N)}(t), Q_2^{ d(N)}(t), \dots, Q_B^{ d(N)}(t))$$
be the system occupancy state, where $Q_i^{ d(N)}(t)$ is the number
of server pools under the JSQ($d(N)$) scheme with $i$~or more active
tasks at time~$t$, $i = 1, \dots, B$. 
A schematic diagram of the $Q_i$-values has been provided in Figure~\ref{figB}.
We occasionally omit the superscript $d(N)$, and replace it by~$N$, to refer
to the $N^{\mathrm{th}}$ system, when the value of $d(N)$ is clear from the context.
In case of a finite buffer size $B < \infty$, when a task is discarded,
we call it an \emph{overflow} event, and we denote by $L^{ d(N)}(t)$ the
total number of overflow events under the JSQ($d(N)$) policy up to time~$t$. 

Throughout we assume that at each arrival epoch the server pools are ordered
in  nondecreasing order of the number of active tasks (ties can be broken arbitrarily), recall Figure~\ref{figB} in Chapter~\ref{chap:introduction}, and whenever we refer to some ordered server pool, it should be understood with respect to this prior ordering, unless mentioned otherwise.

\paragraph{Notation.}
Boldfaced letters will be used to denote vectors.
A sequence of random variables $\big\{X_N\big\}_{N\geq 1}$ is said to be $\Op(g(N))$, or $\op(g(N))$, for some function $g:\N\to\R_+$, if the sequence of scaled random variables $\big\{X_N/g(N)\big\}_{N\geq 1}$ is  a tight sequence, or converges to zero in probability, respectively. 
Whenever we mention `with high probability', it should be understood as `with probability tending to 1 as the underlying scaling parameter tends to infinity'.
For stochastic boundedness of a process we refer to \cite[Definition 5.4]{PTRW07}. 
Also, $f$ will be called `diverging to infinity' if $g(N)\to\infty$ as $N\to\infty$.
For any complete separable metric space $E$, denote by $D_E[0,\infty)$, the set of all $E$-valued c\`adl\`ag (right continuous with left limits exist) processes.
By the symbol `$\dto$' we denote convergence in distribution for real-valued random variables, and  with respect to Skorohod-$J_1$ topology for \emph{c\'adl\'ag} processes.


\subsection{Fluid-limit results}\label{ssec:fluid}

In order to state the fluid-limit results, we first introduce some
useful notation.
Denote the fluid-scaled system occupancy state by
$\qq^{ d(N)}(t) := \QQ^{ d(N)}(t) / N$.
We will denote by $\tilde{S}=\big\{\QQ\in\Z^B:Q_i \leq Q_{i-1} \mbox{ for all } i = 2, \dots, B\big\}$ and $S =\big\{\qq \in [0, 1]^B: q_i \leq q_{i-1} \mbox{ for all } i = 2, \dots, B\big\}$
 the set of all possible unscaled and fluid-scaled occupancy states, respectively.
Further define $S^N:= S\cap\big\{i/N:1\leq i\leq N\big\}^B$ as the space of all fluid-scaled occupancy states of the $N^\mathrm{th}$ system.
We take the following product norm on $S$: for $\mathbf{q}_1=(q_{1,1},q_{1,2},\ldots, q_{1,B})$, $\mathbf{q}_2=(q_{2,1},q_{2,2},\ldots, q_{2,B})\in \R^B$,
$$\rho(\mathbf{q}_1,\mathbf{q}_2):=\sum_{i=1}^{B}\frac{|q_{1,i}-q_{2,i}|\wedge 1}{2^i},$$
and  all the convergence results below will be with respect to product topology.
We often write $\rho(\qq_1,\qq_2)$ as $\norm{\qq_1-\qq_2}$. 
Let $(E,\hat{\rho})$ be a metric space.
We call a function $g:S\to E$  Lipschitz continuous on $S$, if there exists $L>0$, such that for all $x,y\in S,$
$$\hat{\rho}(g(x),g(y))\leq L \hat{\rho}(x,y).$$
For any $\qq \in S$, denote by $m(\qq) = \min\big\{i: q_{i + 1} < 1\big\}$
the minimum number of active tasks among all server pools, with the convention that
$q_{B+1} = 0$ if $B < \infty$.
If $m(\qq)=0$, then define $p_0(m(\qq))=1$ and $p_i(m(\qq))=0$ for all $i=1,2,\ldots$. 
Otherwise, in case $m(\qq)>0$, we distinguish two cases, depending on whether the normalized arrival
rate $\lambda$ is larger than $m(\qq)(1 - q_{m(\qq) + 1})$ or not.
If $\lambda \leq m(\qq) (1 - q_{m(\qq) + 1})$, then define 
$$p_{m(\qq) - 1}(\qq) = 1,\quad\mbox{and}\quad p_i(\qq) = 0\quad\mbox{for all}\quad i \neq m(\qq) - 1.$$
On the other hand, if $\lambda >m(\qq) (1 - q_{ m(\qq) + 1})$,
then 
\begin{equation}
p_{i}(\qq)=
\begin{cases}
m(\qq)(1 - q_{ m(\qq) + 1})/\lambda & \quad\mbox{ for }\quad i=m(\qq)-1,\\
1 - p_{ m(\qq) - 1}(\qq) & \quad\mbox{ for }\quad i=m(\qq),\\
0&\quad \mbox{ otherwise.}
\end{cases}
\end{equation}
Note that the assumption $\lambda \leq B$ ensures that the latter case
cannot occur when $B<\infty$ and $m(\qq) = B$.

\begin{theorem}[{Universality of fluid limit for JSQ($d(N)$) scheme}]
\label{fluidjsqd-mor}
Assume that the initial occupancy state $\qq^{ d(N)}(0)$ converges to $\qq^\infty \in S$ as $N \to \infty$.
For the JSQ($d(N)$) scheme with $d(N)$ diverging to infinity, with probability~1, any subsequence of $\{N\}$ has a further subsequence along which on any finite time interval, the sequence of processes
$\big\{\qq^{ d(N)}(t)\big\}_{t \geq 0}$ converges to some deterministic trajectory $\big\{\qq(t)\big\}_{t \geq 0}$ that
satisfies the system of integral equations
\begin{equation}\label{eq:fluidjsqd-mor}
q_i(t) = q_i(0)+
\lambda\int_0^t p_{i-1}(\qq(s))\dif s - i\int_0^t (q_i(s) - q_{i+1}(s))\dif s, \quad i=1,\ldots, B,
\end{equation}
where $\qq(0) = \qq^\infty$ and the coefficients $p_i(\cdot)$ are
as defined above.
\end{theorem}

The above theorem shows that the fluid-level dynamics do not depend
on the specific growth rate of $d(N)$ as long as $d(N) \to\infty$
as $N \to \infty$.
In particular, the JSQ$(d(N))$ scheme with $d(N) \to\infty$ as $N \to\infty$ exhibits
the same behavior as the ordinary JSQ policy, and thus achieves
fluid-level optimality. This result can be intuitively interpreted as follows. Since $d(N)$ is growing, for large $N$, at an arrival epoch, if the fraction of server pools with the minimum number of active tasks becomes positive, then with high probability at least one of the $d(N)$ selected server pools will be from the ones with the minimum number of active tasks.
This ensures that as long as $d(N)\to\infty$ as $N \to\infty$, the difference in $Q_i$-values between the ordinary JSQ policy and the JSQ$(d(N))$ scheme can not become $O(N)$, yielding fluid-level optimality.

The coefficient $p_i(\qq)$ represents the fraction
of incoming tasks assigned to server pools with exactly~$i$ active
tasks in the fluid-level state $\qq \in S$.
Assuming $m(\qq) < B$, a strictly positive fraction $1 - q_{m(\qq) + 1}$
of the server pools have exactly $m(\qq)$ active tasks.
Since $d(N)\to\infty$ as $N \to\infty$, the fraction of incoming tasks that get assigned
to server pools with $m(\qq) + 1$ or more active tasks is therefore zero:
$p_i(\qq) = 0$ for all $i = m(\qq) + 1, \dots, B - 1$.
Also, tasks at server pools with exactly~$i$ active tasks are completed
at (normalized) rate $i(q_i - q_{i + 1})$, which is zero for all
$i = 1, \dots, m(\qq) - 1$, and hence the fraction of incoming tasks
that get assigned to server pools with $m(\qq) - 2$ or less active
tasks is zero as well: $p_i(\qq) = 0$ for all $i = 0, \dots, m(\qq) - 2$.
This only leaves the fractions $p_{m(\qq) - 1}(\qq)$
and $p_{m(\qq)}(\qq)$ to be determined.
Now observe that the fraction of server pools with exactly $m(\qq) - 1$
active tasks is zero.
However, since tasks at server pools with exactly $m(\qq)$ active tasks
are completed at (normalized) rate $m(\qq) (1 - q_{m(\qq) + 1}) > 0$,
incoming tasks can be assigned to server pools with exactly $m(\qq) - 1$
active tasks at that rate.
We thus need to distinguish between two cases, depending on whether
the normalized arrival rate $\lambda$ is larger than
$m(\qq) (1 - q_{m(\qq) + 1})$ or not.
If $\lambda \leq m(\qq) (1 - q_{m(\qq) + 1})$, then all the incoming tasks
can be assigned to server pools with exactly $m(\qq) - 1$ active tasks,
so that $p_{m(\qq) - 1}(\qq) = 1$ and $p_{m(\qq)}(\qq) = 0$.
On the other hand, if $\lambda > m(\qq )(1 - q_{m(\qq) + 1})$, then not
all incoming tasks can be assigned to server pools with exactly
$m(\qq) - 1$ active tasks, and a positive fraction will be assigned
to server pools with exactly $m(\qq)$ active tasks:
$p_{m(\qq) - 1}(\qq) = m(\qq) (1 - q_{m(\qq) + 1}) / \lambda$
and $p_{m(\qq)}(\qq) = 1 - p_{m(\qq) - 1}(\qq)$.

It is easily verified that the unique fixed point of the differential
equation in Theorem~\ref{fluidjsqd-mor} is given by
\begin{equation}
\label{eq:fixed point}
q_i^\star = \left\{\begin{array}{ll} 1 & i = 1, \dots, K, \\
f & i = K + 1, \\
0 & i = K + 2, \dots, B, \end{array} \right.
\end{equation}
and thus $\sum_{i=1}^B q_i^\star =\lambda$.
This is consistent with the results in Mukhopadhyay
{\em et al.}~\cite{MKMG15,MMG15} and Xie {\em et al.}~\cite{XDLS15}
for fixed~$d$, where taking $d \to \infty$ yields the same fixed point.
However, the results in \cite{MKMG15,MMG15,XDLS15} for fixed~$d$
cannot directly be used to handle joint scalings, and do not yield the
universality of the entire fluid-scaled sample path
for arbitrary initial states as established in Theorem~\ref{fluidjsqd-mor}.

Having obtained the fixed point of the fluid limit, we now establish the interchange of the mean-field 
$(N\to\infty)$ and stationary $(t\to\infty)$ limits. 
The fixed point in~\eqref{eq:fixed point} in conjunction with the interchange of limits result in Proposition~\ref{prop:interchange-mor} below indicates that in stationarity the fraction of servers with at least $K+2$ and at most $K-1$ active tasks is negligible.
Let 
$$\pi^{ d(N)}(\cdot)=\lim_{t\to\infty}\Pro{\qq^{ d(N)}(t)=\cdot}$$ 
be the stationary measure of the occupancy states of the $N^{\mathrm{th}}$ system.  

\begin{proposition}[{Interchange of limits}]
\label{prop:interchange-mor}
Let $d(N)\to\infty$ as $N \to\infty$. Then
the sequence of stationary measures $\big\{\pi^{ d(N)}\big\}_{N\geq 1}$   converges weakly to $\pi^\star$, where $\pi^\star=\delta_{\qq^\star}$ with $\delta_x$ being the Dirac measure concentrated upon $x$, and $\qq^\star$ defined by~\eqref{eq:fixed point}.
\end{proposition}
The above proposition relies on tightness of $\big\{\pi^{d(N)}\big\}_{N\geq 1}$ and the global stability of the fixed point, and is proved in Subsection~\ref{ssec:globstab-mor}.

\subsection{Diffusion-limit results for non-integral \texorpdfstring{$\boldsymbol{\lambda}$}{lambda}}

As it turns out, the diffusion-limit results may be qualitatively
different, depending on whether $f = 0$ or $f > 0$,
and we will distinguish between these two cases accordingly.
Observe that for any assignment scheme, in the absence of overflow events, the total number of active
tasks evolves as the number of jobs in an M/M/$\infty$ system
with arrival rate $\lambda(N)$ and unit service rate,
for which the diffusion limit is well-known~\cite{Robert03}.
For the JSQ$(d(N))$ scheme with $d(N)/(\sqrt{N} \log(N))\to\infty$ as $N\to\infty$, we can
establish, for suitable initial conditions, that the total number of server
pools with $K - 2$ or less and $K + 2$ or more tasks is negligible
on the diffusion scale.
When $f > 0$, the number of server pools with $K - 1$ tasks is negligible
as well, and the dynamics of the number of server pools with $K$
or $K + 1$ tasks can then be derived from the known diffusion limit
of the total number of tasks mentioned above.
In contrast, when $f = 0$, the number of server pools with $K - 1$ tasks
is not negligible on the diffusion scale, and the limiting behavior is
qualitatively different, but can still be characterized.

We first consider the case $f > 0$, and define $f(N):= \lambda(N)-KN.$
Based on the above observations,
we define the following centered and scaled processes:
\begin{equation}
\begin{split}
\bar{Q}_i^{ d(N)}(t) &:= \dfrac{N-Q^{ d(N)}_i(t)}{\sqrt{N}}\geq 0,\quad i \leq K,\\  \\
\bar{Q}_{K+1}^{ d(N)}(t) &:= \dfrac{Q^{ d(N)}_{K+1}(t)-f(N)}{\sqrt{N}}\in\R,\\ \\
\bar{Q}_i^{ d(N)}(t) &:=\frac{Q^{ d(N)}_i(t)}{\sqrt{N}}\geq 0,\quad\text{for}\quad i \geq K+2.
\end{split}
\end{equation}
\begin{theorem}[{Universality of diffusion limit for JSQ($d(N)$) scheme, $f > 0$}]
\label{th:diff pwr of d 1}
If $f > 0$, $\bQ^{ d(N)}_{K+1}(0)\to\bQ_{K+1}\in\R$, $\bQ^{ d(N)}_i(0)\to 0$ for $i\neq K+1$, and $d(N)$ is such that $d(N)/(\sqrt{N} \log(N))\to\infty$ as $N\to\infty$, then the following holds
as $N \to \infty$:
\begin{enumerate}[{\normalfont(i)}]
\item For $i \leq K$, $\big\{\bar{Q}_i^{ d(N)}(t)\big\}_{t \geq 0} \dto
\big\{\bar{Q}_i(t)\big\}_{t \geq 0}$, where $\bar{Q}_i(t) \equiv 0$.
\item $\big\{\bar{Q}_{K+1}^{ d(N)}(t)\big\}_{t \geq 0} \dto 
\big\{\bar{Q}_{K+1}(t)\big\}_{t \geq 0}$, where $\bar{Q}_{K+1}(t)$ is given
by the Ornstein-Uhlenbeck process satisfying the stochastic
differential equation
\begin{equation}
\label{eq:OU process1}
\dif \bar{Q}_{K+1}(t) =
- \bar{Q}_{K+1}(t) \dif t + \sqrt{2 \lambda} \dif W(t),
\end{equation}
where $W(t)$ is the standard Brownian motion.
\item For $i \geq K+2$, $\big\{\bar{Q}_i^{ d(N)}(t)\big\}_{t \geq 0} \dto 
\big\{\bar{Q}_i(t)\big\}_{t \geq 0}$, where $\bar{Q}_i(t) \equiv 0$.
\end{enumerate}

\end{theorem}
Loosely speaking, the above theorem says that, if $f > 0$
and $d(N)/(\sqrt{N} \log(N))\to\infty$ as $N\to\infty$, then over any finite time horizon,
there will only be $o_P(\sqrt{N})$ server pools with fewer than
$K$ or more than $K+1$~active tasks, and $fN+O_P(\sqrt{N})$ server pools with precisely
$K + 1$ active tasks.
Also, as long as $d(N)/(\sqrt{N} \log(N))\to\infty$ as $N\to\infty$, the JSQ$(d(N))$
scheme exhibits the same behavior as the ordinary JSQ policy (i.e., $d(N)=N$),
and thus achieves diffusion-level optimality. 
The result can be heuristically explained as follows. 
When the number of server pools with $K-1$ or less number of active tasks is $\Theta(\sqrt{N})$, the JSQ$(d(N))$ scheme should be able to assign the incoming tasks
with high probability to one of those server pools. 
To be able to select one of the $\Theta(\sqrt{N})$ server pools out of $N$ server pools, $d(N)$ must grow faster than $\sqrt{N}$. 
Now further observe that in any finite time interval there are on average $\Theta(N)$ arrivals, and hence it is not enough to assign the incoming task to the appropriate server pool only once. The number of times that the JSQ$(d(N))$ scheme fails to assign a task to the `appropriate' server pool in any finite time interval, should be $\op(\sqrt{N})$. This gives rise to the additional $\log( N)$ factor in the growth rate of $d(N)$.

\subsection{Diffusion-limit results for integral \texorpdfstring{$\boldsymbol{\lambda}$}{lambda}}
We now turn to the case $f = 0$, and assume that
\begin{equation}
\label{eq:f=0}
\frac{K N - \lambda(N)}{\sqrt{N}} \to \beta \in \R \quad\mbox{ as }\quad N \to \infty,
\end{equation} 
which can be thought of as an analog of the so-called Halfin-Whitt
regime~\cite{HW81}.
As mentioned above, the limiting behavior in this case is
qualitatively different from the case $f > 0$.
Hence, we now consider the following scaled quantities:
\begin{equation}\label{eq:scaling-f=0}
\begin{split}
\hQ_{K-1}^{ d(N)}(t)&:=  \sum_{i=1}^{K-1}\dfrac{ N-Q_i^{ d(N)}(t)}{\sqrt{N}}\geq 0,\\\\
\hQ_K^{ d(N)}(t)&:=\dfrac{ N-Q_K^{ d(N)}(t)}{\sqrt{N}}\geq 0,\\\\
\hQ_i^{ d(N)}(t)&:=\dfrac{Q_i^{ d(N)}(t)}{\sqrt{N}}\geq 0,\quad \mathrm{for}\ i\geq K+1.
\end{split}
\end{equation}

\begin{theorem}[{Universality of diffusion limit for JSQ($d(N)$) scheme, $f = 0$}]
\label{th:diff pwr of d 2}
Suppose there exists $M\geq K+1$, such that $Q^{ d(N)}_{M+1}(0)\equiv 0$, and 
$$(\hQ^{ d(N)}_{K-1}(0),\hQ^{ d(N)}_K(0),\ldots,\hQ^{ d(N)}_M(0))\dto (\hQ_{K-1}(0),\hQ_K(0),\ldots,\hQ_M(0))$$ in $\R^{ M-K+2}$. 
If $f = 0$, $d(N)/(\sqrt{N} \log(N))\to\infty$, Equation~\eqref{eq:f=0}
is satisfied, and $\hQ^{ d(N)}_{K-1}(0)\pto 0$, as $N\to\infty$, then the process 
$$\left\{\big(\hQ^{ d(N)}_{K-1}(t),\hQ^{ d(N)}_K(t),\ldots,\hQ^{ d(N)}_M(t),\hQ^{ d(N)}_{M+1}(t)\big)\right\}_{t\geq 0}$$ 
converges weakly to the process defined as the unique solution to the stochastic integral equation
\begin{equation}\label{eq:OU process2}
\begin{split}
\hQ_K(t) &= \hQ_K(0) + \sqrt{2K} W(t) -
\int_0^t (\hQ_K(s) + K \hQ_{K+1}(s)) \dif s + \beta t + V_1(t), \\
\hQ_{K+1}(t) &= \hQ_{K+1}(0) + V_1(t) - (K + 1) \int_0^t (\hQ_{K+1}(s)-\hQ_{K+2}(s)) \dif s,\\
\hQ_{i}(t) &= \hQ_{i}(0)  - i \int_0^t (\hQ_{i}(s)-\hQ_{i+1}(s)) \dif s, \quad i= K+2.\ldots, M-1,\\
\hQ_{M}(t) &= \hQ_{M}(0) - M \int_0^t \hQ_{M}(s) \dif s,
\end{split}
\end{equation}
$\hQ_{K-1}(t)\equiv 0$, and $\hQ_{M+1}(t)\equiv 0$, where $W(t)$ is the standard Brownian motion, and $V_1(t)$ is the unique
non-decreasing process in $D_{\R_+}[0,\infty)$ satisfying
$$\int_0^t \ind{\hQ_K(s) \geq 0} \dif V_1(s) = 0.$$
\end{theorem}
Unlike the $f > 0$ case, the above theorem says that, if $f = 0$,
then over any finite time horizon, there will be $O_P(\sqrt{N})$
server pools with fewer than $K$ or more than $K$~active tasks,
and hence most of the server pools have precisely $K$~active tasks.


\section{Proof outline}\label{sec:eqiv}
The proofs of the asymptotic results for the JSQ$(d(N))$ scheme in Theorems~\ref{fluidjsqd-mor}, \ref{th:diff pwr of d 1}, and \ref{th:diff pwr of d 2} involve two main components:
\begin{enumerate}[{\normalfont (i)}]
\item deriving the relevant limiting processes for the ordinary JSQ policy,
\item establishing a universality result which shows that the limiting processes for the JSQ$(d(N))$ scheme are `asymptotically equivalent' to those for the ordinary JSQ policy for suitably large~$d(N)$.
\end{enumerate}
For Theorems~\ref{fluidjsqd-mor}, \ref{th:diff pwr of d 1} and \ref{th:diff pwr of d 2}, part (i) will be dealt with in Theorems~\ref{th:genfluid},~\ref{th:diffusion} and~\ref{th: f=0 diffusion-mor}, respectively. 
For all three theorems, part (ii) relies on a notion of asymptotic equivalence between different schemes, which is formalized in the next definition. 
\begin{definition}
Let $\Pi_1$ and $\Pi_2$ be two schemes parameterized by the number of server pools $N$. For any positive function $g:\N\to\R_+$, we say that $\Pi_1$ and $\Pi_2$ are `$g(N)$-alike' if 
there exists a common probability space, such that for any fixed $T\geq 0$, for all $i\geq 1$,
$$\sup_{t\in[0,T]}(g(N))^{-1}|Q_i^{\Pi_1}(t)-Q_i^{\Pi_2}(t)|\pto 0\quad \mathrm{as}\quad N\to\infty.$$
\end{definition}
Intuitively speaking, if two schemes are $g(N)$-alike, then in some sense, the associated system occupancy states are indistinguishable on the $g(N)$-scale. 
For brevity, for two schemes $\Pi_1$ and $\Pi_2$ that are $g(N)$-alike, we will often say that $\Pi_1$ and $\Pi_2$ have the same process-level limits on the $g(N)$-scale.
The next theorem states a sufficient criterion for the JSQ$(d(N))$ scheme and the ordinary JSQ policy to be $g(N)$-alike, and thus, provides the key vehicle in establishing the universality result in part (ii) mentioned above. 
\begin{theorem}\label{th:pwr of d}
Let $g:\N\to\R_+$ be a function diverging to infinity. Then the JSQ policy and the JSQ$(d(N))$ scheme are $g(N)$-alike, with $g(N)\leq N$, if 
\begin{align}\label{eq:fNalike cond1}
\mathrm{(i)}&\quad d(N)\to\infty, \quad\text{for}\quad g(N) = O(N),\\\nonumber\\
\mathrm{(ii)}&\quad d(N)\left(\frac{N}{g(N)}\log\left(\frac{N}{g(N)}\right)\right)^{-1}\to\infty,\quad\text{for}\quad g(N)=o(N).\label{eq:fNalike cond2}
\end{align}
\end{theorem}
Theorem~\ref{th:pwr of d} can be intuitively explained as follows. 
The choice of $d(N)$ should be such that the JSQ$(d(N))$ scheme, at each arrival, with high probability selects one of the server pools with the minimum number of tasks, if the total number of server pools with the minimum number of tasks is of order $g(N)$. 
Moreover, in any finite time interval, the total number of times it fails to do so, should be of order lower than that of $g(N)$. 
These conditions imply that $d(N)$ must diverge if $g(N)=O(N)$, or grow faster than $(N/g(N))\log(N/g(N))$, if $g(N)=o(N)$. 

In order to obtain the fluid and diffusion limits for various schemes, the two main scales that we consider are $g(N)\sim N$ and $g(N)\sim\sqrt{N}$, respectively. 
The next two immediate corollaries of Theorem~\ref{th:pwr of d} will imply that it is enough to investigate the ordinary JSQ policy in various regimes.
\begin{corollary}\label{cor:fluid}
If $d(N)\to\infty$ as $N\to\infty$, then the JSQ$(d(N))$ scheme and the ordinary JSQ policy are $N$-alike.
\end{corollary}
\begin{remark}\label{rem:necessity-fluid}
\normalfont
The growth condition on $d(N)$ in order for the JSQ$(d(N))$ scheme to be $N$-alike to the ordinary JSQ policy, stated in the above corollary, is not only sufficient, but also necessary.
Specifically, if $\liminf_{N\to\infty}d(N)= d<\infty$, then 
consider a subsequence along which the limit of $d(N)$ exists and is uniformly bounded by $d$.
Therefore, one can choose a further subsequence, such that $d(N)=d$ for all $N$ along that subsequence.
Now, from the fluid-limit result for the JSQ$(d)$ scheme~\cite{MKMG15, MMG15},
one can see that it differs from that of the JSQ policy stated in~\eqref{fluidjsqd-mor}, and hence the JSQ($d(N)$) scheme is not $N$-alike to the ordinary JSQ policy.
\end{remark}

\begin{corollary}\label{cor-diff}
If $d(N)/(\sqrt{N}\log (N))\to\infty$ as $N\to\infty$, then the JSQ$(d(N))$ scheme and the ordinary JSQ policy are $\sqrt{N}$-alike.
\end{corollary}

We will prove the universality result in Theorem~\ref{th:pwr of d} in the next section.
The key challenge is that
a direct comparison of the JSQ$(d(N))$ scheme and the ordinary JSQ policy is not straightforward. 
Hence, to compare the JSQ$(d(N))$ scheme with the JSQ policy, we adopt a two-stage approach based on a novel class of schemes, called CJSQ$(n(N))$, as a convenient intermediate scenario. 
Specifically, for some nonnegative integer-valued sequence $\big\{n(N)\big\}_{N\geq 1}$, with $n(N)\leq N$, we introduce a class of schemes named CJSQ($n(N)$), containing all the schemes that always assign the incoming task to one of the  $n(N)+1$ lowest ordered server pools.
Note that when $n(N)=0$, the class only contains the ordinary JSQ policy. 

Just like the JSQ$(d(N))$ scheme, the schemes in the class CJSQ$(n(N))$ may be thought of as ``sloppy'' versions of the JSQ policy, in the sense that tasks are not necessarily assigned to a server pool with the minimum number of active tasks but to one of the $n(N)+1$
lowest ordered server pools, as graphically illustrated in Figure~\ref{fig:sfigCJSQ} in Chapter~\ref{chap:introduction}.
Below we often will not differentiate among the various schemes in the class CJSQ$(n(N))$, and prove a common property possessed by all these schemes. Hence, with a minor abuse of notation, we will often denote a typical assignment scheme in this class by CJSQ($n(N)$).
Note that the JSQ$(d(N))$ scheme is guaranteed to identify the lowest ordered server pool, but only among a randomly sampled subset of $d(N)$ server pools.
In contrast, a scheme in the class in CJSQ$(n(N))$ only guarantees that
one of the $n(N)+1$ lowest ordered server pools is selected, but 
across the entire system of $N$ server pools. 
We will show that for sufficiently small $n(N)$, any scheme from the class CJSQ$(n(N))$ is still `close' to the ordinary JSQ policy in terms of $g(N)$-alikeness as stated in the next proposition.
\begin{proposition}\label{prop: modified JSL}
For any function $g:\N\to\R_+$ diverging to infinity, if 
$$n(N)/ g(N)\to 0\quad \mbox{as}\quad N\to\infty,$$ then the JSQ policy and the CJSQ$(n(N))$ schemes are $g(N)$-alike.
\end{proposition}
In order to prove this proposition, we introduce in Section~\ref{sec:stoch-coupling} a novel stochastic coupling called the T-coupling, to construct a common probability space, and establish the property of $g(N)$-alikeness.\\

Next we compare the CJSQ($n(N)$) schemes with the JSQ($d(N)$) scheme. 
The comparison follows a somewhat similar line of argument as 
in Section~\ref{sec:coupling} in Chapter~\ref{chap:univjsqd}, and involves a JSQ$(n(N),d(N))$ scheme
which is an intermediate blend between the CJSQ$(n(N))$ schemes and the
JSQ$(d(N))$ scheme.
Specifically, the scheme JSQ$(n(N),d(N))$ selects a candidate server pool in the exact same way as JSQ$(d(N))$.
However, it only assigns the task to that server pool if it belongs to 
the $n(N)+1$ lowest ordered ones,
and to a randomly selected server pool among these otherwise.
By construction, the JSQ$(n(N),d(N))$ scheme belongs to the class CJSQ$(n(N))$.

Next consider two T-coupled systems with a JSQ$(d(N))$
and a JSQ$(n(N),d(N))$ scheme.
Assume that at some specific arrival epoch, the incoming task is assigned to the $k^{\mathrm{th}}$ ordered server pool in the system under the JSQ($d(N)$) scheme. 
If $k\in\big\{1,2,\ldots,n(N)+1\big\}$, then the scheme JSQ$(n(N),d(N))$ also assigns the arriving task to the $k^{\mathrm{th}}$ ordered server pool. 
Otherwise it dispatches the arriving task uniformly at random among the first $n(N)+1$ ordered server pools.

We will establish a sufficient criterion on $d(N)$ in order for the JSQ$(d(N))$ scheme and JSQ$(n(N),d(N))$ scheme to be close in terms of $g(N)$-alikeness, as stated in the next proposition.
\begin{proposition}\label{prop: power of d}
Assume that $n(N)/g(N)\to 0$ as $N\to\infty$ for some function $g:\N\to\R_+$ diverging to infinity. The JSQ$(n(N),d(N))$ scheme and the JSQ($d(N)$) scheme are $g(N)$-alike if the following condition holds:
\begin{equation}\label{eq:condition-same}
\frac{n(N)}{N}d(N)-\log\left(\frac{N}{g(N)}\right)\to\infty, \quad\text{as}\quad N\to\infty.
\end{equation}
\end{proposition} 
Finally, Proposition~\ref{prop: power of d} in conjunction with Proposition~\ref{prop: modified JSL} yields Theorem~\ref{th:pwr of d}.
The overall proof strategy as described above, has been schematically represented in Figure~\ref{fig:sfigrelation-inf} in Chapter~\ref{chap:introduction}.

\begin{remark} \normalfont
Note that sampling \emph{without} replacement polls more server pools than \emph{with} replacement, and hence the minimum number of active tasks among the selected server pools is stochastically smaller in the case without replacement.
As a result, for sufficient conditions as in Theorem~\ref{th:pwr of d} it is enough to consider sampling with replacement.
\end{remark}

\section{Universality property}
\label{sec:proof-equiv}

In this section we formalize the proof outlined in the previous section.
In Subsection~\ref{sec:stoch-coupling} we first introduce the T-coupling between any two task assignment schemes.
This coupling is used to derive stochastic inequalities in Subsection~\ref{ssec:stoch-ineq}, stated as Proposition~\ref{prop:stoch-ord2-mor} and Lemma~\ref{lem:majorization-mor}, which
in turn,  are used to prove Propositions~\ref{prop: modified JSL} and~\ref{prop: power of d} and Theorem~\ref{th:pwr of d} in Subsection~\ref{ssec:asymp-eq}.
\subsection{Stochastic coupling}\label{sec:stoch-coupling}
Throughout this subsection we fix $N$, and suppress the superscript $N$ in the notation. Let $Q_{i}^{\Pi_1}(t)$ and $Q_{i}^{\Pi_2}(t)$ denote the number of server pools with at least $i$ active tasks, at time $t$, in two systems following schemes $\Pi_1$ and $\Pi_2$, respectively.
With a slight abuse of terminology, we occasionally use $\Pi_1$ and $\Pi_2$ to refer to systems following schemes $\Pi_1$ and $\Pi_2$, respectively.
To couple the two systems, we synchronize the arrival epochs 
and maintain a single exponential departure clock with instantaneous rate at time $t$ given by $M(t):=\max\left\{\sum_{i=1}^BQ_{i}^{\Pi_1}(t), \sum_{i=1}^BQ_{i}^{\Pi_2}(t)\right\}$.
 We couple the arrivals and departures in the various server pools as follows:

(1) \emph{Arrival:} At each arrival epoch, assign the incoming task in each system to one of the server pools according to the respective schemes.

(2) \emph{Departure:} 
Define
$$H(t):=\sum_{i=1}^{B}\min\left\{Q_{i}^{\Pi_1}(t), Q_{i}^{\Pi_2}(t)\right\}$$
and
$$p(t):=
\begin{cases}
\dfrac{H(t)}{M(t)},&\quad\text{if}\quad M(t)>0,\\
0,&\quad \text{otherwise.}
\end{cases}
$$
At each departure epoch $t_k$ (say), draw a uniform$[0,1]$ random variable $U(t_k)$. 
The departures occur in a coupled way based upon the value of $U(t_k)$. In either of the systems, assign a task index $(i,j)$, if that task is at the $j^{\mathrm{th}}$ position of the $i^{\mathrm{th}}$ ordered server pool. 
Let $\mathcal{A}_1(t)$ and $\mathcal{A}_2(t)$ denote the set of all task-indices present at time $t$ in systems $\Pi_1$ and $\Pi_2$, respectively. 
Color the indices (or tasks) in $\mathcal{A}_1\cap\mathcal{A}_2$, $\mathcal{A}_1\setminus\mathcal{A}_2$ and $\mathcal{A}_2\setminus\mathcal{A}_1$, green, blue and red, respectively,
and note that $|\mathcal{A}_1\cap \mathcal{A}_2|=H(t)$. 
Define a total order on the set of indices as follows:
$(i_1,j_1)<(i_2,j_2)$ if $i_1<i_2$, or $i_1=i_2$ and $j_1<j_2$.  Now, if $U(t_k)\leq p(t_k-)$, then select one green index uniformly at random and remove the corresponding tasks from both systems. Otherwise, if $U(t_k)> p(t_k-)$, then choose one integer $m$, 
uniformly at random from all the integers between 1 and  $M(t)-H(t)=M(t)(1-p(t))$, and remove the tasks corresponding to the $m^{\mathrm{th}}$ smallest (according to the order defined above) red and blue indices in the corresponding systems.
If the number of red (or blue) tasks is less than $m$, then do nothing.


The above coupling has been schematically represented in Figure~\ref{fig:tcoupling} in Chapter~\ref{chap:introduction},
and will henceforth be referred to as T-coupling, where T stands for `task-based'. 
We need to show that, under the T-coupling, the two systems, considered independently, evolve according to their own statistical laws. 
This can be seen in several steps.
Indeed, the T-coupling basically uniformizes the departure rate by the maximum number of tasks present in either of the two systems. 
Then informally speaking, the green region signifies the common portion of tasks, and the red and  blue regions represent the separate contributions. 
Now observe that
\begin{enumerate}[{\normalfont (i)}]
\item The total departure rate from $\Pi_i$ is 
\begin{align*}
&M(t)\left[p(t)+(1-p(t))\frac{|\mathcal{A}_i\setminus\mathcal{A}_{3-i}|}{M(t)-H(t)}\right]
= |\mathcal{A}_1\cap\mathcal{A}_2|+|\mathcal{A}_i\setminus\mathcal{A}_{3-i}|
=|\mathcal{A}_i|,
\end{align*}
$i= 1,2$.
\item Assuming without loss of generality~$|\mathcal{A}_1|\geq |\mathcal{A}_2|$, each task in $\Pi_1$ is equally likely to depart.
\item Each task in $\Pi_2$ within $\mathcal{A}_1\cap\mathcal{A}_2$ and 
each task within $\mathcal{A}_2\setminus\mathcal{A}_1$ is equally likely to depart, and the probabilities of departures are proportional to $|\mathcal{A}_1\cap\mathcal{A}_2|$ and $|\mathcal{A}_2\setminus\mathcal{A}_1|$, respectively.
\end{enumerate}

\subsection{Stochastic inequalities}\label{ssec:stoch-ineq}

Now, as in Chapter~\ref{chap:univjsqd} we define a notion of comparison between two T-coupled systems. 
Two T-coupled systems are said to~\emph{differ in decision} at some arrival epoch, if the index of the ordered server pool joined by the arriving task at that epoch, differs in the two systems.
Denote by $\Delta_{\Pi_1,\Pi_2}(t)$, the cumulative number of times that the two systems~$\Pi_1$ and~$\Pi_2$ differ in decision up to time $t$. 

\begin{proposition}\label{prop:stoch-ord2-mor}
For two T-coupled systems under any two schemes $\Pi_1$ and $\Pi_2$ the following inequality is preserved
\begin{equation}\label{eq:stoch-ord2-mor}
\sum_{i=1}^B \big|Q_i^{\Pi_1}(t)-Q_i^{\Pi_2}(t)\big|\leq 2\Delta_{\Pi_1,\Pi_2}(t)\qquad\forall\ t\geq 0,
\end{equation} 
provided the two systems start from the same occupancy state at time $t=0$.
\end{proposition}
The proof follows a somewhat similar line of argument as in Chapters~\ref{chap:univjsqd} and~\ref{chap:jiq}, but is provided below since the coupling is different here. 
For any scheme $\Pi$, define $I_\Pi(c):=\max\big\{i:Q_i^\Pi\geq N-c+1\big\}$, $c=1,\ldots,N$.\\
\begin{proof}[Proof of Proposition~\ref{prop:stoch-ord2-mor}.]
We use forward induction on event times, i.e., time epochs when either an arrival or a departure takes place. 
Assume that the inequality in~\eqref{eq:stoch-ord2-mor} holds at time epoch $t_0$. 
We 
denote by $\tilde{Q}^{\Pi}$ the updated occupancy state after the next event at time epoch $t_1$, and
 distinguish between two cases depending on whether $t_1$ is an arrival epoch or a departure epoch.

If $t_1$ is an arrival epoch and
 if the systems differ in decision, then observe that the left side of \eqref{eq:stoch-ord2-mor} can increase at most by two. In this case, the right side also increases by two, and the ordering is preserved.
Therefore, it is enough to prove that the right side of \eqref{eq:stoch-ord2-mor} remains unchanged if the two systems do not differ in decision. In that case, assume that both $\Pi_1$ and $\Pi_2$ assign the arriving task to the $k^{\mathrm{th}}$ ordered server pool.
Then
\begin{equation}\label{eq:addition-mor}
\tilde{Q}^{\Pi}_i=
\begin{cases}
Q^{\Pi}_i+1, &\mbox{ for }i=I_{\Pi}(k)+1,\\
Q^{\Pi}_j,&\mbox{ otherwise, }
\end{cases}
\end{equation}
if $I_{\Pi}(k)<B$; otherwise all the $Q_i$-values remain unchanged. 
If $I_{\Pi_1}(k)=I_{\Pi_2}(k)$, then the left side of \eqref{eq:stoch-ord2-mor} clearly remains unchanged. Now, without loss of generality, assume $I_{\Pi_1}(k)<I_{\Pi_2}(k)$. Therefore, 
$$Q_{ I_{\Pi_1}(k)+1}^{\Pi_1}(t_0)< Q_{ I_{\Pi_1}(k)+1}^{\Pi_2}(t_0)\quad \mathrm{and}\quad Q_{ I_{\Pi_2}(k)+1}^{\Pi_1}(t_0)<Q_{ I_{\Pi_2}(k)+1}^{\Pi_2}(t_0).$$
After an arrival, the $ (I_{\Pi_1}(k)+1)^{\mathrm{st}}$ term in the left side of \eqref{eq:stoch-ord2-mor} decreases by one, and the $ (I_{\Pi_2}(k)+1)^{\mathrm{st}}$ term increases by one. Thus the inequality is preserved.

If $t_1$ is a departure epoch, then first consider the case when the departure occurs from the green region. In that case, without loss of generality, assume that a potential departure occurs from the $k^{\mathrm{th}}$ ordered server pool, for some $k\in\big\{1,2,\ldots,N\big\}.$ Also note that a departure in either of the two systems can change at most one of the $Q_i$-values.
Thus
\begin{equation}\label{eq:removal-mor}
\tilde{Q}^{\Pi}_i=
\begin{cases}
Q^{\Pi}_i-1, &\mbox{ for }i=I_{\Pi}(k),\\
Q^{\Pi}_j,&\mbox{ otherwise, }
\end{cases}
\end{equation}
if $I_{\Pi}(k)\geq 1$; otherwise all the $Q_i$-values remain unchanged.

If at time epoch $t_0$, $I_{\Pi_1}(k)=I_{\Pi_2}(k)=I$, then both $Q_{ I^{\Pi_1}}$ and $Q_{ I^{\Pi_2}}$ decrease by one, and hence the left side of \eqref{eq:stoch-ord2-mor} does not change. 

Otherwise, without loss of generality assume that $I_{\Pi_1}(k)<I_{\Pi_2}(k).$ Then observe that 
$$Q_{ I_{\Pi_1}(k)}^{\Pi_1}(t_0)\leq Q_{ I_{\Pi_1}(k)}^{\Pi_2}(t_0)\quad \mathrm{and}\quad Q_{ I_{\Pi_2}(k)}^{\Pi_1}(t_0)<Q_{ I_{\Pi_2}(k)}^{\Pi_2}(t_0).$$ 
Furthermore, after the departure, $Q_{ I_{\Pi_1}(k)}^{\Pi_1}$ decreases by one, therefore $|Q_{ I_{\Pi_1}(k)}^{\Pi_1}- Q_{ I_{\Pi_1}(k)}^{\Pi_2}|$ increases by one, and $Q_{ I_{\Pi_2}(k)}^{\Pi_2}$ decreases by one, thus $|Q_{ I_{\Pi_2}(k)}^{\Pi_1}- Q_{I_{ \Pi_2}(k)}^{\Pi_2}|$ decreases by one. Hence, in total, the left side of \eqref{eq:stoch-ord2-mor} remains the same.
Now if a departure occurs from the blue and/or red region, then for some $i_1$ and/or $i_2$, $(Q_{i_1}^{\Pi_1}-Q_{i_1}^{\Pi_2})^+$ or $(Q_{i_2}^{\Pi_2}-Q_{i_2}^{\Pi_1})^+$ (or both) decreases, and the other terms remain unchanged, and hence the left side clearly decreases or remains unchanged. 
 \end{proof}

In order to compare the JSQ policy with the CJSQ$(n(N))$ schemes, and to prove Proposition~\ref{prop: modified JSL}, we will need the following lemma.

\begin{lemma}\label{lem:majorization-mor}
Let $Q_i^{\Pi_1}(t)$ and $Q_i^{\Pi_2}(t)$ denote the number of server pools with at least $i$ tasks under the JSQ policy and CJSQ$(n(N))$ scheme, respectively.
Then for any $k\in\big\{1,2,\ldots, B\big\}$, 
\begin{equation}\label{eq:majorization-mor}
\left\{\sum_{i=1}^k Q_i^{\Pi_1}(t)-kn(N)\right\}_{t\geq 0}\leq_{st} \left\{\sum_{i=1}^k Q_i^{\Pi_2}(t)\right\}_{t\geq 0}\leq_{st}\left\{\sum_{i=1}^k Q_i^{\Pi_1}(t)\right\}_{t\geq 0}, 
\end{equation}
provided the two systems start from the same occupancy states at $t=0$.
\end{lemma}
In the next two remarks we contrast Lemma~\ref{lem:majorization-mor} and the underlying T-coupling with stochastic dominance properties for the ordinary JSQ policy in the existing literature and the S-coupling technique in Chapter~\ref{chap:univjsqd}, respectively.

\begin{remark}\label{rem:novelty2}\normalfont
The stochastic ordering in Lemma~\ref{lem:majorization-mor} is to be contrasted with the weak majorization results in~\cite{Winston77, towsley, Towsley95, Towsley1992, W78} in the context of the ordinary JSQ policy in the single-server queueing scenario, and in~\cite{STC93,J89, M87, MS91} in the scenario of state-dependent service rates, non-decreasing with the number of active tasks.
In the current infinite-server scenario, the results in~\cite{STC93,J89, M87, MS91} imply that for any non-anticipating scheme~$\Pi$ taking assignment decisions based on the number of active tasks only, for all $t\geq 0$,
\begin{align}\label{eq: towsley-mor}
\sum_{m=1}^\ell X_{(m)}^{\JSQ}(t)&\leq_{st}\sum_{m=1}^\ell X_{(m)}^{\Pi}(t),\mbox{ for } \ell=1,2,\ldots, N,\\\nonumber\\
\left\{L^{\JSQ}(t)\right\}_{t\geq 0}&\leq_{st}\left\{L^{\Pi}(t)\right\}_{t\geq 0},
\end{align}
where $X^{\Pi}_{(m)}(t)$ is the number of tasks in the $m^{\mathrm{th}}$ ordered server pool at time $t$ in the system following scheme $\Pi$ and $L^{\Pi}(t)$ is the total number of overflow events under policy $\Pi$ up to time $t$. 
Observe that $X_{(m)}^{\Pi}$ can be visualized as the $m^{\mathrm{th}}$ largest (rightmost) vertical bar (or stack) in Figure~\ref{figB}.
 Thus~\eqref{eq: towsley-mor} says that the sum of the lengths of the $\ell$ largest \emph{vertical} stacks in a system following any scheme $\Pi$ is stochastically larger than or equal to that following the ordinary JSQ policy for any $\ell=1,2,\ldots,N$. Mathematically, this ordering can be equivalently written as
 \begin{equation}\label{eq:equiv-ord-mor}
 \sum_{i = 1}^{B} \min\big\{\ell, Q_i^{\JSQ}(t)\big\}  \leq_{st}
\sum_{i = 1}^{B} \min\big\{\ell, Q_i^{\Pi}(t)\big\},
 \end{equation}
for all $\ell = 1, \dots, N$.
In contrast, in order to show asymptotic equivalence on various scales, we need to both upper and lower bound the occupancy states of the CJSQ$(n(N))$ schemes in terms of the JSQ policy, and therefore need a much stronger hold on the departure process.
The T-coupling provides us just that, and has several useful properties that are crucial for our proof technique.
For example, Proposition~\ref{prop:stoch-ord2-mor} uses the fact that if two systems are T-coupled, then departures cannot increase the sum of the absolute differences of the $Q_i$-values, which is not true for the coupling considered in the above-mentioned  literature.
The left stochastic ordering in~\eqref{eq:majorization-mor} also does not remain valid in those cases.
Furthermore, observe that the right inequality in~\eqref{eq:majorization-mor} (i.e., $Q_i$'s) implies the stochastic inequality is \emph{reversed} in~\eqref{eq:equiv-ord-mor}, which is counter-intuitive in view of the optimality properties of the ordinary JSQ policy studied in the literature, as mentioned above.
The fundamental distinction between the two coupling techniques is also reflected by the fact that the T-coupling does not allow for arbitrary nondecreasing state-dependent departure rate functions, unlike the couplings in~\cite{STC93,J89, M87, MS91}.
\end{remark}

\begin{remark}\label{rem:contrast2}  \normalfont
As briefly mentioned in the introduction, in the current infinite-server scenario, the departures of the ordered server pools cannot be coupled, mainly since the departure rate at the $m^\mathrm{ th}$ ordered server pool, for some $m = 1,2,\ldots, N$, depends on its number of active tasks.
It is worthwhile to mention that the coupling in this chapter is stronger than that used in Chapter~\ref{chap:univjsqd}.
Observe that due to Lemma~\ref{lem:majorization-mor}, the absolute difference of the occupancy states of the JSQ policy and any scheme from the CJSQ class at any time point can be bounded deterministically (without any terms involving the cumulative number of lost tasks).
It is worth emphasizing that the universality result on some specific scale, stated in Theorem~\ref{th:pwr of d} does not depend on the behavior of the JSQ policy on that scale, whereas in Chapter~\ref{chap:univjsqd} it does, mainly because the upper and lower bounds in Corollary~\ref{cor:bound-ssy} involve tail sums of two different policies.
Also, the bound in the current chapter does not depend upon $t$, and hence, applies in the steady state as well.
Moreover, the coupling in Chapter~\ref{chap:univjsqd} compares the $k$ \emph{highest} horizontal bars, whereas the present chapter compares the $k$ \emph{lowest} horizontal bars.
As a result, the bounds on the occupancy states established in Corollary~\ref{cor:bound-ssy} involve tail sums of the occupancy states of the ordinary JSQ policy,
which necessitates proving convergence of the occupancy states of the ordinary JSQ policy with respect to the $\ell_1$ topology. 
In contrast, the bound we establish in the present chapter, involves only a single component (see equations \eqref{eq:upperocc-mor} and \eqref{eq:lowerocc-mor}), and thus, the convergence with respect to product topology suffices.
\end{remark}

\begin{remark}\label{rem:jap-comp-mor}\normalfont
As mentioned in the introduction, a coupling method is used in Chapter~\ref{chap:jiq} to establish the diffusion limit of the Join-the-Idle Queue (JIQ) policy starting from specific initial occupancy states.
Comparing the JIQ and JSQ policies in that scaling regime was much facilitated when viewed as follows:
(i) If there is an idle server in the system, both JIQ and JSQ perform similarly.
(ii)~Also, when there is no idle server and only $O(\sqrt{N})$ servers with queue length two, JSQ assigns the arriving task to a server with queue length one. 
In that case, since JIQ assigns at random, the probability that the task will land on a server with queue length two and thus JIQ acts differently than JSQ is $O(1/\sqrt{N})$.
Since on any finite time interval the number of times an arrival finds all servers busy is at most $O(\sqrt{N})$, all the arrivals except an $O(1)$ of them are assigned in exactly the same manner in both JIQ and JSQ, which then leads to the same scaling limit for both policies.
Note that in the computation of the expected number of events when JIQ and JSQ performs differently, both the specific initial state condition and the scaling regime were crucial.
In the current chapter the stochastic comparison framework is inherently different.
Here the idea pivots on two key observations: 
(i)~For any scheme, if each arrival is assigned to \emph{approximately} the shortest queue, then the scheme can still retain its optimality on various scales, and
(ii)~For any two schemes, if on any finite time interval not \emph{too many} arrivals are assigned to different ordered servers, then they still have the same scaling limits. 
Combination of the above two ideas provides a much wider coupling framework involving an intermediate class of schemes that enables us to consider arbitrary starting states and different scaling regimes.
In addition, the consideration of the arbitrary starting state will turn out to be crucial in order to extend the fluid-scale universality result to the steady state.
\end{remark}

\begin{proof}[Proof of Lemma~\ref{lem:majorization-mor}.]
Fix any $k\geq 1$. We will use forward induction on the event times, i.e., time epochs when either an arrival or a departure occurs,
and assume the two systems to be T-coupled as described in Section~\ref{sec:stoch-coupling}. 
We suppose that the two inequalities hold at time epoch $t_0$, and will prove that they continue to hold at time epoch $t_1$.

(a) We first prove the left inequality in~\eqref{eq:majorization-mor}. 
We distinguish between two cases depending on whether the next event time $t_1$ is an arrival epoch or a departure epoch.
We first consider the case of an arrival.
Since at each arrival, there can be an increment of size at most one, if $\sum_{i=1}^k Q_i^{\Pi_1}(t_0)-kn(N)< \sum_{i=1}^k Q_i^{\Pi_2}(t_0)$, the inequality holds trivially at time $t_1$. 
Therefore, consider the case when $\sum_{i=1}^k Q_i^{\Pi_1}(t_0)-kn(N)= \sum_{i=1}^k Q_i^{\Pi_2}(t_0)$. 
Now observe that
$$\sum_{i=1}^k Q_i^{\Pi_2}(t_0)=\sum_{i=1}^k Q_i^{\Pi_1}(t_0)-kn(N)\leq kN-kn(N).$$
Hence, $Q_k^{\Pi_2}(t_0)\leq N-n(N)$, which in turn implies that at time $t_1$, $\sum_{i=1}^k Q_i^{\Pi_2}$ increases by 1, and the inequality is preserved. 
We now assume the case of a departure.
Then also if $\sum_{i=1}^k Q_i^{\Pi_1}(t_0)-kn(N)< \sum_{i=1}^k Q_i^{\Pi_2}(t_0)$, the inequality holds trivially at time $t_1$. Otherwise assume $\sum_{i=1}^k Q_i^{\Pi_1}(t_0)-kn(N)= \sum_{i=1}^k Q_i^{\Pi_2}(t_0)$. 
In this case if the departure occurs from the green region in Figure~\ref{fig:tcoupling}, then both $\sum_{i=1}^k Q_i^{\Pi_1}$ and $\sum_{i=1}^k Q_i^{\Pi_2}$ change in a similar fashion (i.e., either decrease by one or remain unchanged). 
Else, if the departure occurs from the red and blue regions, since $\sum_{i=1}^k Q_i^{\Pi_1}\geq \sum_{i=1}^k Q_i^{\Pi_2}$, by virtue of the T-coupling, if $\sum_{i=1}^k Q_i^{\Pi_2}$ decreases by one, then so does $\sum_{i=1}^k Q_i^{\Pi_1}$. 
To see this observe the following:
\begin{equation}
\sum_{i=1}^k Q_i^{\Pi_1}\geq \sum_{i=1}^k Q_i^{\Pi_2}\implies\sum_{i=1}^{k}(Q_i^{\Pi_1}-Q_i^{\Pi_2})^+\geq\sum_{i=1}^{k}(Q_i^{\Pi_2}-Q_i^{\Pi_1})^+.
\end{equation}
Therefore, if $m\leq \sum_{i=1}^{k}(Q_i^{\Pi_2}-Q_i^{\Pi_1})^+$, then $m\leq \sum_{i=1}^{k}(Q_i^{\Pi_1}-Q_i^{\Pi_2})^+$.
 Hence the inequality will be preserved.

(b) We now prove the right inequality in~\eqref{eq:majorization-mor}
and again distinguish between two cases. 
If $t_1$ is an arrival epoch, then
 for a similar reason as above,  we assume that $\sum_{i=1}^k Q_i^{\Pi_2}(t_0)=\sum_{i=1}^k Q_i^{\Pi_1}(t_0)$. 
In this case when a task arrives, if it gets admitted under the CJSQ($n(N)$) scheme and increases $\sum_{i=1}^k Q_i^{\Pi_2}$, then clearly $\sum_{i=1}^k (N-Q_i^{\Pi_1}(t))>0$, and hence the incoming task will increase $\sum_{i=1}^k Q_i^{\Pi_1}$, as well, and the inequality will be preserved. 
If $t_1$ is a departure epoch with $\sum_{i=1}^k Q_i^{\Pi_2}(t_0)=\sum_{i=1}^k Q_i^{\Pi_1}(t_0)$, then by virtue of the T-coupling again, if  $\sum_{i=1}^k Q_i^{\Pi_1}$ decreases by one, then by the argument in (a) above, so does $\sum_{i=1}^k Q_i^{\Pi_2}$, thus preserving the inequality.
\end{proof}

\subsection{Asymptotic equivalence}\label{ssec:asymp-eq}

\begin{proof}[Proof of Proposition~\ref{prop: modified JSL}.]
Using Lemma~\ref{lem:majorization-mor},
there exists a common probability space such that for any $k\geq 1$ we can write
\begin{equation}\label{eq:upperocc-mor}
\begin{split}
Q_k^{\Pi_2}(t)&=\sum_{i=1}^k Q_i^{\Pi_2}(t)-\sum_{i=1}^{k-1} Q_i^{\Pi_2}(t)\\
&\leq \sum_{i=1}^{k} Q_i^{\Pi_1}(t)-\sum_{i=1}^{k-1} Q_i^{\Pi_1}(t)+kn(N)\\
&=Q_k^{\Pi_1}(t)+kn(N).
\end{split}
\end{equation}
Similarly, we can write
\begin{equation}\label{eq:lowerocc-mor}
\begin{split}
Q_k^{\Pi_2}(t)&=\sum_{i=1}^k Q_i^{\Pi_2}(t)-\sum_{i=1}^{k-1} Q_i^{\Pi_2}(t)\\
&\geq \sum_{i=1}^{k} Q_i^{\Pi_1}(t)-kn(N)-\sum_{i=1}^{k-1} Q_i^{\Pi_1}(t)\\
&=Q_k^{\Pi_1}(t)-kn(N).
\end{split}
\end{equation}
Therefore, for all $k\geq 1$, we have, $\sup_t|Q_k^{\Pi_2}(t)-Q_k^{\Pi_1}(t)|\leq kn(N)$, and since $n(N)/ g(N)\to 0$ as $N\to\infty$, the proof is complete.
\end{proof}


\begin{proof}[Proof of Proposition~\ref{prop: power of d}.]
For any $T\geq 0$, let $A^N(T)$ and $\Delta^N(T)$ be the total number
of arrivals to the system and the cumulative number of times that the
JSQ$(d(N))$ scheme and the JSQ$(n(N),d(N))$ scheme differ in decision up to time $T$.
Using Proposition~\ref{prop:stoch-ord2-mor} it suffices to show that
for any $T\geq 0$, $\Delta^N(T)/g(N)\pto 0$ as $N\to\infty$.
Observe that at any arrival epoch, the systems under the JSQ$(d(N))$ and
JSQ$(n(N),d(N))$ schemes will differ in decision only if none of the $n(N)+1$ lowest ordered server pools get selected by the JSQ$(d(N))$ scheme.

Now at the time of an arrival, the probability that the JSQ$(d(N))$ scheme does not select one of the $n(N)+1$ lowest ordered server pools, is given by 
$$p(N)=\left(1-\frac{n(N)+1}{N}\right)^{ d(N)}.$$
Since at each arrival epoch $d(N)$ server pools are selected independently, given $A^N(T)$, $\Delta^N(T)\sim \mbox{Bin}(A_N(T),p(N))$.

Note that, for $T\geq 0$, Markov's inequality yields
$$\Pro{\Delta^N(T)\geq g(N)\given A_N(T)}\leq \frac{\E{\Delta^N(T)}}{g(N)}=\frac{A_N(T)}{g(N)}\left(1-\frac{n(N)+1}{N}\right)^{ d(N)}.$$
Since $\big\{A^N(T)/N\big\}_{N\geq 1}$ is a tight sequence of random variables, in order to ensure that $\Delta^N(T)/g(N)$ converges to zero in probability, it is enough to have 
\begin{equation}\label{eq:nN-order}
\begin{split}
&\frac{N}{g(N)}\left(1-\frac{n(N)+1}{N}\right)^{ d(N)}\to 0\\
\Longleftarrow\hspace{.15cm}&\exp\left(\log\left(\frac{N}{g(N)}\right)-d(N)\frac{n(N)}{N}\right)\to 0\\
\iff& d(N)\frac{n(N)}{N}-\log\left(\frac{N}{g(N)}\right)\to\infty,
\end{split}
\end{equation}
which completes the proof.
\end{proof}

We now use Propositions~\ref{prop: modified JSL} and~\ref{prop: power of d}  to prove Theorem~\ref{th:pwr of d}.\\
\begin{proof}[Proof of Theorem~\ref{th:pwr of d}.]
Fix any $d(N)$ satisfying either~\eqref{eq:fNalike cond1} or~\eqref{eq:fNalike cond2}.
From Propositions~\ref{prop: modified JSL} and~\ref{prop: power of d} observe that it is enough to 
show that there exists an $n(N)$ with $n(N)\to\infty$ and $n(N)/g(N)\to 0$, as $N\to\infty$, such that
$$\frac{n(N)}{N}d(N)-\log\left(\frac{N}{g(N)}\right)\to\infty.$$

(i) If $g(N) = O(N)$, then observe that $\log(N/g(N))$ is $O(1)$. 
Since $d(N)\to\infty$, choosing $n(N) = N/ \log (d(N))$ satisfies the above criteria, and hence part (i) of the theorem is proved.

(ii) Next we obtain a choice of $n(N)$ if $g(N)=o(N)$.
Note that, if
$$h(N):=\frac{d(N)\frac{g(N)}{N}}{\log\left(\frac{N}{g(N)}\right)}\to\infty,\quad \text{as}\quad N\to\infty,$$
then choosing $n(N)= g(N)/\log(h(N))$, it can be seen that as $N\to\infty$, we have $n(N)/g(N)\to 0$ and
\begin{equation}\label{eq:equiv-nN-dN}
\begin{split}
&\frac{d(N)\frac{n(N)}{N}}{\log\left(\frac{N}{g(N)}\right)}=\frac{h(N)}{\log(h(N))}\to\infty\\
\implies& \frac{n(N)}{N}d(N)-\log\left(\frac{N}{g(N)}\right)\to\infty.
\end{split}
\end{equation}
\end{proof}

\section{Fluid limit of JSQ}\label{sec:fluid-mor}
In this section we establish the fluid limit for the ordinary JSQ policy.
In the proof we will leverage the time scale separation technique developed in~\cite{HK94}, suitably extended to an infinite-dimensional space. 
Specifically, note that the rate at which incoming tasks join a server pool
with $i$ active tasks is determined only by the process $\ZZ^N(\cdot)=(Z_1^N(\cdot),\ldots,Z_B^N(\cdot))$, where $Z_i^N(t)=N-Q_i^N(t)$, $i=1,\ldots,B$, represents the number of server pools with fewer than $i$ tasks at time $t$.
Furthermore, in any time interval $[t,t+\varepsilon]$ of length $\varepsilon>0$, the $\ZZ^N(\cdot)$ process experiences $O(\varepsilon N)$ events (arrivals and departures), 
while the $\qq^N(\cdot)$ process can change by only $O(\varepsilon)$ amount.
Therefore, the $\ZZ^N(\cdot)$ process evolves on a much faster time scale than
the $\qq^N(\cdot)$ process.
As a result, in the limit as $N\to\infty$, at each time point $t$, the $\ZZ^N(\cdot)$ process
achieves stationarity depending on the instantaneous value of the $\qq^N(\cdot)$ process, i.e., a separation of time scales takes place. 

In order to illuminate the generic nature of the proof construct,
we will allow for a more general task assignment probability and departure dynamics than described in Section~\ref{sec: model descr-mor}.
Denote by $\bZ_+$ the one-point compactification of the set of nonnegative integers~$\Z_+$, i.e., $\bZ_+=\bZ_+\cup\{\infty\}$.
Equip $\bZ_+$ with the order topology. Denote $G=\bZ_+^B$ equipped with product-topology, and with the Borel $\sigma$-algebra, $\mathcal{G}$.
Let us consider the $G$-valued process $\ZZ^N(s):=(Z_i^N(s))_{i\geq 1}$ as introduced above.
Let $\big\{\mathcal{R}_i\big\}_{1\leq i\leq B}$ be a partition of $G$ such that $\mathcal{R}_i\in\mathcal{G}$.
We assume that a task arriving at (say) $t_k$ is assigned to 
some server pool with $i$ active tasks is given by $p_{i-1}^N(\QQ^N(t_k-))=\ind{\ZZ^N(t_k-)\in\mathcal{R}_{i}}f_i(\qq^N(t_k-))$, where $\ff=(f_1,\ldots,f_B):[0,1]^B\to [0,1]^B$ is Lipschitz continuous, i.e., there exists $C_\ff$, such that for any $\qq_1,\qq_2\in S,$
$$\norm{\ff(\qq_1)-\ff(\qq_2)}\leq C_\ff \norm{\qq_1-\qq_2}.$$
The partition corresponding to the ordinary JSQ policy can be written as
\begin{equation}\label{eq:partition-mor}
\mathcal{R}_i := \big\{(z_1,z_2,\ldots, z_B): z_1=\ldots=z_{i-1}=0<z_i\leq z_{i+1}\leq\ldots\leq z_B\big\},
\end{equation}
with the convention that $Q^N_B$ is always taken to be zero, if $B<\infty$, and $f_i\equiv 1$ for all $i=1,2,\ldots,B$.
The fluid-limit results up to Proposition~\ref{prop:rel compactness-mor} (the relative compactness of the fluid-scaled process) hold true for these general assignment probabilities. 
It is only when proving Theorem~\ref{th:genfluid}, that we need to assume the specific $\big\{\mathcal{R}_i\big\}_{1\leq i\leq B}$ in~\eqref{eq:partition-mor}.
For the departure dynamics, when the system occupancy state is $\QQ^N=(Q_1^N,Q_2^N,\ldots,Q_B^N)$, 
define
the total rate at which departures occur from a server pool with $i$ active tasks by $\mu_i^N(\QQ)$, where $\mmu^N(\QQ)=(\mu_1^N(\QQ),\ldots,\mu_B^N(\QQ))$ will be referred to as the departure rate function.
The departure dynamics described in Section~\ref{sec: model descr-mor} correspond to $\mu_i^N(\QQ)=i(Q_i-Q_{i+1})$ and will be referred to as the infinite-server scenario,
since all active tasks are executed concurrently.
The single-server scenario, where tasks are executed sequentially, corresponds to  the case $\mu_i^N(\QQ)=Q_i-Q_{i+1}$.
\begin{assumption}[{Condition on departure rate function}] \label{assump:mu}
The departure rate function $\mmu^N:\tilde{S}\to[0,\infty)^B$ satisfies the following conditions:
\begin{enumerate}[{\normalfont (a)}]
\item There exists a function $\mmu:S\to[0,\infty)^B$, such that 
$$\lim_{N\to\infty}\sup_{\qq\in S^N}\norm{\frac{1}{N}\mmu^N(\lfloor N\qq\rfloor)-\mmu(\qq)}=0.$$
\item The function $\mmu$ is Lipschitz continuous in $S$, i.e., there exists a constant $C_{\mmu}<\infty$, such that for any $\qq_1,\qq_2\in S$, 
$$\norm{\mmu(\qq_1)-\mmu(\qq_2)}\leq C_{\mmu} \norm{\qq_1-\qq_2}.$$
\item Also, $\mmu^N$ satisfies linear growth constraints in each coordinate, i.e., for all $i\geq 1$, there exists $C_i>0$, such that for all $\qq\in S$,
 $$\mu^N_i(\lfloor N\qq\rfloor)\leq NC_i(1+\norm{\qq}).$$
We will often omit $\lfloor \cdot\rfloor$ in the argument of $\mmu^N$ for notational convenience.
\end{enumerate}
\end{assumption}
Under these assumptions on the departure rate function, we prove the following fluid-limit result for the ordinary JSQ policy.
Recall the definition of $m(\qq)$ in Subsection~\ref{ssec:fluid}.
If $m(\qq)=0$, then define $p_0(m(\qq))=1$ and $p_i(m(\qq))=0$ for all $i=1,2,\ldots$. 
Otherwise, in case $m(\qq)>0$, define
\begin{equation}\label{eq:fluid-gen-mor}
p_{i}(\qq)=
\begin{cases}
\min\big\{\mu_{m(\qq)}(\qq)/\lambda,1\big\} & \quad\mbox{ for }\quad i=m(\qq)-1,\\
1 - p_{ m(\qq) - 1}(\qq) & \quad\mbox{ for }\quad i=m(\qq),\\
0&\quad \mbox{ otherwise.}
\end{cases}
\end{equation}
Note that $p_i(\cdot)$ in~\eqref{eq:fluid-gen-mor} is consistent with the one defined in Subsection~\ref{ssec:fluid} for the proper choice of the departure rate function $\mu_i(\qq)=i(q_i-q_{i+1})$.\\
\begin{theorem}[{Fluid limit of JSQ}]
\label{th:genfluid}
Assume $\qq^N(0)\pto \qq^\infty\in S$ and $\lambda(N)/N\to\lambda>0$ as $N\to\infty$. Further assume that the departure rate function $\mmu^N$ satisfies Assumption~\ref{assump:mu}. Then  with probability~1, any subsequence of $\{N\}$ has a further subsequence along which on any finite time interval, the sequence of processes
$\big\{\qq^N(t)\big\}_{t \geq 0}$ converges to some deterministic trajectory $\big\{\qq(t)\big\}_{t \geq 0}$ that satisfies the system of integral equations
\begin{equation}\label{eq:fluidfinal-mor}
q_i(t) = q_i(0)+\lambda \int_0^t p_{i-1}(\qq(s))\dif s - \int_0^t \mu_i(\qq(s))\dif s, \quad i=1,2,\ldots,B,
\end{equation}
where $\qq(0) = \qq^\infty$ and the coefficients $p_i(\cdot)$ are
defined in~\eqref{eq:fluid-gen-mor}, and may be interpreted as the fractions of incoming tasks
assigned to server pools with exactly $i$ active tasks.
\end{theorem}
We will now verify that the departure rate functions corresponding to the infinite-server and single-server scenarios satisfy the conditions in Assumption~\ref{assump:mu}.\\
\begin{proposition}
The following departure rate functions denoted by $\mmu = (\mu_1, \mu_2,\ldots,\mu_B)$, satisfy the conditions in Assumption~\ref{assump:mu}.
For $\QQ\in \tilde{S}$, and $\qq\in S$,
\begin{enumerate}[{\normalfont (i)}]
\item $\mu_i^N(\QQ)=Q_i-Q_{i+1}$, and $\mu_i(\qq)=q_i-q_{i+1}$, $i\geq 1$.
\item $\mu_i^{N}(\QQ)=i(Q_i-Q_{i+1})$, and $\mu_i(\qq)=i(q_i-q_{i+1})$, $i\geq 1$.
\end{enumerate}
\end{proposition}
\begin{proof}
Observe that if $B<\infty$, then since componentwise $\mu_i$ satisfies all the conditions for all $i\geq 1$, $\mmu$ satisfies the conditions in the product space as well.
Therefore, let us consider the case when $B=\infty$.
In this case observe that,
for both (i) and (ii) Assumption~\ref{assump:mu} (a) is immediate, since $\mmu^N(\lfloor N\qq\rfloor)/N=\mmu(\qq)$ for all $\qq\in S^N.$
Also, the linear growth rate constraint in Assumption~\ref{assump:mu}~(c) is satisfied in both cases by taking $C_i=1$ in~(i) and $C_i=i$ in~(ii).

Now we will show that in both cases $\mmu$ is Lipschitz continuous in $S$.

(i) For $\mu_i(\qq)=q_i-q_{i+1}$, $i\geq 1$, and $\qq_1,\qq_2\in S$,
$$
\norm{\mmu(\qq)}=\sum_{i\geq 1}\frac{|q_{i}-q_{i+1}|}{2^i}
\leq \sum_{i\geq 1}\frac{q_{i}}{2^i}+\sum_{i\geq 1}\frac{q_{i+1}}{2^i}
\leq 2\norm{\qq}.
$$

(ii) Now assume $\mu_i(\qq)=i(q_i-q_{i+1})$, $i\geq 1$. 
Since $\mmu$ is a linear operator on the Banach space (complete normed linear space) $\R^B$, to prove Lipschitz continuity of $\mu$, it is enough to show that $\mu$ is continuous at zero. 
Specifically, we will show that for any sequence $\big\{\qq^n\big\}_{n\geq 1}$, in $\R^B$, $\norm{\qq^n}\to 0$ implies $\norm{\mmu(\qq^n)}\to 0$. 
This would imply that there exists fixed $\kappa>0$, such that whenever $\norm{\qq^n}\leq\kappa$ with $\qq^n\in \R^B$, we have $\norm{\mmu(\qq^n)}<1$. Then due to linearity of $\mmu$, for any $\qq\in \R^B$,
\begin{align*}
\norm{\mmu(\qq)}&=\norm{\frac{\norm{\qq}}{\kappa}\mmu\left(\kappa\frac{\qq}{\norm{\qq}}\right)}\\
&\leq \frac{\norm{\qq}}{\kappa}\norm{\mmu\left(\kappa\frac{\qq}{\norm{\qq}}\right)}\\
&\leq \frac{1}{\kappa}\norm{\qq}.
\end{align*}
To show that $\mmu$ is continuous at $\mathbf{0}\in\R^B$, fix
any $\varepsilon>0$. 
Also, fix an $M>0$, depending upon~$\varepsilon$, such that $\sum_{i> M}1/2^i<\varepsilon/2$. Now, choose $\delta<\varepsilon/(4M)$. 
Then, for any $\qq\in \R^B$ such that $\norm{\qq}<\delta$,
\begin{align*}
\norm{\mmu(\qq)}&=\sum_{i=1}^\infty \frac{i|q_i-q_{i+1}|\wedge 1}{2^i}
\leq \sum_{i=1}^M \frac{i|q_i-q_{i+1}|\wedge 1}{2^i}+\frac{\varepsilon}{2}\\
&\leq M\sum_{i=1}^M \frac{|q_i-q_{i+1}|\wedge 1}{2^i}+\frac{\varepsilon}{2}
\leq 2M\norm{\qq}+\frac{\varepsilon}{2}\leq \varepsilon.
\end{align*}
Hence, $\mmu$ is Lipschitz continuous on $\R^\infty$.
\end{proof}

\subsection{Martingale representation}
In this subsection we construct the martingale representation of the occupancy state process $\QQ^N(\cdot)$.
The component $Q_i^N(t)$, satisfies the identity relation
\begin{align}\label{eq:recursion-mor}
Q_i^N(t)=Q_i^N(0)+A_i^N(t)-D_i^N(t),&\quad\mbox{ for }\quad i=1,\ldots, B,
\end{align}
where
\begin{align*}
A_i^N(t)&=\mbox{ number of arrivals during $[0,t]$ to some server pool with }i-1\mbox{ active tasks,} \\
D_i^N(t)&=\mbox{ number of departures during $[0,t]$ from some server pool with }i\mbox{ active tasks}.
\end{align*}
We can express $A^N_i(t)$ and $D_i^N(t)$ as
\begin{align*}
A^N_i(t) &=  \mathcal{N}_{A,i}\left(\lambda(N)\int_0^t p_{i-1}^N(\QQ^N(s))\dif s\right),\\
D_i^N(t) &=  \mathcal{N}_{D,i}\left(\int_0^t \mu_i^N(\QQ^N(s))\dif s\right),
\end{align*}
where $\mathcal{N}_{A,i}$ and $\mathcal{N}_{D,i}$ are mutually independent unit-rate Poisson processes, $i=1,2,\ldots,B$.
Define the following sigma fields.
\begin{align*}
\mathcal{A}^N_i(t)&:= \sigma\left(A^N_i(s): 0\leq s\leq t\right),\\
\mathcal{D}_i^N(t)&:= \sigma\left(D_i^N(s): 0\leq s\leq t\right),\mbox{ for }i\geq 1,
\end{align*}
and the filtration $\mathbf{F}^N\equiv\big\{\mathcal{F}^N_t:t\geq 0\big\}$ with
\begin{equation}\label{eq:filtration-mor}
\mathcal{F}^N_t:=\bigvee_{i=1}^{\infty} \big[\mathcal{A}_i^N(t)\vee \mathcal{D}_i^N(t)\big]
\end{equation}
augmented by all the null sets. 
Now we have the following martingale decomposition from the classical result in \cite[Lemma~3.2]{PTRW07}.

\begin{proposition}\label{prop:mart-rep-mor}
The following are $\mathbf{F}^N$-martingales, for $i\geq 1$:
\begin{equation}\label{eq:martingales-mor}
\begin{split}
M_{A,i}^N(t)&:=  \mathcal{N}_{A,i}\left(\lambda(N)\int_0^t p_{i-1}^N(\QQ^N(s))\dif s\right)-\lambda(N)\int_0^t p_{i-1}^N(\mathbf{Q}^N(s)) \dif s,\\
M_{D,i}^N(t)&:=\mathcal{N}_{D,i}\left(\int_0^t \mu_i^N(\QQ^N(s))\dif s\right)-\int_0^t \mu_i^N(\QQ^N(s))\dif s,
\end{split}
\end{equation}
with respective compensator and predictable quadratic variation processes given by
\begin{align*}
\langle M_{A,i}^N\rangle(t)&:= \lambda(N)\int_0^t p_{i-1}^N(\mathbf{Q}^N(s-))\dif s,\\
\langle M_{D,i}^N\rangle(t)&:=\int_0^t \mu_i^N(\QQ^N(s))\dif s.
\end{align*}
\end{proposition}

Therefore, we finally have the following martingale representation of the $N^{\mathrm{th}}$ process:
\begin{equation}\label{eq:mart-unscaled-mor}
\begin{split}
Q_i^N(t)&=Q_i^N(0)+\lambda(N)\int_0^t p_{i-1}^N(\mathbf{Q}^N(s))\dif s\\
&-\int_0^t \mu_i^N(\QQ^N(s))\dif s +(M_{A,i}^N(t)-M_{D,i}^N(t)),\quad t\geq 0,\quad i= 1,\ldots,B.
\end{split}
\end{equation}
In the proposition below, we prove that the martingale part vanishes when scaled by $N$. Since convergence in probability in each component implies convergence in probability with respect to the product topology, it is enough to show convergence in each component.

\begin{proposition}\label{prop:mart zero1-mor}
For all $i\geq 1$,
$$\left\{\frac{1}{N}(M_{A,i}^N(t)-M_{D,i}^N(t))\right\}_{t\geq 0}\dto \big\{m(t)\big\}_{t\geq 0}\equiv 0.$$
\end{proposition}
\begin{proof}
Fix any $T\geq 0$, and $i\geq 1$. From Doob's inequality \cite[Theorem 1.9.1.3]{LS89}, we have for any~$\epsilon>0$,
\begin{align*}
\Pro{\sup_{t\in[0,T]}\frac{1}{N}|M_{A,i}^N(t)|\geq \epsilon}&=\Pro{\sup_{t\in[0,T]}|M_{A,i}^N(t)|\geq N\epsilon}\\
&\leq \frac{1}{N^2\epsilon^2}\E{\langle M_{A,i}^N\rangle (T)}\\
&\leq \frac{1}{N^2\epsilon^2}\int_0^T p_{i-1}(\mathbf{Q}^N(s-))\lambda(N)\dif s\\
&\leq \frac{\lambda(N)}{N^2}\cdot\frac{T}{\epsilon^2}\to 0,\mbox{ as }N\to\infty.
\end{align*}
Similarly, for $M_{D,i}^N$, 
\begin{align*}
\Pro{\sup_{t\in[0,T]}\frac{1}{N}|M_{D,i}^N(t)|\geq \epsilon}&=\Pro{\sup_{t\in[0,T]}|M_{D,i}^N(t)|\geq N\epsilon}\\
&\leq \frac{1}{N^2\epsilon^2}\E{\langle M_{D,i}^N\rangle (T)}\\
&\leq \frac{1}{N^2\epsilon^2}\int_0^T \mu_i^N(\QQ^N(s))\dif s\\
&\leq \frac{2L'}{N\epsilon^2}\to 0,\mbox{ as }N\to\infty,
\end{align*}
where the last inequality follows from the linear growth constraint stated in Assumption~\ref{assump:mu}~(c).
Therefore we have uniform convergence over compact sets, and hence with respect to the Skorohod-$J_1$ topology.  
\end{proof}

\subsection{Relative compactness and uniqueness}\label{subsec:fluidlimit}
Now we will first prove the relative compactness of the sequence of fluid-scaled processes.
Recall that we denote all the fluid-scaled quantities by their respective small letters, e.g.~$\mathbf{q}^N(t):=\mathbf{Q}^N(t)/N$, componentwise, i.e., $q_i^N(t):=Q_i^N(t)/N$ for $i\geq 1$. Therefore the martingale representation in~\eqref{eq:mart-unscaled-mor} can be written as
\begin{equation}\label{eq:mart1-mor}
\begin{split}
q_i^N(t)&=q_i^N(0)+\frac{\lambda(N)}{N}\int_0^t p_{i-1}^N(\mathbf{Q}^N(s))\dif s\\
&-\int_0^t \frac{1}{N}\mu_i^N(\QQ^N(s))\dif s +\frac{1}{N}(M_{A,i}^N(t)-M_{D,i}^N(t)),\quad i=1,2,\ldots, B,
\end{split}
\end{equation}
or equivalently,
\begin{equation}\label{eq:martingale rep assumption 2-mor}
\begin{split}
q_i^N(t)&=q_i^N(0)+\frac{\lambda(N)}{N}\int_0^t f_i(\qq^N(s))\ind{\ZZ^N(s)\in\mathcal{R}_{i}}\dif s\\
&-\int_0^t \frac{1}{N}\mu_i^N(\QQ^N(s))\dif s +\frac{1}{N}(M_{A,i}^N(t)-M_{D,i}^N(t)),\quad i=1,2,\ldots, B.
\end{split}
\end{equation}
Now, we consider the Markov process $(\qq^N,\ZZ^N)(\cdot)$ defined on $S\times G$. 
Define a random measure $\alpha^N$ on the measurable space $([0,\infty)\times G, \mathcal{C}\otimes\mathcal{G})$, when $[0,\infty)$ is endowed with the Borel sigma algebra $\mathcal{C}$, by
\begin{equation}
\alpha^N(A_1\times A_2):=\int_{A_1} \ind{\ZZ^N(s)\in A_2}\dif s,
\end{equation}
for $A_1\in\mathcal{C}$ and $A_2\in\mathcal{G}$. 
Then the representation in \eqref{eq:martingale rep assumption 2-mor} can be written in terms of the random measure as,
\begin{equation}\label{eq:martingale rep assumption 2-2-mor}
\begin{split}
q_i^N(t)&=q_i^N(0)+\lambda\int_{[0,t]\times\mathcal{R}_i} f_i(\qq^N(s))\dif\alpha^N\\
&-\int_0^t \frac{1}{N}\mu_i^N(\QQ^N(s))\dif s +\frac{1}{N}(M_{A,i}^N(t)-M_{D,i}^N(t)),\quad i=1,2,\ldots, B.
\end{split}
\end{equation}
Let $\mathfrak{L}$ denote the space of all measures on $[0,\infty)\times G$ satisfying $\gamma([0,t],G)= t$, endowed with the topology corresponding to weak convergence of measures restricted to $[0,t]\times G$ for each $t$.
\begin{proposition}\label{prop:rel compactness-mor}
Assume that $\mathbf{q}^N(0)\dto\mathbf{q}(0)$ as $N\to\infty$, then $\big\{(\mathbf{q}^N(\cdot),\alpha^N)\big\}$ is a relatively compact sequence in $D_{S}[0,\infty)\times\mathfrak{L}$ and the limit $\big\{(\mathbf{q}(\cdot),\alpha)\big\}$ of any convergent subsequence satisfies
\begin{equation}\label{eq:rel compact-mor}
q_i(t)=q_i(0)+\lambda\int_{[0,t]\times\mathcal{R}_i} f_i(\qq(s))\dif\alpha -\int_0^t \mu_i(\qq(s))\dif s,\quad i=1,2,\ldots, B.
\end{equation}
\end{proposition}
\begin{remark} \normalfont
Proposition~\ref{prop:rel compactness-mor} is true even when the function $\ff$ in the assignment probability depends on $N$. In that case the proof will go through by assuming that $\ff^N$ converges uniformly to some Lipschitz-continuous function $\ff$ in the sense of Assumption~\ref{assump:mu}.(a).
\end{remark}
\begin{remark} \normalfont
The relative compactness result in the above proposition holds for an even more general class of assignment probabilities than those considered above.
Since the proof will follow a nearly identical line of arguments, we briefly mention them here. 
Consider a scheme for which the assignment probabilities can be written as
$$p_{i}^N(\QQ^N) = \eta_1 \ind{\ZZ^N\in\mathcal{R}_i}+\eta_2 g_i(\qq^N),\quad i=1,\ldots,B,$$
for some fixed $\eta_1,\eta_2\in[0,1]$, and some Lipschitz continuous function $$\gb =(g_1,g_2,\ldots,g_B):S\to[0,\infty)^B.$$
The above scheme assigns a fixed fraction $\eta_1$ of incoming tasks according to the ordinary JSQ policy, and a fraction $\eta_2$ as some suitable function of the fluid-scaled occupancy states $\gb(\qq)$, for $\qq\in S$.
In practice, the above scheme can handle (two or more) priorities among the incoming tasks, by assigning the high-priority tasks in accordance with the ordinary JSQ policy, and others governed by the JSQ$(d)$ scheme, say.
In that case, the fluid limit in~\eqref{eq:rel compact-mor} will become
\begin{equation}\label{eq:rel compact-gen}
q_i(t)=q_i(0)+\lambda\eta_1\alpha([0,t]\times\mathcal{R}_i)+\eta_2\int_0^t g_i(\qq(s))\dif s -\int_0^t \mu_i(\qq(s))\dif s,
\end{equation}
$i=1,2,\ldots, B.$
\end{remark}
To prove Proposition~\ref{prop:rel compactness-mor}, we will verify the conditions of relative compactness from~\cite{EK2009}. 
Let $(E,r)$ be a complete and separable metric space. For any $x\in D_E[0,\infty)$, $\delta>0$ and $T>0$, define
\begin{equation}\label{eq:mod-continuity-mor}
w'(x,\delta,T)=\inf_{\{t_i\}}\max_i\sup_{s,t\in[t_{i-1},t_i)}r(x(s),x(t)),
\end{equation}
where $\big\{t_i\big\}$ ranges over all partitions of the form $0=t_0<t_1<\ldots<t_{n-1}<T\leq t_n$ with $\min_{1\leq i\leq n}(t_i-t_{i-1})>\delta$ and $n\geq 1$.
 Below we state the conditions for the sake of completeness.
\begin{theorem}[{\cite[Corollary~3.7.4]{EK2009}}]\label{th:from EK-mor}
Let $(E,r)$ be complete and separable, and let $\big\{X_n\big\}_{n\geq 1}$ be a family of processes with sample paths in $D_E[0,\infty)$. Then $\big\{X_n\big\}_{n\geq 1}$ is relatively compact if and only if the following two conditions hold:
\begin{enumerate}[{\normalfont (a)}]
\item For every $\eta>0$ and rational $t\geq 0$, there exists a compact set $\Gamma_{\eta, t}\subset E$ such that $$\varliminf_{n\to\infty}\Pro{X_n(t)\in\Gamma_{\eta, t}}\geq 1-\eta.$$
\item For every $\eta>0$ and $T>0$, there exists $\delta>0$ such that
$$\varlimsup_{n\to\infty}\Pro{w'(X_n,\delta, T)\geq\eta}\leq\eta.$$
\end{enumerate}
\end{theorem}

\begin{proof}[Proof of Proposition~\ref{prop:rel compactness-mor}.]
The proof goes in two steps. We first prove the relative compactness, and then show that the limit satisfies~\eqref{eq:rel compact-mor}.

Observe from \cite[Proposition 3.2.4]{EK2009} that, to prove the relative compactness of the process $\big\{(\mathbf{q}^N(\cdot),\alpha^N)\big\}$, it is enough to prove relative compactness of the individual components.
Note that, from Prohorov's theorem~\cite[Theorem 3.2.2]{EK2009}, $\mathfrak{L}$ is compact, since $G$ is compact. Now, relative compactness of $\alpha^N$ follows from the compactness of $\mathfrak{L}$ under the topology of weak convergence of measures and Prohorov's theorem.

To claim the relative compactness of $\big\{\mathbf{q}^N(\cdot)\big\}$, first observe that $[0,1]^{B}$ is compact with respect to product topology, and $S$ is a closed subset of $[0,1]^B$, and hence $S$ is also compact with respect to product topology. So, the compact containment condition (a) of Theorem~\ref{th:from EK-mor} is satisfied by taking $\Gamma_{\eta,t}\equiv S$.

For condition (b), we will show for each coordinate $i$, that for any $\eta>0$, there exists $\delta>0$, such that for any $t_1,t_2>0$ with $|t_1-t_2|<\delta$,
\begin{equation*}
\varlimsup_{n\to\infty}\Pro{|q^n_i(t_1)-q^n_i(t_2)|\geq \eta}=0.
\end{equation*}
With respect to product topology, this will imply that for any $\eta>0$, there exists $\delta>0$, such that for any $t_1,t_2>0$ with $|t_1-t_2|<\delta$,
\begin{equation*}
\varlimsup_{n\to\infty}\Pro{\norm{q^n(t_1)-q^n(t_2)}\geq \eta}=0,
\end{equation*}
which in turn will imply condition~(b) in Theorem~\ref{th:from EK-mor}.
To see this, observe that for any fixed $\eta>0$ and $T>0$, we can
choose $\delta'>0$ small enough, so that for any fine enough finite partition $0=t_0<t_1<\ldots<t_{n-1}<T\leq t_n$ of $[0,T]$ with $\min_{1\leq i\leq n}(t_i-t_{i-1})>\delta'$ and $\max_{1\leq i\leq n}(t_i-t_{i-1})<\delta$,  
$$\varlimsup_{n\to\infty}\Pro{\norm{q^n(t_i)-q^n(t_{i+1})}\geq \eta}=0$$ 
for all $1\leq i\leq n$.

Now fix any $0\leq t_1<t_2<\infty$, and  $1\leq i\leq B$. Then
\begin{align*}
&|q_i^N(t_1)-q_i^N(t_2)|\\
&\leq \lambda \alpha^N([t_1,t_2]\times\mathcal{R}_i)+\int_{t_1}^{t_2} \frac{1}{N}\mu_i^N(\QQ^N(s))\dif s \\
&+\frac{1}{N}|M_{A,i}^N(t_1)-M_{D,i}^N(t_1)-M_{A,i}^N(t_2)+M_{D,i}^N(t_2)|\\
&\leq\lambda'(t_2-t_1)+\frac{1}{N}|M_{A,i}^N(t_1)-M_{D,i}^N(t_1)-M_{A,i}^N(t_2)+M_{D,i}^N(t_2)|,
\end{align*}
for some $\lambda'\in\R$, using the linear growth constraint of $\mmu^N$ due to Assumption~\ref{assump:mu}(c). 
Now, from Proposition~\ref{prop:mart zero1-mor}, we get, for any $T\geq 0$,
$$\sup_{t\in[0,T]}\frac{1}{N}|M_{A,i}^N(t_1)-M_{D,i}^N(t_1)-M_{A,i}^N(t_2)+M_{D,i}^N(t_2)|\pto 0.$$
To prove that the limit $\big\{(\qq(\cdot),\alpha)\big\}$ of any convergent subsequence satisfies~\eqref{eq:rel compact-mor}, we will use the continuous-mapping theorem~\cite[Section~3.4]{W02}.
Specifically, we will show that the right side of~\eqref{eq:martingale rep assumption 2-2-mor} is a continuous map of suitable arguments.
Let $\big\{\qq(t)\big\}_{t\geq 0}$ and $\big\{\yy(t)\big\}_{t\geq 0}$ be an $S$-valued and an $\R^B$-valued c\`adl\`ag function, respectively. 
Also, let $\alpha$ be a measure on the measurable space $([0,\infty)\times G, \mathcal{C}\otimes\mathcal{G})$. Then for $\qq^0\in S$, define for $i\geq 1$,
$$F_i(\qq,\alpha,\qq^0,\yy)(t):=q_i^0+y_i(t)+\lambda\int_{[0,t]\times\mathcal{R}_i} f_i(\qq(s))\dif\alpha-\int_0^t\mu_i(\qq(s))\dif s.$$
Observe that it is enough to show that $\FF=(F_1,\ldots,F_B)$ is a continuous
operator. Indeed, in that case the right side of~\eqref{eq:martingale rep assumption 2-2-mor} can be written as $\FF(\qq^N,\alpha^N,\qq^N(0),\yy^N)$, where $\yy^N=(y_1^N,\ldots,y_B^N)$ with $y_i^N= (M_{A,i}^N-M_{D,i}^N)/N$, and since each argument converges we will get the convergence to the right side of~\eqref{eq:rel compact-mor}.
Therefore, we now prove the continuity of $\FF$ below. 
In particular assume that the sequence of processes $\big\{(\qq^N,\yy^N)\big\}_{N\geq 1}$ converges to $\big\{(\qq,\yy)\big\}$, for any fixed $t\geq 0$, the measure $\alpha^N([0,t],\cdot)$ on $G$ converges weakly to $\alpha([0,t],\cdot)$, and the sequence of $S$-valued random  variables $\qq^N(0)$
 converges weakly to $\qq(0)$.
 Fix any $T\geq 0$ and $\varepsilon>0$.
 \begin{enumerate}[{\normalfont (i)}]
 \item  Choose $N_1\in\N$, such that for all $N\geq N_1$, 
 $$\sup_{t\in[0,T]}\norm{\qq^N(t)-\qq(t)}<\varepsilon/(4TC_{\mmu}).$$ In that case, observe that
 \begin{align*}
 \sup_{t\in [0,T]}\int_0^t\norm{\mmu(\qq^N(s))-\mmu(\qq(s))}\dif s& \leq T\sup_{t\in [0,T]}\norm{\mmu(\qq^N(t))-\mmu(\qq(t))}\\
 & \leq TC_{\mmu}\sup_{t\in [0,T]}\norm{\qq^N(t))-\qq(t)}<\frac{\varepsilon}{4},
 \end{align*}
 where we have used the Lipschitz continuity of $\mmu$ due to Assumption~\ref{assump:mu}(b).
 \item Choose $N_2\in\N$, such that for all $N\geq N_2$, $$\sup_{t\in[0,T]}\norm{\yy^N(t)-\yy(t)}<\varepsilon/4.$$
 \item   Choose $N_3\in\N$, such that for all $N\geq N_3$,
 $$\sum_{i\geq 1}\frac{\lambda}{2^i}\left|\int_{[0,T]\times\mathcal{R}_i} f_i(\qq^N(s))\dif\alpha^N-\int_{[0,T]\times\mathcal{R}_i} f_i(\qq(s))\dif\alpha\right|<\frac{\varepsilon}{4}.$$
 This can be done as follows: choose $M\in \N$ large enough so that $\sum_{i> M}2^{-i}<\varepsilon/8.$ Now for $i\leq M$, since $\alpha^N([0,T],\cdot)$ converges weakly to $\alpha([0,T],\cdot)$, and $M$ is finite, we can choose $N_3\in\N$ such that 
 \begin{align*}
 &\sum_{i= 1}^M\frac{\lambda}{2^i}\left|\int_{[0,T]\times\mathcal{R}_i} f_i(\qq^N(s))\dif\alpha^N-\int_{[0,T]\times\mathcal{R}_i} f_i(\qq(s))\dif\alpha\right|\\
 \leq & \sum_{i= 1}^M\frac{\lambda}{2^i}\int_{[0,T]\times\mathcal{R}_i} |f_i(\qq^N(s))-f_i(\qq(s))|\dif\alpha^N\\
 &\hspace{4cm}+\sum_{i = 1}^M\frac{\lambda}{2^i}|\alpha^N([0,T]\times \mathcal{R}_i)-\alpha([0,T]\times \mathcal{R}_i)|
 \end{align*}
 \begin{align*}
 \leq &  \sum_{i= 1}^M\frac{\lambda }{2^i}TC_\ff\sup_{s\in [0,T]} \norm{\qq^N(s)-\qq(s)}\\
 &\hspace{4cm}+\sum_{i = 1}^M\frac{\lambda}{2^i}|\alpha^N([0,T]\times \mathcal{R}_i)-\alpha([0,T]\times \mathcal{R}_i)|<\frac{\varepsilon}{4}.
 \end{align*}
 \item  Choose $N_4\in\N$, such that for all $N\geq N_4$, 
 $$\norm{\qq^N(0)-\qq(0)}<\varepsilon/4.$$
 \end{enumerate}
Let $\hat{N}=\max\big\{N_1,N_2,N_3,N_4\big\}$, then for $N\geq \hat{N}$,
\begin{align*}
&\sup_{t\in [0,T]}\norm{\FF(\qq^N,\alpha^N,\qq^N(0),\yy^N)-\FF(\qq,\alpha,\qq(0),\yy)}(t)<\varepsilon.
\end{align*}
Thus the proof of continuity of $\FF$ is complete. 
\end{proof}

To characterize the limit in~\eqref{eq:rel compact-mor}, for any $\qq\in S$, define the Markov process  $\ZZ_{\qq}$ on $G$ as
\begin{equation}\label{eq:slowprocess-mor}
\ZZ_{\qq} \rightarrow 
\begin{cases}
\ZZ_{\qq}+e_i& \quad\mbox{ at rate }\quad \mu_i(\qq)\\
\ZZ_{\qq}-e_i& \quad\mbox{ at rate }\quad \lambda \ind{\ZZ_\qq\in\mathcal{R}_{i}},
\end{cases}
\end{equation}
where $e_i$ is the $i^{\mathrm{th}}$ unit vector, $i=1,\ldots,B$.\\

\begin{proof}[{Proof of Theorem~\ref{th:genfluid}}.]
Having proved the relative compactness in Proposition~\ref{prop:rel compactness-mor},  it follows from analogous arguments as used in the proof of~\cite[Theorem 3]{HK94}, that the limit of any convergent subsequence of the sequence of processes $\big\{\qq^N(t)\big\}_{t\geq 0}$ satisfies
\begin{equation}
q_i(t) = q_i(0)+\lambda \int_0^t \pi_{\qq(s)}(\mathcal{R}_i)\dif s - \int_0^t \mu_i(\qq(s))\dif s, \quad i=1,2,\ldots,B,
\end{equation}
for \emph{some} stationary measure $\pi_{\qq(t)}$ of the Markov process  $\ZZ_{\qq(t)}$ described in~\eqref{eq:slowprocess-mor} satisfying $\pi_{\qq}\big\{\ZZ: Z_i=\infty\big\}=1$ if $q_i<1$. 

Now it remains to show that $\qq(t)$ \emph{uniquely} determines $\pi_{\qq(t)}$, and that $\pi_{\qq(s)}(\mathcal{R}_i)=p_{i-1}(\qq(s))$ as described in~\eqref{eq:fluid-gen-mor}. 
As mentioned earlier, in this proof we will now assume the specific assignment probabilities in~\eqref{eq:partition-mor}, corresponding to the ordinary JSQ policy.
To see this, fix any $\qq=(q_1,\ldots,q_B)\in S$, 
and assume that there exists $m\geq 0$, such that $q_{m+1}<1$ and $q_1=\ldots=q_m=1$,
with the convention that $q_0\equiv 1$ and $q_{B+1}\equiv 0$ if $B<\infty$. In that case,
$$\pi_{\qq}\big(\big\{Z_{m+1}=\infty, Z_{m+2}=\infty,\ldots,Z_B=\infty\big\}\big)=1.$$
Also, 
note that $q_i = 1$ forces $\dif q_i/\dif t \leq 0$, i.e., $\lambda \pi_{\qq}(\mathcal{R}_i) \leq \mu_i(\qq)$ for all $i = 1, \ldots, m$, and in particular $\pi_{\qq}(\mathcal{R}_i) = 0$ for all $i = 1, \ldots, m - 1.$ Thus,
$$\pi_{\qq}\big(\big\{Z_1=0,Z_2=0,\ldots,Z_{m-1}=0\big\}\big)=1.$$

Therefore, $\pi_\qq$ is determined only by the stationary distribution of the $m^{\mathrm{th}}$ component, which can be described as a birth-death process
\begin{equation}\label{eq:bdprocess-mor}
Z \rightarrow 
\begin{cases}
Z+1& \quad\mbox{ at rate }\quad \mu_m(\qq)\\
Z-1& \quad\mbox{ at rate }\quad \lambda\ind{Z>0}
\end{cases}
\end{equation}
and let $\pi^{(m)}$ be its stationary distribution. 
Now it is enough to show that $\pi^{(m)}$ is uniquely determined by $\mu_m(\qq)$. 
First observe that the process on $\bZ$ described in~\eqref{eq:bdprocess-mor} is reducible, and can be decomposed into
two irreducible classes given by $\mathbb{Z}$ and $\{\infty\}$, respectively.
Therefore, if $\pi^{(m)}(Z=\infty)=0$ or $1$, then it is unique. 
Indeed, if $\pi^{(m)}(Z=\infty)=0$, then $Z$ is a birth-death process on $\mathbb{Z}$ only, and hence it has a unique stationary distribution. 
Otherwise, if $\pi^{(m)}(Z=\infty)=1$, then it is trivially unique. 
Now we distinguish between two cases depending upon whether $\mu_m(\qq)\geq \lambda$ or not.

Note that if $\mu_m(\qq)\geq\lambda$, then $\pi^{(m)}(Z\geq k)=1$ for all $k\geq 0$. 
On $\bZ$ this shows that $\pi^{(m)}(Z=\infty)=1$.
Furthermore, if $\mu_m(\qq)<\lambda$, then we will show that $\pi^{(m)}(Z=\infty)=0$.
On the contrary, assume $\pi^{(m)}(Z=\infty)=\varepsilon\in (0,1]$.
Also, let $\hat{\pi}^{(m)}$ be the unique stationary distribution of the birth-death process in~\eqref{eq:bdprocess-mor} on $\mathbb{Z}$.
Therefore, 
$$\pi_\qq(\mathcal{R}_m)=\pi^{(m)}(Z>0)=(1-\varepsilon)\hat{\pi}^{(m)}(Z>0)+\varepsilon = (1-\varepsilon)\frac{\mu_m(\qq)}{\lambda}+\varepsilon.$$
Substituting onto the differential form of the fluid equation~\eqref{eq:fluidfinal-mor} at the given time $t$, we obtain that 
\begin{align*}
\frac{\dif q_m(t)}{\dif t} = \lambda \Big[(1-\varepsilon)\frac{\mu_m(\qq)}{\lambda}+\varepsilon \Big] - \mu_m(\qq) = -\varepsilon\mu_m(\qq) +\lambda \varepsilon >0,
\end{align*}
where the last inequality follows since we are considering the case when $\mu_m(\qq)<\lambda$.
Since $q_m(t)=1$, this leads to a contradiction for any $\varepsilon>0$, and hence it must be the case that $\pi^{(m)}(Z=\infty)=0$. 

Therefore, for all $\qq\in S$, $\pi_\qq$ is uniquely determined by $\qq$. 
Furthermore, we can identify the expression for $\pi_q(\mathcal{R}_i)$ as
\begin{equation}
\pi_\qq(\mathcal{R}_i)=
\begin{cases}
\min\big\{\mu_i(\qq)/\lambda,1\big\}& \quad\mbox{ for }\quad i=m,\\
1- \min\big\{\mu_i(\qq)/\lambda,1\big\} & \quad\mbox{ for }\quad i=m+1,\\
0&\quad \mbox{ otherwise,}
\end{cases}
\end{equation}
and hence $\pi_{\qq(s)}(\mathcal{R}_i)=p_{i-1}(\qq(s))$ as claimed. 
\end{proof}

\subsection{Global stability and interchange of limits}\label{ssec:globstab-mor}

To prove the interchange of limits result stated in Proposition~\ref{prop:interchange-mor}, we will establish the global stability of the fixed point, i.e., all fluid paths converge to the fixed point in~\eqref{eq:fixed point} as $t\to\infty$. 
\begin{lemma}\label{lem:global-stab-mor}
Let $\qq(t)$ be the fluid limit, i.e., the solution of the dynamical system described by the system of integral equations in~\eqref{eq:fluidjsqd-mor}.
For any $\qq^\infty\in S$ with $\sum_{i=1}^Bq_i^\infty<\infty$, if $\qq(0) = \qq^\infty$, then $\qq(t)\to
\qq^\star$ as $t\to\infty$, where $\qq^\star$ is defined as in \eqref{eq:fixed point}.
\end{lemma}
\begin{proof}
The proof follows in three steps: we will first establish that as $t\to\infty$,
$$q_{\leq K}(t):=\sum_{i=1}^Kq_i(t)\to K,$$ 
and then show that $\sum_{i=1}^Bq_{i}(t)\to \lambda$.
Finally using the above two facts we will show that $q_{K+2}(t)\to 0$ as $t\to\infty$, which will complete the proof.

Observe that the rate of change of $q_{\leq K}(t)$ is $\lambda \sum_{i=1}^Kp_{i-1}(\qq(t))-(q_{\leq K}(t)-Kq_{K+1}(t))$.
For any $\varepsilon> 0$, if $q_{\leq K}(t)\leq K-\varepsilon$, then $\sum_{i=1}^Kp_{i-1}(\qq(t))=1$, so that the rate of change is  $\lambda -(q_{\leq K}(t)-Kq_{K+1}(t))\geq \lambda - K + \varepsilon\geq \varepsilon>0$, i.e., positive and bounded away from zero.
Also, observe that $q_{\leq K}(t)$ cannot decrease if $q_{\leq K}(t)\leq K$.
This shows that for all $\varepsilon>0$, there exists a time $t_0 = t_0(\varepsilon, \qq^\infty)<\infty$, such that $q_{\leq K}(t)\geq K -\varepsilon$ for all $t\geq t_0$.
Thus, $\liminf_{t\to\infty} q_{\leq K}(t)\geq K$, and consequently, $q_{\leq K}(t)\to K$, as $t\to\infty$.

Define $y(t):= \sum_{i=1}^Bq_i(t)$ as the total amount of fluid in the system.
Then note that the rate of change of $y(t)$ is given by $\lambda \sum_{i=1}^Bp_{i-1}(\qq(t))-y(t) = \lambda-y(t)$, and therefore, $y(t) = \lambda + \e^{-t}(y(0) - \lambda)$.
Since $y(0)=\sum_{i=1}^Bq_i^\infty<\infty$, this yields that $y(t)\to\lambda$ as $t\to\infty$.

Finally, define $q_{\geq K+1}(t):= \sum_{i=K+1}^Bq_i(t)$ and $q_{\geq K+2}(t):= \sum_{i=K+2}^Bq_i(t)$.
Since $q_{\leq K}(t)\to K$ and $y(t)\to \lambda$, as $t\to\infty$, we obtain $q_{\geq K+1}(t) = y(t) - q_{\leq K}(t)\to \lambda -K= f$.
Consequently, for any $\qq^\infty\in S$ with $\sum_{i=1}^Bq_i^\infty<\infty$, and $\varepsilon>0$, if $\qq(0) = \qq^\infty$, then there exists a time $t_2 = t_2(\qq^\infty, \varepsilon)<\infty$, such that $q_{K+1}(t)\leq f+\varepsilon$ for all $t\geq t_2$.
Choosing $\varepsilon = (1-f)/2$ say, for all $t\geq t_2$, $q_{K+1}(t)<1$, and thus $\sum_{i=1}^{K+1}p_{i-1}(\qq(t))=1$, i.e., $\sum_{i=K+2}^B p_{i-1}(\qq(t))=0$.
Observe that
\begin{align*}
q_{\geq K+2}(t)&= q_{\geq K+2}(t_2)+\lambda\int_{t_2}^t\sum_{i=K+2}^B p_{i-1}(\qq(s))\dif s- \int_{t_2}^tq_{\geq K+2}(s)\dif s\\
&=  q_{\geq K+2}(t_2)- \int_{t_2}^tq_{\geq K+2}(s)\dif s\qquad \mbox{for all}\quad t\geq t_2,
\end{align*}
which implies $q_{\geq K+2}(t)\leq q_{\geq K+2}(t_2)\e^{-(t-t_2)}$.
Since $q_{\geq K+2}(t_2)\leq q_{\geq K+2}(0)+\lambda t_2<\infty$, we obtain that $q_{\geq K+2}(t)$ and thus $q_{K+2}(t)$ converges to 0 as $t\to\infty$.
This completes the proof of global stability of the fixed point. 
\end{proof}

\begin{proof}[Proof of Proposition~\ref{prop:interchange-mor}.]
Observe that $\pi^{ d(N)}$ is defined on $S$, and $S$ is a compact set when endowed with the product topology. 
Prohorov's theorem implies that the sequence of measures $\big\{\pi^{ d(N)}\big\}_{N\geq 1}$ is relatively compact, and hence, has a convergent subsequence. 
Let $\big\{\pi^{ d(N_n)}\big\}_{n\geq 1}$ be a convergent subsequence, with $\big\{N_n\big\}_{n\geq 1}\subseteq\N$, such that $\pi^{ d(N_n)}\dto\hat{\pi}$. We show that $\hat{\pi}$ is unique and equals the measure $\pi^\star = \delta_{\qq^\star}.$

First of all note that if $\qq^{ d(N_n)}(0)\sim\pi^{ d(N_n)}$, then $\qq^{ d(N_n)}(t)\sim\pi^{ d(N_n)}$ for all $t\geq 0$. 
Also, the fact that $\qq^{ d(N_n)}(t)\dto\qq(t)$, and $\pi^{ d(N_n)}\dto\hat{\pi}$, means that $\hat{\pi}$ is an invariant distribution of the deterministic process $\big\{\qq(t)\big\}_{t\geq 0}$. 
If $B<\infty$, then clearly, $\sum_{i=1}^Bq_i(0)<\infty$ with probability 1.
Also, if $B=\infty$, then observe that for any $N\geq 1$, the total number of active tasks in the system under the JSQ$(d(N))$ scheme behaves as that in an M/M/$\infty$ system.
Since $\lambda(N)/N\to \lambda<\infty$, this implies that $\sum_{i=1}^Bq^{d(N)}_i(\infty)\to \lambda<\infty$, where $\qq^{d(N)}(\infty)$ is the steady-state occupancy state of the system under the JSQ$(d(N))$ scheme.
Thus again,  $\sum_{i=1}^Bq_i(0)<\infty$ with probability 1.
This in conjunction with the global stability in Lemma~\ref{lem:global-stab-mor} implies that $\hat{\pi}$ must be the fixed point of the fluid limit.  
Since the latter fixed point is unique and equals $\qq^\star$, we can conclude the desired convergence of the stationary measure. 
\end{proof}

\section{Diffusion limit of JSQ: Non-integral \texorpdfstring{$\boldsymbol{\lambda}$}{lambda}}\label{sec:non-itegral}
In this section we establish the diffusion-scale behavior of the ordinary JSQ policy in the case when $\lambda$ is not an integer, i.e., $f>0$.
Recall that $f(N)=\lambda(N)-KN.$
In this regime, let us define the following centered and scaled processes:
\begin{equation}
\begin{split}
\bar{Q}^N_i(t)&=N-Q^N_i(t)\geq 0\quad \mathrm{for}\quad i\leq K-1,\\\\
\bar{Q}_{K}^N(t)&:=\frac{N-Q_K^N(t)}{\log (N)}\geq 0,\\\\
\bar{Q}_{K+1}^N(t)&:=\frac{Q^N_{K+1}(t)-f(N)}{\sqrt{N}}\in\R,\\\\
\bar{Q}^N_i(t)&:=Q^N_i(t)\geq 0\quad \mathrm{for}\quad i\geq K+2.
\end{split}
\end{equation}
\begin{theorem}[{Diffusion limit for JSQ policy; $f>0$}]\label{th:diffusion}
Assume that $\bar{Q}^N_i(0)\to \bar{Q}_i(0)$ in $\R$, $i\geq 1$, and $\lambda(N)/N\to\lambda>0$ as $N\to\infty$, with $f=\lambda-\lfloor\lambda\rfloor>0$, then 
\begin{enumerate}[{\normalfont(i)}]
\item $\lim_{N\to\infty}\Pro{\sup_{t\in[0,T]}\bar{Q}_{K-1}^N(t)\leq 1}=1$, and $\big\{\bar{Q}^N_i(t)\big\}_{t\geq 0}\dto \big\{\bar{Q}_i(t)\big\}_{t\geq 0}$, where $\bar{Q}_i(t)\equiv 0$, provided that 
$\lim_{N\to\infty}\Pro{\bar{Q}_{K-1}^N(0)\leq 1}=1$, and
$\bar{Q}_i^N(0)\to 0$ for $i\leq K-2$.
\item $\big\{\bar{Q}^N_K(t)\big\}_{t\geq 0}$ is a stochastically bounded sequence of processes in $D_{\R}[0,\infty)$.
\item $\big\{\bar{Q}^N_{K+1}(t)\big\}_{t\geq 0}\dto \big\{\bar{Q}_{K+1}(t)\big\}_{t\geq 0}$, where $\bar{Q}_{K+1}(t)$ is given by the Ornstein-Uhlenbeck process satisfying the following stochastic differential equation:
$$d\bar{Q}_{K+1}(t)=-\bar{Q}_{K+1}(t)dt+\sqrt{2\lambda}dW(t),$$
where $W(t)$ is the standard Brownian motion,
provided that $\bar{Q}_{K+1}^N(0) \to  \bar{Q}_{K+1}(0)$ in $\mathbb{R}$.
\item For $i\geq K+2$, $\big\{\bar{Q}^N_i(t)\big\}_{t\geq 0}\dto \big\{\bar{Q}_i(t)\big\}_{t\geq 0}$, where $\bar{Q}_i(t)\equiv 0$, provided that $\bar{Q}_i^N(0)\to 0$.\\
\end{enumerate}
\end{theorem}
Note that statements~(i) and~(ii) in Theorem~\ref{th:diffusion} imply statement~(i) in Theorem~\ref{th:diff pwr of d 1}, for the JSQ policy, while (iii) and (iv) in Theorem~\ref{th:diffusion} are equivalent with statements (ii) and (iii) in Theorem~\ref{th:diff pwr of d 1}.
In view of the universality result in Corollary~\ref{cor-diff}, it thus suffices to prove Theorem~\ref{th:diffusion}.

The rest of this section is devoted to the proof of Theorem~\ref{th:diffusion}.
From a high level, the idea of the proof is the following. 
Introduce 
\begin{equation}\label{eq:pos-neg}
\begin{split}
Y^N(t):=\sum_{i=1}^B Q_i^N(t),\quad
D^N_+(t):=\sum_{i=1}^{K}(N-Q_i^N(t)),\quad
D^N_-(t):=\sum_{i=K+2}^{B}Q_i^N(t).
\end{split}
\end{equation}
and observe that
\begin{align*}
Q^N_{K+1}(t)+KN &= \sum_{i=1}^B Q_i^N(t)+\sum_{i=1}^{K} (N-Q_i^N(t)) - \sum_{i=K+2}^B Q_i^N(t)\\
&= Y^N(t) + D^N_+(t) -D^N_-(t).
\end{align*}
In Proposition~\ref{prop:positive} we show
that on any finite time interval, the sequence of processes $\big\{D^N_+(t)\big\}_{t\geq 0}$ is $\Op(\log(N))$, which implies
 that the number of server pools with fewer than $K$ active tasks is negligible on $\sqrt{N}$-scale.
 Furthermore, in Proposition~\ref{prop:negative} we prove that
since $\lambda<B$ the number of  tasks that are assigned to  server pools with at least $K+1$ tasks converges to zero in probability
 and hence, for a suitable starting state, $\big\{D^N_-(t)\big\}_{t\geq 0}$ converges 
 to the zero process.
 As we will show, this also means that $Y^N(t)$ behaves with high
 probability as the total number of tasks in an M/M/$\infty$ system. 
Therefore with the help of the following diffusion limit result for the M/M/$\infty$ system in~\cite[Theorem 6.14]{Robert03}, we conclude the proof of statement~(iii) of Theorem~\ref{th:diffusion}.
\begin{theorem}[{\cite[Theorem 6.14]{Robert03}}]
\label{th:robert-book-mmn-mor}
Let $\big\{Y^N_\infty(t)\big\}_{t\geq 0}$ be the total number of tasks in an M/M/$\infty$ system with arrival rate $\lambda (N)$ and unit-mean service time. 
If $(Y^N_\infty(0)-\lambda (N))/\sqrt{N}\to v\in\R$, then the process $\big\{\bar{Y}^N_{\infty}(t)\big\}_{t\geq 0}$, with
$$\bar{Y}^N_{\infty}(t)=\frac{Y^N_\infty(t)-\lambda( N)}{\sqrt{N}},$$
converges weakly to an Ornstein-Uhlenbeck process $\big\{X(t)\big\}_{t\geq 0}$ 
described by the stochastic differential equation
\begin{align*}
X(0)=v,\qquad \dif X(t)& = -X(t)\dif t + \sqrt{2\lambda}\dif W(t).\\
\end{align*}
\end{theorem}
The next two propositions state the asymptotic properties of  $\big\{D^N_+(t)\big\}_{t\geq 0}$ and $\big\{D_-^N(t)\big\}_{t\geq 0}$ mentioned before, which play a crucial role in the proof of Theorem~\ref{th:diffusion}. 
Let $B_{K+1}^N(t)$ be the cumulative number of tasks up to time $t$ that are assigned to some server pool having at least $K+1$ active tasks if $B>K+1$, and that are lost if $B=K+1$.
\begin{proposition}\label{prop:negative}
Under the assumptions of Theorem~\ref{th:diffusion}, for any $T\geq 0$, $B_{K+1}^N(T)\pto 0$, and consequently,
$\sup_{t\in[0,T]}D_-^N(t)\pto 0$ as $N\to\infty,$ 
provided $D_-^N(0)\pto 0$.
\end{proposition}
Informally speaking, the above proposition implies that for large $N$, there will be almost no server pool with $K+2$ or more tasks in any finite time horizon, if the system starts with no server pools with more than $K+1$ tasks. 
The next proposition shows that the number of server pools having fewer than $K$ tasks is of order $\log (N)$ in any finite time horizon.

\begin{proposition}\label{prop:positive}
Under the assumptions of Theorem~\ref{th:diffusion}, the sequence of processes $\left\{D^N_+(t)/\log (N)\right\}_{t\geq 0}$ is stochastically bounded in $D_{\mathbb{R}}[0,\infty)$, provided that the sequence of random variables $\left\{D^N_+(0)/\log (N)\right\}_{N\geq 1}$ is tight.
\end{proposition}
Before providing the proofs of the above two propositions, we first prove Theorem~\ref{th:diffusion} using Propositions~\ref{prop:negative} and~\ref{prop:positive}.
\begin{proof}[Proof of Theorem~\ref{th:diffusion}.]
First observe that (iv) and (ii) immediately follows from Propositions~\ref{prop:negative} and~\ref{prop:positive}, respectively.

To prove (i), fix any $T\geq 0$.
We will show that
\begin{equation}
\lim_{N\to\infty}\Pro{\sup_{t\in[0,T]}\sum_{i=1}^{K-1}\bQ_{i}^N(t)\leq 1}=1.
\end{equation} 
 Since $\bQ^N_{i}\leq 1$ implies that $\bQ^N_{i-1}\leq 1$ for $i= 2,\ldots,K$, this then completes the proof of~(i).
Note that the process $\sum_{i=1}^{K-1}\bQ_{i}^N(\cdot)$ increases by one when there is a departure from some server pool with at most $K-1$ active tasks, and if positive, decreases by one whenever there is an arrival.
Therefore it can be thought of as a
birth-death process with state-dependent instantaneous birth rate $\sum_{i=1}^{K-1}i(Q^N_i(t)-Q^N_{i+1}(t))$, and constant instantaneous death rate $\lambda (N)$. 
Observe that
\begin{align*}
\sum_{i=1}^{K-1}i(Q^N_i(t)-Q^N_{i+1}(t))&=\sum_{i=1}^{K-1}Q^N_i(t)-(K-1)Q^N_{K}(t)\leq (K-1)(N-Q_K^N(t)),
\end{align*}
and due to (ii), we can claim that for \textit{any} nonnegative sequence $\ell(N)$ diverging to infinity,
\begin{equation}
\lim_{N\to\infty}\Pro{\sup_{t\in[0,T]}(N-Q^N_K(t))\leq \ell (N)\log (N)}=1.
\end{equation}
Thus, given $\sup_{t\in[0,T]}(N-Q^N_K(t))\leq \ell (N)\log (N),$ in the interval $[0,T]$, the process $\sum_{i=1}^{K-1}\bQ_{i}^N(t)$ 
is stochastically upper bounded by a birth-and-death process $Z^N(t)$ with birth rate $(K-1)\ell (N)\log (N)$ and death rate $\lambda( N).$
Consequently,
\begin{equation}\label{eq:lNlogN}
\begin{split}
\Pro{\sup_{t\in[0,T]}\sum_{i=1}^{K-1}\bQ_{i}^N(t)> 1}&\leq 
\Pro{\sup_{t\in[0,T]}Z^N(t)> 1}\\
&+ \Pro{\sup_{t\in[0,T]}(N-Q^N_K(t))> \ell (N)\log (N)}.
\end{split}
\end{equation}
Let $\big\{\eta^N(n)\big\}_{n\geq 1}$ denote the discrete uniformized chain of the upper bounding birth-death process. 
Also, let $K_N(t)$ denote the number of jumps taken up to time $t$ by $\big\{\eta^N(n)\big\}_{n\geq 1}$. 
Since the jump rate of the process is $O(N)$, we have for \textit{any} nonnegative sequence $\ell^0(N)$ diverging to infinity, and for any $T\geq 0$, 
$$\lim_{N\to\infty}\Pro{K_N(T)\leq N\ell^0(N)}=1.$$
Considering the Markov chain $\eta^N(\cdot)$, the probability of one birth is bounded from above by 
$$p_{Q_K^N}=\frac{(K-1)\ell (N)\log (N)}{N+(K-1)\ell (N)\log (N)}.$$
Now, $Z^N(\cdot)$ will exceed 1 if and only if there are at least two successive births. Hence,
\begin{equation}\label{eq:split}
\begin{split}
&\Pro{\sup_{t\in[0,T]}Z^N(t)\leq 1}=\Pro{\sup_{n\leq K_N(T)}\eta^N(n)\leq 1}\\
&\hspace{2cm}\geq \Pro{\sup_{n\leq N\ell^0(N)}\eta^N(n)\leq 1}\Pro{K_N(T)\leq N\ell^0(N)}.
\end{split}
\end{equation}
Again we can write the first term of the last inequality above as
\begin{align*}
&\Pro{\sup_{n\leq N\ell^0(N)}\eta^N(n)\leq 1}
\geq \left(1-\left(\frac{(K-1)\ell(N)\log (N)}{N+(K-1)\ell(N)\log (N)}\right)^2\right)^{N\ell^0(N)}.
\end{align*}
If we choose $\ell(N)$ and $\ell^0(N)$, both diverging to infinity, such that $$\ell(N)^2\ell^0(N)\log (N)/ N\to 0\quad \mbox{as} \quad N\to\infty, $$ then the expression on the right side of~\eqref{eq:split} converges to 1, and consequently, the right of~\eqref{eq:lNlogN} converges to~0 (one can see that this choice is always feasible).
 Hence the proof of (i) is complete.

For (iii), recall that $Y^N_\infty(t)$ denotes the total number of tasks in an M/M/$\infty$ system with arrival rate $\lambda(N)$ and exponential service time distribution with unit mean.
Also, Proposition~\ref{prop:negative} implies that under the assumptions of the theorem, in any finite time horizon, with high probability there will be no arrival to a server pool with $K+1$ or more active tasks.
Now observe that since $B\geq K+1$, for any $T\geq 0$,
\begin{align*}
\Pro{\exists\ t\in[0,T]: Y^N(t) \neq Y^N_\infty(t)}
\leq \Pro{\exists\ t\in[0,T]: B^N_{K+1}(t)\geq 1}\to 0,
\end{align*}
as $N\to\infty$.
Propositions~\ref{prop:negative} and~\ref{prop:positive} then yield
\begin{align*}
&\sup_{t\in[0,T]}\frac{1}{\sqrt{N}}\left|Q_{K+1}^N(t)-f(N)-(Y^N_\infty(t)-\lambda (N))\right|\\
= & \sup_{t\in[0,T]}\frac{1}{\sqrt{N}}\Big|\sum_{i=1}^B Q_{i}^N(t)+\sum_{i=1}^K(N- Q_{i}^N(t))-\sum_{i=K+2}^BQ^N_i(t)\\
&\hspace{6cm}-KN-f(N)-(Y^N_\infty(t)-\lambda (N))\Big| \\
= & \sup_{t\in[0,T]}\frac{1}{\sqrt{N}}\left[Y^N(t)-Y^N_\infty(t)
+D_N^+(t)-D^N_-(t)\right]\to 0,
\end{align*}
as $N\to\infty$,
which in conjunction with \cite[Theorem 6.14]{Robert03}, as mentioned earlier, gives the desired diffusion limit. 
\end{proof}

\begin{proof}[Proof of Proposition~\ref{prop:negative}.]
Couple the M/M/$\infty$ system and a system under the ordinary JSQ policy in the natural way, until an overflow event occurs in the latter system.
Fix any $\varepsilon>0$ with $\lambda+\varepsilon<K+1$.
Observe that the event 
$$\left[\sup_{t\in[0,T]}B^N_{K+1}(t)>0\right]$$ will occur only if for some $t'\leq T$,
some arriving task is assigned to a server pool with more than $K$ active tasks,
and in that case, there exists $t''\leq t'$, such that $Y^N(t'')> (\lambda+\varepsilon)N$.
Since, for any $t\in[0,t'']$, $Y^N(t) = Y^N_\infty(t)$, 
we have
\begin{equation}\label{eq:neg_part}
\begin{split}
&\sup_{t\in[0,T]}B^N_{K+1}(t)> 0\\
&\implies\sup_{t''\in[0,t']}Y^N(t'')\geq  (\lambda+\varepsilon)N\\
&\implies\sup_{t''\in[0,t']}Y^N_\infty(t'')\geq  (\lambda+\varepsilon)N\\
&\implies\sup_{t\in[0,T]}(Y^N_\infty(t)-\lambda( N))> \varepsilon N +o(N)\\
&\implies\sup_{t\in[0,T]}\frac{1}{\sqrt{N}}(Y^N_\infty(t)-\lambda (N))> \varepsilon\sqrt{N}+o(\sqrt{N}).
\end{split}
\end{equation}
From Theorem 6.14 of \cite{Robert03}, we know that the process $\big\{(Y^N(t)-\lambda (N))/\sqrt{N}\big\}_{t\geq 0}$ is stochastically bounded. 
Hence, Equation~\eqref{eq:neg_part} yields that $\sup_{t\in[0,T]}B^N_{K+1}(t)$ converges to zero in probability as $N\to\infty$ for any $T\geq 0$.
Consequently, from the assumption of Theorem~\ref{th:diffusion} that $D^N_-(0)\pto 0$, the conclusion $\sup_{t\in[0,T]}D^N_-(t)\pto 0$, is immediate.
\end{proof}

\begin{proof}[Proof of Proposition~\ref{prop:positive}.]
Observe that $\sum_{i=1}^K(N-Q_i^N(\cdot))$ increases by one when there is a departure from some server pool with at most $K$ active tasks, and if positive, decreases by one whenever there is an arrival.
Therefore the process $\big\{D^N_+(t)\big\}_{t\geq 0}$ increases by one at rate 
$\sum_{i=1}^K i(Q_i(t)-Q_{i+1}(t))=\sum_{i=1}^K (Q_i(t)-Q_{K+1}(t))$, and while positive, decreases by one at  constant rate $\lambda(N)$. 
Now, to prove stochastic boundedness of the sequence of processes $\big\{D^N_+(t)/\log (N)\big\}_{t\geq 0}$, we will show that for any fixed $T\geq 0$ and any function $\ell(N)$ diverging to infinity (i.e., such that $\ell(N)\to\infty$ as $N\to\infty$), 
\begin{equation}\label{eq:conv_prob}
\Pro{\sup_{t\in[0,T]}D^N_+(t)>\ell(N)\log (N)}\to 0.
\end{equation}

Let $\big\{X^N(n)\big\}_{n\geq 0}$ be the discrete jump chain, and $K_N(t)$ be the number of jumps before time $t$, of the process $\big\{D_+^N(t)\big\}_{t\geq 0}$. Hence, for any fixed $T\geq 0$,
\begin{equation}\label{eq:condition}
\begin{split}
&\Pro{\sup_{t\in[0,T]}D^N_+(t)>\ell(N)\log (N)}\\
&=\Pro{\sup_{n\leq K_N(T)}X^N(n)>\ell(N)\log (N)}\\
&\leq \Pro{\sup_{n\leq N\ell_0(N)}X^N(n)>\ell(N)\log (N)}\Pro{K_N(T)\leq N\ell_0(N)}\\
&\qquad+\Pro{K_N(T)> N\ell_0(N)},
\end{split}
\end{equation}
for some function $\ell_0(N):\mathbb{N}\to\mathbb{N}$, to be chosen according to Lemma~\ref{lem:discr-walk} below. 
Now, observe that $K_N(T)$ is upper bounded by a Poisson random variable with parameter $\lambda(N)T+\int_0^T \sum_{i=1}^K (Q_i(s)-Q_{K+1}(s))ds$, and $\sum_{i=1}^K (Q_i(s)-Q_{K+1}(s))\leq KN$. 
Hence for any function $\ell_0(N)$ diverging to infinity, we have 
$$\Pro{K_N(T)> N\ell_0(N)}\to 0.$$

To control the first term, it is enough to note that $\sum_{i=1}^K (Q_i(t)-Q_{K+1}(t))\leq KN<\lambda N$. Hence the process $\big\{X^N(n)\big\}_{n\geq 1}$ can be stochastically upper bounded by the process $\big\{\hat{X}^N(n)\big\}_{n\geq 1}$, defined as follows:
\begin{equation}\label{eq:upperboundingbd}
\hat{X}^N(n+1)=
\begin{cases}
\hat{X}^N(n)+1&\mbox{ with prob. }K/(K+\lambda),\\
(\hat{X}^N(n)-1)\vee 0&\mbox{ with prob. }\lambda/(K+\lambda),
\end{cases}
\end{equation}
Therefore, combining Lemma~\ref{lem:discr-walk} below for the above Markov process $\big\{\hat{X}^N(n)\big\}_{n\geq 0}$ with Equation~\eqref{eq:condition} we obtain Equation~\eqref{eq:conv_prob}. Hence the proof is complete. 
\end{proof}
\begin{lemma}\label{lem:discr-walk}
For any function $\ell(N):\mathbb{N}\to\mathbb{N}$, diverging to infinity, there exists another function $\ell_0(N):\mathbb{N}\to\mathbb{N}$, diverging to infinity, such that 
$$\Pro{\sup_{n\leq N\ell_0(N)}\hat{X}^N(n)>\ell(N)\log (N)}\to 0.$$
\end{lemma}
\begin{proof}
We will use a regeneration approach to prove the lemma. Let $p:=K/(K+\lambda)$. 
Note that then $p<q:=1-p$. 
Define the $i^{\mathrm{th}}$ regeneration time $\rho_i$ of the Markov chain as follows: $\rho_0=0$, and $\rho_i:=\min\big\{k>\rho_{i-1}:\hat{X}_k=0\big\}$, for $i\geq 1$. 
Also define, $m_i:=\max\big\{\hat{X}_k: \rho_{i-1}\leq k<\rho_i\big\}$, for $i\geq 1$, and $\xi(n):=\min\big\{i:\rho_i\geq n\big\}$, for $n\geq 1$. 
Now observe that \cite[XIV.2]{feller1}, 
\begin{equation}
\Pro{m_i\geq M}=p\times\frac{\frac{q}{p}-1}{\left(\frac{q}{p}\right)^M-1}\leq a^{-M},
\end{equation}
for some $a>1$, since $q/p>1$. 
Thus the tail of the distribution of the maximum attained in one regeneration period decays exponentially. Recall that, in $n$ steps the Markov chain exhibits $\xi(n)$ regenerations. Hence, for any $\ell_0(N)$ and $\ell(N)$,
\begin{equation}\label{eq:discrete-maximum}
\begin{split}
&\Pro{\sup_{n\leq N\ell_0(N)}\hat{X}^N(n)>\ell(N)\log (N)}=\Pro{\sup_{i\leq \xi(N\ell_0(N))}m_i>\ell(N)\log (N)}\\
&\leq 1-\left(1-a^{-\ell(N)\log (N)}\right)^{\xi(N\ell_0(N))}\leq 1-\left(1-a^{-\ell(N)\log (N)}\right)^{N\ell_0(N)}.
\end{split}
\end{equation}

Now, for given $\ell(N)$, choose $\ell_0(N)$ diverging to infinity, such that 
$$N\ell_0(N) a^{-\ell(N)\log (N)}\to 0\quad \text{as}\quad N\to\infty.$$ 
Since the condition is equivalent to 
$$\log (N)+\log(\ell_0(N))-\ell(N)\log( a)\log (N)\to-\infty,$$ 
it is evident that such a choice of $\ell_0(N)$ is always possible. Hence, for such a choice of $\ell_0(N)$ the probability in Equation~\eqref{eq:discrete-maximum} converges to zero and the proof is complete. 
\end{proof}

\section{Diffusion limit of JSQ: Integral \texorpdfstring{$\boldsymbol{\lambda}$}{lambda}}\label{sec:integral}
In this section we analyze the diffusion-scale behavior of the ordinary JSQ policy when
$\lambda$ is an integer, i.e., $f=0$, and
$$\frac{KN-\lambda(N)}{\sqrt{N}}\to\beta,\quad \text{as}\quad N\to\infty,$$ 
with $\beta\in\R$ being a fixed real number. 
Throughout this section we assume $B=K+1.$
Thus, tasks that arrive when all the server pools have $K+1$ active tasks, are permanently discarded. 
For brevity, define $Z^N_1(t)=\sum_{i=1}^K (N-Q_i^N(t))$ and $Z^N_2(t):=Q_{K+1}^N(t)$. 
Note that $Z^N_1(t)$ corresponds to $D^N_+(t)$ in the previous section.
Also recall \eqref{eq:scaling-f=0}, and define
\begin{equation}
\begin{split}
\zeta_1^N(t)&:=\frac{Z_1^N(t)}{\sqrt{N}}=\hQ_{K-1}^N(t)+\hQ_{K}^N(t)\\\\
\zeta_2^N(t)&:=\frac{Z_2^N(t)}{\sqrt{N}}=\hQ_{K+1}^N(t),
\end{split}
\end{equation}
with $\hQ_{K-1}^N(t)$, $\hQ_{K}^N(t)$, and $\hQ_{K+1}^N(t)$ as in~\eqref{eq:scaling-f=0}.

\begin{theorem}\label{th: f=0 diffusion-mor}
Assume that $(\zeta^N_1(0),\zeta^N_2(0))\to (\zeta_1(0),\zeta_2(0))$ in $\R^2$ as $N\to\infty$. 
Then the two-dimensional process $\big\{(\zeta^N_1(t),\zeta^N_2(t))\big\}_{t\geq 0}$ converges weakly to the process $\big\{(\zeta_1(t),\zeta_2(t))\big\}_{t\geq 0}$ in $D_{\R^2}[0,\infty)$ governed by the stochastic recursion equation:
\begin{align*}
\zeta_1(t)&=\zeta_1(0)+\sqrt{2K}W(t)-\int_0^t(\zeta_1(s)+K\zeta_2(s))\dif s+\beta t+V_1(t),\\
\zeta_2(t)&=\zeta_2(0)+V_1(t)-(K+1)\int_0^t\zeta_2(s)\dif s,
\end{align*}
where $W$ is the standard Brownian motion, and $V_1(t)$ is the unique non-decreasing process in $D_{\R_+}[0,\infty)$ satisfying
\begin{align*}
\int_0^t\ind{\zeta_1(s)\geq 0}\dif V_1(s)=0.\\
\end{align*}
\end{theorem}

\begin{remark}\normalfont
Note that $Y^N(t)-KN = Z_2^N(t) - Z_1^N(t)$. Thus, under the assumption in~\eqref{eq:f=0}, the diffusion limit in Theorem~\ref{th: f=0 diffusion-mor} implies that 
\begin{align*}
\frac{Y^N(\cdot)-\lambda(N)}{\sqrt{N}} = \frac{Y^N(\cdot)-KN}{\sqrt{N}} + \frac{KN-\lambda(N)}{\sqrt{N}}\dto \zeta_2(\cdot)-\zeta_1(\cdot) + \beta .
\end{align*}
Writing $X(t) = \zeta_2(t) -\zeta_1(t) -\beta$, from Theorem~\ref{th: f=0 diffusion-mor}, one can note that the process $\big\{X(t)\big\}_{t\geq 0}$ satisfies 
$$\dif X(t) = -X(t)\dif t - \sqrt{2K}\dif W(t),$$
which is consistent with the diffusion-level behavior of $Y^N(\cdot)$ stated in Theorem~\ref{th:robert-book-mmn-mor}. 
\end{remark}
Next, using the arguments in the proof of Proposition~\ref{prop:positive} one can see that the process
$$\sum_{i=1}^{K-1}\frac{N-Q_i^N(\cdot)}{\sqrt{N}}=\hQ^N_{K-1}(\cdot)\pto 0,$$
provided  $\hQ_{K-1}^N(0)\pto 0$.
Thus, Theorem~\ref{th: f=0 diffusion-mor} yields the diffusion limit for the ordinary JSQ policy in the case $B=K+1$.
The proof for $B>K+1$ then follows from exactly the same arguments as provided in~\cite[Section~5.2]{EG15}.
The idea is that since the process $Q_{K+1}^N(\cdot)$, when scaled by $\sqrt{N}$, is stochastically bounded, the probability that on any finite time interval, it will take value $N$ (or equivalently, all server pools will have at least $K+1$ active tasks) vanishes as $N$ grows large.
Therefore, the dynamics of the limit of $(\hQ_{K+2}^N(\cdot),\ldots,\hQ_M(\cdot))$ becomes deterministic, and the limit of $\hQ^N_{K+1}(\cdot)$ for $B>K+1$ becomes a transformation of the limit of $\hQ^N_{K+1}(\cdot)$ for $B=K+1$, as described in Theorem~\ref{th:diff pwr of d 2}. 
Hence, note that the diffusion limit in Theorem~\ref{th: f=0 diffusion-mor} is equivalent to the one in Theorem~\ref{th:diff pwr of d 2}.
In view of the universality result in Corollary~\ref{cor-diff}, it thus suffices to prove Theorem~\ref{th: f=0 diffusion-mor}.

We will use the reflection argument developed in \cite{EG15} to prove Theorem~\ref{th: f=0 diffusion-mor}. 
Observe that the evolution of $\big\{(Z_1^N(t),Z_2^N(t))\big\}_{t\geq 0}$ can be described by the following stochastic recursion which is explained in detail below.
\begin{equation}\label{eq: f=0 process}
\begin{split}
Z_1^N(t)&=Z_1^N(0)+A_1\left(\int_0^t(KN-Z_1^N(s)-KZ_2^N(s))\dif s\right)
-D_1(\lambda(N) t)+U_1^N(t),\\
Z_2^N(t)&=Z_2^N(0)+U_1^N(t)-D_2\left(\int_0^t(K+1)Z_2^N(t)\dif s\right)-U_2^N(t),
\end{split}
\end{equation}
where $A_1$, $D_1$ and $D_2$ are unit-rate Poisson processes, and
\begin{equation}
\begin{split}
U_1^N(t)&=\int_0^t\ind{Z_1^N(s)=0}\dif D_1(\lambda(N)s),\\
U_2^N(t)&=\int_0^t\ind{Z_2^N(s)=C\sqrt{N}}\dif D_1(\lambda(N)s).
\end{split}
\end{equation}

The components of \eqref{eq: f=0 process} can be explained as follows. 
The process $Z_1(t)$ increases by one when a departure occurs from a server pool with  at most $K$ active tasks, and it decreases by one when an arriving task is assigned to a server pool with at most $K$ active tasks. 
Hence the instantaneous rate of increase at time $s$ is given by 
\begin{align*}
\sum_{i=1}^{K}i(Q^N_i(t)-Q^N_{i+1}(t))&=\sum_{i=1}^{K}Q^N_i(t)-KQ^N_{K+1}(t)\\
&=KN - \sum_{i=1}^{K}(N-Q^N_i(t))-KQ^N_{K+1}(t) \\
&= KN-Z_1^N(t)-KZ_2^N(t),
\end{align*}
and the instantaneous rate of decrease is given by the arrival rate $\lambda(N)$. 
But $Z_1^N$ cannot be negative, and hence the arrivals when $Z_1^N$ is zero, add to $Z_2^N$, and the rate of increase of the $Z_2^N$ process is given by the overflow process $U_1^N$. 
Since $B=K+1$, the rate of decrease of $Z_2^N$ equals the total number of tasks at server pools with exactly $K+1$ tasks, which is given by $(K+1)Z_2^N$. 
This explains the rate in the Poisson process $D_2(\cdot)$. 
Finally, since $Z_2^N$ is upper bounded by $N$, $U_2^N$ is the overflow of the $Z_2^N$ process with $C=\sqrt{N}$, i.e., the number of arrivals to the system when $Z_2^N=N$. 
The existence and uniqueness of the above stochastic recursion can be proved following the arguments in \cite[Section 2]{PTRW07}.

\paragraph{Martingale representation.}
We now introduce the martingale representation for \eqref{eq: f=0 process}, and following similar arguments as in \cite[Subsection 4.3]{EG15}, we obtain the following scaled, square integrable martingales with appropriate filtration:
\begin{equation}\label{eq:martingales1}
\begin{split}
M^N_{1,1}(t)&=\frac{1}{\sqrt{N}}A_1\left(\int_0^t(KN-Z_1^N(s)-KZ_2^N(s))\dif s\right)\\
&\hspace{5cm}-\frac{1}{\sqrt{N}}\int_0^t(KN-Z_1^N(s)-KZ_2^N(s))\dif s,\\
M_{1,2}^N(t)&=\frac{1}{\sqrt{N}}(D_1(\lambda(N) t)-\lambda(N)t),\\
M^N_{2,1}(t)&=\frac{1}{\sqrt{N}}D_2\left(\int_0^t(K+1)Z_2^N(t)\dif s\right)-\frac{K+1}{\sqrt{N}}\int_0^tZ_2^N(s)\dif s,
\end{split}
\end{equation}
with $V_1^N(t):=U_1^N(t)/\sqrt{N}$ and $V_2^N(t):=U_2^N(t)/\sqrt{N}$, and the predictable quadratic variation processes given by
\begin{equation}
\begin{split}
\langle M^N_{1,1}\rangle(t)&=\frac{1}{N}\int_0^t(KN-Z_1^N(s)-KZ_2^N(s))\dif s,\\
\langle M_{1,2}^N\rangle(t)&=\frac{\lambda(N)t}{N},\\
\langle M^N_{2,1}\rangle(t)&=\frac{K+1}{N}\int_0^tZ_2^N(s)\dif s.
\end{split}
\end{equation}
Therefore, we have the following martingale representation for~\eqref{eq: f=0 process}:
\begin{equation}\label{eq:martingale rep}
\begin{split}
\zeta^N_1(t)&=\zeta^N_1(0)+M_{1,1}^N(t)-M_{1,2}^N(t)
-\int_0^t(\zeta_1^N(s)+K\zeta_2^N(s))\dif s\\
&\hspace{5.8cm}+\frac{t(KN-\lambda(N))}{\sqrt{N}}+V_1^N(t),\\
\zeta_2^N(t)&=\zeta_2^N(0)+V_1^N(t)-M_{2,1}^N(t)
-(K+1)\int_0^t\zeta_2^N(s)\dif s -V_2^N(t).
\end{split}
\end{equation}

\paragraph{Convergence of independent martingales.}
We now show the convergence of the martingales defined in~\eqref{eq:martingales1} using the functional central limit theorem.
\begin{lemma}\label{lem:martingale convergence}
As $N\to\infty$,
$$\left\{\left(M^N_{1,1}(t), M^N_{1,2}(t), M^N_{2,1}(t)\right)\right\}_{ t\geq 0}\dto \left\{\left(\sqrt{K} W_1(t),\sqrt{K} W_2(t), 0\right)\right\}_{t\geq 0}$$ in $D_{\R^3}[0,\infty),$  where $W_1$, $W_2$ are independent standard Brownian motions.
\end{lemma}
\begin{proof}
From Theorem~\ref{fluidjsqd-mor} we know that for any fixed $T\geq 0$, 
$$\sup_{t\in [0,T]}Z_1^N(t)/N\pto 0 \quad\text{and}\quad\sup_{t\in [0,T]}Z_2^N(t)/N\pto 0.$$
This yields the following convergence results: 
\begin{equation}
\begin{split}
\langle M^N_{1,1}\rangle(T)&\pto KT,\\
\langle M_{1,2}^N\rangle(T)&\pto\lambda T=KT,\\
\langle M^N_{2,1}\rangle(T)&\pto 0.
\end{split}
\end{equation}
Then, using a random time change, the continuous-mapping theorem and the functional central limit theorem \cite[Theorem 4.2]{PTRW07}, \cite[Lemma 6]{EG15}, we get the convergence of the martingales. 
\end{proof}
Now we use the continuous-mapping theorem to prove the convergence of the processes described in~\eqref{eq:martingale rep}. To proceed in that direction, we need the following proposition, which is analogous to \cite[Lemma 1]{EG15}.

\begin{proposition}\label{th:continuous}
Let $B\in\bar{\R}_+$, $b\in\R^2$, $(y_1,y_2)\in D^2[0,\infty)$, and $(x_1,x_2)\in D^2[0,\infty)$ be defined by the following recursion: for $t\geq 0$,
\begin{equation}\label{eq:recursion1}
\begin{split}
x_1(t)&=b_1+y_1(t)+\int_0^t(-x_1(s)- Kx_2(s))\dif s+u_1(t),\\
x_2(t)&=b_2+y_2(t)+(K+1)\int_0^t(-x_2(s))\dif s+u_1(t)-u_2(t),
\end{split}
\end{equation}
where $u_1$ and $u_2$ are unique non-decreasing functions in $D$, such that
\begin{equation}\label{eq:reflection}
\begin{split}
\int_0^{\infty}\ind{x_1(s)>0}\dif u_1(t) & =0,\\
\int_0^{\infty}\ind{x_2(s)<B}\dif u_2(t) & =0.
\end{split}
\end{equation}
Then, $(x,u)$ is the unique solution to the above system. 
Furthermore, there exist functions $(f,g):(\bar{\R},\R^2, D^2_\R[0,\infty))\to (D^2_\R[0,\infty),D^2_\R[0,\infty))$ with $x=f(B,b,y)$ and $u=g(B,b,y)$, which are continuous when $\bar{\R}_+$ is equipped with the order topology, $D_\R[0,\infty)$ is equipped with the topology of uniform convergence over compact sets, and $(\bar{\R},\R^2, D^2_\R[0,\infty))$ and $(D^2_\R[0,\infty),D^2_\R[0,\infty))$ are equipped with the product topology.
\end{proposition}
The proof of the above proposition follows from similar arguments as described in the proof of \cite[Lemma 1]{EG15}, and hence is omitted.\\
\begin{proof}[Proof of Theorem~\ref{th: f=0 diffusion-mor}.]
Observe that the stochastic recursion equations described by \eqref{eq:martingale rep} fit in the framework of the recursion described by \eqref{eq:recursion1}, by taking $b_i=\zeta_i^N(0)$, $i=1,2$, $C=\sqrt{N}$, $y_1(t)=M^N_{1,1}(t)-M^N_{1,2}(t)+t(KN-\lambda(N))/\sqrt{N}$, and $y_2(t)=-M^N_{2,1}(t)$ for the $N^{\mathrm{th}}$ process.

By the assumptions of the theorem we have $\zeta_i^N(0)\dto \zeta_i(0)$, for $i=1,2$. Also, by Lemma~\ref{lem:martingale convergence}, 
$$\big\{(M^N_{1,1}(t), M^N_{1,2}(t), M^N_{2,1}(t))\big\}_{t\geq 0}\dto \big\{(\sqrt{K} W_1(t), \sqrt{K} W_2(t), 0)\big\}_{t\geq 0}.$$ 
Hence, for the limiting process, $y_1(t)=\sqrt{K} W_1(t)-\sqrt{K} W_2(t)+\beta t\equiv \sqrt{2K}W(t)+\beta t$ and $y_2(t)\equiv 0$.
Finally, using the continuous-mapping theorem we get the desired convergence as in the proof of~\cite[Theorem~2]{EG15}. 
\end{proof}

\section{Performance implications}\label{sec:performance}

\subsection{Evolution of the number of tasks at a tagged server pool}
We now provide some insights into the steady-state dynamics of the number of tasks at a particular server pool in the regime $d(N)\to\infty$ as $N\to\infty$.
Due to exchangeability of the server pools, asymptotically, the dynamics at a particular server pool depends on the system only through the mean-field limit, or the global system state averages.
Based on the fixed point~\eqref{eq:fixed point}, we claim (without proof) that the steady-state dynamics can be described as follows:
\begin{enumerate}[{\normalfont (i)}]
\item If a server pool contains $\lceil\lambda\rceil$ active tasks, then with high probability no further task will be assigned to it. 
\item Similarly, if a departure occurs from a server pool having $K=\lfloor\lambda\rfloor$ active tasks, a task will immediately be assigned to it.
\item Since the total flow of arrivals that join server pools with exactly $K$ active tasks, are
distributed uniformly among all such server pools, each server pool with exactly $K$ active tasks will observe an arrival rate $\lambda p_K(\qq^\star)/(q^\star_K-q^\star_{K+1}) = (K+1)f/(1-f)$.
\item Finally, the rate of departure from a server pool with $K+1$ active tasks is given by $K+1$.
\end{enumerate}
Let $S^{ d(N)}_k(t)$ denote the number of tasks at server pool~$k$ at time $t$ in the $N^{\mathrm{th}}$ system under the JSQ$(d(N))$ scheme.
Combining all the above, provided $d(N)\to\infty$ as $N\to\infty$, the process $\big\{S^{ d(N)}_k(t)\big\}_{t\geq 0}$ converges in distribution to the process $\big\{S(t)\big\}_{t\geq 0}$, described as follows:
\begin{enumerate}[{\normalfont (i)}]
\item If $f>0$, then $\big\{S(t)\big\}_{t\geq 0}$ is a two-state process, taking values $K$ and $K+1$, with transition rate from $K$ to $K+1$ given by $(K+1)f/(1-f)$, and from $K+1$ to $K$ given by $K+1$. 
So the steady-state distribution is $\Pro{S=K}=1-f$, and $\Pro{S=K+1}=f$, i.e., for $i\geq 1$, $\Pro{S=i}=q_i^\star-q_{i+1}^\star$, which agrees with the fixed point~\eqref{eq:fixed point} of the fluid limit.
\item If $f=0$, then $\big\{S(t)\big\}_{t\geq 0}$ is a constant process, taking value $\lambda=K$.
\end{enumerate}

\subsection{Evolution of the number of tasks observed by a tagged task}
To analyze the performance perceived by a particular tagged task with execution time $T$, observe that in steady state the probability that it will join a server pool with $i$ active tasks is given by
$p_i(\qq^\star)=K(1-f)/\lambda$ for $i=K-1$, $(K+1)f/\lambda$ for $i=K$, and 0 otherwise.
In the time interval $[0,T]$, the number of active tasks in the server pool it joins, is again a birth-death process $\big\{\hat{S}(t)\big\}_{0\leq t\leq T}$, whose dynamics is the same as of $\big\{S(t)\big\}_{t\geq 0}$ process conditioned on having one permanent task (i.e., its departure is not allowed).
Therefore, $\big\{\hat{S}(t)\big\}_{0\leq t\leq T}$ can be described as follows:
\begin{enumerate}[{\normalfont (i)}]
\item If $f>0$, then $\big\{\hat{S}(t)\big\}_{0\leq t\leq T}$ is a two-state process, taking values $K$ and $K+1$, with transition rate from $K$ to $K+1$ given by $(K+1)f/(1-f)$, and from $K+1$ to $K$ given by~$K$.
The steady-state distribution of the process is then given by
$\Pro{\hat{S}=K}=K(1-f)/\lambda$, and $\Pro{\hat{S}=K+1}=(K+1)f/\lambda$.
\item If $f=0$, then $\big\{\hat{S}(t)\big\}_{0\leq t\leq T}$ is a constant process, taking value $\lambda=K$.
\end{enumerate}
Observe in both of the above two cases that the initial distribution of $\hat{S}(t)$ coincides with its stationary distribution.
Now, if the performance perceived by the tagged task is measured as a function $h:\N\to \R$ of the number of concurrent tasks, then the relevant performance measure is given by 
\begin{align}
\E{\frac{1}{T}\int_0^T h(\hat{S}(t))\dif t}=\frac{1}{\lambda}((1-f)Kh(K)+f(K+1)h(K+1)),
\end{align}
independent of the execution time $T$.
Notice that if $h(x) = 1/(x+1),$ then the above performance measure becomes the constant $(K(K+2)-f)/((K+1)(K+2)).$

\subsection{Loss probabilities}
We now examine the asymptotic behavior of the loss probability when the buffer capacity at each server is $B<\infty$ and the arrival rate $\lambda(N)$ satisfies~\eqref{eq:f=0} with $K=B$.
We will establish lower and upper bounds, and prove that these asymptotically coincide.
When the buffer capacity $B$ is finite, to characterize the asymptotic steady-state loss probability of the JSQ$(d(N))$ scheme, we bound it from below and above by that of an ordinary and a modified Erlang loss system, respectively.
The lower and upper bounds rely on a stochastic comparison.

Suppose that $Y_1(t)$ and $Y_2(t)$ are two non-explosive, continuous-time Markov processes taking values in a complete separable metric space $E$. 
Let $X_1(t)$ and $X_2(t)$ be two birth-death processes defined on the same probability space, with finite state spaces $\big\{0,1,\ldots,n_1\big\}$ and $\big\{0,1,\ldots,n_2\big\}$, whose birth rates are $f_1(X_1(t), Y_1(t))$ and $f_2(X_2(t), Y_2(t))$, and death rates are $g_1(X_1(t), Y_1(t))$ and $g_2(X_2(t), Y_2(t))$, respectively.
\begin{lemma}\label{lem:bdprocess-stoch-comp}
If $n_1\leq n_2$, and for all $x\in \big\{0,1,\ldots,n_1\big\}$, $f_1(x,y_1)\leq f_2(x,y_2)$ and $g_1(x,y_1)\geq g_2(x,y_2)$, for all $y_1, y_2\in E$, then 
$\big\{X_1(t)\big\}_{t\geq 0}\leq_{st} \big\{X_2(t)\big\}_{t\geq 0}$, provided $X_1(0)\leq_{st} X_2(0)$.
\end{lemma}
\begin{proof}
The proof is fairly straightforward, but we present it briefly for the sake of completeness.
First we suitably couple the two processes, and then
as before, using the forward induction on event times, we show that the inequality holds throughout the sample path.
Define the processes $\big(X_1(\cdot),X_2(\cdot),Y_1(\cdot),Y_2(\cdot)\big)$ on the same probability space.
Due to the assumptions in the theorem, we do not need any condition on the evolution of $Y_1$ and $Y_2$, provided that they are defined on the same probability space.
Maintain two exponential clocks of rate $M_B:=\max\big\{f_1(x_1,y_1),f_2(x_2,y_2)\big\}$ (birth-clock) and $M_D:=\max\big\{g_1(x_1,y_1),g_2(x_2,y_2)\big\}$ (death-clock), respectively.
When the birth-clock rings, draw a single uniform$[0,1]$ random variable $u$ say, and a birth occurs in the $X_1$ process and $X_2$ process if $u\leq f_1(x_1,y_1)/M_B$ and $u\leq f_2(x_2,y_2)/M_B$, respectively.
Couple the deaths also, in a similar fashion.
Note that the processes thus constructed satisfy the relevant statistical laws in terms of the transition rates $f_1(x_1, y_1)$ and $f_2(x_2, y_2)$.

Now under the above coupling we prove the inequality. Assume that the inequality holds at event time $t_0$, and $X_1(t_0) = x_1$ and $X_2(t_0) = x_2$.
Note that if $x_1<x_2$, then trivially the inequality holds at the next event time $t_1$. 
Therefore, without loss of generality, assume $x_1= x_2=x\leq n_1$.
We will distinguish between two cases depending on whether the birth-clock or death-clock rings at time epoch $t_1$.
In the former case, observe that since $f_1(x,y_1)\leq f_2(x,y_2)$ for all $y_1, y_2\in E$, whenever there is a birth in the $X_1$ process, there will be a birth in the $X_2$ process as well.
Thus the inequality is preserved.
Alternatively, if the death-clock rings at time epoch $t_1$, then observe that since $g_1(x,y_1)\geq g_2(x,y_2)$ for all $y_1, y_2\in E$, whenever there is a death in the $X_2$ process, there will be a death in the $X_1$ process as well, and the inequality is preserved. This completes the proof.
\end{proof}

Denote by $\Er(C,\lambda)$ an Erlang loss system with
capacity $C$, load $\lambda$, and exponential service times with unit mean.
We further introduce a modified Erlang loss system $\hEr(n,d)$ with capacity $B(N-n)$, 
and arrival rate $\lambda$, with unit-exponential service times, where a fraction
$$p(n,d):=\left(1-\frac{n+1}{N}\right)^{d},$$ 
of tasks is rejected upfront,
independently of any other processes.
Note that the number of active tasks in the $\hEr(n,d)$ system evolves like an $\Er(B(N-n),\lambda p(n,d))$ system.

Define $C(N):=BN$, $\hat{C}(N):=B (N - n(N))$, and $\hat{\lambda}(N):= \lambda(N)p(n(N),d(N))$. 
Denote the total number of active tasks at time $t$ in the $N^{\mathrm{th}}$ system following the JSQ$(d(N))$ scheme, an $\Er(C(N),\lambda(N))$ system, and an $\hEr(n(N),d(N))$ system by $Y^{ d(N)}(t)$, $Y^N_\Er(t)$, and $Y^N_{\hEr}(t)$, respectively.
Denote the associated steady-state loss probabilities by $L^{ d(N)}$, $L(C,\lambda)$ and $\hat{L}(n,d)$, respectively.

\begin{lemma}\label{lem:stationary-blocking}
For all $N\geq 1$, $d(N)\geq 1$, and $n(N)<N$,
\begin{align*}
&\mathrm{\normalfont(a)}\quad \big\{Y^N_{\hEr}(t)\big\}_{t\geq 0}\leq_{st} \big\{Y^{ d(N)}(t)\big\}_{t\geq 0} \leq_{st} \big\{Y^N_{\Er}(t)\big\}_{t\geq 0},\\\\
&\mathrm{\normalfont(b)}\quad L(C(N),\lambda(N))\leq L^{ d(N)}\leq \hat{L}(n(N),d(N)).
\end{align*}
\end{lemma}
\begin{proof}
(a)
For the lower bound, observe that the rate of increase of the process $Y^{ d(N)}(\cdot)$ is at most that of the process $Y^N_\Er(\cdot)$, and the rate of decrease at any state is the same in both processes.
Thus, Lemma~\ref{lem:bdprocess-stoch-comp} implies that if both systems start from the same occupancy states, then $\big\{Y^{ d(N)}(t)\big\}_{t\geq 0}\leq_{st} \big\{Y^{N}_{\Er}(t)\big\}_{t\geq 0}$.
Consequently, in the steady state, $Y^{ d(N)}(\infty)\leq_{st} Y^N_{\Er}(\infty)$, and invoking Little's law yields 
$$L(C(N),\lambda(N))\leq L^{ d(N)}.$$

For the upper bound, 
first observe that at any arrival, as long as one of the $n(N)$ lowest-ordered server pools is sampled, which occurs with probability $1 - p(n(N), d(N))$, a task can only get lost when the total number of active tasks is at least $B (N - n(N))$.
Thus when the total number of active tasks $Y^{ d(N)}(\cdot)$ in the system under the JSQ$(d(N))$ scheme is $y$, the rate of increase of $Y^{ d(N)}(t)$ is at least $\lambda(N)(1 - p(n(N), d(N)))$ if $y \leq B (N - n(N))$,  and the rate of decrease is given by $y$.
Comparing with the modified Erlang loss system $\hEr(n(N),d(N))$ and using Lemma~\ref{lem:bdprocess-stoch-comp}, we obtain that if $Y^{ d(N)}(0)\geq_{st} Y^N_{\hEr}(0)$, then
$$\big\{Y^{ d(N)}(t)\big\}_{t\geq 0}\geq_{st} \big\{Y^N_{\hEr}(t)\big\}_{t\geq 0}.$$
The proof of the upper bound $L^{ d(N)}\leq \hat{L}(n(N),d(N))$ is then completed by again invoking Little's law.

(b) Little's law implies 
$$L^{ d(N)} = 1 - \frac{1}{\lambda(N)}\lim_{T\to\infty}\int_0^T Y^{ d(N)}(t)\dif t,$$
and similarly for the $\Er(C(N),\lambda(N))$ and $\hEr(n(N),d(N))$ systems. 
 Statement (b) then follows from statement (a). 
\end{proof}

The proposition below states that the limiting loss probability for the JSQ$(d(N))$ scheme vanishes as long as $d(N)\to\infty$. 
\begin{proposition}
 For any $\lambda\leq B$, if $d(N)\to\infty$ as $N \to\infty$, then $L^{ d(N)}\to 0$, as $N\to\infty.$
\end{proposition}
\begin{proof}
From~\eqref{eq:nN-order} and~\eqref{eq:equiv-nN-dN}, we know if $d(N)\to\infty$, then 
there exists $n(N)$ such that as $N\to\infty$, $n(N)/N\to 0$ and $p(n(N),d(N))\to 0$.
For such a choice of $n(N)$, $\lambda(N)/C(N)\to \lambda/B\leq 1$, and $\hat{\lambda}(N)/\hat{C}(N)\to \lambda/B\leq 1$ as $N\to\infty$.
Therefore, using Lemma~\ref{lem:stationary-blocking} and the standard results of the Erlang loss function~\cite{J74}, we complete the proof of the proposition. 
\end{proof}
\begin{remark} \normalfont
Note that in view of the results in~\cite{MKMG15, MMG15} for the JSQ$(d)$ schemes with fixed $d$, following the arguments as in Remark~\ref{rem:necessity-fluid}, the growth condition $d(N)\to\infty$ as $N\to\infty$ is also necessary to achieve an asymptotically zero probability of loss.
\end{remark}

We now further show that  the steady-state loss probability multiplied by $\sqrt{N}$ converges to a non-degenerate 
limit, which is the same as in an $\Er(C(N),\lambda(N))$ system.
The next theorem also establishes that if $d(N)/ (\sqrt{N}\log( N))\to 0$  as $N\to\infty$ and~\eqref{eq:f=0} is satisfied, then the steady-state loss probability is of higher order than $1/\sqrt{N}$. This indicates that  
the growth rate $\sqrt{N}\log (N)$ is not only sufficient but also nearly necessary. 

\begin{theorem}[{Scaled loss probability}]
\label{th:blocking}
Assume that $d(N)/ (\sqrt{N}\log (N))\to \infty$, as $N\to\infty$, and $\lambda(N)$ satisfies~\eqref{eq:f=0} with $K=B$. Then,  
\begin{equation}\label{eq:block-scale}
\lim_{N\to\infty}\sqrt{N}\ L^{ d(N)}=\frac{\phi(\beta)}{\sqrt{B}\Phi(\beta)},
\end{equation}
where $\phi(\cdot)$ and $\Phi(\cdot)$ are the density and distribution function of the standard Normal distribution, respectively.
\end{theorem}
Since the right side of~\eqref{eq:block-scale} corresponds to the 
asymptotic steady-state loss probability in an $\Er(C(N),\lambda(N))$ system~\cite{J74, B76, W84}, we thus conclude that~\eqref{eq:block-scale} is optimal on $\sqrt{N}$-scale in terms of loss probability.\\
\begin{proof}[Proof of Theorem~\ref{th:blocking}.]
 The idea again is to suitably bound the steady-state loss probability of the JSQ$(d(N))$ scheme.
Using Lemma~\ref{lem:stationary-blocking} and~\cite[Chapter 7, Theorem 15 (2)]{B76}, \cite{W84}, we obtain the lower bound as
\begin{equation}\label{eq:block-lower}
\begin{split}
&L^{ d(N)}\geq L(C(N),\lambda(N))\\
\implies & \varliminf_{N\to\infty}\sqrt{N}L^{ d(N)} \geq \varliminf_{N\to\infty} \sqrt{N} L(C(N),\lambda(N))= \frac{\phi(\beta)}{\sqrt{B}\Phi(\beta)}.
\end{split}
\end{equation}
For the upper bound, from~\eqref{eq:nN-order} and~\eqref{eq:equiv-nN-dN}, we know if $d(N)/(\sqrt{N}\log(N))\to\infty$ as $N\to\infty$, then there exists $n(N)$ with $n(N)/\sqrt{N}\to 0$ and 
\begin{equation}\label{eq:blockingupper2}
\sqrt{N}p(n(N),d(N))\to 0,\quad\mbox{as}\quad N\to\infty.
\end{equation}
Take such an $n(N)$. 
Again using~\cite[Chapter 7, Theorem 15 (2)]{B76}, we know that since
 as $N\to\infty$, $\hat{\lambda}(N)/\hat{C}(N)$ converges to one and $\hat{C}(N)/N$ converges to $B$,
\begin{equation}\label{eq:blockingupper3}
\lim_{N\to\infty}\sqrt{N}L(\hat{C}(N),\hat{\lambda}(N))= \frac{\phi(\beta)}{\sqrt{B}\Phi(\beta)}.
\end{equation}
Therefore, Lemma~\ref{lem:stationary-blocking}, and Equations~\eqref{eq:blockingupper2},~\eqref{eq:blockingupper3} yield
\begin{align*}
\varlimsup_{N\to\infty}\sqrt{N}L^{ d(N)} &\leq \varlimsup_{N\to\infty} \sqrt{N} L(C(N),\lambda(N))+\varlimsup_{N\to\infty}\sqrt{N}p(n(N),d(N))\\
&= \frac{\phi(\beta)}{\sqrt{B}\Phi(\beta)}.
\end{align*}
Combination of the lower bound in~\eqref{eq:block-lower} and the above upper bound completes the proof.
\end{proof}

\begin{remark}[{\normalfont Almost necessary condition for growth rate}]
\normalfont
It is worthwhile to mention that when $\lambda =K>0$ and $\lambda(N)$ satisfies~\eqref{eq:f=0}, the growth condition $d(N)/(\sqrt{N}\log(N))\to\infty$, as $N\to\infty$, is nearly necessary in order for the JSQ$(d(N))$ scheme to have the same diffusion limit as the ordinary JSQ policy.
More precisely, if  $d(N)/(\sqrt{N}\log(N))\to 0$ as $N\to\infty$, then the diffusion limit of the JSQ$(d(N))$ scheme differs from the ordinary JSQ policy.
In this remark we briefly sketch the outline of the proof. 
We will assume that the $d(N)$ server pools are chosen with replacement, to avoid cumbersome notation. But the proof technique and the result holds if the server pools are chosen without replacement.

\noindent
Assume on the contrary that as in the ordinary JSQ policy, if the centered and scaled initial occupancy state 
$N^{-1/2}(KN-\sum_{i=1}^KQ_i^{ d(N)}(0))$ is tight, then
 $N^{-1/2}(KN-\sum_{i=1}^KQ_i^{ d(N)}(t))$ is a stochastically bounded process.
 We argue that in this case, for any finite time $t$, the cumulative number of tasks joining a server with $K$ active tasks (or the cumulative number of lost tasks in case $K=B$) $L^{ d(N)}(t)$ does not scale with $\sqrt{N}$, and arrive at a contradiction. 
Indeed, $\big\{L^{ d(N)}(t)\big\}_{t\geq 0}$ admits the following martingale decomposition:
\begin{equation}\label{eq:Lmart-decomp}
L^{ d(N)}(t) = M_L^N(t) +\langle M_L^N\rangle(t),
\end{equation}
where $\big\{M_L^N(t)\big\}_{t\geq 0}$ is a martingale with compensator and predictable quadratic variation process given by 
$$\langle M_L^N\rangle(t)=\lambda(N) \int_0^t \left(Q_K^{ d(N)}(s-)/N\right)^{ d(N)}\dif s.$$
Since $\langle M_L^N\rangle(t)/N\leq \lambda t$, $\big\{M_L^N(t)/\sqrt{N}\big\}_{t\geq 0}$ is stochastically bounded.
We will show that $\langle M_L^N\rangle(t)$ is stochastically unbounded on $\sqrt{N}$-scale. 
From~\eqref{eq:Lmart-decomp}, this will imply that the process $\big\{L^{ d(N)}(t)/\sqrt{N}\big\}_{t\geq 0}$ is stochastically unbounded, which will complete the proof.
Note that 
\begin{align*}
Q_K^{ d(N)}(s)= N- (N-Q_K^{ d(N)}(s))\geq N- \sum_{i=1}^K(N-Q_i^{ d(N)}(s)),
\end{align*}
and hence,
\begin{align*}
 \langle M_L^N\rangle(t)&\geq \lambda(N) \int_0^t \left(1- \frac{1}{N}\sum_{i=1}^K(N-Q_i^{ d(N)}(s))\right)^{ d(N)}\dif s\\
&\geq \lambda(N) t \left(1- \frac{1}{N}\sup_{s\in [0,t]}\sum_{i=1}^K(N-Q_i^{ d(N)}(s))\right)^{ d(N)}.
\end{align*}
For any $T\geq 0$, since $\sup_{t\in[0,T]}\left(KN-\sum_{i=1}^KQ_i^{ d(N)}(t)\right)$ is $\Op(\sqrt{N})$, 
for any function $c(N)$ growing to infinity (to be chosen later), we have with probability tending to 1,
\begin{align*}
&\frac{\lambda(N)T}{\sqrt{N}}\left(1-\frac{1}{N}\sup_{t\in[0,T]}\left(KN-\sum_{i=1}^KQ_i^{ d(N)}(t)\right)\right)^{ d(N)}\\
&\hspace{2cm}\geq \frac{\lambda(N)T}{\sqrt{N}}\left(1-\frac{\sqrt{N}c(N)}{N}\right)^{ d(N)}
\geq  \frac{\lambda(N)T}{\sqrt{N}}\left(1-\frac{c(N)}{\sqrt{N}}\right)^{ d(N)}.
\end{align*}
Now since $d(N)/ \sqrt{N}\log (N)\to 0$ as $N\to\infty$, define $\omega(N):=\sqrt{N}\log (N)/d(N)$, which tends to infinity as $N$ grows large. Choose $c(N)$ such that $c(N)/\omega(N)\to 0$, as $N\to\infty.$ In that case,
\begin{align*}
\frac{\lambda(N)T}{\sqrt{N}}\left(1-\frac{c(N)}{\sqrt{N}}\right)^{ d(N)}
&= T\exp \left[\log (\sqrt{N}-\beta)+\frac{\sqrt{N}\log (N)}{\omega(N)}\log\left(1-\frac{c(N)}{\sqrt{N}}\right)\right]\\
& = T\exp \left[\log (\sqrt{N}-\beta)-\frac{\sqrt{N}\log (N)}{\omega(N)}\frac{c(N)}{\sqrt{N}}\right]\\
&\to\infty\quad\mbox{as }N\to\infty.
\end{align*}
\end{remark}

\section{Conclusion}\label{sec:conclusion-mor}
In this chapter we have investigated asymptotic optimality
properties for JSQ$(d)$ load balancing schemes in large-scale systems.
Specifically, we considered a system of $N$ parallel identical
server pools and a single dispatcher which assigns arriving tasks
to the server pool with the minimum number of tasks among $d(N)$
randomly selected server pools.
We showed that the fluid limit in a regime where the total arrival
rate and number of server pools grow large in proportion coincides
with that for the ordinary JSQ policy ($d(N) = N$) as long as
$d(N) \to \infty$ as $N \to \infty$, however slowly.
We also proved that the diffusion limit in the Halfin-Whitt regime
corresponds to that for the ordinary JSQ policy as long as
$d(N)$ grows faster than $\sqrt{N} \log(N)$, and that the latter
growth rate is in fact nearly necessary.
These results indicate that the optimality of the JSQ policy can be
preserved at the fluid-level and diffusion-level while reducing
the communication overhead by nearly a factor $O(N)$
and $O(\sqrt{N} / \log(N))$, respectively.
In future work we plan to further establish convergence rates and extend
the results to non-exponential service requirement distributions.

The proofs of the asymptotic optimality properties rely on a novel
stochastic coupling construction to bound the difference in the
system occupancy processes between the JSQ policy and a JSQ($d$)
scheme with an arbitrary value of~$d$.
It is worth observing that the coupling construction is two-dimensional in nature, and fundamentally different from the classical coupling approach used for deriving stochastic dominance properties for the ordinary JSQ policy and for establishing universality in the single-server case in Chapter~\ref{chap:univjsqd}.
As it turns out, a direct comparison between the JSQ policy
and a JSQ($d$) scheme is a significant challenge.
Hence, we adopted a two-stage approach based on a novel class
of schemes which always assign the incoming task to one of the
server pools with the $n(N) + 1$ smallest number of tasks.
Just like the JSQ($d(N)$) scheme, these schemes may be thought
of as `sloppy' versions of the JSQ policy.
Indeed, the JSQ($d(N)$) scheme is guaranteed to identify the
server pool with the minimum number of tasks, but only among
a randomly sampled subset of $d(N)$ server pools.
In contrast, the schemes in the above class only guarantee that
one of the $n(N) + 1$ server pools with the smallest number of tasks
is selected, but across the entire system of $N$ server pools.
We showed that the system occupancy processes for an intermediate
blend of these schemes are simultaneously close on a $g(N)$ scale
(e.g.~$g(N) = N$ or $g(N) = \sqrt{N}$) to both the JSQ policy
and the JSQ($d(N)$) scheme for suitably chosen values of $d(N)$
and $n(N)$ as function of $g(N)$.
Based on the latter asymptotic universality, it then sufficed to
establish the fluid and diffusion limits for the ordinary JSQ policy.


%% file: energy1.tex
\begin{abstract}
A fundamental challenge in large-scale service systems 
is to achieve highly efficient server utilization and limit energy consumption, while providing excellent user-perceived performance in the presence of uncertain and time-varying demand patterns. Auto-scaling provides a popular paradigm for automatically adjusting service capacity in response to demand while meeting performance targets, and queue-driven auto-scaling techniques have been widely investigated in the literature. In typical data center architectures and cloud environments however, no centralized queue is maintained, and load balancing algorithms immediately distribute incoming tasks among parallel queues. In these distributed settings with vast numbers of servers, centralized queue-driven auto-scaling techniques involve a substantial communication overhead and major implementation burden, or may not even be viable at all.

Motivated by the above issues, we propose a joint auto-scaling and load balancing scheme which does not require any global queue length information or explicit knowledge of system parameters, and yet provides provably near-optimal service elasticity. We establish the fluid-level dynamics for the proposed scheme in a regime where the total traffic volume and nominal service capacity grow large in proportion. The fluid-limit results show that the proposed scheme achieves asymptotic optimality in terms of user-perceived delay performance as well as energy consumption. Specifically, we prove that both the waiting time of tasks and the relative energy portion consumed by idle servers vanish in the limit. At the same time, the proposed scheme operates in a distributed fashion and involves only constant communication overhead per task. 
Extensive simulation experiments corroborate the fluid-limit results, and demonstrate that the proposed scheme can match the user performance and energy consumption of state-of-the-art approaches that do take full advantage of a centralized queue.
\end{abstract}

\section{Introduction}
\label{sec:intro}
\vspace{.25cm}

In this chapter we propose a joint auto-scaling and load balancing scheme which does not require
any global queue length information or explicit knowledge of system
parameters, and yet achieves near-optimal service elasticity.
The latter property is crucial in reducing the energy consumption in data centers in the presence of variable demand as described in Section~\ref{token}.
For convenience, we focus on a system with just a single dispatcher,
but the proposed scheme naturally extends to scenarios with
multiple dispatchers.

The proposed scheme involves a token-based feedback protocol,
allowing the dispatcher to keep track of idle-on servers in standby mode as well as
servers in idle-off mode and setup mode, as further described below.
Specifically, when a server becomes idle, it sends a message to the
dispatcher to report its status as idle-on.
Once a server has remained continuously idle for more than
an exponentially distributed amount of time with parameter $\mu>0$
(standby period), it turns off, and sends a message to the
dispatcher to change its status to idle-off.

When a task arrives, and there are idle-on servers available,
the dispatcher assigns the task to one of them at random,
and updates the status of the corresponding server to busy accordingly.
Otherwise, the task is assigned to a randomly selected busy server.
In the latter event, if there are any idle-off servers,
the dispatcher instructs one of them at random to start the setup
procedure, and updates the status of the corresponding server from idle-off to setup mode. 
It then takes an exponentially distributed amount of time with
parameter $\nu>0$ (setup period) for the server to become on,
at which point it sends a message to the dispatcher to change its
status from setup mode to idle-on.

Note that tasks are only dispatched to `on' servers (idle or busy), and in no
circumstance assigned to an `off' server (idle-off or setup mode).
Also, a server only sends a (green, say) message when a task
completion leaves its queue empty, and sends at most one (red, say)
message when it turns off after a standby period per green message,
so that at most two messages are generated per task.


In order to analyze the response time performance and energy consumption of the
proposed scheme, we consider a scenario with $N$ homogeneous servers,
and  establish the fluid-level dynamics for the proposed scheme
in a regime where the total task arrival rate and nominal number
of servers grow large in proportion. 
This regime not only offers analytical tractability, but is also highly relevant given the massive numbers of servers in data centers and cloud networks.
The fluid-limit results show that the proposed scheme achieves
asymptotic optimality in terms of response time performance
as well as energy consumption.
Specifically, we prove that for any positive values of~$\mu$
and~$\nu$ both the waiting time incurred by tasks and the relative
energy portion consumed by idle servers vanish in the limit. 
The latter results not only hold for exponential service time distributions, but also extend to a multi-class scenario with phase type service time distributions.
To the best of our knowledge, this is the first scheme to provide auto-scaling capabilities in a setting with distributed queues and achieve near-optimal service elasticity.
Extensive simulation experiments corroborate the fluid-limit results,
and demonstrate that the proposed scheme can match the user
performance and energy consumption of state-of-the-art approaches
that do assume the full benefit of a centralized queue.\\

As mentioned above, centralized queue-driven auto-scaling
mechanisms have been widely considered in the literature~\cite{ALW10, GDHS13, LCBWGWMH12, LLWLA11a, LLWLA11b, LLWA12, LWAT13, PP16, UKIN10, WLT12}.
Under Markovian assumptions, the behavior of these mechanisms can
be described in terms of various incarnations of M/M/N queues
with setup times.
A particularly interesting variant  considered  
by Gandhi \emph{et~al.}~\cite{GDHS13} is referred to as M/M/N/setup/ delayedoff.
In this mechanism, when a server $s$ finishes a service, and finds no immediate waiting task, it waits for an exponentially distributed amount of time with parameter $\mu$.
In the meantime, if a task arrives, then it is immediately assigned to server $s$ (or one of the idle-on servers at random), otherwise server $s$ is turned off.
When a task arrives, if there is no idle-on server, then it selects one of the switched off servers $s'$ say (if any), starts the setup procedure in $s'$, and waits in the queue for service. 
The setup procedure also takes an exponentially distributed amount of time with parameter $\nu$.
During the setup procedure, if some other server completes a service, then the waiting task at the head of the queue is assigned to that server, and the server $s'$ terminates its setup procedure unless there is any task $w$ waiting in the queue that had not started a setup procedure (due to unavailability of idle-off servers at its arrival epoch).
In the latter event, the server continues to be in setup mode for task $w$.
Gandhi \emph{et~al.}~\cite{GDHS13} provide an exact analysis of this model, and observe that this mechanism performs very well in a work-conserving pooled server scenario.
There are several further recent papers which examine on-demand server addition/removal in a somewhat different vein~\cite{PS16, NS16}. 
 Generalizations towards non-stationary arrivals and impatience effects have also been considered recently~\cite{PP16}.

Another related strand of research that starts from the seminal paper~\cite{YDS95} is concerned with scaling the speed of a single processor in order to achieve an optimal trade-off between energy consumption and response time performance. 
In this framework, a stream of tasks having specific deadlines arrive at a processor that either accepts the task and finishes serving it before the deadline, or discards the task at arrival.
The processor can work faster at the cost of producing more heat.
To strike the optimal balance between the revenue earned due to task completions and the energy usage, the server can scale its speed, (possibly) depending on its current load.
Dynamic versions of this speed-scaling scenario have been studied in~\cite{BPS07, B05, C72, WS87, WLT12}.

In case standby periods are infinitely long, idle servers always remain active and the proposed scheme
corresponds to the so-called Join-the-Idle-Queue (JIQ) policy,
as considered in Chapter~\ref{chap:jiq}.
Fluid-limit results described in Section~\ref{ssec:fluidjiq} show that
under Markovian assumptions, the JIQ policy achieves a zero
probability of wait for any fixed subcritical load per server
in a regime where the total number of servers grows large.
Results in Chapter~\ref{chap:jiq} indicate that the JIQ policy exhibits the
same diffusion-limit behavior as the Join-the-Shortest-Queue (JSQ)
strategy, and thus achieves optimality at the diffusion level.
These results show that the JIQ policy provides asymptotically
optimal delay performance while only involving minimal
communication overhead (at most one message per task).
However, in the JIQ policy no servers are ever deactivated,
resulting in a potentially excessive amount of energy wastage.
The scheme that we propose retains the low communication overhead
of the JIQ policy (at most two messages per task) and also
preserves the asymptotic optimality at the fluid level,
in the sense that the waiting time vanishes in the limit of $N\to\infty$.
At same time, however, any surplus idle servers are judiciously
deactivated in our scheme, ensuring that the relative energy
wastage vanishes in the limit as well. \\

The remainder of the chapter is organized as follows.
In Section~\ref{sec:model-sig} we present a detailed model description, and provide
a specification of the proposed scheme.
In Section~\ref{sec:results} we state the main results, and offer an interpretation
and discussion of their ramifications with
the full proof details relegated to Section~\ref{sec:proofs}. 
In Section~\ref{sec:phase-type} we describe how the fluid-limit results extend to phase type service time distributions.
In Section~\ref{sec:simulation} we discuss the simulation experiments that we
conducted to support the analytical results and to benchmark the
proposed scheme against state-of-the-art approaches.
We make a few brief concluding remarks and offer some suggestions
for further research in Section~\ref{sec:conclusion-sig}.

\section{Model description and algorithm specification}\label{sec:model-sig}
Consider a system of $N$~parallel queues with identical servers and a single dispatcher. 
Tasks with unit-mean exponentially distributed service requirements arrive as a Poisson process of rate $\lambda_N(s) = N\lambda(s)$ at time $s\geq 0$, where $\lambda(\cdot)$ is a bounded positive real-valued function, bounded away from zero.
In case of a fixed arrival rate, $\lambda(s)\equiv \lambda$ is assumed to be constant. 
Incoming tasks cannot be queued at the dispatcher, and must immediately and irrevocably be forwarded to one of the servers where they can be queued, possibly subject to a finite buffer capacity limit $B$. 
The service discipline at each server is oblivious to the actual service requirements (e.g., FCFS).
A turned-off server takes an Exp$(\nu)$ time (setup period) to be turned on.

We now introduce a token-based joint auto-scaling and load balancing scheme called TABS (Token-based Auto Balance Scaling).\\

\noindent
\textbf{Algorithm specification. TABS:} 
\begin{itemize}
\item When a server becomes idle, it sends a `green' message to the dispatcher, waits for an $\expn(\mu)$ time (standby period), and turns itself off by sending a `red' message to the dispatcher (the corresponding green message is destroyed).
\item When a task arrives, the dispatcher selects a green message at random if there are any, and assigns the task to the corresponding server (the corresponding green message is replaced by a `yellow' message). 
Otherwise, the task is assigned to an arbitrary busy  server  (and is lost if there is none), and if at that arrival epoch there is a red message at the dispatcher, then it selects one at random, and the setup procedure of the corresponding server is initiated, replacing its red message by an `orange' message.
\item Any server which activates due to the latter event, sends a green message to the dispatcher (the corresponding orange message is replaced), waits for an $\expn(\mu)$ time for a possible assignment of a task, and again turns itself off by sending a red message to the dispatcher.
\end{itemize}
\begin{figure}
\begin{center}
\includegraphics[scale=1]{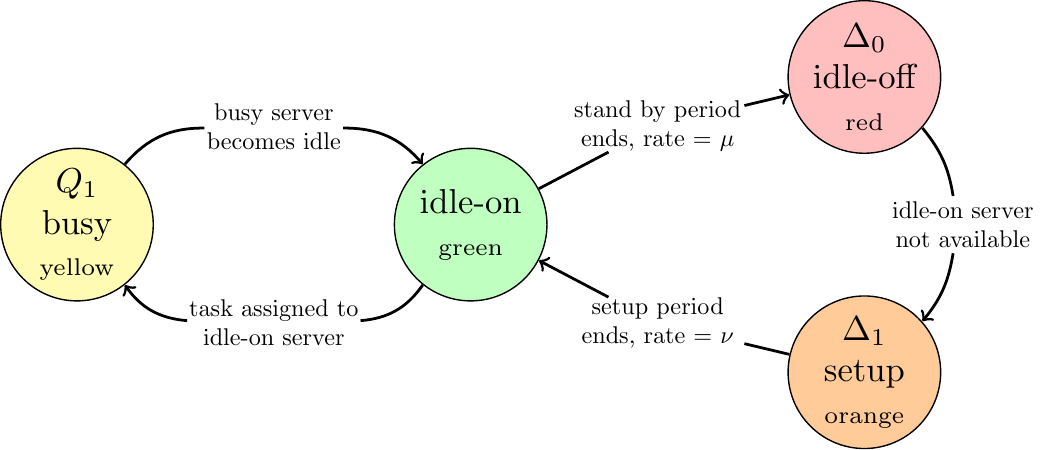}
\end{center}
\caption{Illustration of server on-off decision rules in the TABS scheme, along with message colors and state variables.}
\label{fig:scheme-sig}
\end{figure}
The TABS scheme gives rise to a distributed operation in which  servers are in one of four states (busy, idle-on, idle-off or standby), and advertize their state to the dispatcher via exchange of tokens. Figure~\ref{fig:scheme-sig} illustrates this token-based exchange protocol. 
Note that setup procedures are never aborted and continued even when idle-on servers do become available.  When setup procedures \emph{are} terminated in the latter event, the proposed scheme somewhat resembles the delayed-off scheme considered by Gandhi \emph{et~al.}~\cite{GDHS13} in terms of auto-scaling actions.  This comes however with an extra overhead penalty, without producing any improvement in response time performance or energy consumption in the large-capacity limit, as will be shown later.  

\noindent
\textbf{Notation.}
Let $\mathbf{Q}^N(t) := (Q_1^N(t), Q_2^N(t), \dots, Q_B^N(t))$ denote the system occupancy state, where $Q_i^N(t)$ is the number of servers with queue length greater than or equal to $i$ at time $t$, including the possible task in service.
Also, let $\Delta_0^N(t)$ and $\Delta_1^N(t)$ denote the number of idle-off servers and servers  in setup mode at time $t$, respectively. 
Note that the process $(\mathbf{Q}^N(t),\Delta_0^N(t),\Delta_1^N(t))_{t\geq 0}$ provides a proper state description by virtue of the exchangeablity of the servers and is Markovian.
The exact analysis of the above system becomes complicated due to the strong dependence among the queue length processes of the various servers.
Moreover, the arrival processes at individual servers are not renewal processes, which makes the problem even more challenging.
Thus we resort to an asymptotic analysis, where the task arrival rate and number of servers grow large in proportion.
In the limit the collective system then behaves like a deterministic system, which is amenable to analysis.
The fluid-scaled quantities are denoted by the respective small letters, \emph{viz.}~$q_i^{N}(t):=Q_i^{N}(t)/N$, $\delta_0^N(t) = \Delta_0^N(t)/N$, and $\delta_1^N(t) = \Delta_1^N(t)/N$. For brevity in notation, we will write $\mathbf{q}^N(t) = (q_1^N(t),\dots,q_B^N(t))$ and $\bld{\delta}^N(t) = (\delta_0^N(t),\delta_1^N(t))$. 
Let
$$
E = \Big\{(\bld{q},\dd)\in [0,1]^{B+2}:  q_i\geq q_{i+1},\ \forall i, \  \delta_0+\delta_1+ q_1\leq 1 \Big\},
$$
denote the space of all fluid-scaled occupancy states,
so that $(\mathbf{q}^N(t),\bld{\delta}^N(t))\in E$ for all $t$.
Endow $E$ with the product topology, and the Borel $\sigma$-algebra $\mathcal{E}$, generated by the open sets of $E$.
For stochastic boundedness of a process we refer to \cite[Definition 5.4]{PTRW07}. 
For any complete separable metric space $E$, denote by $D_E[0,\infty)$, the set of all $E$-valued \emph{c\`adl\`ag} (right continuous with left limits exist) processes.
By the symbol `$\dto$' we denote weak convergence for real-valued random variables, and convergence with respect to Skorohod-$J_1$ topology for c\`adl\`ag processes.

\section{Overview of results}
\label{sec:results}
In this section we provide an overview of the main results and discuss their ramifications.  For notational transparency, we focus on the case of exponential service time distributions.  In Section~\ref{sec:phase-type} we show how some of the results extend to phase type service time distributions, at the expense of more complex notation.
\begin{theorem}[{Fluid limit for exponential service time distributions}]
\label{th: fluid}
Assume that $(\qq^N(0),\bld{\delta}^N(0))$ converges to $(\qq^\infty,\bld{\delta}^\infty)\in E$, as $N\to\infty$, where $q_1^\infty>0$. Then the process $\{(\qq^N(t),\dd^N(t))\}_{t\geq 0}$ converges weakly to the deterministic process $\{(\qq(t),\dd(t))\}_{t\geq 0}$ as $N\to\infty$, which satisfies the following integral equations:
\begin{align*}
 q_i(t)&=q_i^\infty +\int_0^t\lambda(s) p_{i-1}(\qq(s),\dd(s),\lambda(s))\dif s
 - \int_0^t(q_i(s)-q_{i+1}(s))\dif s,
 \end{align*}
 for $i= 1,\ldots,B,$ and
 \begin{align*}
\delta_0(t)&=\delta_0^\infty+\mu\int_0^t u(s)\dif s-\xi(t),  \qquad
\delta_1(t)=\delta_1^\infty+\xi(t)-\nu\int_0^t\delta_1(s)\dif s,\nonumber
\end{align*}
where by convention $q_{B+1}(\cdot) \equiv 0$,
and
\begin{align*}
u(t) &= 1- q_1(t) - \delta_0(t) - \delta_1(t),\\
\xi(t) &= \int_0^t\lambda(s)(1-p_0(\qq(s),\dd(s),\lambda(s)))\ind{\delta_0(s)>0}\dif s.
\end{align*}
For any $(\qq,\dd)\in E$, $\lambda>0$, $(p_i(\qq,\dd,\lambda))_{i\geq 0}$ are given by 
\begin{align*}
 p_0(\qq,\dd,\lambda) &= 
\begin{cases}
&1\qquad \text{if}\qquad u=1-q_1-\delta_0-\delta_1>0,\\
&\min \{\lambda^{-1}(\delta_1\nu + q_1-q_2), 1\},\quad\text{otherwise,}
\end{cases}\\
\quad p_i(\qq,\dd,\lambda)  &= (1-p_0(\qq,\dd,\lambda)) (q_{i}-q_{i+1})q_1^{-1},\  i =1,\ldots,B.
\end{align*}
\end{theorem}

We now provide an intuitive explanation of the fluid limit stated above.
The term $u(t)$ corresponds to the asymptotic fraction of idle-on servers in the system at time $t$, and $\xi(t)$ represents the asymptotic cumulative number of server setups (scaled by $N$) that have been initiated during $[0,t]$.
The coefficient $p_i(\qq,\dd,\lambda)$ can be interpreted as the instantaneous fraction of incoming tasks that are assigned to some server with queue length $i$, when the fluid-scaled occupancy state is $(\qq,\dd)$ and the scaled instantaneous arrival rate is $\lambda$.
Observe that as long as $u>0$, there are idle-on servers, and hence all the arriving tasks
will join idle servers. 
This explains that if $u>0$, $p_0(\qq,\dd,\lambda) = 1$ and $p_i(\qq,\dd,\lambda)=0$ for $i=1,\ldots,B-1$.
If $u=0$, then observe that 
servers become idle at rate $q_1-q_2$, and servers in setup mode turn on at rate $\delta_1\nu$.
Thus the  idle-on servers are created at a total rate $\delta_1\nu + q_1-q_2$.
If this rate is larger than the arrival rate $\lambda$, then almost all the arriving tasks can be assigned to idle servers.
Otherwise, only a fraction $(\delta_1\nu + q_1-q_2)/\lambda$
of arriving tasks join idle servers. 
The rest of the tasks are distributed uniformly among busy servers, so a proportion $(q_{i}-q_{i+1})q_1^{-1}$ are assigned to servers having queue length~$i$.
For any $i=1,\ldots,B$, $q_i$ increases when there is an arrival to some server with queue length $i-1$, which occurs at rate $\lambda p_{i-1}(\qq,\dd,\lambda)$, and it decreases when there is a departure from some server with  queue length~$i$, which occurs at rate $q_i-q_{i-1}$. 
Since each idle-on server turns off at rate $\mu$, the fraction of servers in the off mode increases at rate 
$\mu u$.
Observe that if $\delta_0>0$, for each task that cannot be assigned to an idle server, a setup procedure is initiated  at one idle-off server. 
As noted above, $\xi(t)$ captures the (scaled) cumulative number of setup procedures initiated up to time~$t$.
Therefore the fraction of idle-off servers and the fraction of servers in setup mode decreases and increases by $\xi(t)$, respectively, during $[0,t]$.
Finally, since each server in setup mode becomes idle-on at rate $\nu$, the fraction of servers in setup mode decreases at rate $\nu\delta_1$.\\

\noindent
\textbf{Fixed point.}
In case of a constant arrival rate $\lambda(t)\equiv\lambda<1$, the fluid limit in Theorem~\ref{th: fluid} has a unique fixed point:
\begin{equation}\label{eq:fixed point-sig}
\delta_0^*=1-\lambda,\qquad\delta_1^*=0,\qquad q_1^*=\lambda\quad\mbox{and}\quad q_i^*=0,
\end{equation}
for $i=2,\ldots,B.$ 
Indeed, it can be verified that $p_0(\qq^*,\dd^*,\lambda)= 1$ and $u^*=0$ for $(\qq^*,\dd^*)$ given by~\eqref{eq:fixed point-sig} so that the derivatives of $q_i$, $i = 1,\dots,B$, $\delta_0$, and $\delta_1$ become zero, and that these cannot be zero at any other point in $E$.
Note that, at the fixed point, a fraction $\lambda$ of the servers have exactly one task while the remaining fraction have zero tasks, independently of the values of the parameters $\mu$ and $\nu$.

The next proposition states the global stability of the fluid limit, i.e., starting from any point in $E$, the dynamical system  defined by the system of integral equations in Theorem~\ref{th: fluid} converges to the fixed point~\eqref{eq:fixed point-sig} as $t\to\infty$.
\begin{proposition}[{Global stability of the fluid limit}]
\label{prop:glob-stab}
Assume that $(\qq(0), \dd(0)) = (\qq^\infty,\dd^\infty)\in E$. Then  
$$(\qq(t),\dd(t))\to (\qq^*,\dd^*),\quad \mbox{as}\quad t\to\infty,$$
where $(\qq^*,\dd^*)$ is as defined in~\eqref{eq:fixed point-sig}.
\end{proposition}
There are general methods to prove global stability if the evolution of the dynamical system  satisfies some kind of monotonicity property induced by the drift structure~\cite{TX11, Mitzenmacher01}.
Here, it is not straightforward to establish such a monotonicity property, and harder to find a suitable Lyapunov function.
Instead we exploit specific properties of the fluid limit in order to prove the global stability.
Observe that the global stability in particular also establishes the uniqueness of the fixed point above. 
The proof of Proposition~\ref{prop:glob-stab} is presented in Subsection~\ref{app:globstab}.

The global stability can be leveraged to show that the steady-state distribution of the $N^{\mathrm{th}}$ system, for large $N$, can be well approximated by the fixed point of the fluid limit in~\eqref{eq:fixed point-sig}. 
Specifically, in the next proposition, whose proof is provided in Subsection~\ref{app:globstab}, we demonstrate the convergence of the steady-state distributions, and hence the interchange of the large-capacity
($N\to\infty$) and steady-state ($t\to\infty$) limits.
Since the buffer capacity $B$ at each server is finite, for each $N$, the Markov process $(\QQ^N(t), \Delta_0^N(t),\Delta_1^N(t))$ is irreducible, has a finite state space, and thus has a unique steady-state distribution.
Let $\pi^N$ denote the steady-state distribution of the $N^{\mathrm{th}}$ system, i.e.,
$$\pi^{N}(\cdot)=\lim_{t\to\infty}\mathbb{P}\ \big(\qq^{N}(t)=\cdot, \dd^N(t) = \cdot\big).$$ 
\begin{proposition}[{Interchange of limits}]
\label{thm:limit interchange}
As $N\to\infty$,
$\pi^N\dto\pi$, where $\pi$ is given by the Dirac mass concentrated upon $(\qq^*,\dd^*)$ defined in~\eqref{eq:fixed point-sig}.
\end{proposition}

\noindent
\textbf{Performance metrics.}
As mentioned earlier, two key performance metrics are the expected waiting time of tasks $\expt[W^N]$ and energy consumption $\expt[ P^N]$ for the $N^{\mathrm{th}}$ system in steady state.
In order to quantify the energy consumption, we assume that the energy usage of a server is $P_{\full}$ when busy or in set-up mode, $P_{\idle}$ when idle-on, and zero when turned off.
Evidently, for any value of $N$, at least a fraction $\lambda$ of the servers must be busy in order for the system to be stable, and hence $\lambda P_{\full}$ is the minimum mean energy usage per server needed for stability.
We will define $\expt[Z^N]=\expt[P^N]-\lambda P_{\full}$ as the relative energy wastage accordingly.
The next proposition demonstrates that asymptotically the expected waiting time and  energy consumption for the TABS scheme vanish in the limit, for any strictly positive values of $\mu$ and $\nu$.
The key implication is that
the TABS scheme, while only involving constant communication overhead per task, provides performance in a distributed setting that is as good at the fluid level as can possibly be achieved, even in a centralized queue, or with unlimited information exchange.
\begin{proposition}[Asymptotic optimality of TABS scheme]
\label{prop:perform}
In a fixed arrival rate scenario $\lambda(t)\equiv \lambda<1$, for any $\mu>0$, $\nu>0$, as $N\to\infty$,
\begin{enumerate}[{\normalfont (a)}]
\item\ $[$zero mean waiting time$]$ $\expt[W^N]\to 0$,
\item\ $[$zero energy wastage$]$ $\expt[Z^N]\to 0$.
\end{enumerate}
\end{proposition}

\begin{proof}[Proof of Proposition~\ref{prop:perform}]
By Little's law, the mean stationary waiting time $\expt[W^N]$ in the $N^{\mathrm{th}}$ system may be expressed as $(N \lambda)^{-1} \expt[L^N]$, where $L^N = \sum_{i = 2}^{B} Q_i^N$ represents a random variable with the stationary distribution of the total number of waiting tasks in the $N^{\mathrm{th}}$ system.
Thus, $\expt[W^N] = \lambda^{-1} \sum_{i = 2}^{B} \expt[q_i^N]$, where $\qq^N$ is a random vector with the stationary distribution of $\qq^N(t)$ as $t \to \infty$.
Invoking Proposition~\ref{thm:limit interchange} and the fixed point as identified in~\eqref{eq:fixed point-sig}, we obtain that $\expt[W^N] \to \sum_{i = 2}^{B} q_i^* = 0$ as $N \to \infty$.

Denoting by $U^N= N-Q_1^N-\Delta_0^N-\Delta_1^N$ the number of idle-on servers,
the stationary mean energy consumption per server in the $N^{\mathrm{th}}$ system may  be expressed as 
$$\frac{1}{N} \expt[(Q_1^N + \Delta_1^N) P_{\full} + U^N P_{\idle}] = \expt[(q_1^N + \delta_1^N) P_{\full} + u^N P_{\idle}].$$
Applying Proposition~\ref{thm:limit interchange} and the fixed point as identified in~\eqref{eq:fixed point-sig}, we deduce that $\expt[P^N] \to (q_1^* + \delta_1^*) P_{\full} + u^* P_{\idle} = (1 - \delta_0^*) P_{\full} - u^* (P_{\full} - P_{\idle}) = \lambda P_{\full}$ as $N \to \infty$.
This yields that $\expt[Z^N] = \expt[P^N]-\lambda P_{\full}$ converges to 0.
\end{proof}

 The quantitative values of the energy usage and waiting time for finite values of $N$ will be evaluated through extensive simulations in Section~\ref{sec:simulation}.\\

\noindent
\textbf{Comparison to ordinary JIQ policy.}
Consider the fixed arrival rate scenario $\lambda(t) \equiv \lambda$. 
It is worthwhile to observe that the component $\qq$ of the fluid limit in Theorem~\ref{th: fluid} coincides with that for the ordinary JIQ policy where servers always remain on, when the system following the TABS scheme starts with all the servers being idle-on, and $\lambda+\mu<1$. 
To see this, observe that the component $\qq$ depends on $\dd$ only through $(p_{i-1}(\qq,\dd))_{i\geq 1}$. 
Now, $p_0 =1$, $p_i = 0$, for all $i\geq 1$, whenever $q_1+\delta_0+\delta_1<1$, irrespective of the precise values of $(\qq,\dd)$. 
Moreover, starting from the above initial state, $\delta_1$ can increase only when $q_1+\delta_0=1$. 
Therefore, the fluid limit of $\qq$ in Theorem~\ref{th: fluid} and the ordinary JIQ scheme are identical if the system parameters $(\lambda,\mu,\nu)$ are such that $q_1(t)+\delta_0(t) < 1$, for all $t\geq 0$. Let $y(t)  = 1-q_1(t) - \delta_0(t)$. The solutions to the differential equations
\begin{equation*}
 \frac{\dif q_1(t)}{\dif t} = \lambda - q_1(t), \quad  \frac{\dif y(t)}{\dif t} = q_1(t) - \lambda -\mu y(t),
\end{equation*}$y(0)=1$, $q_1(0) = 0$ are given by 
\begin{equation*}
 q_1(t)= \lambda(1-\e^{-t}),\quad y(t) = \frac{\e^{-(1+\mu)t}}{\mu-1}\big(\e^t(\lambda+\mu-1)-\lambda \e^{\mu t}\big).
\end{equation*}
Notice that if $\lambda+\mu<1$, then $y(t) > 0$ for all $t\geq 0$ and thus, $q_1(t)+\delta_0(t) < 1$, for all $t\geq 0$.
The fluid-level optimality of the JIQ scheme was shown in~\cite{Stolyar15,Stolyar17}. 
This observation thus establishes the optimality of the fluid-limit trajectory under the TABS scheme for suitable parameter values in terms of response time performance.
From the energy usage perspective, under the ordinary JIQ policy, since 
the asymptotic steady-state fraction of busy servers ($q_1^*$) and idle-on servers are given by $\lambda$ and $1-\lambda$, respectively, the asymptotic steady-state (scaled) energy usage is given by 
\begin{align*}
\expt[P^{\mathrm{JIQ}}] 
 = \lambda P_{\full} + (1-\lambda) P_{\idle} 
= \lambda P_{\full}(1+ (\lambda^{-1}-1)f),
\end{align*}
where $f = P_{\idle}/P_{\full}$ is the relative energy consumption of an idle server.
Proposition~\ref{prop:perform} implies that the asymptotic steady-state (scaled) energy usage under the TABS scheme is $\lambda P_{\full}.$
Thus the TABS scheme reduces the asymptotic steady-state energy usage by $\lambda P_{\full}(\lambda^{-1}-1)f = (1-\lambda)P_{\idle},$ which amounts to a relative saving of $(\lambda^{-1}-1)f/(1+(\lambda^{-1}-1)f).$ 
In summary, the TABS scheme performs as good as the ordinary JIQ policy in terms of the waiting time and communication overhead while providing a significant energy saving.

\section{Extension to phase type service time distributions}
\label{sec:phase-type}

In this section we extend the fluid-limit results to phase type service time distributions. 
Specifically, the service time of each task is described by a time-homogeneous,  continuous-time Markov process with a finite state space $\{0,1,\dots,K\}$, initial distribution $\bld{r} = (r_i:0\leq i\leq K)$, transition probability matrix  $R = (r_{i,j})$, and the mean sojourn time in state $i$ being $\gamma_i^{-1}$. 
State 0 is an absorbing state, and thus  represents a service completion, while state $j$ is referred to as a type-$j$ service, and is assumed to be transient. 
For convenience, and without loss of generality, it is assumed that $r_{i,i}=0$ for all $i$, and that any incoming task  has a non-zero service time ($r_0~=~0$).  
Consider a time-homogeneous discrete-time Markov chain with the state  space $\{0,1,\dots,K\}$, and transition probability matrix $P=(p_{i,j})$, where $p_{i,j} = r_{i,j}$ for $i\geq 1$, $p_{0,j} = r_j$ $j\geq 1$, and $p_{0,0}=0$. Let $\bld{\eta}= (\eta_0,\dots,\eta_K)$ be the stationary distribution, i.e., $\bld{\eta}$ satisfies
\begin{equation}
 \eta_0r_i + \sum_{j=1}^Kr_{j,i}\eta_j = \eta_i, \quad i\geq 1,\quad \sum_{i=0}^K \eta_i = 1.
\end{equation}
The mean of the phase type service time distribution~\cite{PR00} is $\big(\sum_{i=1}^K\eta_i/\gamma_i\eta_0\big)^{-1}$, and is assumed to be one.

We assume now that the service discipline at each server is not only oblivious of the actual service requirements, but also non-preemptive, and allows at most one task to be served at any given time.
Let $Q_{i,j}^{N}(t)$ denote the number of servers with queue length at least $i$ and providing a type-$j$ service at time~$t$. Thus, $Q_i^N(t) = \sum_{j=1}^KQ_{i,j}^N(t)$. 
Denote the fluid-scaled quantities by $q_{i,j}^N(t) = Q_{i,j}^{N}(t)/N$ and the vector $\qq^N(t) = (q_{i,j}^N(t):1\leq i\leq B , 1\leq j\leq K)$. 
Let $\delta_0^N(t)$ and $\delta^N_1(t)$ be as defined before.
Let 
\begin{align*}
\hat{E} &= \bigg\{\big((q_{i,j})_{1\leq i\leq B, 1\leq j\leq K}, (\delta_0,\delta_1)\big): q_{1,j},\delta_0,\delta_1\in [0,1],
 q_{i+1,j}\leq q_{i,j},\ \forall i,j, \\
 &\hspace{8cm} \delta_0 + \delta_1 + \sum_{j=1}^K q_{1,j}\leq 1
\bigg\}
\end{align*}
denote the space of all fluid-scaled occupancy states,
so that $(\mathbf{q}^N(t),\bld{\delta}^N(t))\in \hat E$ for all $t$, and as before,
endow $\hat{E}$ with the product topology, and the Borel $\sigma$-algebra $\hat{\mathcal{E}}$, generated by the open sets of $\hat E$.
\begin{theorem}[Fluid limit for phase type service time distributions]
\label{th: fluid gen service}
Let $(\qq^N(0),\bld{\delta}^N(0))$ converge to $(\qq^\infty,\bld{\delta}^\infty)\in \hat{E}$, as $N\to\infty$, where $\sum_{j=1}^Kq_{1,j}^{\infty}>0$. Then the sequence of processes $\{\qq^N(t),\dd^N(t)\}_{t\geq 0}$ converges weakly to the  deterministic process $\{\qq(t),\dd(t)\}_{t\geq 0}$, as $N\to\infty$, which satisfies the following integral equations: for $i = 1,\ldots,B$ and $j = 1,\ldots,K$,
\begin{align*}
q_{i,j}(t)&=q_{i,j}^\infty+\int_0^t\lambda(t) p_{i-1,j}(\qq(s),\dd(s),\lambda(s))\dif s\\
&\hspace{2cm} + \int_0^t\sum_{k=1}^K (q_{i,k}(s)-q_{i+1,k}(s))\gamma_kr_{k,j}\dif s   - \gamma_j\int_0^t q_{i,j}(s)\dif s\\
&\hspace{4cm}+ \int_0^t\sum_{k=1}^K (q_{i+1, k}(s) - q_{i+2, k}(s))\gamma_kr_{k,0}r_j\dif s,  \\
\delta_0(t)&=\delta_0^\infty + \mu\int_0^tu(s)\dif s
 -\xi(t), \qquad
\delta_1(t) =\delta_1^\infty + \xi(t)-\nu\int_0^t\delta_1(s)\dif s,
\end{align*}
where by convention $q_{B+1,j}(\cdot) \equiv 0$, $j= 1,\ldots,K$,
and
\begin{align*}
 u(t) &= 1-\sum_{j=1}^Kq_{1,j}(t)-\delta_0(t)-\delta_1(t),\\ 
 \xi(t) &= \int_0^t\lambda(s)\bigg(1-\sum_{j=1}^Kp_{0,j}(\qq(s),\dd(s),\lambda(s))\bigg)\ind{\delta_0(s)>0}\dif s.\nonumber
\end{align*}
For any $(\qq,\dd)\in\hat{E}$, $\lambda>0$,
$p_{0,j}(\qq,\dd,\lambda) = r_j$ if $u=1- \sum_{j=1}^Kq_{1,j}-\delta_0-\delta_1>0$, $j= 1,\ldots,K$, and otherwise
$$p_{0,j}(\qq,\dd,\lambda) =r_j\min \bigg\{\lambda^{-1}\Big(\delta_1\nu + \sum_{j=1}^K(q_{1,j}-q_{2,j})\gamma_jr_{j,0}\Big), 1\bigg\},$$
and for $i = 1,\dots,B$,
$$p_{i,j} (\qq,\dd,\lambda)= \bigg(1-\sum_{j=1}^Kp_{0,j}(\qq,\dd,\lambda)\bigg)\frac{q_{i-1,j}-q_{i,j}}{\sum_{j=1}^Kq_{1,j}}.$$
\end{theorem}
Let us provide a heuristic justification of the fluid limit stated above. 
As in Theorem~\ref{th: fluid}, 
$u(t)$ corresponds to the asymptotic fraction of idle-on servers in the system at time $t$, 
$\xi(t)$ represents the asymptotic cumulative number of server setups (scaled by $N$) that have been initiated during $[0,t]$. The coefficient $p_{i,j}(\qq(t),\dd(t),\lambda(t))$ can be interpreted as the instantaneous fraction of incoming tasks 
that are assigned to a server with queue length $i\geq 1$ and currently providing a type-$j$ service, 
while $p_{0,j}$ specifies the fraction of incoming tasks assigned to idle servers starting with a type-$j$ service. 
The heuristic justification for the $p_{i,j}$ values builds on the same line of reasoning as for Theorem~\ref{th: fluid}. 
As long as there are idle-on servers, i.e., 
$u>0$, incoming tasks are \emph{immediately} assigned to one of those servers, and the initial service type is chosen according to the distribution $\mathbf{r}$. 
Notice that the busy servers and the servers in setup become idle at total rate $\delta_1\nu + \sum_{j=1}^K(q_{1,j}-q_{2,j})\gamma_jr_{j,0}$. 
For the case when $u=0$, we need to distinguish between two cases, depending on whether 
$\delta_1\nu + \sum_{j=1}^K(q_{1,j}-q_{2,j})\gamma_jr_{j,0} > \lambda$ or not.
In the first case, the incoming tasks are again assigned to idle-on servers immediately.  
However, if $\delta_1\nu + \sum_{j=1}^K(q_{1,j}-q_{2,j})\gamma_jr_{j,0} \leq  \lambda$, then only a fraction $\lambda^{-1}(\delta_1\nu + \sum_{j=1}^K(q_{1,j}-q_{2,j})\gamma_jr_{j,0}$ of the incoming tasks are immediately taken into service. 
In both of the above two subcases, the service types of the incoming tasks follow the distribution $\mathbf{r}$. 
This explains the expression for the $p_{0,j}$ values. 
Also, given that an incoming task does not find an  idle-on server, it is assigned to a server that has queue length $i$ and is currently providing a type-$j$ service with probability $\big(\sum_{j=1}^Kq_{1,j}\big)^{-1}(q_{i-1,j}-q_{i,j})$.
This explains the expression for $p_{i,j}$ for $i\geq 1$.
Now, notice that the expressions for $\delta_0$, and $\delta_1$ remain essentially the same as in Theorem~\ref{th: fluid} due to the fact that the dynamics of $\delta_0$ and $\delta_1$ 
depend on $q_{i,j}$'s only through the fraction of incoming tasks that join an idle-on server, which is determined by the coefficients $p_{0,j}(\qq,\bld{\delta},\lambda)$.
Finally,  $q_{i,j}$ decreases if and only if there is a completion of type-$j$ service at a server with queue length at least $i$. 
Here, we have used the fact $r_{i,i} = 0$.
Now, $q_{i,j}$ can increase due to three events: 
(i) assignment of an arriving task, which occurs at rate $\lambda p_{i-1,j}(\qq,\dd,\lambda)$, 
(ii) service completion of some other type, which now requires service of type $j$, and this occurs at rate $\sum_{k}(q_{i,k}-q_{i+1, k})\gamma_k r_{k,j}$, 
(iii) service completion occurs at some server, the task exits from the system, and the next task at that server starts with a type-$j$ service. 
This occurs at rate $\sum_{k = 1}^K (q_{i+1,k} - q_{i+2,k})\gamma_kr_{k,0}r_j$. \\

\noindent
\textbf{Fixed point of the fluid limit.} 
In case of a constant arrival rate $\lambda(t)\equiv\lambda<1$, the unique fixed point of the fluid limit in Theorem~\ref{th: fluid gen service} is given by 
\begin{equation}\label{fix-point-2}
 \delta_0^* = 1-\lambda, \quad \delta_1^* = 0,\quad q_{1,j}^*= \frac{\eta_j}{\eta_0\gamma_j}\lambda,\quad  j=1,\dots,K,
\end{equation}
and $q_{i,j}^*=0$ for all $i=2,\ldots,B$.
Indeed, it can  be verified that the derivatives of $q_{i,j}$, $i=1,\ldots,B$, $j=1,\ldots,K$, $\delta_0$, and $\delta_1$ are zero at $(\qq^*,\dd^*)$ given by~\eqref{fix-point-2}, and that these cannot be zero at any other point in $\hat{E}$. Thus, the fixed point is unique as before.
Notice that in this case also at the fixed point a fraction $\lambda$ of the servers have exactly one task while the remaining fraction have zero tasks,
independent of the values of the parameters $\mu$ and $\nu$, revealing the insensitivity of the asymptotic fluid-scaled steady-state occupancy states to the duration of the standby periods and setup periods.
Further, note that $\sum_{j=1}^K q_{1,j}^* = \lambda$ from the fact that the mean service time is one, irrespective of the initial distribution $\mathbf{r}$, transition probability matrix $R$, and parameters $\gamma_j$. 
Thus the values of $q_1^*,\ldots, q_B^*$ in the fixed point are insensitive in a distributional sense with respect to the service times.
They only depend on the service time distribution through its mean, and higher-order characteristics like variance have no impact on the steady-state performance in the large capacity limit whatsoever.

\section{Simulation experiments}\label{sec:simulation}
In this section we present extensive simulation results to illustrate the fluid-limit results, and to examine the performance of the proposed TABS scheme in terms of mean waiting time and energy consumption, and compare that with existing strategies.\\

\begin{figure}
\begin{center}
$
\begin{array}{c}
\includegraphics[width=90mm]{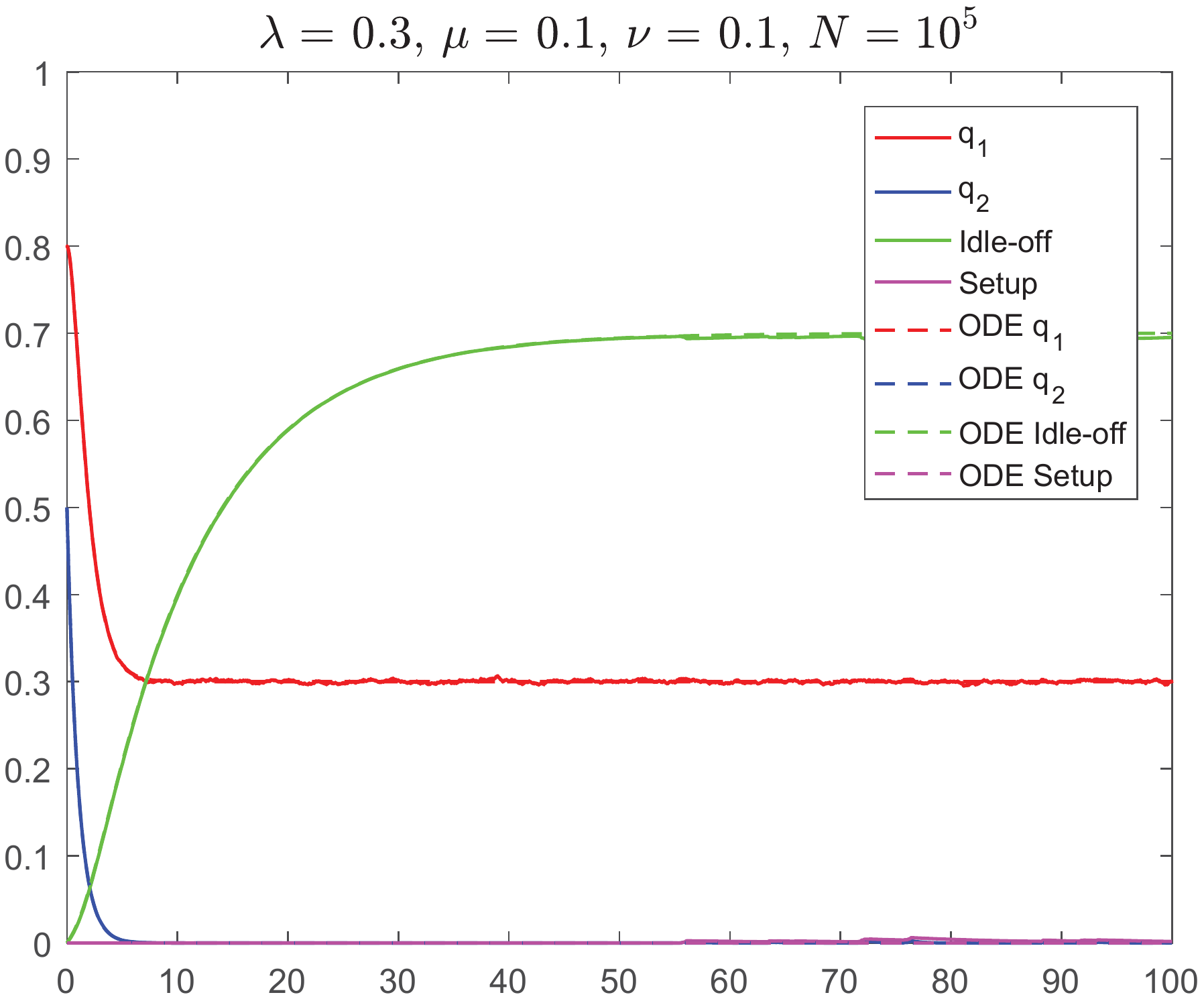}\\
\includegraphics[width=90mm]{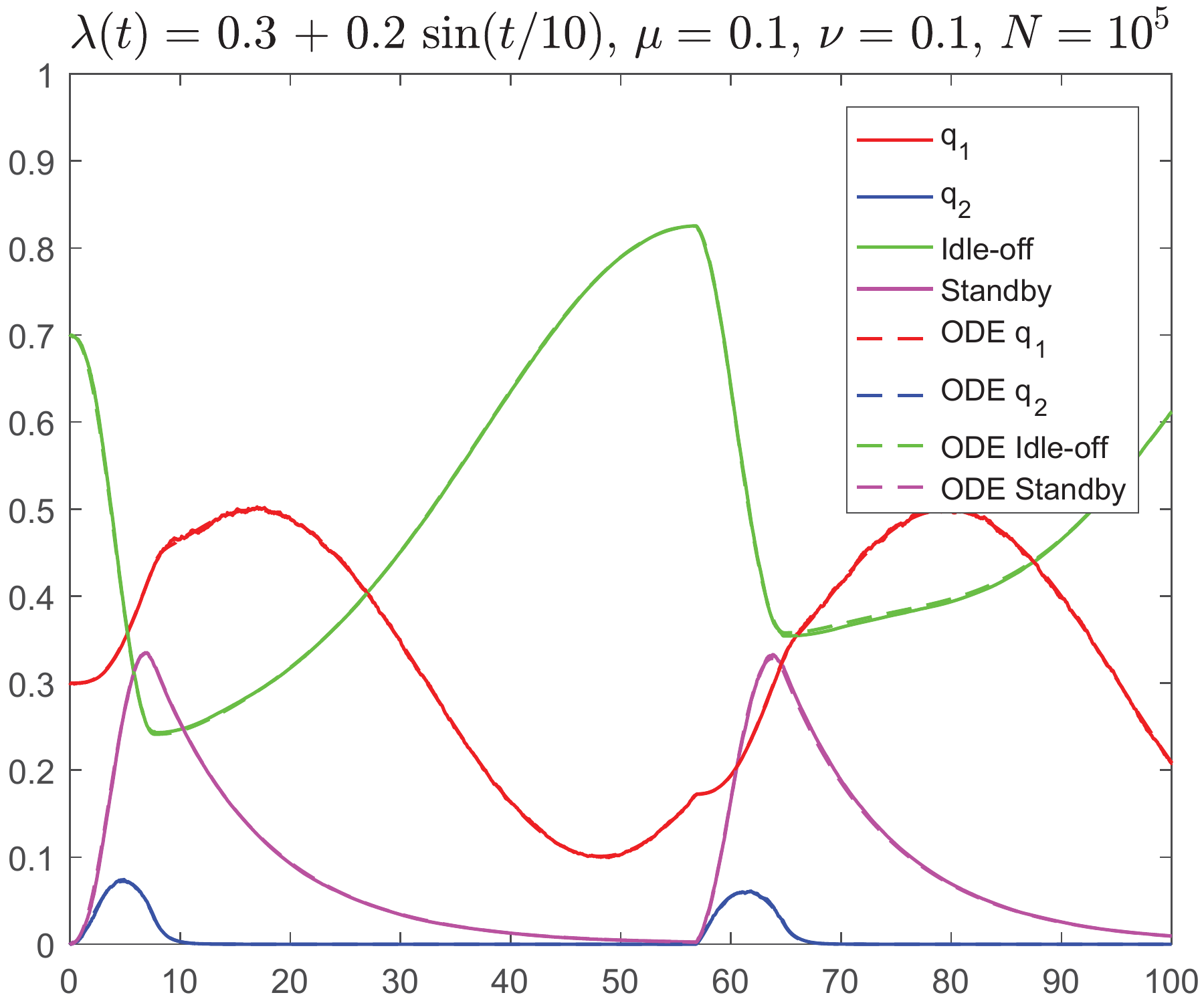}
\end{array}$
\end{center}
\caption{Illustration of the fluid-limit trajectories for $N = 10^5$ servers. The top figure is for constant arrival rate $\lambda(t) \equiv 0.3$, and the bottom figure considers a periodic arrival rate given by $\lambda(t) = 0.3+0.2\sin(t/10)$.}
\label{fig:fluid}
\end{figure}

\begin{figure}
  \begin{center}
    \includegraphics[width=90mm]{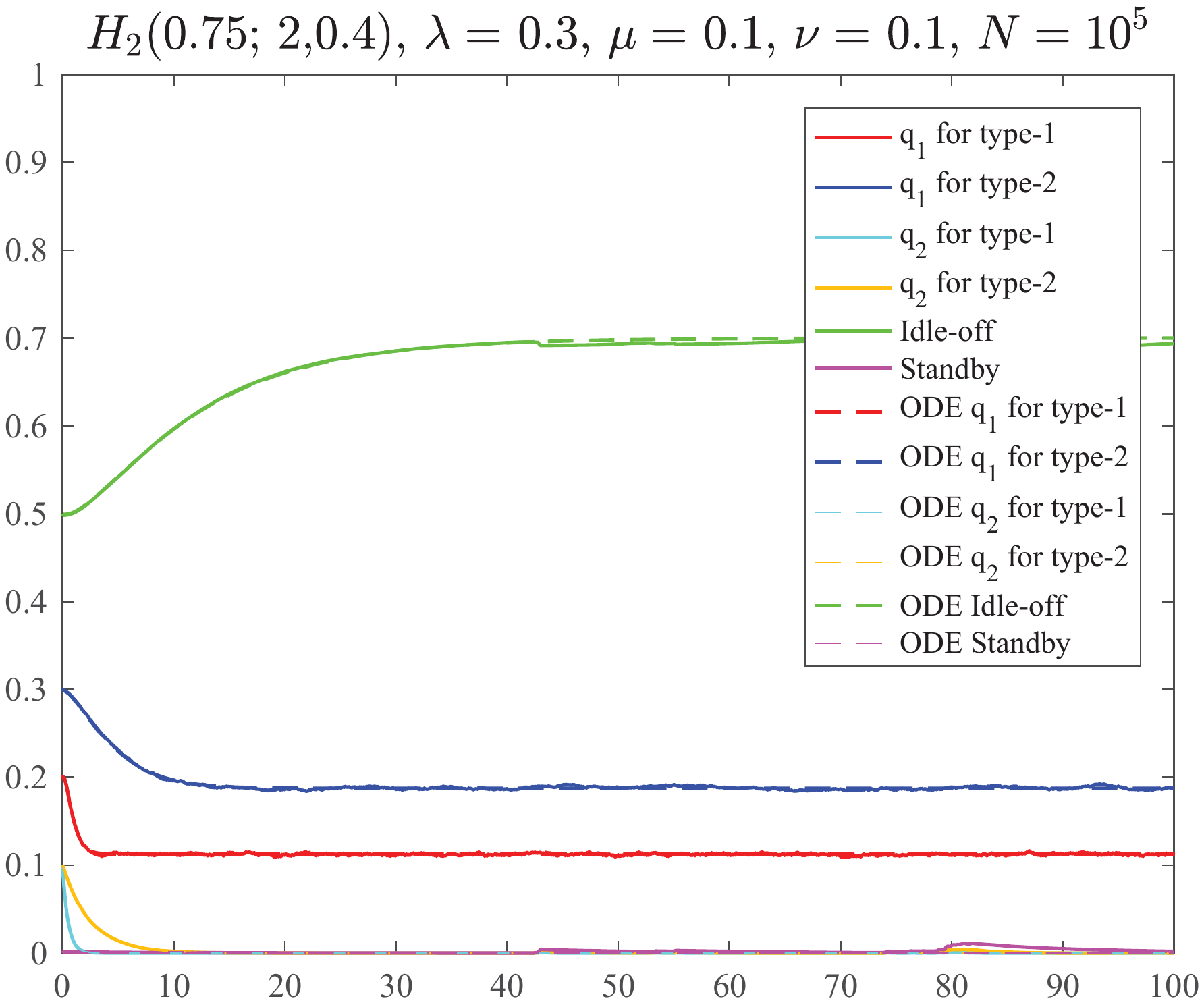}
  \end{center}
  \caption{The figure considers a hyper-exponential service time distribution. An incoming task demands  either type-1 or type-2 service with probabilities $0.75$ and $0.25$, respectively. 
The durations of type-1 and type-2 services are exponentially distributed with parameters 2 and 0.4, respectively, and thus the mean service time is 1.}
\label{fig:hyper}
\end{figure}

\noindent
\textbf{Convergence of sample paths to fluid-limit trajectories.}
The fluid-limit trajectories for the TABS scheme in Theorems~\ref{th: fluid} and~\ref{th: fluid gen service} are illustrated in Figures~\ref{fig:fluid} and~\ref{fig:hyper} for $N=10^5$ servers and three scenarios (constant arrival rate, periodic arrival rate and hyper-exponential service time distribution). 
In all three scenarios the mean standby periods are $\mu^{-1}=10$ and the mean setup periods are $\nu^{-1}=10$. 
In all cases,  the fluid-limit paths and the sample paths obtained from simulation are nearly indistinguishable.
Notice that in case of a time-varying arrival rate the period of fluctuation is only $20\pi\approx 63$ times as long as the mean service time, which is far shorter than what is usually observed in practice. Typically, service times are of sub-second order and variations in the arrival rate occur only over time scales of tens of minutes, if not several hours.
Even in such a challenging scenario, however, the fractions of idle-on servers and those with waiting tasks are negligible. 
In case of the hyper-exponential service time distribution, we note from Figure~\ref{fig:hyper} that the long-term values of $q_1 = q_{1,1} + q_{1,2}, q_{2} = q_{2,1} + q_{2,2}$, $\delta_0$ and $\delta_1$ agree with the corresponding quantities in the top chart for exponential service times.  
This reflects the asymptotic insensitivity in a distributional sense mentioned at the end of Section~\ref{sec:phase-type}, and in particular supports the observation that the proposed TABS scheme achieves asymptotically optimal response time performance and energy consumption for phase type service time distributions as well.\\

\noindent
\textbf{Convergence of steady-state performance metrics to fluid-limit values.}
In order to quantify the energy usage, we will adopt the parameter values from empirical measurements reported in~\cite{GHK12, BH07, GDHS13}. 
A server that is busy or in setup mode, consumes $P_{\full} = 200$ watts,  an idle-on server consumes $P_{\idle} = 140$ watts, and an idle-off servers consumes no energy.
We will consider the normalized energy consumption. 
Thus, the asymptotic steady-state expected normalized energy consumption $\mathbb{E}[P/340]$ is given by $10/17(q_1+\delta_1)+7/17u = 10/17(1-\delta_0)-3/17u$.
Note that the optimal energy usage (with no wastage, i.e., $\delta_0 = 1-\lambda =0.7,$ $\delta_1 = 0$, $q_1=\lambda =0.3$) is given by $3/17$.
Also recall that the asymptotic expected steady-state waiting time is given by 
$\expt[W] = \lambda^{-1} \sum_{i = 2}^{B} q_i$.\\

\begin{figure}
\begin{center}
\includegraphics[width=90mm]{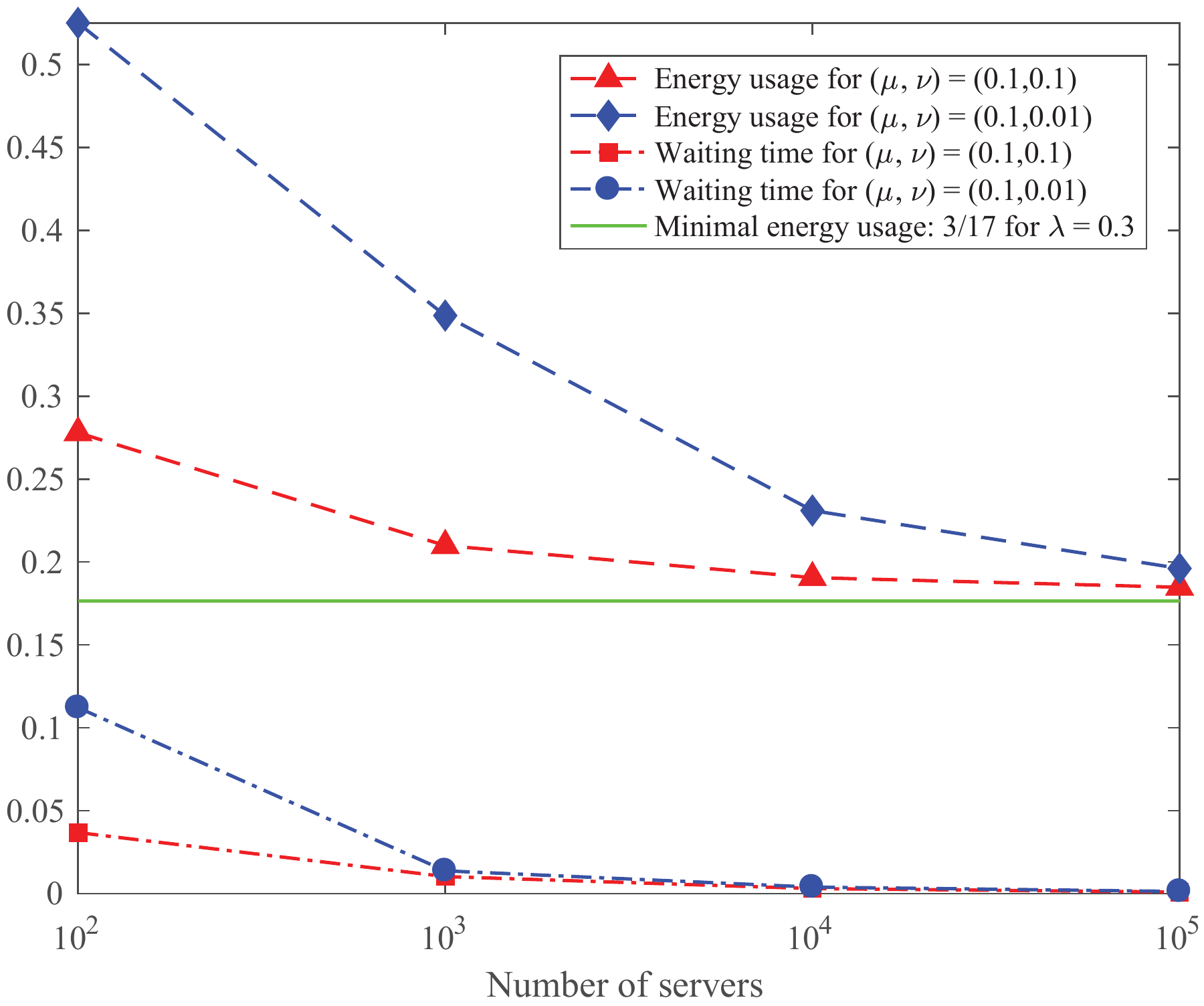}
\end{center}
\caption{Energy usage and mean waiting time for $N=10^2,10^3,10^4,10^5$ servers, mean standby period $\mu^{-1} = 10$, and mean setup periods $\nu^{-1}= 10, 100$.}
\label{fig:increasing-N}
\end{figure}

In Figure~\ref{fig:increasing-N} average values of the performance metrics, taken over time 0 to 250, have been plotted.
 We can clearly observe that both performance metrics approach the asymptotic values associated with the fixed point of the fluid limit as the number of servers grows large.
Comparison of the results for $\nu = 0.01$ and $\nu = 0.1$ shows that the convergence is substantially faster, and the performance correspondingly closer to the asymptotic lower bound, for shorter setup periods.
This is a manifestation of the fact that, even though the fraction of servers in setup mode vanishes in the limit for any value of $\nu$, the actual fraction for a given finite value of $N$ tends to increase with the mean setup period.
This in turn means that in order for the fluid limit values to be approached within a certain margin, the required value of $N$ increases with the mean setup period, as reflected in Figure~\ref{fig:increasing-N}.
 \\

\begin{figure*}
\begin{center}$
\begin{array}{c}
\includegraphics[width=90mm]{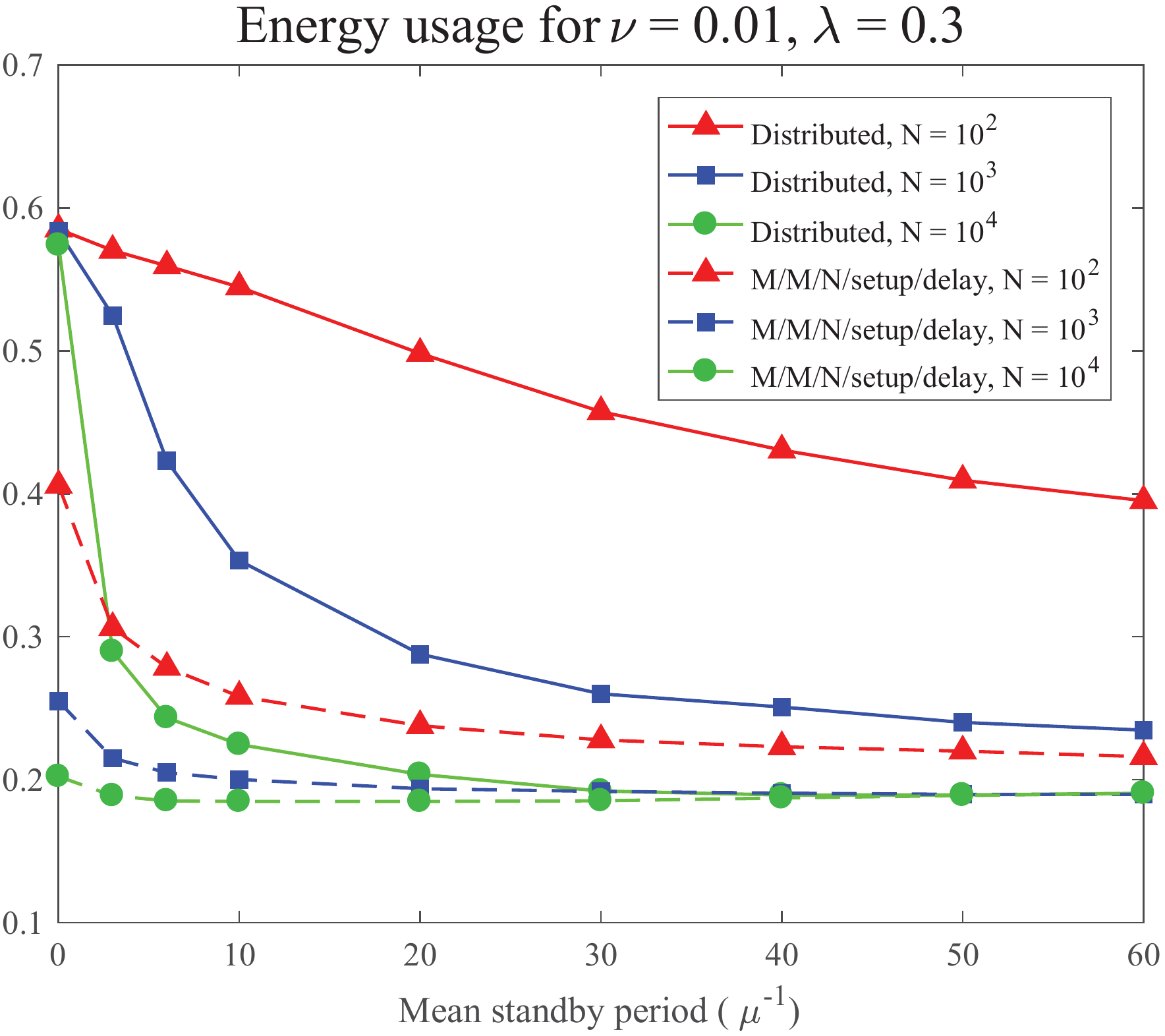}\\
\includegraphics[width=90mm]{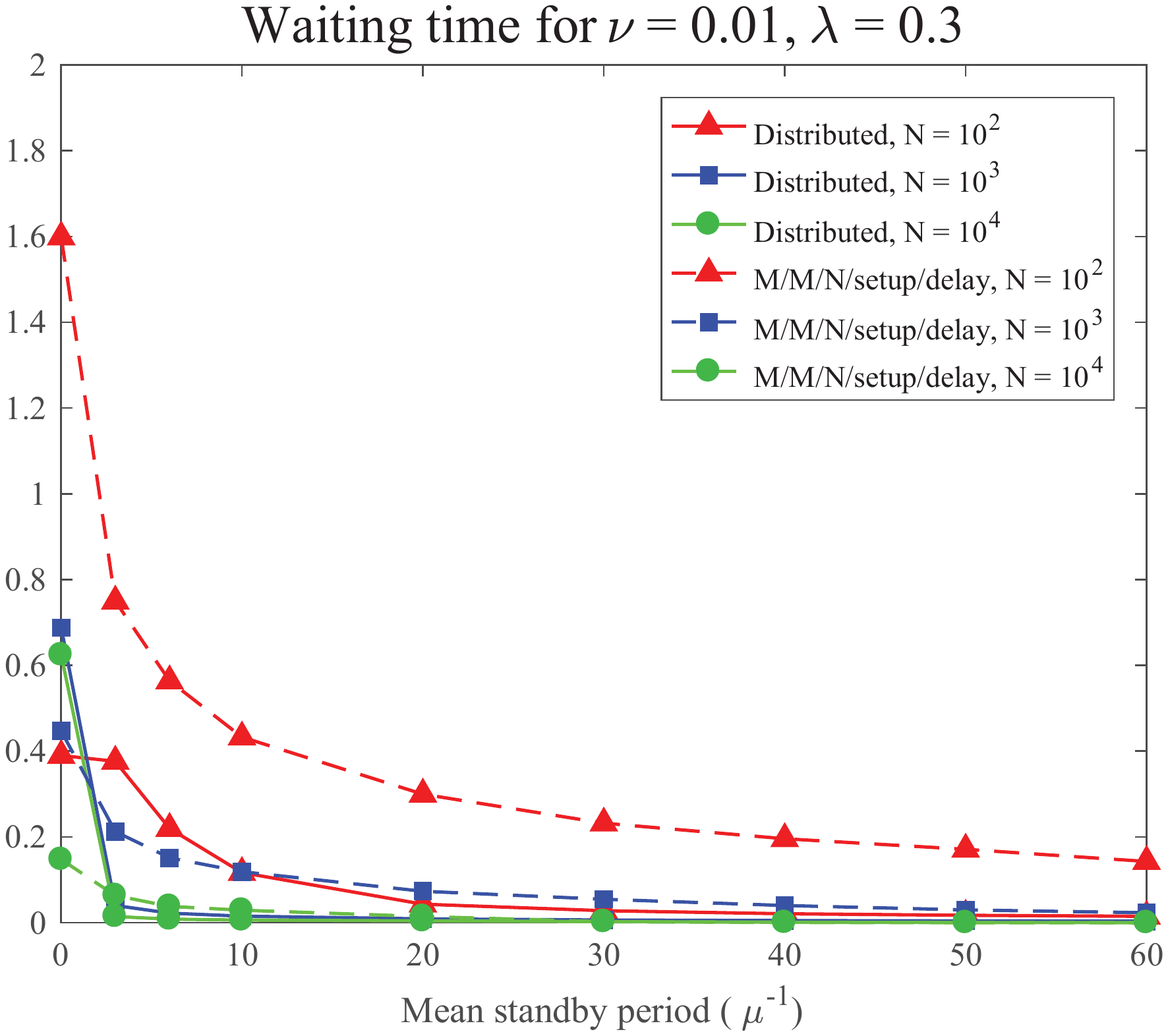}
\end{array}$
\end{center}
\caption{Comparison between TABS and M/M/N/setup/delayedoff schemes as functions  of the mean standby period $\mu^{-1}$ in terms of mean energy consumption and waiting time, for mean setup periods $\nu^{-1} =100$,  $N = 10^2,10^3,10^4$ servers.}
\label{fig:power-sig1}
\end{figure*}

\begin{figure*}
\begin{center}$
\begin{array}{c}
\includegraphics[width=90mm]{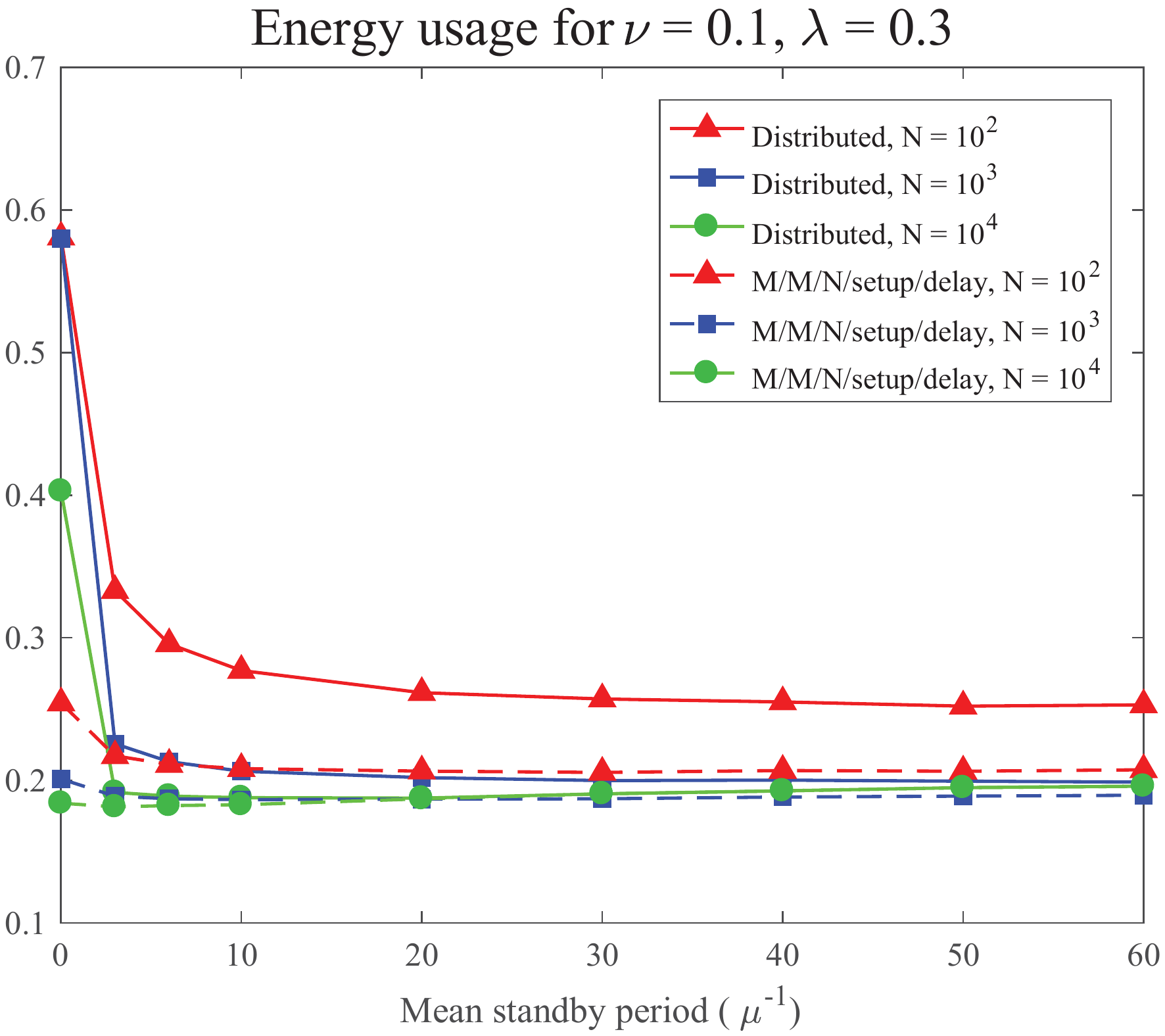}\\
\includegraphics[width=90mm]{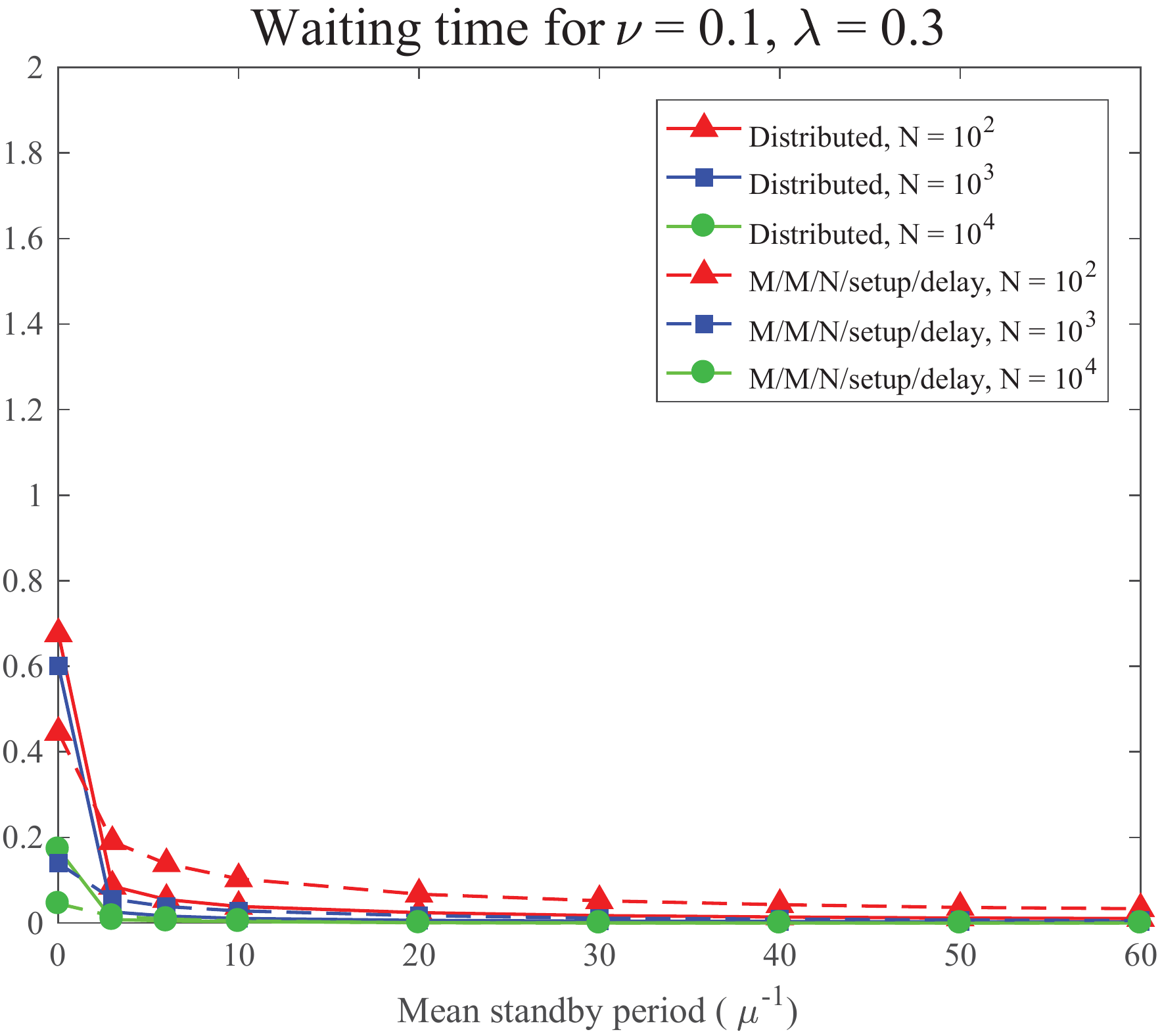}
\end{array}$
\end{center}
\caption{Comparison between TABS and M/M/N/setup/delayedoff schemes as functions  of the mean standby period $\mu^{-1}$ in terms of mean energy consumption and waiting time, for mean setup periods $\nu^{-1} =10$,  $N = 10^2,10^3,10^4$ servers.}
\label{fig:power-sig2}
\end{figure*}

\begin{figure*}
\begin{center}$
\begin{array}{c}
\includegraphics[width=90mm]{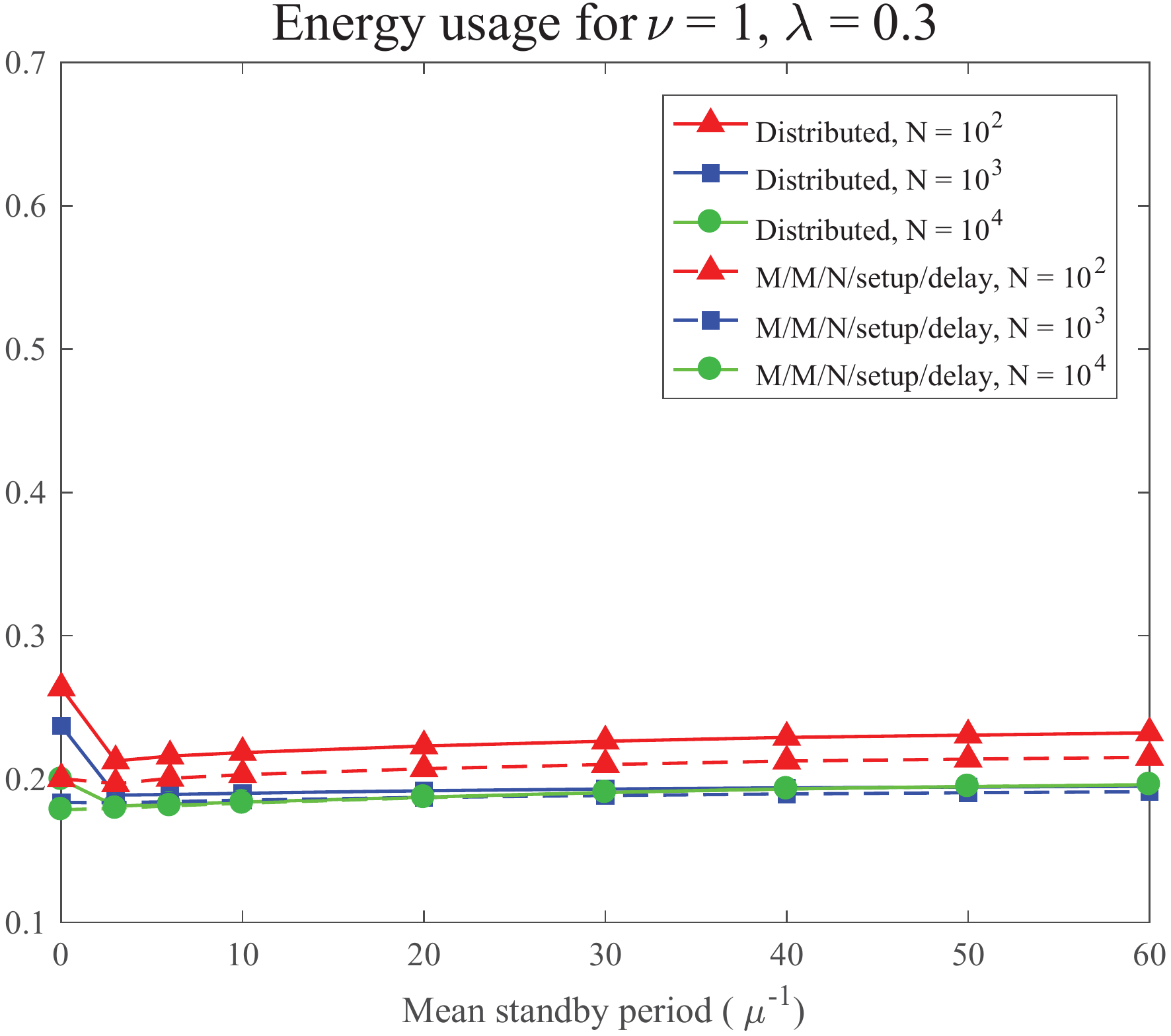}\\
\includegraphics[width=90mm]{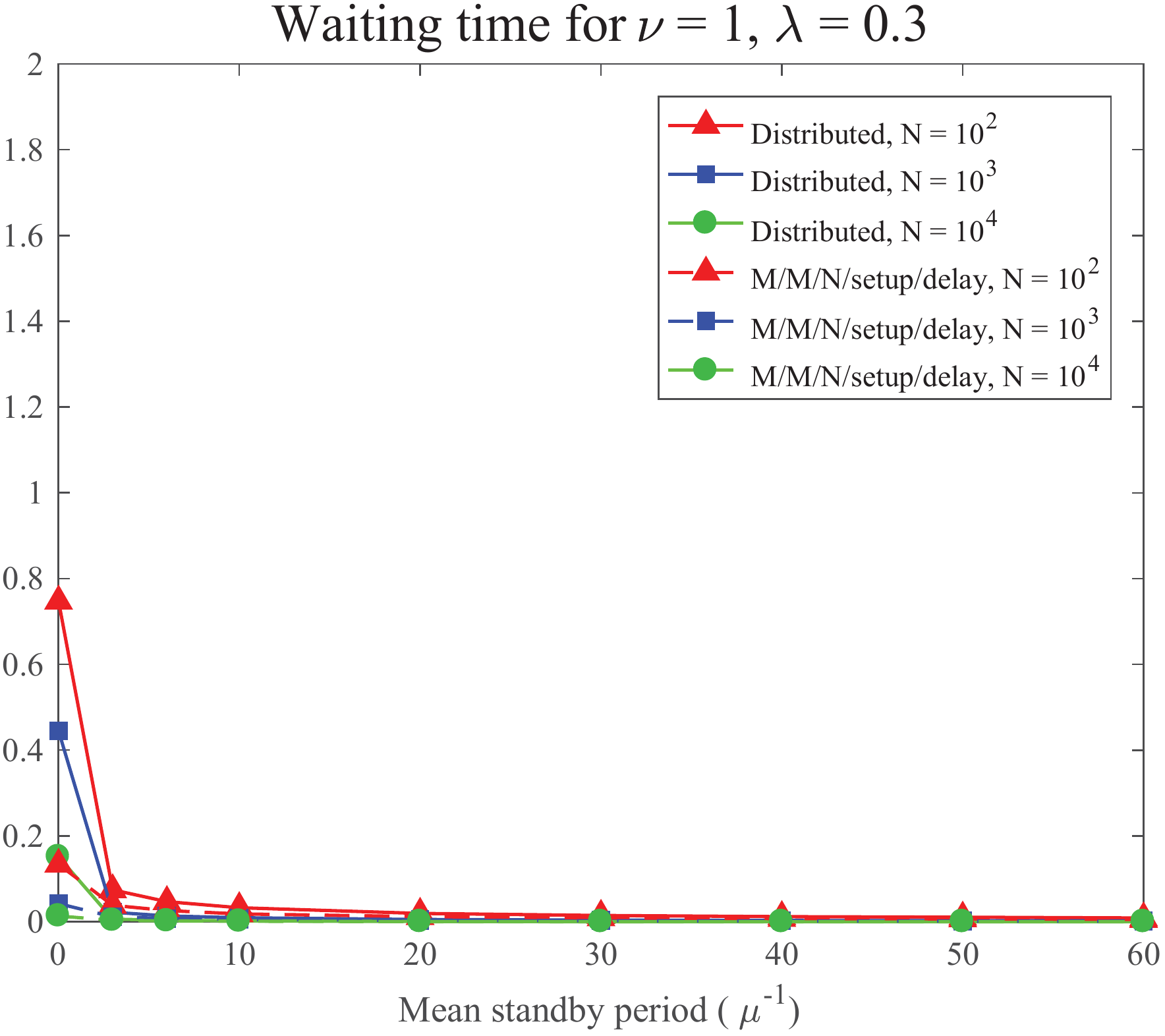}
\end{array}$
\end{center}
\caption{Comparison between TABS and M/M/N/setup/delayedoff schemes as functions  of the mean standby period $\mu^{-1}$ in terms of mean energy consumption and waiting time, for mean setup periods $\nu^{-1} =1$,  $N = 10^2,10^3,10^4$ servers.}
\label{fig:power-sig3}
\end{figure*}

In order to further examine the above observations and also investigate the impact of the mean standby period, we present in Figures~\ref{fig:power-sig1}-\ref{fig:power-sig3} the expected waiting time of tasks $\expt[W^N]$ and energy consumption $\expt[P^N]$ for $\lambda = 0.3$ and various values of $N$ and $\nu$, as a function of the mean standby period $\mu^{-1}$.
The results are based on 100 to 200 independent simulation runs, and we confirmed through  careful inspection that the numbers in fact did not show significant variation across runs.
In order to examine the impact of the load, we have also conducted experiments for $\lambda = 0.9$ which show qualitatively similar results, and hence are omitted.\\

The performance impact of the mean standby period $\mu^{-1}$ appears to be somewhat less pronounced.
Both performance metrics generally tend to improve as the mean standby period increases, although the energy consumption starts to slightly rise when the standby period increases above a certain level in scenarios with extremely short setup periods.
The latter observation may be explained as follows. 
For finite $N$-values, if the standby period is extremely small relative to the setup period, then the servers tend to deactivate too often, and as a result, setup procedures are also initiated too often (which in turn involve a relatively long time to become idle-on). 
Note that the servers in setup mode use $P_{\full}$ while providing no service. 
Thus the energy usage decreases by choosing longer standby periods (smaller~$\mu$).
On the other hand, again for small $N$-values, very long standby periods (smaller $\mu$) are not good either.
The reason in this case is straightforward; the idle-on servers will unnecessarily remain idle for a long time, and thus substantially increase energy usage with very little gain in the performance (reduction in waiting time).\\

As mentioned above, the required value of $N$ for the fluid-limit regime to kick in increases with the mean setup period, and broadly speaking, the asymptotic values are approached within a fairly close margin for $N = 10^3$ servers, except when the setup periods are long or the standby periods are extremely short.  
By implication, for scenarios with $N = 10^3$ or more servers, the TABS scheme delivers near-optimal performance in terms of energy consumption and waiting time, provided the setup periods are not too long and the standby periods are not too short.
It is worth observing that setup periods are basically determined by hardware factors and system constraints, while standby periods are design parameters that can be set in a largely arbitrary fashion.  Based on the above observations, a simple practical guideline is to set standby periods to relatively long values.\\

For smaller numbers of servers, long setup periods, or extremely short standby periods, finite-$N$ effects manifest themselves, and the actual performance metrics will differ from the fluid-limit values.
This does not imply though that the performance of the TABS scheme is necessarily far from optimal, since the absolute lower bound attained in the fluid limit may simply not be achievable by any scheme at all for small $N$ values.
\\

\noindent
\textbf{Comparison with centralized queue-driven strategies.}
To compare the performance in distributed systems under the TABS scheme with that of the corresponding pooled system under the M/M/N/setup/delayedoff mechanism, 
we also present in Figures~\ref{fig:power-sig1}-\ref{fig:power-sig3} the relevant metrics for the latter scenario. 
Quite surprisingly, even for moderate values of the total number of servers $N$, the performance metrics in a non-work-conserving scenario under the TABS scheme are very close to those for the M/M/N/setup/delayedoff mechanism.
Thus, the TABS scheme provides a significant energy saving in distributed systems which is comparable with that in a work-conserving pooled system, while achieving near zero waiting times as well.
In fact, it is interesting to observe that for relatively long setup periods the waiting time in the distributed system under the TABS scheme is even lower than for the M/M/N/setup/delayedoff mechanism!
This can be understood from the dynamics of the two systems as follows.
When an incoming task does not find an idle server, in both systems an idle-off server $s$ (if available) is switched to the setup mode.
By the time $s$ completes the setup procedure and turns idle-on, in the pooled system if a service completion occurs, then the task is assigned to that new idle-on server and the setup procedure of $s$ is discontinued. 
Therefore, when a next arrival occurs, the setup procedure must be initiated again.
As a result, this might cause the effective average waiting time to become higher.
On the other hand, in the distributed system once a setup procedure is initiated, it is completed in any event. 
This explains why for relatively long setup periods the TABS scheme provides a lower waiting time than the M/M/N/setup/delayedoff mechanism.

\section{Proofs}\label{sec:proofs}

\subsection{Fluid convergence}\label{app:conv}

The proof of Theorem~\ref{th: fluid} consists of describing the evolution of the system as a suitable time-changed Poisson process, which can be further decomposed into a martingale part and a drift part.
This formulation can be viewed as a \emph{density-dependent population process} (cf.~\cite[Chapter 11]{EK2009}).
 The martingale fluctuations become negligible on the fluid scale, and the drift terms converge to deterministic limits.
While the convergence of the martingale fluctuations is fairly straightforward to show, the analysis of the drift term is rather involved since the derivative of the drift is not continuous. 
As a result, the classical approaches developed by Kurtz~\cite{EK2009} cannot be applied in the current scenario.
In the literature, these situations have been tackled in various different ways~\cite{HK94, PW13, K92, GG12, TX11, GG10, B16, BG16}.
In particular, we leverage the time-scale separation techniques developed in~\cite{HK94} in order to identify the limits of drift terms.\\

 First, we verify the existence of the coefficients $p_i(\cdot,\cdot,\cdot)$ for all $t\geq 0$, $i = 1,2,\ldots,B$.
From the assumptions of Theorem~\ref{th: fluid}, and the fact that $\lambda(t)$ is bounded away from 0 (by some $\lambda_{\min}$ say), we claim that if $q_1(0)=q_1^\infty>0$, then $q_1(t)>0$ for all $t\geq 0$. 
To see this, it is enough to observe that in the fluid limit the rate of change of $q_1(t)$ is non-negative whenever $q_1(t)<\lambda_{\min}$.
Indeed, if $q_1(t)<\lambda_{\min}$, then 
\begin{align*}
\lambda(t)p_0(\qq(t),\dd(t),\lambda(t))-(q_1(t)-q_2(t))&\geq \min\{\lambda(t) - (q_1(t)-q_2(t)), \delta_1(t)\nu\}\\
&\geq \min\{\lambda_{\min}-q_1(t),\delta_1\nu\}\geq 0,
\end{align*}
and thus the claim follows.
Therefore below we will prove Theorem~\ref{th: fluid} until the time $q_1^N$ hits 0, and the above argument then shows that
if $q_1^N(0)\pto q_1^\infty>0$, then on any finite time interval $[0,T]$, with probability tending to 1, the process $q_1^N(\cdot)$ is bounded away from 0, proving the theorem for any finite time interval.
  \\

Let us introduce the variables $U^N(t) = N-Q_1^N(t)-\Delta_0^N(t)-\Delta_1^N(t)$, $u^N(t) = U^N(t)/N$,
$I_0^N(t) = \ind{U^N(t)>0}$, and $I_1^N(t) = \ind{\Delta_0^N(t)>0}$.
Note that $U^N(t)$ represents the number of idle-on servers at time $t$.\\

\noindent
\textbf{Martingale representation.}
For a unit-rate Poisson process $\big\{\mathcal{N}(t)\big\}_{t\geq 0}$ and a real-valued c\`adl\`ag process $\{A(t)\}_{t\geq 0}$, the random time-change~\cite{PTRW07, EK2009} is the unique process $\big\{\mathcal{N}(\int_0^tA(s)\dif s)\big\}_{t\geq 0}$ such that
\begin{equation}\label{def:timechange}
 \mathcal{N}\Big(\int_0^tA(s)\dif s\Big) - \int_0^tA(s)\dif s \quad\text{ is a martingale.}
\end{equation} 
Thus the evolution of the system is described by  
\begin{eq}\label{eq:poisson-descr}
Q^N_1(t)&=Q^N_1(0)+\mathcal{N}_{ A}\left(\int_0^t (1-I_0^N(s))\lambda_N(s)\dif s\right)\\
&\hspace{5cm}-\mathcal{N}_{ 1, D}\left(\int_0^t(Q^N_1(s)-Q^N_2(s))\dif s\right),\\
Q^N_i(t)&=Q^N_i(0)+\mathcal{N}_{  A}\left(\int_0^t  I_0^N(s)\frac{Q^N_{i-1}(s)-Q^N_i(s)}{Q^N_1(s)}\lambda_N(s)\dif s\right) \\ 
&\hspace{3cm}-\mathcal{N}_{ i, D}\left(\int_0^t (Q^N_i(s)-Q^N_{i+1}(s))\dif s\right),\quad i = 2,\dots,B,\\
\Delta_0^N(t)&=\Delta_0^N(0)+\mathcal{N}_{0}\left(\mu\int_0^t U^N(s)\dif s\right)-\mathcal{N}_{ A}\left(\int_0^tI_0^N(s)I_1^N(s)\lambda_N(s)\dif s\right),\\
\Delta_1^N(t) &= \Delta_1^N(0)+\mathcal{N}_{ A}\left(\int_0^t I_0^N(s)I_1^N(s)\lambda_N(s)\dif s\right)-\mathcal{N}_{1}\left(\nu\int_0^t\Delta_1^N(s)\dif s\right),
\end{eq}
where $\mathcal{N}_{ A}$, $\mathcal{N}_{ i,D}$ for $i= 1,\dots,B$, $\mathcal{N}_{0}$, $\mathcal{N}_{1}$ are independent unit-rate Poisson processes.
Using~\eqref{def:timechange} and \eqref{eq:poisson-descr}, we obtain the martingale representation of the process as 
\begin{eq}\label{eq:mart}
Q^N_1(t)&=Q^N_1(0)+\mathcal{M}_{  A} (t) - \mathcal{M}_{ 1, D}(t) + \int_0^t (1-I_0^N(s))\lambda_N(s)\dif s\\
&\hspace{6cm}-\int_0^t(Q^N_1(s)-Q^N_2(s))\dif s,\\
Q^N_i(t)&=Q^N_i(0)+\mathcal{M}_{ A}(t)-\mathcal{M}_{ i, D}(t)+\int_0^t I_0^N(s)\frac{Q^N_{i-1}(s)-Q^N_i(s)}{Q^N_1(s)}\lambda_N(s)\dif s\\
&\hspace{3cm}-\int_0^t (Q^N_i(s)-Q^N_{i+1}(s))\dif s, \quad i = 2,\dots,B,\
\end{eq}
\begin{eq}
\Delta_0^N(t)&=\Delta_0^N(0)+\mathcal{M}_{0}(t) - \mathcal{M}_{ A}(t) +\mu\int_0^t U^N(s)\dif s -\int_0^tI_0^N(s)I_1^N(s)\lambda_N(s)\dif s,\nonumber\\
\Delta_1^N(t) &= \Delta_1^N(0)+\mathcal{M}_{  A}(t)- \mathcal{M}_{1}(t)+\int_0^t I_0^N(s)I_1^N(s)\lambda_N(s)\dif s-\nu\int_0^t\Delta_1^N(s)\dif s,\nonumber
\end{eq}
where recall that $\mathcal{M}_{A}$, $\mathcal{M}_{0}$, $\mathcal{M}_{1}$, $\mathcal{M}_{i,D}$ for $i = 1,\dots,B$ are square-integrable martingales. 
The fluid-scaled martingale decomposition is thus given by
\begin{align}\label{eq:fluid mart}
q^N_1(t)&=q^N_1(0)+\frac{1}{N}\left(\mathcal{M}_{  A} (t) - \mathcal{M}_{ 1, D}(t)\right) + \int_0^t (1-I_0^N(s))\lambda(s)\dif s\\
&\hspace{6cm}-\int_0^t(q^N_1(s)-q^N_2(s))\dif s,\nonumber\\
q^N_i(t)&=q^N_i(0)+\frac{1}{N}\left(\mathcal{M}_{ A}(t)-\mathcal{M}_{ i, D}(t)\right)+\int_0^t I_0^N(s)\frac{q^N_{i-1}(s)-q^N_i(s)}{q^N_1(s)}\lambda(s)\dif s\\
&\hspace{4cm}-\int_0^t (q^N_i(s)-q^N_{i+1}(s))\dif s, \quad i = 2,\dots,B,\nonumber\\
\delta_0^N(t)&=\delta_0^N(0)+\frac{1}{N}\left(\mathcal{M}_{0}(t) - \mathcal{M}_{ A}(t) \right) +\mu\int_0^t u^N(s)\dif s -\int_0^tI_0^N(s)I_1^N(s)\lambda(s)\dif s,\nonumber\\
\delta_1^N(t) &= \delta_1^N(0)+\frac{1}{N}\left(\mathcal{M}_{  A}(t)- \mathcal{M}_{1}(t)\right)+\int_0^t I_0^N(s)I_1^N(s)\lambda(s)\dif s-\nu\int_0^t\delta_1^N(s)\dif s.
\end{align}

\noindent
\textbf{Random measure representation.}
We will now write the system evolution equation in terms of a suitable random measure.
The transition rates of the process $\{\ZZ^N(t)\}_{t\geq 0}:=\{(\Delta_0^N(t),U^N(t))\}_{t\geq 0}$ are described as follows.
\begin{enumerate}[(i)]
 \item When an idle server turns-off, $\Delta_0^N$ increases by one and $U^N$ decreases by one, and this occurs at rate $N\mu U^N$;
 \item When a server is requested to initiate the setup procedure, $U^N$ must be zero at that epoch. Thus,  $\Delta_0^N$ decreases by one while $U^N$ remains unchanged, and this occurs at rate $\lambda_N(t)\ind{U^N = 0 , \Delta_0^N> 0}$;
 \item When a busy server becomes idle, or a server finishes its setup procedure to become idle-on,  $\Delta_0^N$ remains unchanged while $U^N$ increases by one, and this occurs at rate $N(q_1^N-q_2^N+\nu\delta_1^N)$;
 \item When an arriving task is assigned to an idle-on server, $\Delta_0^N$ remains unchanged while $U^N$ decreases by one, and this occurs at rate $\lambda_N(t)\ind{U^N>0}$.
\end{enumerate}

    Let $\bar{\mathbb{Z}}_+=\mathbb{Z}_+\cup\{\infty\}$ denote the one-point compactification of the set of non-negative integers, equipped with the Euclidean metric, and the Borel $\sigma$-algebra $\mathfrak{B}$, induced by the mapping $f:\bar{\mathbb{Z}}_+\to[0,1]$ given by $f(x)=1/(x+1)$. 
    Let $\mathbf{v}^N(t)$ denote the vector $(\mathbf{q}^N(t),\dd^N(t))$.
  
  Observe that $\{(\mathbf{v}^N(t),\ZZ^N(t))\}_{t\geq 0}$ is a  Markov process defined on $E\times\bar{\mathbb{Z}}_+^2$. 
  Further, equip $[0,\infty)$ with the usual Euclidean metric and the Borel $\sigma$-algebra $\mathfrak{T}$.
We define a random measure $\alpha^N$ on the product space $[0,\infty)\times\bar{\mathbb{Z}}_+^2$ by
\begin{equation}
\alpha^N(A_1\times A_2):=\int_{A_1} \ind{\ZZ^N(s)\in A_2}\dif s,
\end{equation}
for $A_1\in \mathfrak{T}$, $A_2\in\mathfrak{B}$. 
Define
\begin{align*}
\mathcal{R}_1 &= \{(z_1,z_2)\in \Z_+^2: z_2 = 0\},\qquad \mathcal{R}_2 = \{(z_1,z_2)\in \Z_+^2: z_2 = 0, z_1>0\}.
\end{align*}
Note that the process $\{\ZZ^N(t)\}_{t\geq 0}$ determines the system constraints (indicator terms $I_0^N$ and $I_1^N$) in~\eqref{eq:fluid mart}. 
Thus, the process $\big\{(\qq^N(t),\dd^N(t))\big\}_{t\geq 0}$ can be written in terms of the random measure $\alpha^N$~as 
\begin{align}\label{eq:fluid mart2}
q^N_1(t)&=q^N_1(0)+\frac{1}{N}\left(\mathcal{M}_{  A} (t) - \mathcal{M}_{ 1, D}(t)\right) + \int_{[0,t]\times\mathcal{R}_1^c} \lambda(s)\dif \alpha^N\\
&\hspace{6cm}-\int_0^t(q^N_1(s)-q^N_2(s))\dif s,\nonumber\\
q^N_i(t)&=q^N_i(0)+\frac{1}{N}\left(\mathcal{M}_{ A}(t)-\mathcal{M}_{ i, D}(t)\right)+\int_{[0,t]\times\mathcal{R}_1} \frac{q^N_{i-1}(s)-q^N_i(s)}{q^N_1(s)}\lambda(s)\dif \alpha^N\\
&\hspace{3cm}-\int_0^t (q^N_i(s)-q^N_{i+1}(s))\dif s,\quad i = 2,\dots,B,\nonumber\\
\delta_0^N(t)&=\delta_0^N(0)+\frac{1}{N}\left(\mathcal{M}_{0}(t) - \mathcal{M}_{ A}(t) \right) +\mu\int_0^t u^N(s)\dif s -\int_{[0,t]\times\mathcal{R}_2}\lambda(s)\dif \alpha^N,\nonumber\\
\delta_1^N(t) &= \delta_1^N(0)+\frac{1}{N}\left(\mathcal{M}_{  A}(t)- \mathcal{M}_{1}(t)\right)+\int_{[0,t]\times\mathcal{R}_2}\lambda(s)\dif \alpha^N-\nu\int_0^t\delta_1^N(s)\dif s.
\end{align}
We first show that the scaled martingale parts converge to zero in probability, as $N\to\infty$.
\begin{proposition}
\label{prop:mart zero1-sig}
For any $T\geq 0$, $\sup_{t\in [0,T]}|\mathcal{M}_{ i,D}(t)|/N$, for all $i = 1,\dots B$  and $\sup_{t\in [0,T]}|\mathcal{M}_{k}(t)|/N$ for $k = A,0,1$ converge in probability to 0.
\end{proposition}
\begin{proof}
 We only give the proof for $\mathcal{M}_{ A}$ and the other cases can be proved similarly. Fix any $T>0$ and $\eta >0$.
 The proof makes use of the fact that the predictable quadratic variation process of a time-changed Poisson process is given by its compensator~\cite[Lemma 3.2]{PTRW07}.
  Using Doob's Martingale inequality~\cite[Theorem 1.9.1.3]{LS89}, we have
 \begin{align*}
  \mathbb{P}\ \Big(\sup_{t\in [0,T]}\frac{\left| \mathcal{M}_{ A}(t)\right|}{N}>\varepsilon\Big)&\leq \frac{1}{N^2\varepsilon^2}\E{\langle \mathcal{M}_{ A} \rangle_T}
  \leq \frac{NT\sup_{t\in [0,T]}\lambda(t)}{N^2\varepsilon^2}\to 0,
 \end{align*}and the proof follows.
\end{proof}

Let $\mathfrak{L}$ denote the space of all measures $\gamma$ on $[0,\infty)\times \bar{\mathbb{Z}}_+^2$ satisfying $\gamma([0,t]\times\bar{\mathbb{Z}}_+^2)= t$, endowed with the topology corresponding to weak convergence of measures restricted to $[0,t]\times \bar{\mathbb{Z}}_+^2$ for each $t$.  We have the following lemma:
\begin{lemma}[Relative compactness]
\label{lem:rel compactness}
Suppose that $\mathbf{v}^N(0)$ converges weakly to $\mathbf{v}^\infty = (\qq^\infty,\dd^\infty)\in E$ as $N\to\infty$, with $q_1^\infty>0$. Then the sequence of processes $\{(\mathbf{v}^N(\cdot),\alpha^N)\}_{N\geq 1}$ is relatively compact in $D_{E}[0,\infty)\times\mathfrak{L}$ and the limit $(\mathbf{v}(\cdot),\alpha)$ of any convergent subsequence satisfies
\begin{eq}\label{eq:rel compact-sig}
q_1(t)&=q_1^\infty+ \int_{[0,t]\times\mathcal{R}_1^c} \lambda(s)\dif \alpha-\int_0^t(q_1(s)-q_2(s))\dif s\\
q_i(t)&=q_i^\infty+\int_{[0,t]\times\mathcal{R}_1} \frac{q_{i-1}(s)-q_i(s)}{q_1(s)}\lambda(s)\dif \alpha
 -\int_0^t (q_i(s)-q_{i+1}(s))\dif s, \\
 &\hspace{8cm} i = 2,\dots,B,\\
\delta_0(t)&=\delta_0^\infty+\mu\int_0^t u(s)\dif s -\int_{[0,t]\times\mathcal{R}_2}\lambda(s)\dif \alpha \\
\delta_1(t)& = \delta_1^\infty+\int_{[0,t]\times\mathcal{R}_2}\lambda(s)\dif \alpha-\nu\int_0^t\delta_1(s)\dif s,
\end{eq}
with $u(t) = 1-q_1(t)-\delta_0(t)-\delta_1(t).$
\end{lemma}

\noindent
\textbf{
Conditions of relative compactness.}
To prove Lemma~\ref{lem:rel compactness}, we verify the conditions of relative compactness from~\cite[Corollary~3.7.4]{EK2009}. 
Let $(E,r)$ be a complete and separable metric space. For any $x\in D_E[0,\infty)$, $\kappa >0$ and $T>0$, define
\begin{equation}\label{eq:mod-continuity-sig}
w'(x,\kappa,T)=\inf_{\{t_i\}}\max_i\sup_{s,t\in[t_{i-1},t_i)}r(x(s),x(t)),
\end{equation}
where $\{t_i\}$ ranges over all partitions of the form $0=t_0<t_1<\ldots<t_{n-1}<T\leq t_n$ with $\min_{1\leq i\leq n}(t_i-t_{i-1})>\kappa$ and $n\geq 1$.
 Below we state the conditions for the sake of completeness.
\begin{theorem}[{\cite[Corollary~3.7.4]{EK2009}}]\label{th:from EK-sig}
Let $(E,r)$ be complete and separable, and let $\{X_n\}_{n\geq 1}$ be a family of processes with sample paths in $D_E[0,\infty)$. Then $\{X_n\}_{n\geq 1}$ is relatively compact if and only if the following two conditions hold:
\begin{enumerate}[{\normalfont (a)}]
\item For every $\eta>0$ and rational $t\geq 0$, there exists a compact set $\Gamma_{\eta, t}\subset E$ such that $$\varliminf_{n\to\infty}\Pro{X_n(t)\in\Gamma_{\eta, t}}\geq 1-\eta.$$
\item For every $\eta>0$ and $T>0$, there exists $\kappa>0$ such that
$$\varlimsup_{n\to\infty}\Pro{w'(X_n,\kappa, T)\geq\eta}\leq\eta.$$
\end{enumerate}
\end{theorem}

\begin{proof}[Proof of Lemma~\ref{lem:rel compactness}]
 From \cite[Proposition 3.2.4]{EK2009} observe that, to prove the relative compactness of the process $(\vv^N(\cdot),\alpha^N)$, it is enough to prove relative compactness of the individual components.

Let $\mathfrak{L}_t$ denote the collection of measures $\gamma^t$ where $\gamma^t$ is the restriction of $\gamma$ on $[0,t]\times\bar{\mathbb{Z}}_+^2$. Note that, by Prohorov's theorem,  $\mathfrak{L}_t$ is compact, since $\bar{\mathbb{Z}}_+^2$ is compact. The topology on $\mathfrak{L}$ is defined such that any sequence $\{\gamma_N\}_{N\geq 1}$ is relatively compact in $\mathfrak{L}$ if and only if $\{\gamma_N^t\}_{N\geq 1}$ is relatively compact in $\mathfrak{L}_t$ for any $t>0$. Since $\mathfrak{L}_t$ is compact, any sequence $\{\gamma_N\}_{N\geq 1}$ is relatively compact in $\mathfrak{L}$.  Thus, the relative compactness of $\alpha^N$ follows.
To see the relative compactness of $\{\mathbf{v}^N(\cdot)\}_{n\geq 1}$, first observe that $E$ is compact and hence the compact containment condition (a) of Theorem~\ref{th:from EK-sig} is satisfied trivially by taking $\Gamma_{\eta,t}\equiv E$. 

Let $\{\mathbf{M}^N(t)\}_{t\geq 0}$ denote the vector of all the martingale quantities appearing in \eqref{eq:fluid mart2}. Denote by $\|\cdot\|$, the Euclidean norm. For condition (b), we can see that, for any $0\leq t_1<t_2<\infty$,
\begin{equation}\label{mart-norm-ub-sig}
\|\mathbf{v}^N(t_1)-\mathbf{v}^N(t_2)\|\leq C (t_2-t_1)+\frac{1}{N}\|\mathbf{M}^N(t_1)-\mathbf{M}^N(t_2)\|,
\end{equation}
for a sufficiently large constant $C>0$ where we have used $q_i^N\leq 1$, for all $i$, $\lambda(t)$ is bounded, and the fact that $(q_{i-1}^N-q_i^N)/q_1^N \leq 1$.
From Proposition~\ref{prop:mart zero1-sig}, we get, for any $T\geq 0$,
$$\sup_{t\in[0,T]}\frac{1}{N}\|\mathbf{M}^N(t)\|\pto 0.$$
Now, the proof of the relative compactness of $(\mathbf{v}^N(t))_{t\geq 0}$ is complete if we can show that for any $\eta>0$, there exists a $\delta >0$ and a partition $(t_i)_{i\geq 1}$ with $\min_i|t_{i}-t_{i-1}|>\delta$ such that 
\begin{equation}
\varlimsup_{N\to\infty}\mathbb{P}\ \Big(\max_i \sup_{s,t\in [t_{i-1},t_i)}\|\mathbf{v}^N(s)-\mathbf{v}^N(t)\| \geq \eta\Big) < \eta.
\end{equation}
Now, \eqref{mart-norm-ub-sig} implies that, for any partition $(t_i)_{i\geq 1}$,
\begin{align*}
\max_i \sup_{s,t\in [t_{i-1},t_i)} \|\mathbf{v}^N(s)-\mathbf{v}^N(t)\|&\leq C \max_i (t_{i}-t_{i-1})+\zeta_N,
\end{align*}where $\Pro{\zeta_N>\eta/2}<\eta$ for all sufficiently large $N$. Now take $\delta = \eta/4C$ and any partition with $\max_i(t_i-t_{i-1})< \eta/2C$ and $\min_i(t_i-t_{i-1})>\delta$. Now on the event $\{\zeta_N\leq \eta/2\}$,  $$\max_i \sup_{s,t\in [t_{i-1},t_i)}\|\mathbf{v}^N(s)-\mathbf{v}^N(t)\| \leq \eta.$$
Therefore, for all sufficiently large $N$,
\begin{eq}
&\mathbb{P}\bigg(\max_i \sup_{s,t\in [t_{i-1},t_i)}\|\mathbf{v}^N(s)-\mathbf{v}^N(t)\| \geq \eta \bigg)
\leq \Pro{\zeta_N>\eta/2}\leq \eta,
\end{eq}and the proof of the relative compactness of $(\mathbf{v}^N(t))_{t\geq 0}$ is now complete.
The fact that the limit $(\vv,\alpha)$ of any convergent subsequence of $(\vv^N,\alpha^N)$ satisfies~\eqref{eq:rel compact-sig}, follows by applying the continuous-mapping theorem.
\end{proof}

We will now prove the fluid-limit result stated in Theorem~\ref{th: fluid}.

\begin{proof}[{Proof of Theorem~\ref{th: fluid}}]
Using \cite[Theorem 3]{HK94}, we can conclude that the measure $\alpha$ can be represented as
\begin{equation}
\alpha (A_1\times A_2) = \int_{A_1} \pi_{\qq(s),\dd(s)}(A_2)\dif s,
\end{equation}
for measurable subsets $A_1\subset[0,\infty)$, and $A_2\subset \bar{\Z}_+^2$, where for any $(\bld{q},\dd) \in E$, $\pi_{\bld{q},\dd}$ is given by some  stationary distribution of the Markov process with transitions  
\begin{eq}\label{eq:stationary-bd-chain}
 (Z_1,Z_2)\rightarrow 
\begin{cases}
(Z_1,Z_2)+(1,-1)& \mbox{at rate } \mu u\\
(Z_1,Z_2)+(-1,0)& \mbox{at rate } \lambda\ind{Z_2=0,Z_1>0}\\
(Z_1,Z_2)+(0,1)& \mbox{at rate } q_1-q_2+\nu\delta_1\\
(Z_1,Z_2)+(0,-1)& \mbox{at rate } \lambda \ind{Z_2>0},
\end{cases}
\end{eq}
with $u=1-q_1-\delta_0-\delta_1$.
Additionally, the measure $\pi_{\qq,\dd}$ satisfies $\pi_{\qq,\dd}(Z_2 = \infty) = 1$, if $u>0$ and $\pi_{\qq,\dd}(Z_1 = \infty) = 1$ if $\delta_0>0$.
Thus we will show that for any $(\qq,\dd)\in E$, $\pi_{\bld{q},\dd}$ is unique, and that $\pi_{\qq(s),\dd(s)}(\mathcal{R}_1^c)=p_{0}(\qq(s),\dd(s))$ and $\pi_{\qq(s),\dd(s)}(\mathcal{R}_2)=(1-p_{0}(\qq(s),\dd(s)))\ind{\delta_0(s)>0}$ as described in Theorem~\ref{th: fluid} (we have omitted the argument $\lambda(s)$ in $p(\cdot,\cdot,\cdot)$ to avoid cumbersome notation). 
We will verify the uniqueness of the stationary measure $\pi_{\qq,\dd}$ of the Markov process $(Z_1,Z_2)$ subsequently case-by-case.\\

\noindent
\textbf{Case I: $\bld{u>0,$ $\delta_0>0}$.}
In this case, by the definition of $\pi_{\qq,\dd}$ stated above, $\pi_{\qq,\dd} (Z_2 = Z_1 =\infty) = 1$. Thus, $\pi_{\qq,\dd}(\mathcal{R}_1) = \pi_{\qq,\dd}(\mathcal{R}_2) = 0$.\\

\noindent
\textbf{Case II: $\bld{u>0,$ $\delta_0=0}$.} 
Here by definition of $\pi_{\qq,\dd}$,  $\pi_{\qq,\dd} (Z_2  =\infty) = 1$. 
However, if $Z_2 = \infty$, then by~\eqref{eq:stationary-bd-chain},  $Z_1$ increases by one at rate $\mu u$, and decreases at rate 0. Since $\pi_{\qq,\dd}$ is the stationary measure, we also have $\pi_{\qq,\dd} (Z_1  =\infty) = 1$, and thus,  $\pi_{\qq,\dd}(\mathcal{R}_1) = \pi_{\qq,\dd}(\mathcal{R}_2) = 0$.\\

\noindent
\textbf{Case III: $\bld{u=0,$ $\delta_0>0}$.}
 In this case, $\pi_{\qq,\dd} (Z_1  =\infty) = 1$. 
Again note that if $Z_1 = \infty$, then by~\eqref{eq:stationary-bd-chain}, $Z_2$ increases by one at rate $q_1-q_2+\nu\delta_1$, and decreases by one at rate $\lambda\ind{Z_2>0}$. Thus, 
\begin{itemize}
\item if $q_1-q_2+\nu\delta_1\geq \lambda$, then $\pi_{\qq,\dd} (Z_2  =0) = 0$, and consequently, $\pi_{\qq,\dd}(\mathcal{R}_1) = \pi_{\qq,\dd}(\mathcal{R}_2) = 0$,
\item if $q_1-q_2+\nu\delta_1< \lambda$,  then $\pi_{\qq,\dd} (Z_2  = 0) = \lambda^{-1}(q_1-q_2+\nu\delta_1)$, and  $\pi_{\qq,\dd}(\mathcal{R}_1) = \pi_{\qq,\dd}(\mathcal{R}_2) = \lambda^{-1}(q_1-q_2+\nu\delta_1)$.
\end{itemize}

\noindent
\textbf{Case IV: $\bld{u=0,$ $\delta_0=0}$.} 
Observe that in this case, due to physical constraints, it must be that $\pi_{\qq,\dd}(\mathcal{R}_2) = 0$. 
To see this, recall the evolution equation from \eqref{eq:rel compact-sig}. 
Note that $\delta_0(t) = 0$ forces its derivative to be non-negative (since $\delta_0$ is non-negative), and thus $\delta_0'(t)\geq 0$.
Now, $\pi_{\qq(t),\dd(t)}(\mathcal{R}_2) > 0$ implies that $\delta_0'(t)<0$, and  hence, this leads to a contradiction.
Furthermore, $\pi_{\qq,\dd}(Z_2=0,Z_1>0) = 0$ implies that $\pi_{\qq,\dd}(Z_2=0) = \pi_{\qq,\dd}(Z_2=0,Z_1=0)$.  
Again, if $Z_1 = 0$, then by~\eqref{eq:stationary-bd-chain}, $Z_2$ increases by one at rate $q_1-q_2+\nu\delta_1$, and decreases by one at rate $\lambda\ind{Z_2>0}$. 
Thus, an argument similar to Case-III yields that $\pi_{\qq,\dd}(\mathcal{R}_1)  = 0$, if $q_1-q_2+\nu\delta_1\geq \lambda$, and $\pi_{\qq,\dd}(\mathcal{R}_1) = \lambda^{-1}(q_1-q_2+\nu\delta_1)$, if $q_1-q_2+\nu\delta_1<\lambda$.
Combining Cases I-IV, we have 
\begin{align*}
 \pi_{\qq,\dd}(\mathcal{R}_1) &= 1- p_0(\qq,\dd,\lambda), \quad \pi_{\qq,\dd}(\mathcal{R}_2) = \ind{\delta_0>0}\pi_{\qq,\dd}(\mathcal{R}_1),
\end{align*}
and the proof of Theorem~\ref{th: fluid} follows from Lemma~\ref{lem:rel compactness}.
\end{proof}

\begin{proof}[Proof of Theorem~\ref{th: fluid gen service}]The proof of Theorem~\ref{th: fluid gen service} is identical to the proof of Theorem~\ref{th: fluid}, which starts again by  establishing the martingale decomposition for $q_{ij}^N$ of the form 
\begin{align}\label{eq:fluid mart3}
q^N_{1,j}(t)&=q^N_{1,j}(0)+\frac{1}{N}\mathcal{M}_{1,j} (t)  + \int_{[0,t]\times\mathcal{R}_1^c} r_j\lambda(s)\dif \alpha^N - \gamma_j\int_0^t q_{1,j}^N(s)\dif s\\
&\hspace{-.2cm}+ \int_0^t\sum_{k=1}^K (q_{1,k}^N(s)-q_{2,k}^N(s))\gamma_kr_{k,j}\dif s+ \int_0^t\sum_{k=1}^K (q_{2, k}^N(s) - q_{3 k}^N(s))\gamma_kr_{k,0}r_j\dif s \nonumber\\
q^N_{i,j}(t)&=q^N_{i,j}(0)+\frac{1}{N}\mathcal{M}_{ i,j}(t)+\int_{[0,t]\times\mathcal{R}_1} \frac{q^N_{i-1,j}(s)-q^N_{i,j}(s)}{\sum_{j=1}^Kq^N_{1,j}(s)}r_j\lambda(s)\dif \alpha^N\\
&\hspace{2cm}+\int_0^t\sum_{k=1}^K (q_{ik}^N(s)-q_{i+1,k}^N(s))\gamma_kr_{k,j}\dif s - \gamma_j\int_0^t q_{i,j}^N(s)\dif s\nonumber\\
  &\hspace{4cm}+ \int_0^t\sum_{k=1}^K (q_{i+1, k}^N(s) - q_{i+2, k}^N(s))\gamma_kr_{k,0}r_j\dif s\nonumber\\
\delta_0^N(t)&=\delta_0^N(0)+\frac{1}{N}\mathcal{M}_{0}(t) +\mu\int_0^t \bigg(1-\sum_{j=1}^Kq_{1,j}^N(s)-\delta_0^N(s)-\delta_1^N(s)\bigg)\dif s \nonumber \\
&\hspace{7.5cm}-\int_{[0,t]\times\mathcal{R}_2}\lambda(s)\dif s,\\
\delta_1^N(t) &= \delta_1^N(0)+\frac{1}{N}\mathcal{M}_{1}(t)+\int_{[0,t]\times\mathcal{R}_2}\lambda(s)\dif s-\nu\int_0^t\delta_1^N(s)\dif s.
\end{align}
The definitions of the sets $\mathcal{R}_1$, $\mathcal{R}_2$ remain exactly the same. Thus the convergence result Lemma~\ref{lem:rel compactness} holds for $\qq^N=(q_{ij}^N)_{1\leq i\leq B, 1\leq j\leq K}$. The arguments for the time scale separation part remain unchanged as well, except the transition rate $(Z_1,Z_2)\rightarrow (Z_1,Z_2)+(0,1)$ in \eqref{eq:stationary-bd-chain} changes to $\sum_{j=1}^K(q_{1j}-q_{2j})+\nu\delta_1$.
\end{proof}

\subsection{Convergence of stationary distribution}\label{app:globstab}

\begin{proof}[Proof of Proposition~\ref{prop:glob-stab}]
The proof follows in three steps: in Lemma~\ref{lem:q1}, we show that $q_1(t)\to \lambda$ as $t\to\infty$, using this we show in Lemma~\ref{lem:q2} that $q_2(t)\to 0$, and then finally we deduce that $\delta_0(t)\to 1-\lambda$ and $\delta_1(t)\to 0$.\\

\begin{lemma}\label{lem:q1}
$q_1(t)\to \lambda$ as $t\to\infty$.
\end{lemma}
First we will establish that
$q_1(t)\to \lambda$ as $t\to\infty$.
The high-level intuition behind the proof can be described in two steps as follows.

(1) First we prove that $\liminf_{t\to\infty}q_1(t)\geq \lambda$.
Assume the contrary.
Because $q_1(t)$ can be shown to be non-decreasing when $q_1(t)\leq \lambda$, there must exist an $\varepsilon>0$, such that 
\begin{equation}\label{eq:contradict1}
q_1(t)\leq \lambda-\varepsilon\nu, \quad \forall\ t\geq 0.
\end{equation}
If $q_1(t)$ were to remain below $\lambda$ by a non-vanishing margin, then the (scaled) rate $q_1(t) - q_2(t)$ of busy servers turning idle-on would not be high enough to match the (scaled) rate $\lambda$ of incoming jobs.
If there are idle-on servers or sufficiently many servers in setup mode, we can still assign incoming jobs to idle-on servers, but this drives up the fraction of busy servers $q_1(t)$ and cannot continue indefinitely due to \eqref{eq:contradict1}.
This means that we cannot initiate an unbounded number of setup procedures.
Since we cannot continue to have idle-on servers either, this also implies that a non-vanishing fraction of the jobs cannot be assigned to idle servers, and hence we will initiate an unbounded number of setup procedures, hence contradiction. 

(2) Next we show that $\limsup_{t\to\infty} q_1(t)\leq \lambda$.
Suppose not, i.e., assume 
$$\limsup_{t\to\infty}q_1(t) = \lambda+\varepsilon$$ for some $\varepsilon>0$.
Recall that $q_1(t)$ is non-decreasing when $q_1(t)\leq \lambda$.
Hence, there must exist a $t_0$ such that $q_1(t)\geq \lambda$ $\forall\ t\geq t_0$.
If $q_1(t)$ were to get above $\lambda$ by a non-vanishing margin infinitely often, then the cumulative number of departures would exceed the cumulative number of arrivals by an infinite amount, which cannot occur since the (scaled) initial number of tasks is bounded.
\begin{proof}[Proof of Lemma~\ref{lem:q1}]
We first state four useful basic facts based on the fluid limit in Theorem~\ref{th: fluid}.
These are then used to prove Claims~\ref{claim:inf} and~\ref{claim:sup} which together imply Lemma~\ref{lem:q1}.
\begin{fact} \label{fact:nondec}
$q_1(t)$ is nondecreasing if $q_1(t)-q_2(t)\leq\lambda$. 
In particular, if $q_1(t)\leq\lambda$, then $q_1(t)$ is nondecreasing.
\end{fact}
\begin{claimproof}
Note that the rate of change of $q_1(t)$ is determined by $\lambda p_0(\qq(t), \dd(t))-q_1(t) +q_2(t)$. 
So it suffices to show that the latter quantity is non-negative when $q_1(t) - q_2(t)\leq \lambda$.
This follows directly from the fact that 
\begin{equation}\label{eq:lb-p0}
p_0(\qq(t), \dd(t))\geq \min\big\{\lambda^{-1}(\delta_1(t)\nu+q_1(t)-q_2(t)),1\big\}.
\end{equation}
\end{claimproof}
\noindent Define the subset $\mathcal{X}\subseteq E$ as
$$\mathcal{X}:= \Big\{(\qq,\dd)\in E: q_1 + \delta_0 + \delta_1 = 1,  \delta_1\nu+q_1-q_2\leq \lambda\Big\},$$
and denote by  $\indn{\cX}(\qq(s),\dd(s))$  the indicator of the event that $(\qq(s),\dd(s))\in \cX.$
Observe that $q_1(t)$ can be written as
\begin{equation}
\begin{split}
q_1(t)&= q_1(u) +\int_{u}^t\delta_1(s)\nu\indn{\cX}(\qq(s),\dd(s))\dif s \\
&\hspace{3cm}+\int_{u}^t[\lambda - q_1(s) + q_2(s)] \indn{\cX^c}(\qq(s),\dd(s))\dif s.
\end{split}
\end{equation}
The above representation leads to Facts~\ref{fact:2} and \ref{fact:3} stated below.
\begin{fact}\label{fact:2}
\begin{align*}
q_1(t)\geq q_1(u) +\int_{u}^t[\lambda - q_1(s) + q_2(s)] \indn{\cX^c}(\qq(s),\dd(s))\dif s.
\end{align*}
\end{fact}
\begin{fact}\label{fact:3}
\begin{align*}
q_1(t)\geq q_1(u)+\nu\int_{u}^t\delta_1(s)\dif s-(\nu+1)\int_{u}^t\indn{\cX^c}(\qq(s),\dd(s))\dif s.
\end{align*}
\end{fact}

\begin{fact}\label{fact:4}
For all sufficiently small $\varepsilon>0$,
\begin{align*}
\xi(t) &\geq  \int_0^t\Big(\lambda-\frac{\varepsilon \nu}{2}-q_1(s)\Big)\dif s 
- \int_0^t\ind{u(s)>0}\dif s - \int_0^t\ind{\delta_1(s)>\varepsilon/2}\dif s.
\end{align*}
\end{fact}
\begin{claimproof}
Observe that
\begin{align*}
\xi(t)&= \int_0^t\lambda(1-p_0(\qq(s),\dd(s),\lambda))\ind{\delta_0(s)>0}\dif s\\
&\geq \int_{0}^t\lambda(1-p_0(\qq(s),\dd(s),\lambda))\ind{\delta_0(s)>0, u(s) = 0,\delta_1(s)\leq \varepsilon/2}\dif s,
\end{align*}
and on the set $\{s:\delta_0(s)>0, u(s) = 0,\delta_1(s)\leq \varepsilon/2\}$ we have  $p_0(\qq(s),\dd(s),\lambda)\leq \lambda^{-1}(\varepsilon\nu/2+q_1(s))$. Therefore,
\begin{align*}
\xi(t)
&\geq \int_{0}^t\Big(\lambda-\frac{\varepsilon\nu}{2}-q_1(s)\Big)\ind{\delta_0(s)>0, u(s) = 0,\delta_1(s)\leq \varepsilon/2}\dif s.
\end{align*}
Moreover, if $\delta_0(s)=0, u(s) = 0,\delta_1(s)\leq \varepsilon/2,$ then $q_1(s)\geq 1-\varepsilon/2$, and for $\varepsilon<2(1-\lambda)/[1-\nu]^+$ we have $\lambda-\varepsilon\nu/2-q_1(s)<0$. Thus we finally obtain that
\begin{align*}
\xi(t)&\geq \int_{0}^t\Big(\lambda-\frac{\varepsilon \nu}{2}-q_1(s)\Big)\ind{ u(s) = 0,\delta_1(s)\leq \varepsilon/2}\dif s \\
&\geq \int_0^t\Big(\lambda-\frac{\varepsilon\nu}{2}-q_1(s)\Big)\dif s - \int_0^t\ind{u(s)>0}\dif s - \int_0^t\ind{\delta_1(s)>\varepsilon/2}\dif s,
\end{align*}
where the second inequality follows from $\lambda-\varepsilon\nu/2-q_1(s)\leq \lambda <1$.
\end{claimproof}
\noindent
In order to break down the proof of Lemma~\ref{lem:q1}, we will establish the following two claims.
\begin{claim}\label{claim:inf}
$\liminf_{t\to\infty}q_1(t)\geq \lambda$.
\end{claim}
\begin{claimproof}
Assume the contrary. 
Using Fact~\ref{fact:nondec}, $q_1(t)$ is non-decreasing when $q_1(t)\leq \lambda$, and thus there must exist an $\varepsilon>0$, such that 
\begin{equation}\label{eq:contradict2}
q_1(t)\leq \lambda-\varepsilon\nu, \quad \forall\ t\geq 0.
\end{equation}
By Fact~\ref{fact:2} there exist positive constants $K_1,K_2$ (possibly depending on $\varepsilon$) such that $\forall\ t\geq 0$
\begin{equation}\label{eq:K1-choice}
\int_0^t \indn{\cX^c}(\qq(s),\dd(s))\dif s<K_1
\implies \int_0^t\ind{u(s)>0}\dif s<K_1,
\end{equation}
and by Fact~\ref{fact:3} and \eqref{eq:K1-choice}
\begin{equation}\label{eq:K2-choice}
 \int_0^t\delta_1(s)\dif s<K_1 \implies \int_0^t\ind{\delta_1(s)>\frac{\varepsilon}{2}}\dif s<K_2.
\end{equation}
Note that since $\delta_1(t) = \delta_1(0) + \xi(t) - \nu\int_0^t\delta_1(s)\dif s,$
it must be the case that 
$$\limsup_{t\to\infty}\xi(t)<\infty.$$
On the other hand, Fact~\ref{fact:4}, together with \eqref{eq:K1-choice} and~\eqref{eq:K2-choice}, implies that $\xi(t)\to\infty$ as $t\to\infty$, which leads to a contradiction.
\end{claimproof}
\begin{claim}\label{claim:sup}
$\limsup_{t\to\infty} q_1(t)\leq \lambda$.
\end{claim}
\begin{claimproof}
Suppose not, i.e., $\limsup_{t\to\infty}q_1(t) = \lambda+\varepsilon$ for some $\varepsilon>0$.
Because $q_1(t)$ is non-decreasing by Fact~\ref{fact:nondec} when $q_1(t)\leq \lambda$, there must exist a $t_0$ such that $q_1(t)\geq \lambda$ $\forall\ t\geq t_0$.
In that case,
\begin{align*}
\sum_{i=1}^B q_i(t) 
&= \sum_{i=1}^Bq_i(t_0) + \lambda\int_{t_0}^t\sum_{i=1}^B p_{i-1}(\qq(s),\dd(s),\lambda)\dif s - \int_{t_0}^t q_1(s)\dif s\\
& \leq \sum_{i=1}^Bq_i(t_0)  - \int_{t_0}^t [q_1(s) - \lambda]^+\dif s,
\end{align*}
and thus,
$$\int_{t_0}^t [q_1(s) - \lambda]^+\dif s\leq \sum_{i=1}^B q_i(t) - \sum_{i=1}^Bq_i(t_0)<\infty.$$
This provides a contradiction with $\limsup_{t\to\infty} q_1(t) = \lambda +\varepsilon$, since the rate of decrease of $q_1(t)$ is at most~1.
\end{claimproof}
Claims~\ref{claim:inf} and~\ref{claim:sup} together imply Lemma~\ref{lem:q1}.
\end{proof} 

\begin{lemma}\label{lem:q2}
$q_2(t)\to 0$ as $t\to\infty$.
\end{lemma}
Based on the fact that $q_1(t)\to \lambda$ as $t\to\infty$, we now claim that
$q_2(t)\to 0$ as $t\to\infty$.
The high-level idea behind the claim is as follows. 
From the convergence of $q_1(t)$, we know that after a large enough time, $q_1(t)$ will always belong to a very small neighborhood of $\lambda$.
On the other hand, if $q_2(t)$ does not converge to 0, then it must  have a strictly positive limit point. 
In that case, since the rate of decrease of $q_2(t)$ is at most $q_2(t)$, it will be bounded away from 0 for a fixed amount of time infinitely often.
In the meantime, the rate at which busy servers become idle-on will be strictly less than the arrival rate of tasks.
This in turn, will cause $q_1(t)$ to increase substantially compared to the small neighborhood where it is supposed to lie, which leads to a contradiction.
\begin{proof}[Proof of Lemma~\ref{lem:q2}]
Lemma~\ref{lem:q1} implies that for any $M,\varepsilon>0$, there exists finite time $T(\varepsilon,M)$, such that $|q_1(t) - \lambda|\leq \varepsilon/M$ for all $t\geq T(\varepsilon,M).$
We will show that $$\limsup_{t\to\infty} q_2(t) = 0.$$ 
Suppose not, i.e., $q_2(T)>\varepsilon>0$ for some $T>T(\varepsilon,M).$
Since the rate of decrease of $q_2(t)$ is at most $q_2(t)$, it follows that $q_2(t)\geq 9\varepsilon/16$ for all $t\in [T,T+1/2]$, and hence 
\begin{equation}\label{eq:diff-l-q1-q2}
q_1(t)-q_2(t)\leq \lambda + \frac{\varepsilon}{M}-\frac{9\varepsilon}{16}\leq \lambda -\frac{\varepsilon}{2},
\end{equation}
for $M\geq 16$.
Due to Fact~\ref{fact:2},
$$q_1\Big(T+\frac{1}{2}\Big)-q_1(T)\geq \frac{\varepsilon}{2}\int_T^{T+\frac{1}{2}}\indn{\cX^c}(\qq(s),\dd(s))\dif s.$$
Since 
\begin{equation}\label{eq:tempbound}
q_1\Big(T+\frac{1}{2}\Big)-q_1(T)\leq \frac{2\varepsilon}{M},
\end{equation}
it follows that
\begin{equation}\label{eq:indicator-upper2}
\int_T^{T+\frac{1}{2}}\indn{\cX^c}(\qq(s),\dd(s))\dif s\leq \frac{4}{M}.
\end{equation} 
Also, Fact~\ref{fact:3} yields
$$q_1\Big(T+\frac{1}{2}\Big)-q_1(T)\geq \nu\int_T^{T+\frac{1}{2}}\delta_1(s)\dif s - \frac{4\nu(\nu+1)}{M}.$$
Again using~\eqref{eq:tempbound} it follows that
\begin{equation}\label{eq:delta-upper}
\nu\int_T^{T+\frac{1}{2}}\delta_1(s)\dif s\leq \frac{4\nu(\nu+1)+2\varepsilon}{M}\leq \frac{5\nu(\nu+1)}{M},
\end{equation}
for $\varepsilon$ sufficiently smaller than $\nu$.
We will now proceed to show that~\eqref{eq:delta-upper} yields a contradiction.
Notice that
\begin{align*}
\delta_1(t) &= \delta_1(T) + \int_T^t\lambda(1-p_0(\qq(s),\dd(s),\lambda))\ind{\delta_0(s)>0}\dif s - \nu \int_T^t\delta_1(s)\dif s\\
&\geq  \int_T^t (\lambda-q_1(s)+q_2(s))\indn{\cX^c}(\qq(s),\dd(s))\dif s \ind{\delta_0(s)>0}\dif s-2\nu\int_T^t\delta_1(s)\dif s.
\end{align*}
Using \eqref{eq:diff-l-q1-q2}, we obtain for all $t\in [T,T+1/2]$,
\begin{align*}
&\delta_1(t)\geq -2\nu\int_T^t\delta_1(s)\dif s + \frac{\varepsilon}{2}\int_T^t \indn{\cX}(\qq(s),\dd(s))\ind{\delta_0(s)>0}\dif s\\
&\geq -2\nu\int_T^t\delta_1(s)\dif s + (t-T)\frac{\varepsilon}{2} -\frac{\varepsilon}{2}\int_T^t\indn{\cX^c}(\qq(s),\dd(s))\dif s-\frac{\varepsilon}{2}\int_T^t\ind{u(s)=0,\delta_0(s)=0}\dif s\\
&\geq -2\nu\int_T^{T+\frac{1}{2}}\delta_1(s)\dif s-\frac{\varepsilon}{2}\int_T^{T+\frac{1}{2}}\indn{\cX^c}(\qq(s),\dd(s))\dif s  +(t-T)\frac{\varepsilon}{2}\\
&\hspace{7cm}-\frac{\varepsilon}{2}\int_T^{T+\frac{1}{2}}\ind{u(s)=0,\delta_0(s)=0}\dif s,
\end{align*}
and using~\eqref{eq:indicator-upper2} and~\eqref{eq:delta-upper}, it follows that
\begin{align}\label{eq:delta1contr}
\delta_1(t)&\geq -\frac{10\nu(\nu+1)}{M}-\frac{2\varepsilon}{M}+(t-T)\frac{\varepsilon}{2}-\frac{\varepsilon}{2}\int_T^{T+\frac{1}{2}}\ind{\delta_1(s)\geq (1-\lambda-\varepsilon/M)}\dif s.
\end{align}
Furthermore, observe that due to~\eqref{eq:delta-upper}, 
$$\int_T^{T+\frac{1}{2}}\ind{\delta_1(s)\geq (1-\lambda-\varepsilon/M)}\dif s \leq \frac{5(\nu+1)}{M(1-\lambda-\frac{\varepsilon}{M})}\leq \frac{10(\nu+1)}{M(1-\lambda)},$$
for $\varepsilon$ small enough, and thus,~\eqref{eq:delta1contr} yields
 \begin{align*}
\delta_1(t)&\geq -\frac{10\nu(\nu+1)}{M}-\frac{2\varepsilon}{M}+(t-T)\frac{\varepsilon}{2}-\frac{\varepsilon}{2}\int_T^{T+\frac{1}{2}}\ind{\delta_1(s)\geq (1-\lambda-\varepsilon/M)}\dif s\geq \frac{\varepsilon}{16},
\end{align*}
for all $t\in [T+1/4, T+1/2]$, for $M$ sufficiently large.
\end{proof}

Since $q_1(t) - q_2(t)\to \lambda$ and $q_2(t)\to 0$, as $t\to\infty$,
it follows from \eqref{eq:lb-p0} that $p_0(\mathbf{q}(t),\boldsymbol{\delta}(t),\lambda)\to 1$ as $t\to\infty.$
Also, an application of Gronwall's inequality to
$$\delta_1(t) = \delta_1(0)+\int_o^t\lambda(1-p_0(\mathbf{q}(s),\boldsymbol{\delta}(s),\lambda))\dif s- \int_0^t\delta_1(s)\nu\dif s,$$
yields $\delta_1(t)\to 0$ as $t\to\infty$.
Consequently, $\delta_0(t)\to 1-\lambda$ as $t\to\infty$. This completes the proof of Proposition~\ref{prop:glob-stab}.
\end{proof}

\begin{proof}[Proof of Proposition~\ref{thm:limit interchange}]
Note that the proof of the proposition follows from~\cite[Corollary 2]{BenB08}. The arguments are sketched briefly for completeness.

Observe that $\pi^{N}$ is defined on $E$, and $E$ is a compact set. 
Prohorov's theorem implies that $\pi^{N}$ is relatively compact, and hence, has a convergent subsequence. Let $\{\pi^{N_n}\}_{n\geq 1}$ be a convergent subsequence, with $\{N_n\}_{n\geq 1}\subseteq\N$, such that $\pi^{N_n}\dto\hat{\pi}$ as $n\to\infty$. 
We will show that $\hat{\pi}$ is unique and equals the measure~$\pi.$

Notice that if $(\qq^{N_n}(0),\dd^{N_n}(0))\sim\pi^{N_n}$, then  we know $(\qq^{N_n}(t),\dd^{N_n}(t))\sim\pi^{N_n}$ for all $t\geq 0$. 
Also, the process $(\qq^{N_n}(t),\dd^{N_n}(t))_{t\geq 0}$ converges weakly to $\{(\qq(t),\dd(t))\}_{t\geq 0}$, and $\pi^{N_n}\dto\hat{\pi}$ as $n\to\infty$. 
 Thus, $\hat{\pi}$ is an invariant distribution of the deterministic process $\{(\qq(t),\dd(t))\}_{t\geq 0}$.
 This in conjunction with the global stability in Proposition~\ref{prop:glob-stab} implies that $\hat{\pi}$ must be the fixed point of the fluid limit. 
 Since the latter fixed point is unique, we have shown the convergence of the stationary measure.
\end{proof}

\section{Conclusion}\label{sec:conclusion-sig}
Centralized queue-driven auto-scaling techniques do not cover
scenarios where load balancing algorithms immediately distribute
incoming tasks among parallel queues, as typically encountered
in large-scale data centers and cloud networks.
Motivated by these observations, we proposed a joint auto-scaling
and load balancing scheme, which does not require any global queue
length information or explicit knowledge of system parameters.
Fluid-limit results for a large-capacity regime show that the
proposed scheme achieves asymptotic optimality in terms of response
time performance as well as energy consumption.
At the same time, the proposed scheme operates in a distributed
fashion, and involves only a constant communication overhead per task,
ensuring scalability to massive numbers of servers.
This demonstrates that, rather remarkably, ideal response time
performance and minimal energy consumption can be simultaneously
achieved in large-scale distributed systems.

Extensive simulation experiments support the fluid-limit results,
and reveal only a slight trade-off between the mean waiting time
and energy wastage in finite-size systems.
In particular, we observe that suitably long but finite standby
periods yield near-minimal waiting time and energy consumption,
across a wide range of setup durations.
We expect that a non-trivial trade-off between response time
performance and (normalized) energy consumption arises at the
diffusion level, and exploring that conjecture would be
an interesting topic for further research.
It might be worth noting that in the present chapter, we have not taken the communication delay into consideration, and assumed that the message transfer is instantaneous.
This is a reasonable assumption when the communication delay is insignificant relative to the typical duration of the service period of a job.
When the communication delay is non-negligible, one might modify the TABS scheme where a task is discarded if it happens to land on an idle-off server.
In this modified scheme, the asymptotic fraction of lost tasks in steady state should be negligible, since the rate at which idle-on servers are turning off is precisely zero at the fixed point, and it would be useful to further examine the impact of communication delays.


%% file: energy2.tex
\begin{abstract}
We consider the model of a token-based joint auto-scaling and load balancing strategy, proposed in Chapter~\ref{chap:energy1}, which offers efficient scalable implementation, and  asymptotically optimal steady-state delay performance and energy consumption as the number of servers $N\to\infty$. In Chapter~\ref{chap:energy1}, the asymptotic results were obtained \emph{under the assumption that the queues have fixed-size finite buffers}, and therefore the fundamental question of stability with infinite buffers was left open. In this chapter, we address this fundamental stability question. The system stability under the usual subcritical load assumption is not automatic. Moreover, the stability may \emph{not} even hold for all $N$. The key challenge stems from the fact that the process \emph{lacks monotonicity}, which has been the powerful primary tool for establishing stability in load balancing models. We develop a novel method to prove that the subcritically loaded system is stable for \emph{large enough}~$N$, and establish convergence of steady-state distributions to the optimal one, as $N \to \infty$. 
The method advances the state of the art with an induction-based idea that exploits a weak monotonicity property of the model.
This novel method is of independent interest and may have broader applicability.
\end{abstract}

\section{Introduction}
In this chapter we return to the TABS scheme introduced in Chapter~\ref{chap:energy1}.
There we left open a fundamental question: Is the system with a given number $N$ of servers stable under the TABS scheme?
The analysis in Chapter~\ref{chap:energy1} bypasses the issue of stability by assuming that each server in the system has a finite buffer capacity.
Thus, it remains an important open challenge to understand the stability property of the TABS scheme without the finite-buffer restriction. \\

In this chapter we address these stability issues and examine the asymptotic behavior of the system as $N$ becomes large.
Analyzing the stability of the TABS scheme in the infinite-buffer scenario poses a significant challenge, because the stability of the finite-$N$ system, i.e., the system with finite number of $N$ servers under the usual subcritical load assumption is not automatic.
In fact, even under subcritical load, the system may \emph{not} be stable for all $N$ (see Remark~\ref{rem:instability} for details).
Our first main result is that for any fixed subcritical load, the system is stable for \emph{large enough} $N$.
Further, using this large-$N$ stability result in combination with mean-field analysis, we establish convergence of the sequence of steady-state distributions as~$N\to\infty$.

The key challenge in showing large-$N$ stability for systems under the TABS scheme stems from the fact that the occupancy state process lacks monotonicity.
It is well-known that monotonicity is a powerful primary tool for establishing stability of load balancing models~\cite{Stolyar15,Stolyar17, VDK96, BLP12}.
In fact, process monotonicity is used extensively not only for stability analysis and not only in queueing literature -- for example,
many interacting-particle systems' results rely crucially on monotonicity; see e.g.~\cite{L85}. The lack of monotonicity immediately complicates the situation, as for example in ~\cite{FS17, SS17}.
Specifically, when the service time distribution is general, it is the lack of monotonicity that has left open the stability questions for the power-of-d scheme when the system load $\lambda>1/4$~\cite{BLP12}, and for the JIQ scheme when $\lambda>1/2$~\cite{FS17}.
We develop a novel method for proving \emph{large-$N$ stability} for subcritically loaded systems, and use that  to establish the convergence of the sequence of steady-state distributions as $N\to\infty$.
Our method uses an induction-based idea, and relies on a ``weak monotonicity'' property of the model, as further detailed below.
To the best of our knowledge, this is the first time both the \emph{traditional fluid limit} (in the sense of a large starting state) and the \emph{mean-field fluid limit} (when the number of servers grows large) are used in an intricate manner to obtain large-$N$ stability results.

To establish the large-$N$ stability, we actually prove a stronger statement.
We consider an artificial system, where some of the queues are infinite at all times. 
Then, loosely speaking, we prove that the following holds for all sufficiently large $N$:
\emph{If the system with $N$ servers contains $k$ servers with infinite queue lengths, $0\leq k\leq N$, then (i)~The subsystem consisting of the remaining (i.e., finite) queues is stable, and 
(ii)~When this subsystem is in steady state, the average rate at which tasks join the infinite queues is strictly smaller than that at which tasks depart from them.}
Note that the case $k=0$ corresponds to the desired stability result.

The use of backward induction in $k$ facilitates proving the above statement.
For a fixed $N$, first we introduce the notion of a fluid sample path (FSP) for systems where some queues might be infinite.
The base case of the backward induction is when $k=N$, and assuming the statement for $k$, we show that it holds for $k-1$.
We use the classical fluid stability argument (as in~\cite{RS92, S95, D99}) in order to establish stability for the system where the number of infinite queues is $k-1$.
As mentioned above, here the notion of the traditional FSP is needed to be suitably extended to fit to the systems where some servers have infinite queue lengths.
Loosely speaking, for the fluid-stability, the `large queues' behave as `infinite queues' for which the induction statement provides us with the drift estimates.
Also, to calculate the drift of a queue in the fluid limit for fixed but \emph{large enough $N$}, we use the mean-field analysis.
A more detailed heuristic roadmap of the above proof argument is presented in Subsection~\ref{ssec:roadmap}.
This technique is of independent interest, and potentially has a much broader applicability in proving large-$N$ stability for non-monotone systems, where the state-of-the-art results have remained scarce so far.
\\

\noindent
\textbf{Organization of the chapter.}
The rest of the chapter is organized as follows.
In Section~\ref{sec:model-inf} we present a detailed model description,  state the main results, and discuss their ramifications along with discussions of several proof heuristics.
The full proof of the main results is deferred till Section~\ref{sec:proofs-sasha}. 
Section~\ref{sec:ind} introduces an inductive approach to prove the large-$N$ stability result.
We present the proof of the large-scale system (when $N\to\infty$) using mean-field analysis in Section~\ref{sec:mf}.
Finally, we make a few brief concluding remarks in Section~\ref{sec:con}.

\section{Main results}\label{sec:model-inf}

Recall the TABS scheme from Section~\ref{sec:model-sig} in Chapter~\ref{chap:energy1},
with the consideration that the buffer capacity at each server is $B = \infty$.
Also, assume the total arrival rate for the $N$-th system is $\lambda N$ for some fixed $\lambda\in (0,1)$ (that does not vary over time).
It is easy to see that, for any fixed $N$, this process is an irreducible countable-state Markov process. 
Therefore, its positive recurrence, which we refer to as \emph{stability}, is equivalent to ergodicity and to the existence of unique stationary distribution. 
Further, let  $U^N(t)$ denote the number of idle-on servers at time $t$.
We will focus on an asymptotic analysis, where the task arrival rate and the number of servers grow large in proportion.

For the description of the occupancy process we refer to Section~\ref{sec:model-sig} in Chapter~\ref{chap:energy1}.
We emphasize that in this chapter we will use the term \emph{mean-field fluid scaling}, corresponding to the term fluid-scaling in Chapter~\ref{chap:energy1}.
Thus, \emph{mean-field fluid-scaled} quantities are denoted by the respective small letters, viz.~$q_i^{N}(t):=Q_i^{N}(t)/N$, $\delta_0^N(t) = \Delta_0^N(t)/N$, $\delta_1^N(t) = \Delta_1^N(t)/N$, and $u^N(t):= U^N(t)/N$, and
$$
E = \Big\{(\bld{q},\dd)\in [0,1]^\infty:  q_i\geq q_{i+1},\ \forall i, \  \delta_0+\delta_1+ q_1\leq 1 \Big\},
$$
denote the space of all mean-field fluid-scaled occupancy states,
so that the process $(\mathbf{q}^N(t),\bld{\delta}^N(t))$ takes value in $E$ for all $t$.
Endow $E$ with the product topology, and the Borel $\sigma$-algebra $\mathcal{E}$, generated by the open sets of $E$.
Notation for the \emph{conventional fluid-scaled} occupancy states for a fixed~$N$ will be introduced later in Subsection~\ref{ssec:fluid-sasha}.
By the symbol `$\pto$' we denote convergence in probability for real-valued random variables.\\

We now present our first main result: 
\begin{theorem}\label{th:stab}
For any fixed $\mu$, $\nu>0,$ and $\lambda<1$, the system with $N$ servers under the TABS scheme is stable (positive recurrent) for large enough $N$. 
\end{theorem}
Theorem~\ref{th:stab} is proved in Section~\ref{sec:proofs-sasha}.

\begin{remark}\label{rem:instability}\normalfont
It is worthwhile to mention that the `large-$N$' stability in Theorem~\ref{th:stab} is the best one can hope for. 
In fact, for fixed $N$ and $\lambda$, there are values of the parameters $\mu$ and $\nu$ such that the system under the TABS scheme may not be stable.
To elaborate further on this point, consider a system with 2 servers A and B, and $1/2<\lambda<1$.
Let server A start with a large queue, while the initial queue length at server B is small.
In that case, observe that every time the queue length at server B hits 0, with positive probability, it turns idle-off before the next arrival epoch.
Once server B is idle-off, the arrival rate into server A becomes $2\lambda>1$. 
Thus, before server B turns idle-on again, the expected number of tasks that join server A is given by at least $2\lambda/\nu$, while the expected number of departures is $1/\nu$.
Thus the queue length at server A increases by $(2\lambda-1)/\nu$,
which can be very large if $\nu$ is small.
Further note that once server B becomes busy again, both servers receive an arrival rate $\lambda<1$, and hence it is more likely that server B will empty out again, repeating the above scenario.
The situation becomes better as $N$ increases. 
Indeed for large $N$, if `too many' servers are idle-off and `too many' tasks do not find an idle queue to join, the system starts producing servers in setup mode fast enough, and as a result, more and more servers start becoming busy.
The above heuristic has been illustrated in Figures~\ref{fig:instability1}--\ref{fig:instability3} with examples of three scenarios with small, moderate, and large values of $N$, respectively.
\begin{figure}
\begin{center}
\includegraphics[width=90mm]{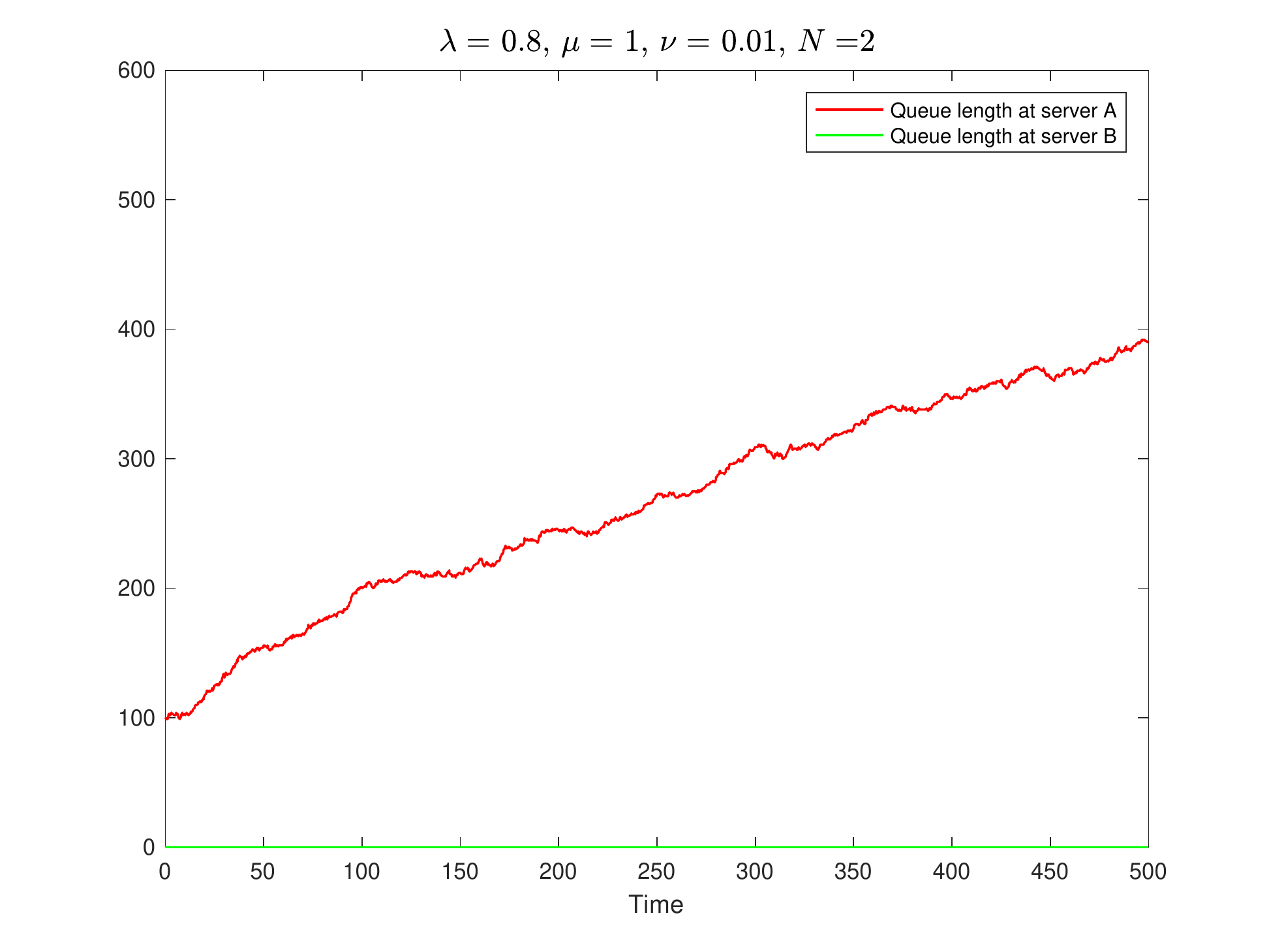}
\end{center}
\caption{Illustration of instability of the TABS scheme for $N=2$ via sample paths of the queue length process. }
\label{fig:instability1}
\end{figure}
\begin{figure}
\begin{center}
\includegraphics[width=90mm]{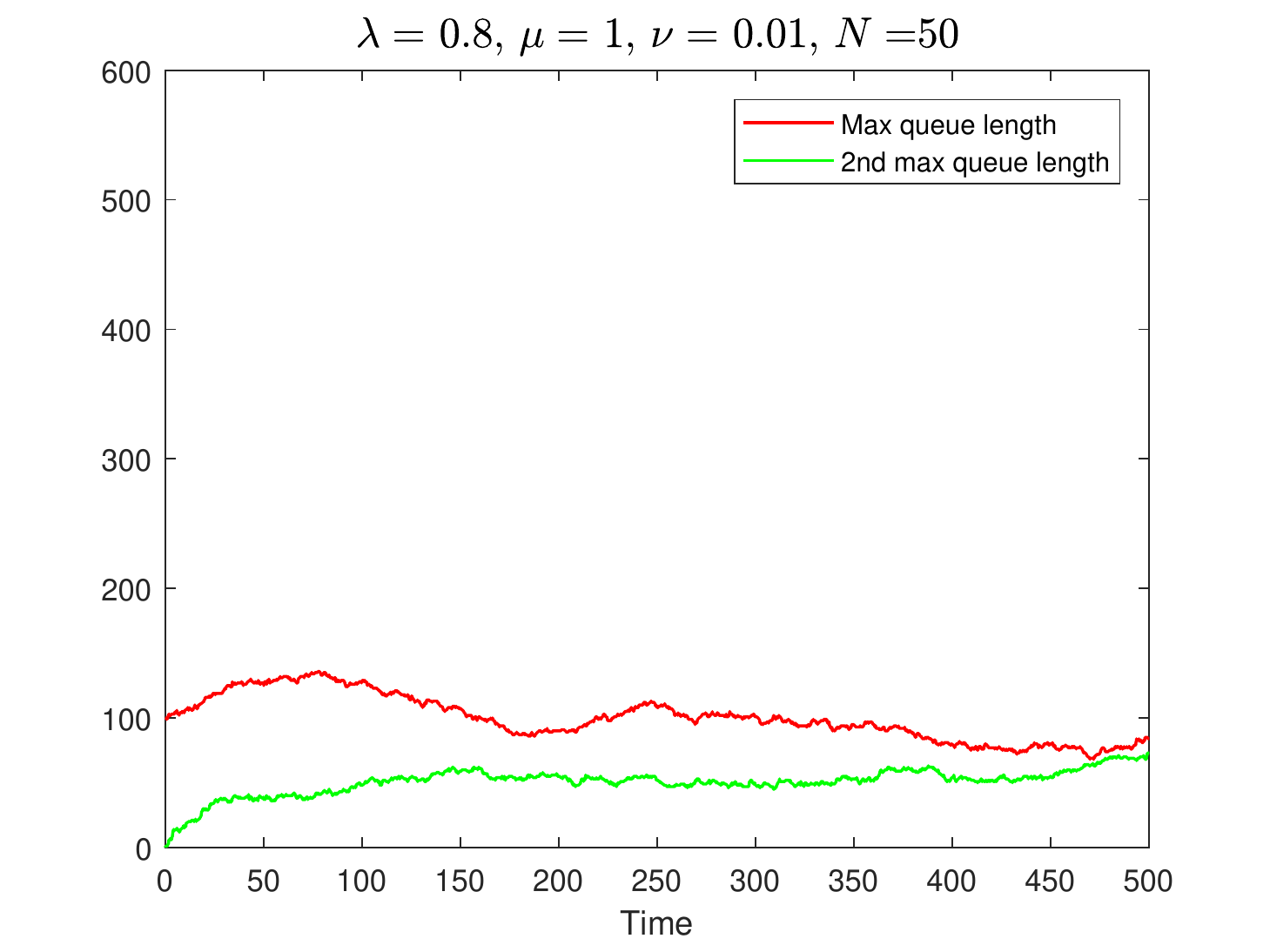}
\end{center}
\caption{Sample paths of the largest and second largest queue length processes in an intermediate system ($N=50$) for the same parameter choices.}
\label{fig:instability2}
\end{figure}
\begin{figure}
\begin{center}
\includegraphics[width=90mm]{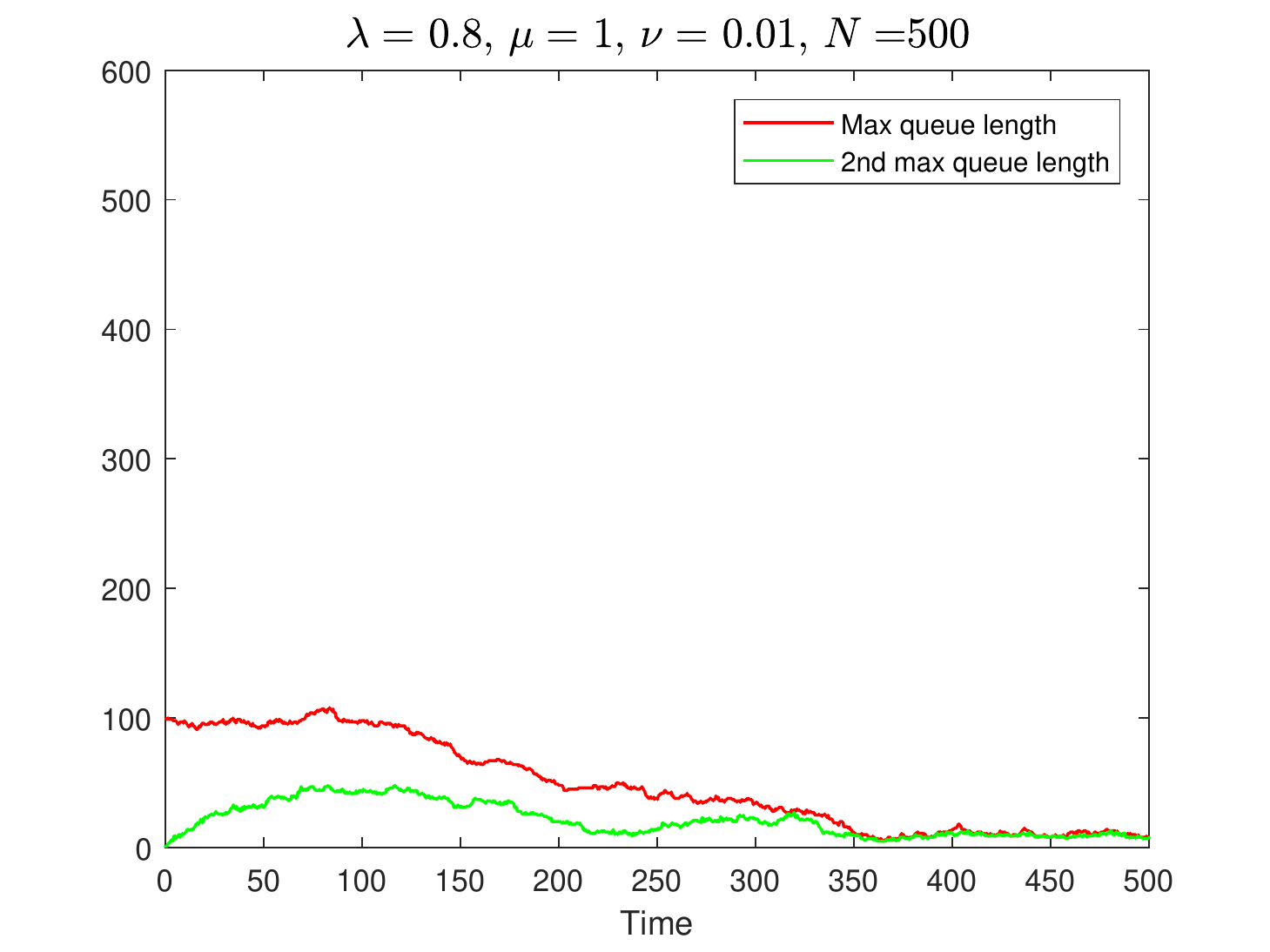}
\end{center}
\caption{The system becomes stable for a large enough number of servers ($N=500$).}
\label{fig:instability3}
\end{figure}
\end{remark}

In the next theorem we will identify the limit of the sequence of stationary distributions of the occupancy processes as $N\to\infty$.
In particular, we will establish that under sub-critical load, for any fixed $\mu$, $\nu>0$, the steady-state occupancy process converges weakly to the unique fixed point.
(For the finite-buffer scenario this was proved in Proposition~\ref{thm:limit interchange} in Chapter~\ref{chap:energy1}.)
Denote by $\qq^N(\infty)$ and $\dd^N(\infty)$ the random values of $\qq^N(t)$ and $\dd^N(t)$ in the steady state, respectively.
\begin{theorem}\label{th:steady-lim}
For any fixed $\mu$, $\nu>0$, and $\lambda<1$, the sequence of steady states
$(\qq^N(\infty), \dd^N(\infty))$ converges weakly to the fixed point $(\qq^\star, \dd^\star)$  as  $N\to\infty$,
where 
\[\delta_0^\star = 1-\lambda\quad\delta_1^\star = 0\quad q_1^\star = \lambda,\quad q_i^\star =0 \quad\mbox{for all}\quad i\geq 2.\]
\end{theorem}
Note that the fixed point $(\qq^\star, \dd^\star)$ is such that the  probability of wait vanishes as $N\to\infty$ and the asymptotic fraction of active servers is the minimum required for stability, and in this sense, the fixed point is optimal.
Thus, Theorem~\ref{th:steady-lim} implies that the TABS scheme provides fluid-level optimality for large-scale systems in terms of delay performance and resource utilization,
while involving only $O(1)$ communication overhead per task.


\section{Proofs of the main results}\label{sec:proofs-sasha}

In Subsection~\ref{ssec:fluid-sasha} we introduce 
the notion of conventional fluid scaling (when the number of servers is fixed) and fluid sample paths (FSP), and state Proposition~\ref{prop:largeN} that implies Theorem~\ref{th:stab} as an immediate corollary.
Subsection~\ref{ssec:large-scale} contains two key results for a sequence of systems with increasing system size, i.e., the number of servers $N\to\infty$, and proves Theorem~\ref{th:steady-lim}.

\subsection{Conventional fluid limit for a system with fixed N}\label{ssec:fluid-sasha}
In this subsection first we will introduce a notion of fluid sample path (FSP) for finite-$N$ systems where some of the queue lengths are infinite.
We emphasize that this is the \emph{conventional fluid limit}, in the sense that the number of servers is fixed, but the time and the queue length at each server are scaled by some parameter that goes to infinity.

Loosely speaking, conventional fluid limits are usually defined as follows: 
For a fixed $N$, consider a sequence of systems with increasing initial norm (total queue length) $R$ say. 
Now scale the queue length process at each server and the time by $R$.
Then any weak limit of this sequence of (space and time) scaled processes is called an FSP.
Observe that this definition is inherently not fit if the system has some servers whose initial queue length is infinite.
Thus we introduce a suitable notion of FSP that does not require the scaled norm of the initial state to be 1.
We now introduce a rigorous notion of FSP for systems with some of the queues being infinite.\\

\noindent
\textbf{Fluid limit of a system with some of the queues being infinite.} 
Consider a system of $N$ servers with indices in $\cN$ (say), among which $k$ servers with indices in $\cK\subseteq\cN$ have infinite queue lengths.
Now consider any sequence of systems indexed by $R$ such that $\sum_{i\in\cN\setminus\cK}X_i^{N,R}(0) <\infty$, and 
\begin{equation}\label{eq:conventional}
x^{N,R}_i(t) := \frac{X_i^{N,R}(Rt)}{R},\quad i\in\cN\setminus\cK
\end{equation}
be the corresponding scaled processes.
For fixed $N$, the scaling in~\eqref{eq:conventional} will henceforth be called the \emph{conventional fluid-scaled} queue length process.
Also, for the $R$-th system, let $A_i^{N,R}(t)$ and $D_i^{N,R}(t)$ denote the cumulative number of arrivals to and departures from server $i$ with $a_i^{N,R}(t):=A_i^{N,R}(Rt)/R$ and $d_i^{N,R}(t):=D_i^{N,R}(Rt)/R$ being the corresponding fluid-scaled processes, $i\in\cN$.
We will often omit the superscript $N$ when it is clear from the context.

Now for any fixed $N$, suppose the (conventional fluid-scaled) initial states converge, i.e., $x^{R}(0) \to x(0)$, for some fixed
$x(0)$ such that $0 \le \sum_{i\in\cN\setminus\cK} x_i(0) < \infty$ and $x_i(0) = \infty$ for $i\in\cK$.
Then a set of uniformly Lipschitz continuous functions $(x_i(t), a_i(t), d_i(t))_{i\in\cN}$ on the time interval $[0,T]$ (where $T$ is possibly infinite) with the convention $x_i(\cdot)\equiv \infty$ for all $i\in \cK$, is called a \emph{fluid sample path} (FSP) starting from $\xx(0)$, if for any subsequence of $\{R\}$ there exists a further subsequence (which we still denote by $\{R\}$) such that with probability~1, along that subsequence the following convergences hold:
\begin{enumerate}[{\normalfont (i)}]
\item For all $i\in \cN$, $a_i^{R}(\cdot)\to a_i(\cdot)$ and $d_i^{R}(\cdot)\to d_i(\cdot)$, uniformly on compact sets.
\item For $i\in\cN\setminus\cK$, $x_i^{R}(\cdot)\to x_i(\cdot)$ uniformly on compact sets.
\end{enumerate}
Note that the above definition is equivalent to convergence in probability to the unique FSP.
For any FSP almost all points (with respect to the Lebesgue measure) are \emph{regular}, i.e., for all $i\in \cN\setminus\cK$, $x_i(t)$ has proper left and right derivatives with respect to $t$, and for all such regular points, \[x_i'(t) = a_i'(t)-d_i'(t).\]

\noindent
\textbf{Infinite queues as part of an FSP.} The arrival and departure functions $a_i(t)$ and $d_i(t)$ are well-defined for each queue, including infinite queues. Of course, the derivative $x'_i(t)$ for an infinite queue makes no direct sense (because an infinite queue remains infinite at all times). However, we adopt a convention that $x'_i(t)=a'_i(t)-d'_i(t)$, for all queues, including the infinite ones. For an FSP, $x'_i(t)$ is sometimes referred to as a ``drift'' of (finite or infinite) queue $i$ at time $t$.\\

We are now in a position to state the key result that establishes the large-$N$ stability of the TABS scheme.
\begin{proposition}\label{prop:largeN}
The following holds for all sufficiently large $N$. 
For each $0\leq k\leq N$, consider a system where $k$ servers with indices in $\cK$ have infinite queues, and the remaining $N-k$ queues are finite.
Then, for each $j=1,2,\ldots,N$, there exists $\varepsilon(j) >0$, such that the following properties hold {\normalfont(}$\varepsilon(j)$ and other constants specified below, also depend on $N${\normalfont)}.
\begin{enumerate}[{\normalfont (1)}]
\item For any $\xx(0)$ such that $0 \le \sum_{i\in\cN\setminus\cK} x_i(0) < \infty$ and $x_i(0) = \infty$ for $i\in\cK$,
there exists $T(k,\xx(0))<\infty$ and a unique FSP on the interval $[0,T(k,\xx(0))]$, which has the following properties:
	\begin{enumerate}[{\normalfont (i)}]
	\item If at a regular point $t$, $\cM(t):=\{i\in \cN: x_i(t)>0\}$ with $|\cM(t)|=m>k$, then $x_i'(t) = -\varepsilon(m)$ for all $i\in \cM(t)$.
	\item For any $i\in\cN\setminus\cK$, if $x_i(t_0)=0$ for some $t_0$, then $x_i(t)=0$ for all $t\geq t_0$.
	\item $T(k,\xx(0))=\inf\ \{t:x_i(t) = 0\mbox{ for all }i\in\cN\setminus\cK\}$.
	\end{enumerate}
\item The subsystem with $N-k$ finite queues is stable.
\item When the subsystem with $N-k$ finite queues is in steady state, the average arrival rate into each of the $k$ servers  having infinite queue lengths is at most $1-\varepsilon(k)$.
\item For any $x(0)$ such that $0 \le \sum_{i\in\cN\setminus\cK} x_i(0) < \infty$ and $x_i(0) = \infty$ for $i\in\cK$,
there exists a unique FSP on the entire interval $[0,\infty)$. 
In $[0,T(k,x(0))]$, it is as described in Statement 1.
Starting from $T(k,x(0))$, all queues in $\cN\setminus\cK$
stay at $0$ and all infinite queues have drift at most $-\varepsilon(k)$.
\end{enumerate}
\end{proposition}
Although Part 2 follows from Part 1, and Part 4 is stronger than Part 1, the statement of Proposition~\ref{prop:largeN} is arranged as it is to facilitate its proof, as we will see in Section~\ref{sec:ind} in detail.
\begin{proof}[Proof of Theorem~\ref{th:stab}]
Note that Theorem~\ref{th:stab} is a special case of Proposition~\ref{prop:largeN} when $k=0$.
\end{proof}

\subsection{Large-scale asymptotics: auxiliary results}
\label{ssec:large-scale}
In this subsection we will state two crucial lemmas that describe asymptotic properties of a sequence of systems as the number of servers $N\to\infty$, \emph{if stability is given}.
Their proofs involve mean-field fluid scaling and limits.
\begin{lemma}\label{lem:expo-bound}
There exist $\varepsilon_1>0$ and $C_q=C_q(\varepsilon_1)>0$, such that the following holds.
Consider any sequence of systems with $N\to\infty$ and $k=k(N)$ infinite queues such that $k(N)/N \to \kappa\in [0,1]$, and assume that each of these systems is stable.
Then for all sufficiently large $N$,
\[\Pro{q_1^N(\infty)<\varepsilon_1}\leq \e^{-C_qN}.\]
\end{lemma}
\begin{lemma}\label{lem:steady-concentration}
Consider any sequence of systems with $N\to\infty$ and $k=k(N)$ infinite queues such that $k(N)/N \to \kappa\in [0,1]$, and assume that each of these systems is stable.
The following statements hold:
\begin{enumerate}[{\normalfont (1)}]
\item If $\kappa \geq 1-\lambda$, then $q_1^N(\infty)\pto 1$ as $N\to\infty$. 
\item If $\kappa < 1-\lambda$, then 
the weak limit of $(\qq^N(\infty),\dd^N(\infty))$ is concentrated at the unique equilibrium point $(\qq^\star(\kappa), \dd^\star(\kappa))$, such that
\begin{align*}
q_1^\star(\kappa) &= \kappa+\lambda,\quad q_2^\star(\kappa) = \kappa,\\
\delta_0^\star(\kappa) &= 1-\lambda-\kappa,\quad\delta_1^\star(\kappa) = 0.
\end{align*}
Consequently,
\begin{equation}\label{eq:foundbusy}
\lim_{N\to\infty}\Pro{Q_1^N(\infty)+\Delta_0^N(\infty)+\Delta_1^N(\infty)=N}=0.
\end{equation}
\end{enumerate}
\end{lemma} 
Lemmas~\ref{lem:expo-bound} and \ref{lem:steady-concentration} are proved in Section~\ref{sec:mf}. 
These results will be used to derive necessary large-$N$ bounds on the expected arrival rate into each of the servers having infinite queue lengths when the system is in steady state.
\begin{remark}\normalfont
It is also worthwhile to note that Lemmas~\ref{lem:expo-bound} and~\ref{lem:steady-concentration} can be thought of as a \emph{weak monotonicity} property of the TABS scheme as mentioned earlier.
Loosely speaking, the weak monotonicity requires that no matter where the system starts, in some fixed time the system arrives at a state with a certain fraction of busy servers.
The purpose of Lemmas~\ref{lem:expo-bound} and~\ref{lem:steady-concentration} is to bound \emph{under the assumption of stability} the expected rate at which tasks arrive to the infinite queues when the subsystem containing the finite queues is in steady state: 
\begin{enumerate}[{\normalfont (i)}]
\item Lemma~\ref{lem:steady-concentration} guarantees high probability bounds on the total number of busy servers, so that with probability tending to 1 as $N\to\infty$, the fraction of busy servers in the whole system is at least $\lambda$ in steady state.
\item However, since the arrival rate is $\lambda N$, when the system has few busy servers (even with an asymptotically vanishing probability), the arrival rate to the infinite servers can become $\Theta(N)$.
Thus we need the exponential bound stated in Lemma~\ref{lem:expo-bound} in order to obtain a bound on the expected rate of arrivals to the infinite queues.
\end{enumerate}
In Subsection~\ref{ssec:part3} we will see that as a consequence of Lemmas~\ref{lem:expo-bound} and~\ref{lem:steady-concentration}, we obtain that for large enough $N$, under the assumption of stability, the steady-state rate at which tasks join an infinite queue is strictly less than 1, and the drift of the infinite queues as defined in Subsection~\ref{ssec:fluid-sasha} becomes strictly negative.
This fact will be used in the proof of Proposition~\ref{prop:largeN}.
\end{remark}

\begin{proof}[Proof of Theorem~\ref{th:steady-lim}]
Note that given the large-$N$ stability property proved in Proposition~\ref{prop:largeN} for $k(N)=0$, and the convergence of 
stationary distributions under the assumption of stability in Lemma~\ref{lem:steady-concentration}, the proof of Theorem~\ref{th:steady-lim} is immediate.
\end{proof}


\section{Proof of Proposition 8.3.1: An inductive approach}\label{sec:ind}
Throughout this section we will prove Proposition~\ref{prop:largeN}.
The proof consists of several steps and uses both a conventional fluid limit and a mean-field fluid scaling and limit in an intricate fashion.
Below we first provide a roadmap of the whole proof argument.

\subsection{Proof idea and the roadmap}\label{ssec:roadmap}
The key idea for the proof of Proposition~\ref{prop:largeN} is to use backward induction in $k$, starting from the base case $k=N$. 
For $k=N$, all the queues are infinite.
In that case, Parts (1) and (2) are vacuously satisfied with the convention $T(N,\xx(0))=0$.
Further observe that the TABS scheme does not differentiate between two large queues (in fact, any two non-empty queues).
Thus, when all queues are infinite, since all servers are always busy, each arriving task is assigned uniformly at random, and each server has an arrival rate $\lambda$ and a departure rate 1.
Thus, it is immediate that the drift of each server is $-(1-\lambda)<0$, and thus, $\varepsilon(N)=1-\lambda$.
This proves (3), and then (4) follows as well.

Now, we discuss the ideas to establish the backward induction step, i.e., assume that Parts (1)--(4) hold for $k\geq k(N)+1$ for some $k(N)\in \{0,1,\ldots, N-1\}$ and verify that the statements hold for $k=k(N)$.
Rigorous proofs to verify Parts (1)--(4) for $k=k(N)$ are presented  in Subsections~\ref{ssec:part1}--\ref{ssec:part4}.
We begin by providing a roadmap of these four subsections.

\paragraph{Part (1).} Recall that we denote by $\cK$ the indices of the servers having infinite queue lengths, and by $\cN$ the set of all server indices.
Denote by $x_{(i)}$ the $i$-th largest component of $\xx$ (ties are broken arbitrarily).
Then for any $\xx$ with $m\in  \{0,1,\ldots, N-1\}$ infinite components, define 
\begin{equation}\label{eq:Tkx}
T(m,\xx) := \frac{x_{(N)}}{\varepsilon(N)}+\sum_{i=1}^{N-m-1}\frac{x_{(N-i)}  -  x_{(N-i+1)}}{\varepsilon(N-i)}
\end{equation}
with the convention that $T(N,\xx)=0$ if all components of $\xx$ are infinite. 
For $k(N)\in \{0, 1,\ldots, N-1\}$, Part (1) is proved with the choice of $T(k, \xx(0))$ as given by~\eqref{eq:Tkx}. 
Indeed, recall that we are at the backward induction step where there are $k(N)$ infinite queues, and we also know from the hypothesis that Parts (1)--(4) hold if there are $k(N)+1$ or larger infinite queues in the system.
Loosely speaking, the idea is that as long as a conventional fluid-scaled queue length $x_j(t)$ at some server $j\in \cN\setminus\cK$ is positive, it can be coupled with a system where the queue length at server $j$ is infinite.
Thus, as long as there is at least one server $j\in \cN\setminus\cK$ with $x_j(t)>0$, the system can be `treated' as a system with at least $k(N)+1$ infinite queues, in which case Part~(4) of the backward induction hypothesis furnishes the drift of each positive component of the FSP (in turn, which is equal to the drift of each infinite queue for the corresponding system).

Now to explain the choice of $T(m,\xx)$ in~\eqref{eq:Tkx}, observe that when all the components of the $N$-dimensional FSP are strictly positive, each component has a negative drift of $-\varepsilon(N)$.
Thus, $x_{(N)}/\varepsilon(N)$ is the time when at least one component of the $N$-dimensional FSP hits 0.
From this time point onwards, each positive component has a  drift of $-\varepsilon(N-1)$, and thus,  $x_{(N)}/\varepsilon(N) + (x_{(N-1)} - x_{(N)})/\varepsilon(N-1)$ is the time when two components hit 0.
Proceeding this way, one can see that at time $T(m, \xx(0))$ all finite positive components of the FSP hit 0.
The above argument is formalized in Subsection~\ref{ssec:part1}.

\paragraph{Part (1) $\boldsymbol{\implies}$ Part (2).}
To prove Part 2, we will use the fluid limit technique of proving stochastic stability as in~\cite{RS92, S95, D99}, see for example \cite[Theorem 4.2]{D99} or \cite[Theorem 7.2]{S95} for a rigorous statement.
Here we need to show that the sum of the non-infinite queues (of an FSP) drains to~0. This is true, because by Part~(1) each positive non-infinite queue will have negative drift. The formal proof is in Subsection~\ref{ssec:part2}.

\paragraph{Part (2) + Lemmas~\ref{lem:expo-bound} and~\ref{lem:steady-concentration} $\boldsymbol{\implies}$ Part (3).}
Note that in the proofs of Parts (1) and (2) we have only used the backward induction hypothesis, and have not imposed any restriction on the value of $N$.
This is the only part where in the proof we use the large-scale asymptotics, in particular, Lemmas~\ref{lem:expo-bound} and~\ref{lem:steady-concentration}.
For that reason, in the statement of Proposition~\ref{prop:largeN} we use ``large-enough $N$''.
The idea here is to use a proof by contradiction.
Suppose Part (3) does not hold for infinitely many values of $N$. 
In that case, it can be argued that there exist a subsequence $\{N\}$ and \emph{some} sequence $\{k(N)\}$ with $k(N)\in \{0,1,\ldots, N-1\}$, such that when the subsystem consisting of $N-k(N)$ finite queues is in the steady state, the average arrival rate into each of the $k(N)$ servers having infinite queue lengths is at least~1, along the subsequence.
Loosely speaking, in that case, Lemmas~\ref{lem:expo-bound} and~\ref{lem:steady-concentration} together imply that for large enough $N$, there are `enough' busy servers, so that the rate of arrival to each infinite queue is strictly smaller than~1, which leads to a contradiction.
Note that we can apply Lemmas~\ref{lem:expo-bound} and~\ref{lem:steady-concentration} here, because Part~(2) ensures the required stability. The rigorous proof is in Subsection~\ref{ssec:part3}.

\paragraph{Parts (2), (3) + Time-scale separation $\boldsymbol{\implies}$ Part (4).}
We assume that Parts (1) -- (3) hold for $k \in \{k(N), k(N)+1, \ldots, N\}$, and we will verify Part (4) for $k=k(N)$.
Observe that it only remains to prove convergence to the FSP on the (scaled) time interval $[T(k,\xx(0)),\infty]$.
For this, observe that it is enough to consider the sequence of systems for which $\xx^R(0)\to\xx(0)$ where $x_i(0)=0$ for all $i\in\cN\setminus\cK$. 
In particular, all that remains to be shown is that the drift of each infinite queue is indeed~$-\varepsilon(k)$.
Recall the conventional fluid scaling and FSP from Subsection~\ref{ssec:fluid-sasha}, and let $R$ be the scaling parameter.
The proof consists of two main parts:
\begin{enumerate}[{\normalfont (i)}]
\item Let us fix any state $z$ of the \emph{unscaled} process. If the sequence of systems is such that $\xx^R(0)\to\xx(0)$ where $x_i(0)=0$ for all $i\in\cN\setminus\cK$, then due to Part (2), for the subsystem consisting of finite queues, the (scaled) hitting time to the (unscaled) state $z$ converges in probability to 0.
Also, since this subsystem is positive recurrent (due to Part (2)), starting from a fixed (unscaled) state~$z$, its expected (unscaled) return time to the state~$z$ is $O(1)$.
This will allow us to split the (unscaled) time line into i.i.d.~renewal cycles of finite expected lengths.
In addition, this also shows that in the scaled time the subsystem of finite queues evolves on a faster time scale and achieves `instantaneous stationarity'.
\item From the above observation we can claim that the number of arrivals to any specific infinite queue can be written as a sum of arrivals in the above defined i.i.d.~renewal cycles.
Using the strong law of large numbers (SLLN) we can then show that in the limit $R\to \infty$, the instantaneous rate of arrival to an specific infinite queue is given by the \emph{average arrival rate} when the subsystem with $N-k$ finite queues is in steady state.
Therefore, Part (3) completes the verification of Part (4).
\end{enumerate}
The above argument is rigorously carried out in Subsection~\ref{ssec:part4}.

\subsection{Coupling with infinite queues to verify Part (1)}\label{ssec:part1}
To prove Part (1), fix any $\xx(0)$ such that $0 \le \sum_{i\in\cN\setminus\cK} x_i(0) < \infty$ and $x_i(0) = \infty$ for $i\in\cK$.
Let $\cK_1\subseteq \cN\setminus\cK$ be the set of server-indices~$i$, such that $x_i(0) > 0$.
We will first show that when $\sum_{i\in\cN\setminus\cK} x_i(0) > 0$ with $|\cM(t)|=m\geq k(N)+1$, then it has a negative drift $-\varepsilon(m)$ for all $i\in\cM(t)$, thus proving Part (1.i). 
Since $\varepsilon(m)$'s are positive, this will then also imply Part (1.ii).
Now assume $\sum_{i\in \cN\setminus\cK}x_i(0)>0$. 
In that case we have that $|\cK_1|=:k_1>0$.
Now consider the sequence of processes $(x_i^R(\cdot), a_i^R(\cdot), d_i^R(\cdot))_{i\in\cN}$ along any subsequence $\{R\}$.
Define the stopping time
$$T^R:=\inf\big\{t: X_i^R(t) = 0 \mbox{ for some }i\in\cK_1\big\},$$
and $\tau^R= T^R/R$.
In the time interval $[0,T^R]$, we will couple this system with a system, let us label it $\Pi$, with $k+k_1$ infinite queues.
Let $(\bx_i^R(\cdot), \ba_i^R(\cdot), \bd_i^R(\cdot))_{i\in\cN}$ be the queue length, arrival, and departure processes corresponding to the system $\Pi$, and assume that $\bx_i^R(0)$ is infinite for $i\in \cK\cup \cK_1$. 
Now couple each arrival to and departure from the $i$-th server in both systems, $i\in\cN$.
Since the scheme does not distinguish among servers with positive queue lengths, observe that up to time $T^R$ both systems evolve according to their own statistical laws.
Also, up to time $T^R$, the queue length processes at the servers in $\cN\setminus (\cK\cup\cK_1)$ in both systems are identical.
Thus, in the (scaled) time interval $[0,\tau^R]$, $a_i^R \equiv \ba_i^R$ and $d_i^R\equiv \bd_i^R$ for all $i\in\cN$, and $x_i^R\equiv \bx_i^R$ for all $i\in \cN\setminus \cK$.
Therefore, using induction hypothesis for systems with $k+k_1\geq k(N)+1$ infinite queues, there exists a subsequence $\{R\}$ along which with probability 1,
$$(\bx_i^R(\cdot), \ba_i^R(\cdot), \bd_i^R(\cdot))_{i\in\cN}\to (\bx_i(\cdot), \ba_i(\cdot), \bd_i(\cdot))_{i\in\cN},$$
where $\bx_i\equiv 0$ for all $i\in \cN\setminus (\cK\cup\cK_1)$, and $\bx_j\equiv \infty$ with $\bx_j'\equiv -\varepsilon(k+k_1)<0$ for all $j\in\cK\cup\cK_1$.
Consequently, in the time interval $[0,\tau]$, along that subsequence with probability 1,
$$(x_i^R(\cdot), a_i^R(\cdot), d_i^R(\cdot))_{i\in\cN}\to (x_i(\cdot), \ba_i(\cdot), \bd_i(\cdot))_{i\in\cN}$$
with $x_i= \bx_i\equiv 0$ for all $i\in\cN\setminus (\cK\cup\cK_1)$ and $x_i'\equiv -\varepsilon(k+k_1)<0$ for all $i\in\cK\cup\cK_1$,
where $\tau = x_{(k+k_1)}/\varepsilon(k+k_1)>0$.
Observe that the above argument can be extended till the time $\sum_{i\in \cN\setminus\cK}x_i(t)$ hits zero.
Furthermore, following the argument as above, this time is given by $T(k(N),\xx(0))$ as given in~\eqref{eq:Tkx}.
This completes the proof of Part 1 (iii).

\subsection{Conventional fluid-limit stability to verify Part (2)}\label{ssec:part2}
As mentioned earlier, we will use the fluid limit technique of proving stochastic stability as in~\cite{RS92, S95, D99} to prove Part (2).
Consider a sequence of initial states with increasing norm $R$, i.e., $\sum_{i\in\cN\setminus\cK}X_i^{R}(0) =R$ and $X_i^R(0)=\infty$ for $i\in\cK$.
Then from Part (1.iii), we know that for any sequence there exists a further subsequence $\{R\}$ along which with probability 1, the fluid-scaled occupancy process $(x_i^R(\cdot))_{i\in\cN}$ converges to the process $(x_i(\cdot))_{i\in\cN}$ 
for which $\sum_{i\in\cN\setminus\cK}x_i(t)$ hits 0 in finite time $T(k(N), \xx(0))$, and stays at 0 afterwards.
This verifies the fluid-limit stability condition in~\cite[Theorem 4.2]{D99} and \cite[Theorem 7.2]{S95}, and thus completes the verification of Part (2).

\subsection{Large-scale asymptotics to verify Part (3)}\label{ssec:part3}
The verification of the backward induction step for Part~(3) uses contradiction.
Namely, assuming that the induction step for Part~(3) does not hold, we will construct a sequence of systems with increasing~$N$, for which we obtain a contradiction using Lemmas~\ref{lem:expo-bound} and~\ref{lem:steady-concentration}.
We note that this is the only part in the proof of Proposition~\ref{prop:largeN}, where we use the large-scale (i.e., $N\to\infty$) asymptotic results.

Observe that we have already argued in Subsection~\ref{ssec:roadmap} that \emph{for all} $N$, Parts (1) -- (4) hold for $k=N$.
Now, if for some $N$, Part (3) does not hold for some $k(N)\in \{0,1,\ldots, N-1\}$ while Parts (1)--(4) hold for all $k\geq k(N)+1$, then from the proofs of Parts (1) and (2), note that Parts (1) and (2) hold for $k=k(N)$ as well.
Consequently, the subsystem with $N-k(N)$ finite queues is stable.
Thus we have the following implication.

\begin{implication}\label{impl:part3not}
Suppose, for infinitely many $N$, the induction step to prove Part~(3) of Proposition~\ref{prop:largeN} does not hold for some $k=k(N)$.
 Then there exists a subsequence of $\{N\}$ $($which we still denote by $\{N\})$ diverging to infinity, such that (i) The system with $k(N)$ infinite queues is stable and (ii) The steady-state arrival rate into each infinite queue is at least 1.
\end{implication}

We will now show that Implication~\ref{impl:part3not} leads to a contradiction -- this will prove Part~(3) of Proposition~\ref{prop:largeN}.
Suppose Implication~\ref{impl:part3not} is true. 
Choose a further subsequence $\{N\}$ along which $k(N)/N$ converges to $\kappa\in [0,1]$.
As in the statement of Lemma~\ref{lem:steady-concentration} we will consider two regimes depending on whether $\kappa\geq 1-\lambda$ or not, and arrive at contradictions in both cases.
Since all the infinite queues are exchangeable, we will use $\sigma$ to denote a typical infinite queue.\\

\noindent
\textbf{Case 1.} 
First consider the case when $\kappa\geq 1-\lambda$.
Note that the expected steady-state instantaneous rate of arrival to $\sigma$ is given by
\begin{equation}\label{eq:inst-rate}
\begin{split}
&\expt\Big(\frac{\lambda N}{Q_1^{N}(\infty)}\ind{Q_1^N(\infty)+\Delta_0^N(\infty)+\Delta_1^N(\infty)=N}\Big)\leq \expt\Big(\frac{\lambda N}{Q_1^{N}(\infty)}\Big)= \expt\Big(\frac{\lambda}{q_1^{N}(\infty)}\Big)+o(1).
\end{split}
\end{equation}
Now observe that for large $N$, $\lambda/q_1^N(\infty)\leq 2\lambda/\kappa$, since $q_1^N(s)\geq \kappa/2>0$.
Further from Lemma~\ref{lem:steady-concentration} we know that $q_1^N(\infty)\pto 1$ as $N\to \infty$.
Consequently, $\expt(\lambda/q_1^N(\infty))\to \lambda$ as $N\to\infty.$
Therefore for large enough $N$,
\begin{equation}\label{eq:arrival1}
\begin{split}
&\expt\Big(\frac{\lambda N}{Q_1^{N}(\infty)}\ind{Q_1^N(\infty)+\Delta_0^N(\infty)+\Delta_1^N(\infty)=N}\Big)
\leq \frac{1+\lambda}{2} = 1- \frac{1-\lambda}{2}<1,
\end{split}
\end{equation}
which is a contradiction to Part (ii) of Implication~\ref{impl:part3not}.\\

\noindent
\textbf{Case 2.} 
In case $\kappa<1-\lambda$, first note that the statement in Part (3) is vacuously satisfied if $k(N) \equiv 0$ for all large enough $N$.
Thus without loss of generality, assume that $k(N)>0$.
Fix $\varepsilon_1$ as in Lemma~\ref{lem:expo-bound}.
In that case \eqref{eq:inst-rate} becomes
\begin{align*}
&\expt\Big(\frac{\lambda N}{Q_1^{N}(\infty)}\ind{Q_1^N(\infty)+\Delta_0^N(\infty)+\Delta_1^N(\infty)=N}\Big) \\
&\leq \expt\Big(\frac{\lambda N}{Q_1^{N}(\infty)}\ind{Q_1^N(\infty)+\Delta_0^N(\infty)+\Delta_1^N(\infty)=N,\ Q_1^N(\infty)\geq \varepsilon_1 N}\Big) \\
&\hspace{6cm}+ \lambda N \Pro{Q_1^{N}(\infty)<\varepsilon_1 N}\\
&\leq \frac{\lambda N}{\varepsilon_1 N}\Pro{Q_1^N(\infty)+\Delta_0^N(\infty)+\Delta_1^N(\infty)=N} 
+ \lambda N \Pro{Q_1^{N}(\infty)<\varepsilon_1 N}.
\end{align*}
Now, due to Part~(2) of Lemma~\ref{lem:steady-concentration}, we know that 
\[\Pro{Q_1^N(\infty)+\Delta_0^N(\infty)+\Delta_1^N(\infty)=N}\to 0,\] 
and furthermore, Lemma~\ref{lem:expo-bound} yields 
\[N \Pro{Q_1^{N}(\infty)<\varepsilon_1 N}\to 0\quad \mbox{as}\quad  N\to\infty.\]
Thus, 
\begin{align}\label{eq:arrival2}
&\expt\Big(\frac{\lambda N}{Q_1^{N}(\infty)}\ind{Q_1^N(\infty)+\Delta_0^N(\infty)+\Delta_1^N(\infty)=N}\Big) \to 0 \quad\mbox{as}\quad N\to\infty.
\end{align}
In particular, for large enough $N$, the expected steady-state arrival rate is bounded away from 1, which is again a contradiction to Part (ii) of Implication~\ref{impl:part3not}. 
This completes the verification of Part (3) of the backward induction hypothesis.

\subsection{Time-scale separation to verify Part (4)}\label{ssec:part4}
Assume Parts (1) -- (3) hold for all $k\in \{k(N), k(N)+1,\ldots, N\}$.
Now consider a system containing $k=k(N)$ infinite queues with indices in $\cK$, and recall the conventional fluid scaling and FSP from Subsection~\ref{ssec:fluid-sasha}.
Also, in this subsection whenever we refer to the process $\{\XX(t)\}_{t\geq 0}$, the components in $\cK$ should be taken to be infinite.

For the queue length vector $\XX$, define the norm $\|\XX\|:=\sum_{i\notin\cK}X_i$ to be the total number of tasks at the finite queues.
Lemmas~\ref{lem:norm-decay} and~\ref{lem:hitting} state two hitting-time results that will be used in verifying Part (4).
\begin{lemma}\label{lem:norm-decay}
For any fixed $\gamma\in (0,1)$, there exists $\tau=\tau(\gamma)$ and $C=C(\gamma)$, such that 
if $\|\XX(0)\|=R\geq C,$ then \[\expt\|X(R\tau)\|\leq (1-\gamma)\|X(0)\|.\]
\end{lemma}
Lemma~\ref{lem:norm-decay} says that if the system starts from an initial state where the total number of tasks in the finite queues is suitably large, then the time it takes until the expected total number of tasks in the finite queues falls below a certain fraction of the initial number, is proportional to itself.
The proof of Lemma~\ref{lem:norm-decay} is fairly straightforward, but is provided below for completeness.
\begin{proof}[Proof of Lemma~\ref{lem:norm-decay}]
Consider a sequence of initial states with an increasing norm, i.e., $\XX^R(0)$ is such that $\|\XX^R(0)\|=R$ where $R^{-1}\XX^R(0)\to\xx(0)$ as $R\to\infty$.
Then from Part 1 we know that as $R\to\infty$, on the time interval $[0,T(m,\xx(0))]$ the process $R^{-1}\XX^R(Rt)$ converges in probability to the unique deterministic process $\xx(t)$ satisfying
\begin{equation}\label{eq:drift}
\sum_{i\in\cN\setminus\cK}x_i'(t)<-c(k(N)) 
\quad \mbox{whenever}\quad \sum_{i\in\cN\setminus\cK}x_i(t)>0,
\end{equation}
where $c(m)= \min\big\{k\varepsilon(k):k(N)+1\leq k\leq N\big\}>0$.
We also know that for any $i\in\cN\setminus\cK$, if $x_i(t_0)=0$ for some $t_0$, then $x_i(t)=0$ for all $t\geq t_0$.
Consequently, since $c(k(N))$ is positive, there exists $\tau= \tau(\gamma)<\infty$, such that 
\[\sup_{\|\xx\|=1}\Big\{\|\xx(\tau)\|: \xx(0) = \xx\in[0,1]^{N-k(N)}\times\{\infty\}^{k(N)}\Big\}<1-\gamma.\]
Now since the expected number of arrivals into the $R$-th system up to time $t$, when scaled by $R$, is $\lambda t$ for any finite $t$, we obtain $\expt(R^{-1}\|X^R(t)\|)\leq 1+\lambda t.$
Therefore, the convergence in probability also implies the convergence in expectation.
Thus for the above choice of $\gamma$,
\[\limsup_{R\to\infty}\expt\Big(\frac{\|\XX^R(R\tau)\|}{R}\Big)<1-\gamma.\]
Hence, there exists $C$ such that for all $R\geq C$,
\[\expt\Big(\frac{\|\XX^R(R\tau)\|}{R}\Big)=\expt\Big(\frac{\|\XX^R(R\tau)\|}{\|\XX^R(0)\|}\Big)\le 1-\gamma.\]
This completes the proof of Lemma~\ref{lem:norm-decay}.
\end{proof}

For any $C>0$, define the set $\cC:=\{\|\XX\|\leq C\}$, and the stopping time $\theta_C := \inf\ \{t: \XX(t)\in\cC\}$.
For large enough $C$, the next lemma bounds the expected hitting time to the fixed set $\cC$ in terms of the norm of the initial state.
\begin{lemma}\label{lem:hitting}
There exists $C, C_1>0$, such that if $\|\XX(0)\|=R\geq C,$ then \[\expt(\theta_C|\XX(0))\leq C_1\|\XX(0)\|.\]
\end{lemma}
\begin{proof}[Proof of Lemma~\ref{lem:hitting}]
Fix any $\gamma\in (0,1)$, and take $\tau = \tau(\gamma)$ and $C=C(\gamma)$ as in Lemma~\ref{lem:norm-decay}.
For $i\geq 1$, define the sequence of random variables $T_i:=\tau \|\XX(T_{i-1})\|$ with the convention that $T_0=0$.
Now consider the discrete-time Markov chain $\{\Phi_i:i\geq 0\}$ adapted to the filtration $\boldsymbol{\cF}=\bigcup_{i\geq 0}\cF_i$, where $\Phi_i= \XX(T_i)$ is the value of the continuous-time Markov process sampled at times $T_i$'s, and $\cF_i=\sigma(\Phi_0,\Phi_1,\ldots,\Phi_i)$ is the sigma field generated by $\{\Phi_0,\Phi_1,\ldots,\Phi_i\}$.
Further, for $i\geq 0$ define the stopping time
$\htheta_C:= \inf\ \{j\geq 0:Z_j\leq C\}$.
Then observe that 
\[\theta_C\leq \sum_{i=1}^{\htheta_C}T_i=:\Psi_C.\]
Also define $\alpha_i = \sum_{j=1}^iT_j$ for $i\geq 1$, and hence $\alpha_{\htheta_C}=\Psi_C.$
Then as a consequence of Dynkin's lemma~\cite[Theorem 11.3.1]{MT93}, using~\cite[Proposition 11.3.2]{MT93} we have  
\[\expt(\theta_C)\leq \expt(\Psi_C)\leq \frac{\tau}{\gamma} \expt\|\XX(0)\|.\]
Choosing $C_1 = \tau/\gamma$ completes the proof.
\end{proof}
Now we have all the ingredients to verify Part (4) of the backward induction hypothesis.
Note that we now look at the sequence of conventional fluid-scaled processes starting at (scaled) time $T(k(N),\xx(0))$.
From the verification of Part (1) we already know that $x_i(t)= 0$ for all $t\geq T(k(N),\xx(0))$, $i\in\cN\setminus\cK$.
Thus, it only remains to show that starting from time $T(k(N),\xx(0))$, the drift of each of the infinite queues is at most $-\varepsilon(k(N))$.
Specifically, we will construct a probability space where the required probability 1 convergence holds.

In order to simplify writing, we assume that the system starts at time 0, and thus it is enough to consider a sequence of initial queue length vectors such that
\[\|\xx^R(0)\|\to 0\quad \mbox{as}\quad R\to\infty,\]
where $R$ is the parameter in the conventional fluid scaling.
Hence, Lemma~\ref{lem:hitting} yields that $R^{-1}\expt(\theta_C|\XX^R(0))\to 0$ as $R\to\infty$.
Consequently, $R^{-1}\theta_C\pto 0.$
Thus, the fluid-scaled time to hit the set~$\cC$ vanishes in probability, which is stated formally in the following claim.
\begin{claim}\label{claim:hittingC}
If the sequence of initial states is such that $\|\xx^R(0)\|\to 0$ as $R\to\infty$, then
$R^{-1}\theta_C\pto 0,$ as $R\to\infty$.
\end{claim}
\noindent
Now pick any (unscaled) state $z\in\cC$, and define the stopping time $\htheta_z$ as 
\[\htheta_z:=\inf\big\{t\geq 0: \XX(t) = z\big\}.\]
Since due to Part (2) of the backward induction hypothesis, the unscaled process $\XX(\cdot)$ is irreducible and positive recurrent, we have the following claim.
\begin{claim}\label{claim:hittingz}
If the sequence of initial states is such that $\xx^R(0)\in\cC$, then $R^{-1}\htheta_z\pto 0$, as $R\to\infty.$
\end{claim}
\noindent
Up to time $\htheta_z$, consider the product topology on the sequence space.
Then Claims~\ref{claim:hittingC} and~\ref{claim:hittingz} yield that for a sequence of initial states such that $\|\xx^R(0)\|\to 0$ as $R\to\infty$, there exists a subsequence $\{R\}$, along which with probability 1, $R^{-1}\htheta_z\to 0$.
Starting from the time $\htheta_z$, along the above subsequence, we construct the sequence of processes $\xx^R(\cdot)$ on the same probability space as follows.\\

\noindent
(1) Define the space of an infinite sequence of i.i.d.~renewal cycles of the unscaled process $\XX(\cdot)$, with the unscaled state $z$ being the renewal state, i.e.,
\[\Big\{\XX^{(i)}(t): 0\leq t\leq \htheta_z^{(i)}, \XX^{(i)}(0)=z\Big\}\] 
for $i=1,2,\ldots$ are i.i.d.~copies, and $\htheta_z^{(i)}$ are also i.i.d.~copies of $\htheta_z$.\\

\noindent
(2) Define the process $\XX^R(\cdot)$ as 
\[\XX^R(Rt)= \sum_{i=1}^\infty \XX^{(i)}\big(Rt- \Theta(i-1)\big)\ind{\Theta(i-1)\leq Rt<\Theta(i)}, \ \text{where}\
\Theta(i):= \sum_{j=1}^i\htheta^{(j)}.\]

\noindent
Let $A(t)$ denote the cumulative number of arrivals up to time $t$ to a fixed server with infinite queue length when the system starts from the state $z$.
Now, in order to calculate the drift of each of the infinite queues, 
observe that cumulative number of arrivals up to time $Rt$ to server  $n\in\cK$ in the $R$-th system can be written as
\begin{align*}
A_n^R(Rt) = \sum_{i=1}^{N_\theta^R} A_n^{(i)}+B_n(t-\Theta(N_\theta^R)), \quad \text{where}\quad N_\theta^R:= \max\{j:\Theta(j)\leq Rt\}.
\end{align*}
 $A_n^{(i)}$'s are i.i.d.~copies of the random variable $A(\htheta_z)$, $B_n(\cdot)$ is distributed as $A(t)$, and $A_n^{(i)}$'s and $B_n(\cdot)$ are independent of the random variable $N_\theta^R$.
 Now, since due to Part (2) of the backward induction hypothesis the subsystem consisting of the finite queues is stable, $\XX(\cdot)$ is irreducible and positive recurrent.
 Thus, we have $\expt(\htheta_z|\XX(0)=z)<\infty$, and hence, with probability 1, 
 \[\frac{N_\theta^R}{R}\to \frac{t}{\expt(\htheta_z|\XX(0)=z)}, \quad\mbox{as}\quad R\to\infty.\]
 Thus, using Part (3) of the backward induction hypothesis, SLLN yields, with probability~1,
 \begin{align*}
 \frac{1}{R}A_n^R(Rt) &= \frac{1}{R}\sum_{i=1}^{N_\theta^R} A_n^{(i)}+\frac{B_n(t-\Theta(N_\theta^R))}{R}\to \ha t, \quad\mbox{as}\quad R\to\infty,
 \end{align*}
 for some $\ha\leq 1-\varepsilon(k(N))$.
 Therefore, in the conventional fluid limit, $a_n(t)\leq  (1-\varepsilon(k(N)))t$. 
 Also, since the departure rate from each of the servers with infinite queue lengths is always~1, it can be seen that in the conventional fluid limit, $d_n(t) = t$, and thus, the drift of the $n$-th infinite queue is given by at most $-\varepsilon(k(N))$.
Combining the probability 1 convergence of the time $\htheta_z$ to 0, and the probability space constructed after time $\htheta_z$, we obtain that along the subsequence $\{R\}$ with probability 1, the fluid-scaled processes converges to a limit where each infinite queue has drift at most $-\varepsilon(k(N))$.
This completes the verification of Part (4), and hence of Proposition~\ref{prop:largeN}.

\section{Mean-field analysis for large-scale asymptotics}\label{sec:mf}
In this section we will analyze the large-$N$ behavior of the system.
In particular, we will prove Lemmas~\ref{lem:expo-bound} and~\ref{lem:steady-concentration}.
The next proposition is a basic mean-field fluid limit result that we need later.
Define 
\[E_\kappa := \Big\{(\bld{q},\dd)\in [0,1]^\infty:  q_i\geq q_{i+1}\geq \kappa,\ \forall i, \  \delta_0+\delta_1+ q_1\leq 1 \Big\}.\]

\begin{proposition}
\label{prop:mf}
Assume $k(N)/N\to\kappa\in [0,1]$ and the sequence of initial states $(\qq^N(0),\bld{\delta}^N(0))$ converge to a fixed $(\qq(0),\bld{\delta}(0))\in E_\kappa$, as $N\to\infty$,  where $q_1(0)>0$. 
Then, 
with probability 1, any subsequence of $\{N\}$ has a further subsequence along which $\{(\qq^N(t),\dd^N(t))\}_{t\geq 0}$ converges, uniformly on compact time intervals, to some deterministic trajectory $\{(\qq(t),\dd(t))\}_{t\geq 0}$ satisfying the following equations:
\begin{align*}
 q_i(t)&=q_i(0) +\int_0^t\lambda p_{i-1}(\qq(s),\dd(s),\lambda)\dif s
- \int_0^t(q_i(s)-q_{i+1}(s))\dif s,\ i\ge 1, \\
\delta_0(t)&=\delta_0(0)+\mu\int_0^t u(s)\dif s-\xi(t),  \\
\delta_1(t)&=\delta_1(0)+\xi(t)-\nu\int_0^t\delta_1(s)\dif s,\nonumber
\end{align*}
where 
\begin{align*}
u(t) &= 1- q_1(t) - \delta_0(t) - \delta_1(t),\\
\xi(t) &= \int_0^t\lambda(1-p_0(\qq(s),\dd(s),\lambda))\ind{\delta_0(s)>0}\dif s.
\end{align*}
For any $(\qq,\dd)\in E$, $\lambda>0$, $(p_i(\qq,\dd,\lambda))_{i\geq 0}$ are given by 
\begin{align*}
 p_0(\qq,\dd,\lambda) &= 
\begin{cases}
&1\qquad \text{if}\qquad u=1-q_1-\delta_0-\delta_1>0,\\
&\min \{\lambda^{-1}(\delta_1\nu + q_1-q_2), 1\},\quad\text{otherwise,}
\end{cases}\\
\quad p_i(\qq,\dd,\lambda)  &= (1-p_0(\qq,\dd,\lambda)) (q_{i}-q_{i+1})q_1^{-1},\  i \ge 1.
\end{align*}
\end{proposition}
This type of result is standard and is obtained using Functional Strong LLN, for example as in~\cite{Stolyar15,Stolyar17, MDBL17}; we omit its proof.
Also, we note that, while Proposition~\ref{prop:mf} is a version of Theorem\ref{th: fluid} in Chapter~\ref{chap:energy1}, it is different in that it is suitably modified for the case of infinite buffers and some queues being infinite, and it states a somewhat different type of convergence, convenient for the use in this chapter.
Define \emph{mean-field fluid sample path} (MFFSP) to be any deterministic trajectory satisfying the properties stated in Proposition~\ref{prop:mf}.

We now provide an intuitive explanation of the mean-field fluid limit stated in Proposition~\ref{prop:mf}.
It is similar to that behind~\cite[Theorem 3.1]{MDBL17}.
The term $u(t)$ corresponds to the asymptotic fraction of idle-on servers in the system at time $t$, and $\xi(t)$ represents the asymptotic cumulative number of server setups (scaled by $N$) that have been initiated during $[0,t]$.
The coefficient $p_i(\qq,\dd,\lambda)$ can be interpreted as the instantaneous fraction of incoming tasks that are assigned to some server with queue length $i$, when the fluid-scaled occupancy state is $(\qq,\dd)$ and the scaled instantaneous arrival rate is $\lambda$.
Observe that as long as $u>0$, there are idle-on servers, and hence all the arriving tasks
will join idle servers. 
This explains that if $u>0$, $p_0(\qq,\dd,\lambda) = 1$ and $p_i(\qq,\dd,\lambda)=0$ for $i=1,2,\ldots$.
If $u=0$, then observe that 
servers become idle at rate $q_1-q_2$, and servers in setup mode turn on at rate $\delta_1\nu$.
Thus the  idle-on servers are created at a total rate $\delta_1\nu + q_1-q_2$.
If this rate is larger than the arrival rate $\lambda$, then almost all the arriving tasks can be assigned to idle servers.
Otherwise, only a fraction $(\delta_1\nu + q_1-q_2)/\lambda$
of arriving tasks join idle servers. 
The rest of the tasks are distributed uniformly among busy servers, so a proportion $(q_{i}-q_{i+1})q_1^{-1}$ are assigned to servers having queue length~$i$.
For any $i=1,2, \ldots$, $q_i$ increases when there is an arrival to some server with queue length $i-1$, which occurs at rate $\lambda p_{i-1}(\qq,\dd,\lambda)$, and it decreases when there is a departure from some server with  queue length~$i$, which occurs at rate $q_i-q_{i-1}$. 
Since each idle-on server turns off at rate $\mu$, the fraction of servers in the off mode increases at rate 
$\mu u$.
Observe that if $\delta_0>0$, for each task that cannot be assigned to an idle server, a setup procedure is initiated  at one idle-off server. 
As noted above, $\xi(t)$ captures the (scaled) cumulative number of setup procedures initiated up to time~$t$.
Therefore the fraction of idle-off servers and the fraction of servers in setup mode decreases and increases by $\xi(t)$, respectively, during $[0,t]$.
Finally, since each server in setup mode becomes idle-on at rate $\nu$, the fraction of servers in setup mode decreases at rate $\nu\delta_1$.

\subsection{Proof of \texorpdfstring{Lemma 3.2}{Lemma~\ref{lem:expo-bound}}}
This subsection is devoted to the proof of Lemma~\ref{lem:expo-bound}.
Within this proof we will use the following terminology. 
Let $A^N$ be an event pertaining to $N$-th system. We will write $\Pro{A^N} = \eta(N)$ to mean the following property: \emph{There exist $C>0$  and $N_1>0$ such that $\Pro{A^N} \le \e^{-C N}$ for all $N \ge N_1$.}
If event $A^N$ depends on some parameter $p$ (say, the process initial state), we say that $\Pro{A^N} = \eta(N)$ \emph{uniformly in $p$} if the 
property holds for common fixed $C>0$  and $N_1>0$.

To prove the lemma, clearly, it suffices to prove that for some fixed $T_0>0$ and $\varepsilon_0>0$
\begin{equation}\label{eq:q-1-lower}
\Pro{q_1^N(T_0)\leq \varepsilon_0} = \eta(N),
\end{equation}
uniformly on the process initial states $(\qq^N(0),\dd^N(0))$. This is what we do in the rest of the proof.

Fix any $T_0>0$; $\varepsilon_0>0$ will be chosen later. We now prove several claims, which rather simply follow from 
the process structure and
basic large-deviations estimates (specifically, Cramer's theorem) -- they will serve as building blocks for the proof argument.

\begin{claim}
\label{claim-busy-decay}
{\normalfont(i)} For any $\varepsilon>0$, uniformly in $\tau\in [0,T_0]$ and uniformly in $q_1^N(0) \ge \varepsilon$, 
\begin{equation}
\label{eq-busy-decay}
\Pro{q_1^N(\tau) \le (\varepsilon/2) \e^{-T_0}} = \eta(N).
\end{equation}
{\normalfont(ii)} For any $\varepsilon>0$, uniformly in $\tau\in [0,T_0]$ and uniformly in $\delta_1^N(0) \ge \varepsilon$, 
\begin{equation}
\label{eq-setup-decay}
\Pro{\delta_1^N(\tau) \le (\varepsilon/2) \e^{-\nu T_0}} = \eta(N).
\end{equation}
\end{claim}

\noindent Indeed, to prove (\ref{eq-busy-decay}), observe that
any busy server at time $t$ stays busy in the interval $[t,t+\tau]$ with probability at least $\e^{-\tau}\ge \e^{-T_0}$. It remains to recall that
$q_1^N(0) \ge \varepsilon$ corresponds to at least $\varepsilon N$ busy servers in the unscaled system and apply Cramer's theorem. 
Statement~(ii) is proved analogously.

\begin{claim}
\label{claim-T1}
For any sufficiently small $T_1>0$, there exists $\varepsilon'_1>0$ such that, uniformly in the initial state $(\qq^N(0),\dd^N(0))$,
\begin{equation}
\label{eq-T1}
\Pro{q_1^N(T_1) + \delta_1^N(T_1) \le \varepsilon'_1} = \eta(N).
\end{equation}
\end{claim}

\noindent Indeed, fix any $T_1>0$ such that $\lambda T_1 \le 1/4$. Suppose first that either $q_1^N(0)\ge 1/4$ or $\delta_1^N(0)\ge 1/4$; uniformly on all such initial conditions, the claim follows by using Claim~\ref{claim-busy-decay}. 
Suppose now that $q_1^N(0)< 1/4$ and $\delta_1^N(0)< 1/4$, and therefore $\delta_0^N(0)+u^N(0) > 1/2$, where recall that $u^N$ is the fraction of idle-on servers. 
The (unscaled) number of new customer arrivals in $[0,T_1]$, denote it by $H[0,T_1]$, is Poisson with mean $\lambda T_1 N$; therefore, 
$$\Pro{| H[0,T_1]/N - \lambda T_1 | \ge (1/2) \lambda T_1} = \eta(N).$$ 
This means that with probability $1-\eta(N)$,  
we have $H[0,T_1]/N < \delta_0^N(0)+u^N(0)$, and therefore
each arrival in $[0,T_1]$ creates either a new busy server or a new setup server; furthermore, each of these newly created busy or setup servers will not change its state until time $T_1$ with probability at least $\e^{-\nu' T_1}$, where $\nu' = \max \{\nu,1\}$.
It remains to choose $\varepsilon'_1 \in (0, (1/4) \lambda T_1 \e^{-\nu' T_1})$ to obtain the claim.

\begin{claim}
\label{claim-T2}
For any $\varepsilon_1>0$ and any $T_2>0$, there exists $\varepsilon'_2>0$ such that, uniformly in $\delta_1^N(0) \ge \varepsilon_1$,
\begin{equation}
\label{eq-T2}
\Pro{q_1^N(T_2) + u^N(T_2) \le \varepsilon'_2} = \eta(N).
\end{equation}
\end{claim}

\noindent Indeed, at time $0$ there are at least $\varepsilon_1 N$ setup servers. Fix any $T_2>0$. In $[0,T_2]$ each of them tuns into an idle-on server with probability at least $1-\e^{-\nu T_2}$; those servers that do turn into idle-on will be either still be idle-on or busy at time $T_2$ with probability at least $\e^{-\nu'' T_2}$, where $\nu'' = \max \{\mu,\nu\}$.
It remains to choose $\varepsilon'_2 \in (0, (1/2) \varepsilon_1 \e^{-\nu'' T_2})$, and apply Cramer's theorem.

\begin{claim}
\label{claim-T3}
For any $\varepsilon_2>0$ and any sufficiently small $T_3>0$, there exists $\varepsilon_3>0$ such that, uniformly in $u^N(0) \ge \varepsilon_2$,
\begin{equation}
\label{eq-T3}
\Pro{q_1^N(T_3) \le \varepsilon_3} = \eta(N).
\end{equation}
\end{claim}

\noindent Indeed, fix $T_3$ small enough so that $\e^{-\mu T_3}>3/4$ and $\lambda T_3 < \varepsilon_2/2$. 
At time $0$ there are at least $\varepsilon_2 N$ idle-on servers; with probability at least $\e^{-\mu T_3}>3/4$ they will still be
idle-on at time $T_3$, \emph{unless} they are taken by a new arrival. The (unscaled) number of new arrivals in $[0,T_3]$, namely
$H[0,T_3]$, is Poisson with mean $\lambda T_3 N$, and therefore $H[0,T_3]/N$ concentrates at $\lambda T_3$:
$\Pro{| H[0,T_3]/N - \lambda T_3 | \ge (1/2) \lambda T_3} = \eta(N)$. We conclude that with probability $1-\eta(n)$ every new arrival
in $[0,T_3]$ will go to an idle-on server and turn it into busy; each of those servers, in turn, will remain busy until $T_3$ with 
probability at least $\e^{-T_3}$. 
It remains to choose $\varepsilon_3 \in (0, (1/4) \lambda T_3 \e^{-T_3})$ to obtain the claim.

With these claims, we are now in a position to conclude the proof of the lemma. Choose small $T_1>0$ and $\varepsilon'_1>0$
as in Claim~\ref{claim-T1}; and then $\varepsilon_1=\varepsilon'_1/2$. 
For the chosen $\varepsilon_1$,
choose small $T_2>0$ and $\varepsilon'_2>0$
as in Claim~\ref{claim-T2}; and then $\varepsilon_2=\varepsilon'_2/2$. Finally, for the chosen $\varepsilon_2$, choose small $T_3>0$ and $\varepsilon_3>0$
as in Claim~\ref{claim-T3}. Note that $T_1, T_2, T_3$ can be taken small enough so that $T'_3\doteq T_1+T_2+T_3 \le T_0$; 
let us also denote
$T'_2=T_1+T_2$. Choose $\varepsilon_0 = (1/2) \min\{\varepsilon_1,\varepsilon_2,\varepsilon_3\} \e^{-T_0}$.

According to Claim~\ref{claim-T1}, with probability $1-\eta(N)$, at time $T_1$ we have either $q_1^N(T_1) \ge \varepsilon_1$ or
$\delta_1^N(T_1) \ge \varepsilon_1$. Conditioned on a state at $T_1$ safisfying $q_1^N(T_1) \ge \varepsilon_1$, we have 
(\ref{eq:q-1-lower}) by applying Claim~\ref{claim-busy-decay}. Therefore, it remains to prove (\ref{eq:q-1-lower}) conditioned on a state at $T_1$ satisfying $\delta_1^N(T_1) \ge \varepsilon_1$. Under this condition at $T_1$, we obtain from Claim~\ref{claim-T2} that, with probability $1-\eta(N)$, at time $T'_2$ we have either $q_1^N(T'_2) \ge \varepsilon_2$ or
$u^N(T'_2) \ge \varepsilon_2$. Then, conditioned on a state at $T'_2$ satisfying $q_1^N(T'_2) \ge \varepsilon_2$, we have 
(\ref{eq:q-1-lower}) by once again applying Claim~\ref{claim-busy-decay}. It now remains to prove (\ref{eq:q-1-lower}) conditioned on a state at $T'_2$ satisfying $u^N(T'_2) \ge \varepsilon_2$. Under this condition at $T'_2$, we obtain from Claim~\ref{claim-T3} that, with probability $1-\eta(N)$, at time $T'_3$ we have $q_1^N(T'_3) \ge \varepsilon_3$; and conditioned on  $q_1^N(T'_3) \ge \varepsilon_3$
at $T'_3$, we have (\ref{eq:q-1-lower}) by, yet again, Claim~\ref{claim-busy-decay}. The proof is complete.

\subsection{Proof of \texorpdfstring{Lemma 3.3}{Lemma~\ref{lem:steady-concentration}}}
In this subsection we will prove Lemma~\ref{lem:steady-concentration}.
Recall that the stability of the subsystem $\cN\setminus\cK$ is assumed,
and hence there exists a unique stationary distribution for each $N$. 
Recall that we denote by $\qq^N(\infty)$ the random value of $\qq^N(t)$ in the steady state. 
We will start by stating a few basic facts about the mean-field limits that will facilitate the proof of Lemma~\ref{lem:steady-concentration}.

Recall the definition of MFFSP from the paragraph after Proposition~\ref{prop:mf}, and that $u(t)=1-q_1(t)-\delta_0(t)-\delta_1(t)$. 
Also, denote by $y_1(t) = q_1(t)-q_2(t)$ and by
$(d^+/dt)$ the right derivative.
\begin{claim}\label{fact:mfl}
For any $\varepsilon>0$ there exists $\alpha>0$, such that any MFFSP with $q_1(0)>0$ satisfies the following properties for all $t\geq 0$:
\begin{enumerate}[{\normalfont (i)}]
\item If $y_1(t) \le \lambda-\varepsilon$ and $u(t)>0$, then $(d^+/dt)q_1(t) \ge \alpha$.
\item If $y_1(t) \le \lambda-\varepsilon$, $u(t)=0$ and $\delta_1(t)\ge \varepsilon$, then $(d^+/dt)q_1(t) \ge \alpha$.
\item If $y_1(t) \le \lambda-\varepsilon$, $u(t)=0$, $\delta_1(t)=0$, and $\delta_0(t)>0$, then $(d^+/dt)\delta_1(t)\ge \varepsilon$.
\end{enumerate}
\end{claim}
\noindent
\begin{proof}
Fix any $\varepsilon>0$. First observe that since $q_1(0)>0$ and due to Proposition~\ref{prop:mf}, $q_1(0)$ is nondecreasing whenever $q_1(t)-q_2(t)\leq \lambda$, we have
$q_1(t)\geq \min\{q_1(0),\lambda\}>0$ for all $t\geq 0$.
Thus, Proposition~\ref{prop:mf} can be applied for all $t\geq 0$, throughout the MFFSP.
Choose $\alpha = \min\{\varepsilon\nu,\ \varepsilon\}$.

For (i), note that if  $y_1(t) \leq \lambda - \varepsilon$ and $u(t)>0$, then 
$$(d^+/dt)q_1(t) = \lambda - (q_1(t) - q_2(t))\geq \varepsilon\geq \alpha.$$

For (ii), note that if $y_1(t) \leq \lambda - \varepsilon$, $u(t)>0$, and $\delta_1(t)\geq \varepsilon$, then due to Proposition~\ref{prop:mf},
\begin{align*}
(d^+/dt)q_1(t) &= \min\big\{(\delta_1(t)\nu + q_1(t) - q_2(t)), \lambda\big\} - (q_1(t) - q_2(t))\\
	&= \min\big\{\delta_1(t)\nu,\quad \lambda - (q_1(t) - q_2(t))\big\}
	\geq \min\big\{\varepsilon\nu,\ \ \varepsilon\big\} = \alpha.
\end{align*}

Finally, for (iii), note that from Proposition~\ref{prop:mf} if $y_1(t) \le \lambda-\varepsilon$, $u(t)=0$, $\delta_1(t)=0$, and $\delta_0(t)>0$, then $(d^+/dt)\delta_1(t) =\lambda-(q_1(t) - q_2(t)) \ge \varepsilon$.
\end{proof}

\noindent
\textbf{Proof of statement (1).}  Note that it is enough to prove the following property of any MFFSP:
\begin{claim}\label{claim:q1case1}
Starting from any state $\qq(0)\in E_\kappa$ with $\kappa\geq 1-\lambda$ and $q_1(0)\in [\kappa, 1)$, along any MFFSP we have $$\lim_{t\to\infty}q_1(t)= 1.$$
\end{claim}
\noindent
Indeed, Claim~\ref{claim:q1case1} implies that under the assumption of stability, asymptotically the stationary distribution of $q_1^N(t)$ must concentrate at $q_1^\star=1$, as $N\to\infty$.

\begin{proof}[Proof of Claim~\ref{claim:q1case1}]
We will prove by contradiction.
Note that for the case under consideration, $q_i(t)\geq \kappa$ for all $i\geq 1$ and $t\geq 0$.
Therefore, throughout the proof of Claim~\ref{claim:q1case1} we can assume $q_1(0)\geq \kappa>0$, and can apply Proposition~\ref{prop:mf} and Claim~\ref{fact:mfl}.

Note that if $q_1(t)< 1$, we have $q_1(t) - q_2(t)<1-\kappa\leq\lambda$, and hence due to Claim~\ref{fact:mfl}, $q_1(t)$ is non-decreasing.
Thus if Claim~\ref{claim:q1case1} does not hold, then there exists an $\varepsilon>0$, such that $q_1(t)\leq 1 - \varepsilon\nu$ for all $t\geq 0$, and hence
\begin{equation}\label{eq:localq1bdd}
q_1(t)-q_2(t)\leq \lambda - \varepsilon\nu \quad\mbox{for all}\quad t\geq 0.
\end{equation}
The high-level proof idea is that if $q_1(t)$ remains below 1 by a non-vanishing amount for all $t\geq 0$, then the (scaled) rate $q_1(t) - q_2(t)$ of busy servers turning idle-on would not be high enough to match the (scaled) rate $\lambda$ of incoming jobs.
If there are idle-on servers (as in Claim~\ref{fact:mfl}.(i)) or sufficiently many servers in setup mode (as in Claim~\ref{fact:mfl}.(ii)), then we can still assign incoming tasks to idle-on servers, but this drives up the fraction of busy servers $q_1(t)$ and cannot continue indefinitely since $q_1(t)\leq 1 - \varepsilon\nu$ for all $t\geq 0$.
This means that we cannot initiate an unbounded number of setup procedures (see Equation~\eqref{eq:contr1}).
At the same time, as argued above, we cannot continue assigning tasks to idle-on servers either. Thus, throughout the MFFSP, a positive fraction of the jobs are assigned to busy servers, which initiates an unbounded (scaled) number of setup procedures, and hence the contradiction. \\

\noindent
Define the subset $\cX_\kappa\subseteq E$ as
$$\cX_\kappa:= \Big\{(\qq,\dd)\in E_\kappa: q_1 + \delta_0 + \delta_1 = 1,  \delta_1\nu+q_1-q_2\leq \lambda\Big\},$$
and denote by  $\indn{\cX_\kappa}(\qq(s),\dd(s))$  the indicator of $(\qq(s),\dd(s))\in \cX_\kappa.$
Observe that due to Proposition~\ref{prop:mf}, $q_1(t)$ can be written as
\begin{equation}\label{eq:q1eqn}
\begin{split}
q_1(t)&= q_1(0) +\int_{0}^t\delta_1(s)\nu\indn{\cX_\kappa}(\qq(s),\dd(s))\dif s \\
&\hspace{3cm}+\int_{0}^t[\lambda - q_1(s) + q_2(s)] \indn{\cX^c_\kappa}(\qq(s),\dd(s))\dif s.
\end{split}
\end{equation}
Thus,
\begin{align*}
q_1(t)\geq q_1(0) +\int_{0}^t[\lambda - q_1(s) + q_2(s)] \indn{\cX^c_\kappa}(\qq(s),\dd(s))\dif s,
\end{align*}
and~\eqref{eq:localq1bdd} yields there exists positive constant $K_1$, which may depend on $\varepsilon$ such that $\forall\ t\geq 0$
\begin{equation}\label{eq:K1-choice-sasha}
\int_0^t \indn{\cX_{\kappa}^c}(\qq(s),\dd(s))\dif s<K_1
\implies \int_0^t\ind{u(s)>0}\dif s<K_1.
\end{equation}
Again from~\eqref{eq:q1eqn} we obtain
\begin{align*}
q_1(t)&\geq  q_1(0) +\int_{0}^t\delta_1(s)\nu\indn{\cX_\kappa}(\qq(s),\dd(s))\dif s\\
&\geq q_1(0)+\nu\int_{0}^t\delta_1(s)\dif s-(\nu+1)\int_{0}^t\indn{\cX_\kappa^c}(\qq(s),\dd(s))\dif s,
\end{align*}
and thus, by \eqref{eq:localq1bdd} and \eqref{eq:K1-choice-sasha}, there exist positive constants $K_2$, $K_2'$  which may depend on $\varepsilon$ such that $\forall\ t\geq 0$
\begin{equation}\label{eq:K2-choice-sasha}
\begin{split}
 \int_0^t\delta_1(s)\dif s<K_1 &\implies \int_0^t\ind{\delta_1(s)>\frac{\varepsilon}{2}}\dif s<K_2,
  \implies \int_0^t\ind{\delta_1(s)>\frac{\varepsilon\nu}{2}}\dif s<K_2'.
 \end{split}
\end{equation}
Consequently, due to Proposition~\ref{prop:mf}, since $\delta_1(t) = \delta_1(0) + \xi(t) - \nu\int_0^t\delta_1(s)\dif s,$
it must be the case that 
\begin{equation}\label{eq:contr1}
\limsup_{t\to\infty}\xi(t)<\infty.
\end{equation}
Furthermore, since $q_1(t)\leq 1-\varepsilon\nu$ for all $t\geq 0$,
\[\ind{\delta_0(t)=0} \leq \ind{u(t)>0}+\ind{\delta_1(t)\geq \frac{\varepsilon\nu}{2}}.\]
Thus, \eqref{eq:K1-choice-sasha} and \eqref{eq:K2-choice-sasha} yield $\forall\ t\geq 0$,
\begin{equation}\label{eq:delta0bdd}
\int_0^t\ind{\delta_0(s)=0}\dif s\leq K_1+K_2'.
\end{equation}
Now from Proposition~\ref{prop:mf} observe that
\begin{align*}
\xi(t)&= \int_0^t\lambda(1-p_0(\qq(s),\dd(s),\lambda))\ind{\delta_0(s)>0}\dif s\\
&\geq \int_{0}^t\lambda(1-p_0(\qq(s),\dd(s),\lambda))\ind{\delta_0(s)>0, u(s) = 0,\delta_1(s)\leq \varepsilon/2}\dif s,
\end{align*}
and on the set $\{s:\delta_0(s)>0, u(s) = 0,\delta_1(s)\leq \varepsilon/2\}$ we have  $p_0(\qq(s),\dd(s),\lambda)\leq \lambda^{-1}(\varepsilon\nu/2+q_1(s)-\kappa)$. Therefore,
\begin{equation}\label{eq:xilowerbdd}
\begin{split}
\xi(t)
&\geq \int_{0}^t\Big(\lambda-\frac{\varepsilon\nu}{2}-q_1(s)+\kappa\Big)\ind{\delta_0(s)>0, u(s) = 0,\delta_1(s)\leq \varepsilon/2}\dif s\\
&\geq \int_0^t\Big(\lambda-\frac{\varepsilon\nu}{2}-q_1(s)+\kappa\Big)\dif s - \int_0^t\ind{\delta_0(s)=0}\dif s 
- \int_0^t\ind{u(s)>0}\dif s\\
&\hspace{6.75cm} - \int_0^t\ind{\delta_1(s)>\varepsilon/2}\dif s,
\end{split}
\end{equation}
where the second inequality is due to the fact that $\lambda-\varepsilon\nu/2-q_1(s)\leq \lambda <1$.
Therefore, since $\lambda+\kappa\geq 1$, we have $\lambda-\varepsilon\nu/2-q_1(s)+\kappa\geq \varepsilon\nu/2>0$, and Equations~\eqref{eq:K1-choice-sasha}, \eqref{eq:K2-choice-sasha}, \eqref{eq:delta0bdd}, and \eqref{eq:xilowerbdd} implies $\liminf_{t\to\infty}\xi(t) = \infty$, which is a contradiction with~\eqref{eq:contr1}.
This completes the proof of Claim~\ref{claim:q1case1}.
\end{proof}

\noindent
\textbf{Proof of statement (2).} First we will establish convergence of $q_1^N(\infty)$ as $N\to\infty$, followed by convergence of $q_2^N(\infty)$, $\delta_0^N(\infty)$, and $\delta_1^N(\infty)$.\\

\noindent
\textbf{Convergence of $q_1^N(\infty)$.} 
First we will show that 
for all $\varepsilon_2>0$, 
\begin{equation}\label{eq:q1limsup}
\limsup_{N\to\infty} \P(q_1^N(\infty) < \kappa+\lambda-\varepsilon_2) = 0.
\end{equation}
Due to Lemma~\ref{lem:expo-bound}, note that any limit of stationary distributions is such that with probability~1, $q_1 \ge \varepsilon_1$ for some fixed $\varepsilon_1>0$.
Therefore, throughout the proof of Part (2) of Lemma~\ref{lem:steady-concentration}, it is enough to consider MFFSP so that $q_1(0)\geq \varepsilon_1$, and Proposition~\ref{prop:mf} and Claim~\ref{fact:mfl} can be used.
Thus, for~\eqref{eq:q1limsup}, it is enough to show that any MFFSP has the following property:
\begin{claim}\label{claim:q1lower}
Starting from any state $\qq(0)\in E_\kappa$ with $q_1(0)\in [\varepsilon_1, \kappa+\lambda)$, along any MFFSP we have 
\begin{equation}\label{eq:q1part2}
\liminf_{t\to\infty}q_1(t)\geq \kappa+\lambda.
\end{equation}
\end{claim}
\noindent
Similar arguments as in the proof of Claim~\ref{claim:q1case1} can be used to prove Claim~\ref{claim:q1lower}, for which we omit the details.
Claim~\ref{claim:q1lower} then implies~\eqref{eq:q1limsup}.

Further, observe that since we have assumed that the system is stable, we have 
\begin{equation}\label{eq:exptbdd}
\lim_{N\to\infty}\expt(q_1^N(\infty)) \leq \kappa+\lambda.
\end{equation}
Fix any $\varepsilon_2'>0$. 
Now for all fixed $M>0$, 
\begin{align*}
&\expt(q_1^N(\infty))-(\kappa+\lambda) \geq \varepsilon_2'\P(q_1^N(\infty) > \kappa+\lambda+\varepsilon_2')- \P\big(q_1^N(\infty) < \kappa+\lambda-\frac{\varepsilon_2'}{M} \big) \\
&\hspace{5cm}-\frac{\varepsilon_2'}{M}\P\big(\kappa+\lambda-\frac{\varepsilon_2'}{M}\leq q_1^N(\infty) \leq \kappa+\lambda+\varepsilon_2'\big),\\
& \geq \varepsilon_2'\P(q_1^N(\infty) > \kappa+\lambda+\varepsilon_2')-\frac{\varepsilon_2'}{M}- \P\big(q_1^N(\infty) < \kappa+\lambda-\frac{\varepsilon_2'}{M} \big),
\end{align*}
and thus, from \eqref{eq:q1limsup} and \eqref{eq:exptbdd} above,
\[\limsup_{N\to\infty} \P(q_1^N(\infty) > \kappa+\lambda+\varepsilon_2') \leq \frac{1}{M}\quad\mbox{for all }M>0\]
which in conjunction with~\eqref{eq:q1limsup} completes the proof of convergence of $q_1^N(\infty)$.\\

\noindent
\textbf{Convergence of $q_2^N(\infty)$.} 
Note that given the above convergence of $q_1^N(\infty)$ to $\kappa+\lambda$ as $N\to\infty$, the following property of the mean-field limit is sufficient to prove that the sequence of stationary distributions $q_2^N(\infty)$ concentrate at $q_2^\star = \kappa$ as $N\to\infty$:
\begin{claim}\label{claim:q2conv}
For any $\varepsilon_3>0$, there exists a fixed $T_0$ and $\varepsilon_4>0$, such that starting from any state $\qq(0)\in E_\kappa$ with $q_1(0)=\lambda+\kappa$ and $q_2(0)> \kappa+\varepsilon_3$, along any MFFSP we have $q_1(T_0) \geq  \kappa+\lambda + \varepsilon_4$.
\end{claim}
\noindent
Indeed, if the sequence of the stationary distributions were such that 
$$\limsup_{N\to\infty} \P\big(q_2^N(\infty)>\kappa+\varepsilon_3\big)>0,$$ then Claim~\ref{claim:q2conv} would imply that $\limsup_{N\to\infty} \P\big(q_1^N(\infty)>\kappa+\lambda+\varepsilon_4/2\big)>0$, which contradicts the convergence of $q_1^N(\infty)$.

\begin{proof}[Proof of Claim~\ref{claim:q2conv}]
We will prove by contradiction.
Note that since $q_1(0)=\lambda+\kappa$ and $q_2(0)>\kappa+\varepsilon_3$, and the rates of change are bounded, in a sufficiently small neighborhood $[0,T_0]$ (depending only on $\varepsilon_3$), we have for all $t\in [0,T_0]$, (i) $q_1(t)\leq \lambda+\kappa+\varepsilon_3/2$, (ii)~$q_2(t)\geq \kappa+\varepsilon_3/2$, and
\[\mathrm{(iii)}~y_1(t)=q_1(t) - q_2(t)\leq \lambda-\frac{\varepsilon_3}{2}.\]
Since due to Claim~\ref{fact:mfl}, $q_1(t)$ is nondecreasing in $[0,T_0]$,
it is enough to produce a subinterval of $[0,T_0]$, where the right-derivative of $q_1(t)$ is bounded away from 0.
Now we will consider two cases:\\

\noindent
\textbf{Case 1:} There exists $t'\in [0,T_0/2]$, such that $u(t')=0$ and $\delta_1(t')\leq\varepsilon_3/2$.
In this case, $\delta_0(t')>0$, and in a sufficiently small time interval almost all points (with respect to Lebesgue measure) are regular for $\delta_1(t)$.
Also, due to Proposition~\ref{prop:mf}, since for $t\geq t'$,
\[\delta_1(t)=\delta_1(t')+\xi(t)-\xi(t')-\nu\int_{t'}^t\delta_1(s)\dif s,\]
with $(d^+/dt)\xi(t)=\lambda-y_1(t)\geq \varepsilon_3/2$ at $t=t'$, we have for sufficiently small $t_1<T_0/4$ (where the choice of $t_1$ does not depend on~$t'$), $\delta_1(t'+t_1)\geq t_1\varepsilon_3/4$.
Also, since the rate of decrease of $\delta_1(t)$ is bounded, there exists $t_2<T_0/4$ (where the choice of $t_2$ does not depend on~$t'$ either), such that, 
\[\delta_1(t)\geq \frac{t_1\varepsilon_3}{8}\quad\mbox{for all}\quad t\in [t'+t_1, t'+t_1+t_2]\subseteq [0,T_0].\] 
Thus, due to Claim~\ref{fact:mfl} there exists $\alpha>0$, such that during the time interval $[t'+t_1, t'+t_1+t_2]$, $(d^+/dt)q_1(t)\geq \min\{\alpha, \nu t_1\varepsilon_3/8\}$. 
Consequently,
\begin{equation}\label{eq:localq1case1}
q_1(T_0)\geq q_1(t'+t_1+t_2)\geq \lambda+\kappa+\min\Big\{\alpha, \frac{\nu t_1\varepsilon_3}{8}\Big\} t_2.
\end{equation}
It is important to note that the choices of $t_1$ and $t_2$ depend only on $\varepsilon_3$ and not on $t'$.

\noindent
\textbf{Case 2:} For all $t\in [0,T_0/2]$, either $u(t)>0$ or $\delta_1(t)>\varepsilon_3/2$.
In this case, due to Claim~\ref{fact:mfl} (i) and (ii), there exists $\alpha>0$, such that $(d^+/dt)q_1(t)>\alpha$ for all $t\in [0,T_0/2]$.
Also, since $q_2(t)$ is non-decreasing in $[0,T_0]$.
we obtain
\begin{equation}\label{eq:localq1case2}
q_1(T_0) \geq q_1\Big(\frac{T_0}{2}\Big)\geq \lambda+\kappa+\frac{\alpha T_0}{2}.
\end{equation}

\noindent
Combining the two cases above, and choosing 
\[\varepsilon_4 = \min\Big\{\min\Big\{\alpha, \frac{\nu t_1\varepsilon_3}{8}\Big\} t_2,  \frac{\alpha T_0}{2}\Big\}>0\]
completes the proof of Claim~\ref{claim:q2conv}.
\end{proof}

\noindent
\textbf{Convergence of $\delta_1^N(\infty)$ and $\delta_0^N(\infty)$.}
Given the convergence of $q_1^N(\infty)$ and $q_2^N(\infty)$, the convergence of $\delta_1^N(\infty)$ and $\delta_0^N(\infty)$ can be seen immediately by observing that the mean-field limit has the following property:
\begin{claim}\label{claim:delta0delta1}
Starting from any state $\qq(0)\in E_\kappa$ with $q_1(0)=\lambda+\kappa$ and $q_2(0)=\kappa$, along any MFFSP $\delta_1(t)\to 0$ and $\delta_0(t)\to 1-\lambda-\kappa$ as $t\to\infty$.
\end{claim}
\noindent
The proof of Claim~\ref{claim:delta0delta1} is immediate from the description of the mean-field limit as in Proposition~\ref{prop:mf}, and hence is omitted.\\

\noindent
The proof of the statement in~\eqref{eq:foundbusy} follows by using the convergence of steady states and the PASTA property.
This completes the proof of Lemma~\ref{lem:steady-concentration}.

\section{Conclusion}\label{sec:con}
In this chapter we studied the stability of systems under the TABS scheme and established large-scale asymptotics of the sequence of steady states.
Understanding stability of stochastic systems is of fundamental importance.
Systems under the TABS scheme, as it turned out, may be unstable for some $N$ even under a sub-critical load assumption.
As in many other cases, the lack of monotonicity makes the stability analysis much more challenging  from a methodological standpoint.
We developed a novel induction-based method and establish that the TABS scheme is stable for all large enough $N$.
The proof technique is of independent interest and potentially has a much broader applicability.
The key model-dependent part of our method is what can be called a weak monotonicity property, which ensures that for large enough $N$, with high probability, no matter where the system starts, in some fixed amount of time, there will be a certain fraction of busy servers.
Both traditional fluid limits (fixed $N$, initial state goes to infinity) and mean-field limits (for a sequence of processes with the number of queues $N \to \infty$) were used in an intricate manner to establish the results.

%% file: networkjsq.tex
\begin{abstract}
We consider a system of $N$~servers inter-connected by some underlying graph topology~$G_N$.  Tasks with unit-mean exponential processing times arrive at the various servers as independent Poisson processes of rate~$\lambda$. Each incoming task is irrevocably assigned to whichever server has the smallest number of tasks among the one where it appears and its neighbors in~$G_N$.  

The above model arises in the context of load balancing in large-scale cloud networks and data centers, and has been extensively investigated in case $G_N$ is a clique.  Since the servers are exchangeable in that case, mean-field limits apply, and in particular it has been proved that for any $\lambda < 1$, the fraction of servers with two or more tasks vanishes in the limit as $N \to \infty$. For an arbitrary graph $G_N$, mean-field techniques break down, complicating the analysis, and the queue length process tends to be worse than for a clique.  Accordingly, a graph $G_N$ is said to be $N$-optimal or $\sqrt{N}$-optimal when the queue length process on $G_N$ is equivalent to that on a clique on an $N$-scale or $\sqrt{N}$-scale, respectively.

We prove that if $G_N$ is an Erd\H{o}s-R\'enyi random graph with average degree $d(N)$, then with high probability it is $N$-optimal and $\sqrt{N}$-optimal if $d(N) \to \infty$ and $d(N) / (\sqrt{N} \log(N)) \to \infty$ as $N \to \infty$, respectively.  This demonstrates that optimality can be maintained at $N$-scale and $\sqrt{N}$-scale while reducing the number of connections by nearly a factor $N$ and $\sqrt{N} / \log(N)$ compared to a clique, provided the topology is suitably random. It is further shown that if $G_N$ contains $\Theta(N)$ bounded-degree nodes, then it cannot be $N$-optimal.  In addition, we establish that an arbitrary graph $G_N$ is $N$-optimal when its minimum degree is $N - o(N)$, and may not be $N$-optimal even when its minimum degree is $c N + o(N)$ for any $0 < c < 1/2$. Simulation experiments are conducted for various scenarios to corroborate the asymptotic results.
\end{abstract}

\section{Introduction}
\label{sec:intro-sig2}
In this chapter we explore the impact of the network topology
on the performance of load-balancing schemes in large-scale systems, as discussed in Section~\ref{networks}.
The chapter is organized as follows.
In Section~\ref{sec:model} we 
present a detailed model description
and introduce some useful notation and preliminaries.
Sufficient and necessary criteria for asymptotic optimality of deterministic graph sequences are developed in Sections~\ref{sec:det} and~\ref{sec:necessary}, respectively.
In Section~\ref{sec:random} we analyze asymptotic optimality of a sequence of random graph topologies.
In Section~\ref{sec:simulation-sig2} we present simulation experiments to support the analytical results, and examine the performance of topologies that are not analytically tractable.
We make a few brief concluding remarks and offer some suggestions
for further research in Section~\ref{sec:conclusion-sig2}.

\paragraph{Notation.} We adopt the usual notations O($\cdot$), o($\cdot$), $\omega(\cdot)$, and $\Omega(\cdot)$ to describe asymptotic comparisons.
For a sequence of probability measures $(\mathbb{P}_N)_{N\geq 1}$, the sequence of events $(\mathcal{E}_N)_{N\geq 1}$ is said to hold with high probability if $\mathbb{P}_N(\mathcal{E}_N)\to 1$ as $N\to\infty$.
Also, for some positive function $f(N):\N\to\R_+$, we write that a sequence of random variables $X_N$ is $\Op(f(N))$ or $\op(f(N))$ if $\{X_N/f(N)\}_{N\geq 1}$ is a tight sequence of random variables or converges to zero as $N\to\infty$, respectively.
The symbols `$\dto$' and `$\pto$' will denote  convergences in distribution and in probability, respectively.

\section{Model description and preliminaries}\label{sec:model}

Let $\{G_N\}_{N\geq 1}$ be a sequence of simple graphs indexed by the number of vertices $N$.
For the $N$-th system with $N$ servers, we assume that the servers are inter-connected by the underlying graph topology $G_N$, where server $i$ is identified with vertex $i$ in $G_N$, $i=1,2,\ldots, N$.
Tasks with unit-mean exponential processing times arrive at the various servers as independent Poisson processes of rate $\lambda$.
 Each server has its own queue with a fixed buffer capacity $b$ (possibly infinite).
 When a task appears at a server $i$, it is immediately assigned to the server with the shortest queue among server $i$ and its neighborhood in $G_N$.
If there are multiple such servers, one of them is chosen uniformly at random.
If $b<\infty$, and server $i$ and all its neighbors have $b$ tasks (including the ones in service), then the newly arrived task is discarded.
The service order at each of the queues is assumed to be oblivious to the actual service times, e.g.~First-Come-First-Served (FCFS).

For $k = 1,\ldots, N$, denote by $X_k(G_N,t)$ the queue length at the $k$-th server at time $t$ (including the one possibly in service), and by $X_{(k)}(G_N,t)$ the queue length at the $k$-th ordered server at time $t$ when the servers are arranged in  nondecreasing order of their queue lengths (ties can be broken in some way that will be evident from the context).
Let $Q_i(G_N, t)$ denote the number of servers with queue length at least $i$ at time $t$, $i = 1, 2,\ldots, b$, and $q_i(G_N,t):=Q_i(G_N,t)/N$ denote the corresponding fractions.
It is important to note that $\{(q_i(G_N, t))_{i\geq 1}\}_{t\geq 0}$ is itself \emph{not} a Markov process, but the joint process $\{(q_i(G_N, t))_{i\geq 1}, (X_{k}(G_N,t))_{k=1}^N\}_{t\geq 0}$ is Markov.
\begin{proposition}
\label{prop:ergod}
For any $\lambda<1$, the joint system occupancy process 
$$\{(q_i(G_N, t))_{i\geq 1}, (X_{k}(G_N,t))_{k=1}^N\}_{t\geq 0}$$ has a unique steady state 
$((q_i(G_N, \infty))_{i\geq 1}, (X_{k}(G_N,\infty))_{k=1}^N)$.
Also, the sequence of marginal random variables $\{(q_i(G_N, \infty))_{i\geq 1}\}_{N\geq 1}$ is tight with respect to the $\ell_1$-topology.
\end{proposition}
\begin{proof}[Proof of Proposition~\ref{prop:ergod}]
Note that if $b<\infty$, the process 
$$\{(q_i(G_N, t))_{i\geq 1}, (X_{k}(G_N,t))_{k=1}^N\}_{t\geq 0}$$ 
is clearly ergodic for all $N\geq 1$.
When $b=\infty$, to prove the ergodicity of the process, first fix any $N\geq 1$ and observe that 
the ergodicity of the queue length processes at the various vertices amounts to proving the ergodicity of the total number of tasks in the system.
Using the S-coupling and Proposition~\ref{prop:det ord} in Chapter~\ref{chap:univjsqd}, we obtain for all $t>0$,
\begin{equation}\label{eq:comparison}
\sum_{i=m}^\infty Q_i(G_N,t) \leq \sum_{i=m}^\infty Q_i(G_N',t),\quad\mbox{for all } m = 1,2, \ldots,
\end{equation}
provided the inequality holds at time $t=0$,
where $G_N'$ is the collection of $N$ isolated vertices.
Thus in particular, the total number of tasks in the system with $G_N$ is upper bounded by that with $G_N'$.
Now the queue length process on $G_N'$ is clearly ergodic since it is the collection of independent subcritical M/M/1 queues.
Next, for the $\ell_1$-tightness of $\{(q_i(G_N, \infty))_{i\geq 1}\}_{N\geq 1}$, we will use the following tightness criterion:
Define 
\begin{equation}\label{eq:s-space}
\mathcal{X} =
\left\{\qq \in [0, 1]^b: q_i \leq q_{i-1} \mbox{ for all } i = 2, \dots, b, \mbox{ and } \sum_{i=1}^b q_i < \infty\right\}
\end{equation}
as the set of all possible fluid-scaled occupancy states equipped with the $\ell_1$-topology.

Recall the criterion for $\ell_1$-tightness stated in Lemma~\ref{lem:tightcond} in Chapter~\ref{chap:univjsqd}.
Note that since $(q_i(G_N, \infty))_{i\geq 1}$ takes value in $[0,1]^\infty$, which is compact with respect to the product topology, Prohorov's theorem implies that $\big\{(q_i(G_N, \infty))_{i\geq 1}\big\}_{N\geq 1}$ is tight with respect to the product topology.
To verify the condition in~\eqref{eq:smalltail}, note that for each $m\geq 1$, Equation~\eqref{eq:comparison} yields
\begin{align*}
&\varlimsup_{N\to\infty}\mathbb{P}\Big(\sum_{i\geq m}q_i(G_N, \infty)>\varepsilon\Big)\\
&\hspace{2cm}\leq \varlimsup_{N\to\infty}\mathbb{P}\Big(\sum_{i\geq m}q_i(G_N', \infty)>\varepsilon\Big)
= (1-\lambda)\sum_{i\geq m}\lambda^i.
\end{align*}
Since $\lambda<1$, taking the limit $k\to\infty$, the right side of the above inequality tends to zero, and hence, the condition in~\eqref{eq:smalltail} is satisfied.
\end{proof}

\noindent
\textbf{ Asymptotic behavior of occupancy processes in cliques.}
We now briefly recall the behavior of the occupancy processes on a clique
as the number of servers $N$ grows large.
Rigorous descriptions of the limiting processes are provided in Theorems~\ref{fluidjsqd-ssy} and~\ref{diffusionjsqd-ssy} in Chapter~\ref{chap:univjsqd}.

The behavior on $N$-scale is observed in terms of 
$q_i(G_N, t) = Q_i(G_N, t)/N$ of servers with queue length
at least~$i$ at time~$t$.
When $\lambda<1$, on any finite time interval,
\begin{equation}\label{eq:fluid-conv}
\big\{(q_1(K_N,t), q_2(K_N,t),\ldots)\big\}_{t\geq 0} \dto \big\{(q_1(t), q_2(t),\ldots)\big\}_{t\geq 0},
\end{equation}
as $N\to\infty$,
where $(q_1(\cdot), q_2(\cdot),\ldots)$ is some deterministic process.
Furthermore, in steady state 
\begin{equation}\label{eq:clique-stn}
q_1(K_N,\infty)\pto \lambda \quad \mbox{and} \quad
q_i(K_N,\infty)\pto 0 \ \mbox{ for all } i = 2, \dots, b,
\end{equation}
as  $N\to\infty$.
Note that $q_1(K_N,\cdot)$ is the fraction of non-empty servers. 
Thus $q_1(K_N,\infty)$ is the steady-state scaled departure rate
which should be equal to the scaled arrival rate $\lambda$.
Surprisingly, however, we observe that the steady-state fraction
of servers with a queue length of two or larger is asymptotically
negligible.

To analyze the behavior on $\sqrt{N}$-scale, we consider
a heavy-traffic scenario (i.e., Halfin-Whitt regime) where the
arrival rate at each server is given by $\lambda(N)/N$ with $\lambda(N)$ satisfying~\eqref{eq:HW}.
In order to describe the behavior in the limit, let 
\[\bar{\QQ}(G_N,t) =
\big(\bar{Q}_1(G_N,t), \bar{Q}_2(G_N,t), \dots, \bar{Q}_b(G_N,t)\big)\]
be a properly centered and scaled version of the occupancy process
$\QQ(G_N,t)$, with 
\begin{equation}\label{eq:HWOcc-sig2}
\bar{Q}_1(G_N,t) = - \frac{N-Q_1(G_N,t)}{ \sqrt{N}},\qquad \bar{Q}_i(G_N,t) =\frac{ Q_i(G_N,t)}{\sqrt{N}},
\end{equation}
$i = 2, \dots, b$.
The reason why $Q_1(\cdot,\cdot)$ is centered around~$N$ while $Q_i(\cdot,\cdot)$,
$i = 2, \dots, b$, are not, is because for $G_N=K_N$, the fraction of servers with a queue length of exactly one tends to one, whereas the fraction of servers with a queue length of two or larger tends to zero as $N\to\infty$, as mentioned above.
As mentioned in Section~\ref{ssec:diffjsq} in Chapter~\ref{chap:introduction}, recent results for $\QQ(K_N,t)$~\cite{EG15} show that from a suitable starting state, 
\begin{equation}\label{eq:diff-conv}
\begin{split}
&\big\{(\bQ_1(K_N,t), \bQ_2(K_N,t),\bQ_3(K_N,t),\ldots)\big\}_{t\geq 0}\dto \big\{(\bQ_1(t), \bQ_2(t),0,\ldots)\big\}_{t\geq 0},
\end{split}
\end{equation}
 as $N\to\infty$, where $(\bQ_1(\cdot), \bQ_2(\cdot))$ is some diffusion process.
A precise description of the limiting diffusion process is provided in Theorem~\ref{diffusionjsqd-ssy} in Chapter~\ref{chap:univjsqd}.
This implies that over any finite time interval,
there will be $O_P(\sqrt{N})$ servers with queue length zero
and $O_P(\sqrt{N})$ servers with a queue length of two or larger,
and hence all but $O_P(\sqrt{N})$ servers have a queue length of exactly one.

\vspace{.25cm}
\noindent
\textbf{ Asymptotic optimality.}
From the stochastic optimality of the JSQ policy as mentioned in Section~\ref{spec}, observe that a clique is an optimal load balancing topology, i.e., the occupancy process is better balanced and smaller (in a majorization sense) than in any other graph topology.
In general the optimality is strict, but it turns out that near-optimality can be achieved asymptotically in a broad class of other graph topologies.
Therefore, we now introduce two notions of \emph{asymptotic optimality}, which will be useful to characterize the performance in large-scale systems. 

\begin{definition}[{Asymptotic optimality}]
A graph sequence $\GG = \{G_N\}_{N\geq 1}$  is called `asymptotically optimal on $N$-scale' or `$N$-optimal', if for any $\lambda<1$, on any finite time interval, the scaled occupancy process $(q_1(G_N,\cdot), q_2(G_N,\cdot),\ldots)$ converges weakly to the process $(q_1(\cdot), q_2(\cdot),\ldots)$ given by~\eqref{eq:fluid-conv}.

Moreover, a graph sequence $\GG = \{G_N\}_{N\geq 1}$  is called `asymptotically optimal on $\sqrt{N}$-scale' or `$\sqrt{N}$-optimal', if for any $\lambda(N)$ satisfying~\eqref{eq:HW}, on any finite time interval, the centered scaled occupancy process   $(\bQ_1(G_N,\cdot), \bQ_2(G_N,\cdot),\ldots)$ as in~\eqref{eq:HWOcc-sig2} converges weakly to the process $(\bQ_1(\cdot), \bQ_2(\cdot),\ldots)$ given by~\eqref{eq:diff-conv}.
\end{definition}
\noindent
Intuitively speaking, if a graph sequence is $N$-optimal or $\sqrt{N}$-optimal, then in some sense, the associated occupancy processes are indistinguishable from those of the sequence of cliques on $N$-scale or $\sqrt{N}$-scale.
In other words, on any finite time interval their occupancy processes can differ from those in cliques by at most $o(N)$ or $o(\sqrt{N})$, respectively. 
For brevity, $N$-scale and $\sqrt{N}$-scale are often referred to as \emph{fluid scale} and \emph{diffusion scale}, respectively.
In particular, due to the $\ell_1$-tightness of the scaled
occupancy processes as stated in Proposition~\ref{prop:ergod},
we obtain that for any $N$-optimal graph sequence $\{G_N\}_{N\geq 1}$,
\begin{equation}
q_1(G_N,\infty)\to \lambda \quad \mbox{and} \quad
q_i(G_N,\infty)\to 0 \ \mbox{ for all } i = 2, \dots, b,
\end{equation}
as  $N\to\infty$,
implying that the stationary fraction of servers with queue length
two or larger and the mean waiting time vanish.


\section{Sufficient criteria for asymptotic optimality}\label{sec:det}

In this section we develop a criterion for asymptotic optimality
of an arbitrary deterministic graph sequence on different scales.
In Section~\ref{sec:random} this criterion will be leveraged to establish optimality
of a sequence of random graphs.

We start by introducing some notation, and two measures
of \emph{well-connectedness}.
Let $G=(V,E)$ be any graph.
For a subset $U\subseteq V$, define $\com(U) := |V\setminus N[U]|$
to be the set of all vertices that do not share an edge with any vertex from $U$,
where $N[U]:=U\cup \{v\in V:\ \exists\ u\in U\mbox{ with }(u,v)\in E\}$.
For any fixed $\varepsilon>0$ define
\begin{equation}\label{def:dis-sig2}
\begin{split}
\dis_1(G,\varepsilon) &:= \sup_{U\subseteq V, |U|\geq \varepsilon |V|}\com(U),\\
\dis_2(G,\varepsilon) &:= \sup_{U\subseteq V, |U|\geq \varepsilon \sqrt{|V|}}\com(U).
\end{split}
\end{equation}

The next theorem provides sufficient conditions for asymptotic
optimality on $N$-scale and $\sqrt{N}$-scale in terms of the above
two well-connectedness measures.

\begin{theorem}\label{th:det-seq-sig2}
For any graph sequence $\GG= \{G_N\}_{N\geq 1}$,
\begin{enumerate}[{\normalfont (i)}]
\item $\GG$ is $N$-optimal if for any $\varepsilon>0$, 
$\dis_1(G_N,\varepsilon)/N\to 0$,  as $ N\to\infty.$
\item $\GG$ is $\sqrt{N}$-optimal if for any $\varepsilon>0$, 
$\dis_2(G_N,\varepsilon)/\sqrt{N}\to 0$,  as $ N\to\infty.$
\end{enumerate}
\end{theorem}

The next corollary is an immediate consequence
of Theorem~\ref{th:det-seq-sig2}.

\begin{corollary}
Let $\GG= \{G_N\}_{N\geq 1}$ be any graph sequence and $d_{\min}(G_N)$ be the minimum degree of $G_N$. Then
\textrm{(i)} If $d_{\min}(G_N) = N-o(N)$, then $\GG$ is $N$-optimal, and 
\textrm{(ii)} If $d_{\min}(G_N) = N-o(\sqrt{N})$, then $\GG$ is $\sqrt{N}$-optimal.
\end{corollary}

The rest of the section is devoted to a discussion of the main
proof arguments for Theorem~\ref{th:det-seq-sig2}, focusing on the
proof of $N$-optimality. 
The proof of $\sqrt{N}$-optimality follows along similar lines.
We establish in Proposition~\ref{prop:cjsq-sig2} that if a system is able
to assign each task to a server in the set $\cS^N(n(N))$ of the
$n(N)$ nodes with shortest queues (ties broken arbitrarily), where $n(N)$ is $o(N)$,
then it is $N$-optimal. 
Since the underlying graph is not a clique however (otherwise there
is nothing to prove), for any $n(N)$ not every arriving task can be
assigned to a server in $\cS^N(n(N))$.
Hence we further prove in Proposition~\ref{prop:stoch-ord-new-sig2}
a stochastic comparison property implying that if on any finite
time interval of length~$t$, the number of tasks $\Delta^N(t)$
that are not assigned to a server in $\cS^N(n(N))$ is $o_P(N)$,
then the system is $N$-optimal as well.
The $N$-optimality can then be concluded when $\Delta^N(t)$ is
$o_P(N)$, which we establish in Proposition~\ref{prop:dis-new-sig2}
under the condition that $\dis_1(G_N,\varepsilon)/N\to 0$ as
$N \to \infty$ as stated in Theorem~\ref{th:det-seq-sig2}.

To further explain the idea described in the above proof outline,
it is useful to adopt a slightly different point of view towards
load balancing processes on graphs.
From a high level, a load balancing process can be thought of as follows:
there are $N$ servers, which are assigned incoming tasks by some scheme.
The assignment scheme can arise from some topological structure as considered in this chapter, in which case we will call it \emph{topological load balancing}, or it can arise from some other property of the occupancy process, in which case we will call it \emph{non-topological load balancing}.
As mentioned earlier, under Markovian assumptions, the JSQ policy or the clique is optimal among the set of all non-anticipating schemes, irrespective of being topological or non-topological.
Also, load balancing on graph topologies other than a clique can be thought of as a `sloppy' version of that on a clique, when each server only has access to partial information on the occupancy state.
Below we first introduce a different type of sloppiness in the task assignment scheme, and show that under a limited amount of sloppiness optimality is retained on a suitable scale.
Next we will construct a scheme which is a hybrid of topological and non-topological schemes, whose behavior is simultaneously close to both the load balancing process on a suitable graph and that on a clique.

\vspace{.25cm}
\noindent
\textbf{ A class of sloppy load balancing schemes.}
Fix some function $n:\N\to\N$, and recall the set $\cS^N(n(N))$ as before.
Consider the class $\CJSQ(n(N))$ where each arriving task is assigned
to one of the servers in $\cS^N(n(N))$. 
It should be emphasized that for any scheme in $\CJSQ(n(N))$, we are not imposing any restrictions on how the ties are broken to select the specific set $\cS^N(n(N))$, or how the incoming task should be assigned to a server in $\cS^N(n(N))$.
The scheme only needs to ensure that the arriving task is assigned to \emph{some} server in $\cS^N(n(N))$ with respect to \emph{some} tie breaking mechanism.
The next proposition provides a sufficient criterion for asymptotic
optimality of any scheme in $\CJSQ(n(N))$.

\begin{proposition}\label{prop:cjsq-sig2}
For $0\leq n(N)<N$, let $\Pi\in\CJSQ(n(N))$ be any scheme. \textrm{(i)} If $n(N)/N\to 0$ as $N\to\infty,$ then $\Pi$ is $N$-optimal, and \textrm{(ii)} If $n(N)/\sqrt{N}\to 0$ as $N\to\infty,$ then $\Pi$ is $\sqrt{N}$-optimal.
\end{proposition}

\begin{proof}
We show that if $n(N)/N\to 0$ and $n(N)/\sqrt{N}\to 0$, then any scheme in the class $\CJSQ(n(N))$ has the same process-level limits on $N$-scale and $\sqrt{N}$-scale, respectively.
This establishes the asymptotic optimality on respective scales. 
The idea is similar to the ones used in the proofs of Theorems~\ref{fluidjsqd-ssy} and~\ref{diffusionjsqd-ssy}.

(i) Define $\bar{N}=N-n(N)$ and $\bar{\lambda}(\bar{N})=\lambda(N)$.
Observe that the MJSQ$(n(N))$ scheme with $N$ servers can be thought of as the clique with $\bar{N}$ servers and arrival rate $\bar{\lambda}(\bar{N})/\bar{N}$ per server.
Also, since $n(N)/N\to 0$,
\begin{align*}
\frac{\bar{\lambda}(\bar{N})}{\bar{N}}=\frac{\lambda(N)}{N-n(N)}\to \lambda\quad \text{as}\quad \bar{N}\to\infty.
\end{align*}
Furthermore, observe that the limit of the scaled occupancy processes in Theorem~\ref{fluidjsqd-ssy} as given by~\eqref{eq:fluid} is characterized by the parameter $\lambda$ only, and hence the fluid limit of the MJSQ$(n(N))$ scheme is the same as that of the clique.

Now, observe from the fluid limit of the occupancy processes of cliques  that if $\lambda< 1$, then for any buffer capacity $b\geq 1$, and any starting state, the fluid-scaled cumulative overflow is negligible, i.e., for any $t\geq 0$, $L^N(t)/N\pto 0$, where $L^N(t)$ is the total number of lost tasks up to time $t$.
Since the above fact is induced by the fluid limit only, the same holds for the MJSQ$(n(N))$ scheme.
Therefore, using the lower and upper bounds in Corollary~\ref{cor:bound-ssy} and the tail bound in Proposition~\ref{prop:stoch-ord}, we complete the proof of (i).

(ii) To show that the MJSQ$(n(N))$ scheme has the same diffusion limit as the occupancy processes of cliques if $n(N)/\sqrt{N}\to 0$ as $N\to\infty$, define $\bar{N}=N-n(N)$ and $\bar{\lambda}(\bar{N})=\lambda(N)$.
As mentioned earlier, the MJSQ$(n(N))$ scheme with $N$ servers can be thought of as the clique with $\bar{N}$ servers and arrival rate $\bar{\lambda}(\bar{N})/\bar{N}$ per server.
Also, since $n(N)/\sqrt{N}\to 0$,
\begin{align*}
\frac{\bar{N}-\bar{\lambda}(\sqrt{\bar{N}})}{\bar{N}}=\frac{N-n(N)-\lambda(N)}{\sqrt{N-n(N)}}\to \beta>0\quad \text{as}\quad \bar{N}\to\infty.
\end{align*}
Furthermore, observe that the diffusion limit of the occupancy processes of cliques in \cite[Theorem~2]{EG15} as given in~\eqref{eq:diffusionjsqd-ssy} is characterized by the parameter $\beta>0$, and hence the diffusion limit of the MJSQ$(n(N))$ scheme is the same as that of the occupancy processes of cliques.

Observe from the diffusion limit of the cliques that if $\beta>0$, then for any buffer capacity $b\geq 2$, and suitable initial state as described in~Theorem~\ref{diffusionjsqd-ssy}, the cumulative overflow is negligible, i.e., for any $t\geq 0$, $L^N(t)\pto 0$.
Indeed observe that if $b\geq 2$, and $\big\{\bQ_2^N(0)\big\}_{N\geq 1}$ is a tight sequence, then the sequence of processes $\big\{\bQ_2^N(t)\big\}_{t\geq 0}$ is stochastically bounded.
Therefore, on any finite time interval, there will be only $\Op(\sqrt{N})$ servers with queue length more than one, whereas, for an overflow event to occur all the $N$ servers must have at least two pending tasks.
Therefore, for any $t\geq 0$,
\begin{align*}
\limsup_{N\to\infty}\Pro{L^N(t)>0}&\leq \limsup_{N\to\infty}\Pro{\sup_{s\in[0,t]}Q^N_2(s)=N}\\
&\leq \limsup_{N\to\infty}\Pro{\sup_{s\in[0,t]}\bQ^N_2(s)=\sqrt{N}}=0.
\end{align*} 
Since the above fact is implied by the diffusion limit only, the same holds for the MJSQ$(n(N))$ scheme.
Therefore, using the lower and upper bounds in Corollary~\ref{cor:bound-ssy}, we complete the proof of~(ii).
\end{proof}

\vspace{.25cm}
\noindent
\textbf{ A bridge between topological and non-topological load balancing.}
For any graph $G_N$ and $n \leq N$, we first construct a scheme called $I(G_N, n)$, which is an intermediate blend between the topological load balancing process on $G_N$ and some kind of non-topological load balancing on $N$ servers.
The choice of $n = n(N)$ will be clear from the context.

To describe the scheme $I(G_N,n)$, first synchronize the arrival
epochs at server~$v$ in both systems, $v = 1,2,\ldots,N$.
Further, the servers in both systems are arranged in  non-decreasing order of the queue lengths, and  the departure epochs at the $k$-th ordered
server in the two systems are synchronized, 
$k = 1,2,\ldots,N$.
When a task arrives at server~$v$ at time~$t$ say, it is assigned
in the graph $G_N$ to a server $v' \in N[v]$ according to its own
statistical law.
For the assignment under the scheme $I(G_N,n)$, first observe that if
\begin{equation}\label{eq:criteria-sig2}
\min_{u\in N[v]}X_u(G_N,t) \leq \max_{u\in \cS(n)}X_u(G_N,t),
\end{equation}
then there exists \emph{some} tie-breaking mechanism for which
$v' \in N[v]$ belongs to $\cS(n)$ under $G_N$.
Pick such an ordering of the servers, and assume that $v'$ is the
$k$-th ordered server in that ordering, for some $k \leq n+1$.
Under $I(G_N,n)$ assign the arriving task to the $k$-th ordered
server (breaking ties arbitrarily in this case).
Otherwise, if \eqref{eq:criteria-sig2} does not hold, then the task is
assigned to one of the $n+1$ servers with minimum queue lengths
under $G_N$ uniformly at random.

Denote by $\Delta^N(I(G_N,n),T)$ the cumulative number of arriving tasks
up to time $T \geq 0$ for which Equation~\eqref{eq:criteria-sig2} is
violated under the above coupling.
The next proposition shows that the load balancing process under
the scheme $I(G_N,n)$ is close to that on the graph $G_N$ in terms
of the random variable $\Delta^N(I(G_N,n),T)$.

\begin{proposition}\label{prop:stoch-ord-new-sig2}
The following inequality is preserved almost surely
\begin{equation}\label{eq:stoch-ord-new-sig2}
\sum_{i=1}^b |Q_i(G_N,t)-Q_i(I(G_N, n),t)|\leq 2\Delta^N(I(G_N,n),t)\ \forall\ t\geq 0,
\end{equation}
provided the two systems start from the same occupancy state at $t=0$.
\end{proposition}
\begin{proof}
With the construction of the scheme $I(G_N,n)$, note that when a task arrives at some vertex $v$ say, the load balancing process on $G_N$ and the scheme $I(G_N,n)$ can differ in decision only if none of the vertices in $\cS(n)$ is a neighbor of $v$, i.e., when Equation~\eqref{eq:criteria-sig2} is not satisfied.
Thus Proposition~\ref{prop:stoch-ord2} completes the proof.
\end{proof}

In order to conclude optimality on $N$-scale or $\sqrt{N}$-scale,
it remains to be shown that for any $T>0$, $\Delta^N(I(G_N,n),T)$ is sufficiently small.
The next proposition provides suitable asymptotic bounds for
$\Delta^N(I(G_N,n),T)$ under the conditions on $\dis_1(G_N,\varepsilon)$
and $\dis_2(G_N,\varepsilon)$ stated in Theorem~\ref{th:det-seq-sig2}.

\begin{proposition}\label{prop:dis-new-sig2}
For any $\varepsilon, T>0$ the following holds.
\begin{enumerate}[{\normalfont (i)}]
\item There exist $\varepsilon'>0$ and $n_{\varepsilon'}(N)$ with $n_{\varepsilon'}(N)/N\to 0$ as $N\to\infty$, such that if $\dis_1(G_N,\varepsilon')/N\to 0$ as $N\to\infty$, then  
\[\Pro{\Delta^N(I(G_N,n_{\varepsilon'}),T)/N>\varepsilon}\to 0.\] 
\item There exist $\varepsilon'>0$ and $m_{\varepsilon'}(N)$ with $m_{\varepsilon'}(N)/\sqrt{N}\to 0$ as $N\to\infty$, such that if $\dis_2(G_N,\varepsilon')/\sqrt{N}\to 0$ as $N\to\infty$, then  
\[\Pro{\Delta^N(I(G_N,m_{\varepsilon'}),T)/\sqrt{N}>\varepsilon}\to 0.\]
\end{enumerate}
\end{proposition}

The proof of Theorem~\ref{th:det-seq-sig2} then readily follows
by combining Propositions~\ref{prop:cjsq-sig2}-\ref{prop:dis-new-sig2}
and observing that the scheme $I(G_N,n)$ belongs to the class
$\CJSQ(n)$ by construction.

\begin{proof}[Proof of Proposition~\ref{prop:dis-new-sig2}]
Fix any $\varepsilon, T>0$ and choose
$\varepsilon' = \varepsilon/(2 \lambda T)$.
With the coupling described above, when a task arrives at some
vertex~$v$ say, Equation~\eqref{eq:criteria-sig2} is violated only if none
of the vertices in $\cS(n_{\varepsilon'}(N))$ is a neighbor of~$v$.
Thus, the total instantaneous rate at which this happens is
\[
\lambda \com(\cS(n_{\varepsilon'}(N),t)) \leq
\lambda \sup_{U \subseteq V_N, |U| \geq n_{\varepsilon'}(N)} \com(U),
\]
irrespective of what this set $\cS^N(n(N))$ actually is.
Therefore, for any fixed $T \geq 0$,
\[
\Delta^N(I(G_N,n_{\varepsilon'}),T) \leq
A\Big(\lambda \sup_{U \subseteq V_N, |U| \geq n_{\varepsilon'}(N)} \com(U)\Big),
\]
where $A(\cdot)$ represents a unit-rate Poisson process.
This can then be leveraged to show that
$\Delta^N(I(G_N,n_{\varepsilon'}),T)$ is small on an $N$-scale
and $\sqrt{N}$-scale, respectively, under the conditions stated
in the proposition, by choosing a suitable $n_{\varepsilon'}$.

Specifically, if $\dis_1(G_N,\varepsilon')/N \to 0$, then there exists
$n_{\varepsilon'}(N)$ with $n_{\varepsilon'}(N)/N \to 0$ such that 
$\dis_1(G_N,\varepsilon') \leq n_{\varepsilon'}(N)$ for all $N \geq 1$,
and hence
$$\sup_{U \subseteq V_N, |U| \geq n_{\varepsilon'}(N)} \com(U) \leq
\varepsilon' N.$$
It then follows that with high probability,
\[
\limsup_{N \to\infty} \frac{1}{N} \Delta^N(I(G_N,n_{\varepsilon'}),T) \leq
\limsup_{N \to\infty} \frac{1}{N} A\Big(\lambda T \varepsilon' N\Big) \leq
2 \lambda T \varepsilon' = \varepsilon.
\]

Likewise, if $\dis_2(G_N,\varepsilon')/\sqrt{N} \to 0$, then there exists
$m_{\varepsilon'}(N)$ with $m_{\varepsilon'}(N)/\sqrt{N}\to 0$ such that
$\dis_2(G_N,\varepsilon') \leq m_{\varepsilon'}(N)$ for all $N \geq 1$,
and hence
$$\sup_{U \subseteq V_N, |U| \geq m_{\varepsilon'}(N)} \com(U) \leq
\varepsilon' \sqrt{N}.$$
It then follows that with high probability,
\begin{align*}
&\limsup_{N\to\infty} \frac{1}{\sqrt{N}} \Delta^N(I(G_N,m_{\varepsilon'}),T) 
\leq
\limsup_{N\to\infty} \frac{1}{\sqrt{N}} A\Big(\lambda T \varepsilon' \sqrt{N}\Big) \leq
2 \lambda T \varepsilon' = \varepsilon.
\end{align*}
\end{proof}

\begin{proof}[Proof of Theorem~\ref{th:det-seq-sig2}]

(i) In order to prove the fluid-level optimality of $G_N$, fix any $\varepsilon>0$. 
Observe from Proposition~\ref{prop:stoch-ord-new-sig2} and Proposition~\ref{prop:dis-new-sig2} (i) that there exists $\varepsilon'>0$ such that with high probability
\begin{align*}
 \sup_{t\in[0,T]}\frac{1}{N}\sum_{i=1}^b|Q_i(G_N,t)-Q_i(I(G_N,n_{\varepsilon'}(N)),t)| 
\leq \frac{2\Delta^N_\varepsilon(T)}{N}\leq \varepsilon.
\end{align*}
Furthermore, since $I(G_N,n_{\varepsilon'}(N))\in\CJSQ(n_{\varepsilon'}(N))$ and $n_{\varepsilon'}(N)/N\to 0$, Proposition~\ref{prop:cjsq-sig2} yields
$$\sup_{t\in[0,T]}\sum_{i=1}^b|q_i(I(G_N,n_{\varepsilon'}(N)),t)-q_i(t)|\pto 0\quad \mathrm{as}\quad N\to\infty.$$
Thus since $\varepsilon>0$ is arbitrary, we obtain that with high probability as $N\to\infty$,
$$\sup_{t\in[0,T]}\sum_{i=1}^b|q_i(G_N,t)-q_i(t)|\leq \varepsilon'',$$
for all $\varepsilon''>0$, which completes the proof of Part (i).\\

(ii) To prove the diffusion-level optimality of $G_N$, again fix any $\varepsilon>0$. 
As in Part (i), using Proposition~\ref{prop:stoch-ord-new-sig2} and Proposition~\ref{prop:dis-new-sig2} (ii), there exists $\varepsilon'>0$
\begin{align*}
\sup_{t\in[0,T]}\frac{1}{\sqrt{N}}\sum_{i=1}^b|Q_i(G_N,t)-Q_i(I(G_N,m_{\varepsilon'}(N)),t)|\leq \frac{\Delta^N_{\varepsilon'}(T)}{\sqrt{N}}\leq \varepsilon.
\end{align*}
Furthermore, since $I(G_N,m_{\varepsilon'}(N))\in\CJSQ(m_{\varepsilon'}(N))$ and $m_{\varepsilon'}(N)/\sqrt{N}\to 0$, Proposition~\ref{prop:cjsq-sig2} yields
\begin{align*}
&\big\{(\bQ_1(I(G_N,m_{\varepsilon'}(N)),t), \bQ_2(I(G_N,m_{\varepsilon'}(N)),t),\ldots)\big\}_{t\geq 0}\\
&\hspace{3cm}\dto \big\{(\bQ_1(t), \bQ_2(t),\ldots)\big\}_{t\geq 0},
\end{align*}
as $N\to\infty$, where the process $(\bQ_1(\cdot), \bQ_2(\cdot),\ldots)$ given by~\eqref{eq:diff-conv}.
Since $\varepsilon>0$ is arbitrary, we thus obtain
$$\big\{(\bQ_1(G_N,t), \bQ_2(G_N,t),\ldots)\big\}_{t\geq 0} \dto \big\{(\bQ_1(t), \bQ_2(t),\ldots)\big\}_{t\geq 0},$$
as $N\to\infty$, which completes the proof of Part (ii).
\end{proof}

\section{Necessary criteria for asymptotic optimality}
\label{sec:necessary}
From the conditions of Theorem~\ref{th:det-seq-sig2} it follows that if for all $\varepsilon>0$, $\dis_1(G_N,\varepsilon)$ and $\dis_2(G_N,\varepsilon)$ are $o(N)$ and $o(\sqrt{N})$, respectively, then
the total number of edges in $G_N$ must be $\omega(N)$ and $\omega(N\sqrt{N})$, respectively.
Theorem~\ref{th:bdd-deg-sig2} below states that the \emph{super-linear} growth rate of the total number of edges is not only sufficient, but also necessary in the sense that any graph with $O(N)$ edges is asymptotically sub-optimal on $N$-scale.

\begin{theorem}
\label{th:bdd-deg-sig2}
Let $\GG= \{G_N\}_{N\geq 1}$ be any graph sequence, such that there exists a fixed integer $M<\infty$ with 
\begin{equation}\label{eq:bdd-deg-sig2}
\limsup_{N\to\infty}\dfrac{\#\big\{v\in V_N:d_v\leq M\big\}}{N}>0,
\end{equation}
where $d_v$ is the degree of the vertex $v$.
Then $\GG$ is sub-optimal on $N$-scale.
\end{theorem}

\begin{proof}
For brevity, denote by $\Xi_N(M)\subseteq V_N$ the set of all vertices with degree at most $M$.
Since $|\Xi_N(M)|/N\leq 1,$ from~\eqref{eq:bdd-deg-sig2} we have a convergent subsequence $\big\{\Xi_{N_n}(M)\big\}_{n\geq 1}$ with $\{N_n\}_{n\geq 1}\subseteq \N$, such that $|\Xi_{N_n}(M)|/N\to \xi>0$, as $N\to\infty$.
For the rest of the proof we will consider the asymptotic statements along this subsequence, and hence omit the subscript $n$.

Let the system start from an occupancy state where all the vertices in $\Xi_N(M)$ are empty.
We will show that in finite time, a positive fraction of vertices in $\Xi_N(M)$ will have at least two tasks. 
This will prove that the fluid limit sample path cannot agree with that of the sequence of cliques, and hence $\{G_N\}_{N\geq 1}$ cannot be $N$-optimal.
The idea of the proof is as follows: If a graph contains $\Theta(N)$ bounded degree vertices, then starting from all empty servers, in any finite time interval there will be $\Theta(N)$ servers $u$ say, for which all the servers in $N[u]$ have at least one task.
For all such servers an arrival at $u$ must produce a server of queue length two.
Thus, it shows that the instantaneous rate at which servers of queue length two are formed is bounded away from zero, and hence $\Theta(N)$ servers of queue length two are produced in finite time. 

Let $u$ be a vertex with degree $M$ or less in $G_N$.
Consider the event $\cE_N(u,t)$ that at time~$t$ all vertices in $N[u]$ have at least one job.
Note that since $M<\infty$ is fixed, for any $t>0$, $\Pro{\cE_N(u,t)}\geq\delta(t)$ for some $\delta(t)>0$, for all $N\geq 1$.
To see this, note that $\delta(t)$ is the probability that before time $t$ there are $M+1$ arrivals at vertex $u$ and no departure has taken place.
Also observe that for two vertices $u,v\in V_N$ with degrees at most $M$, 
\begin{equation}
\Pro{\cE_N(u,t)\cap \cE_N(v,t)}\geq \delta(t)^2.
\end{equation}
Indeed the probability of the event $\cE_N(u,t)\cap \cE_N(v,t)$ can be lower bounded by the probability of the event that before time $t$ there are $M+1$ arrivals at vertex $u$, $M+1$ arrivals at vertex $v$, and no departure has taken place from $N[u]\cup N[v]$.
Thus, at time $t$, the fraction of vertices in $\Xi_N(M)$ for which all the neighboring vertices have at least one task, is lower bounded by $\delta(t)$.
Now the proof is completed by considering the following: let $u$ be a vertex of degree $M<\infty$ for which all the neighbors have at least one task.
Then at such an instance if a task arrives at server $u$, it must be assigned to a server with queue length one, and hence a server with queue length two will be formed.
Therefore the total scaled instantaneous rate at which the number of queue length two is being formed at time $t$ is at least $\lambda \delta(t)>0$, which also gives the total rate of increase of the fraction of vertices with at least two tasks. 
\end{proof}

\vspace{.25cm}
\noindent
\textbf{ Worst-case scenario.}
Next we consider the worst-case scenario.
Theorem~\ref{th:min-deg-negative-sig2} below asserts that a graph sequence can be sub-optimal for some $\lambda<1$ even when the minimum degree  $d_{\min}(G_N)$ is~$\Theta(N)$.

\begin{theorem}
\label{th:min-deg-negative-sig2}
For any $\big\{d(N)\big\}_{N\geq 1}$, such that $d(N)/N\to c$ with $0<c<1/2$, there exists $\lambda<1$, and a graph sequence $\GG= \{G_N\}_{N\geq 1}$ with $d_{\min}(G_N)= d(N)$, such that $\GG$ is sub-optimal on $N$-scale.
\end{theorem}

To construct such a sub-optimal graph sequence, consider a sequence of complete bipartite graphs $G_N = (V_N, E_N)$, with $V_N = A_N \sqcup B_N$ and $|A_N|/N\to c\in (0,1/2)$ as $N\to\infty$.
If this sequence were $N$-optimal, then starting from all empty servers, asymptotically 
the fraction of servers with queue length one would converge to $\lambda$, and the fraction of servers with queue length two or larger should remain zero throughout.
Now note that for large $N$ the rate at which tasks join the empty servers in $A_N$ is given by $(1-c)\lambda$, whereas the rate of empty server generation in $A_N$ is at most $c$.
Choosing $\lambda>c/(1-c)$, one can see that in finite time each server in $A_N$ will have at least one task.
From that time onward with at least instantaneous rate $\lambda(\lambda - c) - c$, servers with queue length two start forming.
The range for $c$ stated in Theorem~\ref{th:min-deg-negative-sig2} is only to ensure that there exists $\lambda<1$ with $\lambda(\lambda - c) - c>0$. 
\begin{proof}[{Proof sketch of Theorem~\ref{th:min-deg-negative-sig2}}]
Fix a $c>0$.
Construct the graph sequence $\big\{G_N\big\}_{N\geq 1}$ as a sequence of complete bipartite graphs with size of one partite set of the $N$-th graph to be $\lceil c N\rceil$, i.e., $V_N = A_N \sqcup B_N$, such that $|A_N| = \lceil c N\rceil$ and $B_N = V_N\setminus A_N$, and the edge set is given by $E_N = \big\{(u,v): u\in A_N, v\in B_N\big\}$.
Note that $d_{\min}(G_N)/N\to c$, as $N\to\infty$.
We will show that for any $0<c<1/2$, there exists $\lambda$, such that $\GG$ is sub-optimal on $N$-scale.

Assume on the contrary that $\GG$ is $N$-optimal.
Denote by $Q_{i,A}^N(t)$ and $Q_{i,B}^N(t)$ the number of vertices with at least $i$ tasks in partite sets $A_N$ and $B_N$, respectively.
Also define $q_{i,A}^N(t) = Q_{i,A}^N(t)/N$ and $q_{i,B}^N(t) = Q_{i,B}^N(t)/N$.
Assume $q_{2,A}^N(0)=0$, for all $N$.
Observe that as long as $c - q_{1,A}^N>0$ by a non-vanishing margin, any external arrival to servers in $B_N$ will be assigned to an empty server in $A_N$ with probability $1-\text{O}(1/N)$.
Similarly, as long as $1-c - q_{1,B}^N>0$ by a non-vanishing margin, any external arrival to servers in $A_N$ will be assigned to an empty server in $B_N$ with probability $1-\text{O}(1/N)$.
Thus one can show that as $N\to\infty$, until $q_{1,A}^N$ hits $c$, the processes $\big\{q_{1,A}^N(t)\big\}$ and $\big\{q_{2,B}^N(t)\big\}$ converges weakly to a deterministic process described by the following set of ODE's:
\begin{equation}\label{eq:fluid-counter}
\begin{split}
q_{1,A}'(t) &= \lambda (1 - c) - q_{1,A}(t),\\
q_{1,B}'(t) &= \lambda c - q_{1,B}(t).
\end{split}
\end{equation}
Since the total scaled arrival rate into the system of $N$ servers is $\lambda$, should the above system follow the fluid-limit trajectory of the occupancy process for a clique, starting from an all-empty state, $q_{1,A}(t) + q_{1,B}(t)$ must approach $\lambda$ as $t\to\infty$, and $q_{i, A}(t)$ and $q_{i,B}(t)$ both remain 0 for all $t\geq 0$, $i\geq 2$.
When $\lambda>c/(1-c)$,~\eqref{eq:fluid-counter} implies that in finite time $q_{1,A}(t)$ hits $c$.
Consequently, $q_{1,B}(t)$ should approach $\lambda - c$ as $t\to \infty$.
Now we claim that when $q_{1,A}(t) =c$, if a task appears at a server $v$ in $B_N$ that has queue length one, then with probability $1-\text{O}(1/N)$, it will be assigned to a server in $A_N$.
To see this, note that at such an arrival if there is an empty server in $A_N$, then the arriving task is clearly assigned to the idle server, otherwise, when there is no empty server in $A_N$, the arriving task is assigned uniformly at random among the vertices in $N[v]$ having queue length one.
Since there are $\Theta(N)$ vertices in $A_N$ with queue length one, the arriving task with probability $1-O(1/N)$ joins a server in $A_N$.
Therefore, the total scaled rate of tasks arriving at the servers in $A_N$ is at least $\lambda(\lambda-c)$, whereas the total scaled rate at which tasks can leave from servers in $A_N$ is at most $c$.
Thus if $\lambda(\lambda-c)>c$, then in finite time, a positive fraction of servers in $A_N$ will have queue length two or larger.
Now observe that
$$\lambda(\lambda-c)>c\implies \lambda> \frac{c+\sqrt{c^2+4c}}{2},$$
and $(c+\sqrt{c^2+4c})/2<1$ for any $c\in (0,1/2)$. This completes the proof of Theorem~\ref{th:min-deg-negative-sig2}.
\end{proof}

\section{Asymptotically optimal random graph topologies}
\label{sec:random}
In this section we use Theorem~\ref{th:det-seq-sig2} to investigate how the load balancing process behaves on random graph topologies. 
Specifically, we aim to understand what types of graphs are asymptotically optimal in the presence of randomness (i.e., in the average case scenario).
Theorem~\ref{th:inhom-sig2} below establishes sufficient conditions for asymptotic optimality of a sequence of inhomogeneous random graphs.
Recall that a graph $G' = (V',E')$ is called a supergraph of $G=(V,E)$ if $V=V'$ and $E\subseteq E'$.

\begin{theorem}
\label{th:inhom-sig2}
Let $\GG= \{G_N\}_{N\geq 1}$ be a graph sequence such that for each $N$, $G_N = (V_N, E_N)$ is a supergraph of the inhomogeneous random graph $G_N'$ where any two vertices $u, v\in V_N$ share an edge with probability $p_{uv}^N$.
\begin{enumerate}[{\normalfont (i)}]
\item If $\inf\ \{p^N_{uv}: u, v\in V_N\}$ is $\omega(1/N)$, then $\GG$ is $N$-optimal.
\item If $\inf\ \{p^N_{uv}: u, v\in V_N\}$ is $\omega(\log(N)/\sqrt{N})$, then $\GG$ is $\sqrt{N}$-optimal.
\end{enumerate}
\end{theorem}

The proof of Theorem~\ref{th:inhom-sig2} relies on Theorem~\ref{th:det-seq-sig2}.
Specifically, if $G_N$ satisfies conditions~(i) and~(ii) in
Theorem~\ref{th:inhom-sig2}, then the corresponding conditions~(i) and~(ii)
in Theorem~\ref{th:det-seq-sig2} hold.
\begin{proof}[{Proof of Theorem~\ref{th:inhom-sig2}}]
In this proof we will verify the conditions stated in Theorem~\ref{th:det-seq-sig2} for fluid and diffusion level optimality.
Fix any $\varepsilon>0$.

(i) Observe that for $G_N = (V_N, E_N)$ as described in Theorem~\ref{th:inhom-sig2} (i), we have $p(N):=\inf\ \{p^N_{uv}: u,v\in V_N\}$ with $Np(N)\to \infty$ as $N\to\infty$. 
For any two subsets $V_1$, $V_2\subseteq V_N$, denote by $E_N(V_1, V_2)$ the number of cross-edges between $V_1$ and $V_2$.
Now, for any function $n:\N\to\N$, 
\begin{equation}\label{eq:fluid-bound}
\begin{split}
&\Pro{\exists\ V_1, V_2\subseteq V_N:~|V_1|\geq\varepsilon N,\ |V_2|\geq n(N), E_N(V_1,V_2) = 0}\\&\\
&= \Pro{\exists\ V_1, V_2\subseteq V_N:|V_1|=\varepsilon N,\ |V_2|= n(N), E_N(V_1,V_2) = 0}\\&\\
&\leq {N(1-\varepsilon) \choose \varepsilon N} {N-2\varepsilon N \choose n(N)}(1-p(N))^{\varepsilon N n(N)}\\&\\
&\lesssim \frac{1}{[\varepsilon^\varepsilon(1-\varepsilon)^{1-\varepsilon}]^N}\times\frac{\big(\frac{N}{n(N)}\big)^{n(N)}}{\big(1-\frac{n(N)}{N(1-\varepsilon)}\big)^{N(1-\varepsilon)}}\times \exp(-\varepsilon Np(N) n(N))\\&\\
&\lesssim \dfrac{\exp(-\varepsilon Np(N)n(N))\times \exp(n(N)\log(N))}{\exp(N\log [\varepsilon^\varepsilon(1-\varepsilon)^{1-\varepsilon}])\exp(-n(N))},
\end{split}
\end{equation}
where the first equality is due to the fact that if there are two  sets of vertices $V_1$ and $V_2$ with  $|V_1|\geq\varepsilon N$ and $|V_2|\geq n(N)$, such that there is no edge between $V_1$ and $V_2$, then the graph must contain two sets $V_1'$ and $V_2'$ of sizes exactly equal to $\varepsilon N$ and $n(N)$, respectively, such that there is no edge between $V_1'$ and $V_2'$, and vice-versa.
Choosing $n(N) = N/\sqrt{Np(N)}$ say, it can be seen that for any $p(N)$ such that $Np(N)\to\infty$ as $N\to\infty$, $n(N)/N\to 0$ and the above probability goes to 0. Therefore for any $\varepsilon, \delta>0$, \eqref{eq:fluid-bound} yields
\begin{align*} 
&\Pro{\dis_1(G_N,\varepsilon)>\delta N}\leq \Pro{\exists\ U\subseteq V_N:\ |U|\geq \varepsilon N \mbox{ and }\com(U)\geq \delta N}\to 0,
\end{align*}
as $N\to\infty$.

(ii) Again, for $G_N = (V_N, E_N)$ as described in Theorem~\ref{th:inhom-sig2} (i), we have $p(N):=\inf\ \{p^N_{uv}: u,v\in V_N\}$ with $Np(N)/(\sqrt{N}\log(N))\to \infty$ as $N\to\infty$.  
Now as in Part (i), for any function $n:\N\to\N$,
\begin{equation}\label{eq:diff-bound}
\begin{split}
&\Pro{\exists\ V_1, V_2\subseteq V_N:|V_1|\geq\varepsilon \sqrt{N},\ |V_2|\geq n(N), E_N(V_1,V_2) = 0}
\\&\\
&= \Pro{\exists\ V_1, V_2\subseteq V_N:|V_1|=\varepsilon \sqrt{N},\ |V_2|= n(N), E_N(V_1,V_2) = 0}
\\&\\
&\leq {N-\varepsilon\sqrt{N} \choose \varepsilon \sqrt{N}} {N-2\varepsilon\sqrt{N} \choose n(N)}(1-p(N))^{\varepsilon\sqrt{N} n(N)}
\\&\\
&\lesssim N^{\varepsilon\sqrt{N}/2}\exp(\varepsilon\sqrt{N}) \times
N^{n(N)}\times \exp(-\varepsilon\sqrt{N}p(N)n(N))\\
&\hspace{1.5cm} \times\exp\Big(\frac{-\varepsilon n(N)}{\sqrt{N}}+n(N)\Big(1-\frac{n(N)}{N-\varepsilon\sqrt{N}}\Big)\Big).
\end{split}
\end{equation}
Choosing $n(N) = \sqrt{N}/\sqrt{\sqrt{N}p(N)/\log(N)}$, it can be seen that $n(N)/\sqrt{N}\to 0$ as $N\to\infty$ and the above probability converges to 0.
 Therefore for any $\varepsilon, \delta>0$, \eqref{eq:diff-bound} yields
\begin{align*} 
&\Pro{\dis_2(G_N,\varepsilon)>\delta \sqrt{N}}
\leq \Pro{\exists\ U\subseteq V_N:\ |U|\geq \varepsilon \sqrt{N} \mbox{ and }\com(U)\geq \delta \sqrt{N}}, 
\end{align*}
which converges to 0 as $N\to\infty$.
This completes the proof of Theorem~\ref{th:inhom-sig2}.
\end{proof}

As an immediate corollary to Theorem~\ref{th:inhom-sig2} we obtain an optimality result for the sequence of Erd\H{o}s-R\'enyi random graphs.

\begin{corollary}\label{cor:errg-sig2}
Let $\GG= \{G_N\}_{N\geq 1}$ be a graph sequence such that for each $N$,
$G_N$ is a supergraph of $\ER_N(p(N))$, and $d(N) = (N-1)p(N)$. Then
{\normalfont (i)}
If $d(N)\to\infty$ as $N\to\infty$, then $\GG$ is $N$-optimal.
{\normalfont (ii)}
If $d(N)/(\sqrt{N}\log(N))\to\infty$ as $N\to\infty$, then $\GG$ is $\sqrt{N}$-optimal.
\end{corollary}

Theorem~\ref{th:det-seq-sig2} can be further leveraged to establish the
optimality of the following sequence of random graphs.
For any $N \geq 1$ and $d(N) \leq N-1$ such that $N d(N)$ is even,
construct the \emph{erased random regular} graph on $N$~vertices
as follows:
Initially, attach $d(N)$ \emph{half-edges} to each vertex. 
Call all such half-edges \emph{unpaired}.
At each step, pick one half-edge arbitrarily, and pair it to another
half-edge uniformly at random among all unpaired half-edges
to form an edge, until all the half-edges have been paired.
This results in a uniform random regular multi-graph with degree
$d(N)$~\cite[Proposition 7.7]{remco-book-1}. 
Now the erased random regular graph is formed by erasing all the
self-loops and multiple edges, which then produces a simple graph.

\begin{theorem}
\label{th:reg-sig2}
Let $\GG= \{G_N\}_{N\geq 1}$ be a sequence of erased random regular
graphs with degree $d(N)$. Then
{\normalfont (i)}
If $d(N)\to\infty$ as $N\to\infty$, then $\GG$ is $N$-optimal.
{\normalfont (ii)}
If $d(N)/(\sqrt{N}\log(N))\to\infty$ as $N\to\infty$, then $\GG$ is $\sqrt{N}$-optimal.
\end{theorem}
\begin{proof}[Proof of Theorem~\ref{th:reg-sig2}]
We will again verify the conditions stated in Theorem~\ref{th:det-seq-sig2} for fluid and diffusion level optimality.
For $k\geq 1$, denote $(2k-1)!! = (2k-1)(2k-3)\cdots3\cdot1.$
Fix any $\varepsilon>0$.

(i) For any function $n:\N\to\N$,
\begin{equation}\label{eq:fluid-bound-reg}
\begin{split}
&\Pro{\exists\ V_1, V_2\subseteq V_N:|V_1|\geq\varepsilon N,\ |V_2|\geq n(N), E_N(V_1,V_2) = 0}\\&\\
&= \Pro{\exists\ V_1, V_2\subseteq V_N: |V_1|=\varepsilon N,\ |V_2|= n(N), E_N(V_1,V_2) = 0}\\&\\
&\leq {N \choose \varepsilon N} {N-\varepsilon N \choose n(N)}\frac{(Nd(N)(1-\varepsilon)-1)!!}{(Nd(N)-1)!!}\times\frac{(Nd(N)-n(N)d(N)-1)!!}{(Nd(N)(1-\varepsilon)-n(N)d(N)-1)!!}\\&\\
&\lesssim \frac{1}{[\varepsilon^\varepsilon(1-\varepsilon)^{1-\varepsilon}]^N}\times\frac{\big(\frac{N}{n(N)}\big)^{n(N)}}{\big(1-\frac{n(N)}{N(1-\varepsilon)}\big)^{N(1-\varepsilon)}}\times \exp(-\varepsilon n(N)d(N))\\&\\
&\lesssim \dfrac{\exp(-\varepsilon d(N)n(N))\times \exp(n(N)\log(N))}{\exp(N\log [\varepsilon^\varepsilon(1-\varepsilon)^{1-\varepsilon}])\exp(-n(N))}.
\end{split}
\end{equation}
Choosing $n(N) = N/\sqrt{d(N)}$ say, it can be seen that for any $p(N)$ such that $d(N)\to\infty$ as $N\to\infty$, $n(N)/N\to 0$ and the above probability goes to 0. Therefore for any $\varepsilon, \delta>0$, \eqref{eq:fluid-bound-reg} yields
\begin{align*} 
&\Pro{\dis_1(G_N,\varepsilon)>\delta N}
\leq \Pro{\exists\ U\subseteq V_N:\ |U|\geq \varepsilon N \mbox{ and }\com(U)\geq \delta N}\to 0,
\end{align*}
 as $N\to\infty$.

(ii) Again, as in Part (i), for any function $n:\N\to\N$,
\begin{equation}\label{eq:diff-bound-reg}
\begin{split}
&\Pro{\exists\ V_1, V_2\subseteq V_N:|V_1|\geq\varepsilon \sqrt{N},\ |V_2|\geq n(N), E_N(V_1,V_2) = 0}
\\&\\
&= \Pro{\exists\ V_1, V_2\subseteq V_N: |V_1|=\varepsilon \sqrt{N},\ |V_2|= n(N) E_N(V_1,V_2) = 0}
\\&\\
&\leq {N \choose \varepsilon \sqrt{N}} {N-\varepsilon\sqrt{N} \choose n(N)}\frac{(Nd(N)-\varepsilon\sqrt{N}d(N)-1)!!}{(Nd(N)-1)!!}\\
&\hspace{4cm}\times\frac{(Nd(N)-n(N)d(N)-1)!!}{(Nd(N)-\varepsilon\sqrt{N}d(N)-n(N)d(N)-1)!!}
\\&\\
&\lesssim \exp\Big(\frac{\varepsilon\sqrt{N}\log(N)}{2}-\frac{n(N)d(N)}{\sqrt{N}}\Big).
\end{split}
\end{equation}
Now, choosing $n(N) = \sqrt{N}/\sqrt{d(N)/(\sqrt{N}\log(N))}$, it can be seen that as $N\to\infty,$ $n(N)/\sqrt{N}\to 0$ and the above probability converges to 0.
 Therefore for any $\varepsilon, \delta>0$, \eqref{eq:diff-bound-reg} yields
\begin{align*} 
&\Pro{\dis_2(G_N,\varepsilon)>\delta \sqrt{N}}
\leq \Pro{\exists\ U\subseteq V_N:\ |U|\geq \varepsilon \sqrt{N} \mbox{ and }\com(U)\geq \delta \sqrt{N}},
\end{align*}
which converges to 0 as $N\to\infty$.
\end{proof}

Note that due to Theorem~\ref{th:bdd-deg-sig2}, we can conclude that the growth rate condition on degrees for $N$-optimality in Corollary~\ref{cor:errg-sig2}~(i) and Theorem~\ref{th:reg-sig2}~(i) is not only sufficient, but necessary as well.
Thus informally speaking, $N$-optimality is achieved under the minimum condition required as long as the underlying topology is suitably random.

\begin{figure}
\begin{center}
\includegraphics[width=100mm]{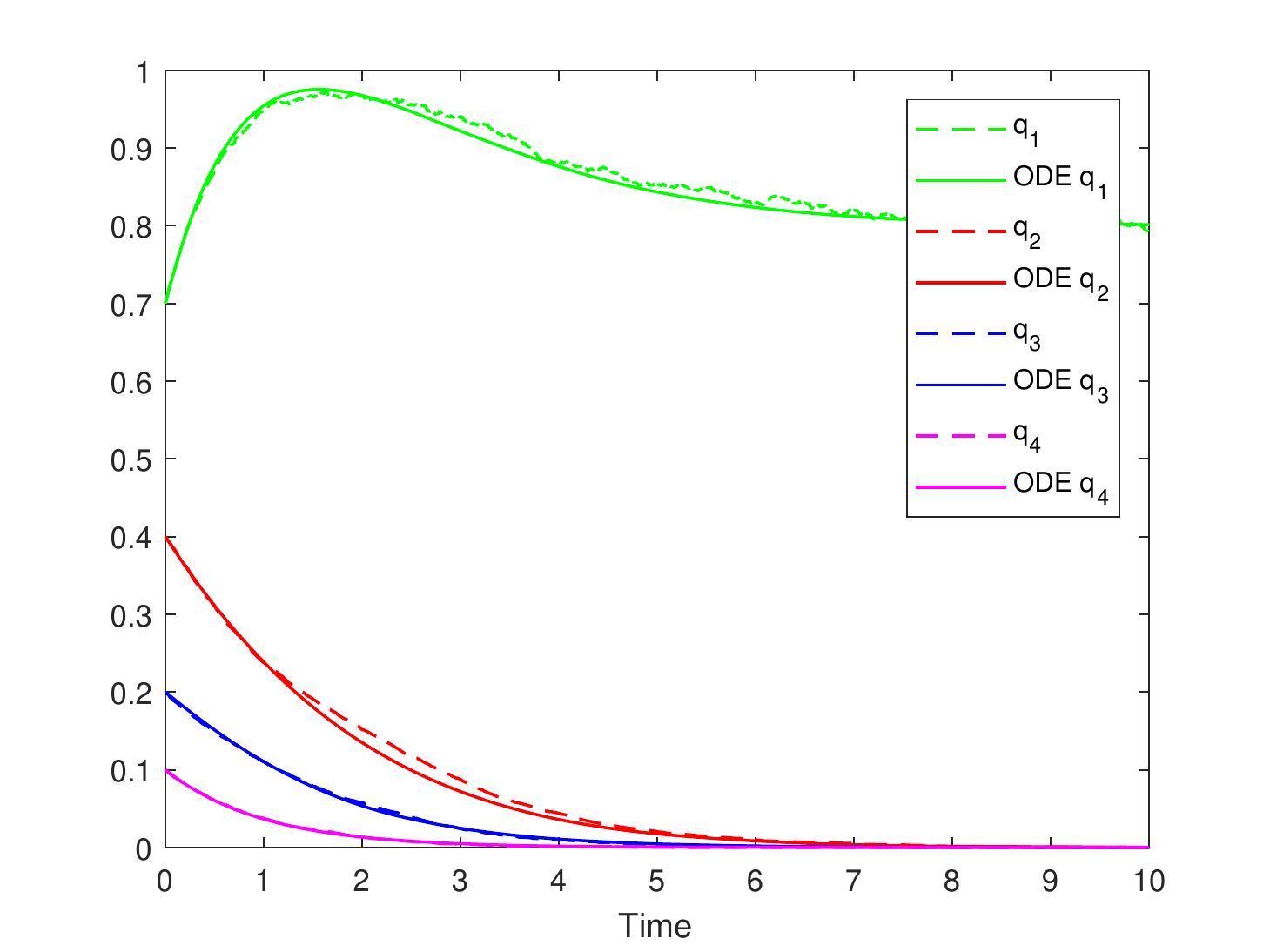}
\end{center}
\caption{Illustration of the fluid-limit trajectories for $\lambda=0.8$ along with a simulation for $N = 10^4$ servers. 
The topology is a single instance of the ERRG on $N=10^4$ vertices with edge probability $1/\sqrt{N}=10^{-2}$, i.e.~the average degree is~100.}
\label{fig:trajectory}
\end{figure}
\begin{figure}
\begin{center}
\includegraphics[width=100mm]{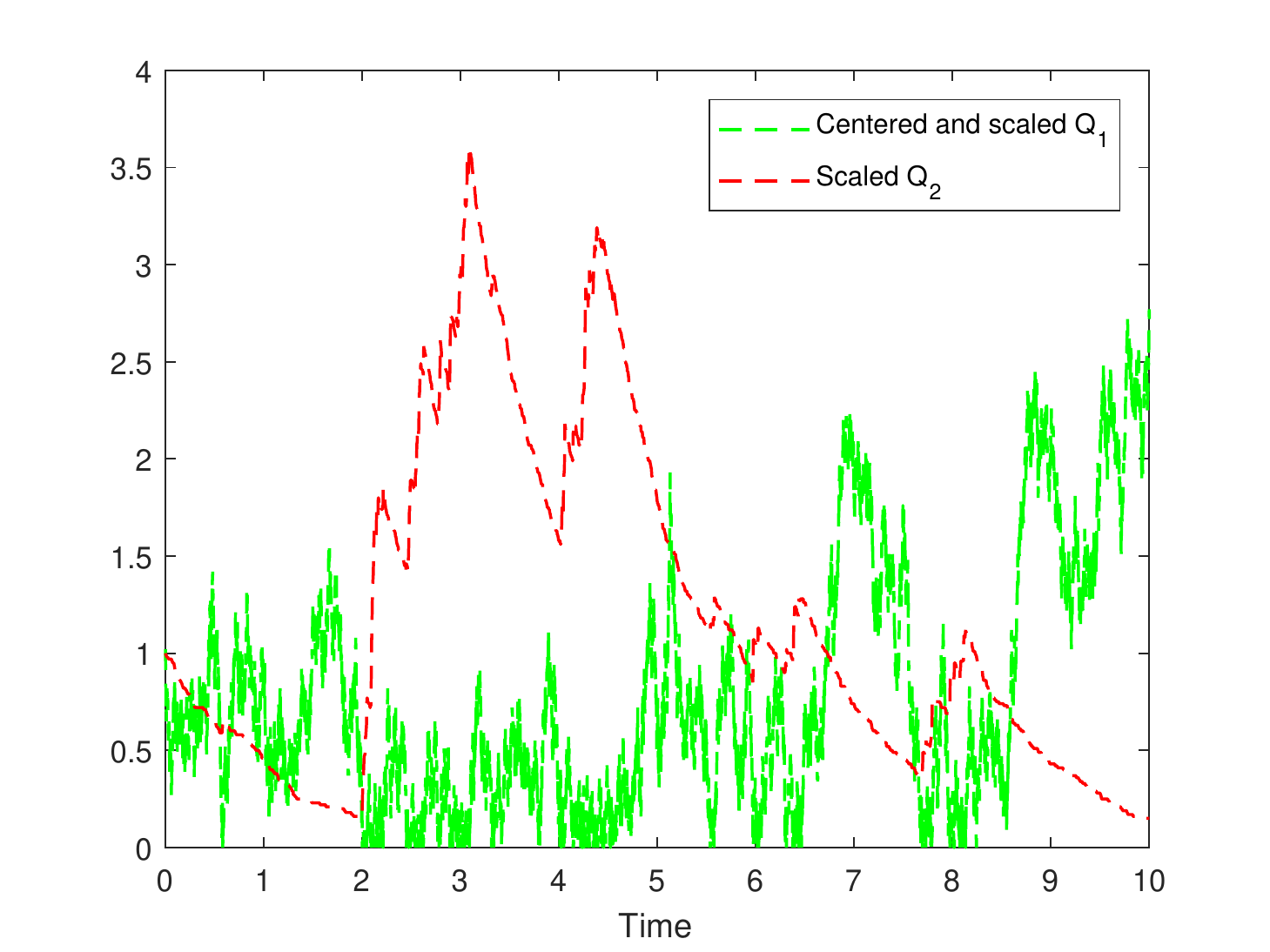}
\end{center}
\caption{Illustration of the diffusion-scaled trajectories in the Halfin-Whitt heavy-traffic regime, for $N = 10^4$ servers and $\lambda(N)= N - \sqrt{N}=9900$.
The topology is a single instance of the ERRG on $N=10^4$ vertices with edge probability $\log(N)^2/\sqrt{N}=0.8483$, i.e.~the average degree is 8483.
}
\label{fig:ht}
\end{figure}
\begin{figure}
\begin{center}
\includegraphics[width=100mm]{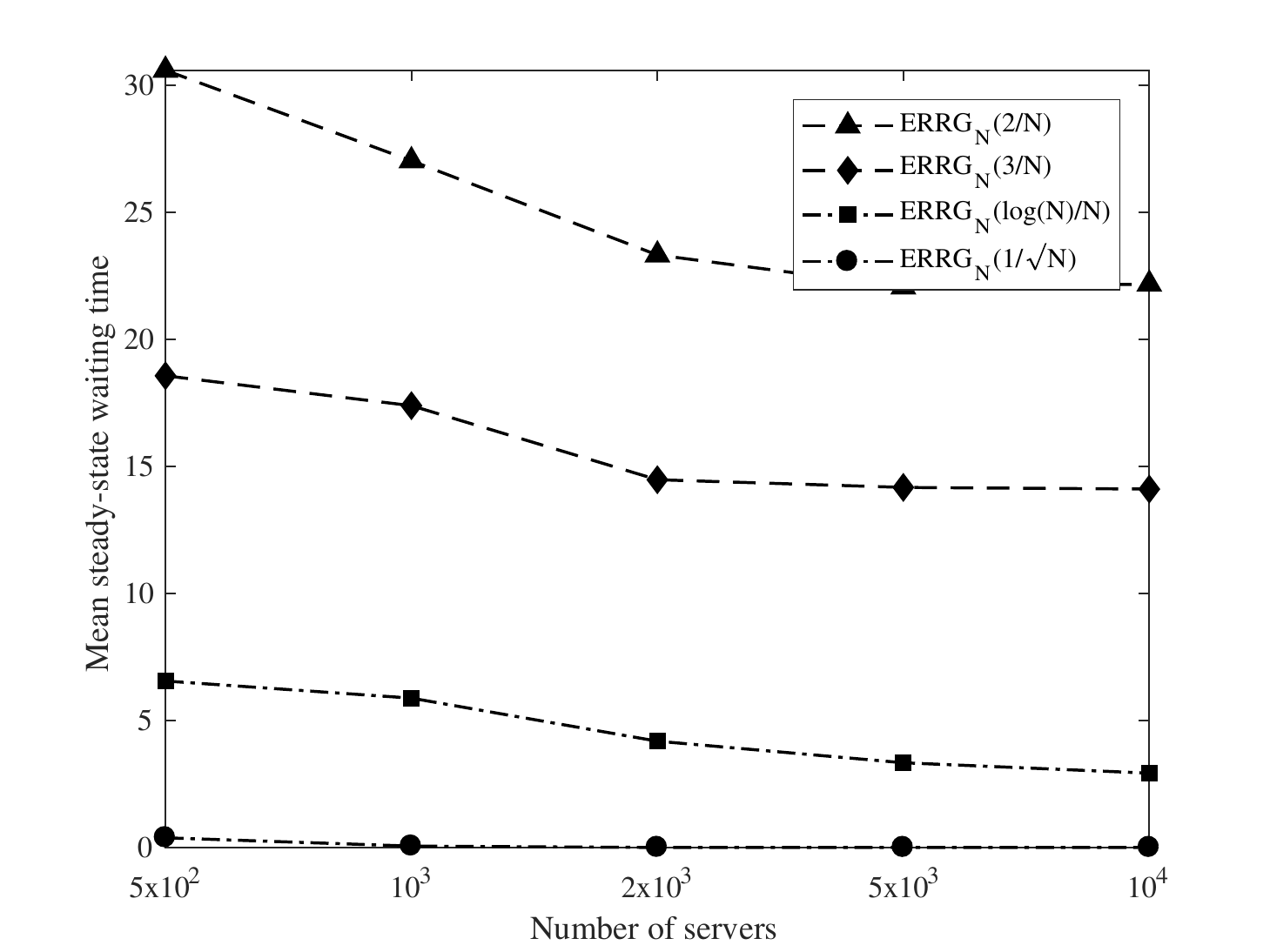}
\end{center}
\caption{Mean steady-state waiting times for $\lambda = 0.9$ and increasing number of servers in ERRG on $N$ vertices with edge probability $c(N)/N$, for $c(N) = 2,3,\log(N),$ and  $\sqrt{N}$.
}
\label{fig:steady-conv}
\end{figure}
\begin{figure}
\begin{center}$
\begin{array}{ccc}
\includegraphics[width=100mm]{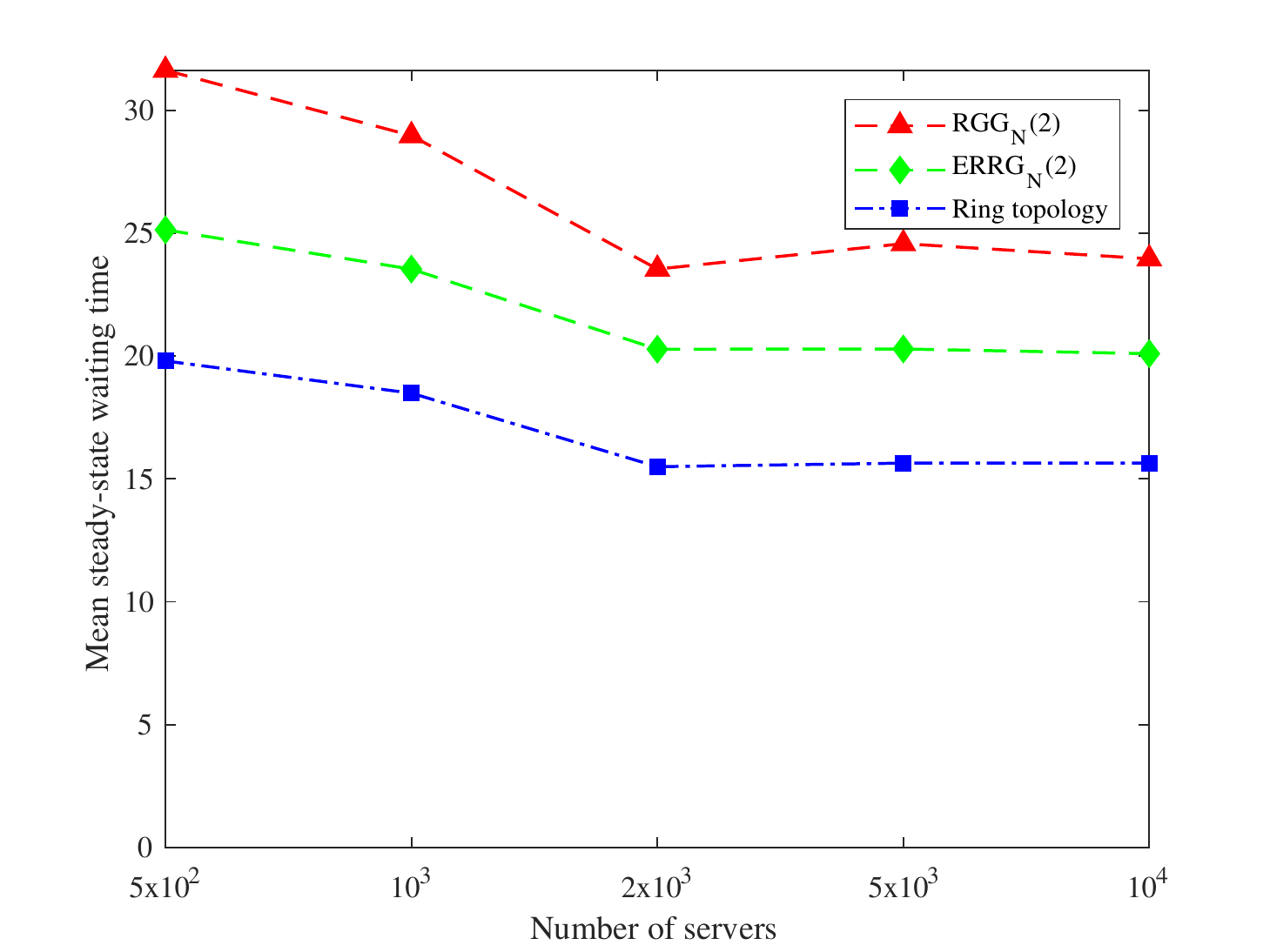}\\
\includegraphics[width=100mm]{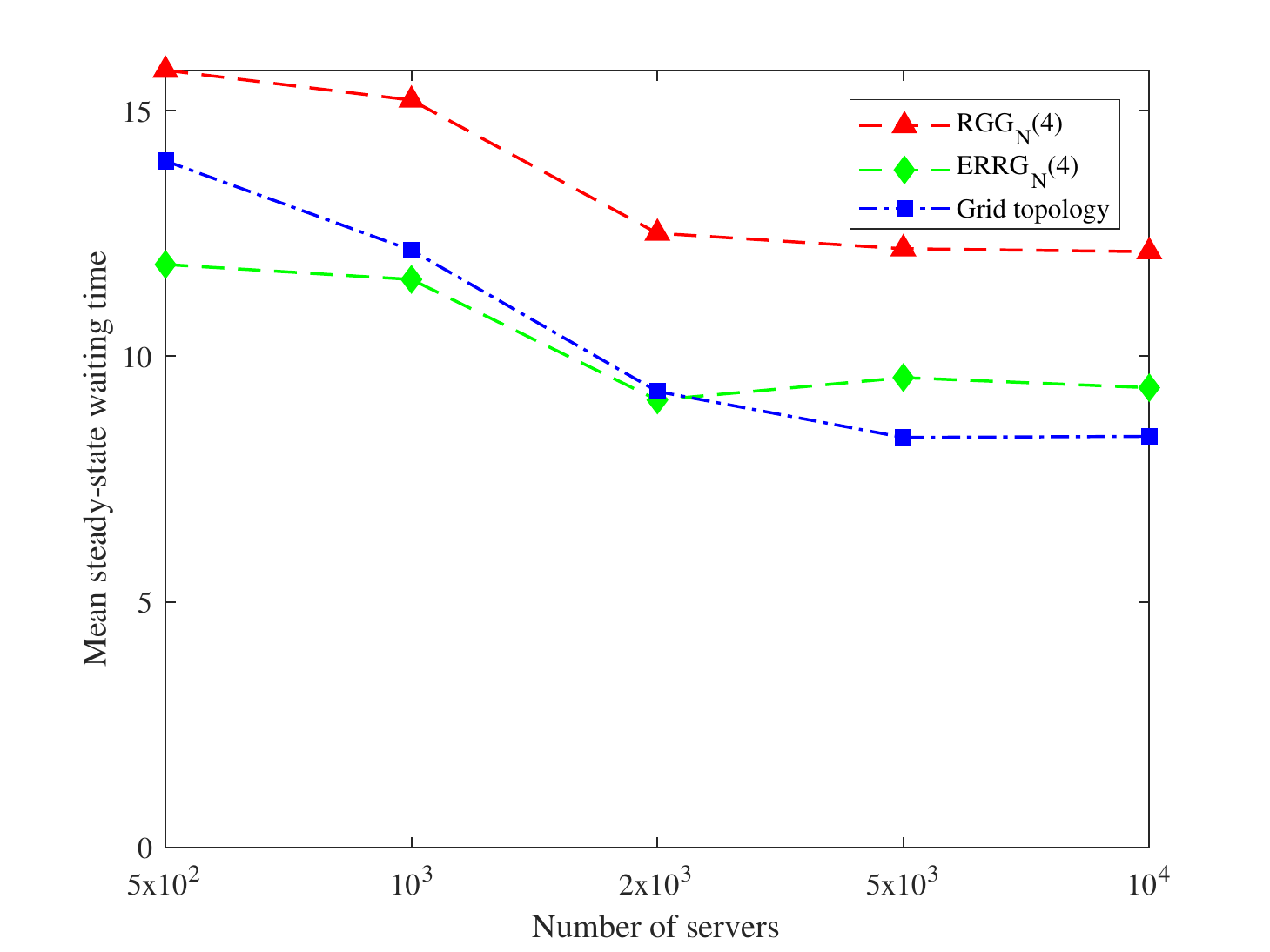}
\end{array}$
\end{center}
\caption{(Top) Performance of the ring topology, and the RGG and ERRG with average degree 2 compared in terms of mean steady-state waiting times. (Bottom) Performance of the grid topology, and the RGG and ERRG with average degree 4 compared in terms of mean steady-state waiting times.}
\label{fig:ring-grid}
\end{figure}
\section{Simulation experiments}
\label{sec:simulation-sig2}
In this section we present extensive simulation results to illustrate the fluid and diffusion-limit results, and compare the performance of various graph topologies in terms of mean waiting times.

\vspace{.25cm}
\noindent
\textbf{ Convergence of sample paths to fluid and diffusion-limit trajectories.}
The fluid-limit trajectory for $\lambda=0.8$ is illustrated in Figure~\ref{fig:trajectory} along with a simulation for $N=10^4$ servers. 
The solid curves represent the case of a clique (i.e.~corresponding to the limit of the occupancy states for the ordinary JSQ policy) as described in Theorem~\ref{fluidjsqd-ssy} in the Chapter~\ref{chap:univjsqd}. The dotted lines correspond to the empirical occupancy process when
the underlying graph topology is a single instance of the Erd\H{o}s-R\'enyi random graph (ERRG) on $N=10^4$ vertices with edge probability $1/\sqrt{N}=10^{-2}$, so the average degree is~100.
Even for a topology much sparser than a clique and  finite $N$-value, the simulated path matches closely with the limiting ODE.
In particular, the above suggests that for a large but finite degree, the behavior may be hard to distinguish from the optimal one for all practical purposes, and there seems to be no prominent effect of graph topologies provided the underlying topology is suitably random.

The diffusion-scaled trajectory has been simulated for $N=10^4$ servers in Figure~\ref{fig:ht}. 
The system load $1-1/\sqrt{N}=0.99$ is quite close to 1. 
The underlying graph topology is taken to be a single instance of the ERRG on $N$ vertices with edge probability $\log(N)^2/\sqrt{N}$.
The green and red curves in Figure~\ref{fig:ht} correspond to the centered and scaled occupancy state processes $-\bQ_1(G_N,\cdot)$ and $\bQ_2(G_N,\cdot)$, respectively.
As stated in Corollary~\ref{cor:errg-sig2}, the centered and diffusion-scaled trajectories can be observed to be recurrent, and the rate of decrease $\bQ_2(G_N,\cdot)$ seems to be proportional to its value --- resembling some properties of the reflected Ornstein-Uhlenbeck process as in the case of a clique (i.e.~the limit of the ordinary JSQ policy) as stated in Theorem~\ref{diffusionjsqd-ssy} in Chapter~\ref{chap:univjsqd}.

\vspace{.25cm}
\noindent
\textbf{ Convergence of steady-state waiting times.} 
Figure~\ref{fig:steady-conv} exhibits convergence of mean steady-state waiting times to their limiting values as $N\to\infty$.
By virtue of Little's law and an interchange of limits argument, note that the asymptotic mean steady-state waiting time can be expressed in terms of the fixed point of the fluid limit as $\lambda^{-1}\sum_{i\geq 2}q_i$.
For each $N$ and average degree $c(N)$ with $c(N)=2$, 3, $\log(N)$, and $\sqrt{N}$, an instance of ERRG on $N$ vertices with average degree $c(N)$ is taken and the time-averaged value of $\lambda^{-1}\sum_{i\geq 2}q_i^N(t)$ is plotted.
The average is taken over the time interval 0 to 200 or 250 depending on the value of $N$.
The figure shows that if the average degree grows with $N$, then the mean steady-state waiting time converges to zero, while it stays bounded away from zero in case the average degree is constant.
It can further be observed that the convergence is notably fast for a higher growth rate of the average degree.

\vspace{.25cm}
\noindent
\textbf{ Effect of the topology in sparse case.}
When the average degree is fixed, the effect of the topology seems to be quite prominent. This has also been observed in prior work~\cite{T98, G15}. Specifically, when comparing graphs with average degree 2, it can be seen in the top chart in Figure~\ref{fig:ring-grid} that the ring topology  has a lower mean steady-state waiting time than random topologies (ERRG or RGG). In case of average degree 4, the (toric) grid topology performs worse for small $N$-values, but the performance improves as $N$ increases. 
There are two crucial effects at play here: 
(i) The regularity in degrees of the vertices: Given a mean degree, higher variability (e.g.~presence of many isolated vertices) is expected to degrade the performance and 
(ii) The locality of the connections: Higher diversity in the connections (i.e., graphs with good expander properties) is expected to improve the performance. 
The RGG has a disadvantage in both these aspects: it contains many isolated vertices and also, its connections are highly localized, and thus its performance is consistently worse in both top and bottom charts in Figure~\ref{fig:ring-grid}. 
The ERRG and the lattice graphs (ring/grid) are good with respect to the degree variability and the connection locality, respectively. However, the presence of many isolated vertices hurts more than the benefit provided by the non-local connections when the average degree is small, as exhibited in Figure~\ref{fig:ring-grid}. In case of higher average degree, the number of isolated vertices in the ERRG is relatively small, and thus the benefit from the non-local connections becomes somewhat prominent for smaller $N$-values. It is therefore worthwhile to note that in case of increasing average degrees, the effect of topology becomes less significant, and so the behavior of random topologies (ERRG, RGG, or random regular graphs) turns out to be as good as the clique.

\begin{figure}
\begin{center}
\includegraphics[width=100mm]{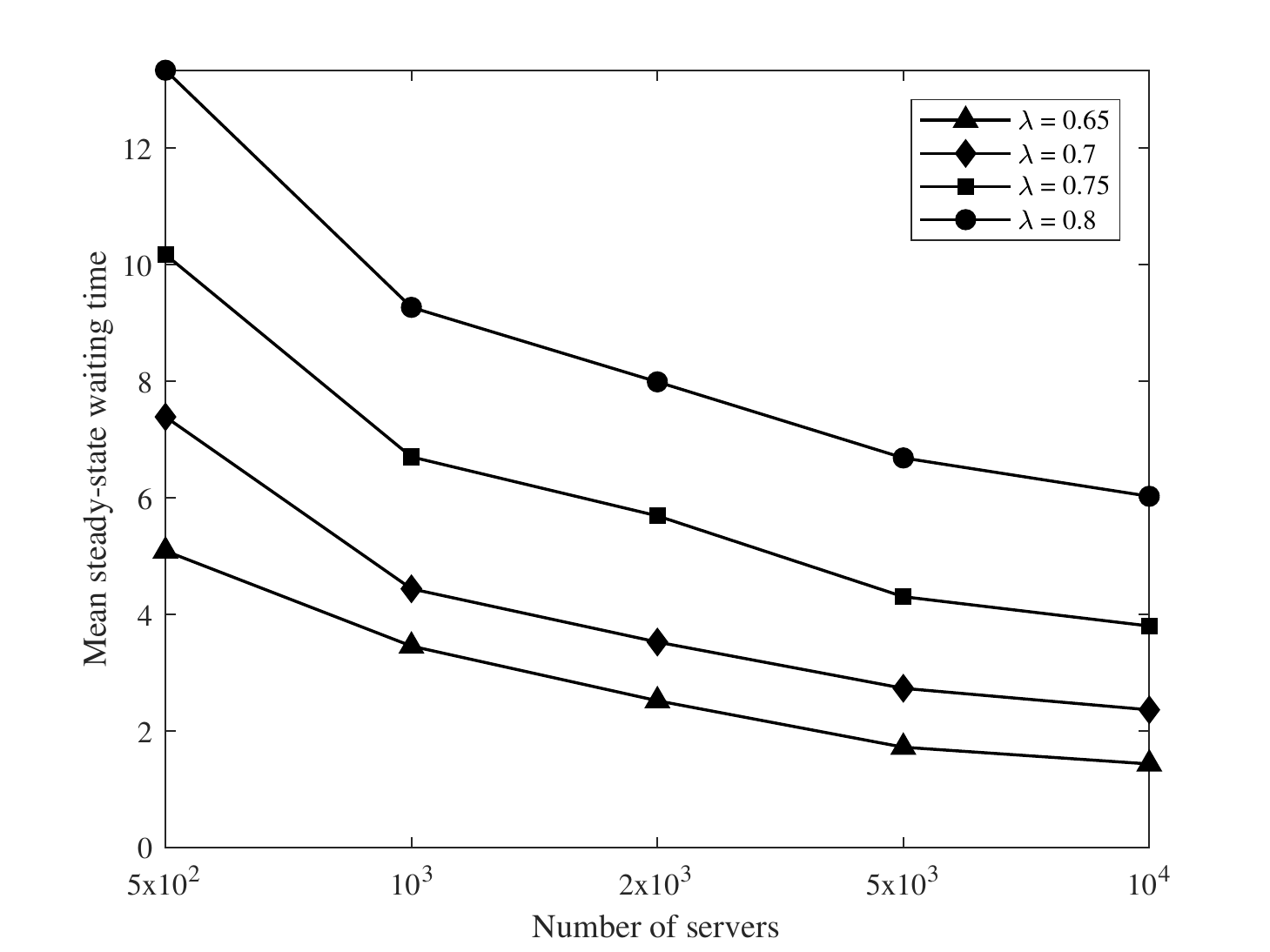}
\end{center}
\caption{Effect of $\lambda$ on the rates of convergence of mean steady-state waiting times.
The underlying topology is an ERRG on $N$ vertices with edge probability $\log(N)/N$, for an increasing number of servers.
}
\label{fig:lambdaeffect}
\end{figure}

\vspace{.25cm}
\noindent
\textbf{ Effect of load on the growth rate of the average degree.}
It is expected that if the system is heavily loaded (i.e., $\lambda$ close to~1), then the rate of convergence of the steady-state measure, and hence that of the mean steady-state waiting time becomes slower.
This can be observed in Figure~\ref{fig:lambdaeffect}.
For moderately loaded systems viz.~$\lambda=0.65$ or $0.75$, the convergence is fast even for topologies that are far from fully connected with average degree as low as $\log(N)$. 

\begin{figure}
\begin{center}
\includegraphics[width=100mm]{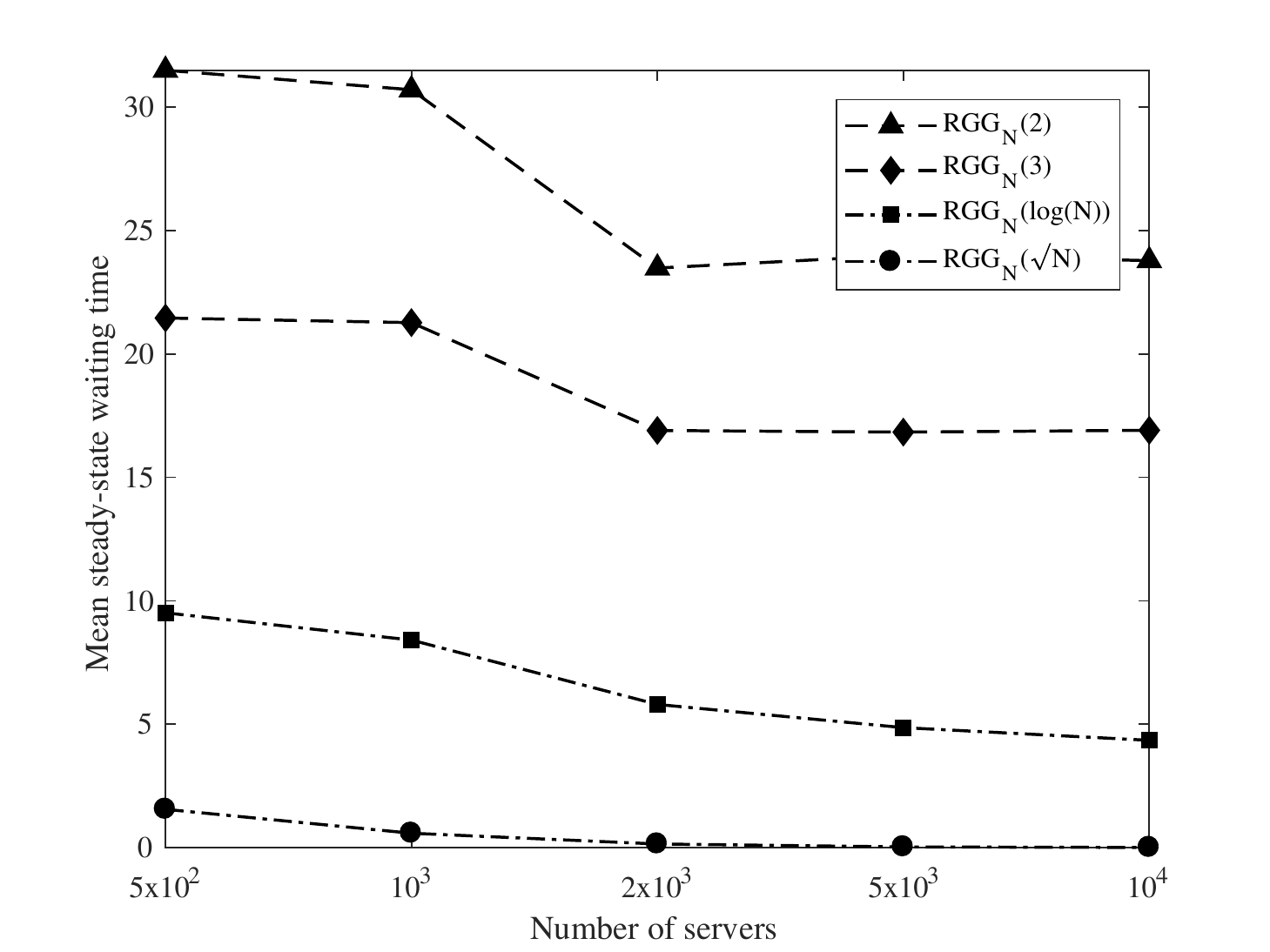}
\end{center}
\caption{Mean steady-state waiting times for $\lambda = 0.8$ and increasing number of servers in RGG on $N$ vertices with average degree $c(N)$, for $c(N) = 2,3,\log(N),$ and  $\sqrt{N}$.
}
\label{fig:steady-conv-rgg}
\end{figure}

\vspace{.25cm}
\noindent
\textbf{ Performance for spatial random network models.}
The conditions stated in Theorem~\ref{th:det-seq-sig2} demand that \emph{any} two \emph{large} portions of the graph share \emph{many} cross edges.
This property is often violated in spatial graph models, where vertices that are closer to each other have a higher tendency to share an edge.
A canonical model for spatial networks is the random geometric graph (RGG), where $N$ vertices correspond to $N$ uniform random locations on $[0,1]^2$ with periodic boundary, and any two vertices share an edge if they are less than a distance $r(N)$ apart.
Note that the average degree in that case is given by $c(N) = (N-1)\pi r(N)^2$.
In other words, for fixed values of $N$ and $c(N)$, the distance $r=r(N)$ scales as $r(N) = \sqrt{c(N)/(\pi N)}$.
To analyze the load balancing process on spatial random graph models, we simulated the processes where the underlying topologies are instances of RGGs on $N$ vertices and average degrees~2,~3, $\log(N)$, and $\sqrt{N}$, and plotted the corresponding mean steady-state waiting times for increasing values of $N$ in Figure~\ref{fig:steady-conv-rgg}.
The surprising resemblance with the ERRG scenario as depicted in Figure~\ref{fig:steady-conv} hints that the asymptotic optimality result can be preserved even under possibly a relaxed set of conditions.
This motivates future study of the asymptotic optimality 
beyond the classes of graphs that we considered.

\section{Conclusion}\label{sec:conclusion-sig2}
We have considered load balancing processes in large-scale systems where the servers are inter-connected
by some graph topology.  
For arbitrary topologies we established sufficient criteria for which the performance is asymptotically similar to that in a clique, and hence optimal on suitable scales.  Leveraging these criteria we showed that unlike fixed-degree scenarios (\emph{viz.}~ring, grid) where the topology has a prominent performance impact, the sensitivity to the topology diminishes in the limit when the average degree grows with the number of servers.
In particular, a wide class of suitably random topologies are provably asymptotically optimal.
In other words, the asymptotic optimality of a clique can be achieved while dramatically reducing the number of connections.
In the context of large-scale data centers,
this translates into significant reductions in communication
overhead and storage capacity, since both are roughly proportional
to the number of connections.

Although a growing average degree is necessary in the sense that any graph with finite average degree is sub-optimal,
it is in no way sufficient. 
Load balancing performance can be provably sub-optimal even when the minimum degree is $cN + o(N)$ with $0 < c < 1/2$.  
What happens for $1/2 < c < 1$ is an open question.
Our proof technique relies heavily on a connectivity property entailing that any two sufficiently large
portions of vertices share a lot of edges.  
This property does not hold however in many networks
with connectivity governed by spatial attributes, such as geometric graphs, although the simulation experiments hint that the family of topologies that are asymptotically optimal is likely to be broader than the ERRG and random regular class as considered in this chapter.  
In future research we aim
to examine asymptotic optimality properties of such spatial network models.


%% file: networkjsqd.tex
\begin{abstract} 
We consider a variation of the supermarket model in which the servers can communicate with their neighbors and where the neighborhood relationships are described in terms of a suitable graph. Tasks with unit-exponential service times arrive at each vertex as independent Poisson processes with rate $\lambda$, and each task is irrevocably assigned to the shortest queue among the one where it first appears and its $d-1$ randomly selected neighbors. This model has been extensively studied when the underlying graph is a clique in which case it reduces to the well known {\em power-of-$d$} scheme. 
We consider settings where the underlying graph need not be a clique and is allowed to be suitably sparse. We show that if the minimum degree approaches infinity (however slowly) as the number of servers $N$ approaches infinity, and the ratio between the maximum degree and the minimum degree in each connected component approaches $1$ uniformly, the occupancy process converges to the same system of ODEs as for the classical supermarket model established in Mitzenmacher~\cite{Mitzenmacher96, Mitzenmacher01} and Vvedenskaya {\em et al.}~\cite{VDK96}. 
In particular, the asymptotic behavior of the occupancy process is insensitive to the precise network topology. We also study the case where the graph sequence is random, with the $N$-th graph given as an Erd\H{o}s-R\'enyi random graph on $N$ vertices with average degree $c(N)$. Annealed convergence of the occupancy process to the same deterministic limit is established under the condition $c(N)\to\infty$, and under a stronger condition $c(N)/ \ln N\to\infty$, convergence (in probability) is shown for almost every realization of the random graph.
\end{abstract}

\section{Introduction}
In this chapter we further explore the impact of the network topology
on the performance of load balancing schemes in large-scale systems, as  discussed in Section~\ref{networks}.
The underlying setup is similar to the one considered in Chapter~\ref{chap:networkjsq}.
The key difference lies in the task assignment strategy: When a task arrives at a server (vertex), it probes a \emph{fixed} number $d$ of its neighbors, in contrast to Chapter~\ref{chap:networkjsq}, where \emph{all} the neighbors are probed.
Thus, the model considered in this chapter can be thought of a network analog of the JSQ($d$) policy, whereas the one in Chapter~\ref{chap:networkjsq} is a network analog of the JSQ policy.
As we will explain below in detail, these changes in the assignment strategy not only make the two systems qualitatively different, but also demand fundamentally different techniques to be developed.\\

We analyze a variation of the supermarket model in which the servers can communicate with their neighbors and where the neighborhood relationships are described in terms of a suitable graph.
Specifically, consider a graph $G_N$ on $N$ vertices, where the vertices represent single-server queues. 
Tasks with unit-exponential service times arrive at each server as independent Poisson processes of rate $\lambda$, and each task is irrevocably assigned to the shortest queue among the one where it first appears and its $d-1$ randomly selected neighbors.

The above model has been extensively investigated in the case where $G_N$ is a clique.
In that case, each task is assigned to the shortest queue among $d\geq 2$ queues selected randomly from the entire system, which is commonly referred to as the `power-of-$d$' or JSQ($d$) scheme (recall Section~\ref{sec:powerofd} in Chapter~\ref{chap:introduction}).
As in Chapter~\ref{chap:networkjsq}, the fundamental challenge in the analysis of load balancing on arbitrary graph topologies is that one cannot reduce the study of the system to that for the state occupancy process $\QQ(\cdot) = (Q_1(\cdot), Q_2(\cdot), \ldots)$ with $Q_i(t)$ being the number of queues with queue length at least $i$ at time $t$, since it is no longer a Markov process. 
In general, one needs to keep track
of the evolution of the number of tasks at each vertex along with the information on neighborhood relationships.
This is a significant obstacle in using tools from classical mean-field analysis for such systems.
Consequently, results for load balancing queuing systems
on general graphs have to date remained scarce.
To the best of our knowledge, this is the first work to study rigorously the limits of the JSQ($d$) occupancy process for non-trivial graph topologies (i.e., other than a clique). 

In Chapter~\ref{chap:networkjsq}, where the tasks are assigned to the shortest queue among \emph{all} the neighbors, we leveraged a stochastic coupling to compare the occupancy process for an arbitrary graph topology with that for the clique, and established that under suitable assumptions on the {\em well-connectedness} of the graph topology, the  occupancy processes and their diffusion-scaled versions
have the same weak limits as for the clique. 
Loosely speaking, for the first convergence, the well-connectedness requires that for any $\varepsilon>0$, the neighborhood of any collection of $\varepsilon N$ vertices contains $N-o(N)$ vertices.
This ensures that on any finite time interval, the fraction of tasks not assigned to servers with the `fluid-scaled minimum queue length' is arbitrarily small.
Thus for large $N$ the occupancy process becomes nearly indistinguishable from that in a clique.
The coupling in Chapter~\ref{chap:networkjsq} is particularly tailored for schemes where on any finite time interval, most of the arrivals are assigned to one of the fluid-scaled shortest queues.
For the setting considered in the current chapter where a fixed number of servers are probed at each arrival,
developing analogous coupling methods  appears to be challenging.
To see this, observe that when all neighbors are probed at arrivals, it is clear that the queue lengths will be better balanced (in the sense of stochastic majorization) for a clique than any other graph topology.
In contrast, for the JSQ($d$) scheme with fixed $d$, even this basic property, namely that the performance of the system will be `optimal' if the topology is a clique, is not clear.
In this chapter, we  take a very different approach, and analyze the evolution of the queue length process at an arbitrary tagged server as the system size becomes large.
The main ingredient is a careful analysis of local occupancy measures associated with the neighborhood of each server and to argue that under suitable conditions their asymptotic behavior is the same for all servers.

Our first result establishes that under fairly mild conditions on the graph topology $G_N$ (diverging minimum degree and a degree regularity condition, see Condition~\ref{cond-reg} and also Remark \ref{rmk:weaker_condition}), for a suitable initial occupancy measure, for any fixed $d\geq 2$, the global occupancy state process for the JSQ($d$) scheme on $G_N$ has the same weak limit as that on a clique, as the number of vertices $N$ becomes large (see Theorem \ref{th:deterministic-d}).
Also, we show that the propagation of chaos property holds, in the sense that the  queue lengths at any finite collection of tagged servers are asymptotically statistically independent, and the 
 queue length process for each server converges in distribution (in the path space) to the corresponding McKean-Vlasov process (see Theorem \ref{th:tagged-d}).
We note that the class of graphs for which the above results hold includes arbitrary $d(N)$-regular graphs, where $d(N)\to\infty$ as $N\to\infty$.
As an immediate consequence of these results, we obtain that the same asymptotic performance of
a JSQ($d$) scheme on  cliques can be achieved by a much sparser graph in which
the number of connections is reduced by almost a factor $N$.  Such a result provides a significant improvement on network connectivity requirements and gives important insights for sparse network design.

When the graph sequence $\{G_N\}_{N\geq 1}$ is random with $G_N$ given as an Erd\H{o}s-R\'enyi random graph (ERRG) with average degree $c(N)$, we establish that for any $c(N)$ that diverges to infinity with $N$, the annealed law of the occupancy process converges weakly to the same limit as in the case of a clique.
For convergence of the quenched law, we require a somewhat more stringent  growth  condition on the average degree.
Specifically, we show that if $c(N)/\log(N)\to\infty$ as $N\to\infty$, then for almost every realization of the random graph the  quenched law of the state occupancy process  converges  to the same  limit as for the case of a clique.
Thus the above results show that the asymptotic performance for cliques  can be achieved by much sparser topologies, even when the connections are random.

In the classical setting of weakly interacting particle systems one considers a collection of
 $N$ stochastic processes on a clique, given as the solution of $N$ coupled stochastic differential equations, where the evolution of any particle at a given time instant depends on its own state and the empirical measure of all particles at that moment (see \cite{Sznitman1989, Kolokoltsov2010,KX99} and references therein).
The asymptotic behavior of the associated state occupancy  measures have been  well studied, including the law of large numbers, propagation of chaos properties, central limit theorems, and large and moderate deviation principles.
However, there is much less work for systems on general graphs except for some recent results for  weakly interacting diffusions on Erd\H{o}s-R\'enyi random graphs. Annealed law of large numbers and central limit theorems for such systems have been established in \cite{BBW16} and a quenched law of large numbers has been shown in \cite{DGL16}.
However these works do not study queueing systems of the form considered here.\\

The rest of the chapter is organized as follows.
In Section~\ref{sec:main-aap2} we present the main results along with some remarks and discussion -- Subsections~\ref{ssec:det} and~\ref{ssec:random} contain the results for sequences of deterministic and random graphs, respectively.
The proofs of  the results in Section \ref{sec:main-aap2} are presented in Section~\ref{sec:proofs-aap2}. 
Finally, we conclude with a discussion of topics for further research in Section~\ref{sec:conclusion-aap2}.

\paragraph{Notation.}
Let $[N] \doteq \{1,\dotsc,N\}$ for $N \in \Nmb$.
For any graph $G_N = (V_N,E_N)$, where $V_N$ is a finite set of vertices and $E_N \subset V_N\times V_N$ is the set of edges,  and $i,j \in V_N$, let
  $\xi_{ij}^N = 1$ if $(i,j)\in E_N$ and $0$ otherwise.
In this chapter, throughout $V_N =[N]$ and $E_N$ will be allowed to be random, in which case $\xi_{ij}^N$ will be random variables. 
Let $\Nmb_0 \doteq \Nmb \cup \{0\}$. For a set $A$, denote by $|A|$ the cardinality. 
For a Polish space $\Smb$, denote by $\Dmb([0,\infty),\Smb)$ the space of right continuous functions with left limits from $[0,\infty)$ to $\Smb$, endowed with the Skorokhod topology. 
For functions $f \colon [0,\infty) \to \Rmb$, let $\|f\|_{*,t} \doteq \sup_{0 \le s \le t} |f(s)|$.
We will use $\kappa,\kappa_1,\kappa_2,\dotsc$ for various non-negative finite constants.
The distribution of an $\Smb$-valued random variable $X$ will be denoted as $\Lmc(X)$.
When the underlying graph is non-random, expectations will be denoted by `$\E$', and when the graphs are random, the notation `$\Ebf$' will be used to denote the expectation (which integrates also over the randomness of the graph topology).


\section{Model description and main results}\label{sec:main-aap2}
Let $\{G_N = (V_N, E_N)\}_{N\geq 1}$ be a sequence of simple graphs where $V_N = [N]$.
The graph $G_N$ corresponds to a system with $N$ servers, where each vertex in the graph  represents a server and edges in the graph define the neighborhood relationships.  
Tasks arrive at the various servers as independent Poisson processes of rate $\lambda$.
 Each server has its own queue with an infinite buffer. Fix $d \in \Nmb$, $d\ge 2$.
 When a task appears at a server $i$, it is immediately assigned to the server with the shortest queue among server $i$ and $d-1$ servers selected uniformly at random from its neighbors in $G_N$.
If there are multiple such servers, one of them is chosen uniformly at random.
Arrivals to any server having less than $d-1$ neighbors in $G_N$ can be assigned in an arbitrary fashion among that server and its neighbors, e.g.
to itself (i.e., without probing the queue length at any other server).
The tasks have independent unit-mean exponentially distributed service times.
The service order at each of the queues is assumed to be oblivious to the actual service time requirements.

Let $X_i^N(t)$ be the number of tasks at the $i$-th server at time instant $t$, and $q^N_j(t)$ be the fraction of servers with queue length at least $j$ in the $N$-th system at time $t$, $i\in [N]$, $j=1,2,\ldots$, namely
\begin{equation}
	q^N_j(t) \doteq \frac{1}{N} \sum_{i=1}^N \sum_{k=j}^{\infty} \one_{\{X^N_i(t)=k\}}, \; t \ge 0, \; j \in \Nmb_0.\label{eq:tailprob}
\end{equation}
Let, $\qq^N(t) \doteq (q^N_i(t))_{i \in \mathbb{N}_0}$.
Then $\qq^N \doteq \{\qq^N(t)\}_{0\le t < \infty}$ is a process with sample paths in $D([0,\infty), S)$ where
$S = \{\qq\in [0,1]^\N: q_0=1, q_i\ge q_{i+1}\ \forall i\in \Nmb_0,\mbox{ and }\sum_{i}q_i<\infty\}$ is equipped with the $\ell_1$-topology.

We will now introduce a convenient representation for the evolution of the queue length processes in the $N$-th system. 
We begin by introducing some  notation. 
For $\xbd = (x_1,\dotsc,x_d) \in \Nmb_0^d$, let $b(\xbd)$ represent the probability that given $d$ servers chosen with queue lengths $\xbd$, the task is sent to the first server in the selection.
Recalling that the task is sent to the shortest queue with ties resolved by selecting at random, the precise definition is as follows:
\begin{equation}
	\label{eq:b}
	b(\xbd) \doteq \sum_{k=1}^d \frac{1}{k} \one_{\displaystyle \{x_1 = \min_{i \in [d]} \{x_i\}, |\mbox{argmin} \{x_i\}|=k \}}.
\end{equation}
Note that 
(i) $b(\xbd)$ is symmetric in $(x_2,\dotsc,x_d)$,
(ii) $b(\xbd) \in [0,1]$, and 
(iii) $b(\xbd)$ is $1$-Lipschitz in $\xbd \in \Nmb_0^d$.
Denote by $D^N_i$ the number of neighbors of a vertex $i$ in $G_N$.
Let $\cN_i$ be iid Poisson processes of rate 1, corresponding to service completions,
and $\Nbar_i$ be iid Poisson random measures on $[0,\infty) \times \Rmb_+$ with intensity $\lambda \, ds \, dy$.
Assume that $\{\cN_i, \Nbar_i\}$ are mutually independent.
Then the evolution of $X_i^N(t)$ can be written as follows:
\begin{equation}
	\label{eq:X_i_n}
	X_i^N(t) = X_i^N(0) - \int_0^t \one_{\{X_i^N(s-)>0\}} \, \cN_i(ds) + \int_{[0,t] \times \Rmb_+} \one_{\{0 \le y \le C_i^N(s-)\}} \, \Nbar_i(ds\,dy),
\end{equation}
where 
\begin{equation}\label{eq:cn}
\begin{split}
	&C_i^N(t)  = \one_{\{D_i^N<d-1\}} \bbar_i^N((X_k^N(t))_{k \in [N]},(\xi_{kl}^N)_{k,l \in [N]}) \\
	&  + \one_{\{D_i^N \ge d-1\}} \sum_{(j_2,\dotsc,j_d) \in \Smc_i^N} \alpha^N(i; j_2, j_3, \ldots, j_d) b(X_i^N(t),X_{j_2}^N(t),\dotsc,X_{j_d}^N(t)) \\
	&  + (d-1) \sum_{(j_2,\dotsc,j_d) \in \Smc_i^N} \one_{\{D_{j_2}^N \ge d-1\}} \alpha^N(j_2; i, j_3, \ldots, j_d) b(X_i^N(t),X_{j_2}^N(t),\dotsc,X_{j_d}^N(t)) \\
	&  + \sum_{j_2 \in [N], j_2 \ne i} \one_{\{D_{j_2}^N<d-1\}}\xi^N_{ij_2} \bbar_{ij_2}^N((X_k^N(t))_{k \in [N]},(\xi^N_{kl})_{k,l \in [N]}), 
	\end{split}
\end{equation}
\begin{equation}
\begin{split}
	\alpha^N(i; j_2, j_3, \ldots, j_d)&:= \frac{\xi^N_{ij_2}\xi^N_{ij_3}\dotsm\xi^N_{ij_d}}{D_i^N(D_i^N-1)\dotsm(D_i^N-d+2)}\\
	\Smc_i^N & := \{ (j_2,\dotsc,j_d) \in [N]^{d-1} : (i,j_2,\dotsc,j_d) \mbox{ are distinct } \}.
\end{split}
\end{equation}
Here $\bbar_i^N$ and $\bbar_{ij}^N$ are measurable functions with 
\begin{equation}
	\label{eq:bbar}
	\bbar_i^N\big((X_k^N(t))_{k \in [N]},(\xi_{kl}^N)_{k,l \in [N]}\big), \bbar_{ij}^N\big((X_k^N(t))_{k \in [N]},(\xi^N_{kl})_{k,l \in [N]}\big) \in [0,D_i^N+1], 
\end{equation}
which define the rules of assigning tasks when $D_i^N<d-1$ or $D_j^N<d-1$, respectively.
The precise form of these functions will not be important in our analysis.
The second term in the expression for $C^N_i(t)$ gives the probability that a task arriving at server $i$ (with $D^N_i\ge d-1$)
is in fact assigned to server $i$ itself, which will happen if server $i$ is one of the queues with minimum queue length among the $d-1$ randomly selected neighbors and itself, and it is the winner of the tie among the servers with minimum queue lengths in the selection. The third term corresponds to the probability that a task arriving at some other server (say $j_2$, with $D^N_{j_2}\ge d-1$) is assigned to server $i$, which will happen if $i$ is a neighbor of $j_2$, server $i$ is among the random selection of $d-1$ neighbors of $j_2$, it is also among the servers with minimum queue length in the selection, and it wins the tie-breaker among the servers with minimum queue length in the selection.


\subsection{Scaling limits for deterministic graph sequences}\label{ssec:det}
In this section we will consider arbitrary deterministic graph sequences, and establish a scaling limit  when the graphs satisfy a certain `regularity' condition  as formulated in Condition~\ref{cond-reg} below.
For any graph $G$, let $d_{\min}(G)$ and $d_{\max}(G)$ denote the minimum and maximum degree, respectively.
\begin{condition}[Regularity of degrees]\label{cond-reg}
The sequence $\{G_N\}_{N\geq 1}$  satisfies the following.
\begin{enumerate}[{\normalfont (i)}]
\item $d_{\min}(G_N)\to\infty$ as $N\to\infty$.
\item $\max_{i\in [N]} \left | \sum_{j\in [N], j\neq i} \frac{\xi_{ji}^N}{D_j^N} - 1\right| \to 0$ as $N\to\infty$. 
\end{enumerate}
\end{condition}

\begin{remark}\normalfont
	\label{rmk:weaker_condition}
	Condition \ref{cond-reg}(ii) holds if for example, $d_{\max}(G_N)/d_{\min}(G_N)\to 1$ as $N\to\infty$, since
	\begin{align*}
		\frac{d_{\min}(G_N)}{d_{\max}(G_N)} \le \frac{D_i^N}{d_{\max}(G_N)} 
		\le \sum_{j \in [N], j \ne i} \frac{\xi_{ji}^N}{D_j^N} \le \frac{D_i^N}{d_{\min}(G_N)} \le \frac{d_{\max}(G_N)}{d_{\min}(G_N)}
	\end{align*}
	for each $i \in [N]$.
	But Condition \ref{cond-reg}(ii) also allows $G_N$ to have degrees of very different orders in different components of the graph. 
For example, if $\{\mathcal{C}^N_k\}_{k\geq 1}$ denote the  connected components of $G_N$, then
Condition~\ref{cond-reg} (ii) is satisfied if
\[\sup_{k\geq 1} \left|\frac{d_{\min}(\mathcal{C}_k^N)}{d_{\max}(\mathcal{C}_k^N)}- 1\right| \to 0 \qquad\mbox{as}\qquad N\to\infty.\]
\end{remark}

Our first  result establishes under Condition~\ref{cond-reg}, the convergence of the occupancy state process $\qq^N$ to the same deterministic limit as for the classical JSQ($d$) policy (i.e. the case when $G_N$ is a clique), as $N\to\infty$.
\begin{theorem}[Convergence of global occupancy states]
\label{th:deterministic-d}
Assume that the sequence of graphs $\{G_N\}_{N\geq 1}$ satisfies Condition~\ref{cond-reg}, and $\{X_i^N(0):i\in [N]\}$ is iid with $\Pro{X_i^N(0) \geq j}=q^{\infty}_j$, $j=1,2,\ldots,$ for some $\qq^{\infty}\in S$. 
Then on any finite time interval, the occupancy state process $\qq^N(\cdot)$ converges weakly with respect to the Skorohod $J_1$-topology to the deterministic limit $\qq(\cdot)$ given by the unique solution to the system of ODEs: 
	\begin{equation}\label{eq:occ-deterministic}
	\frac{dq_i(t)}{dt} = \lambda [(q_{i-1}(t))^d - (q_i(t))^d] - (q_{i}(t)-q_{i+1}(t)),\quad i = 1, 2,\ldots,
	\end{equation}	 
	and $\qq(0) = \qq^\infty$.
\end{theorem}
\begin{remark}\normalfont
We make the following observations.
	\begin{enumerate}[\normalfont (i)]
\item Unique solvability of the system of equations	\eqref{eq:occ-deterministic} is a consequence of Lipschitz continuity of the right side.
Specifically, define the function $\FF(\cdot) = (F_1(\cdot), F_2(\cdot), \ldots)$ on $S$ as 
\[F_i(\qq) = \lambda(q_{i-1}^d - q_i^d) - (q_i - q_{i+1}), \quad i= 1,2, \ldots,\] 
with $\qq\in S$ and $F_i(\qq)$ being the $i$-th component of $F(\qq)$. 
It is easily seen that $F$ is  Lipschitz on $S$ (equipped with the  $\ell_1$-distance).
Standard results then imply that the system of ODEs defined by $\dif\qq(t)/\dif t = \FF(\qq)$ admits a unique solution.

\item The above result shows in particular  that the evolution of the limiting global occupancy process as described by~\eqref{eq:occ-deterministic} coincides with that when the underlying graph is a clique, i.e., when each arriving task can probe any set of $d$ servers.
Thus under Condition~\ref{cond-reg}, the system exhibits the same asymptotic transient performance even when the underlying graph is much sparser.
As an immediate corollary we see that \eqref{eq:occ-deterministic} describes  the limiting system occupancy process associated with {\em arbitrary} $d(N)$-regular graphs as long as $d(N)\to\infty$ as $N\to\infty$.
\end{enumerate}
\end{remark}

\begin{remark}\normalfont\label{rem:compare-20}
Now we contrast Condition \ref{cond-reg} with
the condition stated in Theorem~\ref{th:det-seq} in Chapter~\ref{chap:networkjsq} for the JSQ policy on a graph to behave as that on a clique.
We note that Condition~\ref{cond-reg} relies only on local properties of the graph, and in particular may hold even when, for example, the graph contains several connected components of sizes that grow to infinity with $N$.
In contrast, the condition in Chapter~\ref{chap:networkjsq} requires that any two $\Theta(N)$-sized component must share $\Theta(N)$ cross-edges, which does not hold in many networks with connectivity governed by spatial attributes, such as geometric graphs.
In this sense, Condition~\ref{cond-reg} includes a much broader class of graphs including arbitrary $d(N)$-regular graphs with  $d(N)\to\infty$, as mentioned above.
On the other hand, our condition requires the minimum degree in the graph to diverge to infinity, whereas Theorem~\ref{th:det-seq} allows any $o(N)$ vertices to have bounded degree (or degree zero).
As noted in the introduction, it is easy to see that the queue length process of the JSQ policy on a clique is better balanced (in stochastic majorization sense) than on any other graph. 
This is also reflected by the fact that the sufficient criterion for fluid optimality as developed in Chapter~\ref{chap:networkjsq} is monotone with respect to edge addition.
Specifically, let $\{G_N=(V_N, E_N)\}_{N\geq 1}$ be a graph sequence which satisfies the sufficient criterion in Theorem~\ref{th:det-seq} for the limit of the  occupancy process coincides with that for cliques. 
Then Theorem~\ref{th:det-seq} guarantees that for any graph sequence $\{\bar{G}_N = (V_N, \bar{E}_N)\}_{N\geq 1}$ with $E_N\subseteq\bar{E}_N$, the limit of the  occupancy process also coincides with that for cliques.
The above property is not immediate for systems considered in this chapter since  adding edges arbitrarily may  result in violating Condition~\ref{cond-reg} (ii).
\end{remark}

Our second result gives the joint asymptotic  behavior of queue length processes for any finite collection of servers.
In particular, it shows that the propagation of chaos holds, i.e., the queue length processes for any finite collection of servers are asymptotically statistically independent.
Recall from Section~\ref{sec:main-aap2} the sequence of Poisson processes $\{\cN_i\}$, Poisson random measures $\{\Nbar_i\}$, and the function $b$.

\begin{theorem}[Evolution of tagged servers]\label{th:tagged-d}
Assume that the  sequence of graphs $\{G_N\}_{N\geq 1}$ satisfies Condition~\ref{cond-reg},  and $\{X_i^N(0):i\in [N]\}$ is iid with 
$$\Pro{X_i^N(0) \geq j}=q^{\infty}_j, \quad j=1,2,\ldots,$$ 
for some $\qq^{\infty}\in S$. 
Then the following convergence results hold.
\begin{enumerate}[{\normalfont (i)}]
\item On any finite time interval, the queue length process  $X_i^N(\cdot)$ at server $i$ converges weakly with respect to the Skorohod-$J_1$ topology to the following  McKean-Vlasov process:
\begin{equation}\label{eq:limit-tagged}
\begin{split}
	X_i(t) & = X_i(0) - \int\limits_0^t \one_{\{X_i(s-)>0\}} \, \cN_i(ds) + \int\limits_{[0,t]\times\Rmb_+} \one_{\{0 \le y \le C_i(s-)\}} \, \Nbar_i(ds\,dy), \\
	C_i(t) & = d\int_{\Nmb^{d-1}} b(X_i(t),x_2,\dotsc,x_d) \mu_t(dx_2)\dotsm\mu_t(dx_d) ,
\end{split}
\end{equation}
where $\mu_t = \Lmc(X_i(t))$ and $\mu_0[j,\infty) = q^{\infty}_j$ for $t\ge 0$ and $j \in \Nmb_0$.
\item For any  $m$-tuple $(i_1,\dotsc,i_m)\in \N^m$ with $i_j \ne i_k$ whenever $j \ne k$,
		\begin{equation*}
			\Lmc(X_{i_1}^N(\cdot),\dotsc,X_{i_m}^N(\cdot)) \to \mu^{\otimes m},
		\end{equation*}
		as probability measures on $D([0, \infty): \Nmb_0^m)$ where $\mu$ is the probability law of $X_1(\cdot)$ in part (i).
	\item For any $i\in\N$, the process $\mu^{i,N}$ denoting the occupancy measure process for the neighborhood of the $i$-th server, defined as
	\begin{equation}
	\label{eq:nbdempmzr}
	\mu^{i,N}_t \doteq   \frac{1}{D_i^N+1} \sum_{j \in [N], j \ne i} \xi_{ij}^N \delta_{X_j^N(t)} + 	\frac{1}{D_i^N+1}\delta_{X_i^N(t)}, \; t\ge 0,
	\end{equation}
	converges weakly with respect to the Skorohod $J_1$-topology to the deterministic limit	$\mu$ as in part (i).
\end{enumerate}
\end{theorem}
\begin{remark}\label{rem:rem2}\normalfont
We note the following.
\begin{enumerate}[{\normalfont (i)}]
		\item The existence and uniqueness of solutions to \eqref{eq:limit-tagged} can be proved by standard arguments using the boundedness and Lipschitz property of the functions $b$ and $x\mapsto \one_{\{x>0\}}$ on $\Nmb_0$.
		\item  Using the propagation of chaos property and the fact that $\{X_i(t):i\in [N]\}$ are iid, it follows
		that
		the limit of the global occupancy measure  at any time instant $t$ is in fact the law of $X_i(t)$ for any fixed $i$.
		Therefore,
		$$\mu_t[j,\infty)= \Pro{X_i(t)\geq j} = q_j(t),\;  j  \in \Nmb_0, i \in \Nmb \mbox{ and } t \ge 0.$$
	\end{enumerate}
\end{remark}

\subsection{Scaling limits for random graph sequences}\label{ssec:random}
Next we will consider the scenario when the underlying graph topology is random.
We consider asymptotics of both annealed and quenched laws of the occupancy process and the queue length process at any tagged server.
The following is our main condition in the study of the annealed law.
\begin{condition}[Diverging mean degree]\label{cond:errg1}
 $\{G_N\}_{N\geq 1}$ is a sequence of Erd\H{o}s-R\'enyi random graphs  where any two vertices share an edge with probability $p_N$, and  $Np_N\to\infty$ as $N\to\infty$. 	$\{G_N\}_{N\geq 1}$ is independent of $\{X^N_j(0), \cN_i, \Nbar_i, j \in [N], N \in \Nmb, i \in \Nmb\}$.
\end{condition}


\begin{theorem}[Asymptotics of annealed law]\label{thm:npn_rate}
Assume that the  sequence of graphs $\{G_N\}_{N\geq 1}$ satisfies Condition~\ref{cond:errg1}, and $\{X_i^N(0):i\in [N]\}$ is iid with 
$$\Pro{X_i^N(0) \geq j}=q^{\infty}_j,\quad j=0, 1,2,\ldots,$$ 
for some $\qq^{\infty}\in S$. 
Then the following  hold.
\begin{enumerate}[{\normalfont (i)}]
\item For any $T\in (0,\infty)$
		\begin{equation}
		\label{eq:npn_rate}
		\sup_{N \ge 1} \max_{i \in [N]} \sqrt{Np_N} \Ebf \|X_i^N - X_i\|_{*,T}^2 < \infty,
		\end{equation}
		where $X_i$ is as defined in \eqref{eq:limit-tagged}.
\item For any $m$-tuple $(i_1,\dotsc,i_m)\in \N^m$ with $i_j \ne i_k$ whenever $j \ne k$, 		
\begin{equation*}
			\Lmc(X_{i_1}^N(\cdot),\dotsc,X_{i_m}^N(\cdot)) \to \mu^{\otimes m},
		\end{equation*}
		as probability measures on $D([0, \infty): \Nmb_0^m)$ where $\mu$ is as in Theorem \ref{th:tagged-d}.
	\item For any $i\in\N$, the law of the neighborhood occupancy measure process  defined as in \eqref{eq:nbdempmzr}
	converges weakly in the Skorohod $J_1$-topology to the deterministic limit $\mu$ as in Theorem \ref{th:tagged-d}. 
\end{enumerate}
\end{theorem}
\begin{remark}\normalfont
We make the following observations.
	\begin{enumerate}
\item	In contrast to standard convergence results for weakly interacting diffusions (see e.g.~\cite{Sznitman1989} or \cite{BBW16}), the estimate in \eqref{eq:npn_rate} gives a rate of convergence of $\sqrt{Np_N}$ instead of $Np_N$. The reason for this can be seen from the proof which shows that  the bound for the  quantity $\Ebf \|X_i^N - X_i\|_{*,T}^2$ is controlled by $\Ebf |C_i^N(s)-C_i(s)|$ rather than $\Ebf |C_i^N(s)-C_i(s)|^2$, due to the form of indicator function in the evolution of $X_i^N$ (cf.~\eqref{eq:X_i_n}).
	\item The conditions needed for   Theorem~\ref{thm:npn_rate} should be contrasted with that for Theorems \ref{th:deterministic-d}
	and \ref{th:tagged-d}.
	In particular, for the study of the  annealed law asymptotics we only need information on the average degree rather than on the  maximum and minimum degrees of the graph.
\end{enumerate}
\end{remark}

 We will now consider the asymptotic behavior of the quenched law of the occupancy process. For this we formulate a condition that is stronger than the one used in the study of the annealed asymptotics.
\begin{condition}[Condition for quenched limit]\label{cond:errg2}
 $\{G_N\}_{N\geq 1}$ is a sequence of Erd\H{o}s-R\'enyi random graphs, such that in $G_N$ any two vertices share an edge with probability $p_N$, and  $Np_N/\ln(N)\to\infty$ as $N\to\infty$. $\{G_N\}_{N\geq 1}$ is independent of 
 $$\{X^N_j(0), \cN_i, \Nbar_i, j \in [N], N \in \Nmb, i \in \Nmb\}.$$
\end{condition}
The following theorem provides, under the above condition, the asymptotic behavior of the quenched law.
\begin{theorem}[Asymptotics of quenched law]\label{thm:npn_rate_quench}
Assume that the  sequence of graphs $\{G_N\}_{N\geq 1}$ satisfies Condition~\ref{cond:errg2}, and $\{X_i^N(0):i\in [N]\}$ is iid with 
$$\Pro{X_i^N(0) \geq j}=q^\infty_j, \quad j=0,1,2,\ldots,$$ for some $\qq^\infty\in S$ for all $N$. 
Then the convergence results as stated in Theorems~\ref{th:deterministic-d} and~\ref{th:tagged-d} hold for almost every realization of the random graph sequence. 
\end{theorem}

\section{Proofs}\label{sec:proofs-aap2}

\subsection{Proofs for deterministic graph sequences}
An overview of the proof idea is as follows.
First note that the queue length process at any two vertices can be exactly coupled to evolve identically if the occupancy measure of the corresponding neighborhoods are indistinguishable.
The main step is to show that if the graph sequence satisfies  Condition~\ref{cond-reg}, then the local occupancy measure associated with the neighborhood of every server over any finite time interval converges to the same limit as for the global occupancy measure, which in turn is the same as that when the whole system uses the ordinary JSQ($d$) policy and the graph is a clique. 
This ensures that the rate of arrival (exogenous + forwarded from the neighboring vertices) to a typical server is (asymptotically) the same as that in the clique case. 
Thus, the law of the number of tasks at each server, and consequently the global occupancy measure, converge to the same limit.
For technical convenience we will provide the proof of Theorem~\ref{th:tagged-d} first, and then use that to establish Theorem~\ref{th:deterministic-d}.

We will define the limiting processes $\{(X_i(t))_{i\geq 1}\}_{t\geq 0}$ and the pre-limit processes $\{(X_i^N(t))_{i\geq 1}\}_{t\geq 0}$ on the same probability space by taking the same sequence of Poisson processes $\{\cN_i\}$ and Poisson random measures $\{\Nbar_i\}$ in both cases.
Also, take $X_i^N(0) = X_i(0)$ for all $i\in [N]$, $N\geq 1$.
Using Condition~\ref{cond-reg} we can find a $N_0 \in \Nmb$ such that for all $N\ge N_0$
\begin{equation}
	\label{eq:degree_large_N}
	d_{\min}(G_N) \ge d, \;  \sup_{i\in [N]} \left | \sum_{j\in [N], j\neq i} \frac{\xi_{ji}^N}{D_j} - 1\right| \le \frac{1}{2},\quad \sup_{i\in [N]} \sup_{t \in [0,T]} \left|C_i^N(t)\right| \le 2d.
\end{equation}
For the rest of this section we will assume that $N\ge N_0$ and therefore, in particular, the first and fourth terms in the definition of $C_i^N(s)$ are zero and the indicators in the second and third terms can be replaced by $1$.
We will frequently suppress $N$ in the notation $D^N_i$ and $\xi_{ij}^N$ and write them as $D_i$ and $\xi_{ij}$ respectively. We begin with the following lemma. The proof is given at the end of the subsection.
\begin{lemma}
	\label{lem:prep_quench}
	For $i \in [N]$ and $s\in [0,T]$ let
		\begin{align*}
			U_s  \doteq \E \Big[  \sum_{(j_2,\dotsc,j_d) \in \Smc_i^N} \alpha^N(i; j_2, j_3, \ldots, j_d)    
			   \Big(b(X_i(s),X_{j_2}(s),\dotsc,X_{j_d}(s)) - \frac{C_i(s)}{d} \Big) \Big]^2
		\end{align*}
	and
		\begin{align*}	
			V_s  \doteq \E \Big[ \sum_{(j_2,\dotsc,j_d) \in \Smc_i^N}  \alpha^N(j_2; i, j_3, \ldots, j_d)      \Big(b(X_i(s),X_{j_2}(s),\dotsc,X_{j_d}(s))-\frac{C_i(s)}{d}\Big) \Big]^2.
		\end{align*}
Under the conditions of Theorem~\ref{th:deterministic-d}, there exists $K \in (0,\infty)$ such that for every $s \in [0,T]$ and  $i \in [N]$,
	\begin{align}
		U_s   \le \frac{K}{d_{\min}(G_N)} , \;\;
		V_s   \le \frac{K}{d_{\min}(G_N)} \left( \sum_{j=1,j \ne i}^N   \frac{\xi_{ji}}{D_j} \right)^2 . \label{eq:prep_quench_i}
	\end{align}
\end{lemma}

\begin{proof}[Proof of Theorem~\ref{th:tagged-d}]
Fix any $i\in\N$ and $T>0$.
From \eqref{eq:X_i_n} and \eqref{eq:limit-tagged}, using the Cauchy--Schwarz and Doob's inequalities we have for any fixed $t\in [0,T]$ and $N\ge i$,
	\begin{align}
		\E \left\|X_i^N - X_i\right\|_{*,t}^2 & \le \kappa_1 \E \int_0^t |\one_{\{X_i^N(s)>0\}} - \one_{\{X_i(s)>0\}}|^2 \, ds\\
		&\qquad + \kappa_1 \E \left( \int_0^t |\one_{\{X_i^N(s)>0\}} - \one_{\{X_i(s)>0\}}| \, ds \right)^2 \notag \\
		& \qquad + \kappa_1 \E \int_{[0,t]\times\Rmb_+} |\one_{\{0 \le y \le C_i^N(s)\}} - \one_{\{0 \le y \le C_i(s)\}}|^2 ds\,dy \notag \\
		& \qquad + \kappa_1 \E \left(\int_{[0,t]\times\Rmb_+} |\one_{\{0 \le y \le C_i^N(s)\}} - \one_{\{0 \le y \le C_i(s)\}}| \,ds\,dy\right)^2 \notag \\
		& \le \kappa_1 \int_0^t \E |X_i^N(s) - X_i(s)|^2 \, ds + \kappa_1 \E \left( \int_0^t |X_i^N(s) - X_i(s)| \, ds \right)^2 \notag \\
		& \quad + \kappa_1 \int_0^t \E |C_i^N(s)-C_i(s)| \, ds + \kappa_1 \E \left( \int_0^t |C_i^N(s)-C_i(s)| \, ds \right)^2 \notag \\
		& \le \kappa_2 \int_0^t \E |X_i^N(s) - X_i(s)|^2 \, ds + \kappa_2 \int_0^t \E |C_i^N(s)-C_i(s)| \, ds \label{eq:npn_quench_1}
	\end{align}
for some  $\kappa_1, \kappa_2 \in (0,\infty)$, where in the last line we have used \eqref{eq:degree_large_N} and the fact that $0 \le \frac{C_i(s)}{d} \le 1$.

Now we analyze the difference $|C_i^N(s)-C_i(s)|$ in \eqref{eq:npn_quench_1}. 
Note that by adding and subtracting terms we have
	\begin{equation}
		|C_i^N(s)-C_i(s)| \le |C_i^N(s)-C_i^{N,1}(s)| + |C_i^{N,1}(s)-C_i^{N,2}(s)| + |C_i^{N,2}(s)-C_i(s)|, \label{eq:npn_2}
	\end{equation}
	where
	\begin{align*}
			C_i^{N,1}(s) & =  \sum_{(j_2,\dotsc,j_d) \in \Smc_i^N} \alpha^N(i; j_2, j_3, \ldots, j_d) b(X_i(s),X_{j_2}(s),\dotsc,X_{j_d}(s)) \\
			& \quad + (d-1) \sum_{(j_2,\dotsc,j_d) \in \Smc_i^N}  \alpha^N(j_2; i, j_3, \ldots, j_d) b(X_i(s),X_{j_2}(s),\dotsc,X_{j_d}(s)) 
		\end{align*}
		and
		\begin{align*}
			C_i^{N,2}(s) & =   \sum_{(j_2,\dotsc,j_d) \in \Smc_i^N} \alpha^N(i; j_2, j_3, \ldots, j_d) \frac{C_i(s)}{d} \\
			&\hspace{3.5cm} + (d-1) \sum_{(j_2,\dotsc,j_d) \in \Smc_i^N}  \alpha^N(j_2; i, j_3, \ldots, j_d) \frac{C_i(s)}{d}.
		\end{align*}
We now analyze each term in \eqref{eq:npn_2}.
In particular,  we will use the Lipschitz property of $b$ to handle the term $|C^N_i-C^{N,1}_i|$, and then use the iid property of the $X_i$'s to handle the term $|C^{N,1}_i-C^{N.2}_i|$. 

First consider $|C_i^N(s)-C_i^{N,1}(s)|$.
From the Lipschitz property of $b$ and the definition of $\alpha^N$ we have
	\begin{align*}
		&\E |C_i^N(s)-C_i^{N,1}(s)|\\
		& \le \E \Big[  \sum_{(j_2,\dotsc,j_d) \in \Smc_i^N} \alpha^N(i; j_2, j_3, \ldots, j_d)  \\ 
		& \qquad (|X_i^N(s)-X_i(s)| + |X_{j_2}^N(s)-X_{j_2}(s)| + \dotsb + |X_{j_d}^N(s)-X_{j_d}(s)|) \\
		& \qquad + (d-1) \sum_{(j_2,\dotsc,j_d) \in \Smc_i^N}  \alpha^N(j_2; i, j_3, \ldots, j_d)\\
		& \qquad  (|X_i^N(s)-X_i(s)| + |X_{j_2}^N(s)-X_{j_2}(s)| + \dotsb + |X_{j_d}^N(s)-X_{j_d}(s)|) \Big], \\
		& \le  \max_{j \in [N]} \E |X_j^N(s)-X_j(s)|\Big(d   + (d-1)d  \sum_{j_2 \in [N], j_2 \ne i}  \frac{\xi_{j_2i}}{D_{j_2}}\Big).
	\end{align*}
From \eqref{eq:degree_large_N} we have
	\begin{equation}
		\label{eq:npn_quench_difference_1}
		\E |C_i^N(s)-C_i^{N,1}(s)| \le \kappa_3 \max_{j \in [N]} \E |X_j^N(s)-X_j(s)|
	\end{equation}
for some  $\kappa_3 \in (0,\infty)$.
Next we consider $|C_i^{N,1}(s)-C_i^{N,2}(s)|$.
It follows from Cauchy--Schwarz inequality  that
	\begin{align*}
		&\E |C_i^{N,1}(s)-C_i^{N,2}(s)|^2\\
		& \le 2 \E \Big[  \sum_{(j_2,\dotsc,j_d) \in \Smc_i^N} \alpha^N(i; j_2, j_3, \ldots, j_d)   
		 \Big(b(X_i(s),X_{j_2}(s),\dotsc,X_{j_d}(s)) - \frac{C_i(s)}{d}\Big) \Big]^2  \\
		& \hspace{3cm} + 2(d-1)^2 \E \Big[ \sum_{(j_2,\dotsc,j_d) \in \Smc_i^N}  \alpha^N(j_2; i, j_3, \ldots, j_d)  \\
		  &\hspace{5cm}\times\Big(b(X_i(s),X_{j_2}(s),\dotsc,X_{j_d}(s))-\frac{C_i(s)}{d}\Big) \Big]^2   \\
		 &\le \kappa_4(U_s + V_s).
	\end{align*}
where $U_s,V_s$ are as in Lemma \ref{lem:prep_quench}.
From Lemma~\ref{lem:prep_quench} and \eqref{eq:degree_large_N} we obtain
\begin{eq}
	\label{eq:npn_quench_difference_2}
\big(\E |C_i^{N,1}(s)&-C_i^{N,2}(s)|\big)^2\leq \E |C_i^{N,1}(s)-C_i^{N,2}(s)|^2 \\
&\le \frac{\kappa_5}{d_{\min}(G_N)}  + \frac{\kappa_5 }{d_{\min}(G_N)} \left( \sum_{j=1,j \ne i}^N   \frac{\xi_{ji}}{D_j} \right)^2
\le \frac{\kappa_6}{d_{\min}(G_N)}.
\end{eq}
Finally we consider $|C_i^{N,2}(s)-C_i(s)|$.
Using the  fact that $0 \le \frac{C_i(s)}{d} \le 1$, we have
	\begin{equation}
	\begin{split}
		\label{eq:npn_quench_difference_3}
		\E |C_i^{N,2}(s)-C_i(s)|
		&\le \E \Big[ \frac{(d-1)C_i(s)}{d} \Big| \sum_{j \in [N], j \ne i} \frac{\xi_{ji}}{D_j} -1 \Big| \Big]\\
		& \le (d-1) \Big| \sum_{j \in [N], j \ne i} \frac{\xi_{ji}}{D_j} -1 \Big|.
\end{split}	
	\end{equation}	
	Combining \eqref{eq:npn_quench_1} -- \eqref{eq:npn_quench_difference_3} with the fact that $|X_i^N(s)-X_i(s)| \le |X_i^N(s)-X_i(s)|^2$ yields
	\begin{align*}
		\max_{i \in [N]} \E \left\|X_i^N - X_i\right\|_{*,t}^2 
		& \le \kappa_7 \int_0^t \max_{i \in [N]} \E \left\|X_i^N - X_i\right\|_{*,s}^2 \, ds\\
		&\hspace{1cm}+ \kappa_7 \Big(\frac{1}{(d_{\min}(G_N))^{1/2}} + \max_{i \in [N]} \Big| \sum_{j \in [N], j \ne i} \frac{\xi_{ji}}{D_j} -1 \Big| \Big).
	\end{align*}
Theorem~\ref{th:tagged-d} (i) now follows from Gronwall's lemma and Condition \ref{cond-reg}. 

Given part (i), the proof of the propagation of chaos property as stated in Theorem~\ref{th:tagged-d} (ii) follows from standard arguments (cf. \cite{Sznitman1989}), and hence is omitted.
Also, having established the asymptotic result in Theorem~\ref{th:tagged-d} (i), the proof of convergence of local occupancy measures as stated in Theorem~\ref{th:tagged-d} (iii) can be established using similar arguments as in \cite[Corollary 3.3]{BBW16}.
\end{proof}
We now complete the proof of Theorem \ref{th:deterministic-d}.

\begin{proof}[Proof of Theorem~\ref{th:deterministic-d}]
	From the propagation of chaos property in Theorem \ref{th:tagged-d}(ii), it follows (cf.\ \cite{Sznitman1989}) that $\qq^N(\cdot)$ converges weakly with respect to the Skorohod $J_1$-topology to the deterministic limit $\tilde \qq(\cdot)$ given by $\tilde q_i(t) = \mu_t[i,\infty) = \Pro{X_i(t)\ge i}$ for all $i \in \Nmb_0$ and $t\ge 0$.
	Thus in order to prove the theorem it suffices to show that $\tilde \qq$ satisfies the system of ODEs in \eqref{eq:occ-deterministic}.
	
Define $f_j(x) = \one_{\{x\geq j\}}$, $j=1,2,\ldots$.
Then Equation~\eqref{eq:limit-tagged} yields
\begin{align*}
\E{f_j(X_i(t))} &= \E{f_j(X_i(0))}+\int_0^t \expt{\one_{\{X_i(s)>0\}}(f_j(X_i(s)-1)-f_j(X_i(s)))}\dif s\\
&\hspace{1.5cm}+ \lambda d \int_0^t \int_{\N^{d-1}}\E\Big[ b(X_i(s),x_2,\ldots,x_d)(f_j(X_i(s)+1)\\
&\hspace{3.5cm}-f_j(X_i(s)))\Big]\mu_s(\dif x_2)\ldots \mu_s(\dif x_d)\dif s\\
&=  \E{f_j(X_i(0))}-\int_0^t \expt{f_j(X_i(s))-f_{j+1}(X_i(s))}\dif s\\
&\hspace{1.5cm}+ \lambda d \int_0^t \int_{\N^{d-1}}\E \Big[b(j-1,x_2,\ldots,x_d)(f_{j-1}(X_i(s))\\
&\hspace{3.5cm}-f_j(X_i(s)))\Big]\mu_s(\dif x_2)\ldots \mu_s(\dif x_d)\dif s.
\end{align*}
Since $\E[f_j(X_i(t))]=\tilde q_j(t)$ for $j=1,2,\ldots$, we  obtain
\begin{align*}
\tilde q_j(t) &= \tilde q_j(0) - \int_0^t (\tilde q_j(s) - \tilde q_{j+1}(s))\dif s + \lambda d\int_0^t (\tilde q_{j-1}(s) - \tilde q_j(s))\\
&\hspace{3cm}\times\int_{\N^{d-1}} b(j-1,x_2,\ldots,x_d)\mu_s(\dif x_2)\ldots \mu_s(\dif x_d)\dif s\\
&= \tilde q_j(0) - \int_0^t (\tilde q_j(s) - \tilde q_{j+1}(s))\dif s + \lambda \int_0^t [(\tilde q_{j-1}(s))^d - (\tilde q_j(s))^d]\dif s,
\end{align*}
where the last equality uses the fact that 
$\Pro{X_i(t)\geq j}=\tilde q_j(t)$, $j=1,2,\ldots$.
This shows that $\tilde \qq$ satisfies the system of ODEs in \eqref{eq:occ-deterministic} and completes the proof of Theorem~\ref{th:deterministic-d}.
\end{proof}

\begin{proof}[Proof of Lemma~\ref{lem:prep_quench}]
	We first show the first inequality in \eqref{eq:prep_quench_i}.
	Observe that
	\begin{align}
		U_s & = \sum_{(j_2,\dotsc,j_d) \in \Smc_i^N} \sum_{(k_2,\dotsc,k_d) \in \Smc_i^N} \Big[  \alpha^N(i; j_2, j_3, \ldots, j_d)   \alpha^N(i; k_2, k_3, \ldots, k_d)  \Big] \notag \\
		& \hspace{3cm}\times\E \Big[ \Big(b(X_i(s),X_{j_2}(s),\dotsc,X_{j_d}(s)) - \frac{C_i(s)}{d} \Big)\\
		&\hspace{3cm}\times \Big(b(X_i(s),X_{k_2}(s),\dotsc,X_{k_d}(s)) - \frac{C_i(s)}{d} \Big) \Big]. \notag
	\end{align}
Now observe that since $\{X_i(0):i\in [N]\}$ are iid, we have $\{X_i(s):i\in [N]\}$ are also iid for any fixed $s>0$.
Thus,
	\begin{equation}
		\label{eq:prep_quench_simplify}
		\begin{split}
		&\E \Big[ \Big(b(X_i(s),X_{j_2}(s),\dotsc,X_{j_d}(s)) - \frac{C_i(s)}{d} \Big)\\ &\hspace{3cm}\times\Big(b(X_i(s),X_{k_2}(s),\dotsc,X_{k_d}(s)) - \frac{C_i(s)}{d} \Big) \Big]=0
		\end{split}
	\end{equation}
	when $(i,j_2,k_2,\dotsc,j_d,k_d)$ are distinct.
Therefore, we have
	\begin{align}
		U_s \le \sum  \alpha^N(i; j_2, j_3, \ldots, j_d) \alpha^N(i; k_2, k_3, \ldots, k_d), \label{eq:prep_quench_2}
	\end{align}
	where the summation is taken over
	\begin{equation}
		\label{eq:prep_quench_3}
		\hat \Smc_i^N \doteq \left\{(j_2,\dotsc,j_d) \in \Smc_i^N, (k_2,\dotsc,k_d) \in \Smc_i^N, (j_2,k_2,\dotsc,j_d,k_d) \mbox{ are not distinct}\right\}
	\end{equation}
	and the inequality follows since $0 \le b \le 1$ and $0 \le \frac{C_i(s)}{d} \le 1$.
	Since the total number of combinations in \eqref{eq:prep_quench_3} such that $(\xi_{ij_2}\xi_{ij_3}\dotsm\xi_{ij_d})(\xi_{ik_2}\xi_{ik_3}\dotsm\xi_{ik_d})=1$ is no more than
	\begin{equation}\label{eq:bdonbindiff}
		\left[(d-1)! \binom{D_i}{d-1}\right]^2 - (2d-2)!\binom{D_i}{2d-2} \le \kappa_1 D_i^{2d-3}, 
	\end{equation}
	we can bound \eqref{eq:prep_quench_2} by
	\begin{align*}
	\frac{\kappa_1 D_i^{2d-3}}{D_i^2(D_i-1)^2\dotsm(D_i-d+2)^2} \le \kappa_2  \frac{1}{D_i} \le \frac{\kappa_2}{d_{\min}(G_N)}.
	\end{align*}
	This gives the first bound in \eqref{eq:prep_quench_i}.
	
	Next we show the second bound in \eqref{eq:prep_quench_i}.
	From  \eqref{eq:prep_quench_simplify} it follows from the same argument used for \eqref{eq:prep_quench_2} that 
	\begin{equation}
		V_s \le \sum  \alpha^N(j_2; i, j_3, \ldots, j_d) \alpha^N(k_2; i, k_3, \ldots, k_d), \label{eq:prep_quench_4}
	\end{equation}
	where the summation is taken over \eqref{eq:prep_quench_3}.
	Since for fixed $(j_2,k_2) \in \Smcbar_i$, where
	\begin{equation}\label{eq:smcbar}
	\Smcbar_i\doteq \{ (j,k) \in [N]^2 : j \ne i, k \ne i\},
	\end{equation}
	the total number of combinations in \eqref{eq:prep_quench_3} such that 
	$$(\xi_{j_2i}\xi_{j_2j_3}\dotsm\xi_{j_2j_d})(\xi_{k_2i}\xi_{k_2k_3}\dotsm\xi_{k_2k_d})=1$$ 
	is no more than
	\begin{align}
		& \left[(d-2)! \binom{D_{j_2}-1}{d-2}\right] \left[ (d-2)! \binom{D_{k_2}-1}{d-2} \right] \nonumber\\
		&\hspace{5cm}- \left[(d-2)!\binom{D_{j_2}-2}{d-2}\right] \left[(d-2)!\binom{D_{k_2}-d}{d-2}\right] \nonumber\\
		& \le \kappa_3 (D_{j_2}^{d-3}D_{k_2}^{d-2} + D_{j_2}^{d-2}D_{k_2}^{d-3}), \label{eq:d2d3bd}
	\end{align}
	where the second term in the first line corresponds to choosing distinct $j_3,\dotsc,j_d$ from $D_{j_2}-2$ neighbors (excluding $i,k_2$) of $j_2$ and then choosing distinct $k_3,\dotsc,k_d$ from $D_{k_2}-d$ neighbors (excluding $i,j_2,\dotsc,j_d$) of $k_2$.
	Now, we can bound \eqref{eq:prep_quench_4} by
	\begin{align*}
		 &\sum_{(j_2,k_2) \in \Smcbar_i}  \frac{\kappa_3 (D_{j_2}^{d-3}D_{k_2}^{d-2} + D_{j_2}^{d-2}D_{k_2}^{d-3})\xi_{j_2i}\xi_{k_2i}}{D_{j_2}(D_{j_2}-1)\dotsm(D_{j_2}-d+2)D_{k_2}(D_{k_2}-1)\dotsm(D_{k_2}-d+2)}\\
		& \le \kappa_4 \sum_{(j_2,k_2) \in \Smcbar_i}  \left( \frac{\xi_{j_2i}\xi_{k_2i}}{D_{j_2}^2D_{k_2}} + \frac{\xi_{j_2i}\xi_{k_2i}}{D_{j_2}D_{k_2}^2} \right) \\
		& \le \kappa_4 \frac{2}{d_{\min}(G_N)} \left( \sum_{j=1,j \ne i}^N   \frac{\xi_{ji}}{D_j} \right)^2.
	\end{align*}
	This completes the proof.
\end{proof}

\subsection{Proofs for random graph sequences}
In this section we give the proofs of Theorems \ref{thm:npn_rate} and \ref{thm:npn_rate_quench}.
As in the proof of Theorem~\ref{th:tagged-d}, we will define the limiting processes $(X_i(\cdot))_{i\geq 1}$ and the pre-limit processes $(X_i^N(\cdot))_{i\geq 1}$ on the same probability space by taking identical sequence of Poisson processes $\{N_i\}$ and Poisson random measures $\{\Nbar_i\}$ in both cases. The random graph sequence $\{G_N\}$ will also be given on this common probability space and is taken to be independent of the Poisson processes and Poisson random measures.
Finally, we take $X_i^N(0) = X_i(0)$ for all $i\in [N]$, $N\geq 1$.
Once again, we will frequently suppress $N$ in the notation $D^N_i$ and write it as $D_i$.
We begin with three lemmas that will be used in the proof.
Let for $s\ge 0$
	\begin{align}
		U_s^A & \doteq \Ebf\ \Big[ \one_{\{D^N_i \ge d-1\}} \sum_{(j_2,\dotsc,j_d) \in \Smc_i^N} \alpha^N(i; j_2, j_3, \ldots, j_d)\nonumber\\   
		&\hspace{3.5cm}\Big(b(X_i(s),X_{j_2}(s),\dotsc,X_{j_d}(s)) - \frac{C_i(s)}{d} \Big) \Big]^2 
		\label{eq:prep_new_i} 
	\end{align}
and
	\begin{align}
		V_s^A & \doteq \Ebf\ \Big[ \sum_{(j_2,\dotsc,j_d) \in \Smc_i^N} \one_{\{D^N_{j_2} \ge d-1\}} \alpha^N(j_2; i, j_3, \ldots, j_d)\nonumber\\   
		&\hspace{3.5cm} \Big(b(X_i(s),X_{j_2}(s),\dotsc,X_{j_d}(s))-\frac{C_i(s)}{d}\Big) \Big]^2.
		\label{eq:prep_new_j}
	\end{align}	
	Note that the dependence of $U_s^A$ and $V_s^A$ on $i$ is suppressed in the notation.
The next lemma provides uniform bounds on $U_s^A$ and $V_s^A$.
\begin{lemma}
	\label{lem:prep_new}
Fix $T\ge 0$. Under the conditions of Theorem~\ref{thm:npn_rate}, there exists $\kappa \in (0,\infty)$ such that for every $s \in [0,T]$ and  $i \in [N]$,
	\begin{align}
		U_s^A  \le \frac{\kappa}{Np_N} \qquad\mbox{and}\qquad 
		V_s^A \le \frac{\kappa}{Np_N} + \frac{\kappa}{(Np_N)^2}.\notag 
	\end{align}
\end{lemma}
\noindent
The proof of Lemma~\ref{lem:prep_new} follows along similar lines  as the proof of Lemma~\ref{lem:prep_quench}, however note that the expectations in \eqref{eq:prep_new_i} and \eqref{eq:prep_new_j} are taken also over the randomness of the graph topology, and thus we need additional arguments. 
The proof of Lemma~\ref{lem:prep_new} is provided at the end of this subsection.

The next lemma is taken from \cite{BBW16}.
\begin{lemma}[{\cite[Lemma~5.2]{BBW16}}]
	\label{lem:prep_2}
	Let $G_N$ be an ERRG with connection probability $p_N$. Then
	\begin{equation*}
		\Ebf\ \Big( \sum_{j \in [N], j \ne i} \frac{\xi_{ij}^N}{D_j^N} \one_{\{D_j^N > 0\}} - 1 \Big)^2 \le \frac{4}{N p_N} + 2 e^{-N p_N}, \quad i \in [N],
	\end{equation*}
\end{lemma}

The following lemma provides useful moment bounds on $|X_i^N - X_i|$ and its proof is given at the end of this subsection.
\begin{lemma}
	\label{lem:momentbd}
	Fix $T\ge 0$. Under the conditions of Theorem~\ref{thm:npn_rate},
	\begin{equation*} 
		\sup_{N \ge 1} \max_{i \in [N]} \Ebf \left\|X_i^N - X_i\right\|_{*,T}^4 < \infty.
	\end{equation*}
\end{lemma}

We now present the proof of Theorem \ref{thm:npn_rate}.

\begin{proof}[Proof of Theorem \ref{thm:npn_rate}]

Fix any $i\in\N$ and $T>0$.
From \eqref{eq:X_i_n} and \eqref{eq:limit-tagged}, using Cauchy--Schwarz and Doob's inequalities we have for any fixed $t\in [0,T]$	
\begin{align}
	\Ebf \left\|X_i^N - X_i\right\|_{*,t}^2 & \le 
	\kappa_1 \int_0^t \Ebf |X_i^N(s) - X_i(s)|^2 \, ds + \kappa_1 \int_0^t \Ebf |C_i^N(s)-C_i(s)| \, ds\notag \\
	&\hspace{3cm}+ \kappa_1 \int_0^t \Ebf |C_i^N(s)-C_i(s)|^2 \, ds \label{eq:npn_1}
\end{align}
for some  $\kappa_1\in (0,\infty)$.
Define $C^{N,1}_i(s)$ and $C^{N,2}_i(s)$ by
\begin{align*}
	&C_i^{N,1}(s)  = \one_{\{D_i<d-1\}} \bbar_i((X_k^N(s))_{k \in [N]},(\xi_{kl})_{k,l \in [N]}) \\
	& \quad + \one_{\{D_i \ge d-1\}} \sum_{(j_2,\dotsc,j_d) \in \Smc_i^N} \alpha^N(i; j_2, j_3, \ldots, j_d) b(X_i(s),X_{j_2}(s),\dotsc,X_{j_d}(s)) \\
	& \quad + (d-1) \sum_{(j_2,\dotsc,j_d) \in \Smc_i^N} \one_{\{D_{j_2} \ge d-1\}} \alpha^N(j_2; i, j_3, \ldots, j_d) b(X_i(s),X_{j_2}(s),\dotsc,X_{j_d}(s)) \\
	& \quad + \sum_{j_2 \in [N], j_2 \ne i} \one_{\{D_{j_2<d-1}\}}\xi_{ij_2} \bbar_{ij_2}((X_k^N(s))_{k \in [N]},(\xi_{kl})_{k,l \in [N]})
\end{align*}
and
\begin{align*}
	C_i^{N,2}(s) & = \one_{\{D_i<d-1\}} \bbar_i((X_k^N(s))_{k \in [N]},(\xi_{kl})_{k,l \in [N]}) \\
	& \quad + \one_{\{D_i \ge d-1\}} \sum_{(j_2,\dotsc,j_d) \in \Smc_i^N} \alpha^N(i; j_2, j_3, \ldots, j_d) \frac{C_i(s)}{d} \\
	& \quad + (d-1) \sum_{(j_2,\dotsc,j_d) \in \Smc_i^N} \one_{\{D_{j_2} \ge d-1\}} \alpha^N(j_2; i, j_3, \ldots, j_d) \frac{C_i(s)}{d} \\
	& \quad + \sum_{j_2 \in [N], j_2 \ne i} \one_{\{D_{j_2<d-1}\}}\xi_{ij_2} \bbar_{ij_2}((X_k^N(s))_{k \in [N]},(\xi_{kl})_{k,l \in [N]}).
\end{align*}
By adding and subtracting terms we have ~\eqref{eq:npn_2} and
\begin{equation}
\begin{split}
	|C_i^N(s)-C_i(s)|^2 &\le 3|C_i^N(s)-C_i^{N,1}(s)|^2 + 3|C_i^{N,1}(s)-C_i^{N,2}(s)|^2 \\
	&\hspace{3cm}+ 3|C_i^{N,2}(s)-C_i(s)|^2. \label{eq:npn_3}
	\end{split}
\end{equation}
Here although one has $\Ebf |C_i^N(s)-C_i(s)| \le \left( \Ebf |C_i^N(s)-C_i(s)|^2 \right)^{1/2}$, in order to get the desired rate $\sqrt{Np_N}$ in \eqref{eq:npn_rate}, we have to estimate $\Ebf |C_i^N(s)-C_i(s)|$ more carefully through \eqref{eq:npn_2}.

Let us consider $|C_i^N(s)-C_i^{N,1}(s)|$ and $|C_i^N(s)-C_i^{N,1}(s)|^2$ first.
We claim that for $m=1,2$, there exists some $\kappa_2 \in (0,\infty)$ such that
\begin{align}
	\Ebf |C_i^N(s)-C_i^{N,1}(s)|^m & \le \kappa_2 \Ebf |X_i^N(s)-X_i(s)|^m +\kappa_2 \left(\frac{1}{Np_N} + e^{-Np_N}\right)^{1/2}  \notag \\
	& \quad + \kappa_2 \Ebf \Big[ \one_{\{D_i \ge d-1\}} \sum_{j \in [N], j \ne i} \frac{\xi_{ij}}{D_i} |X_j^N(s)-X_j(s)|^m \Big]. \label{eq:CN_CN1}
\end{align}
To see this, note that from the Lipschitz property of $b$ and the definition of $\Smc_i^N$ we have
	\begin{align*}
		&\Ebf |C_i^N(s)-C_i^{N,1}(s)| \\
		& \le \Ebf \Big[ \one_{\{D_i \ge d-1\}} \sum_{(j_2,\dotsc,j_d) \in \Smc_i^N} \alpha^N(i; j_2, j_3, \ldots, j_d) \notag \\ 
		& \qquad (|X_i^N(s)-X_i(s)| + |X_{j_2}^N(s)-X_{j_2}(s)| + \dotsb + |X_{j_d}^N(s)-X_{j_d}(s)|) \notag \\
		& \qquad + (d-1) \sum_{(j_2,\dotsc,j_d) \in \Smc_i^N} \one_{\{D_{j_2} \ge d-1\}} \alpha^N(j_2; i, j_3, \ldots, j_d) \notag \\
		& \qquad  (|X_i^N(s)-X_i(s)| + |X_{j_2}^N(s)-X_{j_2}(s)| + \dotsb + |X_{j_d}^N(s)-X_{j_d}(s)|) \Big], \notag \\
		& = d\ \Ebf \Big[ \one_{\{D_i \ge d-1\}} \sum_{(j_2,\dotsc,j_d) \in \Smc_i^N} \alpha^N(i; j_2, j_3, \ldots, j_d)   \notag \\ 
		& \qquad  (|X_i^N(s)-X_i(s)| + |X_{j_2}^N(s)-X_{j_2}(s)| + \dotsb + |X_{j_d}^N(s)-X_{j_d}(s)|) \Big] \notag \\
		& \le d\ \Ebf |X_i^N(s)-X_i(s)| + d(d-1) \Ebf \Big[ \one_{\{D_i \ge d-1\}} \sum_{j \in [N], j \ne i} \frac{\xi_{ij}}{D_i} |X_j^N(s)-X_j(s)| \Big], 
	\end{align*}
	where in obtaining the equality we have used the exchangeability property:
	\begin{equation}
	\begin{split}
		& \Lmc(\xi_{ij_2},\xi_{ij_3},\dotsc,\xi_{ij_d},D_i, X_i^N(s),X_i(s),X_{j_2}^N(s),X_{j_2}(s),\\
		&\hspace{5cm} X_{j_3}^N(s),X_{j_3}(s),\dotsc,X_{j_d}^N(s),X_{j_d}(s)) \\
		& = \Lmc(\xi_{j_2i},\xi_{j_2j_3},\dotsc,\xi_{j_2j_d},D_{j_2}, X_{j_2}^N(s),X_{j_2}(s),X_i^N(s),X_i(s),\\
		&\hspace{5cm} X_{j_3}^N(s),X_{j_3}(s),\dotsc,X_{j_d}^N(s),X_{j_d}(s)) \label{eq:exchangeability}
	\end{split}
	\end{equation}	
	for $(j_2,\dotsc,j_d) \in \Smc_i^N$.
	Therefore the claim \eqref{eq:CN_CN1} holds for $m=1$.
	Next we verify \eqref{eq:CN_CN1} when $m=2$.
	Note that
	\begin{align*}
		\Ebf |C_i^N(s)-C_i^{N,1}(s)|^2 & \le 2R_i^{N,1}(s) + 2(d-1)^2R_i^{N,2}(s),
	\end{align*}
	where
	\begin{align*}
		R_i^{N,1}(s) & \doteq \Ebf \Big[ \one_{\{D_i \ge d-1\}} \sum_{(j_2,\dotsc,j_d) \in \Smc_i^N} \alpha^N(i; j_2, j_3, \ldots, j_d) \\
		& \qquad  [b(X_i(s),X_{j_2}(s),\dotsc,X_{j_d}(s)) - b(X_i(s),X_{j_2}(s),\dotsc,X_{j_d}(s))] \Big]^2, \notag \\
		R_i^{N,2}(s) & \doteq \Ebf \Big[ \sum_{(j_2,\dotsc,j_d) \in \Smc_i^N} \one_{\{D_{j_2} \ge d-1\}} \alpha^N(j_2; i, j_3, \ldots, j_d) \\
		& \qquad  [b(X_i(s),X_{j_2}(s),\dotsc,X_{j_d}(s)) - b(X_i(s),X_{j_2}(s),\dotsc,X_{j_d}(s))] \Big]^2. \notag	
	\end{align*}
	From the Lipschitz property of $b$, the definition of $\Smc_i^N$ and Cauchy-Schwarz inequality we have
	\begin{align*}		
		&R_i^{N,1}(s)  \le \Ebf \Big[ \one_{\{D_i \ge d-1\}} \sum_{(j_2,\dotsc,j_d) \in \Smc_i^N} \alpha^N(i; j_2, j_3, \ldots, j_d) \\
		& \qquad\qquad  (|X_i^N(s)-X_i(s)| + |X_{j_2}^N(s)-X_{j_2}(s)| + \dotsb + |X_{j_d}^N(s)-X_{j_d}(s)|) \Big]^2 \notag \\
		& = \Ebf \Big[ \one_{\{D_i \ge d-1\}} \Big( |X_i^N(s)-X_i(s)| +  (d-1) \sum_{j \in [N], j \ne i} \frac{\xi_{ij}}{D_i} |X_j^N(s)-X_j(s)| \Big) \Big]^2 \\
		& \le 2 \Ebf |X_i^N(s)-X_i(s)|^2 + 2(d-1)^2 \Ebf \Big[ \one_{\{D_i \ge d-1\}} \sum_{j \in [N], j \ne i} \frac{\xi_{ij}}{D_i} |X_j^N(s)-X_j(s)|^2 \Big].
	\end{align*}	
	From Cauchy-Schwarz inequality we have
	\begin{align*}
		R_i^{N,2}(s) & \le \Ebf \Big( \Big[ \sum_{(j_2,\dotsc,j_d) \in \Smc_i^N} \one_{\{D_{j_2} \ge d-1\}} \alpha^N(j_2; i, j_3, \ldots, j_d) \Big]\\
		&\times \Big[ \sum_{(j_2,\dotsc,j_d) \in \Smc_i^N} \one_{\{D_{j_2} \ge d-1\}} \alpha^N(j_2; i, j_3, \ldots, j_d) \\
		& \qquad\times  [b(X_i(s),X_{j_2}(s),\dotsc,X_{j_d}(s)) - b(X_i(s),X_{j_2}(s),\dotsc,X_{j_d}(s))]^2 \Big] \Big)
	\end{align*}
	\begin{align*}
		& = \Ebf \Big[ \sum_{(j_2,\dotsc,j_d) \in \Smc_i^N} \one_{\{D_{j_2} \ge d-1\}} \alpha^N(j_2; i, j_3, \ldots, j_d) \\
		& \qquad  [b(X_i(s),X_{j_2}(s),\dotsc,X_{j_d}(s)) - b(X_i(s),X_{j_2}(s),\dotsc,X_{j_d}(s))]^2 \Big] \notag \\
		& \quad + \Ebf \Big( \Big[ \sum_{j \in [N], j \ne i} \one_{\{D_j \ge d-1\}} \frac{\xi_{ji}}{D_j} - 1 \Big]\\
		&\qquad\times \Big[ \sum_{(j_2,\dotsc,j_d) \in \Smc_i^N} \one_{\{D_{j_2} \ge d-1\}} \alpha^N(j_2; i, j_3, \ldots, j_d) \\
		& \qquad  [b(X_i(s),X_{j_2}(s),\dotsc,X_{j_d}(s)) - b(X_i(s),X_{j_2}(s),\dotsc,X_{j_d}(s))]^2 \Big] \Big) \notag \\
		& \doteq R_i^{N,3}(s) + R_i^{N,4}(s),
	\end{align*}
	where the equality follows by adding and subtracting one in the first term.
	From the Lipschitz property of $b$, the definition of $\Smc_i^N$ and the exchangeability property \eqref{eq:exchangeability} we have
	\begin{align*}
		&R_i^{N,3}(s)  \le d^2\ \Ebf \Big[ \sum_{(j_2,\dotsc,j_d) \in \Smc_i^N} \one_{\{D_{j_2} \ge d-1\}} \alpha^N(j_2; i, j_3, \ldots, j_d) \notag \\
		& \qquad\quad \times (|X_i^N(s)-X_i(s)|^2 + |X_{j_2}^N(s)-X_{j_2}(s)|^2 + \dotsb + |X_{j_d}^N(s)-X_{j_d}(s)|^2) \Big] \notag \\
		& = d^2\ \Ebf \Big[ \sum_{(j_2,\dotsc,j_d) \in \Smc_i^N} \one_{\{D_i \ge d-1\}} \alpha^N(i; j_2, j_3, \ldots, j_d) \notag \\
		& \qquad \cdot (|X_i^N(s)-X_i(s)|^2 + |X_{j_2}^N(s)-X_{j_2}(s)|^2 + \dotsb + |X_{j_d}^N(s)-X_{j_d}(s)|^2) \Big] \notag \\
		& \le d^2\ \Ebf |X_i^N(s)-X_i(s)|^2 + d^2(d-1) \Ebf \Big[ \one_{\{D_i \ge d-1\}} \sum_{j \in [N], j \ne i} \frac{\xi_{ij}}{D_i} |X_j^N(s)-X_j(s)|^2 \Big]. 
	\end{align*}
	From the fact that $\|b\|_\infty \le 1$ we have
	\begin{align*}
		R_i^{N,4}(s) & \le \Ebf \Big( \Big| \sum_{j \in [N], j \ne i} \one_{\{D_j > 0\}} \frac{\xi_{ji}}{D_j} - 1 \Big|\\
		&\hspace{2cm}\times \Big[ 4 \sum_{(j_2,\dotsc,j_d) \in \Smc_i^N} \one_{\{D_{j_2} \ge d-1\}} \alpha^N(j_2; i, j_3, \ldots, j_d) \Big] \Big) \notag \\
		& \le 4 \Ebf \Big( \Big| \sum_{j \in [N], j \ne i} \one_{\{D_j > 0\}} \frac{\xi_{ji}}{D_j} - 1 \Big| \sum_{j \in [N], j \ne i} \one_{\{D_j > 0\}} \frac{\xi_{ji}}{D_j} \Big) \notag \\
		& \le \kappa_3 \left( \frac{1}{Np_N} + e^{-Np_N} \right)^{1/2}.
	\end{align*}
	where the last inequality follows from Lemma \ref{lem:prep_2} and Condition \ref{cond:errg1}.
	Combining the above estimates on $R_i^{N,k}(s)$ for $k=1,2,3,4$ gives the claim \eqref{eq:CN_CN1} when $m=2$.
	
	Now using the exchangeability property:
	\begin{align*}
		\Lmc(\xi_{ij}, D_i, X_j^N(s), X_j(s)) & = \Lmc(\xi_{ji}, D_j, X_i^N(s), X_i(s)), \quad i \ne j,
	\end{align*}	
	we have for $m=1,2$,
	\begin{align*}
		& \quad \Ebf \Big[ \one_{\{D_i \ge d-1\}} \sum_{j \in [N], j \ne i} \frac{\xi_{ij}}{D_i} |X_j^N(s)-X_j(s)|^m \Big] \\
		& = \Ebf \Big[ \sum_{j \in [N], j \ne i} \one_{\{D_j \ge d-1\}} \frac{\xi_{ji}}{D_j} |X_i^N(s)-X_i(s)|^m \Big] \\
		& \le \Ebf \Big[ \Big( \sum_{j \in [N], j \ne i} \one_{\{D_j > 0\}}  \frac{\xi_{ji}}{D_j} - 1 \Big) |X_i^N(s)-X_i(s)|^m \Big] + \Ebf |X_i^N(s)-X_i(s)|^m 
	\end{align*}
	\begin{align*}
		& \le \Big[ \Ebf \Big( \sum_{j \in [N], j \ne i} \one_{\{D_j > 0\}}  \frac{\xi_{ji}}{D_j} - 1 \Big)^2 \Ebf |X_i^N(s)-X_i(s)|^{2m} \Big]^{1/2} + \Ebf |X_i^N(s)-X_i(s)|^m \\
		& \le \kappa_4 \Big(\frac{1}{Np_N} + e^{-Np_N}\Big)^{1/2} + \Ebf |X_i^N(s)-X_i(s)|^m, 
	\end{align*}	
	where the second inequality follows from Cauchy-Schwarz inequality and the last line follows from Lemmas \ref{lem:prep_2} and \ref{lem:momentbd}.
	Combining this, \eqref{eq:CN_CN1} with the fact that $|X_i^N(s)-X_i(s)| \le |X_i^N(s)-X_i(s)|^2$ gives
	\begin{equation}
	\begin{split}
		\label{eq:npn_difference_1}
		\Ebf |C_i^N(s)-C_i^{N,1}(s)| &+ \Ebf |C_i^N(s)-C_i^{N,1}(s)|^2 \\
		&\le \kappa_5 \Ebf |X_i^N(s)-X_i(s)|^2 + \kappa_5 \left(\frac{1}{Np_N} + e^{-Np_N}\right)^{1/2}.
		\end{split}
	\end{equation}
	
	Next we consider $|C_i^{N,1}(s)-C_i^{N,2}(s)|^2$.
	From the inequality $(a+b)^2 \le 2a^2+2b^2$, it follows  that
	\begin{equation}
	\begin{split}
		\left( \Ebf |C_i^{N,1}(s)-C_i^{N,2}(s)| \right)^2 &\le \Ebf |C_i^{N,1}(s)-C_i^{N,2}(s)|^2\\
		 &\le 2U_s^A + 2(d-1)^2V_s^A \le \frac{\kappa_6}{Np_N} + \frac{\kappa_6}{(Np_N)^2},\label{eq:npn_difference_2}
	\end{split}
	\end{equation}
	where $U_s^A$ and $V_s^A$ were introduced in \eqref{eq:prep_new_i} and \eqref{eq:prep_new_j} and the last inequality is from Lemma
	\ref{lem:prep_new}.

	Finally we consider $|C_i^{N,2}(s)-C_i(s)|^2$.
	Note that $C_i^{N,2}(s)$ can be rewritten as
	\begin{align*}
		C_i^{N,2}(s) & = \one_{\{D_i<d-1\}} \bbar_i((X_k^N(t))_{k \in [N]},(\xi_{kl})_{k,l \in [N]}) \\
		& \quad + \one_{\{D_i \ge d-1\}} \frac{C_i(s)}{d} + (d-1) \sum_{j \in [N], j \ne i} \one_{\{D_{j} \ge d-1\}} \frac{\xi_{ji}}{D_j} \frac{C_i(s)}{d} \\
		& \quad + \sum_{j \in [N], j \ne i} \one_{\{D_{j}<d-1\}}\xi_{ij} \bbar_{ij}((X_k^N(t))_{k \in [N]},(\xi_{kl})_{k,l \in [N]}).
	\end{align*}
	Using the Cauchy-Schwarz inequality and the fact that $0 \le \frac{C_i(s)}{d} \le 1$, we have
	\begin{align*}
		& \Ebf |C_i^{N,2}(s)-C_i(s)|^2 \notag \\
		& \le 5\Ebf \Big[ \one_{\{D_i<d-1\}} (D_i+1) \Big]^2 + 5\Ebf \Big[ \one_{\{D_i<d-1\}} \frac{C_i(s)}{d} \Big]^2 \\
			&\hspace{2cm} + 5\Ebf \left[(d-1)\sum_{j \in [N], j \ne i} \one_{\{0 < D_{j} < d-1\}} \frac{\xi_{ji}}{D_j} \frac{C_i(s)}{d}\right]^2 \notag \\
			& \hspace{2cm} + 5\Ebf \Big[ \frac{(d-1)C_i(s)}{d} \Big| \sum_{j \in [N], j \ne i} \one_{\{D_{j} > 0\}} \frac{\xi_{ji}}{D_j} -1 \Big| \Big]^2 \\
			&\hspace{2cm} + 5\Ebf \Big[ \sum_{j \in [N], j \ne i} \one_{\{D_{j}<d-1\}}\xi_{ij} (D_i+1) \Big]^2 \notag \\
		& \le 5(d^2+1) \Pbf(D_i < d-1) + 5(d-1)^2\Ebf \left[\sum_{j \in [N], j \ne i} \one_{\{0 < D_{j} < d-1\}} \frac{\xi_{ji}}{D_j}\right]^2 \notag \\
		& \hspace{3cm} + 5(d-1)^2 \Ebf \Big[ \sum_{j \in [N], j \ne i} \one_{\{D_{j} > 0\}} \frac{\xi_{ji}}{D_j} -1 \Big]^2 \\
		&\hspace{3cm}+ 5\Ebf \Big[ \sum_{j \in [N], j \ne i} \one_{\{D_{j}<d-1\}}\xi_{ij} (D_i+1) \Big]^2. 
	\end{align*}	
	For the second and last terms on the right hand side, we have
	\begin{align*}
		&\Ebf \left[\sum_{j \in [N], j \ne i} \one_{\{0 < D_{j} < d-1\}} \frac{\xi_{ji}}{D_j}\right]^2
		 \le \Ebf \left[\sum_{j \in [N], j \ne i} \one_{\{0 < D_{j} < d-1\}} \xi_{ji}\right]^2 \\
		& \le \Ebf \Big[ \sum_{j \in [N], j \ne i} \one_{\{D_{j}<d-1\}}\xi_{ij} (D_i+1) \Big]^2 \\
		& \le \Ebf \Big[ \sum_{j \in [N], j \ne i} \one_{\{D_{j}<d-1\}}\xi_{ij} (D_i+1)^2 \Big] \Big[ \sum_{j \in [N], j \ne i} \xi_{ij} \Big] \\
		& = \sum_{j \in [N], j \ne i} \Ebf \Big[ \one_{\{D_{j}<d-1\}}\xi_{ij} (D_i+1)^2 D_i \Big] \\
		& = \sum_{j \in [N], j \ne i} \Ebf \Big[ \one_{\{D_{j}-\xi_{ij}+1<d-1\}} (D_i-\xi_{ij}+2)^2 (D_i-\xi_{ij}+1) \Big] p_N \\
		& \le \kappa_7 (N-1) \Pbf(D_i<d) (Np_N+1)^3p_N 
	\end{align*}
	where the third inequality follows from Cauchy-Schwarz inequality, the second equality follows by conditioning on $\xi_{ij}=1$, and the last inequality follows from independence and Condition \ref{cond:errg1}.	
	Note that
	\begin{align}
 		\Pbf(D_i<d) 
		& = \sum_{k=0}^{d-1}\binom{N-1}{k} p_N^k (1-p_N)^{N-1-k}  \notag \\
		& \le \kappa_8 (1-p_N)^{N-d} \left[ 1 + Np_N + \dotsb + (Np_N)^{d-1} \right] \notag \\
		& \le \kappa_9 [1+(Np_N)^{d-1}] e^{-(N-d)p_N}.
		\label{eq:CN_bd_2}
	\end{align}	
	Combining above three estimates with Lemma \ref{lem:prep_2} gives
	\begin{eq}
		\label{eq:npn_difference_3}
		\left( \Ebf |C_i^{N,2}(s)-C_i(s)| \right)^2 &\le \Ebf |C_i^{N,2}(s)-C_i(s)|^2\\
		 &\hspace{-1cm}\le \kappa_0 [1+(Np_N)^{d+3}] e^{-Np_N} + \kappa_0 \Big(\frac{1}{Np_N} + e^{-Np_N}\Big).
	\end{eq}	
	Combining \eqref{eq:npn_2}, \eqref{eq:npn_1}, \eqref{eq:npn_3}, \eqref{eq:npn_difference_1}, \eqref{eq:npn_difference_2}, \eqref{eq:npn_difference_3} and Condition \ref{cond:errg1} gives us
	\begin{align*}
		\max_{i \in [N]} \sqrt{Np_N} \Ebf \left\|X_i^N - X_i\right\|_{*,t}^2 & \le \kappa \int_0^t \max_{i \in [N]} \sqrt{Np_N} \Ebf \left\|X_i^N - X_i\right\|_{*,s}^2 \, ds + \kappa.
	\end{align*}
	Part (i) of the theorem now follows from Gronwall's lemma.
	
	The proof of propagation of chaos property as stated in Theorem~\ref{thm:npn_rate} (ii) follows now from standard arguments (cf.\ \cite{Sznitman1989}), and hence is omitted.
	Also, having proved Theorem~\ref{thm:npn_rate} (i), the proof of convergence of local occupancy measures as stated in Theorem~\ref{thm:npn_rate} (iii) can be established using similar arguments as in \cite[Corollary 3.3]{BBW16}.
\end{proof}

We now complete the proof of Theorem \ref{thm:npn_rate_quench}.

\begin{proof}[Proof of Theorem~\ref{thm:npn_rate_quench}]
In order to prove the theorem it suffices, in view of Theorems \ref{th:deterministic-d} and~\ref{th:tagged-d}, to show that  if  $\{G_N\}$ satisfies Condition~\ref{cond:errg2}, then it satisfies Condition~\ref{cond-reg} a.s.

Using the Chernoff inequality (cf. \cite[Theorem 2.4]{CL06}), it follows that for every $x \ge 0$ and $N \in \Nmb$,
\begin{equation*}
	\Pbf(|D_i^N-\Ebf D_i^N| \ge x) \le 2 \exp \left\{ -\frac{x^2}{2\Ebf D_i^N + 2x/3}\right\}.
\end{equation*}
Let $k(N)\doteq Np_N/\ln(N).$ Note that by Condition \ref{cond:errg2}, $k(N) \to \infty$ as $N\to\infty$.
Since $\Ebf D_i^N = (N-1)p_N$ taking $x=x(N) =\ln(N)(k(N))^{3/4}$ in the above expression yields, for some $\kappa_1 \in (0,\infty)$,
\begin{equation}\label{eq:chernoff}
\begin{split}
\Pbf(|D_i^N-Np_N| \ge x(N)) & \le \Pbf(|D_i^N-\Ebf D_i^N| \ge x(N)-p_N) \\
	& \le 2 \exp \Big\{ -\frac{(x(N)-p_N)^2}{2(N-1)p_N + 2(x(N)-p_N)/3}\Big\}\\
	& \le \kappa_1\exp \Big\{ -\kappa_1\frac{(x(N))^2}{Np_N}\Big\},
\end{split}
\end{equation}
for sufficiently large $N$.
Thus
\begin{equation}
	\Pbf\left(\bigcup_{i\in [N]}\left\{|D_i^N-Np_N| \ge x(N)\right\}\right) \le \kappa_1N\exp \Big\{ -\kappa_1\frac{(x(N))^2}{Np_N}\Big\}.
	\label{eq:uninbd}
\end{equation}
From the choice of $x(N)$, we have  $(x(N))^2/[Np_N\ln(N)]\to \infty$, as $N\to\infty$.
Therefore, the right side of~\eqref{eq:uninbd} is summable over $N$.
From  the Borel-Cantelli lemma we conclude a.s., for all sufficiently large~$N$ and all $i \in [N]$
\begin{equation*}
	|D_i^N-Np_N| \le x(N),
\end{equation*}
and therefore for all such $N$
\begin{equation}
		\label{eq:degree_key}
		Np_N - x(N) \le d_{\min}(G_N) \le d_{\max}(G_N) \le Np_N + x(N).
\end{equation}
Finally, observe that
\begin{equation}\label{eq:temp-10.41}
\frac{x(N)}{Np_N} = \frac{\ln(N) (k(N))^{3/4}}{k(N) \ln(N)}= \frac{1}{(k(N))^{1/4}}\to 0 \quad\mbox{as }N\to\infty.
\end{equation}
Combining \eqref{eq:degree_key} and~\eqref{eq:temp-10.41}, $d_{\min}(G_N)  \to \infty$ and
$$\frac{d_{\max}(G_N)-d_{\min}(G_N)}{d_{\min}(G_N)} = \frac{2x(N)}{Np_N-x(N)}\to 0,$$
as $N\to\infty$. 
This together with Remark \ref{rmk:weaker_condition} shows that Condition~\ref{cond-reg} holds for $\{G_N\}$ a.s., completing the proof of Theorem \ref{thm:npn_rate_quench}.
\end{proof}
We now complete the proof of Lemma~\ref{lem:prep_new}. We begin with the following lemma from \cite{BBW16}.

\begin{lemma}[{\cite[Lemma~5.1]{BBW16}}]
	\label{lem:prep_1}
	Let $X$ be a Binomial random variable with number of trials $N$ and probability of success $p$.
	Let $q \doteq 1 - p$.
	Then for each $m \in \Nmb$,
	\begin{equation*}
		\E \left[ \one_{\{X>0\}} \frac{1}{(2X)^m} \right] \le \E \frac{1}{(X+1)^m} \le \frac{m^m}{(N+1)^mp^m}.
	\end{equation*}
\end{lemma}

\begin{proof}[Proof of Lemma~\ref{lem:prep_new}]
As before, we will omit the superscript in the $\xi_{ij}$'s and $D_i$'s for notational convenience.
	We first show \eqref{eq:prep_new_i}.
	From the independence between $\{X_i\}$ and $\{\xi_{ij}\}$ it follows that
	\begin{align}
		U_s^A & = \sum_{\substack{(k_2,\dotsc,k_d) \in \Smc_i^N\\(j_2,\dotsc,j_d) \in \Smc_i^N}}\Ebf \left[ \one_{\{D_i \ge d-1\}} \alpha^N(i; j_2, j_3, \ldots, j_d) \alpha^N(i; k_2, k_3, \ldots, k_d)\right] \notag \\
		& \hspace{3.5cm}\times \Ebf \Big[ \Big(b(X_i(s),X_{j_2}(s),\dotsc,X_{j_d}(s)) - \frac{C_i(s)}{d} \Big)\notag\\
		&\hspace{4.5cm}\times \Big(b(X_i(s),X_{k_2}(s),\dotsc,X_{k_d}(s)) - \frac{C_i(s)}{d} \Big) \Big]. \notag
	\end{align}
	Noting that
	\begin{equation}
	\begin{split}
		\label{eq:prep_new_simplify}
		&\Ebf \Big[ \Big(b(X_i(s),X_{j_2}(s),\dotsc,X_{j_d}(s)) - \frac{C_i(s)}{d} \Big)\\
		&\hspace{3cm}\times \Big(b(X_i(s),X_{k_2}(s),\dotsc,X_{k_d}(s)) - \frac{C_i(s)}{d} \Big) \Big]=0
		\end{split}
	\end{equation}
	when $(i,j_2,k_2,\dotsc,j_d,k_d)$ are distinct, we have
	\begin{align}
		U_s^A & = \sum \Ebf \left[ \one_{\{D_i \ge d-1\}} \alpha^N(i; j_2, j_3, \ldots, j_d) \alpha^N(i; k_2, k_3, \ldots, k_d)\right] \notag \\
		& \hspace{4cm} \Ebf \Big[ \Big(b(X_i(s),X_{j_2}(s),\dotsc,X_{j_d}(s)) - \frac{C_i(s)}{d} \Big)\notag\\ 
		&\hspace{5cm}\Big(b(X_i(s),X_{k_2}(s),\dotsc,X_{k_d}(s)) - \frac{C_i(s)}{d} \Big) \Big] \notag \\
		& \le \Ebf \left[ \sum \one_{\{D_i \ge d-1\}} \alpha^N(i; j_2, j_3, \ldots, j_d) \alpha^N(i; k_2, k_3, \ldots, k_d)\right], \label{eq:prep_new_2}
	\end{align}
	where the summation is taken over the collection $\hat \Smc_i^N$ defined in \eqref{eq:prep_quench_3}
	%
	and the inequality follows since $0 \le b \le 1$ and $0 \le \frac{C_i(s)}{d} \le 1$.
	As noted in \eqref{eq:bdonbindiff},  the total number of combinations in \eqref{eq:prep_quench_3} such that $(\xi_{ij_2}\xi_{ij_3}\dotsm\xi_{ij_d})(\xi_{ik_2}\xi_{ik_3}\dotsm\xi_{ik_d})=1$ is no more than
	$ \kappa_1 D_i^{2d-3}$ and thus
	we can bound \eqref{eq:prep_new_2} by
	\begin{align*}
		\Ebf \left[ \one_{\{D_i \ge d-1\}} \frac{\kappa_1 D_i^{2d-3}}{D_i^2(D_i-1)^2\dotsm(D_i-d+2)^2} \right] \le \kappa_2 \Ebf \left[ \one_{\{D_i > 0\}} \frac{1}{D_i} \right] \le \frac{2\kappa_2}{Np_N},
	\end{align*}
	where the last inequality uses Lemma \ref{lem:prep_1}.
	This gives the first inequality in Lemma \ref{lem:prep_new}.
	
	Next we show the second inequality in Lemma \ref{lem:prep_new}.
	From the independence between $\{X_i\}$ and $\{\xi_{ij}\}$ and \eqref{eq:prep_new_simplify} it follows from the same argument used for \eqref{eq:prep_new_2} that 
	\begin{equation}
		V_s^A \le \Ebf \left[ \sum \one_{\{D_{j_2} \ge d-1\}} \one_{\{D_{k_2} \ge d-1\}} \alpha^N(j_2; i, j_3, \ldots, j_d) \alpha^N(k_2; i, k_3, \ldots, k_d)\right], \label{eq:prep_new_4}
	\end{equation}
	where the summation is taken over $\hat \Smc_i^N$ defined in \eqref{eq:prep_quench_3}.
	As noted in \eqref{eq:d2d3bd},  for fixed $(j_2,k_2) \in \Smcbar_i$ with $\Smcbar_i$ as in~\eqref{eq:smcbar}, the total number of combinations in $\hat \Smc_i^N$  such that 
	\[(\xi_{j_2i}\xi_{j_2j_3}\dotsm\xi_{j_2j_d})(\xi_{k_2i}\xi_{k_2k_3}\dotsm\xi_{k_2k_d})=1\] is no more than
	$\kappa_3 (D_{j_2}^{d-3}D_{k_2}^{d-2} + D_{j_2}^{d-2}D_{k_2}^{d-3})$
	we can bound \eqref{eq:prep_new_4} by
	\begin{align}
		& \Ebf \Big[ \sum_{(j_2,k_2) \in \Smcbar_i} \one_{\{D_{j_2} \ge d-1\}} \one_{\{D_{k_2} \ge d-1\}} \frac{\kappa_3 (D_{j_2}^{d-3}D_{k_2}^{d-2} + D_{j_2}^{d-2}D_{k_2}^{d-3})\xi_{j_2i}\xi_{k_2i}}{D_{j_2}(D_{j_2}-1)\dotsm(D_{j_2}-d+2)}\Big] \notag \\
		&\hspace{5cm}\times\frac{1}{D_{k_2}(D_{k_2}-1)\dotsm(D_{k_2}-d+2)}\\
		& \le \kappa_4 \sum_{(j_2,k_2) \in \Smcbar_i} \Ebf \left[ \one_{\{D_{j_2} \ge d-1\}} \one_{\{D_{k_2} \ge d-1\}} \left( \frac{\xi_{j_2i}\xi_{k_2i}}{D_{j_2}^2D_{k_2}} + \frac{\xi_{j_2i}\xi_{k_2i}}{D_{j_2}D_{k_2}^2} \right) \right] \notag \\
		& = 2 \kappa_4 \sum_{(j,k) \in \Smcbar_i} \Ebf \left[ \one_{\{D_{j} \ge d-1\}} \one_{\{D_{k} \ge d-1\}} \frac{\xi_{ji}\xi_{ki}}{D_{j}^2D_{k}} \right]. \label{eq:prep_new_5}
	\end{align}
	Now for $(j,k) \in \Smcbar_i$ with $j \ne k$, we have
	\begin{align*}
		& \Ebf \left[ \one_{\{D_{j} \ge d-1\}} \one_{\{D_{k} \ge d-1\}} \frac{\xi_{ji}\xi_{ki}}{D_{j}^2D_{k}} \right] \\
		& = \Ebf \left[ \one_{\{\xi_{jk}=1\}} \one_{\{D_{j} \ge d-1\}} \one_{\{D_{k} \ge d-1\}} \frac{\xi_{ji}\xi_{ki}}{D_{j}^2D_{k}} \right] \\
		&\hspace{5cm}+ \Ebf \left[ \one_{\{\xi_{jk}=0\}} \one_{\{D_{j} \ge d-1\}} \one_{\{D_{k} \ge d-1\}} \frac{\xi_{ji}\xi_{ki}}{D_{j}^2D_{k}} \right] \\
		& \le \Ebf \left[ \frac{\xi_{ji}\xi_{ki}}{(D_{j}-\xi_{jk}+1)^2(D_{k}-\xi_{jk}+1)} \right] \\
		&\hspace{3.5cm}+ \Ebf \left[ \one_{\{D_j-\xi_{jk} > 0\}} \one_{\{D_{k}-\xi_{jk} > 0\}} \frac{\xi_{ji}\xi_{ki}}{(D_{j}-\xi_{jk})^2(D_{k}-\xi_{jk})} \right] 
	\end{align*}
	\begin{align*}
		& = \Ebf \left[ \frac{\xi_{ji}}{(D_{j}-\xi_{jk}+1)^2} \right] \Ebf \left[ \frac{\xi_{ki}}{D_{k}-\xi_{jk}+1} \right]\\
		&\hspace{3cm}+ \Ebf \left[ \one_{\{D_j-\xi_{jk} > 0\}} \frac{\xi_{ji}}{(D_{j}-\xi_{jk})^2} \right] \Ebf \left[ \one_{\{D_{k}-\xi_{jk} > 0\}} \frac{\xi_{ki}}{D_{k}-\xi_{jk}} \right],
	\end{align*}
	where the last equality follows from independence between $(\xi_{ji},D_j-\xi_{jk})$ and $(\xi_{ki},D_k-\xi_{jk})$.
	Using exchangeability and Lemma \ref{lem:prep_1} we have
	\begin{align*}
		\Ebf \left[ \frac{\xi_{ji}}{(D_{j}-\xi_{jk}+1)^2} \right] & = \frac{1}{N-2} \sum_{l \in [N], l \ne j,k} \Ebf \left[ \frac{\xi_{jl}}{(D_{j}-\xi_{jk}+1)^2} \right]\\
		& = \frac{1}{N-2} \Ebf \left[ \frac{D_j-\xi_{jk}}{(D_{j}-\xi_{jk}+1)^2} \right] \\
		& \le \frac{1}{N-2} \Ebf \left[ \frac{1}{D_{j}-\xi_{jk}+1} \right] \le \frac{1}{(N-2)(N-1)p_N}.
	\end{align*}
	Similarly one can verify that
	\begin{align*}
		 &\Ebf \left[ \frac{\xi_{ki}}{D_{k}-\xi_{jk}+1} \right] \le \frac{1}{N-2}, \quad
		 \Ebf \left[ \one_{\{D_{k}-\xi_{jk} > 0\}} \frac{\xi_{ki}}{D_{k}-\xi_{jk}} \right] \le \frac{1}{N-2} \\
		 &\Ebf \left[ \one_{\{D_j-\xi_{jk} > 0\}} \frac{\xi_{ji}}{(D_{j}-\xi_{jk})^2} \right] \le \frac{4}{(N-2)(N-1)p_N}.
	\end{align*}
	Combining these gives us 
	\begin{equation*}
		\Ebf \left[ \one_{\{D_{j} \ge d-1\}} \one_{\{D_{k} \ge d-1\}} \frac{\xi_{ji}\xi_{ki}}{D_{j}^2D_{k}} \right] \le \frac{5}{(N-2)^2(N-1)p_N}, \mbox{ when } j \ne k.
	\end{equation*}
	Also note that the summation in \eqref{eq:prep_new_5} when $j=k$ is 
	\begin{align*}
		\sum_{j=1,j \ne i}^N \Ebf \left[ \one_{\{D_{j} \ge d-1\}} \frac{\xi_{ji}}{D_{j}^3} \right] &= \sum_{j=1,j \ne i}^N \Ebf \left[ \one_{\{D_{i} \ge d-1\}} \frac{\xi_{ij}}{D_{i}^3} \right]\\
		& = \Ebf \left[ \one_{\{D_{i} \ge d-1\}} \frac{1}{D_i^2} \right] \le \frac{4}{(Np_N)^2},
	\end{align*}
	where the first equality uses exchangeability and the inequality uses Lemma \ref{lem:prep_1}.
	Combining these two estimates with \eqref{eq:prep_new_5} gives
	\begin{equation*}
		V_s^A \le \kappa_5  \frac{N^2}{(N-2)^2(N-1)p_N} + \kappa_5 \frac{1}{(Np_N)^2}\le \frac{\kappa_6}{Np_N} + \frac{\kappa_6}{(Np_N)^2}
	\end{equation*}
	for some $\kappa_5, \kappa_6 \in (0,\infty)$.
	This completes the proof of Lemma~\ref{lem:prep_new}.
\end{proof}

Finally we complete the proof of Lemma \ref{lem:momentbd}.
\begin{proof}[Proof of Lemma \ref{lem:momentbd}]
	As before, we will omit the superscript in $\xi_{ij}$'s and $D_i$'s for notational convenience.
	Fix $i\in\N$.
	From \eqref{eq:X_i_n} and \eqref{eq:limit-tagged}, using Cauchy--Schwarz and Doob's inequalities we have for any fixed $t\in [0,T]$	
	\begin{align}
		\Ebf \left\|X_i^N - X_i\right\|_{*,t}^4  \le 
		\kappa_1 \int_0^t \Ebf |X_i^N(s) - X_i(s)|^4 \, ds &+ \kappa_1 \int_0^t \Ebf |C_i^N(s)-C_i(s)|^2 \, ds\notag \\
		&\hspace{-1cm}+ \kappa_1 \int_0^t \Ebf |C_i^N(s)-C_i(s)|^4 \, ds. \label{eq:momentbd1}
	\end{align}
	Recall the definition of $C_i^N(s)$ and $C_i(s)$ from~\eqref{eq:cn} and \eqref{eq:limit-tagged}.
	From the bound $\|b\|_\infty \le 1$ and \eqref{eq:bbar}, for $s \in [0,T]$ we have $|C_i(s)| \le d$ and
	\begin{align}
		\Ebf |C_i^N(s)|^4 & \le \Ebf \left| \one_{\{D_i<d-1\}}(D_i+1) + 1 + (d-1) \sum_{j_2 \in [N], j_2 \ne i} \one_{\{D_{j_2} \ge d-1\}} \frac{\xi_{j_2i}}{D_{j_2}} \right. \notag \\
		& \left. \hspace{3cm} + \sum_{j_2 \in [N], j_2 \ne i} \one_{\{D_{j_2} < d-1\}} \xi_{ij_2}(D_i+1) \right|^4 \notag \\
		& \le \kappa_2 + \kappa_2 \Ebf \Big[ \sum_{j_2 \in [N], j_2 \ne i} \one_{\{D_{j_2} \ge d-1\}} \frac{\xi_{j_2i}}{D_{j_2}} \Big]^4 \\
		&\hspace{3cm}+ \kappa_2 \Ebf\Big[\sum_{j_2 \in [N], j_2 \ne i} \one_{\{D_{j_2} < d-1\}} \xi_{ij_2}(D_i+1)\Big]^4. \label{eq:CN_bd}
	\end{align}
	Here the second term on the right hand side can be written as
	\begin{align*}
		& \kappa_2 \Ebf \Big[ \sum_{j_2 \in [N], j_2 \ne i} \one_{\{D_{j_2} \ge d-1, D_i > 0\}} \frac{D_i}{D_{j_2}} \frac{\xi_{j_2i}}{D_i} \Big]^4 \\
		& \le \kappa_2 \Ebf \Big[ \sum_{j_2 \in [N], j_2 \ne i} \one_{\{D_{j_2} \ge d-1, D_i > 0\}} \left(\frac{D_i}{D_{j_2}}\right)^4 \frac{\xi_{j_2i}}{D_i} \Big]\\
		&\hspace{6cm} \Big[ \sum_{j_2 \in [N], j_2 \ne i} \one_{\{D_{j_2} \ge d-1, D_i > 0\}} \frac{\xi_{j_2i}}{D_i} \Big]^3 \\
		& \le \kappa_2 \Ebf \sum_{j_2 \in [N], j_2 \ne i} \one_{\{D_{j_2} \ge d-1\}} \frac{D_i^3 \xi_{j_2i}}{D_{j_2}^4}\\
		& = \kappa_2 \sum_{j_2 \in [N], j_2 \ne i} \Ebf \Big[ \one_{\{D_{j_2}-\xi_{j_2i}+1 \ge d-1\}} \frac{(D_i-\xi_{j_2i}+1)^3 }{(D_{j_2}-\xi_{j_2i}+1)^4} \Big] p_N \\
		& = \kappa_2 \sum_{j_2 \in [N], j_2 \ne i} \Ebf \Big[ \one_{\{D_{j_2}-\xi_{j_2i}+1 \ge d-1\}} \frac{1}{(D_{j_2}-\xi_{j_2i}+1)^4} \Big] \Ebf \Big[ D_i-\xi_{j_2i}+1 \Big]^3 p_N \\
		& \le \kappa_3 (N-1) \frac{1}{(N-1)^4p_N^4}(Np_N+1)^3p_N \le \kappa_4,
	\end{align*}
	where the first inequality uses Holder's inequality, the first equality follows by conditioning on $\xi_{j_2i}=1$, the second equality follows from independence, and the third inequality uses Lemma \ref{lem:prep_1} and moment estimates of binomial random variables.
	Following the similar argument, we can write the last term in \eqref{eq:CN_bd} as
	\begin{align*}
		& \kappa_2 \Ebf\Big[\sum_{j_2 \in [N], j_2 \ne i} \one_{\{D_{j_2} < d-1, D_i > 0\}} D_i(D_i+1) \frac{\xi_{ij_2}}{D_i} \Big]^4 \notag \\
		& \le \kappa_2 \Ebf \Big[ \sum_{j_2 \in [N], j_2 \ne i} \one_{\{D_{j_2} < d-1, D_i > 0\}} D_i^4(D_i+1)^4 \frac{\xi_{ij_2}}{D_i} \Big] \\
		&\hspace{6cm}\Big[ \sum_{j_2 \in [N], j_2 \ne i} \one_{\{D_{j_2} \ge d-1, D_i > 0\}} \frac{\xi_{j_2i}}{D_i} \Big]^3 
	\end{align*}
	\begin{align*}		
		& \le \kappa_2 \Ebf \sum_{j_2 \in [N], j_2 \ne i} \one_{\{D_{j_2} < d-1\}} D_i^3(D_i+1)^4\xi_{ij_2} \notag \\
		& = \kappa_2 \sum_{j_2 \in [N], j_2 \ne i} \Ebf \Big[\one_{\{D_{j_2}-\xi_{ij_2}+1 < d-1\}} (D_i-\xi_{ij_2}+1)^3 (D_i-\xi_{ij_2}+2)^4 \Big]p_N \notag \\
		& = \kappa_2 \sum_{j_2 \in [N], j_2 \ne i} \Ebf \Big[\one_{\{D_{j_2}-\xi_{ij_2}+1 < d-1\}} \Big] \Ebf \Big[ (D_i-\xi_{ij_2}+1)^3 (D_i-\xi_{ij_2}+2)^4 \Big] p_N \notag \\
		& \le \kappa_5 (N-1) \Pbf (D_{i}< d) (Np_N+1)^7 p_N. 
	\end{align*}
	Combining above three estimates with \eqref{eq:CN_bd_2} and using Condition \ref{cond:errg1}, we have
		$\Ebf |C_i^N(s)|^4 \le \kappa_6$.
	It then follows from \eqref{eq:momentbd1} that 
	\begin{equation*}
		\Ebf \left\|X_i^N - X_i\right\|_{*,t}^4 \le \kappa_7 \int_0^t \Ebf \|X_i^N - X_i\|_{*,s}^4 \, ds + \kappa_7.
	\end{equation*}
	The result then follows from Gronwall's inequality.
\end{proof}

\section{Conclusion}\label{sec:conclusion-aap2}
We have considered the JSQ($d$) policy in large-scale systems where the servers communicate with their neighbors and the neighborhood relationships are described in terms of a suitable graph.
We have developed sufficient criteria for arbitrary graph sequences so that asymptotically the evolution of the occupancy process on any finite time interval is indistinguishable from that for the case when the graph is a  clique.
We have also  considered  sequences of Erd\H{o}s-R\'enyi random graphs and established sufficient criteria in terms of the growth rates of the average degree that ensure the annealed and quenched  limit of the occupancy process on any finite time interval to coincide with that in the clique.

The steady-state behavior of the occupancy measure process associated with
 the above graph sequences is an important and challenging  open question.
Steady-state properties of the JSQ($d$) scheme has been well studied in the case of a clique.
For example, in~\cite{Mitzenmacher96, Mitzenmacher01} it is shown that $\pi^N$, the stationary measure of the occupancy process of the $N$-th system, converges in distribution to $\delta_{\qq^*}$, where $\qq^*$ is the unique fixed point of the limiting deterministic dynamical system $\qq(\cdot)$.
Roughly speaking such a result says that the limits  $t\to\infty$ and $N\to\infty$ can be interchanged.
Based on Theorems \ref{th:deterministic-d}, \ref{th:tagged-d}, \ref{thm:npn_rate}, and \ref{thm:npn_rate_quench}, it is natural to conjecture that a similar interchangeability also holds for more general graphs considered in this chapter.
However, the setting here is more complicated, in particular, the occupancy process is not Markov any more.
One may conjecture that with $\pi^N$ replaced by the time asymptotic limit of the law of the occupancy process, the convergence $\pi^N \to \delta_{\qq^*}$ still holds. However,
currently even the existence of such a time asymptotic limit is not clear.


%% file: TT_Summary.tex
{
A fundamental challenge in large-scale networked systems viz.~data centers and cloud networks is to distribute tasks to a pool of servers, using minimal instantaneous state information, while providing excellent delay performance. 
In this thesis we design and analyze load balancing algorithms that aim to achieve a highly efficient distribution of tasks, optimize server utilization, and minimize communication overhead. 
A canonical model of these systems consists of $N$ parallel single-server queues with unit-exponential service rates and a single dispatcher where tasks arrive at rate $\lambda(N)$. 
Motivated by the immense operational scale of data centers and cloud networks, we analyze this model in a large-capacity regime, where both $N$ and $\lambda(N)$ grow large in proportion. 

In Chapter 1 we start with an introduction to the various load balancing algorithms studied in the existing literature. 
We further provide an overview of our contributions along with a high-level description of the stochastic coupling techniques that play an instrumental role in establishing the results.

In Chapter 2 we consider a class of schemes where the dispatcher assigns each arriving task to a server with the shortest queue among $d(N)$ randomly selected servers $(1 \leq d(N) \leq N)$. 
This load balancing strategy is referred to as a JSQ($d(N)$) scheme, marking that it subsumes the Join-the-Shortest Queue (JSQ) policy as a crucial special case for $d(N) = N$. 
The JSQ policy exhibits several strong optimality properties, but at the cost of a prohibitively high communication burden of $N$ per task. 
In contrast, a random assignment policy ($d(N)=1$) has no communication overhead but has much worse delay performance. Leveraging a novel stochastic coupling construction, we first show that asymptotic optimality can be achieved on fluid and diffusion scale while reducing the overhead by nearly a factor O$(N)$ and O$(\sqrt{N}/\log N)$, respectively, as the number of servers $N$ becomes large.  
 
In Chapter 3 we consider the Join-the-Idle queue (JIQ) strategy, which assigns the incoming tasks to idle servers, if any, and to a server selected uniformly at random otherwise. 
Exploiting the memory at the dispatcher as a further dimension, JIQ can be implemented in a token-based manner with only O$(1)$ communication overhead per task. We establish that the JIQ strategy achieves diffusion-level optimality, and thus comparing with the results presented in Chapter 2, the JIQ strategy outperforms the JSQ$(d(N))$ schemes in terms of delay and communication. 

In Chapter 4 we analyze the steady-state diffusion process that comes as the weak limit of the appropriately scaled occupancy states of a system under the JSQ policy as the number of servers grows large.
From a methodological point of view, the diffusion process under consideration goes beyond the state-of-the-art techniques in the study of the steady state of diffusion processes.
Exploiting a technique involving the theory of regenerative processes, we establish precise tail asymptotics of the stationary distribution and scaling of extrema of the process on large time intervals. Our results imply that the asymptotic steady-state scaled number of servers with queue length two or larger exhibits an exponential tail, whereas that for the number of idle servers turns out to be Gaussian.

The asymptotic optimality results of Chapter 2 are extended in Chapter 5 to an infinite-server scenario where the single-server queues are replaced by server pools. 
As it turns out, due the intrinsic difference in the dynamics, a fundamentally different coupling argument is required to establish asymptotic optimality results.
 
Along with the delay-communication trade-off, in Chapter 6 we further consider the issue of energy consumption, which has increasingly become a concern in large-scale data centers in recent years. Specifically, we provide the first token-based joint auto-scaling and load balancing algorithm that has a distributed scalable implementation and yet achieves asymptotic optimality on fluid scale in all three aspects of the delay-communication-energy trade-off.

Chapter 7 extends the asymptotic optimality result of Chapter 6 in the case when the servers have an infinite buffer capacity.
As it turns out, establishing asymptotic optimality results in this case is fundamentally more challenging, since the system stability under the usual subcritical load assumption is not automatic.
We develop a novel method to prove this stability for any subcritical load, and establish convergence of steady-state distributions to the optimal one, as the system size grows large.
 
The analysis becomes much more complicated when the servers are interconnected by some graph topology. The topology can arise, for example, due to selective placement of data files. Due to the lack of a tractable Markovian state description, analysis of such systems had remained intractable for topologies sparser than a clique. We exploit a novel stochastic coupling framework in Chapter 8 and an asymptotic decoupling method in Chapter 9 to establish that for various load balancing algorithms the asymptotic behavior of a clique can be achieved in the large-capacity regime while dramatically reducing the number of connections.}

%% file: TT_cv.tex
Debankur Mukherjee was born on  16 September 1991 in Hooghly, India.
After finishing high school in 2009 at Hooghly Collegiate School in Chinsurah, India, he studied B.Sc.~with statistics major at Ramakrishna Mission Residential College, Narendrapur (affiliated with the University of Calcutta) in Kolkata, India. 
In 2012 he joined the Indian Statistical Institute, where he obtained his master's degree with mathematical statistics and probability specialization. 
In August 2014 he started a Ph.D.~project at Eindhoven University of Technology (TU/e) in the Stochastics Operations Research (SOR) group under the supervision of Sem Borst and Johan van Leeuwaarden. 

His Ph.D.~research spans the areas of probability theory and stochastic networks, at the interface of stochastic processes and computer science, with applications in queueing theory, performance analysis, random graphs, and randomized algorithms. 
His key contributions towards the field of load balancing and scheduling in large-scale stochastic networks include developing novel coupling techniques and establishing non-classical stochastic process limit theorems to study the delay-communi-cation-energy trade-off in large networks viz.~data centers and cloud networks. 
Debankur also introduced a stochastic comparison framework to study mean-field limits of processes on networks and examined the impact of the network topology on the performance of load balancing schemes in large-scale systems, for which he received the Best Student Paper Award at ACM SIGMETRICS 2018. 
The results of this research are collected in this thesis.

In July 2018 Debankur joined Brown University as a Prager assistant professor in the Division of Applied Mathematics.
In August 2019 he will join Georgia Institute of Technology as a tenure-track assistant professor in the H.~Milton Stewart School of Industrial and Systems Engineering.